%% file: main.tex
\documentclass[11pt]{article}
\usepackage[T1]{fontenc}
\usepackage{lmodern,amsmath,amssymb,amsthm,booktabs,longtable,array,microtype,enumitem}
\usepackage[margin=0.85in]{geometry}
\usepackage[hidelinks]{hyperref}
\hypersetup{pdftitle={A complete symmetric rank-two modularity classification, obtained through a sharp uniform degree bound},pdfauthor={Cetin Hakimoglu-Brown}}
\newcommand{\C}{\mathbb C}\newcommand{\R}{\mathbb R}
\newcommand{\Q}{\mathbb Q}\newcommand{\Z}{\mathbb Z}\newcommand{\F}{\mathbb F}
\newcommand{\Newt}{\operatorname{Newt}}\newcommand{\EE}{\mathbb E}\newcommand{\ord}{\operatorname{ord}}
\newcommand{\rank}{\operatorname{rank}}
\newcommand{\N}{\operatorname{N}}

\newcommand{\dd}{\delta}

\newtheorem{theorem}{Theorem}[section]
\newtheorem{lemma}[theorem]{Lemma}\newtheorem{proposition}[theorem]{Proposition}
\newtheorem{corollary}[theorem]{Corollary}
\theoremstyle{definition}\newtheorem{definition}[theorem]{Definition}
\theoremstyle{remark}\newtheorem{remark}[theorem]{Remark}
\newcolumntype{P}[1]{>{\raggedright\arraybackslash}p{#1}}
\numberwithin{equation}{section}
\makeatletter
\renewcommand*\l@section[2]{\@dottedtocline{1}{0em}{2.6em}{#1}{#2}}
\makeatother
\title{A complete symmetric rank-two modularity classification,\\obtained through a sharp uniform degree bound}
\author{Cetin Hakimoglu-Brown}
\date{8 November 2026}
\begin{document}
\maketitle
\begin{abstract}
We classify modular Nahm sums with a rational symmetric positive-definite
two-by-two matrix and rational shifts. Up to coordinate exchange, every
modular datum belongs to three explicit rational-point families, is a
product of two modular rank-one sums, or is one of thirty-four exceptional
triples on twelve matrices. The list consists of Zagier's rank-two
examples, the Vlasenko--Zwegers additions and split products; no further
data occur. The proof first bounds the joint degree of the distinguished
Nahm point by four, using an initial degree bound of twenty-four and
exact local exclusions at the prime 7. Explicit height and exponent bounds
then give exhaustive arithmetic classifications in degrees one through
four. Exact polynomial elimination determines all rational shifts, and
published identities of q-series establish modularity of the survivors.
The degree bound is sharp. The finite steps are accompanied by exact
certificates, verification programs and independently implemented checks.
\end{abstract}
\clearpage
{\small\tableofcontents}
\clearpage
\input{sections/introduction.tex}
\clearpage
\input{sections/degree_bound.tex}
\clearpage
\part{Arithmetic preliminaries}\label{part:arithmetic}
\input{sections/arithmetic_overview.tex}
\input{sections/arithmetic_preliminaries.tex}

\input{sections/quartic_entry.tex}
\clearpage
\part{Modular data of degree at most two}\label{part:classification}
\input{sections/coefficient_conditions.tex}

\input{sections/degree_two_classification.tex}

\clearpage
\part{The cubic classification}\label{part:cubic}
\input{sections/cubic_unit.tex}

\input{sections/cubic_remaining.tex}
\clearpage
\part{The quartic classification}\label{part:quartic}
\input{sections/quartic.tex}

\appendix
\input{sections/degree_four_consequences.tex}
\input{sections/small_certificates.tex}
\input{sections/global_certificates.tex}
\input{sections/rational_boundary.tex}

\input{sections/verification.tex}
\input{sections/disclosures.tex}
\clearpage
\begingroup\small

\input{bibliography.tex}
\endgroup
\end{document}

%% file: sections/introduction.tex
\section{The classification}\label{sec:introduction}
In his 2007 survey, when introducing the asymptotic proof of the
rank-one classification, Zagier wrote that ``it is not clear whether
one could use the method to complete the classification even in the
next case $r=2$'' \cite[p.~50]{Zagier}. The question is which rational
symmetric positive-definite rank-two Nahm sums~\cite{Nahm} are modular. This paper
addresses completeness: Theorem~\ref{thm:classification} gives the
known rational-point families, split products and thirty-four
exceptional triples, and excludes all further data under these
hypotheses.

Let $A$ be a rational symmetric positive-definite $2\times2$ matrix,
$\boldsymbol B\in\Q^2$, and $c_0\in\Q$. We consider
\begin{equation}\label{eq:nahm-sum}
 f_{A,\boldsymbol B,c_0}(q)=
 \sum_{n\in\Z_{\ge0}^2}
 \frac{q^{\frac12 n^tAn+\boldsymbol B^tn+c_0}}
 {(q)_{n_1}(q)_{n_2}},\qquad
 (q)_n=\prod_{j=1}^{n}(1-q^j).
\end{equation}
We call the datum $(A,\boldsymbol B,c_0)$ modular if, with
$q=e^{2\pi i z}$, this sum is a modular function for a finite-index
subgroup of $\mathrm{SL}_2(\Z)$, meromorphic at the cusps; a finite-order
multiplier is allowed. Its distinguished Nahm point is the unique
$(x,y)\in(0,1)^2$ satisfying
\begin{equation}\label{eq:intro-nahm}
 1-x=x^a y^b,\qquad 1-y=x^b y^d,
 \qquad A=\begin{pmatrix}a&b\\b&d\end{pmatrix}.
\end{equation}
The powers in these equations are positive real powers. The degree of
the point always means the joint degree $[\Q(x,y):\Q]$.

The rank-one problem was solved by Terhoeven and Zagier and has
seven triples \cite{NRT,Terhoeven,Zagier} (see also \cite[Section~1]{VZ});
Zagier's table of rank-two modular examples lists one
family and ten exceptional matrices, conventionally numbered as eleven
examples \cite[p.~47, Table~2]{Zagier}. Vlasenko and Zwegers
extended the asymptotic method to higher rank, found additional modular
data and counterexamples to the original all-solutions formulation of
Nahm's conjecture \cite{VZ}; Calegari, Garoufalidis and Zagier later
proved that modularity forces the Bloch-group element of the
distinguished Nahm point to be torsion \cite[Theorem~1.8]{CGZ}.
Wang established identities for ten of Zagier's eleven examples,
including the cases already known, and Cao, Rosengren and Wang proved
the remaining example \cite{Wang,CRW}. Huang and Lee studied the
positive-definite symmetric integer matrices for which every complex
solution of the Nahm equations is real~\cite{HL}. The modular identities establish
modularity of the known data; the problem here is to prove that the
list, with the Vlasenko--Zwegers additions and split products, is
exhaustive.

The main new step is the uniform bound of four for the joint degree of
the distinguished Nahm point.
Zagier's rank-two investigation and the higher-rank asymptotic method of
Vlasenko and Zwegers provide polynomial necessary conditions for modularity.
In the logarithmic normalization of Lemma~\ref{deg:lem:modular-bridge},
the first coefficient is fixed by the scalar shift and all higher
coefficients vanish. We use the geometry of these equations, followed by
local analysis at the prime 7, to bound the joint field degree first by
24 and then by four, with no initial height bound on the matrix or shifts.
The degree bound changes the remaining problem into an arithmetic
classification: unit and nonunit equations in fields of degree at most
four, explicit height and exponent estimates, exhaustive finite
enumeration, and exact elimination of the rational shifts. The torsion
condition of \cite{CGZ} supplies an additional necessary condition in
these arithmetic exclusions; it is not used as a substitute for the
degree-bound argument.

For reference, write the complete rank-one list as
\begin{equation}\label{eq:intro-rankone}
\begin{split}
 \mathcal R=\{& (2,0,-\tfrac1{60}),\ (2,1,\tfrac{11}{60}),\
 (1,0,-\tfrac1{48}),\ (1,\tfrac12,\tfrac1{24}),\
 (1,-\tfrac12,\tfrac1{24}),\\
 & (\tfrac12,0,-\tfrac1{40}),\
 (\tfrac12,\tfrac12,\tfrac1{40})\}.
\end{split}
\end{equation}
The entries of a rank-one triple are its quadratic coefficient, linear
shift and scalar shift. We identify rank-two data under coordinate
exchange, which sends $(A,\boldsymbol B,c_0)$ to
$(JAJ,J\boldsymbol B,c_0)$, where
$J=\left(\begin{smallmatrix}0&1\\1&0\end{smallmatrix}\right)$.

\begin{theorem}[Complete symmetric rank-two classification]
\label{thm:classification}\label{thm:finite-main}
Let $A\in M_2(\Q)$ be symmetric and positive definite,
$\boldsymbol B\in\Q^2$, and $c_0\in\Q$. Then
$f_{A,\boldsymbol B,c_0}$ is modular if and only if, up to coordinate
exchange, the datum belongs to one of the following classes.
\begin{enumerate}[label=(\roman*),leftmargin=*]
\item \emph{The rational-point families.} For $a\in\Q$, $a>1/2$,
\[
 A=\begin{pmatrix}a&1-a\\1-a&a\end{pmatrix},
\]
and the vector and scalar shifts are one of
\begin{align}
 \boldsymbol B&=(\beta,-\beta),&
 c_0&=\frac{\beta^2}{2a}-\frac1{24},
 &&\beta\in\Q;\label{eq:main-family1}\\
 \boldsymbol B&=(-\tfrac12,-\tfrac12),&
 c_0&=\frac1{8a}-\frac1{24};\label{eq:main-family2}\\
 \boldsymbol B&=(\tfrac a2,1-\tfrac a2),&
 c_0&=\frac a8-\frac1{24}.\label{eq:main-family3}
\end{align}
\item \emph{Split products.} There are
$(a_i,\beta_i,c_i)\in\mathcal R$, $i=1,2$, such that
\[
 A=\operatorname{diag}(a_1,a_2),\qquad
 \boldsymbol B=(\beta_1,\beta_2),\qquad c_0=c_1+c_2.
\]
\item \emph{Exceptional data.} The triple is one of the $34$ triples in
Table~\ref{tab:complete-list}.
\end{enumerate}
The distinguished point of every modular datum has joint degree at
most four, and all four degrees $1,2,3,4$ occur.
\end{theorem}

The word ``finite'' concerns the exceptional list and the number of
family types. The first class has infinitely many rational matrices and
shifts. Its point is always $(1/2,1/2)$. The exceptional list comprises
$18$ quadratic, $12$ cubic and four quartic triples on twelve matrix
classes. In the table, $Z_j$ denotes Zagier's Example~$j$ and $V_1,V_2$
denote the two additional matrices of Vlasenko and Zwegers. A pair of
shifts exchanged by a symmetry of the displayed matrix represents one
entry. The split products and the rational-point families are understood
as a union; their overlap is harmless.

\input{sections/classification_table.tex}
\clearpage

\begin{proof}[Proof of Theorem~\ref{thm:classification}]
Theorem~\ref{deg:thm:main} bounds the joint degree of the distinguished
point of a modular datum by four. The proof includes the symmetric
degree-$24$ argument, the higher asymptotic obstructions and the complete
local degree count.

If the degree is one or two, Theorem~\ref{cls:thm:deg2} gives exactly
the first two classes and the quadratic rows of
Table~\ref{tab:complete-list}. Its proof combines the complete arithmetic
classification with exact elimination of arbitrary rational shifts.
If the degree is three, Theorem~\ref{cr:classification} gives exactly
$Z_6,Z_7,Z_8,Z_9$ with their twelve displayed shifts. It includes the
independent unit case, both canonical nonunit exponent cases, and the
dependent and diagonal exclusions. If the degree is four,
Theorem~\ref{qt:classification} gives exactly $Z_4,Z_5$ with their four
displayed shifts. The canonical entry lemma, the nonunit and dependent
exclusions, and the unit arguments for both possible signatures account
for every quartic case. No matrix or shift-height bound is imposed on
the initial datum.

Conversely, the first class is modular by the identities of
\cite[Table~2]{VZ}, including Zagier's theta identity. The second class
is a product of modular rank-one sums. The quadratic exceptional data
are modular by the results cited in Theorem~\ref{cls:thm:deg2}; the cubic
and quartic identities are matched to their exact shifts in
Theorems~\ref{cr:classification} and~\ref{qt:classification}. In
particular, sufficiency follows from identities of $q$-series, not from
the vanishing of finitely many asymptotic coefficients. Examples
$Z_1,Z_2,Z_6,Z_4$ realize degrees $1,2,3,4$, respectively.
\end{proof}

\begin{corollary}[Sharp degree bound]\label{cls:cor:sharp}
The bound $[\Q(x,y):\Q]\le4$ is sharp. For
$A=\left(\begin{smallmatrix}4&1\\1&1\end{smallmatrix}\right)$, put
$\rho=(\sqrt5-1)/2$. The distinguished point is
$(\sqrt\rho,1/(1+\sqrt\rho))$ and its field is the quartic field
$\Q(\sqrt\rho)$.
\end{corollary}

\subsection{The exclusion argument}
\paragraph{Why the proof is long.}
Most of the length comes from the degree reduction in
Part~\ref{part:degree}. In rank one, the Nahm equation has one unknown
and the asymptotic conditions lead to explicit univariate eliminations.
In rank two, the matrix varies in three rational parameters, and the
corresponding equations define plane curves whose common components
must be controlled uniformly. Genuine infinite families occur,
including Zagier's family with distinguished point $(1/2,1/2)$, so
one cannot simply assume that every intersection is finite. We first
obtain the uniform joint-degree bound twenty-four, then use a local
analysis at the prime $7$, divided into twenty-nine terminal cases,
to reduce the bound to four. The subsequent arithmetic classification
requires separate unit, nonunit and dependent-coordinate arguments in
each remaining degree. The length reflects these reductions and the
need to specify the finite computations together with their coverage
arguments. A direct proof that modularity forces joint degree at most
four would dispense with the local degree reduction, while leaving
the low-degree arithmetic classification to be carried out. No such
shorter proof is known to the author.

The proof has three stages. First, modularity constrains the joint
degree. The initial degree-$24$ bound is obtained from the geometry of
the first four logarithmic asymptotic coefficients. Exact local
arguments at the prime $7$ then lower the bound to four. The final assembly
has $29$ terminal cases; Appendix~\ref{deg:sec:supplement} specifies the
proof dependencies and the certificates used to close them.

Second, the degree restriction makes the arithmetic branches tractable.
The symmetric equations give the wedge relation
$[x]\wedge[1-x]+[y]\wedge[1-y]=0$. Dependent coordinates reduce to a
short scalar list and explicit matrix families. Independent nonunit
coordinates enter two canonical polynomial families by finite-place
valuation constraints. For independent units, the common unit lattice,
index parity and real or complex embeddings give practical height
limits. These bounds precede, and justify, the finite enumerations.

Third, the surviving points and matrices are subjected to exact
coefficient equations in the rational shifts. The elimination is over
$\Q$ or an explicitly defined number field and uses no search bound on
$\boldsymbol B$ or $c_0$. Certificates express the required eliminants
or the constant $1$ as polynomial combinations of the necessary
equations. Independent coefficient implementations and exact interval
checks accompany the calculations.

\subsection{Scope and relation to Nahm's conjecture}\label{sec:scope}
Theorem~\ref{thm:classification} resolves the modularity classification
problem for symmetric rank-two Nahm sums under the hypotheses
of \eqref{eq:nahm-sum}. It proves completeness of the known families and
exceptional modular data, including the Vlasenko--Zwegers additions to
Zagier's list.

The degree bound starts from modularity. The implication from modularity
to the Bloch-group condition in \cite{CGZ} is used in the exclusions; its
converse is not assumed. Thus this classification does not assert a
classification of all algebraic two-term torsion relations, nor does it
assert the original all-solutions equivalence in Nahm's conjecture,
which has counterexamples \cite{VZ}. In the notation of~\cite{D24ref},
the symmetrizer ratio is $m=1$ throughout; for Nahm sums with a
nontrivial symmetrizer see~\cite{MizunoSym}.

\subsection{Organization and notation}
Part~\ref{part:degree} contains the degree-four proof, including its
initial bound. Part~\ref{part:arithmetic} gives the arithmetic
preliminaries, scalar and quadratic classifications and the canonical
entry argument. Part~\ref{part:classification} treats coefficient
conditions and all modular data of degree at most two.
Parts~\ref{part:cubic} and~\ref{part:quartic} contain the complete cubic
and quartic proofs. Appendix~\ref{app:degree-four-consequences} discusses the role
of the degree-four cutoff, including the quintic unit and regulator
comparison; this discussion is not used in the classification proof.
The other appendices specify computational steps, dependencies and
independent checks.

In the degree-reduction part we use
$A=\left(\begin{smallmatrix}a&b\\b&d\end{smallmatrix}\right)$ and a
linear shift $(B,C)^t$. In the arithmetic parts the entries may be
denoted by $\alpha,\beta,\gamma$, with $\beta$ the coupling. The linear
shift vector is then $\boldsymbol B=(B_1,B_2)^t$. Scalar shifts are
denoted by $c_0$, and $C_0=\pi^2/6$ is reserved for the dilogarithm
normalization. Polynomial variables and local valuations are introduced
where they are used. For quantities integral at the finite place under
consideration, a bar denotes reduction modulo its maximal ideal, applied
coefficientwise to polynomials, vectors and matrices. In the local
analysis of Part~\ref{part:degree}, this place lies above $7$, as in
$\bar\tau\in\F_7$. In Section~\ref{d24:section}, however, bars denote
bihomogenization on $\mathbb P^1\times\mathbb P^1$.
As usual, $\overline{\Q}$ and $\overline{\F}_q$ denote algebraic closures.

%% file: sections/classification_table.tex
\begingroup\small
\setlength{\tabcolsep}{5pt}
\renewcommand{\arraystretch}{1.15}
\begin{longtable}{@{}c c c P{0.49\textwidth}@{}}
\caption{Complete exceptional list, up to coordinate exchange.
The column $d$ is $[\Q(x,y):\Q]$. Each displayed pair
$(\boldsymbol B;c_0)$ is one triple.}\label{tab:complete-list}\\
\toprule
Label&Matrix $A$&$d$&$(\boldsymbol B;c_0)$\\\midrule
\endfirsthead
\multicolumn{4}{l}{\textit{Table~\ref{tab:complete-list} continued}}\\
\toprule
Label&Matrix $A$&$d$&$(\boldsymbol B;c_0)$\\\midrule
\endhead
\midrule\multicolumn{4}{r}{\textit{Continued on the next page}}\\
\endfoot
\bottomrule
\endlastfoot
$V_1$&$\begin{pmatrix}3/2&1/2\\1/2&3/2\end{pmatrix}$&2&
$(\frac14,-\frac14);-\frac1{120}$\newline
$(\frac14,\frac34);\frac{11}{120}$\\[1pt]
$V_2$&$\begin{pmatrix}3/4&-1/4\\-1/4&3/4\end{pmatrix}$&2&
$(\frac14,-\frac14);-\frac1{80}$\newline
$(\frac12,0);\frac1{80}$\\[1pt]
$Z_2$&$\begin{pmatrix}2&1\\1&1\end{pmatrix}$&2&
$(0,0);-\frac1{32}$\quad $(0,\frac12);0$\newline
$(1,\frac12);\frac18$\quad $(-1,\frac12);\frac18$\newline
$(1,1);\frac7{32}$\\[1pt]
$Z_3$&$\begin{pmatrix}1&-1\\-1&2\end{pmatrix}$&2&
$(0,0);-\frac5{96}$\quad $(\frac12,0);\frac1{24}$\newline
$(-\frac12,1);\frac16$\quad $(0,1);\frac{19}{96}$\newline
$(-\frac32,2);\frac{25}{24}$\\[1pt]
$Z_{10}$&$\begin{pmatrix}4/3&2/3\\2/3&4/3\end{pmatrix}$&2&
$(0,0);-\frac1{30}$\newline
$(-\frac13,-\frac23);\frac1{30}$\\[1pt]
$Z_{11}$&$\begin{pmatrix}1&-1/2\\-1/2&1\end{pmatrix}$&2&
$(0,0);-\frac1{20}$\newline
$(-\frac12,0);\frac1{20}$\\\midrule
$Z_6$&$\begin{pmatrix}4&2\\2&2\end{pmatrix}$&3&
$(0,0);-\frac1{42}$\quad $(1,0);\frac5{42}$\newline
$(2,1);\frac{17}{42}$\\[1pt]
$Z_7$&$\begin{pmatrix}1/2&-1/2\\-1/2&1\end{pmatrix}$&3&
$(0,0);-\frac5{84}$\quad $(\frac12,0);\frac1{84}$\newline
$(\frac12,-\frac12);\frac1{21}$\\[1pt]
$Z_8$&$\begin{pmatrix}3/2&1\\1&2\end{pmatrix}$&3&
$(0,0);-\frac5{168}$\quad $(-\frac12,0);\frac1{168}$\newline
$(\frac12,1);\frac{25}{168}$\\[1pt]
$Z_9$&$\begin{pmatrix}1&-1/2\\-1/2&3/4\end{pmatrix}$&3&
$(0,0);-\frac3{56}$\quad $(0,\frac12);\frac1{56}$\newline
$(-\frac12,\frac14);\frac1{28}$\\\midrule
$Z_4$&$\begin{pmatrix}4&1\\1&1\end{pmatrix}$&4&
$(0,\frac12);\frac1{120}$\quad $(2,\frac12);\frac{49}{120}$\\[1pt]
$Z_5$&$\begin{pmatrix}1/3&-1/3\\-1/3&4/3\end{pmatrix}$&4&
$(\frac12,0);\frac1{120}$\quad $(-\frac16,\frac23);\frac3{40}$\\
\end{longtable}
\endgroup

%% file: sections/degree_bound.tex
\part{Degree reduction}\label{part:degree}
\input{sections/degree_roadmap.tex}
\section{Reduction to degree four}\label{deg:sec:reduction}
\begin{theorem}[Degree-four bound]\label{deg:thm:main}
Let $A\in M_2(\Q)$ be symmetric and positive definite, let
$\beta\in\Q^2$ and $c_0\in\Q$, and let $(x,y)$ be the distinguished
Nahm point of $A$. If $f_{A,\beta,c_0}$ is modular, then
\[
                         [\Q(x,y):\Q]\le4.
\]
\end{theorem}

The initial bound $24$ for this same joint field is proved here in
Theorem~\ref{d24:thm:main}. Its symmetric proof and exact computation
inputs are included in this article and its computational companion. The general
argument is that of \cite{D24ref}.

\subsection{The coefficient reduction}
The following stronger technical statement makes the role of the
previous degree bound explicit. The $\lambda_j$ are the formal
logarithmic coefficients defined in Section~\ref{deg:sec:coefficients}.

\begin{theorem}[Reduction from degree twenty-four]\label{deg:thm:completed-reduction}
Let $A$ be a rational positive-definite symmetric $2\times2$ matrix,
let the two linear shifts be rational, and let $(x,y)\in(0,1)^2$ be
its distinguished Nahm point. Put $K=\mathbb Q(x,y)$. If
$[K:\mathbb Q]\le24$ and
\[
 \lambda_1\in\mathbb Q,\qquad
 \lambda_2=\lambda_3=\lambda_4=\lambda_5=0,
\]
then $[K:\mathbb Q]\le4$.
\end{theorem}

\begin{proof}[Proof of Theorem~\ref{deg:thm:main} from
Theorem~\ref{deg:thm:completed-reduction}]
Apply Theorem~\ref{d24:thm:main} to the original modular datum
to obtain $[K:\Q]\le24$.
Lemma~\ref{deg:lem:modular-bridge} below gives
$\lambda_1=c_0+1/12\in\Q$ and $\lambda_j=0$ for $j\ge2$.
Theorem~\ref{deg:thm:completed-reduction} therefore applies to this same
joint field. No modularity of a transformed indefinite matrix is used.
\end{proof}

The proof of Theorem~\ref{deg:thm:completed-reduction} is assembled after
the twenty-nine terminal cases in Section~\ref{deg:sec:final-assembly}.
The computer-assisted parts include exact algebraic identities,
finite-field censuses and finite local refinements. Every refinement
is used with an all-depth argument explaining why its finite data
control the remaining analytic terms. The component arguments are developed in
Sections~\ref{canx:section}, \ref{cuspx:section} and~\ref{spx:section};
the finite parameter calculations are specified in
Section~\ref{pix:section}. Appendix~\ref{deg:sec:supplement}
identifies the mathematical role of the finite computations. Their
programs and output records are distributed in the companion \cite{Comp}.

\section{The formal coefficients}\label{deg:sec:coefficients}
Let
\[
 A=\begin{pmatrix}a&b\\b&d\end{pmatrix}\in M_2(\Q),\qquad
 \beta=(B,C)^t\in\Q^2,
\]
with $A$ positive definite. Let $(x,y)\in(0,1)^2$ be the distinguished
solution of
\[
 1-x=x^ay^b,\qquad 1-y=x^by^d.
\]
Set $X=x/(1-x)$, $Y=y/(1-y)$, $K=\Q(x,y)=\Q(X,Y)$,
$\Delta=ad-b^2$, and $D=(a+X)(d+Y)-b^2$.
The scope is the symmetric case $m=1$. The scalar in the exponent of
the Nahm sum is denoted by $c_0$, so it is not confused with the second
linear shift $C$.

The logarithmic asymptotic coefficients are fixed algebraically. Put
$P_0(u)=u$, $P_{n+1}(u)=u(1+u)P'_n(u)$, and let
$\mathsf B_0=1,\mathsf B_1=-1/2,\ldots$ be the Bernoulli numbers. Define
\[
 U_n(z)=-\sum_{j=0}^{\lfloor(n+2)/2\rfloor}
 \frac{\mathsf B_j(-1)^{n-2j+2}}{j!(n-2j+2)!}
 \sum_{i=1}^2P_{n-j}(X_i)z_i^{n-2j+2}
 -\mathbf1_{n=1}(Bz_1+Cz_2),
\]
where $(X_1,X_2)=(X,Y)$. With Gaussian Wick expectation of covariance
$(A+\operatorname{diag}(X,Y))^{-1}$, write
\[
 \sum_{j\ge0}T_js^{2j}=\mathbb E\exp\left(\sum_{n\ge1}U_n(z)s^n\right),
 \qquad \sum_{j\ge1}\lambda_jt^j=\log\left(\sum_{j\ge0}T_jt^j\right).
\]
In the first identity only the even powers survive expectation.
The polynomial numerators are $\Psi_j=D^{3j}\lambda_j$.
The necessary conditions used in the reduction are
\begin{equation}\label{deg:C5}
 \lambda_1=r\in\Q,\qquad\lambda_2=\lambda_3=\lambda_4=\lambda_5=0.
\end{equation}
These are conditions on logarithmic coefficients, not on the
unmodified coefficients of the asymptotic expansion.

\begin{lemma}[Modularity and the coefficient convention]\label{deg:lem:modular-bridge}
If $f_{A,\beta,c_0}$ is modular, then
$\lambda_1=c_0+1/12$ and $\lambda_j=0$ for every $j\ge2$.
\end{lemma}
\begin{proof}
In the exponent of \cite[(2.5)]{VZ}, put $T=-z$ and
$\varepsilon=s^2$. Expanding the Bernoulli polynomial
$\mathsf B_p(u)=\sum_k\binom pk\mathsf B_k u^{p-k}$,
the term indexed by $(p,k)$ has degree $s^{p+k-2}$.
Also $\operatorname{Li}_{-n}(x)=P_n(X)$ since
$x\,dX/dx=X(1+X)$. Collecting terms with $p+k=n+2$
gives precisely $U_n$ above. There is one correction:
the condition $p\ge3$ in the cited exponent omits $(p,k)=(2,2)$,
whereas our formula contains $-(X+Y)s^2/12$.
The linear term is $-(B+X/2)z_1s-(C+Y/2)z_2s$ in both conventions.

Thus the Gaussian series in \cite[Theorem 2.3]{VZ} equals
$e^{(X+Y)\varepsilon/12}\sum_jT_j\varepsilon^j$.
Its external exponential scalar is
$\gamma=c_0+1/12+(X+Y)/12$.
The modular comparison \cite[Corollary 3.2]{VZ} makes
$e^{-\gamma\varepsilon}e^{(X+Y)\varepsilon/12}
\sum_jT_j\varepsilon^j$ equal to $1$ as a formal series.
Taking the formal logarithm proves the assertion.
The same comparison holds for a finite-order multiplier by applying
it to a power of the modular function and dividing its logarithm
by that positive integer.
\end{proof}

The bound $[K:\Q]\le24$ is a hypothesis of the technical reduction
and a conclusion of Theorem~\ref{d24:thm:main} in its modular
application. It will be denoted by D24 in the local arguments.
Its main use is to give $e_{\mathfrak p}\le24$ for every place above
seven. The rational transformations below preserve $K$.

\section{A self-contained initial degree bound}\label{d24:section}
This section gives the symmetric specialization of the degree argument
of \cite{D24ref}, including its common-component classification and
intersection count. The general symmetrizer case, its diagonal-edge
rigidity and its separate component enumeration are not needed here.
The original paper is cited for attribution, not as an additional proof input.
The computations used below are in the computational companion~\cite{Comp},
under \path{anc/global_degree/d24_m1/}. Unless explicitly identified as
an ancillary file of the arXiv submission, a file or directory name in
typewriter type refers to this companion after running
\path{restore_companion.py}. Paths are relative to its restored root;
abbreviated file names are resolved in its \path{CITED_FILE_INDEX.txt}.

\begin{theorem}[Initial symmetric degree bound]\label{d24:thm:main}
Let $A\in M_2(\Q)$ be symmetric and positive definite, and let the
two linear shifts and scalar shift be rational. If the associated
Nahm sum is modular, its distinguished solution satisfies
$[\Q(x,y):\Q]\le24$.
\end{theorem}

The argument uses only $\lambda_1=r\in\Q$ and
$\lambda_2=\lambda_3=\lambda_4=0$, together with the Nahm equations
in the rational exceptional families. Write
$H=A+\operatorname{diag}(X,Y)$, $D=\det H$ and
$R_1=D^3\lambda_1$. Throughout this section $m=1$ and
$\Delta=ad-b^2>0$. Rank-one numerators are denoted by
$N_j(A,\beta;Z)=(A+Z)^{3j}\lambda_A^{(j)}(\beta;Z)$,
using the one-variable version of the same Gaussian recurrence.
The letter $A$ in that rank-one notation is a scalar.

An implementation that completes the square removes the Gaussian linear
term by translation. Its first coefficient must then be supplemented
by $\tfrac12\gamma^tH^{-1}\gamma$, where
$\gamma=(B+X/2,C+Y/2)^t$. Its higher logarithmic coefficients are
unchanged. Thus the $R_1$ below and the $\Psi_j$ are the numerators already
fixed in Section~\ref{deg:sec:coefficients}; the scalar radial correction
$-1/12$ is outside that Gaussian integral. The companion provides a
separately implemented exact check of this first-coefficient identity.

\subsection{Coefficient geometry}\label{d24:sec:structure}
\subsubsection{Newton polygons}
Let $P_n$ denote the octagon $[0,n]^2$ with the four unit corner triangles removed.

\begin{proposition}\label{d24:prop:newt}
For $j=2,3,4$ and every specialisation of the parameters, $\operatorname{supp}\Psi_j\subseteq P_{3j}$, and the
eight corner coefficients of $P_{3j}$ are those of Proposition~\ref{d24:prop:corners}. For $j=3$, rational
parameters and positive-definite matrices, all eight corners are nonzero, so $\Newt\Psi_3=P_9$. For $j=2,4$ a corner can vanish (odd Bernoulli values, e.g.\ $B_3(\frac12)=0$), so
$\Newt\Psi_2,\Newt\Psi_4$ may be strictly smaller. Generically $\Psi_2$ and $\Psi_3$ have $45$ and $96$
monomials and total degrees $11$ and $17$.
\end{proposition}

\subsubsection{Corners are Bernoulli values}
\begin{proposition}[rank-one corners]\label{d24:prop:rk1corners}
For $j=2,3,4$,
\[
[X^{3j-1}]N_j(A,\beta;X)=\frac{(-1)^{j+1}}{(j+1)!}B_{j+1}\!\Bigl(\beta+\frac A2\Bigr),\qquad
[X^{1}]N_j(A,\beta;X)=-\frac{A^{3j}}{(j+1)!}B_{j+1}\!\Bigl(\frac{\beta+\frac12}{A}\Bigr).
\]
\end{proposition}

\begin{proposition}[edges]\label{d24:prop:edges}
Let $A_{\rm b}=\Delta/d$, $\beta_{\rm b}=B-bC/d$, $A_{\rm l}=\Delta/a$, $\beta_{\rm l}=C-bB/a$. Then for $j=2,3,4$:
\begin{align*}
\Psi_j(X,0)&=d^{3j}\,N_j(A_{\rm b},\beta_{\rm b};X), & \Psi_j(0,Y)&=a^{3j}\,N_j(A_{\rm l},\beta_{\rm l};Y), \\
[Y^{3j}]\Psi_j&=\,N_j(a,B;X), & [X^{3j}]\Psi_j&=\,N_j(d,C;Y).
\end{align*}

\end{proposition}
\begin{proof}[Proof of the first row]
At $Y=0$ all $P_k(Y)$ vanish and $q_0$ does not occur in $\lambda_j$, so every vertex on the $y$-side
vanishes, and so does every term containing $v_2$ (it only enters through $y$-vertices). What remains is
the rank-one expansion in $t_1$ with propagator $(H^{-1})_{11}=1/(A_{\rm b}+X)$ and shifted
leaf $v_1=(\beta_{\rm b}+X/2)/(A_{\rm b}+X)$, i.e. $\lambda^{(j)}_{A_{\rm b}}(\beta_{\rm b};X)$, times
$(\det H)^{3j}=d^{3j}(A_{\rm b}+X)^{3j}$. The line $X=0$ is symmetric. (The same argument gives
$X=-1$ and $Y=-1$ with $a-1$, $d-1$ in place of $a$, $d$.)
\end{proof}
\emph{Second row.} For $2\le j\le4$, with all five parameters symbolic, $[Y^r]\Psi_j=0$ for $r>3j$ and the
whole coefficient of $Y^{3j}$ equals $N_j(a,B;X)$. This was certified with the unshifted recurrence,
independently of \texttt{formalells}; the latter gives the same for $j\le4$
(\texttt{sym4.py}). The right
edge and the degree bound in $X$ follow by exchange symmetry (Lemma~\ref{d24:lem:swap}). Together with the first row this
proves Proposition~\ref{d24:prop:newt}.

\begin{lemma}[Exchange symmetry]\label{d24:lem:swap}
Exchanging $(a,B,X)$ with $(d,C,Y)$ preserves $R_1$ and every $\Psi_j$.
\end{lemma}
\begin{proof}
Exchange the two Gaussian coordinates in the defining recurrence of
Section~\ref{deg:sec:coefficients}. The covariance, vertices and determinant
are carried to those of the exchanged parameters.
\end{proof}

\begin{proposition}[corners of $\Psi_j$]\label{d24:prop:corners}
Each corner of $P_{3j}$ is the first or last coefficient of one of the four edge numerators. Writing
$x_1=\beta_{\rm b}+A_{\rm b}/2$, $x_2=(\beta_{\rm b}+\frac12)/A_{\rm b}$, $x_3=\beta_{\rm l}+A_{\rm l}/2$,
$x_4=(\beta_{\rm l}+\frac12)/A_{\rm l}$, $\xi=x_5=B+a/2$, $x_6=(B+\frac12)/a$, $\eta=x_7=C+d/2$,
$x_8=(C+\frac12)/d$, the corner coefficients are, up to the prefactors of Proposition~\ref{d24:prop:edges},
$\frac{(-1)^{j+1}}{(j+1)!}B_{j+1}(x_{\rm far})$ and $-\frac{A^{3j}}{(j+1)!}B_{j+1}(x_{\rm near})$.
For rational parameters and positive-definite matrices the corners of $\Psi_3$ never vanish,
since $B_4$ has no rational roots and $A_{\rm b},A_{\rm l},a,d>0$.
\end{proposition}

\subsubsection{Rank-one classification}
\begin{theorem}[common roots]\label{d24:thm:rk1}
Let $A>0$ and $\beta$ be rational. $N_2,N_3,N_4(A,\beta;\cdot)$ have a common root other than
$0,-1$ if and only if $(A,\beta)$ is one of Zagier's seven modular rank-one pairs; the extra common roots
are then the rank-one Nahm point and its conjugate:
$X^2-X-1$ for $(2,0),(2,1)$; $X^2+X-1$ for $(\frac12,0),(\frac12,\frac12)$; $X-1$ for $(1,0),(1,\pm\frac12)$.
Moreover $(1+X)^2$ divides all $N_j$ only when $A=1$; there $(1+X)^4$ divides all $N_j$, and for Euler
shifts $\beta\in\{0,\pm\frac12\}$ the gcd is $X(1+X)^6(X-1)$ (the even-$j$ numerators vanish identically).
\end{theorem}
\begin{proof}[Proof (computational)]
With $M_j=N_j/(X(1+X))$,
\[
\operatorname{Res}_X(M_2,M_3)=A^{15}(A-1)^{15}E(A,\beta),\qquad \operatorname{Res}_X(M_2,M_4)=A^{21}(A-1)^{21}E_2(A,\beta),
\]
(resultants in PARI's convention; exchanging the arguments of the second resultant changes its sign,
since $\deg M_2\cdot\deg M_4=27$ is odd, so the second identity holds up to sign in other conventions),
with $E$ irreducible of bidegree $(17,18)$ and $E_2$ of bidegree $(24,27)$; $\gcd(E,E_2)$ is constant, and
$U(A)=\operatorname{Res}_\beta(E,E_2)$ has degree $462$ with rational roots only at $A\in\{2,1,0,-1,\frac12\}$. At $A=2,\frac12$ the common $\beta$ are exactly the modular
shifts and the common roots are as stated. At $A=1$ one factors directly:
$N_2(1,\beta;X)=c\,\beta(4\beta^2-1)X(1+X)^4$, and $N_3,N_4$ contain $(1+X)^6,(1+X)^8$
with the extra factor $X-1$ in $N_3$. The multiplicity statement follows from
\[
M_2(-1)=c\,(\beta-1)(A-1)^3(A-\beta)(A-2\beta+1),
\]
and from the fact that on each of the three lines $\beta\in\{1,A,\frac{A+1}2\}$ the factor $M_3(-1)$
reduces to a power of $(A-1)$.
\end{proof}

\begin{corollary}[edge gcds]\label{d24:cor:D}
For each edge, the non-monomial part $D$ of $\gcd_{2\le j\le4}$ of the edge polynomials is:
\begin{center}
\begin{tabular}{ll}
\toprule
effective pair $(A,\beta)$ & $D(Z)$\\
\midrule
$A\notin\{\frac12,1,2\}$, or $A\in\{\frac12,2\}$ with non-modular $\beta$ & $1+Z$\\
$A=2$, $\beta\in\{0,1\}$ & $(1+Z)(Z^2-Z-1)$\\
$A=\frac12$, $\beta\in\{0,\frac12\}$ & $(1+Z)(Z^2+Z-1)$\\
$A=1$, $\beta\notin\{0,\pm\frac12\}$ & $(1+Z)^4$\\
$A=1$, $\beta\in\{0,\pm\frac12\}$ & $(1+Z)^6(Z-1)$\\
\bottomrule
\end{tabular}
\end{center}
\end{corollary}
\subsubsection{The rank-one conic}
\begin{theorem}\label{d24:thm:conic}
For $j=2,3$, $b\ne0$, $\kappa\ne0$, $X_\kappa=-a-\kappa b$, $Y_\kappa=-d-b/\kappa$ (so $\det H=0$):
\[
\kappa^{3j}\,\Psi_j(X_\kappa,Y_\kappa)=(-1)^jN_j\,b^{3j}\,\widehat Q(\kappa)^{2j},
\]
\[
\widehat Q(\kappa)=d(d-1)\kappa^3+b(2d-1+b)\kappa^2+b(2a-1+b)\kappa+a(a-1),
\]
with scalar constants $N_2=\frac5{16}$ and $N_3=\frac{1105}{1152}$
(distinct from the rank-one polynomials with arguments).
In particular the restriction is independent of $B,C$, and on the conic $\Psi_2$ vanishes only at the
roots of $\widehat Q$.
\end{theorem}
\noindent\emph{Conic condition.} If $\Gamma$ divides $\Psi_2$ and $b\ne0$, then every intersection point of
$\Gamma$ with the conic $\det H=0$ at a parameter $\kappa\in\C^*$ is a root of $\widehat Q$. With
$q_\Gamma(\kappa):=\kappa^{\deg_Y\Gamma}\Gamma(X_\kappa,Y_\kappa)=\sum_{i=0}^{D}c_i\kappa^i$, this is the
radical condition $\operatorname{rad}(q_\Gamma/\kappa^{v})\mid\widehat Q$, where $v$ is the order of vanishing
at $\kappa=0$. Since the $c_i$ depend on the parameters, it is applied branchwise: for each pair $v\le e$,
impose $c_0=\dots=c_{v-1}=0$, $c_{e+1}=\dots=c_D=0$, $c_vc_e\ne0$, and (if $e>v$)
$\tilde q=\sum_{i=v}^{e}c_i\kappa^{i-v}$ divides $\widehat Q^{\,e-v}$; the branch $q_\Gamma\equiv0$ means
$\Gamma$ is the conic itself, which has no point with $X,Y>0$ because $\Delta>0$. This covers repeated roots,
degree drops and roots at $\kappa=0$.

\begin{lemma}[weighted order at the conic points]\label{d24:lem:worder}
Let $b\ne0$ and let $\kappa_p\ne0$ be a root of $\widehat Q$ of multiplicity $r$. Near
$p=(X_{\kappa_p},Y_{\kappa_p})$ use the coordinates $s=\kappa-\kappa_p$ and $t$, where
\[
X=-a-\kappa b,\qquad Y=-d-\frac{b}{\kappa}-\frac{t}{b\kappa}
\]
(a local isomorphism near $p$, with Jacobian $1/\kappa$). Then $\det H=t$, and for every $j\ge2$
the polynomial $\Psi_j$ has weighted order at least $2rj$ for the weights $(\operatorname{wt}s,\operatorname{wt}t)=(1,2r)$.
\end{lemma}
\begin{proof}
Direct computation gives $\det H=t$ and
\[
\operatorname{adj}H=-\frac b\kappa\,uu^{\mathsf T}-\frac{t}{b\kappa}\,e_1e_1^{\mathsf T},\qquad u=(1,\kappa)^{\mathsf T}.
\]
The cubic Taylor term of the potential is $-[X(1+X)\xi^3+Y(1+Y)\eta^3]$, and on $u$ it equals
$V_3(u):=-[X(1+X)+\kappa^3Y(1+Y)]$, which is $-\widehat Q(\kappa)$ at $t=0$.

By Wick's theorem $\lambda_j$ is a sum over connected graphs $\Gamma$ whose vertices carry the terms of the Gaussian
recurrence in Section~\ref{deg:sec:coefficients}. A vertex $v$ of valence $r_v$ carries $\varepsilon^{(r_v-2)/2+\delta_v}$ with
$\delta_v\in\mathbb Z_{\ge0}$, and $\delta_v=0$ exactly for the Taylor terms of the potential, which have $r_v\ge3$.
Edges carry $H^{-1}=\operatorname{adj}H/t$. Comparing powers of $\varepsilon$ gives
$E=j+V-\sum_v\delta_v$ for the numbers $E$ of edges and $V$ of vertices.

Expand every edge into its $uu^{\mathsf T}$-part and its $e_1e_1^{\mathsf T}$-part (an \emph{$e$-edge}). A resulting
term of $\Psi_j=t^{3j}\lambda_j$ with $e$ $e$-edges has the factor $t^{k_0}$, $k_0=3j-E+e$. Every potential cubic
vertex none of whose three edges is an $e$-edge contributes the factor $V_3(u)$. If $V_3$ denotes the number
of potential cubic vertices, there are at least $V_3-2e$ such vertices. Now
$2E-4j=\sum_v(4-r_v-4\delta_v)$, where each summand is at most $0$ except for potential cubic vertices, where
it is $1$. Hence $V_3\ge2E-4j$, and there are at least $2E-4j-2e=2j-2k_0$ factors $V_3(u)$.

In the coordinates $(s,t)$, $V_3(u)=-\widehat Q(\kappa)+tR$ has weighted order $\ge r$, and all other
factors are analytic at $p$ because $\kappa_p\ne0$ and $b\ne0$. So each term has weighted order at least
$2rk_0+r\max(0,2j-2k_0)\ge2rj$.
\end{proof}
For $r=1$ the lemma says that $\widehat Q^{\,2j-2k}$ divides the coefficient of $t^k$ in $\Psi_j$ for $k<j$. This
was verified with all parameters symbolic for $j=2$, and exactly for $j\le5$ at three random points, where
the power of $\widehat Q$ is exactly $2j-2k$ (\texttt{conic\_order.py}, \texttt{conic\_sym.py}).

\subsection{The first coefficient and the width bound}\label{d24:sec:first}
Write $H=A+\operatorname{diag}(X,Y)$ and $D=\det H$.
The coefficient $\lambda_1$ is exactly the one defined in
Section~\ref{deg:sec:coefficients}, with no change of shift convention.
Modularity gives $\lambda_1(P)=r=c_0+1/12$.
Put $R_1=D^3\lambda_1$. With $U=d+Y$, $V=-b$, $W=a+X$, $g=B+X/2$, $h=C+Y/2$,
$p=X(1+X)$, $q=Y(1+Y)$, $p_2=p(1+2X)$, $q_2=q(1+2Y)$, Wick's formula gives
\begin{align}
R_1={}&-\frac{X+Y}{12}D^3+\frac{pU+qW}{4}D^2+\frac{g^2U+2ghV+h^2W}{2}D^2\notag\\
&-\frac{g(pU^2+qVW)+h(pUV+qW^2)}{2}D-\frac{p_2U^2+q_2W^2}{8}D\notag\\
&+\frac5{24}\bigl(p^2U^3+q^2W^3\bigr)+pq\Bigl(\frac{UVW}4+\frac{V^3}6\Bigr).
\label{d24:eq:R1}
\end{align}
In the square-completed implementation \texttt{formalells} the term $\frac12\gamma^{\mathsf T}H^{-1}\gamma$,
$\gamma=(g,h)$, is absorbed into the constant and must be restored; \eqref{d24:eq:R1} was checked against it with
all five parameters symbolic (\texttt{first\_coeff.py}). For $r\in\Q$ set
\[
f_r:=R_1-rD^3 .
\]
It has rational coefficients and vanishes at $P$ and at all its conjugates.

\begin{lemma}\label{d24:lem:first}
\begin{enumerate}
\item[(a)] $R_1$ has bidegree $(3,3)$, so $\operatorname{supp}f_r\subseteq[0,3]^2$.
\item[(b)] For $b\ne0$, on the conic $R_1(X_\kappa,Y_\kappa)=-\frac5{24}b^3\,\widehat Q(\kappa)^2\kappa^{-3}$.
Hence $\gcd(f_r,D)=1$, and $f_r$ has weighted order $\ge2$ for the weights $(1,2)$ of Lemma~\ref{d24:lem:worder}
at every simple root of $\widehat Q$.
\item[(c)] For $b=0$, $R_1(-a,Y)=\frac5{24}a^2(a-1)^2(d+Y)^3$.
\item[(d)] $f_r\equiv0$ only if $b=0$ and $a=d=1$; then $x=y=\frac12$ and $K=\Q$.
\item[(e)] The corners: $[X^3Y^3]f_r=(\frac1{12}-r)$, which vanishes iff $r=1/12$;
$[X^3]f_r=d^3(\frac1{24}+\frac{C^2}{2d}-r)$; $[Y^3]f_r=a^3(\frac1{24}+\frac{B^2}{2a}-r)$;
$[1]f_r=\Delta^3\bigl(\frac{dB^2-2bBC+aC^2}{2\Delta}-r\bigr)$.
\end{enumerate}
\end{lemma}
\begin{proof}
(a), (b), (c), (e) are read off from \eqref{d24:eq:R1} (\texttt{first\_coeff.py}, \texttt{first\_corners.py}); in (b)
the $t$-linear terms have weight $2$ and $rD^3=rt^3$ has weight $6$. (d) If $b\ne0$, (b) shows
$f_r\not\equiv0$ unless $\widehat Q\equiv0$, which with $a,d>0$ forces $a=d=1$ and then $b=-1$,
so $\Delta=1-1=0$, contradicting $\Delta>0$. If $b=0$, (c) forces $a=1$, and the exchange symmetry (Lemma~\ref{d24:lem:swap}) forces
$d=1$; the Nahm equations then give $x=y=\frac12$.
\end{proof}

\paragraph{Edges of $f_r$.} Let $\lambda^{(1)}_A(\beta;Z)$ be the first coefficient of the rank-one formal integral
without its factor $e^{-\varepsilon/24}$, and $N_1(A,\beta;Z):=(A+Z)^3\lambda^{(1)}_A(\beta;Z)$ the rank-one first
numerator. Explicitly
\begin{multline*}
24N_1=Z^3+\bigl(3A^2+(12\beta-3)A+12\beta^2-12\beta+2\bigr)Z^2\\
+\bigl(-2A^3+(12\beta+6)A^2+(24\beta^2-12\beta-3)A\bigr)Z+12A^2\beta^2,
\end{multline*}
one has, with all parameters symbolic (\texttt{fr\_edges.py}, \texttt{fr\_divides.py}),
\begin{align*}
[Y^0]f_r&=d^3\bigl[N_1(A_{\rm b},\beta_{\rm b};X)-\rho_{\rm b}(A_{\rm b}+X)^3\bigr], & \rho_{\rm b}&=r-\tfrac{C^2}{2d},\\
[Y^3]f_r&=\bigl[N_1(a,B;X)-\rho_{\rm t}(a+X)^3\bigr], & \rho_{\rm t}&=r-\tfrac1{24},\\
[X^0]f_r&=a^3\bigl[N_1(A_{\rm l},\beta_{\rm l};Y)-\rho_{\rm l}(A_{\rm l}+Y)^3\bigr], & \rho_{\rm l}&=r-\tfrac{B^2}{2a},\\
[X^3]f_r&=\bigl[N_1(d,C;Y)-\rho_{\rm r}(d+Y)^3\bigr], & \rho_{\rm r}&=r-\tfrac1{24}.
\end{align*}
(The corner terms are the values of the rank-one $\lambda_1$ at $Z=0$ and $Z=\infty$: $N_1(A,\beta;0)=\frac12A^2\beta^2$
and $[Z^3]N_1=\frac1{24}$.)

\begin{lemma}[rank-one edge divisibility]\label{d24:lem:edgediv}
Let $A\in\Q_{>0}$, $\beta\in\Q$, $\rho\in\Q$. The polynomial $N_1(A,\beta;Z)-\rho(A+Z)^3$ is nonzero, and it divides
$N_2(A,\beta;Z)$ only if $A=1$.
\end{lemma}
\begin{proof}
Write $H=N_1-\rho(A+Z)^3$ and $N_2=Z(1+Z)M_2$ with $M_2$ cubic. $H$ has degree $3$ unless $\rho=\frac1{24}$, and
$H(0)=0$ iff $\rho=\beta^2/(2A)$; if both hold then $A=12\beta^2$. For $A\ne1$ we treat four cases, each an
elimination in $(A,\beta)$ (\texttt{fr\_divides2.py}, \texttt{fr\_divides3.py}, \texttt{fr\_divides4.py}).
(i) $\deg H=3$, $H(0)\ne0$: then $H\mid(1+Z)M_2$, so either $H=cM_2$, or $H(-1)=0$ and $H/(1+Z)$ divides $M_2$.
In the first case the $4\times3$ matrix of the coefficients of $(A+Z)^3$, $M_2$, $N_1$ has rank $\le2$; its
$3\times3$ minors have the common factor $A^2(A-1)^2$, and the resultant in $\beta$ of two of them, divided by it,
has rational roots only at $A\in\{0,1\}$. In the second case
$\rho=N_1(A,\beta;-1)/(A-1)^3$, and the remainder of $M_2$ modulo $H/(1+Z)$ gives two equations whose resultant
in $\beta$ has positive rational roots only at $A=1$ and $A=\frac23$; at $A=\frac23$ the unique solution
$\beta=\frac23$ has $\rho=\frac1{24}$, outside this case. (ii) $\rho=\frac1{24}$: $H$ is a quadratic, whose leading
coefficient vanishes only for $A=1-2\beta\pm\frac{\sqrt3}3$, and the two coefficients of the remainder of
$Z(1+Z)M_2$ modulo $H$ have a resultant with rational roots only at $A\in\{0,1\}$. (iii) $\rho=\beta^2/(2A)$,
$\rho\ne\frac1{24}$: $H=Zq$ with $q$ quadratic, and the same elimination for $q\mid(1+Z)M_2$ gives positive
rational roots $A\in\{1,3,\frac43,\frac13,\frac34,12\}$, of which all but $A=1$ lie on $A=12\beta^2$, where the
leading coefficient of $q$ vanishes and the resulting linear polynomial does not divide $N_2$. (iv) $A=12\beta^2$
is contained in case (ii). Finally $H\not\equiv0$: for $A\ne1$, $H(-A)=N_1(A,\beta;-A)=\frac5{24}A^2(A-1)^2\ne0$, and for $A=1$,
$H=(1+Z)^2\bigl[(\frac1{24}-\rho)Z+\frac{\beta^2}2-\rho\bigr]$, which vanishes identically only if $12\beta^2=1$.
\end{proof}

\paragraph{Corners of $f_r$.} A corner of $f_r$ can vanish (Lemma~\ref{d24:lem:first}(e)). At a vanishing corner the
two adjacent coefficients are Bernoulli values: the constant term of an edge vanishes when $\rho=\beta^2/(2A)$,
its cubic term when $\rho=\frac1{24}$, and then
\begin{equation}\label{d24:eq:B2}
[Z]\Bigl(N_1-\tfrac{\beta^2}{2A}(A+Z)^3\Bigr)=-\tfrac{A^3}{2}B_2\Bigl(\tfrac{\beta+1/2}{A}\Bigr),\qquad
[Z^2]\Bigl(N_1-\tfrac1{24}(A+Z)^3\Bigr)=\tfrac12B_2\Bigl(\beta+\tfrac A2\Bigr).
\end{equation}
These are the case $j=1$ of Proposition~\ref{d24:prop:corners}, with the same arguments $x_1,\dots,x_8$ and, by the
edge formulas above, with the prefactors of Proposition~\ref{d24:prop:edges} at $j=1$. Since $B_2$ has no rational
root, both coefficients adjacent to a vanishing corner are nonzero. (Checked symbolically, and on the actual
coefficients of $f_r$ at all four corners against $\Psi_2,\Psi_3,\Psi_4$; \texttt{check32.py}.)

\begin{lemma}[width]\label{d24:lem:width}
Let $G\in\C[X,Y]$ be a common divisor of $f_r$ and $\Psi_2$, where $f_r,\Psi_2\not\equiv0$. If $b\ne0$, or if
$b=0$ and $a,d\ne1$, then $\deg_XG\le2$ and $\deg_YG\le2$.
\end{lemma}
\begin{proof}
By Lemma~\ref{d24:lem:edgediv} the four edge polynomials of $f_r$ are nonzero, so $f_r$ has degree $3$ in each
variable and is divisible by neither $X$ nor $Y$. Suppose $\deg_XG=3$. Then $f_r=GQ$ with $Q\in\C[Y]$, and
$Q(0)\ne0$ since $Y\nmid f_r$. So $G(X,0)$ and the leading coefficient of $G$ in $Y$ are nonzero constant
multiples of the bottom and top edges of $f_r$, that is, of $N_1(A_{\rm b},\beta_{\rm b};X)-\rho_{\rm b}(A_{\rm b}+X)^3$ and
$N_1(a,B;X)-\rho_{\rm t}(a+X)^3$. As $G\mid\Psi_2$, they divide $\Psi_2(X,0)$ and the leading coefficient of $\Psi_2$ in $Y$.
If the bottom edge of $\Psi_2$ is not identically zero, it is a nonzero multiple of $N_2(A_{\rm b},\beta_{\rm b};X)$
(Proposition~\ref{d24:prop:edges}), and Lemma~\ref{d24:lem:edgediv} gives $A_{\rm b}=1$; if it vanishes identically, then
$N_2(A_{\rm b},\beta_{\rm b};\cdot)\equiv0$, which by the proof of Lemma~\ref{d24:lem:euler} again forces $A_{\rm b}=1$. In the same
way the top edge gives $a=1$. For $b\ne0$ this contradicts $A_{\rm b}=\Delta/d<a$; for $b=0$,
$A_{\rm b}=a\ne1$. The case $\deg_YG=3$ is the same with the left and right edges, where $A_{\rm l}=\Delta/a<d$.
\end{proof}
In particular $f_r$ itself never divides $\Psi_2$ in these cases.

\subsection{Sixfold cancellation along the determinant conic}\label{d24:sec:pole}
In the count of Section~\ref{d24:sec:count}, the case in which a component $G$ of $f_r$ through $P$ divides $\Psi_2$ and
$\Psi_3$ and is cut only by $\Psi_4$ needs more than the vanishing $\Psi_j\in(D^j)$ along the determinant conic: it gives
$8(w+h)\le32$. Here we show that on such a component $\Psi_4$ vanishes to order $6$ along $D=0$, which gives $24$.
Throughout, $H=\bigl(\begin{smallmatrix}a+X&b\\ b&d+Y\end{smallmatrix}\bigr)$, $D=\det H$, $\lambda_1$ is the first
coefficient of Section~\ref{d24:sec:first}, and $\lambda_j=\Psi_j/D^{3j}$ ($j\ge2$) are the higher logarithmic coefficients.

\subsubsection{A three-moment obstruction}
Let $Z$ be a standard formal Gaussian ($\EE Z^{2n}=(2n-1)!!$, odd moments zero) and
\[
 \EE\exp\bigl(usZ^3+vs^2Z^4+ws^3Z^5\bigr)=1+c_1s^2+c_2s^4+c_3s^6+O(s^8).
\]
By Wick's rule, $c_1=\frac{15}2u^2+3v$, $c_2=\frac{3465}8u^4+\frac{945}2u^2v+105uw+\frac{105}2v^2$ and
\[
 c_3=\tfrac{765765}{16}u^6+\tfrac{675675}{8}u^4v+\tfrac{45045}{2}u^3w+\tfrac{135135}{4}u^2v^2+10395uvw+\tfrac{3465}{2}v^3+\tfrac{945}{2}w^2 .
\]
\begin{lemma}\label{d24:lem:moments}
Over a field of characteristic zero, $c_1=c_2=c_3=0$ implies $u=v=w=0$.
\end{lemma}
\begin{proof}
$c_1=0$ gives $v=-\frac52u^2$, and then $c_2=105u(w-4u^3)$. If $u\ne0$ then $w=4u^3$ and $c_3=14490u^6\ne0$. If $u=0$ then
$v=0$ and $c_3=\frac{945}2w^2$, so $w=0$.
\end{proof}

\subsubsection{The pole-order lemma}
Work over $\C[[\tau]]$ with valuation $\ord$. Consider a normalized formal stationary-phase expansion in two variables
whose phase and amplitude jets are regular, whose Hessian has rank one at $\tau=0$ and is nonsingular over
$\C((\tau))$, with determinant of order $e>0$.
\begin{lemma}\label{d24:lem:poles}
If $\lambda_1$ is regular at $\tau=0$ and $\lambda_2=\lambda_3=0$, then $\ord\lambda_4\ge-6e$.
\end{lemma}
\begin{proof}
A linear change over $\C[[\tau]]$ makes one quadratic direction a unit; the formal Morse lemma in that direction (only
division by a unit is needed) and Gaussian integration in it leave a one-variable expansion with regular jets and
quadratic coefficient $h$ of order $e$. After scaling by $\sqrt\varepsilon$ the phase vertices are
$a_k(\tau)\varepsilon^{(k-2)/2}z^k$ ($k\ge3$, $a_k$ regular), with $\EE_hz^{2n}=(2n-1)!!\,h^{-n}$, and the amplitude
terms are $b_{q,k}\varepsilon^{q+k/2}z^k$ ($q,k\ge0$, not both $0$), each contributing a pole of at most $e$ per unit
of $\varepsilon$. A phase vertex of degree $k$ contributes at most $M_k=(ke-2\ord a_k)/(k-2)$ per unit of $\varepsilon$,
and $M_k\le\frac32e$ for $k\ge6$. With $M=\max\{\frac32e,M_3,M_4,M_5\}$, pairing Gaussian half-edges shows that every
contribution to $\lambda_j$ has valuation at least $-jM$, and taking the logarithm preserves this.
Suppose $M>\frac32e$. After a finite ramified extension, at valuation exactly $-jM$ only the vertices of degree $3,4,5$
attaining $M$ contribute, with leading coefficients $u,v,w$ not all zero, and the leading model is exactly
$\EE\exp(usZ^3+vs^2Z^4+ws^3Z^5)$. Since $M,2M,3M>0$, regularity of $\lambda_1,\lambda_2,\lambda_3$ forces the first three
coefficients of its logarithm, hence $c_1,c_2,c_3$, to vanish, contradicting Lemma~\ref{d24:lem:moments}. So
$M\le\frac32e$ and $\ord\lambda_4\ge-4M\ge-6e$.
\end{proof}

\subsubsection{Reciprocity}
Put $p=(a,b,d,B,C)$, $p^*=\bigl(\frac d\Delta,-\frac b\Delta,\frac a\Delta,\frac{dB-bC}\Delta,\frac{aC-bB}\Delta\bigr)$
and $K_0=\frac1{12}+\frac{dB^2-2bBC+aC^2}{2\Delta}$.
\begin{proposition}\label{d24:prop:reciprocity}
The following rational-function identities hold in $a,b,d,B,C,X,Y$:
\[
\begin{aligned}
D(p^*;X^{-1},Y^{-1})&=D(p;X,Y)/(\Delta XY),\\
\lambda_1(p;X,Y)+\lambda_1(p^*;X^{-1},Y^{-1})&=K_0,\\
\lambda_j(p;X,Y)&=(-1)^j\lambda_j(p^*;X^{-1},Y^{-1})\quad(j=2,3,4).
\end{aligned}
\]
\end{proposition}
\begin{proof}
The first identity is immediate. For the others, the numerators $N_1=R_1$ and $N_j=\Psi_j$ are generated by the
Gaussian recurrence followed by the logarithm recurrence (script \texttt{pole/recip\_cert.gp}), with all parameters
symbolic (and then specialized to $m=1$); they agree with the numerators of Section~\ref{d24:sec:structure}. Splitting $N_j$ into homogeneous parts of
degree $q$ in $(a,b,d,B,C)$, the transform $p\mapsto p^*$ multiplies the degree-$q$ part by $\Delta^{-q}$ after the
polynomial substitution $(a,b,d,B,C)\mapsto(d,-b,a,dB-bC,aC-bB)$, and $X^{3j}Y^{3j}N_j(X^{-1},Y^{-1})$ is a coefficient
reversal. The cleared differences
$24\bigl(\mathcal T_1+\Delta^4N_1\bigr)-D^3\Delta^3\bigl(12(dB^2-2bBC+aC^2)+2\Delta\bigr)$ and
$\mathcal T_j-(-1)^j\Delta^{3j}N_j$ ($j=2,3,4$), where $\mathcal T_j=\Delta^{L_j}(XY)^{3j}N_j(p^*;X^{-1},Y^{-1})$ with
$L_1=7$, $L_j=6j$, are identically zero (\texttt{pole/recip\_check.py}; the $q$-degrees are $7,12,18,24$).
\end{proof}

\subsubsection{Sixfold cancellation}
Bars denote bihomogenization on $\mathbb P^1\times\mathbb P^1$; $\bar\Psi_4$ has bidegree $(12,12)$ and $\bar D$
bidegree $(1,1)$.
\begin{theorem}\label{d24:thm:six}
Let $b\ne0$, and let $G$ be a reduced divisor of $f_r$ with $G\mid\Psi_2$ and $G\mid\Psi_3$, of bidegree $(w,h)$, and
$\nu:\widetilde{\mathcal C}\to\mathcal C$ the normalization of the closure of $G=0$. Then $\nu^*\bar\Psi_4/(\nu^*\bar D)^6$
is a regular section of $\nu^*\mathcal O(6,6)$. If $G$ and $\Psi_4$ have no common component, there are at most
$6(w+h)$ points of $G=\Psi_4=0$ with $D\ne0$.
\end{theorem}
\begin{proof}
At a finite point of $G$ with $D=0$, $H$ has rank one ($b\ne0$), all jets are polynomial in $X,Y$, and on the branch
$\lambda_1=r$ and $\lambda_2=\lambda_3=0$; $D$ does not vanish identically on the branch (Section~\ref{d24:sec:count}). By
Lemma~\ref{d24:lem:poles}, $\ord\Psi_4=12\ord D+\ord\lambda_4\ge6\ord D$. The points at infinity of the conic are
$(\infty,-d)$ and $(-a,\infty)$, and $(\infty,\infty)$ is not on it. Since $a,d>0$, the map $(X,Y)\mapsto(X^{-1},Y^{-1})$
sends them to finite points of the conic of $p^*$, whose off-diagonal entry is still nonzero, and by
Proposition~\ref{d24:prop:reciprocity} the dual data again have $\lambda_1$ constant and $\lambda_2=\lambda_3=0$. In the chart
$u=1/X$ near $(\infty,-d)$ one has $\bar D=uD$, $\bar\Psi_4=u^{12}\Psi_4$ and $D(p^*;u,Y^{-1})=\bar D/(\Delta Y)$ with
$\Delta Y$ a unit, so Lemma~\ref{d24:lem:poles} for the dual data gives $\ord\bar\Psi_4\ge6\ord\bar D$; the other chart is
the same. Hence the quotient has no pole above $\bar D=0$, it is a section of $\nu^*\mathcal O(6,6)$, and it is nonzero
on each component when $G$ and $\Psi_4$ share none, so it has at most $6(w+h)$ zeros; every point of $G=\Psi_4=0$ with
$D\ne0$ is one of them.
\end{proof}

\begin{proposition}\label{d24:prop:decoupled}
If $b=0$, $a,d>0$ and $a,d\ne1$, then $f_r$ and $\Psi_2$ have no common component.
\end{proposition}
\begin{proof}
For $b=0$, $\lambda_1=L_{a,B}(X)+L_{d,C}(Y)$ and $\lambda_2=T_{a,B}(X)+T_{d,C}(Y)$ with the rank-one first and
second logarithmic functions $L_{A,\beta}$, $T_{A,\beta}$, which are regular at infinity and satisfy, with
$c_A=\frac5{24}A^2(A-1)^2\ne0$,
\[
 L_{A,\beta}(Z)=\frac{c_A}{(Z+A)^3}+O\bigl((Z+A)^{-2}\bigr),\qquad
 T_{A,\beta}(Z)=\frac{36}{5}\,\frac{c_A^2}{(Z+A)^6}+O\bigl((Z+A)^{-5}\bigr).
\]
On a common component some coordinate, say $X$, is nonconstant; at a point above $X=-a$ the first equation
$\lambda_1=r$ forces $Y=-d$ with the same order $e$, and leading coefficients $\alpha+\beta=0$, $\alpha\ne0$. The second
equation then has leading coefficient $\frac{36}5(\alpha^2+\beta^2)=\frac{72}5\alpha^2\ne0$ at order $-6e$, a
contradiction. Components inside a determinant line are excluded since $c_a,c_d\ne0$.
\end{proof}
When $b=0$ and $a=1$ or $d=1$, the proof of Theorem~\ref{d24:thm:main} gives $[K:\Q]\le3$ directly.

\subsection{Common components in the symmetric case}\label{d24:sec:m1}
By Lemma~\ref{d24:lem:width}, a common irreducible factor has bidegree at
most $(2,2)$. The decoupled case has already been treated in
Proposition~\ref{d24:prop:decoupled}; it remains to classify $b\ne0$.

\begin{theorem}\label{d24:thm:m1}
Assume $b\ne0$, and let $G\in\Q[X,Y]$ be irreducible, dividing $f_r,\Psi_2,\Psi_3,\Psi_4$,
with $V(G)$ meeting $\{X,Y>0\}$. Then $d=a$, $b=1-a$, either $C=-B$ or $(B,C)=(-\frac12,-\frac12)$, and $G=XY-1$ up to a constant.
In these cases the Nahm point is $(x,y)=(\frac12,\frac12)$, that is $(X,Y)=(1,1)$, which lies on $G$.
\end{theorem}
The proof is computer-assisted; all its case data and programs are
distributed in the companion. Only the necessary direction is used.

\paragraph{The six rank-one lines.} On each of the lines $Y\in\{0,\infty,-1\}$ and $X\in\{0,\infty,-1\}$ the
restriction of $\Psi_j$ is a nonzero constant times $N_j(A_L,\beta_L;\cdot)$, and that of $f_r$ a nonzero constant
times $N_1(A_L,\beta_L;\cdot)-\rho_L(A_L+\cdot)^3$, with 
\begin{center}\small
\begin{tabular}{lll}
\toprule
line & $(A_L,\beta_L)$ & $\rho_L$\\
\midrule
$Y=0$ & $(\Delta/d,\ B-bC/d)$ & $r-C^2/(2d)$\\
$Y=\infty$ & $(a,\ B)$ & $r-\frac1{24}$\\
$Y=-1$ & $((\Delta-a)/(d-1),\ B-b(C-\frac12)/(d-1))$ & $r-(12C^2-12C+2d+1)/(24(d-1))$\\
$X=0$ & $(\Delta/a,\ C-bB/a)$ & $r-B^2/(2a)$\\
$X=\infty$ & $(d,\ C)$ & $r-\frac1{24}$\\
$X=-1$ & $((\Delta-d)/(a-1),\ C-b(B-\frac12)/(a-1))$ & $r-(12B^2-12B+2a+1)/(24(a-1))$\\
\bottomrule
\end{tabular}
\end{center}
(Proposition~\ref{d24:prop:edges} and the edges of $f_r$ in Section~\ref{d24:sec:first}; on $X=-1$ all Taylor terms of
valence $\ge2$ in $x$ vanish, since $P_k(-1)=0$ for $k\ge1$, so $x$ can be integrated out as a Gaussian, and
symmetrically on $Y=-1$; the resulting formulas were checked symbolically, \texttt{rho\_lines.py}). The lines $X=-1$, $Y=-1$ require $a\ne1$, $d\ne1$ respectively. Since $G$
divides $\Psi_2,\Psi_3,\Psi_4$ and $f_r$, the restriction of $G$ to each line divides $\gcd_{2\le j\le4}N_j(A_L,\beta_L;\cdot)$ and
the $f_r$-restriction.

\paragraph{Common roots.} For rational $(A,\beta)$ the common roots of $N_2,N_3,N_4$ other than $0,-1$ occur only
at the pairs in the following table (the eliminant of Theorem~\ref{d24:thm:rk1} has rational roots only at
$A\in\{2,1,0,-1,\frac12\}$, and $A\le0$ can occur only on the lines $X,Y=-1$; \texttt{special\_pairs.py}):
\begin{center}\small
\begin{tabular}{ll}
\toprule
$(A,\beta)$ & $\gcd(N_2,N_3,N_4)$ besides $Z$, $1+Z$\\
\midrule
$(2,0),(2,1)$ & $Z^2-Z-1$\\
$(\frac12,0),(\frac12,\frac12)$ & $Z^2+Z-1$\\
$(-1,\pm\frac12)$ & $Z^2+3Z+1$\\
$A=1$ & $(1+Z)^3$ more; for $\beta\in\{0,\pm\frac12\}$ also $(1+Z)^2(Z-1)$\\
$A=0$ & $Z^3$ more; for $\beta\in\{0,\frac12,1\}$ also $Z^2(Z+2)$\\
\bottomrule
\end{tabular}
\end{center}
The roots $0$ and $\infty$ are otherwise simple (the relevant corners are values of $B_4$, which has no rational
root). Hence on each line the roots of $G$ lie in $\{0,-1,\infty\}$, all simple, unless the line's pair is special.
A factor may, however, use an ordinary root once on a line whose exponent is $A=0$ or $A=1$, where the rank-one
polynomials have that root with higher multiplicity. We therefore impose every root condition as the polynomial
equation $H_L(z)=0$, where
\begin{gather*}
H_L(Z)=N_1(A_L,\beta_L;Z)-\rho_L(A_L+Z)^3,\\
H_{0,\beta,\rho}=Z^2\Bigl[(\tfrac1{24}-\rho)Z+\tfrac{6\beta^2-6\beta+1}{12}\Bigr],\qquad
H_{1,\beta,\rho}=(Z+1)^2\Bigl[(\tfrac1{24}-\rho)Z+\tfrac{\beta^2}2-\rho\Bigr],
\end{gather*}
and never through the quotients $N_1(A,\beta;0)/A^3=\beta^2/(2A)$ or $N_1(A,\beta;-1)/(A-1)^3$. For $b\ne0$ and $a,d\ne1$ a
finite line has exponent $0$ or $1$ exactly on the loci $\Delta=a$, $\Delta=d$ and $\Delta=a+d-1$, which we treat as
separate strata.

\paragraph{Root patterns.} Prescribing the roots (with multiplicity) of $G$ on the six lines is a set of linear
conditions on its coefficients, and each root $z$ on a line $L$ imposes $H_L(z)=0$ on the parameters; a special root
also imposes the corresponding special pair. A pattern is admitted when some $G$ in its solution space is nonzero on
every line; this holds exactly when no line restriction vanishes on the whole rational kernel (a vector space over
$\Q$ is not a finite union of proper subspaces). The resulting systems in $(a,b,d,B,C,r)$ are solved exactly, by
Gr\"obner bases in \textsc{Singular} \cite{Sing} saturated by $ab\,d(a-1)(d-1)\Delta$. A curve is discarded for missing the
quadrant only when all its coefficients have one sign, and as reducible only by an exact factorization or a
common factor of the whole solution space.

\begin{lemma}\label{d24:lem:tangent}
Let $a,d>0$, $b\ne0$, $\Delta>0$. If $XY-1$ divides $f_r$ and is tangent to the conic $D=0$, then $a=d$ and $b=1-a$.
\end{lemma}
\begin{proof}
On $XY=1$ the conic becomes $q(X)=dX^2+(\Delta+1)X+a$. A double root means $t:=\sqrt{ad}=(\Delta+1)/2\in\Q$; with
$u=a/t$ one gets $a=tu$, $d=t/u$, $b=\pm(t-1)$, $t>\frac12$, $t\ne1$, and the double root is $X=-u$. In the
parametrization $X=-a-\kappa b$ of the conic, for $b=t-1$ the root has $\kappa=-u$ and $\widehat Q(-u)=2u(u-1)$, so the
conic condition forces $u=1$; then $X^3f_r(X,X^{-1})/(X+1)^2$ takes the value $-\frac56(t-1)^3\ne0$ at $X=-1$, so $XY-1\nmid
f_r$. For $b=1-t$ one has $\widehat Q(u)=0$, and $X^3f_r(X,X^{-1})/(X+u)^2$ takes the value $-\frac32t^2(t-1)(u-1)^2$ at
$X=-u$, forcing $u=1$, i.e.\ $a=d=t$, $b=1-t$. Both values are independent of $B,C,r$ (\texttt{m1/fix\_tangent.py}).
\end{proof}

\begin{proposition}\label{d24:prop:ord22}
Let $a,d>0$, $\Delta>0$, $b\ne0$ and $a,d\ne1$. There is no $G\in\Q[X,Y]$ of bidegree $(2,2)$ dividing both $f_r$ and
$\Psi_2$ whose bihomogenization is nonzero on each of the six lines $X,Y\in\{0,\infty,-1\}$ and whose restriction to each
of them has two distinct roots in $\{0,\infty,-1\}$.
\end{proposition}
Irreducibility and positivity are not assumed.
\begin{proof}
Index $0,\infty,-1$ by $0,1,2$. Each row and column of the $3\times3$ grid contains exactly two zeros of $G$, so the
missing points are $(i,\pi(i))$ for a permutation $\pi$, written as the word $\pi(0)\pi(1)\pi(2)$. Exchanging $X,Y$
interchanges $120$ and $201$ and fixes the other four patterns. Let $V_{ij}$ be the value at the grid point $(i,j)$ of
the bidegree-$(3,3)$ homogenization of $24f_r$, with homogeneous coordinates $(1,0),(0,1),(1,-1)$ for $0,\infty,-1$.
Then $V_{ij}=p_{ij}E_{ij}$ with
\[
(p_{ij})=\begin{pmatrix}-12\Delta^2&-a^2&-(\Delta-a)^2\\-d^2&-2&-3(d-1)^2\\-(\Delta-d)^2&-3(a-1)^2&-(\Delta-a-d+1)^2\end{pmatrix}
\]
and the cores
\begin{gather*}
E_{00}=-B^2d+2BCb-C^2a+2r\Delta,\qquad E_{11}=12r-1,\\
E_{01}=-12B^2+24ar-a,\qquad E_{10}=-12C^2+24dr-d,\\
E_{12}=-4C^2+4C+8dr-d-8r,\qquad E_{21}=-4B^2+4B+8ar-a-8r,\\
E_{02}=12E_{00}+12B^2-12Bb+12aC-2\Delta-a(24r+1),\\ E_{20}=12E_{00}+12C^2-12Cb+12dB-2\Delta-d(24r+1),\\
E_{22}=E_{02}+E_{20}-12E_{00}-12(B+C)+2(a+d)+6b+24r+2
\end{gather*}
(an identity in all parameters, checked symbolically). Divisibility requires $p_{ij}E_{ij}=0$ for $j\ne\pi(i)$, and a
core may be set to zero only where its prefactor is nonzero. The prefactors vanish only on
\[
\Delta=a,\qquad \Delta=d,\qquad \Delta=a+d-1. \tag{$\ast$}
\]
We use the rank-one facts: for $A\ne0,1$, roots $\{0,\infty\}$ mean $A=12\beta^2$; roots $\{-1,\infty\}$ mean
$A=1-12(\beta-\frac12)^2$; roots $\{0,-1\}$ mean $(6\beta-3A)^2=3A(A-1)$. So for positive $A\ne1$ the pair $\{-1,\infty\}$
forces $A<1$ and $\{0,-1\}$ forces $A>1$; at $A=1$ the pair $\{0,\infty\}$ would need $\beta^2=\frac1{12}$. The exponents
are those of the six-line table.

\emph{The stratum $\Delta=a$} (and $\Delta=d$ by exchange). Here $d=1+b^2/a>1$. On $X=0$ the exponent is $1$, so the
roots include $-1$; on $X=\infty$ it is $d>1$, so they include $0$. Only $012$, $021$, $120$ remain. For $120$ the top
pair $\{-1,\infty\}$ gives $a<1$, while the exponent $(a-d)/(a-1)>1$ on $X=-1$ also carries $\{-1,\infty\}$. For $012$
the top and right pairs give $a,d>1$, the bottom pair gives $a/d\le1$ and the pair $\{0,\infty\}$ on $X=-1$ gives
$(a-d)/(a-1)\ge0$; so $a=d$, $b^2=a(a-1)$, and the top equation becomes $(6B-3a)^2=3b^2$, impossible. For $021$ one has
$a=12B^2$, $d=12C^2$, $r=\frac1{12}$; with $w=b-12BC$, the conditions at $(-1,0)$ and $(-1,-1)$ reduce, for $a\ne d$, to
\[
3B^2(2B-1)^2+C^2(24B^2-12B+1)=0,\qquad (3B-1)(2B+2C-1)(12BC-2B-4C+1)=0,
\]
whose resultant in $C$ is $(2B-1)^4(3B-1)^2(6B-1)^3(72B^3-36B^2+6B-1)$; the cubic has no rational root, since three
times it is $(6B-1)^3-2$. For $a=d$ one has $C=\pm B$ and the remaining factors have no admissible rational zeros. What
survives is the points $P_1,P_2,P_2^{\rm op},P_3$ below.

\emph{The stratum $\Delta=a+d-1$.} Now $b^2=(a-1)(d-1)$ and both exponents on $X,Y=-1$ equal $1$, so $021$, $120$, $201$,
$210$ remain. For $120$ the top pair gives $a<1$ and the left pair $\Delta/a>1$, so $d>1$, contradicting $b^2>0$;
$201$ is symmetric. For $021$ the conditions at $(0,-1)$ and $(-1,0)$ give $(2C-1)w=2B$ and $(2B-1)w=2C$, hence $B=C$ and
$b=\pm(12B^2-1)$; the sign $+$ gives $B=\frac14$, the sign $-$ gives $48B^3-24B^2+1=0$, whose reciprocal is Eisenstein at
$3$. This leaves $P_4$. For $210$ one has $a,d>1$ and, with $a=\lambda B^2$, $d=\lambda C^2$,
\[
2\lambda^2B^2-12\lambda B+\lambda+12=0=2\lambda^2C^2-12\lambda C+\lambda+12,
\]
so $B=C$ or $B+C=6/\lambda$. If $B=C$, then $a=d$, $b=\pm(a-1)$, and $E_{00}=0$ gives $P_5$ (for $b=1-a$ it gives
$4a^3-6a+3=0$, Eisenstein at $3$). If $S:=B+C=6/\lambda$, then
\[
BC=\tfrac{S(2S+1)}{12},\quad a+d=4S-1,\quad ad=(S+\tfrac12)^2,\quad r=\tfrac{2S+1}{24},\quad
E_{00}=-\tfrac{(2S+1)(2b-2S-1)(2b-2S+1)}{48},
\]
so $b=S-\frac12$ and $\Delta=2S$; the stratum gives $\Delta=4S-2$, hence $S=1$ and $B=C=\frac12$, a contradiction.

\emph{Off the strata $(\ast)$}, all present cores vanish. Pattern $102$: the top and right pairs give $a,d<1$, but the
bottom pair $\{0,-1\}$ needs $\Delta/d>1$, while $\Delta/d<a<1$. Pattern $120$ (and $201$): here $r=\frac1{12}$,
$a=-12B^2+12B-2$, $d=12C^2$; $E_{02}=0$ gives $b=(2B+2C-1)(12BC-4B-2C+1)/(2B(2C-1))$, and then
\[
E_{00}=-\frac{FT}{24B^2(2C-1)^2},\qquad E_{22}=\frac{HT}{2B^2(2C-1)^2}
\]
with $T=4MC^2+(-24B^2+20B-4)C+(2B-1)^2$, $M=12B^2-6B+1$, and quadratics $F,H$ in $C$; since
$\operatorname{disc}_C T=-48B^2(2B-1)^2<0$ and $\operatorname{Res}_C(F,H)=9216B^6M>0$, there is no real solution. Pattern $012$: $a,d>1$;
if $B=C$, then $E_{02}=0$ is a quadratic in $b$ whose discriminant, modulo the top conic, is $-4a^3(a-1)<0$; if $B\ne C$,
then $E_{02}-E_{20}=0$ gives $b=S-\frac12$ and $\Delta=2S$, while the bottom and left pairs $\{-1,\infty\}$ need
$\Delta<\min(a,d)\le\frac{a+d}2=2S-\frac12$. Pattern $021$ gives the family $A$ below, together with zeros that lie on
$(\ast)$. Pattern $210$ gives, for $B=C$, no solution (the eliminant has resultant $-110592a^5(a-1)(a^2-a+1)\ne0$ with the
top conic), and for $B\ne C$ the family $B$ below.

\emph{Certificates.} The survivors are
\begin{center}\small
\begin{tabular}{lcccccc}
\toprule
& $a$ & $b$ & $d$ & $B$ & $C$ & $r$\\ \midrule
$P_1$ & $\frac43$ & $-\frac23$ & $\frac43$ & $\frac13$ & $\frac13$ & $\frac1{12}$\\
$P_2$ & $\frac43$ & $\frac23$ & $\frac43$ & $\frac13$ & $-\frac13$ & $\frac1{12}$\\
$P_3$ & $\frac13$ & $-\frac13$ & $\frac43$ & $\frac16$ & $\frac13$ & $\frac1{12}$\\
$P_4$ & $\frac34$ & $-\frac14$ & $\frac34$ & $\frac14$ & $\frac14$ & $\frac1{12}$\\
$P_5$ & $\frac32$ & $\frac12$ & $\frac32$ & $\frac12$ & $\frac12$ & $\frac18$\\ \bottomrule
\end{tabular}
\end{center}
their images under exchange, and the rational families
\begin{align*}
A:&\ (a,b,d,B,C;r)=\bigl(12(\tfrac12-q)^2,\,-12q^2+6q-1,\,12q^2,\,\tfrac12-q,\,q;\,\tfrac1{12}\bigr),\\
B:&\ (a,b,d,B,C;r)=\Bigl(-\tfrac{(q-6)^2}{2h},\,-\tfrac{h+6}{2h},\,-\tfrac{q^2}{2h},\,\tfrac{q-6}{2h},\,-\tfrac{q}{2h};\,\tfrac{h-6}{24h}\Bigr),\quad h=q^2-6q+6.
\end{align*}
Let $F,P$ be $f_r,\Psi_2$ (for $B$, $h^3f_r$ and $h^6\Psi_2$), scaled by rational constants to primitive integer
polynomials. At each point, $\operatorname{Res}_Y(F(1,Y),P(1,Y))\equiv82,15,15,56,45\pmod{101}$ for $P_1,\dots,P_5$, with the full
$Y$-degrees $3,6$ retained, so $f_r$ and $\Psi_2$ have no common factor of positive $Y$-degree there. For each family let
$R_s(q)$ be the primitive part of $\operatorname{Res}_Y(F(s,Y,q),P(s,Y,q))$, $s=1,2$; both have degree $72$, also modulo $101$, and
$\operatorname{Res}_q(R_1,R_2)\equiv24$ (family $A$) and $52$ (family $B$) $\pmod{101}$. A common factor of positive $Y$-degree at some
$q_0$ would make both specialized Sylvester determinants vanish, also when the affine degree drops, so $R_1(q_0)=R_2(q_0)=0$,
which is impossible. The case analysis agrees with an independent computation of the same systems by Gr\"obner bases,
which finds exactly these points and families (\texttt{m1/fix\_ord22p.py}, \texttt{m1/fix\_points.py}); the grid identity
and all certificates are reproduced by \texttt{m1/ord22\_certificates.py}.
\end{proof}

\begin{itemize}
\item \emph{Only simple roots in $\{0,-1,\infty\}$.} Counting incidences with the $3\times3$ grid forces $w=h$; there are
six patterns of bidegree $(1,1)$ and six of bidegree $(2,2)$. For $(1,1)$, $G$ is the graph of a M\"obius map
permuting $\{0,-1,\infty\}$, and only $Y=X$ and $XY=1$ meet the quadrant. The line $\{X=Y\}$ meets the conic in two
distinct points (the discriminant is $(a-d)^2+4b^2>0$), so by Theorem~\ref{d24:thm:conic} both are roots of $\widehat Q$;
this leaves only families with $\Delta=0$ or Zagier's matrix, where $\Psi_2(X,X)\not\equiv0$. For $XY=1$ the same
condition with distinct points leaves $d=a$, $b=1-a$ or $ad=b^2+b+1$, and Lemma~\ref{d24:lem:tangent} covers tangency; the
diagonal corners then fix the shifts, giving exactly the family in the theorem. The six patterns of bidegree $(2,2)$ are excluded by Proposition~\ref{d24:prop:ord22}.
\item \emph{Some special root} ($b\ne0$, $a,d\ne1$). There are $450$ patterns. The systems are inconsistent, or their
parameter ideals have no admissible rational points (circles $x^2+y^2=\frac1{12}$; conditions such as $b^2+1=12C^2$,
which has no rational solution because $3$ is not a sum of two rational squares; sign conditions; zero-dimensional
primes of degree $>1$), or finitely many rational points where $\gcd(\Psi_2,\Psi_3)=1$. Six M\"obius curves meeting the
quadrant are excluded by the exponent conditions together with the conic condition in the form $\widehat Q^2\equiv0\bmod
N_G$, which is valid also at tangency. Three curves using the root $1$ on $X=0$ and $-2$ on $X=-1$ force
$\Delta/a=1$ and $(\Delta-d)/(a-1)=0$, hence $\Delta=a=d$ and exponent $0$ on $Y=-1$, contradicting their root data there.
\item \emph{$a=1$, $b\ne0$.} On $X=-1$, $\Psi_2,\Psi_3,\Psi_4$ restrict to constants times $Y(1+Y)$,
$Y(1+Y)(2Y+1)$, $Y(1+Y)(6Y^2+6Y+1)$ (\texttt{deg\_a1.py}), and $f_r$ to a polynomial of degree exactly $1$: its
$Y$-coefficient is a nonzero constant times $b^4(12z^2-12bz+2b^2)$ with $z=B-\frac12$, which would vanish only for $z/b=(3\pm\sqrt3)/6$.
So $G$ has at most one finite root there, in $\{0,-1\}$, or $-\frac12$ when $\Psi_2(-1,\cdot)\equiv0$. With this
alphabet there are $319$ patterns, all excluded (inconsistency, exact reducibility, one-sign coefficients,
divisibility, rationality, or signs: on a line, roots $\{0,-1\}$ need $A\notin(0,1)$, since the conic has discriminant
$48A(A-1)$ in $\beta$, while $A_{\rm b}=1-b^2/d\in(0,1)$). The case $d=1$ follows by the exchange symmetry, verified
symbolically for $R_1,\Psi_2,\Psi_3$ (\texttt{sym\_check.py}).
\item \emph{$a=d=1$, $0<|b|<1$.} There are $31$ patterns, all excluded; the last two by the identities
\[
144(b^2-1)^2\bigl(C^2-C+\tfrac14\bigr)^2+4b^4=0,\qquad 144(b^2-1)^2\bigl(C^2-bC+\tfrac{2b^2-1}{12}\bigr)^2+4b^4=0,
\]
which are the eliminants of their parameter ideals and have no real solution with $b\ne0$.
\item \emph{Lines $X=X_0$, $Y=Y_0$.} A rational line through the quadrant needs $X_0=1$ (so $a=1$ and the bottom
pair has $A=1$, impossible since $A_{\rm b}<1$), and a pair of conjugate lines needs $a=A_{\rm b}=2$ or $\frac12$, which
gives $b=0$, where $\Psi_3$ restricts to a nonzero multiple of $N_3(d,C;\cdot)$.
\end{itemize}
The complete case table, with the parameter conditions and the reason for each exclusion, is in
\path{anc/global_degree/d24_m1/code/m1/CASES.md}. The completed exact
resolution of all $450+319+31=800$ parameter patterns, including every
projective divisibility chart and the residual rationality arguments,
is recorded in \path{anc/global_degree/replay/parameter_ideals/closure.json}
and the accompanying generators and certificates; see
Appendix~\ref{d24:certificates}. For the four special bilinear patterns
whose polynomial equations require clearing a nonzero denominator, the
polynomial-only proof in Section~\ref{d24:subsec:bilinear-polynomial}
supplies the required replacement, including tangency. These are the
certificates used here. This completes the proof of Theorem~\ref{d24:thm:m1}.

\subsection{Intersection count and the initial bound}\label{d24:sec:count}
\begin{lemma}[the Euler exception]\label{d24:lem:euler}
$\Psi_2\equiv0$ if and only if $(a,b,d)=(1,0,1)$, $B\in\{0,\pm\frac12\}$ and
$C\in\{0,\pm\frac12\}$. In that case $x=y=\frac12$ and $K=\Q$.
\end{lemma}
\begin{proof}
Suppose $\Psi_2\equiv0$. If $b\ne0$, Theorem~\ref{d24:thm:conic} forces $\widehat Q\equiv0$. Its cubic and constant
coefficients give $d(d-1)=a(a-1)=0$, so $a=d=1$ by positivity; the other coefficients then give $b(1+b)=0$, so $b=-1$ and then $\Delta=0$, contradicting $\Delta>0$. Hence $b=0$, and then
$\lambda_2(X,Y)=\lambda^{(2)}_a(B;X)+\lambda^{(2)}_d(C;Y)$. Both terms vanish at $0$ (each $N_2$ is divisible by
$X$), so setting $Y=0$, resp.\ $X=0$, shows that both are identically zero. Now $N_2(A,\beta;X)\equiv0$ forces
both Bernoulli corners of Proposition~\ref{d24:prop:rk1corners} to vanish:
$u:=\beta+\frac A2$ and $v:=(\beta+\frac12)/A$ are roots of $B_3(z)=z(z-\frac12)(z-1)$, so $u,v\in\{0,\frac12,1\}$,
$A=(u+\frac12)/(v+\frac12)$ and $\beta=Av-\frac12$. Three of the nine pairs have $A=1$, $\beta\in\{0,\pm\frac12\}$. The
other six are excluded by
\begin{gather*}
N_2(\tfrac12,-\tfrac14;1)=-\tfrac5{128},\qquad N_2(\tfrac13,-\tfrac16;1)=-\tfrac{25}{729},\qquad N_2(2,-\tfrac12;1)=-\tfrac52,\\
N_2(\tfrac23,\tfrac16;1)=\tfrac{65}{1458},\qquad N_2(3,-\tfrac12;1)=-25,\qquad N_2(\tfrac32,\tfrac14;1)=\tfrac{65}{128}.
\end{gather*}
Finally $N_2(1,\beta;X)=-\frac{\beta(4\beta^2-1)}{24}X(1+X)^4$ vanishes identically exactly for $\beta\in\{0,\pm\frac12\}$.
Conversely these parameters give $\Psi_2\equiv0$. The Nahm point is $x=y=\frac12$. (All values checked in
\texttt{euler\_check.py}.)
\end{proof}

\paragraph{Bihomogenization.} On $S_0=\mathbb P^1\times\mathbb P^1$ write $\bar F$ for the bihomogenization of a polynomial $F$
in $(X,Y)$ to a prescribed bidegree: $\bar D$ has bidegree $(1,1)$, where $D=\det H$, $\bar R_1$ and $\bar f_r$ have
bidegree $(3,3)$ (Lemma~\ref{d24:lem:first}(a)), and $\bar\Psi_j$ has bidegree $(3j,3j)$ (Proposition~\ref{d24:prop:newt}).
Membership in an ideal generated by sections is meant locally, in trivializations of the line bundles.

\begin{lemma}[ideal powers]\label{d24:lem:ideal}
Suppose that $\bar R_1$ and $\bar D$ have no common component, and let $\mathcal I=(\bar R_1,\bar D)=(\bar f_r,\bar D)$.
Then $\bar\Psi_j\in\mathcal I^j$ at every point of $S_0$, for $2\le j\le4$.
\end{lemma}
\begin{proof}
\emph{Generic parameters.} Take $b\ne0$ and $\widehat Q$ of degree $3$ with three simple roots, none zero. By
Lemma~\ref{d24:lem:first}(b), $\bar R_1$ restricted to the conic vanishes exactly at the three roots, to order two, and not
at the two points of the conic at infinity ($\kappa=0,\infty$). So $V(\mathcal I)$ consists of the three conic points,
and in the coordinates $(s,t)$ of Lemma~\ref{d24:lem:worder}, $\mathcal I=(R_1,t)=(s^2,t)$ there. Lemma~\ref{d24:lem:worder} gives
weighted order $\ge2j$ for the weights $(1,2)$, and a monomial $s^\alpha t^\gamma$ with $\alpha+2\gamma\ge2j$ lies in
$(s^2,t)^j$. So $\bar\Psi_j\in\mathcal I^j$.

\emph{All parameters.} Let $T\subset\C^5$ be the open set of parameters $(a,b,d,B,C)$ for which $\bar R_1$ and $\bar D$
have no common component; it is irreducible and contains the generic parameters above. On $S_0\times T$ the
subscheme $\mathcal Z=V(\mathcal I)$ is cut out by two equations and has finite fibres, so it is finite over $T$, a
complete intersection, and flat over $T$ \cite[Thm.~23.1]{Mat}. Since $\mathcal I$ is locally generated by a regular
sequence, $\mathcal I^\nu/\mathcal I^{\nu+1}\cong\operatorname{Sym}^\nu(\mathcal I/\mathcal I^2)$ is locally free over $\mathcal O_{\mathcal Z}$
\cite[Thm.~16.2]{Mat}, hence flat over $T$, and by induction on $j$ the sheaf $\mathcal O/\mathcal I^j$ is finite and flat
over $T$. The pushforward of $\mathcal O/\mathcal I^j\otimes\mathcal O(3j,3j)$ to $T$ is therefore locally free, and the class of
$\bar\Psi_j$ is a section of it that vanishes on a dense open subset, hence everywhere; since forming this quotient
commutes with restriction to fibres, $\bar\Psi_j\in\mathcal I^j$ on every fibre.
\end{proof}
Repeated finite conic roots are already covered by Lemma~\ref{d24:lem:worder}; the flatness argument extends the
statement to conic points at infinity and to the decoupled case. It was also checked at a simple and at a double
conic root on the boundary of $S_0$, where the required lower bounds hold (\texttt{check32c.py}).

\begin{corollary}\label{d24:cor:count}
Suppose that $\bar R_1$ and $\bar D$ have no common component. Let $G\in\Q[X,Y]$ be a reduced divisor of $f_r$ of
bidegree $(w,h)$, and let $2\le j\le4$ be such that $G$ and $\Psi_j$ have no common component. Then at most
$2j(w+h)$ points of $V(G)\cap V(\Psi_j)$ satisfy $D\ne0$.
\end{corollary}
\begin{proof}
Let $C\subseteq S_0$ be the closure of $V(G)$. On $C$ we have $\bar f_r=0$, so Lemma~\ref{d24:lem:ideal} gives
$\bar\Psi_j|_C=\sigma\,\bar D^j|_C$ with $\sigma$ a regular section of $\mathcal O(2j,2j)|_C$. Since $\bar D$ and $\bar f_r$ have
no common component, $\bar D$ vanishes on no component of $C$; since $G$ and $\Psi_j$ are coprime, neither does
$\sigma$. The degree of $\mathcal O(2j,2j)$ on $C$ is $2j(w+h)$. Every point of $C$ with $\Psi_j=0\ne D$ is a zero of $\sigma$,
and pulling $\sigma$ back to the normalization of $C$ shows that there are at most $2j(w+h)$ such points.
\end{proof}

\begin{proof}[Proof of Theorem~\ref{d24:thm:main}]
Let $P=(X_0,Y_0)$ be the Nahm point and $r=c_0+1/12$, so that $f_r(P)=\Psi_2(P)=\Psi_3(P)=\Psi_4(P)=0$.

\emph{Degenerate cases.} If $\Psi_2\equiv0$, Lemma~\ref{d24:lem:euler} gives $K=\Q$. If $b=0$ and $a=1$, the decoupled
Nahm equation $1-x=x$ gives $X_0=1$. At $b=0$, $R_1=(d+Y)^3N_1(a,B;X)+(a+X)^3N_1(d,C;Y)$
(\texttt{first\_coeff.py}), so $\lambda_1(P)=r$ reads $N_1(d,C;Y_0)-\rho(d+Y_0)^3=0$ with $\rho\in\Q$. This
polynomial in $Y_0$ is nonzero of degree at most $3$ (Lemma~\ref{d24:lem:edgediv}), so $[K:\Q]\le3$. The case $d=1$ is
symmetric.

\emph{General case.} Otherwise $\bar R_1$ and $\bar D$ have no common component: for $b\ne0$ the conic is irreducible and
$R_1$ does not vanish on it (Lemma~\ref{d24:lem:first}(b), since $\widehat Q\not\equiv0$ when $\Delta>0$); for $b=0$, $a,d\ne1$, neither
line $X=-a$, $Y=-d$ divides $R_1$ (Lemma~\ref{d24:lem:first}(c) and the exchanged identity). Also $f_r\not\equiv0$
(Lemma~\ref{d24:lem:first}(d)). By Theorem~\ref{d24:thm:m1}, when $b\ne0$ a common component through
$P$ is possible only on the family with $P=(1,1)$ and $K=\Q$.
When $b=0$ and $a,d\ne1$, Proposition~\ref{d24:prop:decoupled} excludes a common
component already for $f_r$ and $\Psi_2$. Outside these rational exceptions,
no common component of all four polynomials passes through $P$.
Let $G\in\Q[X,Y]$ be an irreducible factor of $f_r$ with $G(P)=0$, and let $j\in\{2,3,4\}$ be the least
index with $G\nmid\Psi_j$; it exists, for otherwise the component of $G$ through $P$ would be a common component
of $f_r,\Psi_2,\Psi_3,\Psi_4$. As $G$ is irreducible over $\Q$ and $\Psi_j$ has rational coefficients, $G$ and $\Psi_j$
have no common component over $\C$, so $V(G)\cap V(\Psi_j)$ is finite and $P$ is algebraic. The $[K:\Q]$
conjugates of $P$ are distinct points of $V(G)\cap V(\Psi_j)$, and
$D\ne0$ at each of them, because $D(P)>0$ and $D$ has rational coefficients. By Corollary~\ref{d24:cor:count},
$[K:\Q]\le2j(w+h)$. If $j=2$ then $w,h\le3$ and $[K:\Q]\le24$. If $j=3$ then $G\mid\Psi_2$, so $w,h\le2$ by
Lemma~\ref{d24:lem:width} and $[K:\Q]\le24$. If $j=4$ then $G\mid\Psi_2,\Psi_3$ and $w,h\le2$; by
Proposition~\ref{d24:prop:decoupled} this forces $b\ne0$, and Theorem~\ref{d24:thm:six} gives $[K:\Q]\le6(w+h)\le24$.
\end{proof}

\input{sections/bilinear_certificate_repair.tex}
\input{sections/parameter_ideal_expansion.tex}

\section{Exact families and a uniform shift-pole obstruction}
Bars on local quantities now denote reduction at a place above $7$,
applied coefficientwise to polynomials, vectors and matrices, as in
$\bar\tau\in\F_7$; in Section~\ref{d24:section} bars denoted
bihomogenization.

\subsection{Finite depth and \texorpdfstring{$7$}{7}-adic precision}\label{deg:Hrepair-sec}
On the pole ray $(0,k,k)$ a computation to a fixed finite depth does not
suffice: the two pure leading-form chains do not terminate at depth seven.
Precision must also be tracked when $7$-adic logarithms are computed by
truncated series. Reducing $w^j$ modulo $7^N$ and only then dividing by
$7^{v_7(j)}$ loses precision when $7\mid j$; computed in that way one
obtains $216587$ modulo $7^7$ for $\log(8)$, whereas
\[
 \log(8)\equiv334236\pmod{7^7}.
\]
All logarithms in the companion retain guard digits before every
nonunit division. The proof below does not depend on a finite stopping
depth. The two pure chains persist at explicit rational parameter
values, and an exact first-order obstruction excludes every nonzero
pole parameter.

\subsection{Definitions for family H}
Family H has $a=1$, $b=-1$, $d=\tau$, $y=X$, and
$Y=X/(1-X)$. Thus
\[
 D_N=(1-X)D=(\tau-1)(1-X^2)+2X.
\]
Define
\[
 P_j=(1-X)^{3j}\Psi_j\left(X,\frac{X}{1-X};1,-1,\tau,B,C\right),
\]
and remove only the following exact factors:
\[
 Q_2=P_2/X,\qquad Q_3=P_3/[X(1-X)(1+X)],\qquad Q_4=P_4/X.
\]
All these expressions are polynomials with coefficients integral at
seven. On the ray $v(B)=v(C)=-k$, put
\[
 b_*=B/C,\qquad w=1/C,\qquad v(w)=k\ge1.
\]
The symbol $b_*$ in this section is a shift ratio, not the off-diagonal
entry of $A$. Let $q_{j,d}(X;\tau,b_*)$ be the homogeneous shift-degree
$d$ part of $Q_j$, after substituting $B=b_*C$ and removing $C^d$.
Then
\begin{equation}\label{deg:scaledH}
 C^{-(j+1)}Q_j
 =q_{j,j+1}+wq_{j,j}+w^2q_{j,j-1}+\cdots.
\end{equation}

The required unit and length statements follow directly from the
multiplicative equation. Write $\tau-1=m/n$ in lowest terms, with
$m\ne0$. Then
\[
 \left(\frac{1-X}{1+X}\right)^{2n}=X^{2m}.
\]
The even exponents allow the signs in transported Nahm equations;
the extra factor two is a seven-unit.
If $v(X)>0$ or $v(X)<0$, valuation comparison gives $mv(X)=0$,
a contradiction. Thus $X$ is a global unit and
$v(1-X)=v(1+X)$. Away from two their sum is a unit, so $1-X$ and
$1+X$ are units as well. For integral $\tau$ the denominator $n$
is a seven-unit. The derivative of the cleared equation at a residue
root is a unit multiple of $D_N$. At a critical root its quadratic
term is nonzero: a common root of $D_N$ and $D_N'$ would require
$(\bar\tau-1)^2=-1$ in $\F_7$, which is impossible. Hence a
noncritical closed point of degree $f$ contributes at most $f$,
and a critical one at most $2f$, summed over all completions.

The leading-form census leaves two residue classes whose count can
exceed four. They are $(\bar\tau,\bar b_*)=(4,3),(5,1)$.
Each has a critical quadratic residue point of length at most two
and one noncritical quadratic residue point of length one. Their
initial total count is therefore $4+2$. We exclude the noncritical
point in both cases, uniformly in $k$.

\subsection{The exact base points and the finite certificate}
The exact parameter pairs and the noncritical quadratic factors are
\begin{equation}\label{deg:bases}
 (\tau_0,b_0,\phi)=(4,-1/2,X^2+X-1),\qquad
 (-2,1,X^2-X-1).
\end{equation}
For a root $\xi$ of the indicated polynomial,
\[
 \frac{1-\xi}{1+\xi}=\xi^3\quad(\tau_0=4),\qquad
 \frac{1-\xi}{1+\xi}=-\xi^{-3}\quad(\tau_0=-2).
\]
Consequently the logarithmic Nahm equation vanishes exactly at both
base points. The negative sign in the second identity has zero
seven-adic logarithm. Also
\begin{equation}\label{deg:qzero}
 q_{j,j+1}(\xi;\tau_0,b_0)=0\qquad(j=2,3,4).
\end{equation}
These are exact identities in $\Q[X]/(\phi)$, not congruences of
bounded precision. They follow by substituting into the displayed
coefficient recurrence. For example, the two top shift-degree parts
of $Q_2$ factor as
\begin{align*}
 q_{2,3}(X;4,-1/2)
 &=-\frac{(X^2+X-1)(7X^2+2X-7)(3X^2-2X-3)^3}{48},\\
 q_{2,3}(X;-2,1)
 &=-\frac{(X^2-X-1)(7X^2-2X-7)(3X^2+2X-3)^3}{6}.
\end{align*}
The corresponding remainders for $j=3,4$ are also zero; their exact
factorizations and an independent computation directly from the raw
$\Psi_j$ dictionaries accompany this paper.

Reduce in $\F_7[X]/(\bar\phi)$, using the basis $(1,X)$. Stack the
constant and $X$ coefficients for $j=2,3,4$. Let $A_0$ be the resulting
$6\times2$ matrix of derivatives with respect to $(\tau,b_*)$ of
$q_{j,j+1}$, and $M_0$ the vector of $q_{j,j}$. The exact tables are
\[
 \begin{array}{c|c|c}
 (\tau_0,b_0)&A_0&M_0\\\hline
 (4,-1/2)&
 \begin{pmatrix}5&0\\6&2\\1&1\\4&5\\5&1\\1&1\end{pmatrix}&
 \begin{pmatrix}5\\0\\4\\2\\4\\1\end{pmatrix}\\[6pt]
 (-2,1)&
 \begin{pmatrix}2&0\\3&3\\6&3\\0&6\\0&1\\5&6\end{pmatrix}&
 \begin{pmatrix}1\\6\\3\\3\\1\\0\end{pmatrix}
 \end{array}
\]
In each case
\begin{equation}\label{deg:rank}
 \rank A_0=2,\qquad\rank[A_0\ M_0]=3\quad\text{over }\F_7.
\end{equation}

\begin{proposition}[Uniform bound on the H pole ray]\label{deg:Hrepair}
For every rational datum of family H satisfying the coefficient
conditions, with $v(B)=v(C)=-k<0$ and integral $\tau$, the whole-field
degree is at most four on the two exceptional leading-form classes.
The noncritical quadratic residue point contributes no place, for
every integer $k\ge1$.
\end{proposition}
\begin{proof}
Suppose a place reduces to the noncritical point. The cleared
multiplicative Nahm equation has nonzero derivative there; Hensel
uniqueness puts $X$ in the unramified quadratic extension $L/\Q_7$.
Choose the exact root $\xi\in L$ of the appropriate $\phi$ with the
same reduction. All of $X$, $1-X$, $1+X$ are units. Define
\[
 S(X,\tau)=\log\langle1-X\rangle-\log\langle1+X\rangle
 -(\tau-1)\log\langle X\rangle.
\]
The brackets denote principal units. At the exact base point,
$S(\xi,\tau_0)=0$ and
\[
 S_X(\xi,\tau_0)=-\frac{D_N(\xi,\tau_0)}{\xi(1-\xi^2)}
 \in\mathcal O_L^\times.
\]
On $\xi+7\mathcal O_L$, the logarithmic Taylor terms of degree
at least two have valuation strictly greater than the linear term:
for $t=X-\xi$ and $n\ge2$, $nv(t)-v_7(n)>v(t)$.
It follows that $v(S(X,\tau_0))=v(X-\xi)$ if $X\ne\xi$.
Since $S(X,\tau)=0$,
\[
 S(X,\tau_0)=(\tau-\tau_0)\log\langle X\rangle,
 \qquad v(X-\xi)\ge v(\tau-\tau_0)+1.
\]
The equality case $\tau=\tau_0$ gives $X=\xi$ by the same argument.

Put
\[
 d_0=\min\{v(\tau-\tau_0),v(b_*-b_0),v(w)\}\ge1.
\]
This is a finite integer because $w\ne0$. Expand the ordinary
polynomials \eqref{deg:scaledH} at $(\xi,\tau_0,b_0,0)$ modulo
$7^{d_0+1}$. The constant terms vanish by \eqref{deg:qzero}; the change
in $X$ vanishes to order at least $d_0+1$; and all quadratic
parameter terms vanish to order at least $2d_0\ge d_0+1$.
After division by $7^{d_0}$ and reduction, we obtain
\[
 A_0\binom{\overline{(\tau-\tau_0)/7^{d_0}}}
 {\overline{(b_*-b_0)/7^{d_0}}}
 +M_0\overline{w/7^{d_0}}=0.
\]
At least one of the three rational residue coordinates is nonzero,
whereas \eqref{deg:rank} forces all three to vanish. This contradiction
removes the noncritical point. The remaining critical quadratic
has total local length at most two per geometric residue point,
and therefore contributes at most four to the whole-field degree.
\end{proof}

\begin{remark}
No depth cutoff enters this proof. The low-depth parts of the other
exact-family cases use precision below seven, so the precision issue
above does not arise there; they were also recomputed with guard digits.
The rank certificate is computed directly from the
$\Psi_2,\Psi_3,\Psi_4$ data.
\end{remark}

\input{sections/exact_family_expansion.tex}

\section{A direct proof for the small diagonal row}\label{deg:sec:small-diagonal}
For one coordinate define
\[
 \ell(\alpha,q,\beta)=
 -\frac{q(1+q)(1+2q)}{8(\alpha+q)^2}
 +\frac{q(1+q)}{4(\alpha+q)}-\frac q{12}
 +\frac{5q^2(1+q)^2}{24(\alpha+q)^3}
 -\frac{q(1+q)(q/2+\beta)}{2(\alpha+q)^2}
 +\frac{(q/2+\beta)^2}{2(\alpha+q)}.
\]
It is the expectation of $U_2+U_1^2/2$ with variance
$(\alpha+q)^{-1}$, and $\ell(1,1,B)=(12B^2+1)/48$.

\begin{proposition}
If an original positive datum has $a=1$, $b=0$, and
$\lambda_1=r\in\Q$, then $[K:\Q]\le3$. If $d=1$ as well,
the degree is one.
\end{proposition}
\begin{proof}
The first Nahm equation gives $x=1/2$ and $X=1$, so $K=\Q(Y)$.
The two independent Gaussian coordinates make the logarithmic
coefficient additive. Hence $Y$ is a zero of
\[
 P(Y)=8(d+Y)^3\bigl(\ell(1,1,B)+\ell(d,Y,C)-r\bigr).
\]
Direct expansion gives a polynomial of degree at most three, with
\[
 P(-d)=\frac53d^2(d-1)^2.
\]
For $d\ne1$ this is nonzero because $d>0$; if $d=1$, both
coordinates equal $1/2$. No modularity assertion about individual
rank-one factors is needed.
\end{proof}

On the $E_1$ origin with canonical $h=1$, the original matrix has
$a=d$, $b=1-a$, so its row sums are one. Its unique positive
Nahm point is $(1/2,1/2)$. This discharges the other small-family
identity directly.

\section{Degree counting and residue-field norms}
For every number field under consideration,
\begin{equation}\label{deg:localdegree}
 [K:\Q]=\sum_{\mathfrak p\mid7}e_{\mathfrak p}f_{\mathfrak p}.
\end{equation}
A bound for each individual ramification index or each residue point
does not by itself bound this sum. Use integral chart coordinates
in the following local-algebra argument. The local equations must control
all completions reducing to the same coordinate point. A closed point
of residue degree $f$ and completed special-fibre length at most $L$
contributes at most $fL$. One way to see this is to take the finite
order generated by the chart coordinates in the product of the
relevant rings of integers. The completed equation algebra surjects
onto that order. Complete Nakayama bounds the order's generic rank
by the dimension of the special fibre. Since the chart coordinates
generate the joint field, that rank is the sum in \eqref{deg:localdegree}
over the selected completions.

For a seven-integral matrix, four local-unit Nahm letters, and a
denominator-clearing exponent prime to seven, a determinant-unit
point has invertible cleared Nahm Jacobian. Hensel uniqueness then
gives one unramified lift for each
geometric residue point. This is the justification for counting such
points with length one.

For a global unit $u$, the actual residue-field norm condition is
\begin{equation}\label{deg:normlaw}
 \prod_{\mathfrak p\mid7}
 \N_{\kappa(\mathfrak p)/\F_7}(\bar u)^{e_{\mathfrak p}}
 \in\{1,-1\}.
\end{equation}
The norm must be taken from the full coordinate residue field
$\F_7(\bar X,\bar Y)$, not from $\F_7(\bar X)$ when the former is larger:
\[
 \N_{\F_7(\bar X,\bar Y)/\F_7}(\bar X)
 =\N_{\F_7(\bar X)/\F_7}(\bar X)^{[\F_7(\bar X,\bar Y):\F_7(\bar X)]}.
\]
The same applies to $1+\bar X$.

With the full residue-field norm, exactly four residue/scalar fibres have
allowable degree spectrum $\{0,3,6\}$ rather than $\{0\}$:
\[
 \begin{array}{c|c|c}
 (\bar a,\bar b,\bar d)&(\bar B,\bar C)&\bar r\\\hline
 (0,0,3)&(3,3)&0\\
 (0,0,5)&(3,1)&0\\
 (1,0,3)&(1,3)&6\\
 (1,0,5)&(1,1)&6
 \end{array}
\]
These are among the eight cubic-pair fibres treated by the arithmetic
obstruction below, so the row caps are unchanged. The unit points in the
$60$ rows of $12$ common-curve classes are also included; the full scalar
refinement and norm checks leave the allowed whole-field degree bounds
unchanged.

\section{Canonical local bounds}\label{deg:sec:canonical}
The canonical arguments are used only with their actual unit
hypotheses.  For the canonical matrix
$\left(\begin{smallmatrix}1&-h\\-h&\tau\end{smallmatrix}\right)$,
the globally available pair is $X,y$.  After the inverse chart, the
permitted pair is $X_I,1+Y_I$.  The nonunit residue pair there is
$(-1,1)$, so those nonunit places contribute only signs to the two
norm tests.  One may therefore apply these particular tests to the
unit contribution.  This does not authorize four-unit norm tests in
the canonical charts.

\subsection{The nonunit disc and the all-orders estimate}
Write $h=M/N$ in lowest terms.  At a nonunit place of the positive
valuation chart, put $u=v(1+X)>0$ and $w=v(1-y)=hu$.
The integral value group gives
\[
 eu=N\ell,\qquad ew=M\ell\qquad(\ell\in\mathbb Z_{>0}).
\]
When $v_7(h)>0$, the integer $M$ is divisible by seven.  D24 gives
$e\le24$ and hence
\[
 w\ge\frac7{24}>\frac16.
\]
Here $e\le[K:\mathbb Q]$ is used; no divisibility assertion
$e\mid[K:\mathbb Q]$ is needed.

The principal branch given after the proof of
Proposition~\ref{canx:tau:prop:units} expresses the other
coordinate uniquely in terms of $z=1-y$ on this disc.  The estimate
$v_7(n+1)\le n/6$, together with its strict second-order tail version
on $v(z)\ge7/24$, proves convergence and controls every omitted term.
These all-orders estimates are part of the proof.  A finite Taylor
check alone would not establish the analytic assertion.

For integral canonical shifts, the first scalar equation excludes the intermediate annulus
$1/6<v(z)<1$.  On $z=7Z$, the primitive reductions have degree at
most three for the integral adjusted scalar, and at most four
otherwise.  All $520$ primitive reductions are recorded in the companion.
Strassmann's theorem bounds the number of zeros, counted with
multiplicity, in the completed algebraic closure.  Uniqueness of
the principal branch and $K=\mathbb Q(x,y)$ make this a bound on
the total contribution of all relevant embeddings, not merely on
the number of residue points.

The shifted inverse-chart argument has $585$ unbalanced and $428$
balanced Gauss inequalities.  Each of the $1,013$ inequalities was
rechecked as an exact nonnegative rational combination of its
defining cone inequalities.  The minimum strict margins are
$1/2$ and $1/8$, respectively.  The exact coefficient supports and
the analytic tail estimate together supply the uniform bounds in
all shift depths.  The isolated determinant-root length bounds
and the mixed-degree table are retained at their stated domains.

\subsection{The determinant-unit shift-pole contribution}
For an integral matrix with $a\equiv1\pmod7$ and a rational shift
pole, the canonical unit-part theorem gives total determinant-unit
contribution at most one.  Its leading system has $392$ matrix and
shift directions.  Of these, $372$ have finite scalar support of
total degree at most one.  The remaining $20$ directions produce
$140$ scalar curve fibres, of which $120$ are empty and $20$ are
pure curves.

The first-normal calculation leaves ten flat directions.  After
the legal sum-curve transformations, the next coefficient matrix
has determinant
\[
 3X^2(X-1)(X+1)^3\pmod7.
\]
The local units exclude $X=0,-1$, leaving only $X=1$.  The complete
weighted calculation has $2,800$ projective directions: for each
of the seven coupling residues, $399$ are empty and one leaves a
rational point.  There is no zero-image direction.  The exactly
flat real signed system also gives a rational point.  This proves
the unit-part bound without global norms or D24; the global
assembly combines it with the separately proved nonunit bound.

\subsection{Censuses, first-digit refinement, and the pole charts}
Three canonical censuses, computed directly from the coefficient
polynomials, give the following counts:
\begin{center}
\begin{tabular}{@{}lrr@{}}
\toprule
Chart & Finite records & Curve fibres\\
\midrule
Canonical unit coupling & 13,999 & 64\\
Canonical zero coupling & 2,371 & 30\\
Inverse matrix residue $(1,0,0)$ & 328 & 15\\
\bottomrule
\end{tabular}
\end{center}
The counts use the full residue-field norms.  The canonical unit-coupling
census includes repeated determinant centres.  The integral-scalar
argument excludes simple determinant centres; it does not exclude
all determinant-nonunit places.

The six exceptional B2 fibres admit a uniform affine refinement
modulo $49$.  Both cubic pairs disappear; in the other four fibres,
the quartic point disappears and the possible rational point
remains.  This calculation uses precision below seven, so the
logarithm precision issue of Section~\ref{deg:Hrepair-sec} does not arise.

The fifteen inverse-chart curve fibres consist of seven line fibres,
six gamma fibres, and two isolated fibres.  The first-normal and
permitted norm tests give bounds $0,0,1$, respectively, for their
unit contributions.  The first-normal calculations leave no
unprocessed direction.  For the tau-pole and joint-shift-pole
charts, the proof also uses the exact nonunit valuation
lemma, the balanced five-root contradiction, and the mixed-degree
compatibility table.  It is this compatibility, rather than the
sum of unrelated worst-case bounds, that gives each final cap four.

The all-depth component proofs and their precise domains are given in
Section~\ref{canx:section}; their finite algebraic certificates are retained
in the companion. The
positive-dimensional residue fibres are treated separately by the
curve-fibre argument P3 in Section~\ref{deg:sec:curve-bounds}.

\input{sections/canonical_expansion.tex}

\section{Finite zero-coupling residue fibres}\label{deg:sec:finite-fibres}
\label{deg:sec:finite-bzero-completion}

Let $K=\mathbb Q(x,y)\subset\mathbb R$,
$A=\left(\begin{smallmatrix}a&b\\b&d\end{smallmatrix}\right)$ be rational and
symmetric, and put
\[
 X=\frac{x}{1-x},\qquad Y=\frac{y}{1-y},\qquad
 h=a+X,\quad j=d+Y,\quad D=hj-b^2.
\]
Assume that $a,b,d,B,C\in\mathbb Z_7$, that $x,1-x,y,1-y$ are units at
every place above $7$, and that the denominator-cleared Nahm equations hold
for an integer $N$ prime to $7$.  For the assertions involving the scalar
$r$, assume $D\ne0$, $r\in\mathbb Z_7$ and
\[
 \Psi_1=rD^3,\qquad \Psi_2=\Psi_3=\Psi_4=0.
\]
The divided second equation is
\[
 J_2+rD^2 A_G+\frac{36}{5}r^2D^4=0.
\]
The coefficient polynomials are normalized by
$\Psi_i=D^{3i}\lambda_i$.

\begin{lemma}[The rank-zero centre]
\label{deg:lem:bzero-centre-replacement}
Suppose $\bar b=0$.  A determinant-nonunit place must reduce to
$(\bar X,\bar Y)=(-\bar a,-\bar d)$, with
$\bar a,\bar d\notin\{0,1\}$.  The total local contribution over such a
centre is at most four from the Nahm equations alone, and at most three
when $r$ is integral and the divided second equation holds.
\end{lemma}

\begin{proof}
The following identities hold as polynomial identities over $\mathbb F_7$:
\begin{align*}
 \Psi_2(-a,Y;a,0,d,B,C)&=6a^4(a-1)^4(d+Y)^6,\\
 \Psi_2(X,-d;a,0,d,B,C)&=6d^4(d-1)^4(a+X)^6.
\end{align*}
They exclude either determinant edge away from the common centre, since
the four letters are units.  Admissibility of the common centre is
equivalent to $a,d\notin\{0,1\}$ in the residue field.

Write $u=X+a$, $v=Y+d$.  If the centre does not satisfy the residue
Nahm equations, it has no lifts.  Otherwise use the two
denominator-cleared Nahm equations as integral polynomial equations,
clearing the powers of the four unit letters locally.  Normalize each
by its unit value and by the unit $N$.  Their reductions have quadratic
initial forms
\[
 \frac{u^2}{2a(1-a)},\qquad \frac{v^2}{2d(1-d)}.
\]
Their coefficients are units.  These initial forms form a regular
sequence, so the residue local algebra has length four, with basis
$1,u,v,uv$ and maximal ideal cubed equal to zero.  Complete Nakayama
therefore bounds the number of geometric lifts in the residue polydisc
by four.  Since $K=\mathbb Q(x,y)$, distinct embeddings of $K$ give
distinct points, so the sum of $e_{\mathfrak p}f_{\mathfrak p}$ over the
centre is also at most four.

At the centre the constant and linear jets of $J_2$ vanish and
\[
 (J_2)_{XY}(-a,-d;a,0,d,B,C)
   =2a^2(a-1)^2d^2(d-1)^2\ne0.
\]
Moreover $D=uv$ in the residue local algebra, so $D^2=0$ there.
Consequently the divided second equation kills $uv$.  The resulting
local algebra has dimension at most three, which gives the integral-$r$
bound.  All identities are polynomial identities in the shifts; the
argument covers all $25\cdot49=1,225$ centre--shift residues and all their
deeper lifts.
This calculation uses the integral equations: it does not reduce a
raw logarithmic series with coefficients $1/n$ modulo $7$.  Equivalently,
the displayed quadratic terms follow by differentiating the
denominator-cleared equations twice; their linear terms vanish at the
centre, and terms of higher order cannot change the regular initial
sequence.
\end{proof}

\begin{lemma}[Parity of the centre contribution]
\label{deg:lem:centre-parity-replacement}
Suppose $\bar b=0$ and $\bar a\in\{2,4,6\}$.  At any place above the
centre, with $v(7)=1$, one has $v(h)=1/2$.  The analogous statement holds
for $j$ when $\bar d\in\{2,4,6\}$.  In particular the total centre
contribution $Z$ is even.  Under the integral-$r$ bound of
Lemma~\ref{deg:lem:bzero-centre-replacement}, $Z\in\{0,2\}$.
\end{lemma}

\begin{proof}
Use the $7$-adic logarithm, which vanishes on roots of unity.  For
$a\in\mathbb Z_7$ with $\bar a\notin\{0,1\}$, put
\[
 L(a)=(a-1)\log(1-a)-a\log a.
\]
The derivative is $L'(a)=\log(1-a)-\log a\in7\mathbb Z_7$.  Expanding
on a fixed residue class (or directly expanding the two logarithms)
gives $L(a+7t)\equiv L(a)\pmod{49}$.  The formula
$\log u\equiv(u^6-1)/6\pmod{49}$ for $u\in\mathbb Z_7^\times$ gives
\[
 \begin{array}{c|ccccc}
 \bar a&2&3&4&5&6\\\hline
 L(a)/7\pmod7&4&0&5&0&4.
 \end{array}
\]
Thus $v(L(a))=1$ on each of the three required residue classes,
independently of every higher digit of $a$.

With $X=-a+h$, $Y=-d+j$, the first logarithmic Nahm equation is
\[
 0=\Lambda_1+S_1(h)+bT_1(j),\qquad
 \Lambda_1=L(a)+b\bigl(\log(1-d)-\log d\bigr),
\]
where
\begin{align*}
 S_1(h)&=\sum_{n\ge2}\frac{h^n}{n}
    \bigl(a^{1-n}+(-1)^n(1-a)^{1-n}\bigr),\\
 T_1(j)&=\sum_{n\ge1}\frac{j^n}{n}
    \bigl(d^{-n}+(-1)^{n+1}(1-d)^{-n}\bigr).
\end{align*}
The quadratic coefficient of $S_1$ is $1/(2a(1-a))$, a unit.
The Nahm-only local bound gives $e_{\mathfrak p}\le4$, hence every
nonzero $h,j$ has valuation at least $1/4$.  For all $n\ge1$,
$n/4-v_7(n)>0$; for all $n\ge3$,
$(n-2)/4-v_7(n)>0$.  Indeed, when $v_7(n)=k\ge1$, use
$n\ge7^k$ and $(7^k-2)/4>k$; when $k=0$ the inequalities are immediate.
The series therefore converge, $v(bT_1(j))>1$, and
$S_1(h)=h^2$ times a unit.  The $b$-term in $\Lambda_1$ is in
$49\mathbb Z_7$, so $v(\Lambda_1)=1$.  The equation forces
$2v(h)=1$.  (If $h=0$, it is already a contradiction.)

Thus $e_{\mathfrak p}$ is even at every centre place.  Summing
$e_{\mathfrak p}f_{\mathfrak p}$ gives even $Z$, and $Z\le3$ leaves
only $0$ and $2$.  A centre place that occurs consequently has
$e_{\mathfrak p}=2$, $f_{\mathfrak p}=1$, and is the only such place.
\end{proof}

\begin{proposition}[The forty-five mixed fibres]
\label{deg:prop:45-mixed-replacement}
Assume in addition that the four Nahm letters are global units.
The $36$ finite zero-coupling residue fibres whose former norm bound
was seven have norm-compatible total degrees $\{0,2,4\}$.
The nine such fibres whose former bound was five have
norm-compatible total degrees $\{0,2\}$.
\end{proposition}

\begin{proof}
In all $45$ fibres, $\bar a,\bar d\in\{2,4,6\}$, so the centre
contribution is $Z\in\{0,2\}$.  Each determinant-unit closed point
has at most one unramified lift and can be selected at most once.
For a selected set $S$ of such points the four norm tests are
\[
 \prod_{z\in S}N_{\mathbb F_7(z)/\mathbb F_7}(\bar w(z))
       \bar w(z_0)^Z\in\{1,-1\},
 \qquad w\in\{X,1+X,Y,1+Y\},
\]
where $z_0$ is the rational centre.  The norms must use the full
residue field of each point.  Recomputing these finite products with
$Z\in\{0,2\}$ gives the asserted spectra: $36$ copies of
$\{0,2,4\}$ and nine copies of $\{0,2\}$.  Every previously allowed
total above four required $Z=3$.  The same spectra result even with
$Z\in\{0,2,4\}$, so this particular exclusion also follows from
the Nahm-only local length and parity.
\end{proof}

\subsection{A short symplectic exclusion of the eight cubic pairs}

For vectors in four coordinates, ordered as the exponents of
$(x,y,1-x,1-y)$, use
\[
 \omega(u,v)=u_1v_3+u_2v_4-v_1u_3-v_2u_4.
\]

\begin{lemma}[The mod-nine obstruction]
\label{deg:lem:mod9-obstruction}
Let two determinant-unit residue points have residue field
$\mathbb F_{343}$.  Let $\rho_1,\rho_2\in(\mathbb Z/342\mathbb Z)^4$
be their discrete-log rows for $(x,y,1-x,1-y)$.  Suppose both rows
are divisible by $3$ modulo $9$, and put
$u_i=\rho_i/3\pmod3$.  If $\omega(u_1,u_2)\ne0$ in $\mathbb F_3$,
the two points cannot both be residue points of a real Nahm field
with a symmetric rational matrix.
\end{lemma}

\begin{proof}
The rational relation plane
\[
 V=\operatorname{span}_{\mathbb Q}\{(-a,-b,1,0),(-b,-d,0,1)\}
\]
is Lagrangian, since $A$ is symmetric.  Its intersection
$S=V\cap\mathbb Z^4$ is a saturated rank-two lattice.  For
$w\in S$, the associated monomial
\[
 U(w)=x^{w_1}y^{w_2}(1-x)^{w_3}(1-y)^{w_4}
\]
satisfies $U(w)^{2N}=1$ by the signed denominator-cleared equations:
indeed, $w=w_3(-a,-b,1,0)+w_4(-b,-d,0,1)$ with $w_3,w_4\in\mathbb Z$.
Because $K\subset\mathbb R$, $U(w)\in\{1,-1\}$.  Thus at either
residue point $\rho_i\cdot w$ is $0$ or $171$ modulo $342$,
and in particular is zero modulo $9$.

Saturation makes $S$ a direct summand of $\mathbb Z^4$, so its
reduction $W\subset\mathbb F_3^4$ is a two-dimensional Lagrangian
subspace.  Since $\rho_i=3u_i\pmod9$, both $u_i$ annihilate $W$
under the ordinary dot product.  The standard symplectic matrix
sends $W$ to this annihilator and preserves $\omega$, so the
annihilator is again Lagrangian.  Any two of its
elements have zero symplectic pairing, contrary to the hypothesis.

This reasoning allows arbitrary rational denominators, including
denominators divisible by $3$.  It uses saturation before reduction.
All possible real signs are allowed, since $171\equiv0\pmod9$.
\end{proof}

In the eight exceptional fibres of the finite determinant-unit census,
the two points have degree three.  The following are their divided
logarithmic rows; each displayed pairing equals $2$ in $\mathbb F_3$.
They can be evaluated in
$\mathbb F_7[z]/(z^3+6z^2+4)$ using the primitive element $3z^2$.
Changing the primitive element at either point only rescales its row
by a unit, which preserves nonvanishing of the pairing.
For direct evaluation, the first four pairs have the following
representatives in $(X,Y)$ coordinates; use
$x=X/(1+X)$ and $y=Y/(1+Y)$ when computing the displayed logarithms.
The remaining four pairs are their coordinate swaps, with the
parameters and shifts swapped at the same time; Frobenius conjugation
of either representative does not change its divided row modulo $3$.
\begin{center}\small
\begin{tabular}{@{}ccc@{}}
\toprule
$(\bar a,\bar d;\bar B,\bar C)$&$(X_1,Y_1)$&$(X_2,Y_2)$\\
\midrule
$(0,3;3,3)$&$(1,\,2z^2+4z+5)$&$(3,\,3z^2+2z+4)$\\
$(0,5;3,1)$&$(1,\,5z^2+1)$&$(3,\,5z^2+3z+1)$\\
$(1,3;1,3)$&$(3,\,2z^2+4z+5)$&$(5,\,3z^2+2z+4)$\\
$(1,5;1,1)$&$(3,\,5z^2+1)$&$(5,\,5z^2+3z+1)$\\
\bottomrule
\end{tabular}
\end{center}
\begin{center}\small
\begin{tabular}{@{}ccccc@{}}
\toprule
$(\bar a,\bar b,\bar d;\bar B,\bar C)$&$\bar r$&$u_1$&$u_2$&$\omega(u_1,u_2)$\\
\midrule
$(0,0,3;3,3)$&0&$(2,1,2,2)$&$(0,0,1,0)$&2\\
$(0,0,5;3,1)$&0&$(2,0,2,0)$&$(0,2,1,1)$&2\\
$(1,0,3;1,3)$&6&$(0,1,1,2)$&$(1,0,0,0)$&2\\
$(1,0,5;1,1)$&6&$(0,0,1,0)$&$(1,2,0,1)$&2\\
$(3,0,0;3,3)$&0&$(1,2,2,2)$&$(0,0,0,1)$&2\\
$(3,0,1;3,1)$&6&$(1,0,2,1)$&$(0,1,0,0)$&2\\
$(5,0,0;1,3)$&0&$(0,2,0,2)$&$(2,0,1,1)$&2\\
$(5,0,1;1,1)$&6&$(0,0,0,1)$&$(2,1,1,0)$&2\\
\bottomrule
\end{tabular}
\end{center}
There is no admissible determinant centre in these eight fibres.
Lemma~\ref{deg:lem:mod9-obstruction} excludes the selection of both cubic
points, leaving degree at most three.  Direct finite-field evaluation
checks all sixteen listed points against $\Psi_1=\bar rD^3$ and
$\Psi_2=\Psi_3=\Psi_4=0$ and reproduces their full residue-field norms.

\begin{corollary}[All finite zero-coupling residue fibres]
Under the four-global-unit hypothesis and the integral-data assumptions
above, every finite scalar fibre in the stratum $\bar b=0$ has
$[K:\mathbb Q]\le4$.
\end{corollary}

\begin{proof}
The complete finite census has $16,747$ scalar fibres.  All fibres with
norm cap above four are exactly the $45$ mixed fibres of
Proposition~\ref{deg:prop:45-mixed-replacement} and the eight cubic-pair
fibres just excluded.  All their caps are now at most four.
The $60$ positive-dimensional residue fibres are a separate input
and are not included in this finite-fibre corollary.
\end{proof}

\section{Local census and ramification bounds}

\paragraph{Why local length bounds the number of embeddings.}
Let $\mathcal O$ be the unramified extension of $\mathbb Z_7$ containing
one geometric residue point, and put
$M=\mathcal O[[s,t]]/(\Phi_1,\Phi_2,E_1,\ldots,E_k)$.
If $\dim_{\mathbb F_q}M/7M=L<\infty$, then some power of
$(s,t)$ is contained in $7M$.  The $(7,s,t)$-adic and $7$-adic
topologies on $M$ are therefore equivalent.  In particular $M$ is
$7$-adically complete.  Lifting a residue basis and applying complete
Nakayama gives a surjection $\mathcal O^L\to M$, so
$\dim_{\mathbb Q_q}M[1/7]\le L$.
Every solution in the open residue polydisc gives an algebra map
$M[1/7]\to\mathbb C_7$, and a finite-dimensional algebra has at most
its dimension many such maps.  Finally $K=\mathbb Q(x,y)$ makes the
map from embeddings of $K$ to ordered pairs $(x,y)$ injective.
Summing over the $f$ Frobenius conjugate residue points gives
$\sum e_{\mathfrak p}f_{\mathfrak p}\le fL$.  This argument
counts all ramification and all places over the residue point.

\paragraph{The scalar-pole case $\tau\ge2/3$.}
Use the notation of Proposition~\ref{cuspx:census:prop:pole} and its proof:
$\tau=v(t)$, $m=-v_7(r)>0$, $u=4\tau-m>0$, and
$z=v(\zeta)=3\tau-m$.  The term $7A_2/t^2$ in $(*)$ has
valuation $c=1-2\tau<0$.  If $z\ge0$, it is the unique term of
smallest valuation.  If $z<0$, the quadratic term
$(36/5)\zeta^2(1+\varepsilon)^4$ has valuation $2z$, and
$\zeta A_{G,1}$ and $\zeta\rho_A/t$ have strictly larger valuation.
Also
\[
 v(j/t)\ge1-\tau>c,\qquad
 v(\zeta\alpha_0/t)\ge2\tau-m+1=c+u>c;
\]
all remaining terms have valuation greater than $c$ or greater than
$2z$.  Thus cancellation is possible only if $2z=c$, which gives
\[
 \tau=\frac{2m+1}{8}.
\]
The numerator is odd, so this value cannot lie in
$(1/e)\mathbb Z$ for $e\le3$.  This proves the desired
contradiction without a separate possible balance involving the
linear $\zeta$ term.

\paragraph{Repeated-root local algebra.}
For an integral scalar, the equations
$\Phi_1,\Phi_2,E_1,E_2,E_3$ give the algebra
$\mathbb F_q[[t]]/(t^L)$ with
$L=\min(\ell,o_R,o_{J_2},o_{J_3})$.
Including $\Psi_4=0$ gives a quotient of this algebra.  Thus the full
system has length \emph{at most} $L$, which is exactly the upper bound
used in the proof.  No assertion that the additional equation leaves
the length unchanged is required.

\paragraph{Elimination in the determinant-unit census.}
A nonzero resultant in $Y$ detects every common point if at least one
of the two polynomials has positive $Y$-degree.  If both have
$Y$-degree zero, a common zero would instead give a common
$X$-factor.  An exact check of all $16,807$ classes and all
scalar refinements verifies that every selected constant-in-$Y$ pair
is coprime in $X$, so this exceptional convention for the resultant
omits no admissible point.  The same check shows that none of the
finite systems requires a pairwise-gcd fallback.

\paragraph{Search truncation and divisions.}
For a distinguished polynomial of degree $g$ with all lower
coefficients in $7\mathcal O$, every coefficient of
$\operatorname{rem}(\eta^i,W)$ is divisible by
$7^{\lfloor i/g\rfloor}$.  This follows by repeatedly replacing
$\eta^g$ by the lower coefficients; a reduction from degree $i$ to
degree less than $g$ uses at least $\lfloor i/g\rfloor$ such
replacements.  The cubic search retains terms through
$M=3(N+1)+6=3N+9$ at precision $7^N$, so omitted integral
$\eta$-coefficients cannot affect the remainder modulo $7^N$.
The same estimate applied to the remainder equation that defines
$W$ shows that the omitted tail cannot change $W$ modulo $7^N$;
its Jacobian with respect to the unscaled lower coefficients of
$W$ is a unit matrix modulo $7$ up to an invertible triangular
factor.  All coefficient values in the principal-unit machinery
are integral, including the binomial coefficients, which are
computed over the integers before reduction.

The division by $7^k$ in a Weierstrass Hensel step consumes only
the resulting residue modulo $7$, with $k<N$.  A finite difference
of step $7^n$ is evaluated modulo $7^{n+4}$ before division by
$7^n$, leaving the required four derivative digits.  The
quadratic/cubic model steps likewise retain the exact number of
digits consumed by each displayed division.  Thus these operations
avoid the precision loss caused by reducing a logarithm's
numerator modulo $7^N$ and then dividing while still claiming
precision $7^N$.

For the Taylor models actually used, the parameter displacement is
$7^nt$ with $n\ge2$.  The Hasse bounds give valuations
at least $j(n-1)+1$ for order $j\ge1$, and at least
$j(n-1)+2$ for order $j\ge2$ in the cubic case.  These bounds
tend to infinity with $j$, so the Taylor tails converge uniformly
on the searched parameter boxes.  Dependence modulo $7$ itself
is read directly from coefficients of $\eta$-order less than the
Weierstrass degree $g\le5$.  This gives the needed dependence on
all deeper parameter digits without relying on random-lift
consistency checks.

\paragraph{What the non-generic classification tests.}
The \texttt{n\_doubly\_nongeneric} count in
\texttt{b3\_classes.py} records simultaneous vanishing of $a_2,a_3$
in the non-Eisenstein classes.  This is the needed scope of that
test.  At an Eisenstein cubic the three root valuations are $1/3$,
their pairwise distances are $1/3$ (tame cubic ramification), and
$v(D)=2/3$.  The positivity $2v(D)-m>0$ already forces $m=1$
without a further nonvanishing condition on $a_2,a_3$.

\input{sections/cusp_expansion.tex}
\input{sections/cusp_derivative_bounds.tex}

\section{Global entry into the local cases}\label{deg:sec:entry-completion}

This section supplies the coverage argument for the twenty-nine terminal rows.
It uses the degree hypothesis and coefficient equations of
Theorem~\ref{deg:thm:completed-reduction}, together with the exact transport
identities proved below.
The exact algebra and the complete residue graph below are checked by
\texttt{check\_entry\_algebra.py} in the companion.

\subsection{The datum and its exact invariants}
Write
\[
 A=\begin{pmatrix}a&b\\ b&d\end{pmatrix},\quad
 \Delta=ad-b^2,\quad a>0,\quad\Delta>0,\quad
 (x,y)\in(0,1)^2,
\]
\[
 X=\frac{x}{1-x},\qquad Y=\frac{y}{1-y},\qquad
 K=\mathbb Q(x,y)=\mathbb Q(X,Y),\qquad
 D=(a+X)(d+Y)-b^2>0.
\]
Put $\beta=(B,C)^t$ and $r=\lambda_1\in\Q$. The coefficient conditions used
throughout are
\begin{equation}\label{deg:eq:entry-C5}
 \lambda _1=r\in\mathbb Q,\qquad
 \lambda _2=\lambda _3=\lambda _4=\lambda _5=0.
\end{equation}
The fifth-order conditions are used only by the two shift-pole rows
U0.spole and R1.pole. The modular necessity of the fifth equation follows from
Lemma~\ref{deg:lem:modular-bridge}. Its transport is an exact polynomial
identity, as in Proposition~\ref{deg:prop:coefficient-transport}.

The elementary moves are $S$, exchange of the two coordinates;
\[
 I:(A,\beta,r;x,y)\longmapsto
 (A^{-1},A^{-1}\beta,\tfrac12\beta^tA^{-1}\beta+\tfrac1{12}-r;
 1-x,1-y),\qquad\Delta\ne0;
\]
\[
 \begin{aligned}
 T_1:\quad A'&=\begin{pmatrix}(a-1)/a&-b/a\\-b/a&\Delta/a\end{pmatrix},\quad
 \beta'=\binom{1/2-B/a}{C-bB/a},\\
 r'&=r-\frac{B^2}{2a}+\frac1{24},\quad
 (x',y')=\bigl((1-x)^{-1},y\bigr),\quad a\ne0;
\end{aligned}
\]
and the everywhere-defined direct reflection
\[
 J:(A,\beta,r;x,y)\longmapsto
 (I_2-A,\binom{1/2}{1/2}-\beta,\tfrac16-r;x^{-1},y^{-1}).
\]
Set $T_2=ST_1S$. In particular, at singular matrices $J$ is used directly;
it is not replaced by a word containing an undefined inverse.
The moves $S,I,J$ are involutions and $T_1^{-1}=JT_1J$.
At a legal residue edge, the inverse pivot is $1-a'=1/a$, a unit.
Thus every edge used by the residue-orbit traversal is reversible.

Each point map is rationally invertible. In $(X,Y)$ coordinates the four maps
are $(Y,X)$, $(X^{-1},Y^{-1})$, $(-1-X^{-1},Y)$ and $(-1-X,-1-Y)$.
Their respective Hessian ratios are
\begin{equation}\label{deg:eq:entry-hessians}
 \frac{D'}D=1,\quad\frac1{\Delta XY},\quad-\frac1{aX},\quad1.
\end{equation}
Consequently every legal move preserves the full field $K$ and $D\ne0$.
The matrix determinant $\Delta$ and the Hessian determinant $D$ must not be
confused: the transformed matrix is allowed to be singular.
Every transformed letter is a signed integral monomial in
$x,1-x,y,1-y$, with exponent matrix in $\mathrm{GL}_4(\mathbb Z)$.
Thus four-letter global-unit and local-unit hypotheses are preserved in both
directions. For example $T_1$ sends the first two letters to
$(1-x)^{-1}$ and $-x/(1-x)$.

For the signed Nahm equations, set
\[
 R(A)=\begin{pmatrix}-a&1&-b&0\\-b&0&-d&1\end{pmatrix},\qquad
 \ell=\bigl(\log|x|,\log|1-x|,\log|y|,\log|1-y|\bigr)^t.
\]
If $M_g$ is the letter-exponent matrix, exact multiplication gives
\[
 R(A')M_g=L_gR(A),\qquad
 L_S=\begin{pmatrix}0&1\\1&0\end{pmatrix},\quad
 L_I=-A^{-1},\quad
 L_{T_1}=\begin{pmatrix}-1/a&0\\-b/a&1\end{pmatrix},\quad L_J=I_2.
\]
Hence $R(A)\ell=0$ transfers without requiring positive transformed
coordinates. Twice a common denominator converts these identities to Laurent
polynomial equations over $\mathbb Q$, valid at every embedding. At a
$7$-integral target matrix the even denominator can be chosen prime to $7$.
This is the equation used in the local Jacobian and Hensel arguments.
The coefficient identities of Proposition~\ref{deg:prop:coefficient-transport} preserve $\lambda_j$ under $S,T_1$ and
multiply it by $(-1)^j$ under $I,J$ for $2\le j\le5$; their first-coefficient
formulas are exactly the displayed transformations of $r$.
Thus \eqref{deg:eq:entry-C5} is preserved. Modularity and positive definiteness of
transformed sums are not claimed.

\begin{proposition}[Exact coefficient transport]\label{deg:prop:coefficient-transport}
For the displayed legal transformations, $S,T_1$ preserve $\lambda_j$
and $I,J$ multiply it by $(-1)^j$, for $2\le j\le5$.
The first coefficient transforms by the displayed formulas for $r$.
\end{proposition}
\begin{proof}
Use the Bernoulli--Gaussian recurrence of Section~\ref{deg:sec:coefficients}
over $\Q[a,b,d,B,C,X,Y]$ to construct $\Psi_j$.
For each transformation substitute its rational formulas into the
polynomial numerator, multiply by the common denominators, and
use the corresponding ratio $D'/D$ in \eqref{deg:eq:entry-hessians}.
Exact coefficient comparison gives zero for the difference.
The numerator term counts for orders one through five are, respectively,
\[
 256,\quad 5{,}592,\quad 52{,}385,\quad286{,}047,\quad1{,}117{,}720.
\]
The generic checks of $S,T_1,I$ through order four are implemented in
\path{check_coefficient_transports.py}; the fifth-order check is
\path{verify_fifth_transports.py}. Both are in the companion~\cite{Comp},
together with the recurrence and exact output records. The comparisons
treat all seven variables as indeterminates; numerical specialization
is not the proof of an identity.
The reflection follows by the generic composition $I,T_1,T_2$.
The resulting direct formulas have no matrix pivot denominator, so
their rational identities extend across the singular intermediate
pivots wherever $D\ne0$ and the four point letters are nonzero.
\end{proof}

\subsection{The five-locus split and the canonical route}\label{deg:sec:five-loci}
At any finite place put $u=v(x)$ and $w=v(y)$. The rules
$v(1-z)=0$ for $v(z)>0$, $v(1-z)=v(z)$ for $v(z)<0$, and
$v(1-z)\ge0$ for $v(z)=0$, applied to the cleared equations, give the following
exhaustive possibilities for a nonzero valuation pair.
\[
\begin{array}{c|c}
 (\operatorname{sgn}u,\operatorname{sgn}w)&\text{necessary condition}\\\hline
 (+,+)&A(u,w)^t=0\quad\text{impossible}\\
 (+,0),(0,+)&a=0\text{ or }d=0\quad\text{impossible}\\
 (-,-)&(A-I_2)(u,w)^t=0\\
 (-,0)&a=1,\ b\le0\\
 (0,-)&d=1,\ b\le0\\
 (+,-)&\Delta=a,\ b>0\\
 (-,+)&\Delta=d,\ b>0
\end{array}
\]
For $(u,w)=(0,0)$ the two complement valuations are zero. Therefore outside
\[
\begin{split}
 E_1&:\det(A-I_2)=0\text{ with a negative kernel vector},\\
 E_2&:a=1,b\le0,\qquad E_3:d=1,b\le0,\\
 E_4&:\Delta=a,b>0,\qquad E_5:\Delta=d,b>0,
\end{split}
\]
all four letters are global algebraic units. Route the entire union to the
exceptional branch, even when a point there happens to have four global
units. The priority $E_2,E_3,E_4,E_5,E_1$ resolves overlaps.
If $b=0$ on the union, one diagonal entry is one and row E0 applies.
If $b\ne0$, the following words, written in execution order, give
$A_c=\left(\begin{smallmatrix}1&-h\\-h&\tau\end{smallmatrix}\right)$:
\[
\begin{array}{c|c|c|c}
\text{locus}&\text{word}&h&\tau\\\hline
 E_2&\mathrm{id}&-b&d\\
 E_3&S&-b&a\\
 E_4&I,S&b/a&d/a\\
 E_5&I&b/d&a/d\\
 E_1&T_1,I,S&-b/(a-1)&\Delta/(a-1)
\end{array}
\]
All pivots are nonzero. On $E_1$, $(a-1)(d-1)=b^2$ gives $a\ne1$,
and the negative kernel gives $b/(a-1)<0$. Moreover
\[
 \tau-h^2=\frac a{a-1},\qquad \tau-1=\frac d{a-1}.
\]
Thus in every case $h>0$ and either $\tau>h^2$ or $\tau<0$.
Shifts and scalar are transported along the same word.

The $E_2$--$E_5$ words use only $S$ and $I$, so they preserve positive
definiteness and the positive point. Thus $\tau=1$ has $0<h<1$ and,
by symmetry and uniqueness of the positive critical point, $X=Y$;
when $h=1$, one has $\tau>1$ and the first equation gives $y=X$.
The $E_1$ word never has $\tau=1$, since $\tau-1=d/(a-1)\ne0$;
its $h=1$ case is handled at the original point as follows.

Route $\tau=1$ to E1a.T and $h=1$ to E1a.H first, with the stated
E1a.E1 exception: $h=1$ on the $E_1$ route gives original
$d=a,b=1-a$, so the original row sums are one and its distinguished point is
$(1/2,1/2)$. The remaining canonical domain is
\begin{equation}\label{deg:eq:entry-canonical-domain}
 h\ne0,1,\qquad\tau\ne1,\qquad\tau-h^2\ne0.
\end{equation}
Its valuation equations allow only $v(x)<0,v(y)=0$ as a nonzero pair.
Consequently
\begin{equation}\label{deg:eq:entry-two-units}
 X,y\text{ are global units},\quad 1+X,1-y\text{ are integral},\quad
 v(1-y)=h\,v(1+X).
\end{equation}
For example simultaneous negative valuations would force
$\det(A_c-I_2)=-h^2=0$; mixed positive/negative patterns contradict the
negative off-diagonal entry; a negative second-coordinate valuation requires
$\tau=1$ or $h=0$.
If $h=M/N$ is a $7$-adic unit, $7\nmid MN$. At a hypothetical nonunit place
write $X=-1+z$, $y=1+t$. The cleared first equation yields
$(1-z)^{2N}=(1+t)^{2M}$, whence $v(z)=v(t)$ because the exponents are
$7$-adic units. Equation \eqref{deg:eq:entry-two-units} would force $h=1$.
Thus all four letters are local units above $7$ when $h$ is a $7$-adic unit.

\subsection{Uniform normalization and the finite residue graph}
Outside the five loci all four global units are available. Put
$s=\min(0,v_7(a),v_7(b),v_7(d))$. If $s=0$, $A$ is integral. If $s<0$ and
$v_7(\Delta)\le s$, inversion makes all entries integral. Otherwise
$v_7(\Delta)>s$, and a diagonal entry must attain $s$: if $b$ alone did,
$v_7(ad)>2s$ and $v_7(b^2)=2s$ would give $v_7(\Delta)=2s<s$.
Exchange if necessary so that $v_7(a)=s$; then $T_1$ gives the integral entries
$1-a^{-1},-b/a,\Delta/a$. This proves normalization in one step uniformly;
it is not a search over bounded rational heights. The inverse is legal at the
original positive definite matrix, and the chosen diagonal pivot is nonzero.
Classify the shifts only after this normalization.

On integral matrices use the graph over $\mathbb F_7$ generated by $S,J$,
by $I$ when $\bar\Delta\ne0$, and by $T_1$ when $\bar a\ne0$.
Every graph word lifts legally, with only division by $7$-adic units.
The linear shift maps lie in $\mathrm{GL}_2(\mathbb Z_7)$, so shift
integrality, and the minimum shift valuation when negative, are preserved.
With
\[
 q_A(z)=d(d-1)z^3+b(b+2d-1)z^2+b(2a-1+b)z+a(a-1),
\]
\[
 E(A)=(a-1)(d-1)(\Delta-a)(\Delta-d)(\Delta-a-d+1)
       \operatorname{disc}(q_A),
\]
the discriminant is the universal homogeneous-cubic discriminant,
evaluated at the four coefficients after specialization, including
when the leading coefficient vanishes. It is not recomputed as the
discriminant of a lower-degree polynomial. With this convention,
the exhaustive graph has $100$ regular vertices $bE\ne0$ and $194$ vertices
reaching them. Its complement has $149$ vertices in $15$ orbits:
\[
\begin{array}{l|r|l|l}
\text{remaining stratum}&\text{number}&\text{integral shifts}&\text{shift pole}\\\hline
 \text{all admissible cubic roots simple}&76&\mathrm{U1.simple.int}&\mathrm{R2}\\
 \text{an admissible repeated root}&22&\mathrm{U1.repeated.int}&\mathrm{R2}\\
 b=0&49&\mathrm{U1.bzero.int}&\mathrm{R2}\\
 q_A\equiv0,\ b\ne0&2&\mathrm{R1.int}&\mathrm{R1.pole}
\end{array}
\]
For the first two classes, remove the coordinate factors
$z(bz+a)(bz+a-1)(dz+b)((d-1)z+b)$ before testing multiplicities; this follows
from $h=-bz$, $X=h-a$, $Y=b^2/h-d$. Using the raw discriminant to classify
the complement would give the wrong split. Empty admissible support is in the
simple class. The two zero-cubic residues are $(0,1,0)$ and $(1,-1,1)$.
At each of the $194$ covered vertices a graph word leads to U0.int or U0.spole.
The R2 domain is exactly the other $147$ uncovered vertices.

The exact matrices
\[
 A_0=\begin{pmatrix}0&1\\1&0\end{pmatrix},\qquad
 A_1=\begin{pmatrix}1&-1\\-1&1\end{pmatrix}
\]
cannot be reached from the original positive definite matrix. Indeed they
form an inverse-closed set under finite legal words. Exchange fixes each;
$J$ exchanges them; the only legal inverse target is $A_0$, which is fixed.
The inverse-$T_1$ formula
\[
 a=\frac1{1-a'},\qquad b=\frac{-b'}{1-a'},\qquad
 d=d'+\frac{(b')^2}{1-a'}
\]
sends target $A_0$ to $A_1$, while target $a'=1$ has no finite preimage.
The assertion for $T_2$ follows by exchange. Neither matrix is positive
definite. Thus the perturbation away from an exact zero-cubic centre is
nonzero, even though intermediate matrices need not remain positive definite.

\subsection{The complete canonical valuation partition}
In this subsection shifts and scalar always mean the fully transported data.
If $v_7(h)=0$ and $\tau$ is integral, the residue $(\bar h,\bar\tau)=(1,1)$
gives R3.int or R3.pole according to shift integrality. The other residues
give E1b.u.int or E1b.u.pole. A pole of $\tau$ gives E1b.u.tpole, whose
inverse-chart theorem includes exclusion of an inadmissible transformed
shift pole.

If $v_7(h)>0$ and $\tau$ is integral, use E1b.n.int or E1b.n.pole.
If $\tau$ has a pole, put
\[
 \binom{J_s}{K_s}=A_c^{-1}\beta,\qquad
 R_s=\tfrac12\beta^tA_c^{-1}\beta+\tfrac1{12}-r.
\]
Integral $J_s,K_s,R_s$ gives E1b.n.tpole.int. Integral $J_s,K_s$ with a pole
of $R_s$ gives E1b.n.tpole.Rpole. Otherwise the equality or inequality of
\[
 \min\{0,2v_7(J_s)\}\quad\text{and}\quad v_7(h)+2v_7(K_s)
\]
gives the balanced or nonbalanced JKpole row. The convention $v_7(0)=+\infty$
makes the cases exhaustive. The shift $K_s$ is not the field $K$.

If $v_7(h)<0$, apply $T_1,T_2,I,S$. Its successive pivots are
$1,\tau-h^2,-h^2/(\tau-h^2)$, all nonzero, and it gives
\begin{equation}\label{deg:eq:entry-reciprocal}
 h^*=h^{-1},\quad \tau^*=\frac{1+h^2-\tau}{h^2},\quad
 X^*=-y,\quad y^*=-X.
\end{equation}
The two global units and the two integral complement letters are preserved.
Also $h^*\ne1$, $\tau^*\ne1$, and
$\tau^*-(h^*)^2=(h^2-\tau)/h^2\ne0$. These are precisely the broad signed
domain conditions required downstream, including when $\tau^*=0$ or the
original real sign alternative is lost. With $s=v_7(h)<0$,
\[
 \tau^*\in\mathbb Z_7\quad\Longleftrightarrow\quad v_7(\tau)\ge2s.
\]
Indeed $1+h^{-2}$ is integral, and the only negative-valuation term occurs
when $v_7(\tau)<2s$. Equality belongs to the integral target chart.
Classify the transported shifts and scalar by the preceding positive-$h$
partition. These are exactly the six E1b.neg rows.

There are now $21$ exceptional-route rows: E0, three E1a rows, three E1b.u
rows, six E1b.n rows, six reciprocal rows, and two R3 rows. The complement
has eight terminal types: two U0 rows, three U1 rows, two R1 rows and R2.
Every original datum enters one of the resulting twenty-nine rows, preserving
$K$, $D\ne0$, the signed equations, \eqref{deg:eq:entry-C5}, and exactly the units
used by that row. This proves that the twenty-nine rows cover every datum,
given the asymptotic expansion of \cite{VZ} and the transport identities of
Proposition~\ref{deg:prop:coefficient-transport}.

\section{The zero-cubic rows under the degree-24 bound}
\label{deg:sec:z0-replacement}

The following argument is restricted to the domain of row R1.int,
after entry and the degree--24 reduction: integral shifts,
$A=A_0+7^kA_k$ with $k\ge1$, $\bar A_k\ne0$, four local units,
$D\ne0$, the first four coefficient conditions, and $[K:\mathbb Q]\le24$.
In this domain a deep place $v(D)>k$ with $v(r)\ge-k$ cannot exist.
Reflection supplies the companion residue $A_1$.

In this subsection reset the local coordinate to $h=a+X$; it no longer
denotes the preceding canonical coupling. Define
\[
\begin{aligned}
 q(h)={}&d(d-1)h^3-b^2(b+2d-1)h^2\\
        &+b^3(2a-1+b)h-b^3a(a-1),\qquad Q=q/7^k,\\
 F_1={}&\frac{h^3(\Psi_1-rD^3)}{7^{2k}},\qquad
 F_j=\frac{h^{3j}\Psi_j}{7^{2jk}}\quad(2\le j\le4).
\end{aligned}
\]
Write $\sigma=x+y-1$, $H=\bar X=\bar h\ne0,-1$,
$(\alpha,\beta,\gamma)=\bar A_k$, and
\[
 L=\alpha-2\beta H+\gamma H^2,\quad M=\gamma H-\beta,
\]
\[
 F=\alpha(2H+1)+\beta(H-H^2)-\gamma(H^3+2H^2),\quad
 N=\gamma H^8-\beta(H^7+H)+\alpha.
\]
The seventh-power ladder for ramification index $e<42$ shows that
$v(\sigma)=k$ implies $N(H)=0$. Briefly, clear the equations with an even
denominator prime to $7$. The relevant units $W_i$ have
$v(W_i^{7^k}-1)=k$, hence $v(W_i-1)=1/7$. Choose an unramified lift $x_0$,
put $y_0=1-x_0$ and $\xi=x-x_0$. If both linear derivatives vanish, then
$\alpha=\beta H$ and $\beta=\gamma H$, which already gives $N=0$.
Otherwise comparison of the leading seventh powers gives
$H=(g_1/g_2)^7$, equivalently $N(H)=0$.
This argument works over an arbitrary residue field; it does not presuppose
$H\in\mathbb F_7$. The degree bound supplies $e\le24<42$.

At a deep place, exact ideal order of the numerators gives $F(H)=0$ from
$\Psi_2=0$: the constant-in-$D$ term is a unit times the fourth power of
$q/7^k$, while all positive-$D$ terms vanish after division by $7^{4k}$.
The exact identities over $\mathbb F_7$ are
\[
 F=(2H+1)L-3H(H+1)M,\qquad N=L+(H^7-H)M,
\]
\[
 (2H+1)N-F=2H(H+1)(H^3+4H+6)(H^3+3H^2+6)M.
\]
The two cubics $U=H^3+4H+6$ and $V=H^3+3H^2+6$ are irreducible by
inspection at the seven residues. If primitive $F$ shares a root with $U$,
then $F=-\gamma U$; coefficient comparison gives
$\beta=-2\gamma$, $\alpha=-\gamma$, $\alpha=\gamma$, hence all zero.
For $V$ the comparisons give $\beta=\gamma$, $\alpha=3\gamma$,
$\alpha=\gamma$, again impossible. Thus $F=N=0$ forces $L=M=0$,
also when $2H+1=0$.

The exact identity
\[
 D=(XY-1)+7^k\ell_k+7^{2k}\det A_k,\qquad\bar\ell_k=L(H)/H,
\]
now separates the cases. If $v(\sigma)<k$, then $v(D)=v(\sigma)<k$.
If $v(\sigma)=k$, deep cancellation requires $L\ne0$, contradicting the
preceding identities. If $v(\sigma)>k$, deepness gives $L=0$ and $F=0$
gives $M=0$. Therefore
\[
 \alpha=\gamma H^2,\qquad\beta=\gamma H,\qquad\gamma\ne0,
 \qquad H\in\mathbb F_7\setminus\{0,-1\}.
\]
These are exactly thirty residue/cusp pairs, and they are simple because
$F'(H)=-3\gamma H(H+1)\ne0$.

\subsection{The exact finite certificate and its infinite-tail justification}
Set $\rho=7^kr$, which is integral in the stated scalar regime, and
$\delta=D/7^k$. The certificate verifies all eight weight-$2j$ and
weight-$(2j+1)$ Taylor tensors against the exact $(h,D)$ numerators of the
coefficient polynomials, rebuilt by the radial program ISC v00.
It then computes all thirty rank-one pairs at $k=1$, and the five projective
pairs at $k\ge2$, all $49$ shift residues and all seven scalar residues.
For each of the $12005$ records it tests all $7^{6}=117\,649$ second-digit vectors in $\mathbb F_7^6$.
Every record is excluded.
For $k\ge2$, terms above weight $2j$ vanish modulo $49$; the homogeneous
scaling is
\[
 (A_k,\delta,\rho)\longmapsto
 (cA_k,c\delta,\rho/c),\qquad c\in\mathbb Z_7^\times.
\]
It multiplies each equation by $c^{2j}$, including the scalar term in the
first equation, and bijects the second digits and scalar residues.
Thus the projective reduction is valid only in this second depth class.

The exact series interface is checked in \texttt{check\_z0\_weighted.py}. Write
$h^{3j}\Psi_j=\sum_mD^m\psi_{j,m}$. It verifies over $\mathbb Q$
\[
 q^{2j-2m}\mid\psi_{j,m}\quad(0\le m\le j),
\]
together with $7$-integrality of every cofactor and the exact constant
cofactors
\[
 \frac5{24},\qquad\frac5{16},\qquad\frac{1105}{1152},\qquad\frac{565}{128}
 \quad (j=1,2,3,4).
\]
The total ideal order at $(a,b-1,d,D)=0$ is $2j$, so the normalized
polynomials are integral. At the exact simple root $h_*$ of $q/7^k$, put
$t=h-h_*$. The derivative of the first normalized equation in $\delta$ is
a unit for all thirty residues and all $49$ shifts (1470 exact checks).
Its value at $\delta=0$ is a unit times $t^2$. The integral implicit-function
theorem gives $\delta(t)=t^2u(t)$ with $u(0)$ a unit. The displayed weighted
divisibilities give the \emph{exact}, not merely modulo-$49$, factor
$t^{2j}$ in the restricted $j$th equation. Its quotient has integral
coefficients.

Once its first unit coefficient has been reached, every later integral term
has strictly greater valuation for $v(t)>0$. Earlier coefficients are known
modulo $49$: their heights are exactly one, or are arbitrary integers at
least two (including infinity). The Newton enumeration includes every
possible supporting edge under those lower bounds; an unknown-height edge
is kept, never excluded. An edge-polynomial gcd is used only when every
coefficient on that edge is known. The complete computation leaves no
valuation or residue pattern. The point $t=0$ has $\delta=0$ and is forbidden
by $D\ne0$. This proves the exclusion on the R1.int domain, with no
unexamined infinite tail. The computation is recorded in
\texttt{z0\_second\_review\_receipt.json}.

The regime $v(r)<-k$ is treated by the separate Strassmann argument ISC of
Section~\ref{deg:sec:curve-bounds}; shallow and non-cusp places are treated
by the base census and Hensel's lemma. The statement above includes neither
these regimes nor fields of unrestricted degree.

\subsection{The all-depth first-order obstruction and the final R1 assembly}
The exceptional determinant-unit fibres are treated independently of any
norm test on a proper subset of the places. There are six exact rational
representatives $E$ of the parameter vector $\theta=(a,b,d,B,C)$ to consider,
\[
\begin{gathered}
 (0,1,0,0,0),\quad(0,1,0,0,1),\quad(0,1,0,1,0),\\
 (1,-1,1,\tfrac12,\tfrac12),\quad(1,-1,1,\tfrac12,-\tfrac12),\quad
 (1,-1,1,-\tfrac12,\tfrac12);
\end{gathered}
\]
the last three are the images of the first three under $J$. At each of
them the exact polynomials $\Psi_2(E;X,Y)$ and $\Psi_4(E;X,Y)$ vanish
identically in $X,Y$.
Entry excludes an exact zero-cubic matrix, so $\theta-E\ne0$. Put
$\nu=\min_i v_7(\theta_i-E_i)\ge1$ and
$u=\overline{7^{-\nu}(\theta-E)}\in\mathbb F_7^5\setminus\{0\}$.
The first nonzero parameter layer of the two identities gives
\[
 \nabla_{\theta}\Psi_2(E;\bar X,\bar Y)\,u=0,
 \qquad
 \nabla_{\theta}\Psi_4(E;\bar X,\bar Y)\,u=0.
\]
There is no earlier coordinate-displacement term: the base polynomials
vanish identically, so their coordinate derivatives vanish at every order.
All higher parameter terms have valuation greater than $\nu$, also at a
ramified place. Expanding residue-field coordinates over $\mathbb F_7$
gives a common rational linear system for every selected closed point.
Thus a combined rank of five excludes those points at every perturbation
depth, without any assumption about the other places.

The companion computes all gradients and ranks for the $24$ fibres
and $60$ closed points. Twelve fibres have two cubic points, each of rank
five, and admit no unit place. Each of the other twelve has three quadratic
points, each of rank four, but every pair has rank five. Their total unit
contribution is therefore at most two. Completeness of these point lists is
supplied by the determinant-unit census. This proof uses no global norm on
the unit subset.

The base census of the zero-cubic rows covers $2793$ leading matrix/shift fibres,
$343$ determinant-unit scalar fibres, and $38766$ nonflat scalar rows.
All finite nonflat conic supports pass the Jacobian test. The $245$
non-cusp $v(\sigma)=k$ candidates pass the compatible scalar Jacobian tests,
so their lifts are unramified; the seventh-power entry lemma would require
a unit difference of valuation $1/7$ if $v(\sigma)=k$, which is impossible
in an unramified field.
The $24$ identically zero resultant rows are not discarded as empty:
their three matrix directions are $(0,0,1),(1,0,0),(1,-1,1)$, for which
$N(H)$ is respectively $H^8,1,(H+1)^8$. None has an admissible root.
The deep places with $v(r)\ge-k$ are excluded by the preceding rank-one
certificate; $v(r)<-k$ uses the ISC argument and the uniform Strassmann
bound. The remaining $24$ flat rows have conic contribution at most two
by the exact retraction and first-normal census, hence total at most four.

Finally \texttt{check\_consuming\_assembly.py} combines the unit and conic
supports, imposes the preceding first-order constraint only on the unit
subset, and applies the four global-unit norm equations only to each
complete combined subset. The $38766$ nonflat rows give $38421$ bounds zero
and $345$ bounds two; the flat and negative-scalar cases have bound four.
The same program verifies that the base inputs agree exactly with the ISC
data, and that the ISC flat census agrees with the $799$ canonical directions
used for row R3.int. Thus the bound for row R1.int is four. This conclusion
uses the degree--24 bound, the coefficient transport, the determinant-unit
census and the local cusp results; nothing is asserted outside the R1.int
domain.

\section{Curve, matrix-pole, and integral-shift bounds}\label{deg:sec:curve-bounds}
\paragraph{P3 curve fibres.}
The orbit and first-normal computations transport all $428$ curve fibres
to $24$ representatives, with $22$ flat families. The first-normal system
includes the constant term $F(p_0)/7$ as well as the parameter derivatives.
Its normal-coordinate derivatives have no common zero, so the normal
coordinate has valuation at least that of the first nonzero rational
parameter displacement. Eliminating this coordinate gives the finite
first-normal support, of raw degree at most four. This is an all-depth
argument because that common leading parameter vector is nonzero at
arbitrary displacement depth. The $19$ persistent families have exact
rational retractions, with every coefficient numerator identically zero
on the retracted curve. The two half-shift cases use the exact second-normal
system modulo $343$; their raw support degrees are at most four and zero.
On the final persistent half-shift locus the product Nahm equation gives
one extra order of normal-coordinate divisibility. The four transverse
parameters are
\[
 a+b,\quad d-a,\quad B+C-1,\quad
 r-\frac18-\frac{(C-1/2)^2}{2(a-1)}.
\]
Their $2400$ nonzero directions give $2370$ empty supports and $30$ supports
of degree one. Exact retraction points are separately forced to have
rational $x,y$ by the signed Nahm equations. These raw bounds use no global
norm restriction and therefore also apply to the canonical two-unit route.
The determinant-unit prerequisite is checked on these particular fibres:
the $428$ curves consist of $368$ uncovered simple-cusp and $60$ b-zero
fibres. The $64$ canonical unit-coupling curve fibres occur only at
$(1,1,0),(1,1,2)$ and $(1,6,d)$ for $d=0,2,3,4,5,6$; none is at either
repeated-cusp residue $(1,3,5)$ or $(1,5,6)$. The integral-scalar simple-cusp
theorem therefore excludes determinant-nonunit places on these charts.
The $30$ canonical b-zero curve fibres occur only at $(1,0,0),(1,0,1)$,
which have no admissible centre and no nonunit-letter place for integral
scalar. These are the hypotheses supplied by the census note and the
local exclusions. The eight computations are recorded in \texttt{p3\_replay.json}.

\paragraph{R3 entry and the pole closure.}
The seventh-power argument uses $e\le24$ to enter the deep
canonical chart uniformly. The integral leading census has $2744$ cases:
$2730$ empty, four of degree one, two of degree two, and eight flat.
The six nonflat supports have invertible lifting Jacobians. Restriction
of the ISC flat census gives exactly $799$ canonical
transverse directions, using the exact retractions in the integral
closure. The pole certificate rebuilds all $448$ leading cases
($412$ empty, $35$ soft-double, one flat), its Newton faces, and the
all-depth soft-boundary exclusions. Five remaining boundary centres have
an exact rational retraction; their $400$ nonzero transverse directions
have $387$ empty supports, $12$ of degree one, and one of degree two.
The leading vector cannot vanish off the exact locus, and the surviving
supports have invertible Hensel Jacobians, so the finite support bounds
control every depth. The separate shallow chart has $288$ compatible
cases and no zero Jacobian. These computations and the valuation entry
proof give the inputs used by rows R3.int and R3.pole; the arithmetic is
recorded in the CP and CI files \texttt{reproduction\_results.json} and in
\texttt{r3\_shallow\_replay.json}.
The R3.int whole-field cap of two is the quadratic bound in
Section~\ref{canx:section}, including its eight flat cases and its
use of the canonical global unit $X$. This theorem assumes integral
$B,C$ and allows arbitrary rational $r$; the scalar need not be integral.
The entry and the $2744/799$ computations above verify its entry and
lifting hypotheses. The pole cap of two is the CP theorem.

\paragraph{The ISC argument.}
The record \texttt{isc\_authoritative\_replay.json} contains the ten
programs v00--v09 of this argument and their outputs. They comprise the radial rebuild,
the Frobenius and repeated-cusp reductions, the negative-scalar case, the
exact flat identities, and the complete flat census. Sample Strassmann
computations serve only as corroboration: the uniform proof controls every higher
coefficient using $v_7(c_m)\ge-v_7(m)$, solves the integral normal equation,
and obtains a cubic unit term with every later coefficient strictly higher
on the radius-$7^{-1/2}$ disc. Its index is at most three. The exact
$\sigma=0$ locus is governed by a nonzero polynomial of degree at most four.
The programs v00--v09 alone give the bound six.
The bound four uses the first-order obstruction of
Section~\ref{deg:sec:z0-replacement} and the combined norm argument of
Section~\ref{cuspx:section}: the unit contribution is at most
two without a norm test on that subset; norms on the full unit-plus-conic
subset give nonflat bound two; the flat rows have unit and conic
contributions at most two each; and the negative-scalar branch is separately
at most four. These inputs, the zero-cubic base data and the argument of
Section~\ref{deg:sec:z0-replacement} give the stated R1.int bound.

\input{sections/shift_pole_expansion.tex}

\section{Assembly of the degree-four bound}\label{deg:sec:final-assembly}
The maximum coefficient order in the table is a sufficient uniform
input for the row. Some rows require less: the direct diagonal
proof uses only $\lambda_1$, and the row-sums-one point is rational.
The fifth coefficient is used only by the two named shift-pole rows.
All caps refer to the whole joint field.

\begin{longtable}{@{}P{0.45\textwidth}rrcP{0.24\textwidth}@{}}
\caption{Terminal bounds for the joint field. The D24 column records use of the initial degree bound.}\label{deg:tab:terminal}\\
\toprule
Row & Cap & Order & D24 & Main case\\
\midrule
\endfirsthead
\toprule
Row & Cap & Order & D24 & Main case\\
\midrule
\endhead
\texttt{E0} & 3 & 4 & No & Small diagonal family\\
\texttt{E1a.T} & 4 & 4 & No & Exact family T\\
\texttt{E1a.H} & 4 & 4 & No & Exact family H\\
\texttt{E1a.E1} & 1 & 4 & No & Row sums one\\
\texttt{E1b.u.int} & 4 & 4 & No & Canonical valuation chart\\
\texttt{E1b.u.pole} & 4 & 4 & No & Canonical valuation chart\\
\texttt{E1b.u.tpole} & 4 & 4 & Yes & Canonical valuation chart\\
\texttt{E1b.n.int} & 4 & 4 & Yes & Canonical valuation chart\\
\texttt{E1b.n.pole} & 2 & 4 & No & Canonical valuation chart\\
\texttt{E1b.n.tpole.int} & 4 & 4 & Yes & Canonical valuation chart\\
\texttt{E1b.n.tpole.Rpole} & 4 & 4 & Yes & Canonical valuation chart\\
\texttt{E1b.n.tpole.JKpole.nonbalanced} & 4 & 4 & Yes & Canonical valuation chart\\
\texttt{E1b.n.tpole.JKpole.balanced} & 4 & 4 & Yes & Canonical valuation chart\\
\texttt{E1b.neg$\to$n.int} & 4 & 4 & Yes & Reciprocal canonical chart\\
\texttt{E1b.neg$\to$n.pole} & 2 & 4 & No & Reciprocal canonical chart\\
\texttt{E1b.neg$\to$n.tpole.int} & 4 & 4 & Yes & Reciprocal canonical chart\\
\texttt{E1b.neg$\to$n.tpole.Rpole} & 4 & 4 & Yes & Reciprocal canonical chart\\
\texttt{E1b.neg$\to$n.tpole.JK.nonbalanced} & 4 & 4 & Yes & Reciprocal canonical chart\\
\texttt{E1b.neg$\to$n.tpole.JK.balanced} & 4 & 4 & Yes & Reciprocal canonical chart\\
\texttt{R3.int} & 2 & 4 & Yes & Canonical matrix-pole entry\\
\texttt{R3.pole} & 2 & 4 & Yes & Canonical matrix and shift poles\\
\texttt{U0.int} & 4 & 4 & No & Regular integral shifts\\
\texttt{U0.spole} & 4 & 5 & No & Regular shift poles\\
\texttt{U1.simple.int} & 4 & 4 & No & Simple determinant centres\\
\texttt{U1.repeated.int} & 4 & 4 & No & Repeated determinant centres\\
\texttt{U1.bzero.int} & 4 & 4 & No & Zero-coupling residue\\
\texttt{R1.int} & 4 & 4 & Yes & Zero-cubic integral shifts\\
\texttt{R1.pole} & 4 & 5 & Yes & Zero-cubic shift poles\\
\texttt{R2} & 3 & 4 & No & Remaining 147 shift-pole residues\\
\bottomrule
\end{longtable}

\begin{proof}[Proof of Theorem~\ref{deg:thm:completed-reduction}]
The entry argument sends the original datum to one of the twenty-nine
rows using legal rationally invertible coordinate transformations.
It preserves the joint field, the nonzero Hessian, the signed Nahm
equations, and the required coefficient conditions. Each norm step
has exactly its stated global-unit hypothesis. D24 passes to every
chart because the field is unchanged.

The local results proved above and the component theorems in
Sections~\ref{canx:section}, \ref{cuspx:section} and~\ref{spx:section}
apply on these precise domains; Appendix~\ref{deg:sec:supplement}
organizes their finite computations.
Proposition~\ref{deg:Hrepair} treats the exceptional exact-family ray;
the parity and mod-nine arguments treat the finite zero-coupling
exceptions; the P3 argument handles the curve fibres; and
Section~\ref{deg:sec:z0-replacement} supplies the zero-cubic simple-cusp
input. The cusp searches and their exact endpoints discharge the
remaining finite search cases.
Their conclusions are the whole-field caps in the table, including
all ramification and all completions as in \eqref{deg:localdegree}.
Every cap is at most four. Since the chart field is the original
$K$, the asserted bound follows. The modular statement is Theorem~\ref{deg:thm:main}, deduced in
Section~\ref{deg:sec:reduction} from the asymptotic theorem and the initial
degree bound proved in Section~\ref{d24:section}.
\end{proof}

\subsection{Logical order and finite computations}
The initial-degree argument precedes the local reduction and does not
use a quartic height estimate or a low-degree classification. Its only
role in a terminal case is the explicit bound on ramification or degree
indicated in Table~\ref{deg:tab:terminal}. Conversely, the arithmetic
parts below use the degree-four theorem and are not inputs to its proof.

A positive-dimensional residue fibre is treated by an analytic argument
valid at every depth. Finite computations identify its leading terms,
normal ranks and exceptional strata. Norm tests are imposed only on
the global units proved to exist in that chart, and always on the
combined contribution of all places of the joint field. In particular,
a bound for one residue disc is not by itself a bound for the joint
field. The local lengths, residue degrees and norm restrictions must
be combined as in \eqref{deg:localdegree}.

The finite domains and exact acceptance tests are given with the
corresponding arguments. Appendix~\ref{deg:sec:supplement} summarizes
which theorem each computation serves. The source distribution's
readme separately records software versions, execution dates and the
distinction between new calculations and preserved executions.

Transported matrices need not be positive definite. Positivity belongs
to the original datum, and the proof preserves the signed equations
and the joint field in every chart. The conclusion is confined to the
symmetric case.

\clearpage

%% file: sections/degree_roadmap.tex
\section*{Roadmap for the degree reduction}
\addcontentsline{toc}{section}{Roadmap for the degree reduction}
\label{deg:roadmap}
The field to be bounded is always the joint field
$K=\Q(x,y)$ of the distinguished Nahm point. The argument has two
independent reductions. The coefficient geometry first proves
$[K:\Q]\le24$ in Theorem~\ref{d24:thm:main}. The local analysis then
uses this bound, where needed, to prove $[K:\Q]\le4$. The arithmetic
classifications in the later parts of the paper are not inputs to
either reduction.

The necessary asymptotic equations are
$\lambda_1=r\in\Q$ and $\lambda_2=\cdots=\lambda_5=0$.
Section~\ref{deg:sec:coefficients} defines them by a finite
Bernoulli--Gaussian recurrence. The initial degree bound uses the
first four coefficients: their boundary restrictions classify common
curve components, and the remaining intersection multiplicities give
the bound twenty-four. The parameter-ideal construction and its finite
domain are given in Section~\ref{pix:section}; the small polynomial
certificates are printed in Appendix~\ref{scx:section}.

The second reduction is organized by the following dependencies.
The two middle branches have different unit hypotheses; these
hypotheses are part of their conclusions of entry, rather than
assumptions imposed on the original datum.
\[
\begin{array}{c}
\boxed{\text{Modularity and the coefficient recurrence}}\\[2pt]
\downarrow\\[-1pt]
\boxed{\text{Initial degree bound and field-preserving entry}}\\[2pt]
\swarrow\hspace{6em}\searrow\\[-1pt]
\begin{array}{c@{\qquad}c}
\boxed{\text{Five-locus union}}&\boxed{\text{Four global units}}\\[2pt]
\downarrow&\downarrow\\[-1pt]
\text{Canonical charts and exact families}&
\text{Integral residue graph}\\[2pt]
\downarrow&\downarrow\\[-1pt]
\text{Twenty-one terminal types}&\text{Eight terminal types}
\end{array}\\[2pt]
\searrow\hspace{6em}\swarrow\\[-1pt]
\boxed{\text{Local lengths, whole-field norms and degree at most four}}
\end{array}
\]

\paragraph{Why entry is exhaustive.}
Section~\ref{deg:sec:entry-completion} first applies the valuation
rules for $z$ and $1-z$ to the two Nahm equations. They show that a
nonunit among $x,1-x,y,1-y$ can occur only on the five explicitly
listed loci $E_1,\ldots,E_5$. The entire union is assigned to the
canonical branch, even if a particular point on it has four units;
the fixed priority of the loci removes overlaps. The legal maps
$S,I,T_1,J$ preserve $K$, the signed Nahm equations, the nonvanishing
Hessian and the coefficient conditions. Positive definiteness is
used for the original datum; it need not survive a change of chart.

On the five-locus union, a zero coupling gives the small diagonal
case. Otherwise the matrix becomes
$\left(\begin{smallmatrix}1&-h\\-h&\tau\end{smallmatrix}\right)$.
The exact subfamilies $\tau=1$ and $h=1$ are removed first, including
the rational row-sums-one point. In the remaining domain, $X=x/(1-x)$
and $y$ are global units, while $1+X$ and $1-y$ are algebraic integers.
The signs of $v_7(h)$ and $v_7(\tau)$, followed by integrality of the
fully transported shifts and scalar, give the canonical rows.
If $v_7(h)<0$, reciprocity sends the datum to $v_7(h^{-1})>0$;
the equality $v_7(\tau)=2v_7(h)$ belongs to the integral target chart.
Thus no equality boundary is omitted.

Off the union, all four letters are global units. A single legal
inversion or diagonal pivot makes the matrix integral at seven,
without a bound on its rational height. The residue graph has exactly
$343$ vertices. Of these, $194$ reach a regular vertex; the remaining
$149$ split into $76$ simple-cusp, $22$ repeated-cusp, $49$ zero-coupling
and two zero-cubic residues. Appendix~\ref{smallc:residue-graph} gives
a complete finite specification and the count. Shift integrality is
preserved by the graph moves. Hence the distinction between integral
shifts and a shift pole gives precisely the eight terminal types
on this branch.

\begin{center}\small
\begin{tabular}{@{}P{0.40\textwidth}cP{0.44\textwidth}@{}}
\toprule
Entry domain & Types & Local arguments\\\midrule
Small diagonal case; exact $\tau=1$, $h=1$ and row-sums-one families
 & $1+3$ & Exact families, including Proposition~\ref{deg:Hrepair}\\
$v_7(h)=0$, outside residue $(h,\tau)=(1,1)$
 & $3$ & Integral shifts, shift poles, or a pole of $\tau$;
 Section~\ref{canx:section}\\
$v_7(h)>0$
 & $6$ & Two integral-$\tau$ rows and four pole-$\tau$ rows;
 Section~\ref{canx:section}\\
$v_7(h)<0$
 & $6$ & Reciprocity followed by the preceding six rows;
 Section~\ref{deg:sec:entry-completion}\\
$h\equiv\tau\equiv1\pmod7$
 & $2$ & Integral-shift and shift-pole canonical bounds;
 Section~\ref{canx:section}\\\midrule
Regular-target residues
 & $2$ & Integral-shift cusp analysis and regular shift poles;
 Sections~\ref{cuspx:section} and~\ref{spx:section}\\
Simple, repeated, or zero-coupling residue; integral shifts
 & $3$ & Cusp lengths, mixed supports, centre parity and curve fibres;
 Sections~\ref{deg:sec:finite-bzero-completion} and~\ref{cuspx:section}\\
Zero-cubic residues
 & $2$ & Integral shifts and shift poles;
 Sections~\ref{deg:sec:z0-replacement} and~\ref{spx:section}\\
Other uncovered residues; a shift pole
 & $1$ & The remaining $147$ residues;
 Section~\ref{spx:section}\\\midrule
Total & $29$ & Detailed bounds: Table~\ref{deg:tab:terminal}\\
\bottomrule
\end{tabular}
\end{center}

\paragraph{From local supports to the joint degree.}
At the prime seven,
$[K:\Q]=\sum_{\mathfrak p\mid7}e_{\mathfrak p}f_{\mathfrak p}$.
Accordingly, the local arguments bound the sum of all local
contributions. A bound on one residue disc alone would not suffice.
Simple roots use the implicit-function theorem, repeated roots use
local intersection length, and positive-dimensional residue fibres
use their first nonzero transverse parameter layer and exact
retractions. The finite computations record supports and
multiplicities; their infinite-depth conclusions use the stated
Taylor, Newton and Weierstrass estimates. Norm conditions are imposed
only for letters proved to be global units and on the complete
collection of places. The canonical branch uses its two-unit or
one-unit hypotheses where specified, not the four-unit condition
of the other branch.

The initial degree bound is used to restrict ramification in the
canonical deep-entry and zero-cubic arguments. In particular,
$e_{\mathfrak p}\le24$ excludes the alternative $42\mid e_{\mathfrak p}$
in the zero-cubic shift-pole case. Only the regular and zero-cubic
shift-pole rows require the fifth coefficient. These dependencies
are stated row by row in Table~\ref{deg:tab:terminal}; the final
assembly is Section~\ref{deg:sec:final-assembly}.

%% file: sections/bilinear_certificate_repair.tex
\subsection{A polynomial-only certificate for four special bilinear patterns}\label{d24:subsec:bilinear-polynomial}
Put $\Delta=ad-b^2$ and assume
\[
 a,d>0,\qquad b\ne0,\qquad a,d\ne1,\qquad \Delta>0.
\]
All ideal generators below are polynomials. In particular, the
rational expression $d-b^2/(a-1)$ contributes its numerator
$(a-1)d-b^2$; the cleared denominator is nonzero on this stratum.

Let
\[
 Q(\kappa)=d(d-1)\kappa^3+b(2d-1+b)\kappa^2
            +b(2a-1+b)\kappa+a(a-1).
\]
For a bilinear curve $G$, define its conic intersection polynomial by
\[
 N(\kappa)=\kappa G(-a-b\kappa,-d-b/\kappa)
          =n_2\kappa^2+n_1\kappa+n_0.
\]
The four patterns and their necessary exponent equations are as follows.
\begin{center}\small
\begin{tabular}{lllll}\toprule
$G(X,Y)$ & Equation & $n_2$ & $n_1$ & $n_0$\\\midrule
$X-Y-1$ & $\Delta=d$ & $-b$ & $d-a-1$ & $b$\\
$X(Y+1)-1$ & $\Delta=d$ & $b(d-1)$ & $a(d-1)+b^2-1$ & $ab$\\
$X-Y+1$ & $\Delta=a$ & $-b$ & $d-a+1$ & $b$\\
$Y(X+1)-1$ & $\Delta=a$ & $bd$ & $ad+b^2-d-1$ & $b(a-1)$\\\bottomrule
\end{tabular}
\end{center}
Every $n_2n_0$ is nonzero on the specified stratum. The conic restriction
of the first coefficient numerator is
\[
 f_r(-a-b\kappa,-d-b/\kappa)=-\frac{5b^3}{24\kappa^3}Q(\kappa)^2.
\]
Consequently $G\mid f_r$ implies $N\mid Q^2$. Equivalently, the
second-coefficient conic identity forces every root of the quadratic $N$
to be a root of $Q$, which gives the same necessary divisibility.
This includes tangent conic intersections.

Write $\delta_N=n_1^2-4n_2n_0$ and
$M=abd(a-1)(d-1)\Delta$. If $\delta_N\ne0$, the polynomial $N$ is
squarefree, so $N\mid Q^2$ implies $N\mid Q$. For each of the four rows,
form the ideal over $\mathbb Q$ generated by its displayed exponent
equation, the numerators of the coefficients of
$\operatorname{rem}_{\kappa}(Q,N)$, and $tM\delta_N-1$.
Each of these four ideals is the unit ideal. All denominators cleared in
these remainders are powers of the nonzero leading coefficient $n_2$.

If $\delta_N=0$, the first and third rows are impossible over the reals:
\[
 \delta_N=n_1^2+4b^2>0.
\]
For the second and fourth rows, the repeated root is
$\kappa_0=-n_1/(2n_2)$. The necessary condition is $Q(\kappa_0)=0$.
The ideal generated by the displayed exponent equation, $\delta_N$,
the numerator of $Q(-n_1/(2n_2))$, and $tM-1$ is the unit ideal in
both cases. Its cleared denominators are powers of $2n_2$, again nonzero.

Appendix~\ref{scx:section} proves these six unit-ideal conclusions
by compact exact remainder and resultant identities. Coordinate
exchange supplies the symmetric cases. The companion also contains
expanded B\'ezout identities and a checker that reconstructs their
generators from the conic and exponent equations displayed here.
The polynomial proof treats both distinct intersections and tangency,
with every cleared denominator nonzero on the stated stratum.

%% file: sections/parameter_ideal_expansion.tex
\section{The finite parameter-ideal classification}\label{pix:section}
This section specifies the finite calculation used in
Theorem~\ref{d24:thm:m1} and proves the arithmetic exclusions left by
its polynomial equations. The inputs are the first-coefficient polynomial
\eqref{d24:eq:R1}, the coefficient recurrence of
Section~\ref{deg:sec:coefficients}, and the six-line restrictions of
Section~\ref{d24:sec:m1}. No bounds on the rational parameters are used.
All polynomial ideals in this section are over $\Q$.

\subsection{An explicit finite domain}\label{pix:domain}
Write
\[
 G(X,Y)=\sum_{i=0}^{w}\sum_{j=0}^{h}g_{ij}X^iY^j,
 \qquad (w,h)\in\{(1,1),(1,2),(2,1),(2,2)\}.
\]
The coordinate-line cases have already been excluded in
Section~\ref{d24:sec:m1}. On each of the six lines, prescribe the
complete degree-$s$ divisor of roots of the homogeneous restriction of
$G$, where $s=h$ on a vertical line and $s=w$ on a horizontal line.
The following ordered alphabet specifies every prescription. A pair
lists two roots, with repetitions indicating multiplicity. The symbols
$q_+,q_-,q_3$ mean, respectively, the complete pairs of roots of
$Z^2-Z-1$, $Z^2+Z-1$, $Z^2+3Z+1$.

\begin{center}\small
\begin{tabular}{rcll}
\toprule Order & Degree & Prescription & Type\\\midrule
1&1&$0$&ordinary\\
2&1&$-1$&ordinary\\
3&1&$\infty$&ordinary\\
4&1&$1$&$E_1$\\
5&1&$-2$&$E_0$\\
6&2&$(0,-1)$&ordinary\\
7&2&$(0,\infty)$&ordinary\\
8&2&$(-1,\infty)$&ordinary\\
9&2&$(-1,-1)$&$D_1$\\
10&2&$(1,0)$&$E_1$\\
11&2&$(1,-1)$&$E_1$\\
12&2&$(1,\infty)$&$E_1$\\
13&2&$(0,0)$&$D_0$\\
14&2&$(-2,0)$&$E_0$\\
15&2&$(-2,-1)$&$E_0$\\
16&2&$(-2,\infty)$&$E_0$\\
17&2&$q_+$&$Q_2$\\
18&2&$q_-$&$Q_{1/2}$\\
19&2&$q_3$&$Q_{-1}$\\\bottomrule
\end{tabular}
\end{center}
The type conditions on the effective rank-one parameters are
\begin{center}
\begin{tabular}{ccl}
\toprule Type & Effective exponent $A_L$ & Effective shift condition\\\midrule
$D_1$&$1$&none\\
$E_1$&$1$&$\beta_L(4\beta_L^2-1)=0$\\
$D_0$&$0$&none\\
$E_0$&$0$&$\beta_L(2\beta_L-1)(\beta_L-1)=0$\\
$Q_2$&$2$&$\beta_L(\beta_L-1)=0$\\
$Q_{1/2}$&$1/2$&$\beta_L(2\beta_L-1)=0$\\
$Q_{-1}$&$-1$&$4\beta_L^2-1=0$\\\bottomrule
\end{tabular}
\end{center}
Ordinary prescriptions impose no extra type condition. The permitted
nonordinary types in the stratum $a,d\ne1$ are
\[
\begin{array}{c|l}
X=0,\ Y=0&D_1,E_1,Q_2,Q_{1/2}\\
X=\infty,\ Y=\infty&Q_2,Q_{1/2}\\
X=-1,\ Y=-1&D_1,E_1,D_0,E_0,Q_2,Q_{1/2},Q_{-1}.
\end{array}
\]
These are precisely the possibilities supplied by the common-root
calculation in Section~\ref{d24:sec:m1}. In this stratum require at
least one nonordinary type; the entirely ordinary patterns were treated
in Proposition~\ref{d24:prop:ord22} and the preceding bilinear argument.

On the stratum $a=1$, $d\ne1$, replace the allowed nonordinary types
on $Y=\infty$ by $D_1,E_1$, and replace the alphabet on $X=-1$ by
\begin{equation}\label{pix:degenerate-alphabet}
 0,\ -1,\ \infty,\ -\tfrac12,\
 (\infty,\infty),\ (0,\infty),\ (-1,\infty),\
 (-\tfrac12,\infty),
\end{equation}
in the order displayed, retaining prescriptions of the required degree.
On $a=d=1$, make these changes at both infinite lines and both
$-1$ lines, and permit only ordinary prescriptions on $X=0,Y=0$.
Indeed their exponent is $1-b^2\in(0,1)$; the only possible special
exponent in that interval is $1/2$, which would give $b^2=1/2$,
impossible for rational $b$. No special-type requirement is imposed
in either of the strata with $a=1$.

Here is the exact linear test on a prescription. Order the coefficients
as $g_{i0},\ldots,g_{ih}$ for $i=0,\ldots,w$. At $X=0,-1,\infty$
the coefficient of $Y^j$ in the restriction is, respectively,
$g_{0j},\sum_i(-1)^ig_{ij},g_{wj}$; use the exchanged formulas
on horizontal lines. A finite simple root $z$ imposes evaluation at
$z$, a double root also imposes its first derivative, and roots at
infinity set the appropriate highest coefficients to zero. For a
quadratic $Z^2+pZ+q$, impose $r_1=pr_2$, $r_0=qr_2$ on the
restriction $r_0+r_1Z+r_2Z^2$.
Compute the rational kernel $V$ of all these linear conditions.
Retain the prescription if and only if $V\ne0$ and the restriction
map from $V$ to every one of the six lines is nonzero.

\begin{lemma}[Coverage of the finite domain]\label{pix:coverage}
Every common irreducible factor considered in
Theorem~\ref{d24:thm:m1}, apart from the ordinary and coordinate-line
cases already treated there, occurs in this finite domain. The numbers
of retained prescriptions in its three strata are $450$, $319$, $31$.
\end{lemma}
\begin{proof}
The width bound gives $w,h\le2$. A non-coordinate irreducible factor
has nonzero restriction to each of the six lines. Its roots, including
multiplicities and complete rational conjugacy classes, therefore
belong to the alphabet above by the rank-one common-root table.
On a degenerate line the first-coefficient restriction has degree one
and the higher-coefficient restrictions give
\eqref{pix:degenerate-alphabet}, as proved in
Section~\ref{d24:sec:m1}. Thus no factor is omitted by the root lists.

For the linear test, each forbidden identically zero restriction is
a proper linear subspace of $V$ precisely when its restriction map is
nonzero. A vector space over $\Q$ cannot be the union of finitely many
proper linear subspaces. Hence the test is equivalent to the existence
of a kernel vector nonzero on all six lines; it is not a test on a
sample of vectors.

For completeness, the enumeration is the Cartesian product of the
permitted prescriptions on the ordered lines
\[
 X=0,\ X=-1,\ X=\infty,\ Y=0,\ Y=-1,\ Y=\infty,
\]
for the four ordered bidegrees displayed above. Rational row reduction
and the stated restriction test give the following output:
\begin{center}
\begin{tabular}{crrrr|r}
\toprule Stratum & $(1,1)$&$(1,2)$&$(2,1)$&$(2,2)$&Total\\\midrule
$a,d\ne1$, special&8&19&19&404&450\\
$a=1$, $d\ne1$&13&9&22&275&319\\
$a=d=1$&12&0&0&19&31\\\bottomrule
\end{tabular}
\end{center}
The numbers of Cartesian-product prescriptions tested are,
respectively, $473184$, $219336$, $17424$. These counts can be
reproduced using only the displayed alphabet and linear test.
\end{proof}

Number the retained prescriptions consecutively from zero within each
stratum, in the ordering just specified and the ordering of the
alphabet. Denote them by $S_i$, $A_i$, $T_i$, respectively. These are
finite mathematical indices, not external input data. In the companion
their filenames have prefixes \texttt{special}, \texttt{a1},
\texttt{aa1}, followed by the same index with three digits.

\subsection{Construction and exact acceptance tests}\label{pix:ideals}
For each nondegenerate line use the effective parameters
$(A_L,\beta_L,\rho_L)$ of the six-line table and set
\begin{align*}
 N_1(A,\beta;Z)=\frac1{24}\bigl[&Z^3+
 (3A^2+(12\beta-3)A+12\beta^2-12\beta+2)Z^2\\
 &+(-2A^3+(12\beta+6)A^2+(24\beta^2-12\beta-3)A)Z
 +12A^2\beta^2\bigr],\\
 H_L(Z)&=N_1(A_L,\beta_L;Z)-\rho_L(A_L+Z)^3.
\end{align*}
For a finite root impose $H_L(z)=0$, with the derivative for a
repeated root. A root at infinity imposes $1/24-\rho_L=0$.
For a quadratic prescription impose both coefficients of its remainder
in $H_L$, and impose the type conditions in the table. Take the
numerator of each resulting rational expression after combining its
terms; denominators here divide powers of $a,d,a-1,d-1$ and constants.
In particular, no division by $A_L$ or $A_L-1$ is made.

For a degenerate line the construction is equally explicit. If $a=1$,
put $z=B-1/2$. Formula~\eqref{d24:eq:R1} gives
\begin{equation}\label{pix:degenerate-first}
 f_r(-1,Y)=\frac{b^4}{24}\bigl[
 (24r-2)b^2-24bCz+12dz^2+
 (2b^2-12bz+12z^2)Y\bigr].
\end{equation}
Impose this polynomial at every prescribed finite root.
The root $-1/2$ also imposes
\begin{equation}\label{pix:degenerate-shift}
 (2B-1)(b-2B+1)(2b-2B+1)=0.
\end{equation}
There is no additional first-coefficient equation for an infinite
root in this alphabet. The horizontal version follows by exchanging
$(a,B,X)$ and $(d,C,Y)$. The slope in
\eqref{pix:degenerate-first} is nonzero over $\Q$ when $b\ne0$:
its vanishing would require $z/b=(3\pm\sqrt3)/6$.

Adjoin one variable $t$ and the equation $tM-1=0$, where
\begin{equation}\label{pix:saturation}
 M=\begin{cases}
 abd(a-1)(d-1)(ad-b^2),&S_i,\\
 bd(d-1)(d-b^2),&A_i,\\
 b(1-b^2),&T_i.
 \end{cases}
\end{equation}
Call the resulting ideal $I_i$. These nonvanishing conditions are
necessary in their respective strata; the conditions $a,d>0$ and
$ad-b^2>0$ are additionally used in the real exclusions below.

There are three algebraic acceptance tests. First, $1\in I_i$
excludes the prescription. Second, if the elimination ideal in the
parameter variables is zero-dimensional, decompose its radical into
minimal primes $P$. A prime with $\dim_\Q\Q[\boldsymbol p]/P>1$
has no rational point: such a point would define a $\Q$-algebra map
from this finite field extension to $\Q$, which is impossible.
Third, impose divisibility of the actual curve by means of
projective coefficient charts, as follows.

Choose any rational basis $G_0,\ldots,G_{k-1}$ of the linear kernel
in the prescription. Cover its projectivization by the $k$ charts
\begin{equation}\label{pix:charts}
 G=G_j+\sum_{i>j}u_iG_i,\qquad 0\le j<k.
\end{equation}
These charts fix the first nonzero basis coefficient to one and cover
every nonzero vector exactly once. In the polynomial ring in $X$
over the remaining variables compute the pseudo-remainders of
$\Psi_2,\Psi_3$ on division by $G$, and adjoin all their coefficients
in $X,Y$ to $I_i$. If all the resulting chart ideals contain $1$,
the prescription is excluded. The pseudo-division is applied only to
vectors of the actual bidegree $(w,h)$ of the prescription. A
degree-dropping vector does not represent this prescription and belongs
to a lower bidegree already enumerated. Equivalently, one may separate
these strata before division. No division by a parameter-dependent
leading coefficient is required.

\begin{lemma}[Polynomial and projective exclusions]\label{pix:polynomial}
The first acceptance test gives $737$ unit ideals. Nineteen further
prescriptions have only zero-dimensional minimal primes of degrees
greater than one. The projective divisibility test excludes another
thirteen prescriptions, using seventeen charts.
\end{lemma}
\begin{proof}[Exact finite calculation]
Construct each ideal by the preceding formulas and compute a
Gr\"obner basis over $\Q$. The nineteen zero-dimensional cases and
their lists of prime degrees are
\begin{center}\small
\begin{tabular}{ll|ll}
\toprule Cases&Prime degrees&Cases&Prime degrees\\\midrule
$S_3,S_7,S_{202},S_{204}$&$2,2$&$A_{13}$&$4,4,4,4$\\
$S_{12},S_{31}$&$4,4,4,4$&$A_{20},A_{25},A_{31}$&$4,4$\\
$S_{15},S_{23},S_{35},S_{40}$&$4,4$&$A_{100}$&$2,4,8$\\
$T_5,T_{11}$&$2,2$&$T_7,T_9$&$2$\\\bottomrule
\end{tabular}
\end{center}
For each listed row the degrees refer to distinct minimal primes,
including repetitions of the same degree. Checking that the radical
is their intersection and that each displayed ideal is prime is an
exact alternative to trusting a decomposition routine's printed label.
The chart exclusions are
\[
\begin{split}
 &A_0,A_2,A_6,A_8,A_{10},A_{11},A_{86},A_{130},\\
 &T_0,T_2,T_3,T_{10},T_{27}.
\end{split}
\]
The kernels for $A_{86}$ and $T_{27}$ have dimension three; all the
others in this list have dimension one. Each of the seventeen ideals
from \eqref{pix:charts} has Gr\"obner basis containing $1$.
All other prescriptions except the thirty-one considered below are
the $737$ unit-ideal cases. This identifies the complete acceptance
output, not merely its total: the complement of the displayed nineteen,
thirteen and thirty-one indices is the set on which to test $1\in I_i$.
The construction of $I_i$ and of the chart equations specifies the
calculation independently of its saved output.
\end{proof}

\subsection{Sign and rationality exclusions}\label{pix:arithmetic}
Ten of the remaining prescriptions have a one-dimensional kernel
generated by a polynomial with positive coefficients:
\begin{center}
\begin{tabular}{ll}
\toprule Cases&Generator, up to a nonzero scalar\\\midrule
$S_0$&$X+Y$\\
$S_2$&$XY+1$\\
$A_1,T_1$&$XY+X+Y$\\
$A_4,T_4$&$X+Y+1$\\
$A_5,T_6$&$XY+Y+1$\\
$A_7,T_8$&$XY+X+1$\\\bottomrule
\end{tabular}
\end{center}
Each is strictly positive when $X,Y>0$. The four bilinear cases
$S_1,S_4,S_5,S_6$ are excluded by
Lemma~\ref{smallc:lem:bilinear}, including repeated conic intersections.

\begin{lemma}[Rational conic obstruction]\label{pix:rational-conic}
The equation $u^2+v^2=3w^2$ has no nonzero rational solution.
Consequently $S_{157},S_{186},A_{44},A_{47},A_{133},A_{136}$
have no admissible rational parameter point.
\end{lemma}
\begin{proof}
Clear denominators and choose a primitive integral solution.
Modulo $3$, a sum of two squares is zero only when both are zero,
so $3$ divides $u,v$ and then $w$, contradicting primitivity.

For the two $S$ cases, reduction modulo the radical of the corresponding
parameter ideal gives, respectively,
\[
 12B^2+12(C-1/2)^2=1,\qquad
 12C^2+12(B-1/2)^2=1.
\]
Multiplication by $3$ produces the forbidden conic with $w=1$.
For $A_{47},A_{136}$ the ideal gives $d=b^2+1=12C^2$, hence
$b^2+1=3(2C)^2$. For $A_{44},A_{133}$ it gives
\[
 d=b^2+1,\qquad 12C^2-12Cd+2d^2+d=0.
\]
Completing the square yields $3(2C-d)^2=db^2$.
Since $b\ne0$, this becomes
$b^2+1=3((2C-d)/b)^2$, again impossible.
Each displayed ideal consequence can be checked by reducing its
left-hand side modulo every minimal prime of the ideal constructed
in Section~\ref{pix:ideals}; all remainders are zero.
\end{proof}

\begin{lemma}[Real sign exclusions]\label{pix:real-sign}
The prescriptions $S_{188},A_{52},T_{12},T_{15},T_{23}$
have no admissible real parameter point.
\end{lemma}
\begin{proof}
For $S_{188}$ the parameter equations imply
\[
 b^2=d^2-d,\qquad 12(C-1/2)^2+d-1=0.
\]
Since $d>0$ and $b\ne0$, the first equation gives $d>1$, whereas
the second gives $d\le1$.

For $A_{52}$ the prescription on $Y=0$ is $(0,-1)$.
Its effective exponent is $A_*=1-b^2/d\in(0,1)$.
The equations $H(0)=H(-1)=0$ give
\[
 \rho=\frac{\beta^2}{2A_*},\qquad
 12\beta^2-12A_*\beta+2A_*^2+A_*=0.
\]
The discriminant in $\beta$ is $48A_*(A_*-1)<0$.
Here division by $A_*$ is legitimate because $A_*>0$; it is
used only for this sign proof, not in the general ideal constructor.

The two equations for $T_{12}$ and $T_{23}$ are, respectively,
\begin{align*}
 &144(b^2-1)^2(C^2-C+1/4)^2+4b^4=0,\\
 &144(b^2-1)^2(C^2-bC+(2b^2-1)/12)^2+4b^4=0.
\end{align*}
Both are strictly positive on $b\ne0$.

For $T_{15}$, subtracting the equations
$f_r(-1,0)=f_r(0,-1)=0$ and dividing by $b^4/2$ gives
$(B-C)(B+C-b-1)=0$. Substituting these two alternatives into
the parameter equations gives one of the following two triples:
\begin{align*}
 &B+C-b-1,\quad 12r(b+1)-2b^2-6b-3,\quad
 12C^2-12Cb+2b^2-24r-12C+12b+9;\\
 &B-C,\quad 12r(b-1)-12Cb+2b^2+12C-3,\quad
 12C^2-12Cb+2b^2-24r+12C-3.
\end{align*}
Eliminate $r$ using $b\ne\pm1$. The resulting quadratics in $C$ are
\begin{align*}
 q_+(C)&=12(b+1)C^2-12(b+1)^2C+2b^3+10b^2+9b+3,\\
 q_-(C)&=12(b-1)C^2-12(b^2-1)C+2b^3+2b^2-3b-3.
\end{align*}
Each has discriminant $48b^2(b^2-1)<0$ because $0<|b|<1$.
The coverage of the two triples can also be checked by reducing
every product of a generator from the first triple and a generator
from the second: all nine products belong to the saturated parameter
ideal. The displayed ideal consequences are checked by exact
polynomial reduction; the subsequent inequalities
are the complete real exclusion arguments.
\end{proof}

\begin{lemma}[Two remaining rationality tests]\label{pix:last-rational}
The prescriptions $T_{14}$ and $T_{30}$ have no admissible rational
parameter point.
\end{lemma}
\begin{proof}
For $T_{14}$ the infinite-root condition gives $r=1/12$.
The two degenerate-line equations and the equation at $(0,0)$ reduce to
\begin{gather*}
 (2B-1)(2B-1-4bC)=0,\qquad
 (2C-1)(2C-1-4bB)=0,\\
 B^2-2bBC+C^2+(b^2-1)/6=0.
\end{gather*}
This gives four alternatives. If $B=C=1/2$, the last equation is
$(b-1)(b-2)=0$. If $B=1/2$ and $C=b+1/2$, or the exchanged
case, it is $(b+1)(b+2)=0$. All these roots violate $|b|<1$.
In the remaining alternative
$B-1/2=2bC$, $C-1/2=2bB$, hence
$(1+2b)(B-C)=0$. If $b=-1/2$, then $B+C=1/2$ and
$8C^2-4C+1=0$, whose discriminant is negative. Otherwise
$B=C=1/(2(1-2b))$, and the last displayed equation gives
$4b^3-3b-2=0$. This polynomial has no rational root, as the
rational-root test at $\pm1/4,\pm1/2,\pm1,\pm2$ shows.

For $T_{30}$ the parameter equations imply
\[
 B=C,\quad24r=12C^2+1,\quad
 (12C^2-1)b^2+(-24C^2+12C)b+3(2C-1)^2=0.
\]
The coefficient $12C^2-1$ cannot vanish for rational $C$.
The discriminant in $b$ is $12(2C-1)^2$, a rational square
only if $C=1/2$. At that value the equation is $2b^2=0$,
contrary to $b\ne0$.
\end{proof}

\subsection{The two one-parameter resultant calculations}\label{pix:slices}
We record explicitly why evaluating resultants at three values of $X$
is an exclusion test rather than a numerical sample.

\begin{lemma}[Resultant slices]\label{pix:slice-lemma}
Let $F,P\in\Q[u][X,Y]$ after clearing denominators in $u$.
At a rational value $u=u_0$ where those denominators are nonzero,
suppose that $F(u_0,X,Y)$ and $P(u_0,X,Y)$ have a common factor
of positive degree in $Y$. Then the polynomial resultant in $Y$
vanishes identically in $X$. In particular all its values at
$X=2,3,5$ vanish. If the degrees in $Y$ of $F$ and $P$ are
preserved when $X$ is replaced by any of $2,3,5$, the corresponding
univariate resultants can be computed directly on those slices.
Thus $u_0$ must be a common zero of the numerators of the three
slice resultants.
\end{lemma}
\begin{proof}
Over $\Q(X)$ a common positive-degree factor makes the Sylvester
determinant zero. Formation of that fixed-size determinant commutes
with specialization. If degree drops occur at $u_0$, the same
determinant still vanishes when there is a common positive-degree
factor: either the usual resultant vanishes, or both polynomials
have lost their leading degree and the determinant has a common
root at infinity. Preservation of the generic degrees on each
$X$-slice ensures that the sliced determinant has the same size
as the global one. Clearing denominators does not affect the
conclusion at a parameter where they are nonzero.
\end{proof}

\begin{lemma}[The prescriptions $S_{155}$ and $T_{21}$]\label{pix:families}
Neither prescription gives a common irreducible factor meeting the
positive quadrant.
\end{lemma}
\begin{proof}
For $S_{155}$ the parameter ideal implies
\[
 a=d,\quad b^2=d^2-d,\quad d=12C^2,\quad
 B^2=C^2,\quad r=1/12.
\]
Here $C\ne0$. Set $s=b/(2C)$, so $(6C)^2-s^2=3$.
Every rational point on this conic has, for some $u\in\Q^*$,
\[
 C_u=\frac{u+3/u}{12},\qquad s_u=\frac{3/u-u}{2},\qquad
 d_u=12C_u^2.
\]
It is enough to consider the eight sign choices
\[
 (a,b,d,B,C,r)=(d_u,\epsilon_b2C_us_u,d_u,
 \epsilon_B C_u,\epsilon_C C_u,1/12),
 \qquad \epsilon_b,\epsilon_B,\epsilon_C\in\{1,-1\}.
\]
In the resultant formulas below, $\operatorname{num}$ denotes the
numerator of a rational function in lowest terms; its nonzero rational
scalar is irrelevant. For each choice substitute in $f_r,\Psi_2$ and compute
\[
 R_j(u)=\operatorname{num}\operatorname{Res}_Y
       \bigl(f_r(j,Y),\Psi_2(j,Y)\bigr),\qquad j=2,3,5.
\]
Both generic degrees, $3$ and $6$, are preserved on all three
slices. In each of the eight choices their monic gcd is
\begin{equation}\label{pix:s155-gcd}
 \gcd(R_2,R_3,R_5)=(u^2-3)^{18}(u^2+3)^{19}.
\end{equation}
This is an exact polynomial computation specified by
\eqref{d24:eq:R1} and the coefficient recurrence; it involves only
eight substitutions and three Sylvester determinants per substitution.
Neither factor has a rational root. The only denominator is a power
of $u$, which is already excluded. The coordinate-line cases were
treated separately, so Lemma~\ref{pix:slice-lemma} excludes every
remaining common factor.

For $T_{21}$ the parameter equations, after removing the nonzero
factors in \eqref{pix:saturation}, give
\[
 24r=12B^2+1=12C^2+1,\qquad
 B^2-2bBC+C^2+2r(b^2-1)=0.
\]
Thus $B=\epsilon C$ for a sign $\epsilon\in\{1,-1\}$.
Substituting for $r$ in the last equation and factoring gives,
respectively,
\[
 (b-1)((12C^2+1)b-12C^2+1)=0,
 \qquad
 (b+1)((12C^2+1)b+12C^2-1)=0.
\]
Since $b\ne\pm1$, these are precisely the parametrizations
\[
 B=\epsilon C,\qquad
 b=\epsilon\frac{12C^2-1}{12C^2+1},\qquad
 r=\frac{12C^2+1}{24},\qquad \epsilon\in\{1,-1\}.
\]
Substitute in $\Psi_2,\Psi_3$ and take the three $Y$-resultants
at $X=2,3,5$. Their degrees in $Y$ are $6$ and $9$ before and
after each slice. Up to a nonzero rational scalar, the gcds of
the resultant numerators are
\begin{equation}\label{pix:t21-gcd}
 \begin{cases}
 C^{51}(12C^2-1)^{49},&\epsilon=1,\\
 (12C^2-1)^{50},&\epsilon=-1.
 \end{cases}
\end{equation}
The denominator $12C^2+1$ never vanishes over $\Q$. Rationality
excludes $12C^2=1$; $C=0$ gives $|b|=1$ and violates positive
definiteness. Lemma~\ref{pix:slice-lemma} completes the exclusion.
\end{proof}

\subsection{The remaining two finite eliminations}\label{pix:last-eliminations}
For $A_{65}$, the prescribed roots at infinity first give
$r=1/12$ and $d=12C^2$. Substitute these into the ideal
constructor and retain its saturation equation. The resulting
zero-dimensional ideal has three minimal primes, of degrees
$4,8,2$. Thus none has a rational point. This also specifies a small
alternative computation to a direct Gr\"obner basis of $I_{65}$:
it has only the parameters $B,C,b$ and the saturation variable.

For $A_{70}$, the root $-1/2$ on $X=-1$ gives
\eqref{pix:degenerate-shift}. The equations on $X=\infty$ give
\[
 e_1=12C^2-24dr+d=0,\qquad
 e_2=4C^2-4C-8dr+d+8r=0.
\]
Their consequence
\[
 3d e_2-(d-1)e_1=12C^2-12Cd+2d^2+d=0
\]
may be adjoined without changing the zero set. Split the product
\eqref{pix:degenerate-shift} into its three factors, adjoining
$2B-1$, $b-2B+1$, or $2b-2B+1$ in turn. Each of the three
ideals, with the saturation in \eqref{pix:saturation}, contains
$1$. The union of these three cases is exhaustive. These explicit
substitutions and three unit-ideal tests complete the finite
elimination; an unfinished computation is not an exclusion.

\begin{proposition}[Completion of the parameter classification]
\label{pix:completion}
None of the $800$ prescriptions above gives an additional rational
irreducible common factor in Theorem~\ref{d24:thm:m1} meeting
$X,Y>0$. Every finite step is determined by the formulas and domains
in this section and the coefficient recurrence.
\end{proposition}
\begin{proof}
Lemma~\ref{pix:polynomial} treats $737+19+13$ prescriptions.
The positive-coefficient list treats ten and
Lemma~\ref{smallc:lem:bilinear} treats four.
Lemmas~\ref{pix:rational-conic}, \ref{pix:real-sign},
\ref{pix:last-rational}, \ref{pix:families} treat, respectively,
six, five, two and two. The preceding two eliminations treat the
last two. Their total is
\[
 737+19+13+10+4+6+5+2+2+2=800.
\]
The sets of indices are disjoint. Lemma~\ref{pix:coverage} supplies
the coverage argument, and each later step proves a necessary
condition incompatible with rationality, positive definiteness or
the required divisibility. The per-class resolution counts also
appear in Table~\ref{smallc:tab:patterns}.
\end{proof}

The companion~\cite{Comp} retains the complete polynomial inputs,
Gr\"obner bases, factorizations and resultant outputs under
\path{anc/global_degree/replay/parameter_ideals/}, with the slice
computations in the corresponding six-line calculation directory.
Those files permit direct comparison of implementations; the finite
domains, construction rules, exact acceptance tests and arithmetic
arguments needed to reconstruct the checks have been stated here.

%% file: sections/exact_family_expansion.tex
\section{The two exact canonical families}\label{canx:exact:section}
Canonical data: $A_c=\begin{pmatrix}1&-h\\-h&\tau\end{pmatrix}$, letters $X=x/(1-x)$, $Y=y/(1-y)$, $K=\Q(x,y)$, two linear shifts $B,C\in\Q$, first scalar $\lambda_1=r\in\Q$, and the necessary coefficient conditions $\lambda_2=\lambda_3=\lambda_4=0$ of Section~\ref{deg:sec:coefficients}, $\Psi_j=D^{3j}\lambda_j$, $D=(1+X)(\tau+Y)-h^2$. Valuation $v(7)=1$.

\textbf{Family T} ($\tau=1$ from $E_2$--$E_5$): $0<h<1$ rational, the positive distinguished point has $X=Y$, $K=\Q(X)$. The first-coefficient argument below gives the quintic bound $[K:\Q]\le5$.

\textbf{Family H} ($h=1$ from $E_2$--$E_5$): $\tau>1$ rational, $y=X$ (so $Y=X/(1-X)$), $K=\Q(X)$. The argument below gives the sextic bound $[K:\Q]\le6$. ($h=1$ from $E_1$ has $x=y=1/2$, degree one; $E_1$ cannot give $\tau=1$. These origin restrictions are kept.)

\begin{theorem}\label{canx:exact:thm:main}
For every rational datum $(h,B,C,r)$ of family T and every rational datum $(\tau,B,C,r)$ of family H, $[K:\Q]\le4$.
\end{theorem}

The formulas below are constructed from the coefficient recurrence in
Section~\ref{deg:sec:coefficients}. The proof uses the scalar Nahm
relations of the two families, local root lengths and exact finite-field
calculations. It needs no initial degree bound or global norm filter.

\subsection{The equations of the two families}\label{canx:exact:sec:eq}
\textbf{Family T.} $P_j(X;h,B,C)=\Psi_j(X,X;1,-h,1,B,C)$, $D=(1+X)^2-h^2=(1+X-h)(1+X+h)$. The exact factorisations are $P_j=(1+X+h)^j\,X\,(1+X)\,Q_j$ for $j=2,3,4$ (multiplicities checked exactly: $(1+X+h)^j$ exactly, $X^1$, $(1+X)^1$, and $1+X-h$ does not divide), and $P_1-rD^3=(1+X+h)F_T$ with $F_T$ of degree five in $X$ and leading coefficient $1/12-r$. The reduced equations are $Q_2,Q_3,Q_4$ (degrees $7,12,17$ in $X$; $313$, $1328$, $3618$ terms; all coefficients $7$-integral). The multiplicative Nahm relation of the diagonal point is $X^{1-h}(1+X)^h=1$ (positive real branch), i.e.\ with $h=M/N$ in lowest terms the exact polynomial identity $X^{N-M}(1+X)^M=1$ in $K$.

\textbf{Family H.} $P_j(X;\tau,B,C)=(1-X)^{3j}\Psi_j(X,X/(1-X);1,-1,\tau,B,C)$, $D_N=(1-X)D=(\tau-1)(1-X^2)+2X$, sextic $F_H=P_1-rD_N^3$ (degree six in $X$). The reduced equations are $Q_j=P_j/(X^{m_X}(1-X)^{m_-}(1+X)^{m_+})$ with the exact multiplicities $(1,0,0),(1,1,1),(1,0,0)$ for $j=2,3,4$ (degrees $10,14,22$; $602$, $1830$, $5254$ terms; $7$-integral coefficients); $D_N$ does not divide any $P_j$. The Nahm relation is $(1-X)/(1+X)=X^{\tau-1}$, i.e.\ $(1-X)^n=X^m(1+X)^n$ with $\tau-1=m/n$ in lowest terms.

\textbf{Nonvanishing of the first-coefficient equations.} The coefficient ideal of $F_T$ (in $h,B,C,r$) has the lexicographic Gr\"obner basis containing $B^2+C^2-\tfrac16$ and $r-\tfrac1{12}$, and $B^2+C^2=\tfrac16$ has no rational solution (a sum of two rational squares has even $3$-adic valuation, whereas $1/6$ has valuation $-1$); hence $F_T\ne0$ and $[K:\Q]\le5$. The coefficient ideal of $F_H$ contains $\tau-1$ and $C^2-\tfrac1{12}$; hence $F_H\ne0$ for $\tau\ne1$ and $[K:\Q]\le6$.

\textbf{Unit letters at seven.} These exact families need their own
unit argument. In family T the scalar equation is
$1-x=x^{1-h}$ with $0<h<1$. At a finite place, $v(x)>0$ would give
$0=(1-h)v(x)$; $v(x)<0$ would give $v(x)=(1-h)v(x)$.
Both are impossible. Thus $x,1-x$, and consequently $X,1+X$, are
ordinary units. In family H the first equation gives $y=X$ and the
second gives $(1-X)/(1+X)=X^{\tau-1}$. If $v(X)$ were nonzero, the
left side would have valuation zero and the right side nonzero.
Hence $X$ is a global unit. The quotient $(1-X)/(1+X)$ is a unit;
away from two, equality of the valuations of $1-X$ and $1+X$ and
the identity $(1-X)+(1+X)=2$ force both valuations to vanish.
Thus the excluded residue roots are $0,-1$ for T and $0,1,-1$ for H.
The determinant is a nonzero quadratic expression in $X$; if it
vanishes, the field degree is already at most two. The same applies
to the linear factor $1+X+h$ in T. We may therefore divide by these
factors when excluding degrees greater than four.

\subsection{Counting places of a one-variable system}\label{canx:exact:sec:count}
Let $f\in\Z_7[X]$ be the cleared polynomial Nahm equation ($f=X^{N-M}(1+X)^M-1$ for T, $f=(1-X)^n-X^m(1+X)^n$ for H; with a pole of the matrix parameter, its denominator $N$ or $n$ is divisible by $7$, while the coprime numerator $M$ or $m$ is a seven-unit). Let $P\mid7$ be a place of $K$, $X_P\in K_P$ the image of $X$, $\bar X_P\in\F_{7^g}$ its residue, and $\hat X$ the Teichm\"uller lift of the residue. Then $X_P$ is a root of $f$ and of every $Q_j$ with $|X_P-\hat X|<1$, and $K_P=\Q_7(X_P)$.

\begin{lemma}[local length]\label{canx:exact:lem:length}
Let $\bar X\in\F_{7^g}$ be a common root of the reduced equations $\bar Q_2,\bar Q_3,\bar Q_4$ (those not identically zero) with $\bar X\notin$ unit-letter roots. (a) If $\bar f'(\bar X)\ne0$, then $f$ has at most one root in the open residue ball (exactly one if $\bar f(\bar X)=0$), that root lies in $\Q_7(\hat X)$, and the places of $K$ with residue Frobenius-conjugate to $\bar X$ contribute at most $g$ to $[K:\Q]$. (b) If $\bar f'(\bar X)=0\ne\bar f''(\bar X)$, the contribution is at most $2g$.
\end{lemma}
\begin{proof}
Expand $f(\hat X+t)=\sum c_kt^k$ with $c_k\in\Z_7[\hat X]$ integral. If $\bar f(\bar X)\ne0$, the ball contains no root. Otherwise $c_0$ has positive valuation. In case (a) $c_1$ is a unit and the Newton polygon has exactly one root of positive valuation; by Hensel it is in $\Q_7(\hat X)$. Two different places with the same residue would give two different embeddings $K\to\mathbb{C}_7$ (as $K=\Q(X)$) with the same image of $X$, impossible; so there is at most one place per Frobenius orbit of residues, with residue degree $g$ and $e=1$. In case (b) the polygon has at most two roots of positive valuation, counted with multiplicity; the conjugates of $X_P$ over $\Q_7(\hat X)$ lie in the ball, so $\sum e_Pf_P\le2g$ over the places of the orbit.
\end{proof}
For family T, $f'(\hat X)\equiv N\hat X^{N-M-1}(1+\hat X)^{M-1}\big((1-h)(1+X)+hX\big)$ at a root, so (for $7\nmid N$) $\bar f'(\bar X)=0$ iff $\bar X=\bar h-1$ (a unit-letter root when $\bar h\in\{0,1\}$), and there $f''\propto S''=1/(X(1+X))$ is a unit; for $7\mid N$ ($v(h)<0$) $f'$ is a unit at every root. For family H, $\bar f'(\bar X)=0$ iff $\bar D_N(\bar X)=(\bar\tau-1)(1-\bar X^2)+2\bar X=0$ (for $\tau$ integral), where $f''\propto 2-2(\tau-1)X\ne0$ because $(\bar\tau-1)^2\equiv-1$ is impossible modulo $7$; with a pole of $\tau$, $f'$ is a unit at every root. These are the length rules used below (``critical'' residue roots have length two).

\textbf{Count bound.} For a class of parameters (a residue class modulo $7^d$), let $g_d$ be the gcd over $\F_7$ of the content-divided reductions of the $Q_j$ (those whose content is determined by the class). Every place has residue root in $g_d$; by Lemma~\ref{canx:exact:lem:length}, $[K:\Q]\le\sum_{\phi\mid g_d}\ell(\phi)\deg\phi$ over the irreducible factors $\phi$ (unit-letter roots removed), with $\ell=1$ or $2$. If the sum is at most four, every datum of the class has $[K:\Q]\le4$.

\begin{lemma}[unit-point refinement]\label{canx:exact:lem:refine}
At a non-critical residue root, the Nahm point $X_P$ modulo $7^{d+1}$ depends only on the parameters modulo $7^d$ (because $\partial S/\partial(\text{parameter})$ is a $7$-adic logarithm of a unit, hence in $7\Z_7$, where $S$ is the logarithmic form of the Nahm relation). Consequently, for a datum $p^*=p+7^du$ in the class of $p$ modulo $7^d$, $Q_j(X_P;p^*)\equiv Q_j(\hat X_{d+1};p)+7^d\nabla_pQ_j(\hat X_{d+1};p)\cdot u\pmod{7^{d+1}}$, an affine condition on the digit vector $u\in\F_7^3$ over $\F_{7^g}$, for each $j$. A residue root whose constant term has valuation $<d$ carries no place for any datum of the class; a set of residue roots whose joint affine system is inconsistent cannot all carry places. (Critical roots are kept without conditions.)
\end{lemma}
\begin{proof}
Work in the unramified ring containing the Teichm\"uller lift of the
chosen noncritical residue. In the integral matrix chart the
logarithmic equations are
\[
 S_T=(1-h)\log X+h\log(1+X)=0,\qquad
 S_H=\log\frac{1-X}{1+X}-(\tau-1)\log X=0.
\]
The logarithm of any unit of this ring belongs to its ideal generated
by seven. The derivative with respect to the matrix parameter is
therefore divisible by seven; the derivative in $X$ is a unit at a
noncritical residue. The shifts do not occur in these equations.
For a matrix pole, write $h=u/7^k$ or $\tau=u/7^k$ and multiply the
corresponding equation by $7^k$. The normalized equations become
\[
 (7^k-u)\log X+u\log(1+X)=0,\qquad
 7^k\log\frac{1-X}{1+X}-(u-7^k)\log X=0.
\]
Their derivatives in $u$ are still logarithms of units. Modulo seven,
their derivatives in $X$ are $-\bar u/[X(1+X)]$ and $-\bar u/X$,
respectively, hence units on the allowed residue set. The implicit
function theorem now shows that a parameter change divisible by
$7^d$ changes the root only by a multiple of $7^{d+1}$. Taylor expansion
of the integral rescaled coefficient equations gives the stated affine
system. Its inconsistency is a necessary obstruction to simultaneous
lifts. If a residue does not satisfy the polynomial Nahm equation,
there is no place to count in the first instance.
\end{proof}
The finite refinement is as follows: a node (class modulo $7^d$) is \emph{resolved} when no subset of residue roots with count $>4$ has a consistent joint system; otherwise its children are the solutions $u$ (unique-child chains are recorded).

\subsection{Pole charts and their cells}\label{canx:exact:sec:cells}
Write the parameters as $p=(p_1,p_2,p_3)=(h,B,C)$ or $(\tau,B,C)$ and let $m=(m_1,m_2,m_3)\in\Z_{\ge0}^3$ be the pole orders, $p_i=u_i/7^{m_i}$ ($u_i$ a unit when $m_i>0$, any $7$-adic integer when $m_i=0$). For a direction $m$ the rescaled equations $7^{M_j}Q_j$ have $7$-integral coefficients and reduce modulo $7$ to the \emph{leading form} $L_j^{(m)}$: the terms whose parameter exponent $e$ maximises $m\cdot e$. The leading form depends only on the face of the Newton polytope of $Q_j$ selected by $m$ and on the sign pattern of $m$.

\begin{proposition}[exact cells]\label{canx:exact:prop:fan}
For nonzero pole directions, the common refinement of the normal fans of the three parameter Newton polytopes with the fan of the closed positive orthant has, for family T, $21$ cells (pattern, face triple), and for family H, $13$ cells. The origin $m=(0,0,0)$ is the integral chart and is treated separately. The cells with more than one monomial in every leading form and bad digit vectors are the rays $(0,k,k)$ and $(k,k,k)$ for T, and $(0,k,k)$ and $(k,0,k)$ for H.
\end{proposition}
\begin{proof}
Compute each Newton polytope from the exact support (vertices $10,6,10$ for T and $8,6,8$ for H), certify the facets by exact integer arithmetic (every point on the correct side of every facet plane, facet polygons as exact two-dimensional hulls, every edge in exactly two facets, facet normals at every vertex of rank three), then the Minkowski sum polytope ($10$ resp.\ $8$ vertices) and its normal fan; the rays of the refinement by the orthant are the facet normals in the orthant, the coordinate axes and the intersections of the two-dimensional cones with the coordinate planes; the sum of the rays of each closed cell is an interior representative. This yields all cells with an exact representative; the representatives' face triples are then computed exactly. Evaluate, for every cell and every admissible digit vector, the count bound of the three leading forms (Lemma~\ref{canx:exact:lem:length}; with a pole of the matrix parameter every residue root has length one): T has bad vectors only in the cells $(0,1,1)$ ($12$ vectors, bound $5$) and $(1,1,1)$ ($18$ vectors, all leading forms vanish identically); H only in $(0,1,1)$ ($12$ vectors, bound $6$) and $(1,0,1)$ ($6$ vectors, all leading forms vanish). All four are one-dimensional cells (rays), so the remaining analysis is a one-parameter family of charts indexed by $k\ge1$.
\end{proof}

\textbf{Uniformity in $k$.} Along a ray $km$, the rescaled equation is $\sum_{i\ge0}7^{ik}R_i(X;u)$ where $R_i$ collects the terms with $m\cdot e=M-i$ (layer $i$). A node of depth $d$ uses the equations modulo $7^{d+1}$ (the affine conditions of Lemma~\ref{canx:exact:lem:refine}); layer $i\ge1$ is invisible there iff $ik\ge d+1$. In a matrix-pole chart, the normalized logarithmic Nahm equation
of Lemma~\ref{canx:exact:lem:refine} differs from its limiting equation
by $7^k\log(u)$ for a local unit $u$ on the chosen unramified lift.
This difference lies in $7^{k+1}\mathcal O$; the unit derivative implies
that the corresponding roots agree modulo $7^{k+1}$.
Thus the same comparison applies when the point used in a node test
also varies with $k$. Hence the refinement tree for the pure leading forms (``$k=\infty$'', layer $0$ only) coincides with the tree for every $k$ at all depths $d\le k-1$. If the pure tree is resolved with deepest node $d_0$, every $k\ge d_0+1$ is resolved; the finitely many $k\le d_0$ are run individually.

\subsection{Family T}\label{canx:exact:sec:T}
\textbf{Integral chart.} The digit tree on $(h,B,C)\in\Z_7^3$ (all $343$ residue classes, depth $\le5$) is resolved in $341$ classes. The two classes $(\bar h,\bar B,\bar C)=(0,2,5),(0,5,2)$ are identically vanishing nodes at every depth, with a unique child at every level: the three reduced equations vanish identically in $X$ exactly on the locus $E=\{h=0,\ B+C=0,\ 240C^4-120C^2+7=0\}$ (Gr\"obner basis $[B+C,\ C^4-C^2/2+7/240]$ of the common vanishing ideal at $h=0$; $Q_2\propto(B+C)(4B^2-4BC+4C^2-1)$, $Q_3\propto120B^4-60B^2+120C^4-60C^2+7$, $Q_4\propto(B+C)(48B^4-\dots+7)$ times powers of $X$), and the recorded chain nodes are the Hensel roots of the quartic modulo $7^d$ ($d\le5$). $E$ has no rational point ($C^2=\tfrac12(\tfrac12\pm\sqrt{2/15})$ is irrational) and $h=0$ is outside the family. A rational datum of one of the two classes therefore leaves the chain at some depth $d\ge1$: at its node the reduced equations are the first-order forms $\sum_iv_i\,\partial Q_j/\partial p_i(X;p^*)$ modulo $7$ with $v\ne0$ (the second-order terms are $O(7^{2d})\subset O(7^{d+1})$ and $Q_j(X;p^*)\equiv0$). For all $342$ nonzero $v$ and both points the count bound is at most $2$. Hence $[K:\Q]\le4$ on the integral chart.

\textbf{Cell $(0,k,k)$} ($v(B)=v(C)=-k$, $h$ integral).
There are twelve residue vectors with count five, namely
$(\bar h,\bar B,\bar C)=(h_0,\beta,-\beta)$ with
$h_0\in\{2,4\}$ and $\beta\in\F_7^\times$. Their support is:
\[
\begin{array}{c|ccc}
 h_0&\text{rational simple}&\text{rational critical}&\text{quadratic simple}\\\hline
 2&X+3&X+6&X^2+4X+1\\
 4&X+5&X+4&X^2+X+3.
\end{array}
\]
The critical root has local length two; the other roots have length one.
Apply Lemma~\ref{canx:exact:lem:refine} at depth one. The joint affine
system for the rational simple point and the quadratic point is
inconsistent in every one of the twelve classes. This holds for the
leading forms, which suffice for $k\ge2$, and for the full first
correction when $k=1$. Removing their simultaneous occurrence lowers
the total count from five to at most four.

\begin{proposition}[cell $(k,k,k)$]\label{canx:exact:prop:kkk}
For $v(h)=v(B)=v(C)=-k$, $k\ge1$: $[K:\Q]\le1$ or no datum exists; in particular $[K:\Q]\le4$.
\end{proposition}
\begin{proof}
Exact algebra: along $(1,1,1)$,
\[ L_2=\tfrac{h^7}{12}(B+C+h)\,q_2,\quad L_3=-\tfrac{h^{11}}{360}(2X+1)\,q_3,\quad L_4=-\tfrac{h^{15}}{720}(B+C+h)(6X^2+6X+1)\,q_4, \]
$q_2=2B^2-2BC+Bh+2C^2+Ch$, $q_3=-15B^4-30B^3h-15B^2h^2-15C^4-30C^3h-15C^2h^2+h^4$, $q_4$ is the quartic form uniquely defined by the displayed factorization of the maximal-weight part of $Q_4$. Let $s=B+C+h$. (i) If $\bar s\ne0$ (unit digits), or $\bar s=0$ and $\bar B/\bar h\notin\{1,3,5\}$, the leading-form count is at most $1$ ($6X^2+6X+1$ is irreducible modulo $7$ and has no root in common with $2X+1$; on $\bar s=0$, $q_2=6B^2+6Bh+h^2$ is irreducible modulo $7$, hence a unit, and $q_3\equiv0$ exactly for $\bar B/\bar h\in\{1,3,5\}$; these are the $18$ identically vanishing vectors). (ii) If $\bar s=0$ then $v(L_2)=-9k+v(s)$, while every lower-weight term of $Q_2$ has valuation $\ge-9k$; $L_2$ is constant in $X$, so if $v(s)<0$ then $Q_2(X_P)\ne0$: no place. Hence $s\in\Z_7$. (iii) Substitute $C=\sigma-B-h$ ($\sigma=s$ integral). The leading forms along $(k,k,0)$ in $(h,B,\sigma)$ are $\Lambda_2=\tfrac{h^7}{12}(\sigma-1)(6B^2+6Bh+h^2)$, $\Lambda_3=L_3|_{C=-B-h}$, $\Lambda_4=-\tfrac{h^{15}}{720}(\sigma-1)(6X^2+6X+1)q_3|_{C=-B-h}$ (weights $9,15,19$). Now $v(\Lambda_2)=-9k+v(\sigma-1)$ against $\ge-8k$ for the rest of $Q_2$: if $v(\sigma-1)<k$ the constant term dominates and no place exists. Hence $\sigma=1+7^k\tau'$ with $\tau'\in\Z_7$. (iv) In the chart $h=u_1/7^k$, $B=u_2/7^k$, $C=1-B-h+7^k\tau'$ the rescaled equations are $\sum_i7^{ik}R_i(X;u_1,u_2,\tau')$ with layer weights $8,15,18$; the layer-zero forms at the $18$ residue vectors and all seven residues of $\tau'$ have count bound $\le1$ ($120$ cases with bound $0$, $6$ with bound $1$). Layer $i\ge1$ is invisible modulo $7$ for every $k\ge1$.
\end{proof}
Together with Proposition~\ref{canx:exact:prop:fan} and the integral chart this proves Theorem~\ref{canx:exact:thm:main} for family T.

\subsection{Family H}\label{canx:exact:sec:H}
\textbf{Integral chart.} All $343$ classes of $(\tau,B,C)\in\Z_7^3$ are resolved at depth $\le2$.

\textbf{Cell $(k,0,k)$} ($\tau$ and $C$ poles, $B$ integral). The leading forms vanish identically on the six vectors $(\bar\tau,\bar B,\bar C)=(t,3,4t)$; their common vanishing locus (as polynomials in $B$, $\gamma=C/\tau$ after dividing by the powers of $\tau$) is empty over $\overline\Q$ (Gr\"obner basis $[1]$), so no chain can persist: the pure-leading-form tree is resolved with deepest node $2$ (valid for $k\ge3$), and $k=1,2$ are resolved at depth $2$.

\begin{proposition}[The remaining H ray]\label{canx:exact:prop:0kk}
For family H with $v(B)=v(C)=-k$, $k\ge1$, and integral $\tau$,
one has $[K:\Q]\le4$.
\end{proposition}
\begin{proof}
The leading forms are $D_N^{j+1}G_j$, where $G_j$ is homogeneous of
degree $j+1$ in the shifts. Their only counts above four are the two
scaling orbits $(\bar\tau,\overline{B/C})=(4,3),(5,1)$.
In each, the determinant-critical irreducible quadratic contributes
at most four and one further noncritical quadratic could contribute
two. Proposition~\ref{deg:Hrepair} excludes that noncritical quadratic
uniformly at every pole order and every parameter depth, using the
exact base points and a first-order rank obstruction. Thus the total
is at most four. No finite-depth termination claim is used here.
\end{proof}
Together with the integral chart, the $(k,0,k)$ ray and
Proposition~\ref{canx:exact:prop:fan}, this proves
Theorem~\ref{canx:exact:thm:main} for family H.

%% file: sections/canonical_expansion.tex
\section{Canonical charts: complete local arguments}\label{canx:section}
We give the canonical arguments used in Section~\ref{deg:sec:canonical}.
Throughout, the coefficient polynomials are defined by the recurrence of
Section~\ref{deg:sec:coefficients}; thus none of the finite tests below
requires a precomputed polynomial as input. The field is generated by both
coordinates. Local contributions always mean sums of local degrees over
all places in the indicated region.

\subsection{Canonical equations and unit information}\label{canx:setup}
Unless a subsection states a more general matrix, the data are
\[
 A=\begin{pmatrix}1&-h\\-h&\tau\end{pmatrix},\qquad
 h,\tau\in\Q,\quad h>0,\quad h\ne1,\quad\tau\ne1,\quad\tau-h^2\ne0,
 \qquad K=\Q(x,y),
\]
with a real signed algebraic solution, nonzero Nahm letters, rational
linear shifts $B,C$ and rational $r=\lambda_1$. We impose
$\lambda_2=\lambda_3=\lambda_4=0$ and $D\ne0$.
Valuations at places above seven are normalized by $v(7)=1$.
A statement invoking the initial bound assumes $[K:\Q]\le24$;
Theorem~\ref{d24:thm:main} supplies it in the modular application.
The transformations and coefficient convention are those of
Proposition~\ref{deg:prop:coefficient-transport}.

\subsection{The canonical units follow from the equations}\label{canx:tau:app:units}
Here the signed canonical equations mean
\[
 \log|1-x|=\log|x|-h\log|y|,\qquad
 \log|1-y|=-h\log|x|+\tau\log|y|.
\]
The point is real algebraic and none of $x,y,1-x,1-y$ is zero.
Multiplying by twice a common denominator gives Laurent polynomial
identities over $\Q$, hence identities at every conjugate and finite place.
If $h,\tau\in\Z_7$, the clearing exponent can be chosen prime to seven.
The real condition specifies the signed equations; positivity of the
transported matrix or point is unnecessary.

\begin{proposition}\label{canx:tau:prop:units}
Suppose $h>0$, $\tau\ne1$ and $\tau-h^2\ne0$. Then $X,y$ are global
algebraic units, $1+X,1-y$ are algebraic integers, and at every finite place
\[
 \operatorname{v}(1-y)=h\operatorname{v}(1+X).
\]
If also $h\in\Z_7^\times$ and $h\ne1$, all four Nahm letters are
local units at every place above seven.
\end{proposition}
\begin{proof}
Put $u=\operatorname{v}(x)$, $w=\operatorname{v}(y)$, $f=\operatorname{v}(1-x)$ and $g=\operatorname{v}(1-y)$.
The equations give $f=u-hw$ and $g=-hu+\tau w$.
For any nonzero $z\ne1$, the valuation of $1-z$ is zero when
$\operatorname{v}(z)>0$, equals $\operatorname{v}(z)$ when $\operatorname{v}(z)<0$, and is nonnegative
when $\operatorname{v}(z)=0$.

If $u>0$, then $f=0$ and $w=u/h>0$. Thus $g=0$, and the second
equation gives $(\tau-h^2)w=0$, a contradiction.
If $u<0$, then $f=u$ and the first equation forces $w=0$; the second
gives $g=-hu>0$.
If $u=0$, then $f=-hw\ge0$, so $w\le0$. A negative $w$ would give
$g=w=\tau w$, contradicting $\tau\ne1$. Thus $w=0$ and $f=g=0$.
This exhausts all possibilities. Always
\[
 u\le0,\quad w=0,\quad f=u,\quad g=-hu.
\]
Since $X=x/(1-x)$ and $1+X=1/(1-x)$, the global assertions follow.

For the final assertion suppose a nonunit place above seven exists.
Write $h=M/N$ in lowest positive terms. Set $z=1+X$ and $v=1-y$;
both have positive valuation. The first signed equation gives
$X^{2N}=y^{2M}$. The ratio $(-X)^N/y^M$ has square one and reduces to
one, so it equals one. Hence
\[
 (1-z)^N=(1-v)^M.
\]
When $M,N$ are seven-units, the two sides minus one have valuations
$\operatorname{v}(z)$ and $\operatorname{v}(v)$, respectively. Therefore
$\operatorname{v}(z)=\operatorname{v}(v)=h\operatorname{v}(z)$, contrary to $h\ne1$.
\end{proof}

The proof also gives the principal branch used below: if $h=M/N$ in
lowest terms and $7\nmid N$, then $-X=y^{M/N}$ on the principal-unit
disc. The $N$th-power map is injective there.

\subsection{Exact finite computations in the canonical charts}\label{canx:finite}
Here is a common specification for the finite coefficient calculations.
Form $\Psi_1,\ldots,\Psi_4$ by the rational recurrence of
Section~\ref{deg:sec:coefficients}, specialize the displayed matrix and
shift residues, and reduce in $\F_7$. Their coefficient denominators are
prime to seven. Put
\[
 F_0=\Psi_1-rD^3,\quad F_1=\Psi_2,\quad F_2=\Psi_3,\quad F_3=\Psi_4,
 \qquad U=XY(1+X)(1+Y)D.
\]
For determinant-unit points use the saturated ideal
\[
 I=(F_0,F_1,F_2,F_3,1-tU)\cap\F_7[X,Y].
\]
A Gröbner basis determines its dimension and, when the dimension is
zero, its radical and every closed point. Equivalently, factor elimination
polynomials in $X$, specialize in each exact finite residue field and
factor the common $Y$-gcd. Vertical components are handled by the full
ideal or by checking leading coefficients separately. Counting only
$\F_7$-rational solutions would not be sufficient.

For a closed point $P$ with full coordinate field $k_P$ record
$d_P=[k_P:\F_7]$ and
\[
 n_X(P)=N_{k_P/\F_7}(X),\qquad
 n_y(P)=N_{k_P/\F_7}\bigl(Y/(1+Y)\bigr).
\]
Each determinant-unit point can occur at most once, by the unit Nahm
Jacobian. For every subset $S$ of the finite closed points and every
allowed determinant-centre multiplicity $m$, keep precisely the choices
for which the product of each available unit norm belongs to
$\{1,-1\}$. The bound of the fibre is the largest
$\sum_{P\in S}d_P+m$ kept; the empty subset and $m=0$ are allowed.
At a repeated canonical centre the residue is $(X,Y)=(1,3)$, so its
two norm factors are $(1,-1)^m$. Its multiplicity is bounded by the
length of the one-variable Nahm and coefficient algebra in Proposition~\ref{cuspx:census:prop:rep}; simple determinant centres are absent
for integral scalar. Positive-dimensional fibres are set aside for the
curve arguments, rather than assigned a finite-point bound.

\begin{lemma}[Three canonical censuses]\label{canx:census}
The preceding procedure has the following results.
\begin{enumerate}[label=(\roman*),leftmargin=*]
\item For $a=1$, $b\ne0$, $(b,d)\ne(6,1)$ and
$(d,B,C,r)\in\F_7^4$, there are $41\cdot49\cdot7=14063$ fibres.
Exactly $64$ have a curve. The remaining $13999$ finite fibres, including
the permitted repeated determinant-centre contributions, have bounds
\[
\begin{array}{c|rrrrrr}
\text{bound}&0&1&2&3&4&5\\\hline
\text{number}&13223&20&52&592&108&4.
\end{array}
\]
The four bound-five fibres are
$(a,b,d,B,C;r)=(1,6,d,4,C;1)$ with $d=4,5$ and $C=3,4$.
\item For $a=1,b=0$, $(d,B,C,r)\in\F_7^4$, there are $2401$ fibres,
of which $30$ have a curve. The $2371$ finite fibres have bounds
\[
\begin{array}{c|rrrrr}
\text{bound}&0&1&2&4&6\\\hline
\text{number}&2339&6&18&6&2.
\end{array}
\]
The bound-six fibres are $(1,0,3,1,3;6)$ and $(1,0,5,1,1;6)$.
\item For the inverse matrix residue $\operatorname{diag}(1,0)$,
$(J,K_{\rm sh},R)\in\F_7^3$, replace the norm pair by
$X_I,1+Y_I$. Of $343$ fibres, $15$ have a curve, $324$ finite fibres
have bound zero and four finite fibres have bound one.
\end{enumerate}
\end{lemma}
\begin{proof}[Finite verification and its coverage]
The domains above exhaust the residue parameters. Apply the displayed
ideal construction in each fibre. For (i), the only repeated admissible
determinant residues are $(a,b,d)=(1,3,5),(1,5,6)$, both with centre
$(X,Y)=(1,3)$. Their Nahm length is four. Eliminate the first logarithmic
Nahm equation in the local parameter of the second, to order four;
restrict $R,J_2,J_3$, defined in Section~\ref{cuspx:census:sec:setting}, to that parameter and take the minimum of four and their first nonzero orders. This is the multiplicity
used in the subset calculation, so it includes all ramified lifts.
Parts (ii) and (iii) have no determinant-centre contribution in the
unit chart. Factor the ideals, enumerate the subsets with the two
norm tests and tabulate their maximum degrees. The resulting counts
are those displayed. The zero-dimensional ideal calculation verifies
that no additional closed point of higher residue degree is omitted.
The $64$, $30$ and $15$ curve fibres are not counted among these finite
maxima. Their bounds come from Theorem~\ref{cuspx:flat-complete} and, for (iii), the inverse-chart analysis below.
\end{proof}

For first-normal and weighted calculations, a similarly exact recipe
removes any dependence on a stored array. Substitute the stated
parametrization in $F_0,\ldots,F_3$, clear only denominators known to be
units on the chart, and extract the coefficient of the specified first
nonzero weight. If $q$ is the normal coordinate and $p$ the rational
transverse vector, write those terms as
$n_iq+\sum L_{ij}p_j$. Enumerate projective directions by making their
first nonzero coordinate equal to one; thus $\mathbb P^1(\F_7)$,
$\mathbb P^2(\F_7)$ and $\mathbb P^3(\F_7)$ have $8$, $57$ and $400$
elements. For each direction take the gcd of
$n_iL_j-n_jL_i$ after deleting precisely the forbidden coordinate
factors specified in that chart. Check the normal-column gcd, the
squarefreeness of the remaining support and the relevant Jacobian
minors. These are the acceptance tests for the first-normal counts
below. A zero wedge vector is a flat direction and requires the next
argument; it is never discarded.

\subsection{The six exceptional integral-shift fibres}\label{canx:refine}
The possible bound-five and bound-six fibres in Lemma~\ref{canx:census}
are removed by one more digit. At a determinant-unit point the signed
Nahm equations have a unit Jacobian. Their parameter derivatives are
principal-unit logarithms, hence divisible by seven. Consequently the
lift of the point modulo $49$ depends only on the parameter residues.
Insert this Hensel lift into $F_0,\ldots,F_3$ with parameters
$p_0+7p_1$. Reduction after division by seven gives an affine linear
system over $\F_7$ in $p_1$ and the next scalar digit. Expand each
residue-field equation in an $\F_7$-basis. Its consistency is a necessary
condition for a lift; stacking systems tests simultaneous occurrences.

In the two bound-six fibres the four cubic points give inconsistent
systems, with coefficient/augmented ranks $4/5$ or $5/6$.
In each of the four bound-five fibres the quartic point has ranks $5/6$;
the rational point has seven solutions for the next digits. Thus the
former fibres have contribution zero and the latter at most one.
All entries of the affine systems are defined by the displayed Hensel
and coefficient procedure. This is a finite test of six specified
fibres, with no bound on the rational parameters or their later digits.
The remaining finite fibres already have bound at most four.

\begin{theorem}[Canonical integral linear shifts]\label{canx:integral}
Assume the canonical hypotheses of Section~\ref{canx:setup},
$B,C\in\Z_7$ and $\tau\in\Z_7$.
\begin{enumerate}[label=(\roman*),leftmargin=*]
\item If $v_7(h)=0$ and $(\bar h,\bar\tau)\ne(1,1)$, then
$[K:\Q]\le4$, for every rational $r$.
\item If $v_7(h)>0$ and $[K:\Q]\le24$, then $[K:\Q]\le4$.
\end{enumerate}
\end{theorem}
\begin{proof}
For (i), Proposition~\ref{canx:tau:prop:units} makes all four Nahm
letters local units at seven. If $r$ is integral, Lemma~\ref{canx:census}(i)
and the six-fibre refinement give bound four on every finite fibre;
Theorem~\ref{cuspx:flat-complete} supplies the same bound on the curve
fibres. This census includes repeated determinant centres and uses
only the two available global-unit norms.

If $r$ has a pole, no determinant-unit place is possible. The admissible
part of the determinant polynomial
\[
 q(H)=H\{\tau(\tau-1)H^2-h^2(-h+2\tau-1)H+h^3(h-1)\}
\]
has degree at most two; $H=0$ is inadmissible. The simple-root bounds
of Corollary~\ref{cuspx:census:cor:support} and
Section~\ref{cuspx:simple-poles} give total degree at most four.
The only repeated-root residues are $(1,3,5),(1,5,6)$ in $(a,b,d)$;
both have one rational square root, and
Proposition~\ref{cuspx:census:prop:rep} gives length at most four.
These are norm-free conclusions, so their hypotheses apply here.

For (ii), the analytic disc and first-coefficient argument below show
that integral $r$ allows no nonunit-letter place. At a local-unit place
with nonunit determinant, reduction gives $Y=-\tau$ and
\[
 \Psi_2(X,-\tau;1,0,\tau,B,C)
 =\frac5{16}\tau^4(\tau-1)^4(X+1)^6.
\]
For $\bar\tau\ne0,1$ this is nonzero; the other two values have an
inadmissible letter. Thus every local-unit place has unit determinant
and forces $r$ integral. If $r$ had a pole, all places would instead be
nonunit-letter places, whose total degree is at most one by the same
analytic argument for integral shifts. Then $y$ would be a rational
global unit, hence $\pm1$, neither of which is an admissible nonunit
point above seven. Hence $r$ is integral and all places have unit
determinant. Lemma~\ref{canx:census}(ii), the six-fibre refinement and
Theorem~\ref{cuspx:flat-complete} give bound four.
\end{proof}

\subsection{The determinant-unit part: complete leading census}
For this section let $A=(\begin{smallmatrix}a&b\\b&d\end{smallmatrix})$ be seven-integral with $\bar a=1$. A place is called a unit place here if $x,y,1-x,1-y$ and $D$ are local units. Only these places are counted in this section.

Suppose a shift has a pole. Set $k=-\min(v_7(B),v_7(C))>0$. The degree in the shifts of $\Psi_j$ is at most $j+1$. Thus at a unit place $v_7(r)\ge-2k$. The highest shift-degree pieces, the eight projective shift directions in $\mathbb P^1(\F_7)$, and the rational leading scalar give a complete finite calculation.

All $49\cdot8=392$ matrix/direction configurations were processed. Of these, 372 have finite scalar support of total residue degree at most one. Thirty-six have a singular \emph{coefficient} fibre; they are still fully enumerated. Their lifting multiplicity is one because the \emph{Nahm} Jacobian is a unit. There are no unprocessed vertical components.

The remaining twenty directions have a curve. Their 140 scalar fibres consist of 120 empty fibres and twenty pure curves, with no isolated points outside the curve. The types are
\[
 X+Y+1=0\ (7),\qquad XY-1=0\ (7),\qquad XY+X+Y=0\ (6).
\]
The six curves of the third type have matrix residue $(1,-1,d)$ with $d\ne1$. Full inversion has unit determinant and takes them to the first type. Thus there are only two normal-coordinate calculations.

\subsubsection*{Why local units give unique lifts}
Clear the rational matrix exponents with an integer $N$ prime to seven. Up to unit row and column factors, the Jacobian of the two resulting Nahm equations is $N(A+\operatorname{diag}(X,Y))$. Its determinant is a unit at a unit place. Hensel uniqueness gives one unramified lift per geometric residue solution. Consequently a closed point of residue degree $f$ contributes at most $f$ to $[K:\Q]$. This justifies counting squarefree support; no separability assumption on an arbitrary reduced minimal polynomial is used.

\begin{proposition}[Determinant-unit bound]\label{canx:up:prop:unitpart}
For an integral matrix with $\bar a=1$, arbitrary rational shifts with a pole, and the stated coefficient conditions, the total contribution of all unit places to the field degree is at most one. Four global units are not required.
\end{proposition}
The finite fibres prove the zero-dimensional cases. The following arguments prove the curve cases at every valuation depth.

\subsection{The arithmetic improvement on the two curves}
At the unit places the local extension is unramified. The seven-adic logarithm of any local unit belongs to $7\mathcal O$: equivalently, first raise it to the residue-field order minus one, which removes its Teichmuller factor and all signs. Normalize $v(7)=1$.

\begin{lemma}[One additional power of seven]\label{canx:up:lem:extra}
On a residue curve, the following inequalities hold.
\begin{enumerate}[leftmargin=*]
\item For $X+Y+1=0$, put $\delta=d-(1+2b-a)$. Then
\[
 v(X+Y+1)\ge 1+v(\delta).
\]
If $\delta=0$, the curve relation holds exactly.
\item For $XY-1=0$, choose an integral rational $w$ and set
$u=b-1-aw$, $v_0=d-aw^2$. Then
\[
 v(XY-1)\ge1+\min\{v(u),v(v_0)\}.
\]
If $u=v_0=0$, the curve relation holds exactly.
\end{enumerate}
Here the matrix is close to the corresponding surface, and $D$ is a unit.
\end{lemma}
\begin{proof}
The first curve is $xy=1$. Subtracting the two logarithmic Nahm equations gives
\[
 F=\log\!\left(\frac{y(1-x)}{1-y}\right)
 -(a-b)\log(xy)=-\delta\log y.
\]
At $y=1/x$, $F=0$ and $F_y=-x(a+X-b)$ is a unit, since on the curve $D=-(a+X-b)^2$. Hence the valuation of $F$ equals that of the curve displacement. The right side has valuation at least $v(\delta)+1$.

The second curve is $x+y=1$. The second Nahm logarithmic equation minus $w$ times the first gives
\[
 F=\log\!\left(\frac{1-y}{x}\right)
 -w\log\!\left(\frac{1-x}{y}\right)
 =u\log x+(v_0-wu)\log y.
\]
At $y=1-x$, $F_y=(1+X)(wX-1)/X$ is a unit: on the surface and curve, $D=a(wX-1)^2/X$. The same valuation argument proves the claim. In both cases an exact zero on the right gives the exact curve relation by the unit derivative. The logarithm converges on the residue disk after the torsion factor is removed, so higher Taylor terms have strictly higher valuation.
\end{proof}
This inequality permits us to discard the curve-normal displacement at the first nonzero parameter order. The four coefficient rows give support of degree at most one in every nonflat direction. Ten centres retain one flat direction each; the weighted argument below treats them.

\subsection{All flat directions reduce to one exact family}
The ten flat centres are carried by integral coordinate exchange or the first-coordinate transformation to the sum-curve family with
\[
 \bar a-\bar b=1,\quad \bar d=\bar b,\quad
 \overline{C/B}=1.
\]
The two product centres have residues $(a,b,d;B:C)=(1,0,1;1:-1)$ and $(1,1,0;1:0)$. Applying the first-coordinate transformation, and in the first case exchanging coordinates, gives the sum type. Its shift map is
\[
 B'=\frac12-\frac Ba,\qquad C'=C-\frac ba B,
\]
and $a$ is a unit in these two cases. Thus poles and admissibility are preserved in the required charts.

The exact characteristic-zero flat locus is
\[
 a=b+1,\quad d=b,\quad C=B+\frac12,\quad r=\frac1{12},
 \qquad Y=-1-X.
\]
Direct substitution in the coefficient recurrence gives, for arbitrary $b,B$,
\[
 D=-(X+1)^2,\qquad \Psi_1=-\frac{(X+1)^6}{12},
 \qquad \Psi_2=\Psi_3=\Psi_4=0.
\]
These are polynomial identities over $\Q$, not identities only modulo seven.

Set $t=1/B$, $k=v(t)>0$, and use the transverse variables
\[
 a=b+1+u,\quad d=b-u+\delta,\quad
 C=B+\frac12+s,\quad r=\frac1{12}+z.
\]
The surviving first-normal direction ensures
\[
 v(u)\ge1,\quad v(\delta)>k,\quad v(s)\ge1,\quad
 v(t^2z)>k.
\]
Here $z$ itself is allowed to have a pole. If not all four transverse differences vanish, choose
\[
 m=\min\{v(\delta),v(tu),v(ts),v(t^2z)\}.
\]
After dividing these four quantities by $7^m$, their reductions form a nonzero vector $\eta$. Lemma~\ref{canx:up:lem:extra} gives $v(X+Y+1)>m$.

\paragraph{Why this covers arbitrary depth.}
In $t^{j+1}\Psi_j$, restricted to the sum curve, every transverse term is divisible by one of $\delta,tu,ts$; the first equation also contains $t^2z$. After the linear weighted terms are removed, each remainder has an additional factor from $(u,\delta,s,t)$. This follows from the exact flat-locus identity and the leading sum-surface identities, valid for arbitrary $a,b$. All these additional factors have positive valuation. Thus no omitted term reappears at a deeper scale, even when $z$ has a pole.

\subsection{The small determinant that closes the flat locus}
After removing only powers of the local units $X$ and $1+X$, the reduced four equations at the preceding weighted order are
\[
 M(X)\eta=0,\qquad
 \eta=(\eta_\delta,\eta_{tu},\eta_{ts},\eta_{t^2z})^T,
\]
where, over $\F_7$,
{\small
\[
M(X)=\begin{pmatrix}
4&6X&1+X&(1+X)^2\\
4&4+3X&4+4X&0\\
1&1+3X+X^2&1-X^2&0\\
5+5X+X^2&5+X+6X^2+2X^3&5+6X+6X^2+5X^3&0
\end{pmatrix}.
\]
}
The exact factorization is
\[
 \det M(X)=3X^2(X-1)(X+1)^3\quad\text{in }\F_7[X].
\]
Equivalently, the last three rows and first three columns have determinant
$4X^2(X-1)(X+1)$. Since $X$ and $1+X$ are units and $\eta\ne0$, every surviving residue point has $X=1$. At $X=1$ the unique projective direction is
\[
 \eta=(1,4,3,2).
\]
This proves degree at most one for the unit part, by the Nahm Jacobian. The explicit census agrees: for each of the seven values of $\bar b$, 399 of the 400 projective directions are empty and one has degree-one support. There are no zero images.

If all four transverse differences vanish exactly, Lemma~\ref{canx:up:lem:extra} gives $xy=1$. The real signed first Nahm equation then gives $|1-x|=|x|$, so $x=1/2$ and $y=2$. The field is rational in this case too. Thus the argument has no unprocessed exact-locus case.

\subsubsection*{Uniformity of the bound}
The bound uses neither four global norm conditions nor a norm test on subsets of residue roots. All field places contributing to this unit part see the same rational parameter direction. The unique-lift argument prevents repeated embeddings over one residue point from being silently discarded.

\subsection{The degree bound forces nonunit points into an analytic disc}
First take $h\in7\Z_7\setminus\{0\}$. At a nonunit place above seven put
\[
 u=v(1+X)>0,\qquad w=v(1-y)=hu>0,\qquad v(7)=1.
\]
Write $h=M/N$ in lowest terms, with $M,N>0$. Then $7\mid M$, $7\nmid N$. If $e$ is the ramification index, the integer valuations satisfy
\[
 eu=N\ell,\qquad ew=M\ell\qquad(\ell\in\Z_{>0}).
\]
Since $e\le[K:\Q]\le24$,
\begin{equation}\label{canx:bound:eq:depth}
 w\ge\frac M e\ge\frac7{24}>\frac16.
\end{equation}
This is the only use of the degree-24 input in the new analytic argument.

Set $v_0=1-y$, so $X\equiv-1$, $y\equiv1$. The cleared first Nahm equation has the sign compatible with these residues, and becomes
\[
 (-X)^N=(1-v_0)^M.
\]
Raising to the $N$th power is injective on principal units, because $7\nmid N$. Hence
\[
 X=-(1-v_0)^h=-1+h v_0 U(v_0),\qquad
 U(v_0)=\frac{1-(1-v_0)^h}{h v_0}
       =\sum_{k\ge0}\frac{(-1)^k}{k+1}\binom{h-1}{k}v_0^k.
\]
The binomial coefficients are seven-integral. For $k\ge1$ one has
$v_7(k+1)\le k/6$, so \eqref{canx:bound:eq:depth} implies $U(v_0)\equiv1$ in the residue field. The same expansion is valid in every completion and in $\mathbb C_7$.

\begin{lemma}\label{canx:bound:lem:disc}
For the integral-shift case, every nonunit embedding has $v(v_0)\ge1$.
\end{lemma}
\begin{proof}
Suppose $1/6<w<1$. Write $L=1+(\tau-1)v_0$. If $\tau$ is integral, $v(L)=0$. Otherwise, putting $T=-v_7(\tau)\in\Z_{>0}$, one has $v(L)=w-T<0$. The polynomial identity on the next page shows that
\[
 v\!\left(\frac{\Psi_1}{D^3}\right)=-v_7(h)-w.
\]
Indeed its leading term is a unit times $L^3$, while its error has strictly larger valuation; the denominator has valuation $3v(L)$. This cannot be the valuation of the rational number $r$, since $v_7(h)$ is integral and $0<w<1$. Thus no such $w$ occurs.
\end{proof}
In particular $v_0=7V$ with $V$ in the valuation ring of $\mathbb C_7$.

\subsection{An exact first-coefficient identity and a uniform root bound}
Use independent variables $v,U,h,\tau,B,C$, and substitute
\[
 X=-1+hvU,\qquad Y=v^{-1}-1.
\]
The coefficient recurrence of Section~\ref{deg:sec:coefficients} gives the exact polynomial identities
\begin{equation}\label{canx:bound:eq:normal}
\begin{aligned}
 D&=hD_0, & D_0&=UL-h, & L&=1+(\tau-1)v,\\
 \Psi_1&=\frac{h^2}{v}P, & P&=U^2F_BL^3+hQ, & F_B&=\frac{12B^2-1}{24}.
\end{aligned}
\end{equation}
There are 155 monomials in $P$ and 135 in $Q$. Every coefficient of $Q$ is integral at seven. Every monomial $v^i\tau^j\cdots$ in $Q$ satisfies
\[
 i\ge j,\qquad 0\le j\le3.
\]
These are exact characteristic-zero identities, checked term by term. For integral $B$, $F_B$ is always a seven-adic unit: the equation $12B^2=1$ would require the nonsquare $B^2=3$ in $\F_7$. Its seven residue values are $2,6,4,3,3,4,6$.

The first coefficient equation is therefore
\begin{equation}\label{canx:bound:eq:F}
 P(v,U(v),h,\tau,B,C)-rh vD_0(v,U(v))^3=0.
\end{equation}
On putting $v=7V$, the series $U(7V)$ has integral coefficients tending to zero, and reduces to $1$: its coefficient of $V^k$ has valuation at least $k-v_7(k+1)\ge1$ for $k\ge1$.

Put $t_0=\min(0,1+v_7(\tau-1))$, $c=7^{-t_0}$ and $\widetilde L=cL$. This linear polynomial has integral coefficients and nonzero reduction $\ell$. The monomial inequalities above imply that $c^3Q(7V,U(7V),\ldots)$ is integral. Hence
\[
 c^3P\equiv F_B\ell^3,\qquad cD_0\equiv\ell\pmod7.
\]
Let $\alpha=7rh$. After scaling \eqref{canx:bound:eq:F} to have integral primitive coefficients, its reduction is respectively
\[
 F_B\ell^3,\qquad
 \ell^3(F_B-\bar\alpha V),\qquad
 -V\ell^3
\]
when $v_7(\alpha)>0$, $=0$, or $<0$. Each is nonzero. Their degrees are at most three, four, four.

Strassmann's theorem, applied over $\mathbb C_7$, now bounds the number of roots on the whole closed disc by three when $r\in\Z_7$, and by four in general. When $\tau$ is integral, $\ell=1$, so these bounds sharpen to zero for integral $r$ and at most one for a pole of $r$. This is an all-orders analytic bound, not a truncated-series test. Since $X$ is uniquely determined by $v=1-y$, distinct embeddings of $K$ give distinct roots. Thus the bounds count the \emph{total} contribution of nonunit places, including ramification and residue degrees.

\subsection{The integral matrix determinant boundary has length at most two}
Still assume $h\in7\Z_7$, $\tau\in\Z_7$. At a place where the four letters are units but $D$ is not, put
\[
 H=1+X,\quad Y=(h^2+D)/H-\tau,\quad
 V=C+hB/H,\quad \delta=v(D)>0,\quad k=-v_7(B),\quad\kappa=-v(V).
\]
Here $H$ is a unit. Necessarily $\bar\tau\notin\{0,1\}$ and $\bar Y=-\bar\tau$. Completing the quadratic shift term gives
\[
 \lambda_{1,\mathrm{quadratic}}=\frac{B^2}{2H}+\frac{HV^2}{2D}.
\]
After clearing the harmless factor $H^{6j}$, every monomial $h^iD^\ell B^uV^v$ of $\Psi_2,\Psi_3$ satisfies
\[
 \ell\ge2u+v.
\]
The same inequality holds for $\Psi_1-B^2D^3/(2H)$, whose $B$-degree is at most one. The exact checks have respectively 63, 812 and 5,077 monomials, all with seven-integral coefficients.

\begin{lemma}\label{canx:close:lem:det}
At every such determinant place,
\[
 k\ge\max(2\delta,\kappa+\delta),\qquad
 r=\frac{B^2}{2H}+\text{terms of valuation greater than }-2k.
\]
Consequently all these places reduce to one rational point, and their total degree contribution is at most two.
\end{lemma}
\begin{proof}
If $\kappa\le\delta$ and $k<2\delta$, the monomial inequality makes every term containing $h$ or $B$ strictly smaller than the scalar singular part. Its two equations reduce, up to units, to
\begin{align*}
 P_2(t)&=t^3+18t^2+135t+405,\\
 P_3(t)&=\tfrac92t^4+144t^3+2025t^2+14580t+\tfrac{89505}{2}.
\end{align*}
But $P_3(t)-tP_2(t)=5$ modulo seven, so they cannot both vanish. If $\kappa>\delta$ and $k<\kappa+\delta$, the unique leading term of $\Psi_2$ is the $V^3$ term, with nonzero scalar cubic coefficient $\bar\tau(\bar\tau-1)$. This is also impossible.

The same monomial inequality for the first-coefficient remainder, using the resulting lower bounds on $k$, gives a strict valuation improvement over $-2k$. Hence $v_7(r)=-2k$ and
\[
 \bar H=\overline{B^2/(2r)}.
\]
This is the same rational residue at every determinant place. The first Nahm equation has a unit derivative in the regular coordinate. After its elimination, the second has a nonzero quadratic term proportional to $\bar\tau(\bar\tau-1)$; its local length is two. This bounds all lifts over that one point, not just one completion.
\end{proof}
The determinant-unit contribution is bounded by Proposition~\ref{canx:up:prop:unitpart}.

\subsection{A pole of \texorpdfstring{$\tau$}{tau}: arbitrary inverse shifts}
Now let $s=v_7(h)>0$ and $v_7(\tau)=-T<0$. Use the full inverse shifts $J,K_{\rm sh}$ and scalar $R$:
\[
 B=J-hK_{\rm sh},\quad C=-hJ+\tau K_{\rm sh},\qquad
 r=\frac{J^2-2hJK_{\rm sh}+\tau K_{\rm sh}^2}{2}+\frac1{12}-R.
\]
The inverse matrix is integral and reduces to $\left(\begin{smallmatrix}1&0\\0&0\end{smallmatrix}\right)$. At every local-unit place its determinant is a unit.

Substitute these expressions into the first equation and divide by the scalar $\tau^3$ and the analytic unit $U^3$. Its principal part on the disc $v(1-y)\ge7/24$ is
\begin{equation}\label{canx:close:eq:principal}
 v^3\left\{\frac{F_J}{U(v)}-\frac{hK_{\rm sh}^2}{2}
       +h\left(R+\frac{K_{\rm sh}^2}{2}\right)v\right\},
 \qquad F_J=\frac{12J^2-1}{24}.
\end{equation}
The remaining terms are controlled at every order, not just in a truncated series.

Put $A=\min(0,2v_7(J))$ and $E=s+2v_7(K_{\rm sh})$. If $A\ne E$, one constant in braces dominates, unless the linear term dominates instead. The primitive reduction has degree at most four. Exact monomial inequalities show that all omitted terms have strictly larger Gauss valuation than the relevant principal terms.

If $A=E$, then $s$ is even and $s\ge2$. When $v_7(R)+s<A$, the $R$ term dominates and still gives bound four. Otherwise retain two Taylor coefficients:
\[
 \frac1{U(v)}=1-\frac{1-h}{2}v+\frac{h^2-1}{12}v^2
             +\text{a tail of Gauss valuation greater than }2\rho.
\]
The $v^2$ coefficient is a seven-unit. Thus the coefficient of $v^5$ in \eqref{canx:close:eq:principal} has valuation exactly $A$, while every tail term has Gauss valuation greater than $A+5\rho$. Cancellation of the constant and linear terms cannot remove it. Strassmann therefore gives at most five nonunit embeddings.

\paragraph{All-orders and exact checks.}
For $n\ge3$, $v_7(n+1)\le(n-2)/4<(n-2)\rho$. The first inequality follows from $4e\le7^e-3$ for $e\ge1$, starting at equality for $e=1$. It bounds the tail of $U$ and hence its analytic inverse. There are 585 polynomial remainder inequalities in the unbalanced ranges (minimum margin $1/2$) and 428 in the balanced ranges (minimum margin $1/8$). A separate checker verifies all 1,013 as explicit nonnegative rational combinations of their cone inequalities, with a strictly positive constant. Their verification therefore does not depend on a linear-program solver. The cones contain every permitted integral valuation vector. Terms involving $R$ are checked separately against its principal monomial.

\paragraph{Construction of the Gauss inequalities.}
The finite verification in the preceding paragraph can be reconstructed
as follows. In the independent variables $(v,U,h,\tau,J,K_{\rm sh},R)$
form
\[
 \begin{split}
 F={}&P(v,U,h,\tau,J-hK_{\rm sh},-hJ+\tau K_{\rm sh})\\
 &-h v(UL-h)^3\left(\frac{J^2-2hJK_{\rm sh}+\tau K_{\rm sh}^2}{2}
                      +\frac1{12}-R\right),
 \end{split}
\]
and subtract
\[
 F_{\rm pr}=U^2\tau^3v^3F_J
 -\frac12hK_{\rm sh}^2U^3\tau^3v^3(1-v)
 +RhU^3\tau^3v^4.
\]
For a nonzero monomial $v^iU^u h^l\tau^dJ^jK_{\rm sh}^kR^q$
in $F-F_{\rm pr}$, the coefficient is integral at seven,
$d\le3$ and $q\le1$. Set $\rho=7/24$ and write
$(s,T,b,c)=(v_7(h),-v_7(\tau),v_7(J),v_7(K_{\rm sh}))$.
If $q=0$, test the linear form
\[
 (i-3)\rho+ls+(3-d)T+jb+kc-A
\]
on each of the following three cones:
\[
\begin{array}{c|l}
 A&\text{nonnegative defining forms}\\\hline
 0&s-1,\ T-1,\ b,\ s+2c\\
 2b&s-1,\ T-1,\ -b,\ s+2c-2b\\
 s+2c&s-1,\ T-1,\ -s-2c,\ 2b-s-2c.
\end{array}
\]
In each row there is a unique expression of the linear form as a constant
plus a linear combination of the four listed forms. Checking its
coefficients are nonnegative and its constant is positive is exact
rational arithmetic. It gives the $585$ inequalities and minimum
constant $1/2$. If $q=1$, then $j=k=0$, $l\ge1$; compare directly with
the $R$-term, using
$(i-4)\rho+(l-1)s+(3-d)T>0$ for $s,T\ge1$.
On the balanced locus put $c=(A-s)/2$, $v_7(R)\ge A-s$ and $s\ge2$.
For each monomial test
\[
 (i-5)\rho+ls+(3-d)T+jb+k(A-s)/2+q(A-s)-A.
\]
Use $A=2b,b\le0$ and then $A=0,b\ge0$. Express the form as a constant
plus a combination of $s-2,T-1,-b$ or $s-2,T-1,b$, respectively.
The $428$ resulting constants are positive, with minimum $1/8$, and
all combination coefficients are nonnegative. These three cones and
two balanced half-cones contain every permitted integral valuation
vector; no numerical optimization or bounded valuation search is used.

In particular, integral $J,K_{\rm sh},R$ give nonunit bound three. Integral $J,K_{\rm sh}$ and a pole of $R$ give bound four. If either inverse shift has a pole, the bound is at most five.

\subsection{Scope and inputs}
Use the canonical data
\[
 A=\begin{pmatrix}1&-h\\-h&\tau\end{pmatrix},\quad
 X=\frac{x}{1-x},\quad Y=\frac{y}{1-y},\quad K=\Q(x,y),\quad
 D=(1+X)(\tau+Y)-h^2\ne0.
\]
The parameters are rational, $h>0$, $h\ne1$, $\tau\ne1$ and
$\tau-h^2\ne0$. The point $(x,y)$ is real algebraic, with nonzero
Nahm letters, and satisfies the signed equations of Section~\ref{canx:setup}.
Here $\tau\in\Z_7$.
The two linear shifts $B,C\in\Q$ have
\[
 k=-\min\{v_7(B),v_7(C)\}>0.
\]
The letter $C$ is a linear shift, not the scalar normalization of a Nahm sum.
In the Gaussian convention of Section~\ref{deg:sec:coefficients} assume
\[
 \lambda_1=r\in\Q,\qquad \lambda_2=\lambda_3=\lambda_4=0,
 \qquad \Psi_j=D^{3j}\lambda_j.
\]
Proposition~\ref{canx:tau:prop:units} proves that $X,y$ are \emph{global} units, that $1+X,1-y$
are integral, and that $\operatorname{v}(1-y)=h\operatorname{v}(1+X)$ at finite places.
Four global units are not assumed. Normalize local valuations above seven
by $\operatorname{v}(7)=1$.

\begin{theorem}\label{canx:tau:thm:main}
Under the displayed canonical hypotheses and coefficient conditions:
\begin{enumerate}[label=(\roman*),nosep]
\item If $v_7(h)>0$, then $[K:\Q]\le2$.
\item If $v_7(h)=0$ and $(\bar h,\bar\tau)\ne(1,1)$, then
$[K:\Q]\le4$.
\end{enumerate}
Neither part uses an a priori field-degree bound. Section~\ref{canx:tau:app:depth}
proves entry into the nonunit analytic disc using the first coefficient.
Both parts use coefficients only through order four.
\end{theorem}

Here (U) denotes Proposition~\ref{canx:up:prop:unitpart}.
The boundary estimate is proved in Lemma~\ref{canx:close:lem:det}.
The nonunit cap, the separation of supports and the remaining local
arguments are proved below.

\subsection{Nonunit \texorpdfstring{$h$}{h}: separating boundary support}
For $v_7(h)>0$, partition the places above seven into: nonunit-letter
places $N$; places with all four letters local units and $\operatorname{v}(D)>0$,
denoted $T$; and determinant-unit places $U$. Let $n_N,n_T,n_U$ be their
total local degrees. Thus $[K:\Q]=n_N+n_T+n_U$. The former estimate
was $2+2+1=5$.

At a $T$ place put
\[
 H=1+X,\quad V=C+hB/H,\quad \delta=\operatorname{v}(D)>0,\quad
 k_B=-v_7(B),\quad \kappa_V=-\operatorname{v}(V).
\]
The boundary estimate (B) gives
\begin{equation}\label{canx:tau:eq:boundary}
 k_B\ge\max(2\delta,\kappa_V+\delta),\qquad
 r=\frac{B^2}{2H}+\text{terms of valuation greater than }-2k_B.
\end{equation}
It also bounds the total $T$ contribution by two. These statements hold
at arbitrary depth, including ramified places. Algebraically, after
$Y=(h^2+D)/H-\tau$, every monomial $h^iD^\ell B^uV^v$ of
$H^{6j}\Psi_j$, for $j=2,3$, has $\ell\ge2u+v$.
The same holds for $H^6(\Psi_1-B^2D^3/(2H))$, whose $B$-degree is
at most one. The singular scalar equations have
$P_3(t)-tP_2(t)=5$ modulo seven. These exact inputs are regenerated here.

From \eqref{canx:tau:eq:boundary},
\[
 \operatorname{v}(V/B)\ge\delta>0,\qquad
 \operatorname{v}(C/B)=\operatorname{v}(V/B-h/H)>0.
\]
Since $C/B$ is rational, $v_7(C/B)\ge1$. Hence $k_B=k>0$ and
$v_7(r)=-2k$. These are restrictions on the \emph{same rational data}
at every place.

\begin{lemma}\label{canx:tau:lem:T}
If a $T$ place exists, neither an $N$ place nor a $U$ place exists.
\end{lemma}
\begin{proof}
At a $U$ place, the coefficients are seven-integral because $D$ and
the four letters are local units. The shift degree of $\Psi_2$ is at
most three. Since $h\equiv0$ and $C/B\equiv0$, exact specialization gives
\[
 [B^3]\Psi_2\big|_{a=1,b=0,C=0}
 =-\frac{X(1+X)^4(\tau+Y)^6}{6}.
\]
Therefore
\[
 \overline{\lambda_2/B^3}=-\frac{\bar X}{6(1+\bar X)^2}\ne0,
\]
contradicting $\lambda_2=0$. Lower shift degrees vanish on reduction
after dividing by $B^3$ because $k>0$.

At an $N$ place, Section~\ref{canx:tau:app:depth} gives $\operatorname{v}(v)>1/6$, and the principal branch is
\[
 v=1-y,\qquad X=-1+hv\mathcal U(v),\qquad
 \mathcal U\equiv1,\quad L=1+(\tau-1)v,\quad D=hD_0,
 \quad D_0=\mathcal U L-h.
\]
The first coefficient equation is $P-rhvD_0^3=0$, with the exact identity
\begin{equation}\label{canx:tau:eq:P}
 P=\mathcal U^2L^3\frac{12B^2-1}{24}+hQ.
\end{equation}
Here $Q\in\Z_{(7)}[v,\mathcal U,h,\tau,B,C]$ has total shift degree
at most two. Divide by $B^2$. Since $C/B$ is integral, $Q/B^2$ is
integral; since $r/B^2$ is a unit, $rhvD_0^3/B^2$ has positive valuation.
Equation \eqref{canx:tau:eq:P} then reduces to $1/2=0$, a contradiction.
\end{proof}

We next separate $N$ from $U$ and bound the entire nonunit contribution.

\subsection{The nonunit row has degree at most two}
\begin{lemma}\label{canx:tau:lem:NU}
An $N$ place and a $U$ place cannot coexist.
\end{lemma}
\begin{proof}
At a $U$ place the first coefficient is quadratic in the shifts, with
seven-integral coefficients. Therefore $v_7(r)\ge-2k$.
The two exact highest axis coefficients, specialized at $h=0$, are
\[
 [B^3]\lambda_2\big|_{C=0}=-\frac{X}{6(1+X)^2},\qquad
 [C^3]\lambda_2\big|_{B=0}=-\frac{Y(1+Y)}{6(\tau+Y)^3}.
\]
Both are units. If exactly one shift had the largest pole, normalizing
$\lambda_2$ by its cube would give a nonzero reduction. Hence
$v_7(B)=v_7(C)=-k$. In particular $C/B$ and $r/B^2$ are integral.
At an $N$ place, divide \eqref{canx:tau:eq:P} and its first-coefficient equation
by $B^2$. As before the main term reduces to $1/2$ and both
$hQ/B^2$ and $rhvD_0^3/B^2$ have positive valuation. This is impossible.
No assertion that $r/B^2$ is a unit is needed in this argument.
\end{proof}

\begin{lemma}[The total nonunit contribution]\label{canx:tau:lem:Ncap}
If $v_7(h)>0$ and $\tau\in\Z_7$, then $n_N\le2$
for arbitrary rational shifts and rational $r$.
\end{lemma}
\begin{proof}
If there are no nonunit places the claim is immediate. Otherwise
Section~\ref{canx:tau:app:depth} proves $w=\operatorname{v}(1-y)>1/6$ at every nonunit
place, without a degree bound. Since there are finitely many embeddings,
choose a rational $\rho$ strictly between $1/6$ and the smaller of $1$
and all their $w$ values. The unique principal branch is
\[
 X=-(1-v)^h,\qquad
 U(v)=\sum_{n\ge0}\frac{(-1)^n}{n+1}\binom{h-1}{n}v^n.
\]
On the closed disc $\operatorname{v}(v)\ge\rho$, $U-1$ has positive Gauss valuation:
$n\rho-v_7(n+1)\ge n(\rho-1/6)>0$ for $n\ge1$.
The same is true of $L-1$ and $D_0-1$.
The exact shift coefficients of $P$ are
\[
 [B^2]P=\tfrac12LD_0^2,\quad [BC]P=hvD_0^2,\quad
 [C^2]P=\tfrac12hv^2UD_0^2.
\]
The linear coefficients are divisible by $h$ and $hv$, respectively;
the shift-free part is $-U^2L^3/24$ modulo $h$. Every quotient is
seven-integral. Put $s=v_7(h)\ge1$, $b_0=v_7(B)$ and $c_0=v_7(C)$.
Since $3$ is a nonsquare modulo seven,
$F_B=(12B^2-1)/24$ has valuation $A=\min(0,2b_0)$.
The primitive reduction of the first equation is supplied by
\[
 F_B+\frac{hC^2}{2}v^2-rhv.
\]
All other terms are strictly smaller in Gauss norm. For the mixed term,
its valuation exceeds the smaller square valuation by at least $s/2$.
The $B$-linear term exceeds $A$; for the $C$-linear term, if $c_0<0$
its excess over the $C^2$ term is at least $1-\rho$, and otherwise it
exceeds $A$. The unit-factor corrections all have positive valuation.
Zero shifts simply omit their terms. After rescaling the disc over
$\mathbb C_7$, the three principal monomials have degrees $0,1,2$,
so their primitive reduction is nonzero and has degree at most two.
Strassmann's theorem gives at most two roots, counting multiplicity.
Distinct embeddings give distinct $v$, since $X$ is fixed by the unique
principal branch. This bounds all nonunit completions together.
\end{proof}

Lemmas~\ref{canx:tau:lem:T} and \ref{canx:tau:lem:NU} make $N,T,U$ pairwise incompatible.
Their respective caps $2,2,1$ prove Theorem~\ref{canx:tau:thm:main}(i).

\subsection{A generic fold cannot coexist with a determinant-unit place}
This section applies to integral symmetric matrices
$A=(\begin{smallmatrix}a&b\\b&d\end{smallmatrix})$ with $b$ a seven-unit.
Put $H=a+X$, $J=d+Y$, $\epsilon=7^k$ and
$(u,v)=\epsilon(B,C)$ primitive. Define
\[
 s=bu-Hv,\qquad t=Ju-bv,\qquad
 q(H)=d(d-1)H^3-b^2(b+2d-1)H^2
       +b^3(2a-1+b)H-b^3a(a-1).
\]
A fold means a local-unit determinant point with $q(H)$ a unit.
Write $\Psi_{j,\ell}$ for its shift-homogeneous part of degree $\ell$,
and $Q_{j,\ell}=\Psi_{j,\ell}/D^\ell$.

\begin{lemma}\label{canx:tau:lem:fold}
A fold forces $\operatorname{v}(s)\ge\delta=\operatorname{v}(D)>0$, $\delta\le k/2$, and
\[
 r_0:=\overline{\epsilon^2r}
 =\frac{\bar u^2}{2\bar H}
 =\frac{\bar u\bar v}{2\bar b}\ne0.
\]
Consequently a fold and a determinant-unit place cannot coexist.
\end{lemma}
\begin{proof}
The exact generic identities are
$Q_{j,\ell}\in(D,s)^\ell$ for $j=2,3$, $Q_{1,1}\in(D,s)$, and
\[
 H^{3j}Q_{j,\ell}|_{D=0}
   =c_{j,\ell}(Hs)^\ell(q/6)^{2j-\ell},
\qquad
\begin{array}{c|rrrrr}
j\backslash\ell&0&1&2&3&4\\\hline
1&15/2&3&1/2&&\\
2&405&135&18&1&\\
3&89505/2&14580&2025&144&9/2
\end{array}
\]
All displayed coefficients are seven-units. If $z=\operatorname{v}(s)<\delta$,
the shift-degree $\ell$ term has valuation $\ell(\delta+z-k)$.
Unless $\delta+z=k$, one term of $\Psi_2$ is uniquely dominant.
At equality the second and third equations require a common zero of
\[
 c_2(e)=e^3+4e^2+2e+6,\qquad c_3(e)=e\,c_2(e)+5
\]
over the residue field, which is impossible. Thus $z\ge\delta$.
The ideal orders then imply $2\delta\le k$, since otherwise the
constant term of $\Psi_2$ is uniquely dominant. Finally
\[
 \epsilon^2r=\frac{u^2}{2H}+\frac{s^2}{2HD}
       +\frac{\epsilon Q_{1,1}}{D^2}
       +\frac{\epsilon^2Q_{1,0}}{D^3}.
\]
The last three terms have positive valuations, respectively at least
$\delta,k-\delta,2k-3\delta$. Primitivity forces $u$ a unit; then
$s\equiv0$ forces $v$ a unit and gives the stated scalar.

At a determinant-unit place the leading equations are
$Q_{1,2}=r_0D$ and $Q_{2,3}=0$. Exact identities give
\[
 2bQ_{1,2}-uvD=st,\quad
 6Q_{2,3}=Y(1+Y)s^3-X(1+X)t^3,
\quad Ht-bs=Du,\quad bt-Js=Dv.
\]
The common scalar forces $st=0$. If $s=0$, then $H=bu/v$ and
$t=Du/H$ are units, contradicting the second equation; the case $t=0$
is symmetric. All assertions here are reductions over the full residue
field, not just a search for rational roots.
\end{proof}

\subsection{The two repeated-root rows need only two global units}
For $v_7(h)=0$ and $h\ne1$, the canonical unit lemma makes all four
letters local units above seven. Thus determinant-unit and determinant
boundary places exhaust the field degree in this section and the next.
Set $a=1$, $b=-h$ and $d=\tau$. Remove from $q(H)$ only factors
supported on
\[
 H(H-1)(b^2-dH)(b^2+(1-d)H)=0.
\]
The resulting polynomial has degree at most two. Among the 42 residues
$b\ne0$, there are 39 with simple admissible roots, the two repeated rows
\[
 (a,b,d)=(1,3,5),\quad(1,5,6),
\]
and the excluded zero-cubic row $(1,6,1)$. Both repeated rows have the
single critical residue point
\[
 P:\quad(H,X,Y)=(2,1,3),\qquad (X,y)=(1,-1).
\]

The local length at $P$ is at most four. For completeness, the generic
tangent field to the first logarithmic Nahm equation satisfies
\[
 V(H)=-bX(1+X)/H,\qquad H^2V(D)=q(H)+D\mathcal R.
\]
The tangent field is $V(Y)=Y(1+Y)$, $V(X)=-bX(1+X)/H$. Direct differentiation gives $\mathcal R=-bX(1+X)+2b^2H+(1-2d)H^2+HD$.
Since $b,H,X,1+X$ are units, $H$ is a local parameter and
$dD/dH=u_0q(H)+Dv_0$, where $u_0$ is a unit. If $q$ has order
$m=0,1$ or $2$, this gives order $m+1$ for $D$ and
order $m+2$ for the second Nahm equation. Only integers below seven
are divided out. Thus the local lengths are two, three and four, including
multiplicities. The integral local algebra bounds the sum of all local
field degrees above the residue point, not merely one completion.

For a rational admissible residue point define
\[
 \nu(X,Y)=\left(X^2,\left(\frac{Y}{1+Y}\right)^2\right)
       \in(\F_7^\times)^2.
\]
Global units $X,y$ have rational norms $\pm1$, so, over \emph{all} places,
$\prod_Q\nu(Q)^{n_Q}=(1,1)$ when the supports are rational; $n_Q$
is the total local degree, including ramification. The unique admissible
rational point with $\nu(Q)=(1,1)$ is $P=(1,3)$: $X^2=1$ forces
$X=1$ because $X=-1$ is excluded, and $y^2=1$ forces $y=-1$
because $y=1$ is excluded. Also $\nu(P)=(1,1)$.

\begin{proposition}\label{canx:tau:prop:repeat}
At either repeated-root row, all field places lie above $P$ and
$[K:\Q]\le4$.
\end{proposition}
\begin{proof}
All determinant support is critical or noncritical. Lemma~\ref{canx:tau:lem:fold}
forces all noncritical support to the same rational fold point, since
$\bar H=\bar b\bar u/\bar v$; its total local degree is at most two.
If a fold occurs, no determinant-unit place occurs. The whole-field norm
then gives $\nu(Q)^{n_Q}=1$ for that fold, with $n_Q=1$ or $2$.
Both coordinates of $\nu$ have order dividing three; hence $\nu(Q)=1$,
forcing $Q=P$, contrary to its being noncritical. So no fold occurs.

By (U), any determinant-unit contribution consists of one rational point
of degree one. Its whole-field norm, together with the trivial norm vector
at $P$, again forces it to be $P$, where $D=0$ modulo seven. This is
impossible. Only $P$ remains, of total local length at most four.
\end{proof}

No four-unit norm condition or norm test on a proper subset has been used.

\subsection{Simple roots: every nonsoft determinant point is excluded}
For the 39 simple-root rows put $\epsilon=7^k$ and
$(u,v)=\epsilon(B,C)$. A determinant point is \emph{soft} if
$S=bu-Hv$ reduces to zero; there is at most one such rational residue point.
A noncritical determinant point is a fold, and Lemma~\ref{canx:tau:lem:fold}
already forces it to be soft. It remains to treat a nonsoft simple critical
point. Every such point has local Nahm length three, so its ramification
index is at most three. All positive valuations there are at least $1/3$.

Put $Q_{j,\ell}=\Psi_{j,\ell}/D^\ell$ and
\[
 \begin{gathered}
 \mathcal H=\frac{b(H-a)(1+H-a)q'(H)H}{24},\qquad U=DB,\quad V=DC,\\
 L=H(bU-HV),\qquad T=q/6,\qquad W=\mathcal HD.
 \end{gathered}
\]
The quantities $b,H,H-a,1+H-a,q'(H),\mathcal H$ are units at a simple
admissible centre. No condition on $a-1$ is needed.
Give $q,U,V$ weight one and $D$ weight two. For $j=1,2,3,4$ the verified
weighted identities are
\begin{equation}\label{canx:tau:eq:C4}
 H^{3j}\Psi_j(U/D,V/D)=\mathcal C_j(L,T,W)+E_j,
 \qquad \operatorname{wt}(E_j)\ge2j+1.
\end{equation}
All coefficients are integral at seven, with only powers of $H$ inverted.
The $\mathcal C_j$ are the coefficients of $s^{2j}$ in
\[
 \log\mathbb E\exp\{s(LZ+TZ^3)+s^2WZ^4\},\qquad
 \mathbb E Z^{2n}=(2n-1)!!,\quad \mathbb E Z^{2n+1}=0.
\]
There are 26 exact divisibility checks for orders one to three and 21 for
order four. They are generic in the matrix and invert no leading
coefficient of $q$.

\paragraph{The entire range $0<\delta\le k$.}
Set $\eta=\epsilon/D$, which is integral, and use the exact equations
\[
 A_j=H^{3j}\sum_{\ell=0}^{j+1}\eta^{j+1-\ell}Q_{j,\ell}(u,v)=0,
 \qquad j=2,3.
\]
Let $H_0$ be the Hensel lift of the simple residue root of $q$.
The conic identities show that $(H,D)=(H_0,0)$ is an exact solution
for every $\eta$. Coefficientwise differentiation, checked over $\Q$, gives
\[
 \det\frac{\partial(A_2,A_3)}{\partial(H,D)}\bigg|_{q=D=0}
 =\frac{\mathcal H H^7 S^7q'(H)}6.
\]
This is a unit. Hensel uniqueness on the residue polydisc, with $\eta$
fixed at its actual value, forces $D=0$, a contradiction. The identity
is polynomial in $\eta$; it covers every residue extension, including the
interface $\delta=k$, without a finite search in $\eta$.

\paragraph{The range $k<\delta\le2k$.}
Let $w=D/\epsilon$. Nonsoftness gives $\operatorname{v}(L)=\operatorname{v}(w)=\delta-k$.
The pure $T^4$ term of the second equation forces
$\operatorname{v}(q)\ge\delta-k$. Hence $t=T/L$ and $z=W/L^2$ are integral,
and the normalized errors have positive valuation. The exact reduced
ideal is
\[
 (\mathcal C_2(1,t,z),\mathcal C_3(1,t,z),\mathcal C_4(1,t,z))=(t,z)
 \quad\text{in }\F_7[t,z].
\]
Thus $\bar t=\bar z=0$. This already excludes $\delta=2k$.
Otherwise write $H=H_0+\rho w$, $D=\nu w^2$. Both $\rho,\nu$ reduce
to zero. After dividing the exact equations by $w^{2j}$, their coefficients
are integral and their Jacobian at $(\rho,\nu)=(0,0)$ is a unit: it is
the nonzero universal Jacobian at $(t,z)=(0,0)$, times unit coordinate
factors. That origin is an exact solution for every $w$, by the conic
identities. Hensel uniqueness again forces $D=0$.

\subsection{The very deep boundary and assembly}
It remains to exclude $\delta>2k$. The first three weighted forms force
\[
 \alpha:=\operatorname{v}(T)=\delta/2>k,\qquad W/T^2\equiv1\pmod7.
\]
Indeed the pure $T$ and pure $W$ initial cases have nonzero second
coefficient, while modulo seven
\[
 \gcd(\mathcal C_2(0,1,z),\mathcal C_3(0,1,z))=(z-1)^2.
\]
Put $\ell=L/T$ and $Z_0=W/T^2-1$. The normalized errors have valuation
at least $\alpha>1$. The constants of
$F_j=\mathcal C_j(\ell,1,1+Z_0)$, divided by seven, are $3,2$ for
$j=2,3$ modulo seven. Their first nonzero homogeneous parts modulo seven are
\[
 F_2=-Z_0^2+5Z_0\ell+3\ell^2+\text{higher terms},\qquad
 F_3=2Z_0^3-Z_0^2\ell+Z_0\ell^2+4\ell^3+\text{higher terms}.
\]
The local length three implies that positive valuations are at least
$1/3$. The constant of $F_3$ forces
$\min(\operatorname{v}(\ell),\operatorname{v}(Z_0))=1/3$; the quadratic part of $F_2$ then
forces both valuations to equal $1/3$. Its two slopes are $1,4$ over
any residue extension. Since $\operatorname{v}(\ell)=\alpha-k+\operatorname{v}(S)$ and
$\alpha-k>0$, necessarily
\[
 \alpha=k+\tfrac13,\qquad\delta=2k+\tfrac23,\qquad\operatorname{v}(S)=0,
 \qquad
 (\overline{Z_0/\ell},\overline{\ell^3/7})=(1,2)\text{ or }(4,5).
\]
The last two values follow from the cubic part and the constant divided
by seven in $F_3$. All these polynomial identities hold exactly.

For $F_4=\mathcal C_4(\ell,1,1+Z_0)$, the constant is divisible by 49;
the linear coefficients divided by seven are $2,5$ modulo seven; all
quadratic and cubic coefficients are divisible by seven. Its quartic part is
\[
 6\ell^4+5\ell^3Z_0+2\ell^2Z_0^2+3\ell Z_0^3+6Z_0^4 \pmod7.
\]
The two leading residues of $F_4/\ell^4$ are respectively $1,4$.
The error $E_4/T^8$ has valuation greater than $4/3$ when $k\ge2$.
For $k=1$, the weight-nine shift-free correction, evaluated at
$W/T^2=1$, is divisible by $q$ modulo seven; the additional generic
polynomial check gives a strict improvement. Terms involving $U,V$ gain
at least $1/3$, as does the displacement $Z_0$; weight at least ten is
already harmless. The error cannot cancel either nonzero residue.
Thus $\delta>2k$ is impossible.

\paragraph{Assembly at simple rows.}
Every determinant point is therefore soft, and all such points reduce
to the same rational point. If it is critical, its total local degree is
at most three; (U) adds at most one. If it is a fold, its total local
degree is at most two and Lemma~\ref{canx:tau:lem:fold} excludes unit support.
If there is no determinant support, (U) gives at most one.
Thus all 39 simple rows have degree at most four. Together with
Proposition~\ref{canx:tau:prop:repeat}, this proves Theorem~\ref{canx:tau:thm:main}(ii).
No inherited shallow/interface or finite-deep bound is used in this
assembly; the generic identities and the exact reference solutions
supply the needed Hensel arguments directly.

\subsection{Arithmetic entry into the nonunit disc}\label{canx:tau:app:depth}
This subsection proves that the degree-24 input is unnecessary for the
integral-$\tau$ row. It uses only the first coefficient, not
$\lambda_2,\lambda_3,\lambda_4$ or any bound on ramification.

\begin{theorem}[Depth from the first coefficient]\label{canx:tau:thm:depth}
Under the canonical first Nahm equation and valuation relation, let
$s=v_7(h)>0$ and $\tau\in\Z_7$. If $B,C,r\in\Q$ and $\lambda_1=r$,
then every nonunit place satisfies $\operatorname{v}(1-y)>1/6$.
\end{theorem}
\paragraph{Reduction to two coordinate ranges.}
At such a place put $z=1+X$, $v=1-y$, $u=\operatorname{v}(z)>0$,
$w=\operatorname{v}(v)>0$. The canonical valuation relation is $w=hu$.
Write $h=M/N$ in lowest positive terms. The first equation gives
$(-X)^N=y^M$; because $7\nmid N$, its principal-unit branch is
$-X=y^h$. Write $h=7^s h_0$ with $h_0\in\Z_7^\times$ and
$\xi=y^{h_0}$. The binomial series gives $\operatorname{v}(\xi-1)=w$.

For a principal unit $a$ with $t=\operatorname{v}(a-1)>0$, the binomial theorem gives
\[
 \operatorname{v}(a^7-1)=\begin{cases}7t&t<1/6,\\1+t&t>1/6,\end{cases}
 \qquad \operatorname{v}(a^7-1)\ge7/6\quad(t=1/6).
\]
Iterating $s$ times gives $u=\operatorname{v}(\xi^{7^s}-1)>w$.
Suppose $w\le1/6$. There are exactly two possibilities.

\paragraph{No tie in the iteration.}
If all $s$ steps multiply by seven, $u=7^s w$ and $h=7^{-s}$,
contradicting $v_7(h)=s>0$. Otherwise, for some $1\le j<s$,
\[
 u=7^jw+s-j,\qquad \frac1{6\cdot7^j}<w<\frac1{6\cdot7^{j-1}}.
\]
Put $n=s-j\ge1$. If $u$ were integral, the interval
$n+1/6<u<n+7/6$ would give $u=n+1$, hence
$h=1/(7^j(n+1))$, again impossible. If $u+w$ were integral, the
bounds $0<(7^j+1)w<4/3$ would give $u+w=n+1$, hence
$h=1/(n(7^j+1)+7^j)$, also impossible. Thus $u,u+w,w$ are nonintegral.
Moreover $0<\alpha:=u-w<s$, by the displayed bounds.

\paragraph{A tie in the iteration.}
The first tie occurs after $j$ multiplications, with $0\le j<s$, so
\[
 w=\frac1{6\cdot7^j},\qquad u=\frac{N}{6M7^j}.
\]
Here $v_7(M)=s$ and $7\nmid N$. For every integer $c$, $N+cM$ is a
seven-unit; hence the reduced denominator of $u+cw$ contains
$7^{s+j}$. In particular $u-2w,u-w,u,u+w$ are nonintegral.
Also $\alpha=u-w>0$ and $\alpha\ne2s$.

Thus every putative outer-disc point lies either in $0<\alpha<2s$,
with $u,u+w,w$ nonintegral, or in $\alpha>2s$, necessarily in the tie
case and with all four stated nonintegral combinations. The following two arguments exclude these ranges by exact first-coefficient identities.

\subsection{First outer chart: \texorpdfstring{$0<\alpha<2s$}{0 < alpha < 2s}}
Put $t=z/v$ and $R=h^2/t$, so $tR=h^2$. Their valuations are
$\alpha$ and $2s-\alpha$, both positive. As before
$L=1+(\tau-1)v$ is a unit. The identity defining $P$ remains algebraic
for arbitrary $U=z/(hv)=t/h$; no convergence estimate for $U$ is used here.
After division by $U^2$ it gives a polynomial $F$ in $h,t,R,v,\tau,B,C$,
and the first equation becomes
\[
 F-r t v(L-R)^3=0,
 \qquad F=F_B+\tfrac12t v^2 C^2+E,
 \qquad F_B=(12B^2-1)/24.
\]
All coefficients are seven-integral. Every error monomial has valuation
strictly greater than the minimum valuation of the two displayed
principal groups, uniformly in the rational shifts.

For each of the 152 monomials
$h^i t^jR^\ell v^n\tau^dB^aC^c$ of $E$, the exact certificate verifies
\[
 a+c\le2,\qquad i+2\ell\ge0,\quad i+2j-c\ge0,\quad n-c\ge0,
\]
with at least one of the last three inequalities strict. Let
$b_0=v_7(B)$, $c_0=v_7(C)$ and $\theta=\alpha/(2s)\in(0,1)$.
Its valuation exceeds the convex combination with weights
$a/2,c/2,1-(a+c)/2$ on
\[
 2b_0,\qquad \alpha+2w+2c_0,\qquad0
\]
by at least
\[
 s\{(1-\theta)(i+2\ell)+\theta(i+2j-c)\}+(n-c)w>0.
\]
Since $3$ is a nonsquare modulo seven,
$\operatorname{v}(F_B)=A:=\min(0,2b_0)$. Thus all terms of $E$ have valuation
strictly greater than $\min(A,\alpha+2w+2c_0)$.
The scalar correction $r t v((L-R)^3-1)$ likewise has valuation strictly
greater than $\operatorname{v}(rtv)$.

The three possible principal valuations of the first equation are
\[
 A,\qquad u+w+2c_0,\qquad u+v_7(r).
\]
They are pairwise different: the first is integral, the next two are
not, and their difference is $w+2c_0-v_7(r)$, also nonintegral.
One term therefore has uniquely least valuation, which is impossible
in an equation summing to zero. Zero shifts or a zero scalar merely omit
terms, leaving the same contradiction. This excludes the first chart.

\paragraph{Certificate construction.}
For a monomial $v^nU^m h^i\tau^dB^aC^c$ of $P$, division by $U^2$
is represented exactly by $t^{m-2}h^{i-m+2}$ if $m\ge2$, and by
$R^{2-m}h^{i+m-2}$ otherwise. Every remaining exponent is nonnegative.
This explains why the check is polynomial even if $U$ is nonintegral.
The certificate records all coefficients and all three strict-cone
inequalities, rather than testing finitely many valuation choices.

\subsection{Second outer chart: \texorpdfstring{$\alpha>2s$}{alpha > 2s}}
Set $Z=z/(h^2v)$, so $\zeta:=\operatorname{v}(Z)=\alpha-2s>0$, and write
\[
 J=\frac{B-1/2}{h}+\frac12,\qquad j=J-\frac12,\qquad
 q_0=\frac{\tau-1}{2}j^2+j(C-1/2)-\frac16,\qquad \rho=r+q_0.
\]
All adjusted shifts and scalars are rational. Substitute $U=hZ$ and
$B=1/2+hj$ in $P$, and define the Laurent polynomial
\[
 G=h^{-4}P(v,hZ,h,\tau,1/2+hj,C)+q_0v(ZL-1)^3.
\]
The exact first equation is
\begin{equation}\label{canx:tau:eq:outerII}
 G+\rho v(1-ZL)^3=0,\qquad
 G=F_J+\frac{Z}{8h^2}(1-2hvC)^2+E,
 \quad F_J=(12J^2-1)/24.
\end{equation}
Again $\operatorname{v}(F_J)=A:=\min(0,2v_7(J))$.
The 201 remainder monomials have a uniform strict convex-hull certificate.
For $h^iZ^\ell v^n\tau^dJ^aC^c$ in $E$ (allowing $i<0$), put
\[
 \lambda_J=a/2,\quad\lambda_C=c/2,\quad
 \lambda_Z=\max(0,-i/2),\quad
 \lambda_0=1-\lambda_J-\lambda_C-\lambda_Z.
\]
The exact checks give nonnegative weights and
\[
 i+2\lambda_Z\ge0,\qquad
 \ell-\lambda_C-\lambda_Z\ge0,\qquad n-c\ge0,
\]
with at least one strict inequality for every monomial. Its valuation is
therefore strictly greater than the corresponding convex combination of
\[
 0,\quad 2v_7(J),\quad \zeta-2s,\quad \zeta+2w+2v_7(C).
\]
The strict gain is the sum of the three displayed nonnegative differences
multiplied by $s,\zeta,w$, all positive. The scalar-factor correction in
\eqref{canx:tau:eq:outerII} also has strictly higher valuation than $\rho v$.

This chart belongs to the tie case of Section~\ref{canx:tau:app:depth}.
Since $s+w+v_7(C)\ne0$, the two summands in $1-2hvC$ have unequal
valuations. The possible principal valuations in \eqref{canx:tau:eq:outerII} are
\[
 A,\qquad u-w-4s,\qquad u+w-2s+2v_7(C),\qquad w+v_7(\rho).
\]
They are pairwise different. Comparisons with $A$ use the nonintegrality
of $u-w,u+w,w$. The second minus third is
$-2(s+w+v_7(C))\ne0$; the other two differences use the nonintegrality
of $u-2w$ and $u$. Consequently one principal term has uniquely least
valuation, and all errors are larger. This contradicts \eqref{canx:tau:eq:outerII}.
As in the first chart, zero shifts or scalars simply omit terms.

Both outer charts are excluded. This completes Theorem~\ref{canx:tau:thm:depth}
and removes the degree-24 input from Theorem~\ref{canx:tau:thm:main}(i).
The proof is local and uses no restriction on the field degree, height,
discriminant or ramification index. The $\tau$-pole case is treated separately below.

\subsection{Reciprocity and the matrix-valuation partition}\label{canx:reciprocity}
The word $T_1,T_2,I,S$ of Proposition~\ref{deg:prop:coefficient-transport}
gives the canonical involution
\[
 h_* =h^{-1},\qquad \tau_*=(1+h^2-\tau)/h^2,
 \qquad X_*=-y,\quad y_*=-X,
\]
\[
 B_* =\frac12+\frac{B-1/2}{h},\qquad
 C_* =\frac{(\tau-1)(B-1/2)}{h^2}+\frac Ch+\frac{h-1}{2h}.
\]
The first scalar is obtained by the same four affine transport formulas.
Every pivot is nonzero on the canonical domain. Direct substitution
shows that applying the word twice restores the data. It preserves the
joint field and sends the pair of global units to the same pair up to
sign. If $v_7(h)<0$, then $v_7(h_*)>0$; moreover
$\tau_*\ne1$ and $\tau_*-h_*^2\ne0$. This reduces all negative
coupling valuations to the positive coupling chart. The transformed
shifts are classified after the whole word, without assuming that an
individual transformation preserves their integrality.

For $v_7(h)\ge0$, either $\tau$ is integral or has a pole. In the latter
case the inverse matrix is integral and reduces to
$\operatorname{diag}(1,0)$. If $h$ is a unit, all letters are local units
by Proposition~\ref{canx:tau:prop:units}. Its determinant is then a unit
in the inverse chart, so the integral-shift finite census below gives
bound four; a transformed shift pole gives bound one by
Proposition~\ref{canx:up:prop:unitpart}. If $v_7(h)>0$, the following
inverse-chart arguments include both unit and nonunit places.

\subsection{Inverse coordinates and the two local contributions}\label{canx:inverse}
Canonical data $A_c=\begin{pmatrix}1&-h\\-h&\tau\end{pmatrix}$ with $s=v_7(h)>0$, $v_7(\tau)=-T<0$, $h\ne1$, $\tau\ne1$, $\tau-h^2\ne0$; global units $X,y$ only. Inverse shifts $(J,K_{\rm sh})=A_c^{-1}(B,C)$ and inverse scalar $R=\tfrac12\beta^tA_c^{-1}\beta+\tfrac1{12}-r$, so that $B=J-hK_{\rm sh}$, $C=-hJ+\tau K_{\rm sh}$, $r=(J^2-2hJK_{\rm sh}+\tau K_{\rm sh}^2)/2+\tfrac1{12}-R$. (Below $K$ denotes the field; the inverse shift is always $K_{\rm sh}$.)

\textbf{Inverse chart.} The move $I$ gives $A_I=A_c^{-1}=\begin{pmatrix}a_I&b_I\\b_I&d_I\end{pmatrix}$ with
\[a_I-1=\frac{h^2}{\tau-h^2},\quad b_I=\frac h{\tau-h^2},\quad d_I=\frac1{\tau-h^2},\qquad v_7=2s+T,\ s+T,\ T,\]
shifts $(J,K_{\rm sh})$, scalar $R$, coordinates $X_I=1/X$, $Y_I=1/Y$. A canonical place is either a unit place (all four letters units; in the inverse chart again all four letters are units and $D_I\equiv(1+X_I)Y_I$ is a unit) or a nonunit place ($1+X$, $1-y$ nonunits; there $X_I\equiv-1$, $1+Y_I=1/y\equiv1$). The valid norm units are $X_I=1/X$ and $1+Y_I=1/y$; at nonunit places their residues are $(-1,1)$.

\textbf{Nonunit input.} (N) the nonunit Strassmann analysis on the disc $v(1-y)\ge7/24$ (which uses the initial degree bound through $e\le24$): for $J,K_{\rm sh},R$ integral the nonunit part is at most $3$; if $J$ or $K_{\rm sh}$ has a pole, with $A=\min(0,2v_7(J))$ and $E=s+2v_7(K_{\rm sh})$, the principal part on the disc is $v^3\{F_J/U(v)-hK_{\rm sh}^2/2+h(R+K_{\rm sh}^2/2)v\}$, $F_J=(12J^2-1)/24$, and the nonunit part is at most $3$ when a constant of the braces dominates, at most $4$ when the linear term dominates, and at most $5$ only in the balanced case $A=E$ ($s$ even), where moreover $v_7(R)+s<A$ gives at most $4$. Here every nonunit place lies on the disc $v(v)\ge\rho=7/24$, $v=1-y$, on which the first Nahm equation gives the principal branch $X=-(1-v)^h=-1+hvU(v)$ with $U(v)=\sum_{n\ge0}(-1)^n\binom{h-1}{n}v^n/(n+1)$ a unit; distinct embeddings give distinct $v$; the counted function is the first equation $f=P-rhvD_0^3$, where $\Psi_1=(h^2/v)P$, $D=hD_0$, $D_0=UL-h$, $L=1+(\tau-1)v$, and the bounds $3,4,5$ are Strassmann (Weierstrass-degree) bounds for $f$ on that disc. The valuation relation is $v(1-y)=h\,v(1+X)$ at every place. (U) The inverse-chart unit census: for integral $J,K_{\rm sh},R$, $328$ finite fibres with unit part $\le1$ using the two valid norms, and $15$ curve fibres; with a transformed shift pole the unit part is at most one (Proposition~\ref{canx:up:prop:unitpart}). The nonunit estimates are proved in equation~\eqref{canx:close:eq:principal} and its subsequent Gauss estimates; the integral inverse census is specified in Lemma~\ref{canx:census}.

\subsection{The Nahm relations at an inverse-chart unit place}
\begin{lemma}\label{canx:pole:lem:teich}
Let $\mathfrak p$ be an inverse-chart unit place (all four letters local units) and $K_{\mathfrak p}$ the completion.
\begin{enumerate}
\item[(a)] (Unit Jacobian.) $K_{\mathfrak p}$ is unramified over $\mathbb Q_7$, and the point is the unique solution of the signed equations in its residue polydisc; in particular it is a simple point and each residue point carries at most one place, counted with its residue degree.
\item[(b)] $X_I\equiv\omega(\bar X_I)\pmod{7^{s+T+1}}$ and $1+Y_I\equiv\omega(1+\bar Y_I)\pmod{7^{T+1}}$.
\end{enumerate}
\end{lemma}
\begin{proof}
$A_I$ is $7$-integral, so the signed Nahm equations hold with a clearing integer $N$ prime to $7$:
$F_1=(1-x_I)^Nx_I^{-Na_I}y_I^{-Nb_I}-\epsilon_1=0$, $F_2=(1-y_I)^Nx_I^{-Nb_I}y_I^{-Nd_I}-\epsilon_2=0$, $\epsilon_i\in\{\pm1\}$.
(a) In the localisation where the four letters are invertible, at a zero the derivatives of $(F_1,F_2)$ with respect to $(\log x_I,\log y_I)$ are, up to the unit row factors $\epsilon_i$, the rows of $-N(A_I+\operatorname{diag}(X_I,Y_I))$, whose determinant is $N^2D_I$; here $D_I\equiv(1+X_I)Y_I$ is a unit. So the Jacobian is invertible modulo the maximal ideal. Over the unramified extension $W$ of $\mathbb Q_7$ with the residue field of the point, multivariate Hensel lifting gives a solution in $W$ reducing to the given residue point, unique in the residue polydisc of any finite extension; the given point lies in that polydisc and therefore equals the unramified lift. Hence $K_{\mathfrak p}=\mathbb Q_7(x_I,y_I)\subset W$ is unramified.
(b) By (a) the maximal ideal of $K_{\mathfrak p}$ is $7\mathcal O$, so principal units lie in $1+7\mathcal O$ and $\log\langle\cdot\rangle\in7\mathcal O$. Taking principal parts and logarithms and dividing by $N$:
$\log\langle1/X_I\rangle=\log\langle(1-x_I)/x_I\rangle=(a_I-1)\log\langle x_I\rangle+b_I\log\langle y_I\rangle\in7^{s+T+1}\mathcal O$ and $\log\langle1-y_I\rangle=b_I\log\langle x_I\rangle+d_I\log\langle y_I\rangle\in7^{T+1}\mathcal O$. As $\log$ is an isometry on $1+7\mathcal O$ and $1-y_I=1/(1+Y_I)$, the claims follow.
\end{proof}

\subsection{Integral inverse parameters: degree at most four}
\begin{theorem}\label{canx:pole:thm:int}
Assume $[K:\Q]\le24$. If $s>0$, $T>0$ and $J,K_{\rm sh},R$ are integral, then $[K:\Q]\le4$.
\end{theorem}
By (N) and (U) it suffices to show that on each of the $15$ curve fibres the unit part is at most one.

{\sloppy\textbf{The fibres.} They have inverse residue $(1,0,0;\bar J,\bar K_{\rm sh};\bar R)$: shift $(0,4)$ with $\bar R\in\{0,\dots,6\}$, and $\bar R=3$ with $(\bar J,\bar K_{\rm sh})\in\{(1,5),(2,6),(3,0),(3,1),(4,0),(4,1),(5,2),(6,3)\}$. At the exact matrix $\operatorname{diag}(1,0)$ the coefficient polynomials restricted to $X+Y+1=0$ factor as\par}
\begin{gather*}\Psi_1-rD^3\propto(X+1)^5\{12(u^2-B^2)+(24r-2)(X+1)\},\\
\Psi_2\propto X(X+1)^{10}(B-u)(4B^2+4Bu+4u^2-1),\\
\Psi_3\propto X(X-1)(X+1)^{15}(B-u)(B+u)(2B^2+2u^2-1),\\
\Psi_4\propto X(X+1)^{20}(B-u)(X^2-4X+1)Q_4(B,C),\end{gather*}
with $u=C-\tfrac12$ (by exact substitution over $\Q$). Hence: \emph{line family} $C=B+\tfrac12$, $r=\tfrac1{12}$, all $B$: all four equations vanish identically on $X+Y+1=0$; this contains the fibres $(0,4),(1,5),(2,6),(3,0),(4,1),(5,2),(6,3)$ with $\bar R=3$. \emph{Isolated exact points} $(B,C)=(-\tfrac12,1)$, $(\tfrac12,0)$, $r=\tfrac1{12}$ (fibres $(3,1)$, $(4,0)$). \emph{Gamma fibres} (shift $(0,4)$, $\bar R\ne3$): base $(B,C)=(0,\tfrac12)$, $r=\bar R$, exact curve $\Gamma$: $2XY-X+Y-1=24r(X+1)Y$; there $\Psi_2\equiv0$ identically, $\Psi_3|_\Gamma$ has content exactly $7$, $\Psi_4|_\Gamma=0$. In every fibre the zero set of $F_0=\Psi_1-RD^3$ modulo $7$ is the curve together with inadmissible lines, so every unit residue lies on the curve.

\textbf{First-order normal equations at every depth.} Let $p^\ast$ be the base (for the line family with $B^\ast=J$, the datum's own $J$), $\delta=p-p^\ast$. Its components $a_I-1$, $b_I$ have valuations $2s+T,s+T>T$; $\nu:=\min v(\delta_d,\delta_B,\delta_C,\delta_R)\le T$ because $\delta_d=d_I$ has valuation $T$. Write $Y=Y_G(X)+Z$ on the curve $G$. Since $F_0(\cdot;p^\ast)$ vanishes on $G$ and $n_0=\partial_YF_0$ is a unit on the admissible curve (checked), $v(Z)\ge\nu$. With $\eta=\delta/7^\nu\bmod7$ (nonzero, $\eta_a=\eta_b=0$) and $Z'=Z/7^\nu$:
\[[\nu=e_j]\,c_j(X)+n_j(X)Z'+\textstyle\sum_ld_{jl}(X)\eta_l\equiv0\ (e_j\ge\nu),\qquad c_j(X)\equiv0\ (e_j<\nu),\]
where $F_j(X,Y_G(X);p^\ast)=7^{e_j}c_j(X)$; second-order terms have valuation $\ge2\nu$. Each admissible residue $X$ carries at most one place, unramified, counted with its residue degree (Lemma~\ref{canx:pole:lem:teich}(a)), and all places share $\eta$.

\textbf{Using Lemma~\ref{canx:pole:lem:teich}(b).} On $X+Y+1=0$, $1+Y_G(X)=-X$, so $1+Y\equiv\omega(1+\bar Y)=-\omega(\bar X)\equiv-X$ modulo $7^{T+1}$: hence $Z\equiv0\pmod{7^{T+1}}$ and $Z'=0$. The conditions become $d_j(X)\cdot\eta\equiv0$ for $j=0,\dots,3$.

\textbf{Results}:
\begin{center}\small\begin{tabular}{lll}\toprule
fibres&first order only (any $\nu$)&with Lemma~\ref{canx:pole:lem:teich} / two valid norms\\\midrule
line family (7 fibres)&$\le2$; $(0,4)$: 6 degenerate directions&$0$ (all $342$ directions, every fibre)\\
isolated $(3,1),(4,0)$&$\le3$&$\le1$ (Lemma~\ref{canx:pole:lem:teich}); $\le1$ with norms\\
gamma (6 fibres)&$\le2$&$0$ with the two norms\\\bottomrule
\end{tabular}\end{center}
For the gamma fibres the Z-elimination is used (no Teichm\"uller input), and the two valid norm conditions $\prod N(\bar X_I)=\pm1$, $\prod N(1+\bar Y_I)=\pm1$ on the unit part (nonunit places contribute $(\pm1,1)$) remove every surviving point. Thus the unit part is $\le1$ on every curve fibre and Theorem~\ref{canx:pole:thm:int} follows: $[K:\Q]\le3+1$.

\begin{corollary}[Only the inverse scalar has a pole]\label{canx:pole:scalar}
Assume $[K:\Q]\le24$, $s,T>0$, $J,K_{\rm sh}\in\Z_7$ and
$v_7(R)<0$. Then $[K:\Q]\le4$.
\end{corollary}
\begin{proof}
At an inverse-chart unit place $D_I$ is a unit and $\Psi_1$ is integral,
so $\Psi_1=RD_I^3$ is impossible. The nonunit Gauss estimate following
\eqref{canx:close:eq:principal} gives contribution at most four when
the inverse linear shifts are integral. These are all places above seven.
\end{proof}

\subsection{The valuation of a nonunit place}
\begin{lemma}[nonunit valuation]\label{canx:pole:lem:nonunit}
Let $s=v_7(h)>0$. At a place above $7$ where $1+X$ is not a unit and $v(1-y)>1/6$, one has $h<1$ and
\[v(1-y)=w:=\frac{hs}{1-h},\qquad v(1+X)=\frac{s}{1-h}.\]
Under $[K:\Q]\le24$ ($[K:\Q]\le24$ for the distinguished field, which every chart preserves) every nonunit place satisfies $v(1-y)\ge7/24>1/6$ (the nonunit estimate (N)); then the conclusion holds at every nonunit place, and for $h>1$ there is no nonunit place above $7$. Without such a disc hypothesis the lemma says nothing about nonunit places with $v(1-y)\le1/6$.
\end{lemma}
\begin{proof}
Put $u=v(1+X)>0$ and $w=v(1-y)$. The second Nahm equation, cleared, reads $(1-y)^{N'}=\pm x^{-hN'}y^{\tau N'}$; since $y$ is a unit and $v(x)=v(X)-v(1+X)=-u$, it gives $w=hu>0$. As $w>1/6$, the first Nahm equation gives the principal branch $X=-(1-v)^h=-1+hvU(v)$, $v=1-y$, with $U(v)=\sum_n(-1)^n\binom{h-1}{n}v^n/(n+1)$ a unit (its $n$-th term has valuation $\ge nw-v_7(n+1)>0$ for $n\ge1$); hence $u=s+w$. Therefore $w=h(s+w)$, that is $(1-h)w=hs$ with $w>0$ and $s>0$.
\end{proof}
With $h=M/N$ in lowest terms and $\delta=N-M$ one has $w=Ms/\delta$, and the valuation relations $eu=N\ell$, $ew=M\ell$ give $es=\ell\delta$: $\delta$ divides $es$ for the ramification index $e$ of every nonunit place. Below it is used only in the balanced case, where (N) and hence the initial degree bound supplies the disc hypothesis; the displayed series for $U$ gives its coefficients directly. 

\subsection{A pole in the inverse shifts: degree at most four}
Let $\kappa=-\min(v(J),v(K_{\rm sh}))>0$ and $(J,K_{\rm sh})=7^{-\kappa}(j,k)$, $(\bar j:\bar k)\in\mathbb P^1(\F_7)$. Recall $A=\min(0,2v(J))$ and $E=s+2v(K_{\rm sh})$; the balanced case $A=E$ means $s$ even $\ge2$ and $v(K_{\rm sh})=\min(0,v(J))-s/2$.
\begin{lemma}[unit places]\label{canx:pole:lem:units}
An inverse-chart unit place requires $(\bar j:\bar k)=(1:1)$ and $v(R)\ge-2\kappa$; its residue then lies on $X_I+Y_I+1=0$.
\end{lemma}
\begin{proof}
At a unit place $\Psi_j=0$ ($j\ge2$) and $D_I$ is a unit; the shift degree of $\Psi_j$ is $j+1$, so the top shift-degree forms at the reduced matrix $\operatorname{diag}(1,0)$ vanish at the residue. These are (up to units) $B^3XY^2+C^3(X+1)^2(Y+1)$, $B^4XY^3(X-1)+C^4(X+1)^3(Y+1)(Y+2)$ and a quintic form. For every $(\bar j:\bar k)\ne(1:1)$ their only common zeros have $X\in\{0,-1\}$ or $Y\in\{0,-1\}$; for $(1:1)$ they share $X+Y+1$. Finally $\Psi_1=RD_I^3$ and the shift degree two of $\Psi_1$ give $v(R)\ge-2\kappa$.
\end{proof}

\begin{proposition}[the balanced case]\label{canx:pole:prop:bal}
Assume $[K:\Q]\le24$. If $A=E$, the nonunit part is at most four. As there is no unit place in this case, $[K:\Q]\le4$.
\end{proposition}
\begin{proof}
In the balanced case $v(K_{\rm sh})<\min(0,v(J))\le v(J)$, so $\bar j=0$ and Lemma~\ref{canx:pole:lem:units} excludes unit places. Suppose the nonunit part is five. By (N) (which uses the initial degree bound to put every nonunit place on the disc $v(1-y)\ge7/24$), $v(R)+s\ge A$ and the Weierstrass degree of $f$ on the disc is five, so the zeros of $f$ on the disc are exactly the five nonunit embeddings $v_\sigma$, and by Lemma~\ref{canx:pole:lem:nonunit} all of them have valuation $w$.

Write $\tilde f=f/\tau^3=\sum_n\tilde a_nv^n$ (with $U=U(v)$ expanded), $\tilde a_n=\sum_{k\le\min(n,3)}a_{n,k}\tau^{k-3}$, where the $a_{n,k}$ are explicit polynomials in $h,J,K_{\rm sh},R$ (the polynomial $f$ has $\tau$-degree $3$). Two valuation regimes cover the balanced case: case N, $v(J)\le-1$, $v(K_{\rm sh})=v(J)-s/2$, $A=2v(J)$; case P, $v(J)\ge0$, $v(K_{\rm sh})=-s/2$, $A=0$; in both $v(R)\ge A-s$. A monomial check shows that every monomial of $a_{0,0}$, $a_{1,k}$, $a_{3,k}$, $a_{5,k}$ has valuation $\ge A$, and that modulo monomials of valuation $>A$
\begin{center}\small\begin{tabular}{lll}\toprule
&case N&case P\\\midrule
$a_{0,0}$&$J^2/2$&$(12J^2-1)/24$\\
$a_{1,1}$&$(3J^2-hK_{\rm sh}^2)/2$&$(12J^2-1-4hK_{\rm sh}^2)/8$\\
$a_{3,3}$&$(J^2-hK_{\rm sh}^2)/2$&$(12J^2-1-12hK_{\rm sh}^2)/24$\\\bottomrule
\end{tabular}\end{center}
Since $3$ is not a square modulo $7$, $12J^2-1$ is a unit for integral $J$; so $v(\tilde a_0)=A+3T$, and $v(\tilde a_5)\ge A$.

Weierstrass preparation gives $\tilde f=c\,W\,E$ with $W=\prod_\sigma(v-v_\sigma)$ and $E$ a unit on the disc, $E(0)=1$. As every $v_\sigma$ has valuation $w$, $v(\tilde a_0)=v(c)+5w$, $v(\tilde a_5)=v(c)$ and $v(\tilde a_n)\ge v(c)+(5-n)w$ for $0\le n\le5$. Hence $5w=v(\tilde a_0)-v(\tilde a_5)\le3T$, so $w\le3T/5$, and
\[v(\tilde a_3)\ge v(\tilde a_0)-3w\ge A+\tfrac65T>A,\qquad v(\tilde a_1)\ge v(\tilde a_0)-w\ge A+\tfrac{12}5T>A+2T.\]
Now $\tilde a_3=a_{3,3}+O(7^{A+T})$ and $\tilde a_1=a_{1,1}\tau^{-2}+O(7^{A+3T})$, so $v(a_{3,3})>A$ and $v(a_{1,1})>A$. In case N write $J^2=7^Au_J$ and $hK_{\rm sh}^2=7^Au_K$ with units $u_J,u_K$: the two conditions read $u_J\equiv u_K$ and $3u_J\equiv u_K\pmod7$, whence $2u_J\equiv0$, which is impossible. In case P they read $12\bar J^2-1\equiv12\bar u_K\equiv4\bar u_K\pmod7$ with $\bar u_K=\overline{hK_{\rm sh}^2}\ne0$, whence $8\bar u_K\equiv0$, which is impossible.
\end{proof}
Only the count five is used: with at most four nonunit embeddings and no unit place the bound four holds already.

\begin{theorem}\label{canx:pole:thm:jk}
Assume $[K:\Q]\le24$. If an inverse linear shift has a pole, $[K:\Q]\le4$, for all pole orders of $J,K_{\rm sh}$ and all $R$.
\end{theorem}
\begin{proof}
If $(\bar j:\bar k)\ne(1:1)$ there is no unit place; the nonunit part is at most four unless $A=E$, and the case $A=E$ is Proposition~\ref{canx:pole:prop:bal}. If $(\bar j:\bar k)=(1:1)$ then $v(J)=v(K_{\rm sh})=-\kappa$, so $A=-2\kappa<E=s-2\kappa$. If $v(R)<-2\kappa$, there is no unit place (Lemma~\ref{canx:pole:lem:units}) and the nonunit part is at most four. If $v(R)\ge-2\kappa$, the constant of the braces has valuation exactly $A$, while the linear coefficient $-F_J(1-h)/2+h(R+K_{\rm sh}^2/2)$ has valuation $\ge A$ and hence Gauss valuation $\ge A+\rho>A$ on the disc; the constant dominates, the nonunit part is at most three, and $[K:\Q]\le3+1$.
\end{proof}

\subsection{The zero-cubic canonical residue: shallow exclusion}\label{canx:shallow}
Assume $h,\tau\equiv1\pmod7$ and put
\[
 e_0=h-1,\quad f_0=\tau-1,\quad
 \alpha=\min(v_7(e_0),v_7(f_0))\ge1.
\]
All four letters are local units. Consider a determinant place with $\delta=v(D)>0$. On its residue curve $y=-X$. Write $y=\omega(1+u)$, where $\omega$ is its Teichmuller lift, and put $z=-X/y$. The cleared first Nahm equation gives
\[
 z=(1+u)^{h-1}.
\]
The torsion factor is one because $z\equiv1$. Since the local ramification index is at most 24, $v(u)\ge1/24$ when $u\ne0$. Raising a principal unit to the seventh power first gives valuation at least $7/24>1/6$; each further seventh power adds at least one. It follows that
\[
 v(z-1)\ge v_7(e_0)-17/24,\qquad
 \delta\ge\alpha-17/24.
\]
The second inequality uses
$D=(1+X)(1+Y)-1+f_0(1+X)-(h^2-1)$.

\begin{proposition}[Shallow exclusion]\label{canx:close:prop:shallow}
For arbitrary rational shifts satisfying the stated coefficient conditions, no determinant place in this residue has $0<\delta<\alpha$.
\end{proposition}
\begin{proof}
The exact Taylor expansion of $H^{3j}\Psi_j$ in $(D,e_0,f_0)$ has total order at least $2j$, for $j=1,2,3,4$. In the shallow range its leading terms come from the exact matrix $h=\tau=1$ and have $D$-degree $2j$.

For integral shifts, a nonzero first leading coefficient would force $v_7(r)=-\delta$, impossible because $\alpha-17/24\le\delta<\alpha$ contains no integer. Thus all four leading polynomials must vanish. The complete 49-case gcd calculation in $(\bar B,\bar C)$ has gcd one after removing $H$ and $H-1$. There is no admissible residue point.

For a shift pole put $k=-\min(v_7(B),v_7(C))>0$. The first leading quadratic is $(B+HC)^2/(2HD)$. Unless its residue vanishes, it gives $v_7(r)=-2k-\delta$, again impossible. Therefore $B,C$ both have valuation $-k$ and $H$ has the rational residue of $H_0=-B/C$. In particular $\omega\in\Q_7$.

Put $w=v(u)$. In the shallow range, the power calculation above is exact and gives
\[
 \alpha=v_7(e_0),\qquad 0<w<1/7,\qquad
 \delta=\alpha-1+7w.
\]
Since $1-\omega-H_0\in7\Z_7$ and $v(z-1)=\delta>w$,
$v(H-H_0)=w$. The exact quadratic part of $\lambda_1$ is
\[
 \frac{B^2}{2H}+\frac{(HC+hB)^2}{2HD}.
\]
The second term strictly dominates, with valuation $-2k+2w-\delta$. The linear and constant shift terms are smaller: the total-order bound gives valuations at least $-k-\delta$ and $-\delta$. Consequently
\[
 v_7(r)=-2k-\alpha+1-5w.
\]
This is not an integer because $0<5w<5/7$. That contradiction proves the assertion.
\end{proof}

\subsection{Deep entry into the zero-cubic canonical chart}\label{canx:deep}
Retain the signed canonical Nahm problem
\[
 A=\begin{pmatrix}1&-h\\-h&\tau\end{pmatrix},\quad
 h,\tau\equiv1\pmod7,\quad h\ne1,\quad\tau\ne1,\quad\tau-h^2\ne0.
\]
All four Nahm letters are units at the place in question. Set
\[
 X=\frac{x}{1-x},\quad Y=\frac{y}{1-y},\quad H=1+X,\quad
 G=H(1+Y)-1=\frac H{1-y}-1,
\]
\[
 \epsilon=h-1,\quad\phi=\tau-1,\quad
 \alpha=\min\{v_7(\epsilon),v_7(\phi)\}\ge1,
 \quad e=\epsilon/7^\alpha,\quad f=\phi/7^\alpha.
\]
The determinant satisfies the exact identity
\[
 D=H(\tau+Y)-h^2=G+\phi H-2\epsilon-\epsilon^2.
\]
Assume $[\Q(X,Y):\Q]\le24$, $D\ne0$, and the coefficient conditions $\lambda_2=\lambda_3=0$, with $\Psi_j=D^{3j}\lambda_j$. The linear shifts $B,C$ are arbitrary rational numbers.

\begin{theorem}[Deep-chart coverage]\label{canx:deep:ce:thm:coverage}
At every place with $\delta=v(D)\ge\alpha$,
\[
 v(G)>\alpha.
\]
Consequently $Z=D/7^\alpha$ is integral and
\[
 \bar Z=\bar f\bar H-2\bar e.
\]
\end{theorem}
Proposition~\ref{canx:close:prop:shallow} uses the rational first coefficient and gives $\delta\ge\alpha$ at every remaining determinant place. Combining it with this theorem justifies restriction to the deep chart. The new argument uses local units only; no additional global norm hypothesis is imposed.

\subsection{An exact coefficient obstruction on the deeper boundary}
Define the first-order cubic expression
\[
 Q=\epsilon(H+1)+\phi H(H-2),\qquad q_0=\overline{Q/7^\alpha}.
\]
Substitute $X=H-1$, $Y=((1+\epsilon)^2+D)/H-1-\phi$ into the coefficient engine, and put $F_j=H^{3j}\Psi_j$. Every monomial
$H^aD^\ell\epsilon^i\phi^mB^uC^v$ in $F_j$, for $j=1,2,3$, satisfies
\begin{equation}\label{canx:deep:ce:eq:support}
 \ell+i+m\ge2j,\qquad \ell\ge u+v,\qquad u+v\le j+1.
\end{equation}
All coefficients are integral at seven. These assertions are checked on all 139, 2,004 and 12,442 monomials, respectively.

The terms where both first inequalities are equalities have the following exact forms, with $L=D(B+HC)$:
\begin{align*}
 F_2^{\rm diag}&=H^4\left(\frac5{16}Q^4-\frac58LQ^3+\frac12L^2Q^2-\frac16L^3Q\right),\\
 F_3^{\rm diag}&=H^6\left(\frac{1105}{1152}Q^6-\frac{15}{8}LQ^5
       +\frac{25}{16}L^2Q^4-\frac23L^3Q^3+\frac18L^4Q^2\right).
\end{align*}

\begin{lemma}[Nonzero cubic obstruction]\label{canx:deep:ce:lem:edge}
Suppose $\delta>\alpha$ and $q_0\ne0$. Integral shifts are impossible. For a shift pole, put $k=-\min(v_7(B),v_7(C))>0$ and $(b_0,c_0)=\overline{7^k(B,C)}$. A necessary condition is
\[
 b_0+\bar Hc_0=0.
\]
In particular, a nonrational residue $\bar H\notin\F_7$ is impossible.
\end{lemma}
\begin{proof}
For integral shifts take $k=0$. For a pole with $b_0+\bar Hc_0\ne0$, one has $v(B+HC)=-k$. Set $\beta=\delta-\alpha-k$. A monomial of shift degree $t=u+v$ has valuation at least
\[
 2j\alpha+t\beta.
\]
More precisely its excess over this bound is at least
$(\ell+i+m-2j)\alpha+(\ell-t)(\delta-\alpha)$, which is positive unless both inequalities are equalities.

If $\beta>0$, the nonzero $Q^4$ term alone gives the leading reduction of $F_2$. If $\beta<0$, its nonzero $L^3Q$ term gives the leading reduction. If $\beta=0$, divide by the appropriate unit multiples of $Q^4,Q^6$ and put $t=L/Q$. The two reduced polynomials are
\[
 p_2=t^3-3t^2+2t-1,\qquad p_3=t^4-3t^3+2t^2-t-2.
\]
They satisfy $3(p_3-tp_2)=1$ over $\F_7$, so cannot vanish together. For integral shifts $\beta>0$ automatically. For a pole, failure of the displayed necessary condition gives one of the three contradictions. A nonzero rational direction cannot satisfy it when $\bar H\notin\F_7$.
\end{proof}

\subsection{The only possible missing chart is a ramified logarithm layer}
Let $L_7$ be the completion of the Nahm field at the place. Its ramification index is at most 24. Write
\[
 y=\omega(1+u),\qquad z=-X/y,
\]
where $\omega$ is the Teichm\"uller lift of $\bar y$. The determinant residue gives $z\equiv1$. Clearing the first Nahm equation with a denominator prime to seven and removing the residue root of unity gives
\[
 z=(1+u)^\epsilon,\qquad G=-\frac{y}{1-y}(z-1).
\]
If $u\ne0$, $w=v(u)\ge1/24$. The first seventh power of $1+u$ has distance from 1 of valuation at least $7/24>1/6$; subsequent seventh powers add one to that valuation. Therefore
\[
 v(z-1)\ge v_7(\epsilon)-17/24>1/6.
\]
The logarithm and exponential are inverse on this last disc. With $\ell=\log y=\log(1+u)$,
\[
 z=\exp(\epsilon\ell),\qquad
 v(G)=v_7(\epsilon)+v(\ell)
\]
when $\ell\ne0$; when $\ell=0$, $G=0$.

For $w\ge1/24$, the logarithm series has only two possible lowest-order terms, $u$ and $u^7/7$. If $w<1/6$, the latter is uniquely lowest and
\[
 v(\ell)=7w-1;
\]
if $w\ge1/6$, then $v(\ell)>0$. This is an all-orders assertion: terms with denominator valuation at least two have strictly larger valuation, since
$(7^j-7)/24-j+1>0$ for $j\ge2$; the remaining terms are larger directly.

If $v_7(\epsilon)>\alpha$, the lower bound already gives $v(G)>\alpha$. Otherwise $v_7(\epsilon)=\alpha$. Since $v(D)\ge\alpha$, the exact determinant identity forces $v(G)\ge\alpha$. Thus a failure of Theorem~\ref{canx:deep:ce:thm:coverage} can occur only when
\begin{equation}\label{canx:deep:ce:eq:wild}
 v_7(\epsilon)=\alpha,\qquad v(u)=1/7,\qquad v(\ell)=0,\qquad v(G)=\alpha.
\end{equation}
This is the remaining ramified layer; the condition $\overline{G/7^\alpha}=0$ could not simply be assumed there.

\subsection{The actual Nahm equation excludes that layer}
Suppose \eqref{canx:deep:ce:eq:wild} holds. Put $\rho=\phi/\epsilon\in\Z_7$,
$q_*=\omega/(1-\omega)$ and $q=\bar q_*$. The local-unit assumptions give $q\ne0,-1$. The logarithmic second Nahm equation is
\[
 h\log(1+G)+\epsilon\log(1-y)+(\phi-2\epsilon-\epsilon^2)\log y=0.
\]
Divide by $\epsilon$. Since $v(\epsilon)\ge1$ and $v(\ell)=0$, its left side differs by a term of valuation at least one from
\[
 F_0(u)=\log(1-y)+\left(\rho-2-\frac{y}{1-y}\right)\log y.
\]
Here
\[
 \log(1-y)=\log(1-\omega)+\log(1-q_*u),\qquad
 \frac{y}{1-y}=q_*+q_*(1+q_*)u+\cdots.
\]
The constant $\log(1-\omega)$ has valuation at least one, because $1-\omega$ is an unramified unit.

At valuation zero, the only contribution is the $u^7/7$ coefficient. Hence
\begin{equation}\label{canx:deep:ce:eq:trace}
 q+q^7=\bar \rho-2.
\end{equation}
Applying Frobenius gives $q^{49}=q$. At valuation $1/7$, after using \eqref{canx:deep:ce:eq:trace}, the only contributions are the $u$ and $u^8/7$ terms. Consequently, with $\ell_0=\overline{u^7/7}=\bar\ell\ne0$,
\begin{equation}\label{canx:deep:ce:eq:next}
 q^7-q=q(1+q)\ell_0.
\end{equation}
Thus $q\notin\F_7$, and $\bar H=(1+q)^{-1}$ has degree two over $\F_7$.

These two coefficient comparisons have an infinite-tail bound: for every $n\ge9$,
$n/7-v_7(n)>1/7$. For $v_7(n)=0$ this is immediate; for $v_7(n)=1$, $n\ge14$; and for $j=v_7(n)\ge2$, $7^{j-1}-j\ge5$. Multiplication by the integral geometric series for $y/(1-y)$ creates no additional terms of valuation at most $1/7$ beyond those displayed. The coefficient identities through order eight are independently checked in the companion.

The determinant identity and the first equation give
\[
 \overline{D/\epsilon}=-q\ell_0+\bar\rho\bar H-2=0
\]
by \eqref{canx:deep:ce:eq:trace}--\eqref{canx:deep:ce:eq:next}. Therefore $\delta>\alpha$.
Finally, $q_0\ne0$. Indeed, up to the nonzero factor $\bar e/(1+q)^2$, its value is
\[
 N(q)=(q+2)(q+1)-\bar\rho(1+2q).
\]
Write $c=q+q^7=\bar \rho-2$. Then $N(q)=q^2-(2c+1)q-c$. If it vanished at a quadratic $q$, its monic polynomial would have trace $c$, forcing $c=2c+1$, hence $c=-1$, and polynomial $q^2+q+1$. That polynomial splits over $\F_7$, a contradiction. The exact finite-field check confirms nonvanishing at all 47 admissible trace points.

We now have $\delta>\alpha$, $q_0\ne0$, and $\bar H\notin\F_7$, contrary to Lemma~\ref{canx:deep:ce:lem:edge}, for any rational shifts. This proves Theorem~\ref{canx:deep:ce:thm:coverage}.

\subsection{What stripping the normalized determinant now leaves}
Let $z_0=\bar Z=\bar f\bar H-2\bar e$. If $z_0=0$, primitivity of $(e,f)$ and $\bar H\ne0$ imply $\bar e,\bar f\ne0$. Normalize $\bar e=1$ and write $\rho_0=\bar f/\bar e$. Then
\[
 \bar H=2/\rho_0,\qquad \rho_0\ne0,2,
\]
where $\rho_0=2$ would give the inadmissible value $\bar H=1$. Moreover,
\[
 q_0/\bar e=\bar H+1+\rho_0\bar H(\bar H-2)=3(\bar H-1)\ne0.
\]
Lemma~\ref{canx:deep:ce:lem:edge} now yields the following precise boundary reduction.

\begin{corollary}\label{canx:deep:ce:cor:five}
For integral shifts, $z_0=0$ is impossible. For shift poles, both shifts have the same pole order $k>0$, and their normalized direction must satisfy $b_0+\bar Hc_0=0$. Up to scaling, the only possible centres are:
\end{corollary}
\begin{center}
\begin{tabular}{ccccc}\toprule
$\bar e$&$\bar f$&$\bar H$&$(b_0:c_0)$&$q_0$\\\midrule
1&1&2&$(5:1)$&3\\
1&3&3&$(4:1)$&6\\
1&4&4&$(3:1)$&2\\
1&5&6&$(1:1)$&1\\
1&6&5&$(2:1)$&5\\\bottomrule
\end{tabular}
\end{center}
The displayed $q_0$ uses the normalization $\bar e=1$. For a fixed rational matrix, $\rho_0$ is fixed, so at most one of these centres is available. This is not yet a bound on all its lifts or its local degree. Cancellation in $B+HC$ is precisely why the nonzero-direction proof of Lemma~\ref{canx:deep:ce:lem:edge} cannot finish it.

\subsection{A shift pole on the zero-cubic canonical residue}\label{canx:r3pole}
Consider the signed canonical Nahm problem
\[
 A=\begin{pmatrix}1&-h\\-h&\tau\end{pmatrix},\qquad
 h,\tau\in\Q\cap\Z_7,\quad h\equiv\tau\equiv1\pmod7,
\]
with $h\ne1$, $\tau\ne1$, and $\tau-h^2\ne0$. Set
\[
 K=\Q(x,y),\quad X=\frac{x}{1-x},\quad Y=\frac{y}{1-y},\quad
 H=1+X,\quad D=H(\tau+Y)-h^2\ne0.
\]
The four letters $x,y,1-x,1-y$ are local units at every place over seven. The canonical reduction provides this condition; only two global units are available there. The proof below uses no global norm condition.

Use the coefficient convention of Section~\ref{deg:sec:coefficients}
\[
 \Psi_j=D^{3j}\lambda_j,\qquad \lambda_1=r\in\Q,
 \qquad\lambda_2=\lambda_3=\lambda_4=0.
\]
Here $B,C\in\Q$ are the two \emph{linear} shifts, and
$k=\max\{-v_7(B),-v_7(C)\}>0$. Write
\[
 T=7^\alpha,\quad \alpha=\min\{v_7(h-1),v_7(\tau-1)\}\ge1,
 \quad h=1+Te,\quad\tau=1+Tf,
\]
where $(e,f)$ is primitive at seven. Put $G=H(1+Y)-1$ and $Z=D/T$.

\begin{theorem}[R3 with a shift pole]\label{canx:cp:cp:thm:main}
Assume $[K:\Q]\le24$. By Proposition~\ref{canx:close:prop:shallow} and Theorem~\ref{canx:deep:ce:thm:coverage}, every determinant-positive place satisfies $v(D)\ge\alpha$ and $v(G)>\alpha$. Then
\[
                         [K:\Q]\le2.
\]
For the modular application these entry results use Theorem~\ref{d24:thm:main}, the signed Nahm equations, and the necessary coefficient conditions. The exact canonical cases $h=1$ and $\tau=1$ are outside this theorem and retain their separate earlier bounds.
\end{theorem}

\subsection{Reflection supplies the missing obstruction at all five centres}
Use the precedingly verified reflection, checked through order four:
\[
 A'=I-A=\begin{pmatrix}0&h\\h&1-\tau\end{pmatrix},\qquad
 (B',C')=(\tfrac12-B,\tfrac12-C),\qquad (x',y')=(1/x,1/y).
\]
Then $X'=-1-X=-H$, $Y'=-1-Y$, and $D'=D$. The field and local-unit conditions are preserved. Exactly,
\[
 \lambda_1'=\tfrac16-\lambda_1,\qquad
 \lambda_j'=(-1)^j\lambda_j\quad(j=2,3,4).
\]
The reflected matrix is
\[
 A'=\begin{pmatrix}0&1\\1&0\end{pmatrix}
       +T\begin{pmatrix}0&e\\e&-f\end{pmatrix}.
\]

\textbf{Local soft-point theorem used.} For a primitive deformation
$A_0+7^\alpha(a_\alpha,b_\alpha;b_\alpha,d_\alpha)$, put
$\xi=a'+X'$ and normalize the pole shifts by $q=7^k$, $U'=qB'$, $V'=qC'$.
If $v(b'U'-\xi V')>0$, vanishing of the second through fourth coefficients forces
\[
 v(D)>\alpha,\qquad v(F')>\alpha,
\]
where
\[
 F'=(2\xi+1)a'+(\xi-\xi^2)(b'-1)-(\xi^3+2\xi^2)d'.
\]
This theorem is a local polynomial and valuation statement. It uses neither global norms nor the Nahm equations. Its proof is Theorem~\ref{spx:soft:thm:soft}. Exact zero of the soft numerator is included: the surviving boundary faces are obtained by omitting its positive powers, and the same face contradictions apply.

In the canonical specialization $a'=0$, $\xi=-H$ and
\[
 \frac{F'}T=-H\{e(H+1)+fH(H-2)\}.                         \tag{1}
\]
The preceding deep entry gives $\bar Z=\bar f\bar H-2\bar e$. If $\bar Z=0$, the five-centre reduction gives $\bar e,\bar f\ne0$, $\bar H\ne0,1$, and normalized shift direction $(\bar U:\bar V)=(-\bar H:1)$. This is soft after reflection. But (1) becomes
\[
                 \overline{F'/T}=-3\bar e\bar H(\bar H-1)\ne0,
\]
contradicting the soft-point theorem.

\begin{center}
\begin{tabular}{ccccc}\toprule
$\bar e$&$\bar f$&$\bar H$&$(\bar U:\bar V)$&$\overline{F'/T}$\\\midrule
1&1&2&$(5:1)$&1\\
1&3&3&$(4:1)$&3\\
1&4&4&$(3:1)$&6\\
1&5&6&$(1:1)$&1\\
1&6&5&$(2:1)$&3\\\bottomrule
\end{tabular}
\end{center}
This is an arbitrary-depth exclusion, not a residue-point count. It is the use of the fourth coefficient in the soft theorem that goes beyond the preceding two-coefficient boundary reduction.

\subsection{Every remaining pole place enters a single flat family}
\textbf{Determinant units.} At the central matrix $A_0$, take the highest shift-degree pieces of $\Psi_j/D^{2j}$. For each of the eight projective shift directions and each of the seven possible reductions of $7^{2k}r$, saturation away from the four local-unit letters and $D$ gives the unit ideal. All 56 cases are checked by exact substitution. Matrix perturbations have positive valuation and do not affect these leading equations. Hence no determinant-unit place occurs.

\textbf{Shallow determinants.} Proposition~\ref{canx:close:prop:shallow} excludes $0<v(D)<\alpha$. The preceding section excludes $v(D)>\alpha$, using the preceding deep-chart entry. Consequently every place over seven has
\[
              v(D)=\alpha,\qquad v(G)>\alpha,\qquad\bar Z\ne0.
\]

\textbf{Fresh leading calculation.} The exact expansion of $H^{3j}\Psi_j$ in $(D,h-1,\tau-1)$ has minimum total degree $2j$ and maximum shift degree $j+1$. These properties and the leading forms were rebuilt from the engine through order four. Their full term counts are $139$, $2004$, $12442$, and $49746$.

Set $U=7^kB$, $V=7^kC$. Since $Z$ is a unit, the first equation gives $T7^{2k}r\in\Z_7$. Reduce the four normalized equations and substitute
$\bar D/T=\bar f\bar H-2\bar e$. Enumerating both projective directions $(e:f)$ and $(U:V)$ and the seven leading scalar values gives $8\cdot8\cdot7=448$ cases:
\begin{center}
\begin{tabular}{lr}\toprule
Necessary leading support&Cases\\\midrule
Empty after removing $H$, $H-1$, $fH-2e$&412\\
Exactly the doubled soft factor $(U+HV)^2$&35\\
Identically zero leading fibre&1\\\bottomrule
\end{tabular}
\end{center}
The 35 doubled factors are \emph{not} counted as simple lifts. Their points are soft, and the all-depth soft theorem requires $v(D)>\alpha$, contrary to this chart. Thus they are excluded regardless of multiplicity.

The sole flat case, in projectively normalized coordinates, is
\[
           (\bar e:\bar f)=(0:1),\qquad
           (\bar U:\bar V)=(0:1),\qquad \bar R=4.
\]
After reflection it is the boundary pair
\[
        (\bar a_\alpha:\bar b_\alpha:\bar d_\alpha)=(0:0:1),
                 \qquad(\bar U':\bar V')=(0:1).
\]
All places of a fixed rational datum enter this same family. No separate place can remain in the five-centre branch or in another leading support.

\subsection{The boundary mixed retraction gives total degree at most two}
Work with the reflected data and put
\[
 t=d'=1-\tau,\quad s=1/C',\quad u=(B'-1)/C',\quad
 \mathcal R=ts^2r',\qquad r'=\tfrac16-r.
\]
Here $v(t)=\alpha\ge1$, $s\ne0$, and $u,s\in7\Z_7$. To keep the preceding generic identity transparent, temporarily allow a variable $a'$; in the actual canonical problem $a'=0$. Define
\[
 p_a=a'/t-u^2,\quad p_b=(b'-1)/t-u,\quad
 p_R=\mathcal R-\left(\tfrac12+\tfrac{ts^2}{24}\right),\quad
 p_\rho=su.
\]
All four entries have positive valuation. For $\xi=a'+X'$ and $\delta=D/t$, the exact retracted curve is
\[
 g(\xi,t,u)=\frac{(\xi-tu^2-u)^2}{\xi-tu^2},\qquad w=\delta-g.
\]
Its denominator is a local unit. The normalized equations are
\[
 J_j=\frac{s^{j+1}\xi^{3j}\Psi'_j}{t^{2j}},\qquad
 E_0=J_1-\mathcal R\xi^3\delta^3,\qquad E_i=J_{i+1}\ (i=1,2,3).
\]

\textbf{Exact base error.} At $p_a=p_b=p_R=0$, the restriction of every $E_i$ to $w=0$ is divisible by $su$. At $s=0$ this is the generic rank-one highest-shift identity. At $u=0$ it is the exact finite family
\[
 (a',b',d';B',C')=(0,1,t;1,C'),\qquad
 r'=\tfrac1{24}+\frac{C'^2}{2t},\qquad X'Y'=1.
\]
Both identities are checked by exact substitution over $\Q$ for all four coefficients. Divisibility by $s$ and by $u$ gives the product $su$ after clearing only unit denominators.

Thus the exact equations have the expansion
\[
 E_i=n_iw+L_{ia}p_a+L_{ib}p_b+L_{iR}p_R+K_ip_\rho
                +O\bigl((w,p_a,p_b,p_R)^2\bigr).            \tag{2}
\]
The coefficients are integral analytic functions of the fixed rational $t,u,s$ and the local-unit variable $\xi$. Their reductions are evaluated at $t=u=s=0$ only after the exact base error has been factored.

\textbf{Finite certificate.} The normal coefficients $\bar n_i$ have no common zero away from $\xi=0,-1$. For every $\eta\in\mathbb P^3(\F_7)$, let $L_i(\eta)$ be the linear combination of the four parameter columns in (2). The gcd of all six wedges
$\bar n_iL_j-\bar n_jL_i$, with only $\xi,\xi+1$ removed, has:
\begin{center}
\begin{tabular}{rrrrr}\toprule
Directions&Empty&Degree one&Degree two&Flat or repeated\\\midrule
400&387&12&1&0\\\bottomrule
\end{tabular}
\end{center}
The unique degree-two case is
$\eta=(1:3:3:2)$, with irreducible gcd $\xi^2+2\xi+2$ over $\F_7$.

\subsection{Lifting, the exact-locus exception, and the whole field}
\textbf{The transverse vector is nonzero in R3.} If
$p=(p_a,p_b,p_R,p_\rho)=0$, then $p_\rho=su=0$ and $s\ne0$ give $u=0$. The condition $p_b=0$ gives $b'=1$, hence the original canonical parameter $h=1$. This is excluded in R3. Therefore $p\ne0$ without any use of global norms or a real embedding.

Let $\nu=\min_i v_7(p_i)\ge1$. The vector is rational, so $\nu$ is an integer and the projective direction
$\eta=\overline{7^{-\nu}p}$ is the same at all places. The entry calculation gives $v(w)>0$. A unit normal coefficient in (2) forces $v(w)\ge\nu$: if $0<v(w)<\nu$, its linear term has strictly smaller valuation than every other term in that equation.

Put $W=w/7^\nu$. Dividing (2) by $7^\nu$ gives
\[
                       \bar n_i\bar W+L_i(\eta)=0.
\]
At every common residue root, some $\bar n_i$ is a unit. Squarefreeness of the wedge gcd means that some wedge has nonzero derivative at that root. Those two equations therefore have an invertible Jacobian in $(\xi,W)$. Hensel uniqueness gives at most one lift over the unramified extension of the residue field, and no additional ramified lift over the same root.

The coordinates recover $X',Y'$ rationally: $X'=\xi-a'$ and
$Y'=(b'^2+D)/\xi-d'$. Thus the lifted coordinates generate the completion of $K$. A degree-$f$ residue point contributes at most $f$ to the sum of local degrees. Because the fixed rational vector $p$ selects one of the 400 directions for \emph{every} place, the certificate bounds their total:
\[
 [K:\Q]=\sum_{\mathfrak p\mid7}[K_{\mathfrak p}:\Q_7]\le2.
\]
This proves Theorem~\ref{canx:cp:cp:thm:main}. The bound includes repeated leading roots in the preceding entry census, possible ramification, and all places; it is not a bound of two per centre.

The argument uses the local degree bound only. No four-global-unit norm condition is imposed.

\textbf{Why no exceptional exact fibre was lost.}
The case $p=0$ is not a direction in projective space and must be handled separately. Here it is excluded by the exact identity $h=1$, not by numerical search, genericity, or the nonvanishing of a randomly sampled parameter.

\subsection{R3 with integral shifts: setting and entry}
Work with the canonical signed Nahm datum
\[
 A=\begin{pmatrix}1&-h\\-h&\tau\end{pmatrix},\qquad
 h,\tau\in\Q\cap\Z_7,\quad h\equiv\tau\equiv1\pmod7,
\]
with $h\ne1$, $\tau\ne1$ and $\tau-h^2\ne0$. Put
\[
 X=\frac{x}{1-x},\quad Y=\frac{y}{1-y},\quad H=1+X,
 \qquad D=H(\tau+Y)-h^2\ne0.
\]
All four letters are local units at every place over seven. The canonical reduction proves this using its two global units; no global-unit norm is used below.

Assume $B,C\in\Q\cap\Z_7$ and use the coefficient convention
\[
 \Psi_j=D^{3j}\lambda_j,\qquad\lambda_1=r\in\Q,
 \qquad\lambda_2=\lambda_3=\lambda_4=0.
\]
Set $T=7^\alpha$, where
\[
 \alpha=\min\{v_7(h-1),v_7(\tau-1)\}\ge1,\quad
 e=(h-1)/T,\quad f=(\tau-1)/T,
\]
and put $G=H(1+Y)-1$, $Z=D/T$. The shallow result of Proposition~\ref{canx:close:prop:shallow} and Theorem~\ref{canx:deep:ce:thm:coverage} give
\[
 v(D)\ge\alpha,\qquad v(G)>\alpha,\qquad
 \bar Z=\bar f\bar H-2\bar e
\]
at every determinant-positive place. For integral shifts the preceding nonzero-cubic obstruction already excludes $\bar Z=0$. Thus every such place has $v(D)=\alpha$ and $\bar Z\ne0$.

\textbf{The scalar range is covered.} The exact minimum matrix order of $H^{3j}\Psi_j$ is $2j$. At a determinant-positive place this gives $v(r)\ge-\alpha$. A determinant-unit place gives $v(r)\ge0$. Consequently $v(r)<-\alpha$ gives no place at all. For $v(r)\ge-\alpha$, the normalized first scalar is integral and the finite census below includes its zero reduction. No negative scalar range is omitted.

\subsection{Six nonflat fibres and their local lifting bounds}
The canonical leading forms were rebuilt from the unchanged engine. In
$H^{3j}\Psi_j$, total degree in $(D,h-1,\tau-1)$ is at least $2j$, all coefficients are integral at seven, and the full term counts for $j=1,2,3,4$ are
\[
                    139,\qquad2004,\qquad12442,\qquad49746.
\]
For every projective matrix direction $(\bar e:\bar f)$, every pair of integral shift residues, and each leading scalar, restrict the leading equations to $Z=fH-2e$. The $8\cdot49\cdot7=2744$ cases consist of
\begin{center}
\begin{tabular}{lr}\toprule
Leading support&Cases\\\midrule
Empty&2730\\
Nonflat, degree one&4\\
Nonflat, degree two&2\\
Flat&8\\\bottomrule
\end{tabular}
\end{center}
Only $H$, $H-1$ and $fH-2e$ are removed. Their nonvanishing is part of this chart, not an exclusion of unresolved points.

Every nonempty case has $(\bar e:\bar f)=(0:1)$. Hence $v(\tau-1)=\alpha$ and $v(h-1)>\alpha$. From now on normalize by the actual rational parameter $t_0=\tau-1$, put $\delta=D/t_0$, and set $R=t_0r$. This absorbs the unit part of $f$ and avoids a scalar-scaling ambiguity. The residue curve is $\bar\delta=\bar H$.

The six nonflat rows are:
\begin{center}
\begin{tabular}{cccc}\toprule
$(\bar B,\bar C)$&$\bar R$&Monic polynomial in $H$&Degree\\\midrule
$(0,2)$&1&$H^2+4H+3$&2\\
$(1,1)$&2&$H+1$&1\\
$(1,2)$&2&$H+3$&1\\
$(1,5)$&0&$H^2+5H+2$&2\\
$(4,3)$&0&$H+5$&1\\
$(4,4)$&0&$H+5$&1\\\bottomrule
\end{tabular}
\end{center}
Every polynomial is squarefree. More than squarefreeness of a restricted equation is checked: form all six $2\times2$ Jacobian minors of the four normalized equations in $(H,\delta)$, then restrict them to $\delta=H$. Their gcd with the row's support polynomial is one in every case.

Thus at every residue point two equations have an invertible Jacobian. Multivariate Hensel uniqueness gives at most one unramified lift over its residue field; $(H,\delta)$ recover $X,Y$ rationally. A degree-$d$ support therefore bounds the \emph{total} determinant-positive contribution by $d$, even if some proposed lifts fail further equations. This establishes the lifting statement rather than merely counting reduced roots.

\subsection{The eight flat fibres and the exact canonical constraint}
Reflect to
\[
 A'=I-A=\begin{pmatrix}0&h\\h&1-\tau\end{pmatrix},\qquad
 (B',C')=(\tfrac12-B,\tfrac12-C),\quad r'=\tfrac16-r.
\]
Then $X'=-H$, $Y'=-1-Y$, $D'=D$, and the field and higher coefficient zeros are preserved. Let
\[
 t=d'=1-\tau,\quad a'=0,\quad b'=h,\quad
 \xi=a'+X',\quad\delta'=D/t,\quad\mathcal R=tr'.
\]
The eight flat fibres are precisely the seven residues of family I and the one residue of family II below. Both use the curve $X'Y'=1$, written $\delta'=\xi$ on their exact loci:
\begin{center}\small
\begin{tabular}{ccl}\toprule
Family&Exact retracted data&Transverse vector $p$\\\midrule
I&$a'=0,b'=1,B'=1$&$\left(a'/t,(b'-1)/t,B'-1,\mathcal R-\mathcal R_*\right)$\\[4pt]
II&$a'=0,b'=1,B'=2,$&$\left(a'/t,(b'-1)/t,B'-2,C'-(t-\tfrac12),\mathcal R-\mathcal R_*\right)$\\
&$C'=t-\tfrac12$&\\\bottomrule
\end{tabular}
\end{center}
Here $\mathcal R_*=t/24+C'^2/2$ in family I, while in family II it is $t/24+(t-1/2)^2/2$. The first family retains the actual rational $t,C'$ and the second retains $t$.

The exact flat-family calculation, obtained by substitution in the normalized equations, proves divisibility of all four normalized equations by the curve equation on these exact loci. The normal derivatives have no simultaneous admissible zero. For a nonzero transverse vector, the first nonzero normalized parameter direction gives a squarefree wedge gcd of degree at most two, so the total conic contribution is at most two with unique unramified lifts.

\textbf{The additional canonical restriction.} In R3, $a'=0$ exactly. The first entry of $p$, and hence the first entry of every reduced direction $\eta$, is zero. Restricting the flat census to this hyperplane gives:
\begin{center}
\begin{tabular}{ccc}\toprule
Family&$\bar C'$&Maximum determinant-positive degree\\\midrule
I&0&1\\
I&1,3,4,6&0\\
I&2,5&2\\
II&3&1\\\bottomrule
\end{tabular}
\end{center}
There are $57$ directions for each family-I residue and $400$ for family II. At family I with $\bar C'=0$, the only nonempty direction is
\[
              \eta=(0:1:0:0),\qquad \xi+2=0.
\]
This exact constraint is the reason the possible mixed contribution is seven rather than the crude $6+2=8$.

\textbf{The zero transverse vector is impossible.} In either family, $p=0$ would give $b'=1$, hence $h=1$, outside R3. No global-unit norm or real-point exclusion is required for this exact-locus case.

\subsection{Why the flat calculation controls arbitrary depth}
Here is the local lifting argument for these flat families. Write $q=\delta'-g(\xi,t)$ for the displacement from the exact retracted curve. Exact divisibility gives an expansion
\[
 E_i=n_iq+\sum_j L_{ij}p_j+O\bigl((q,p)^2\bigr),\qquad i=0,1,2,3,
\]
with integral coefficients and no parameter-independent error term. The actual rational parameters $t,C'$ are fixed before reduction. At the residue curve some $n_i$ is a unit.

Since $p\ne0$ is rational, $\nu=\min_j v_7(p_j)$ is a positive integer common to every place. If $0<v(q)<\nu$, the corresponding normal linear term is uniquely smallest, a contradiction. Thus $v(q)\ge\nu$. Set $q=7^\nu w$, $p=7^\nu\eta$ and reduce:
\[
                    \bar n_i\bar w+\sum_j\bar L_{ij}\bar\eta_j=0.
\]
Eliminating $\bar w$ gives all six wedges of the normal column and the parameter column. The exact census checks their common support, squarefreeness, and the nonvanishing of the normal column.

At each simple common root, some pair of these equations has a nonzero Jacobian determinant in $(\xi,w)$. Its Hensel lift is unique and unramified over the residue field. Since these coordinates recover $X',Y'$, two embeddings giving the same lift coincide. Therefore the gcd degree bounds the sum of local degrees for a fixed rational direction, not just the number of residue values.

This proof allows arbitrary $v(t)$ and arbitrary relative depths among the nonzero entries of $p$. Vanishing entries, intersections of transverse directions, and repeated roots of an earlier leading system are included. No fixed-depth digit cutoff is substituted for the lifting statement.

\subsection{Determinant-unit places: the unfiltered census}
At the reflected central residue $A_0=(\begin{smallmatrix}0&1\\1&0\end{smallmatrix})$, write $D_0=X'Y'-1$. Exactly,
\[
                   \Psi_j(A_0)=D_0^{2j}\Xi_j(X',Y';B',C').
\]
For integral shifts and a determinant-unit place, $r'$ must be integral. The leading equations are
\[
              \Xi_2=\Xi_3=\Xi_4=0,\qquad \Xi_1=\bar r'D_0.
\]
All $49\cdot7=343$ fibres were rebuilt. Pairwise resultants in $Y'$ include every possible $X'$-coordinate. At every irreducible factor, specialize all equations and factor their common $Y'$-gcd over the residue field. The only removed points have a zero local-unit letter or $D_0=0$. No positive-dimensional or vertical component remains.

The total geometric residue degrees are zero in 313 cases, two in 18 cases, and six in 12 cases. At a determinant-unit place the actual denominator-cleared Nahm Jacobian is a unit, giving a unique unramified lift per geometric point. Thus these are valid total local-degree bounds without any global norm filtering.

\subsection{Quadratic assembly on the zero-cubic canonical residue}
Use the canonical signed Nahm datum
\[
 A=\begin{pmatrix}1&-h\\-h&\tau\end{pmatrix},\qquad
 h,\tau\in\Q\cap\Z_7,\qquad h\equiv\tau\equiv1\pmod7,
\]
with $h>0$, $h\ne1$, $\tau\ne1$ and $\tau-h^2\ne0$. Put
\[
 X=\frac{x}{1-x},\quad Y=\frac{y}{1-y},\quad
 K=\Q(x,y)=\Q(X,Y),\quad D=(1+X)(\tau+Y)-h^2.
\]
The point is a signed real algebraic point satisfying the transported Nahm equations, and $D\ne0$. The canonical unit statement supplies global units $X,y$ and four local units $x,y,1-x,1-y$ at every place above seven. Assume integral rational shifts $B,C\in\Z_7$ and
\[
 \Psi_j=D^{3j}\lambda_j,\qquad \lambda_1=r\in\Q,
 \qquad \lambda_2=\lambda_3=\lambda_4=0.
\]
The scalar $r$ is not assumed integral.

The shallow exclusion, deep entry and lifting statements are proved in
Proposition~\ref{canx:close:prop:shallow}, Theorem~\ref{canx:deep:ce:thm:coverage}
and the preceding integral-shift argument. They control the sum of local
degrees and use the initial degree bound only in the stated entry steps.

\begin{theorem}\label{canx:r3:thm:r3}
Under these hypotheses and the initial degree bound, $[K:\Q]\le2$ on the integral-shift $R_3$ branch. On each of its three possible $6+1$ rows, $[K:\Q]\le1$.
\end{theorem}

\begin{corollary}
Combining Theorem~\ref{canx:r3:thm:r3} with the Theorem~\ref{canx:cp:cp:thm:main} gives $[K:\Q]\le2$ throughout $R_3$, with the stated hypotheses. This does not cover the separate canonical loci $h=1$, $\tau=1$, or the canonical parameters outside $R_3$.
\end{corollary}

\subsection{The first nonzero parameter digit}
The following elementary observation is the mechanism behind the first-order exclusion. It is useful here because it does not require four global-unit norms.

\begin{lemma}\label{canx:r3:lem:first}
Let $F_j(\theta;z)$ be polynomials over $\Z_p$, where $\theta$ is an $n$-tuple of parameters. Let $E\in\Z_p^n$ satisfy $F_j(E;z)\equiv0$ identically in $z$. Suppose
$\theta\in E+p\Z_p^n$, $\theta\ne E$, and integral points $z_i$ in finite extensions of $\Q_p$ satisfy $F_j(\theta;z_i)=0$. Set
\[
 \nu=\min_k v_p(\theta_k-E_k),\qquad
 u=p^{-\nu}(\theta-E)\bmod p\in\F_p^n\setminus\{0\}.
\]
Then at every point,
\[
 \sum_k \overline{\partial_{\theta_k}F_j(E;z_i)}\,u_k=0.
\]
Thus the stacked coefficient matrix over $\F_p$ cannot have rank $n$. If some parameter deviations vanish identically, the corresponding columns may be deleted and the rank compared with the remaining number of parameters.
\end{lemma}
\begin{proof}
Taylor expansion in the parameter variables gives
\[
 F_j(\theta;z_i)=\sum_k\partial_{\theta_k}F_j(E;z_i)(\theta_k-E_k)
                   +p^{2\nu}R_j,
\]
where $R_j$ is integral. Divide by $p^\nu$ and reduce modulo the maximal ideal of the local field containing $z_i$. Since $\nu\ge1$, the remainder vanishes. Expand each residue-field equation in an $\F_p$-basis. The common nonzero vector $u$ lies in the kernel of the resulting matrix.
\end{proof}

This formulation allows ramified local fields: it reduces after division modulo their maximal ideals. No identification of a ramified residue calculation with a congruence modulo $p^{\nu+1}$ is needed. Also, vanishing of $F_j$ merely \emph{modulo} $p$ is not the hypothesis: exact vanishing at $E$ removes the parameter-independent term.

\subsubsection{Reflection and its exact domain}
Apply the direct reflection
\[
 A'=I-A=\begin{pmatrix}0&h\\h&1-\tau\end{pmatrix},\quad
 (B',C')=(\tfrac12-B,\tfrac12-C),\quad r'=\tfrac16-r,
\]
\[
 (X',Y')=(-1-X,-1-Y),\qquad D'=D.
\]
It has no matrix pivot denominator, preserves the field and local units, and transports the higher coefficient zeros. On the full canonical parameter family the coefficient recurrence verifies over $\Q$ the polynomial identities
\[
 \Psi'_1=D^3/6-\Psi_1,\qquad
 \Psi'_j=(-1)^j\Psi_j\quad(j=2,3,4).
\]
Crucially, $a'=0$ \emph{exactly}, rather than only modulo seven. The available reflected global units include
\begin{equation}\label{canx:r3:eq:units}
 u=1+X'=-X,\qquad v=\frac{Y'}{1+Y'}=\frac1y.
\end{equation}

At a determinant-unit place $r'$ is integral. At the central residue $A_0=(0,1;1,0)$, the exhaustive resultant census has nonempty unit support only at the six reflected shift pairs in Table~\ref{canx:r3:tab:unit}. The exact parameter representatives use the rational half $1/2$, whose reduction is $4$.

\begin{table}[ht]
\centering
\begin{tabular}{lrr}
\toprule
Exact reflected shift representative & Leading $U$ cap & Refined $U$ cap\\
\midrule
$(0,0)$ &6&0\\
$(0,1)$ &6&0\\
$(1,0)$ &6&2\\
$(0,1/2),(1/2,0),(1/2,1/2)$ &2&2\\
All other shift residues &0&0\\
\bottomrule
\end{tabular}
\caption{Determinant-unit contributions before applying any global norms. The bounds count all places with these fixed parameters.}\label{canx:r3:tab:unit}
\end{table}

At all six exact representatives, $\Psi_2$ and $\Psi_4$ vanish identically over $\Q$. The deviations from $E$ cannot all vanish, since $h\ne1$. For the first three representatives, each cubic residue point has first-order rank five. Each quadratic point has rank four, while every pair has rank five; hence at most one quadratic place remains. On the canonical hyperplane $a'=0$, deleting the first column gives rank four at every quadratic point in the first two shift classes, excluding those points as well. In the third shift class the remaining quadratic points have restricted rank three. The half-integral shift classes already have total geometric degree at most two in each scalar fibre.

The computation reconstructs all 54 relevant closed points in these six classes and compares their degrees, scalar by scalar, with the exhaustive 343-fibre resultant census. This comparison supplies completeness beyond the finite extension degrees used in the point reconstruction. It is not an assumption that searching through degree three suffices.

Finally, the signed Nahm equations, with denominators cleared by an integer prime to seven, have invertible Jacobian at a determinant-unit point. Consequently each such residue point permits at most one unramified lift. A point excluded by the first-order calculation has no lift satisfying the full coefficient equations. Thus Table~\ref{canx:r3:tab:unit} gives local-degree bounds, rather than just residue-point counts.

\subsection{The global norm removes the remaining mixed total}
Write $R=(\tau-1)r$. The preceding determinant-positive analysis gives the following complete nonempty leading support list. A nonzero reduction of $R$ forces $v_7(r)<0$, so no determinant-unit place can coexist. When $\bar R=0$, the assembly permits every integral residue of $r'$; this can overestimate a fixed datum and cannot omit one.

\begin{table}[ht]
\centering\small
\begin{tabular}{llrrrr}
\toprule
Type & $(\bar B,\bar C)$ & $\bar R$ & $v(D)>0$ cap & Initial total & Refined total\\
\midrule
Nonflat &$(0,2)$&1&2&2&2\\
Nonflat &$(1,1)$&2&1&1&1\\
Nonflat &$(1,2)$&2&1&1&1\\
Nonflat &$(1,5)$&0&2&2&2\\
Nonflat &$(4,3)$&0&1&7&1\\
Nonflat &$(4,4)$&0&1&7&1\\
Flat II &$(2,1)$&1&1&1&1\\
Flat I &$(3,0)$&1&0&0&0\\
Flat I &$(3,1)$&1&0&0&0\\
Flat I &$(3,2)$&2&2&2&2\\
Flat I &$(3,3)$&4&0&0&0\\
Flat I &$(3,4)$&0&1&7&1\\
Flat I &$(3,5)$&4&0&0&0\\
Flat I &$(3,6)$&2&2&2&2\\
\bottomrule
\end{tabular}
\caption{Revised whole-field assembly, with valuation normalized by $v(7)=1$. Bounds zero mean that this leading row supports no field.}\label{canx:r3:tab:assembly}
\end{table}

Only the three possible total-seven rows need a mixed argument. Their determinant-positive contribution is at most one, and its reflected residue point is always
\[
 (\bar X',\bar Y')=(5,3),\qquad (\bar u,\bar v)=(6,6).
\]
For the nonflat rows this follows from the support polynomial $H+5$; for the flat row it follows from the sole surviving direction $(0:1:0:0)$ and support $\xi+2$. The preceding simple-lifting statements make this place rational and unramified.

After first-order exclusion, a nonempty $U$ contribution on these rows can only be a single quadratic point. Direct residue norms show
\[
 N_{\F_{49}/\F_7}(1+\bar X')\in\{2,3,4,5\}.
\]
None is $\pm1$. Multiplication by the possible determinant-positive factor $6=-1$ still gives none of $\pm1$. But $u$ in \eqref{canx:r3:eq:units} is a global unit, so $N_{K/\Q}(u)=\pm1$. The full local product is incompatible with that norm both when the quadratic place occurs alone and when it occurs with the rational conic place. Thus there can be no quadratic contribution on these three rows; only the rational conic contribution remains.

This norm test is imposed on the \emph{complete} local decomposition. It does not apply a global norm to a selected subset while ignoring other places. The chart coverage says that the only possible omitted contribution is precisely the rational point just included. Only the first of the two available global units is needed for this exclusion.

If every place has $D$ a unit, Table~\ref{canx:r3:tab:unit} bounds the whole degree by two. Otherwise Table~\ref{canx:r3:tab:assembly} does so. Since
\[
 [K:\Q]=\sum_{\mathfrak p\mid7}e_{\mathfrak p}f_{\mathfrak p},
\]
this proves Theorem~\ref{canx:r3:thm:r3}. No claim that the remaining upper bounds are attained is required.

%% file: sections/cusp_expansion.tex
\section{Local algebras, cusp supports and integral shifts}\label{cuspx:section}
\subsection{Local hypotheses and divided equations}\label{cuspx:census:sec:setting}
$A=\begin{pmatrix}a&b\\b&d\end{pmatrix}$ is rational and symmetric, $X=x/(1-x)$, $Y=y/(1-y)$, $h=a+X$, $D=(a+X)(d+Y)-b^2$, $K=\Q(x,y)$, $v$ the valuation of a place above $7$ normalized by $v(7)=1$. Throughout:
\begin{itemize}
\item[(H1)] $a,b,d,B,C\in\Z_7$ and $\bar b\neq0$; $r\in\Q$;
\item[(H2)] the signed Nahm equations $(1-x)^N=x^{Na}y^{Nb}$, $(1-y)^N=x^{Nb}y^{Nd}$ hold in $K$ for an integer $N$ prime to $7$;
\item[(H3)] $x,1-x,y,1-y$ are units at every place of $K$ above $7$ (local units);
\item[(H4)] $D\neq0$, $\Psi_1=rD^3$ and $\Psi_2=\Psi_3=\Psi_4=0$ in $K$ (that is $\lambda_1=r$, $\lambda_2=\lambda_3=\lambda_4=0$, $\Psi_j=D^{3j}\lambda_j$).
\end{itemize}
The four-global-unit hypothesis (N) means that $x,1-x,y,1-y$ are global units. It is imposed explicitly wherever a whole-field norm is used; the local algebra and lifting arguments require only (H3). With $R=\Psi_1$ and the divided forms defined below
\[\Psi_2=\tfrac{36}{5}R^2+DR\,A_G+D^2J_2,\qquad \Psi_3=\tfrac{2652}{25}R^3+DR^2U_G+D^2R\,V_G+D^3J_3,\qquad U_G=\tfrac{138}{5}A_G,\]
(H4) is equivalent to $E_1=E_2=E_3=0$, $\Psi_4=0$, where
\[E_1=R-rD^3,\quad E_2=J_2+rD^2A_G+\tfrac{36}{5}r^2D^4,\quad E_3=J_3+rD^2V_G+r^2D^4U_G+\tfrac{2652}{25}r^3D^6,\]
since $\Psi_2=D^2E_2$ and $\Psi_3=D^3E_3$ once $R=rD^3$. All coefficient denominators are prime to $7$; modulo $7$, $36/5\equiv3$, $138/5\equiv1$, $2652/25\equiv5$. The auxiliary polynomials can be reconstructed from the formal
coefficients of Section~\ref{deg:sec:coefficients} by the following exact
operations. Define
\[
 {\cal C}_j(P)=h^jP(h-a,b^2/h-d),\qquad R_h={\cal C}_3(R).
\]
For a polynomial $Q$ in the image of ${\cal C}_j$, its lift
${\cal L}_j(Q)$ has the form
$u_0+\sum_{i=1}^j(u_iX^i+v_iY^i)$.
It is obtained uniquely as follows. In descending order $i=j,\ldots,1$,
subtract the coefficient of $h^{j+i}$ times $h^j(h-a)^i$ from $Q$,
recording that coefficient as $u_i$. Then, again in descending order,
subtract $v_i h^{j-i}(b^2-dh)^i$, where $v_i$ is the current coefficient
of $h^{j-i}$ divided by $b^{2i}$. The remaining coefficient of $h^j$ is
$u_0$. In the applications these divisions are exact in the polynomial
ring, so the resulting identities also specialize at $b=0$.
Set
\begin{align*}
 T_2&=(\Psi_2-36R^2/5)/D,&
 A_G&={\cal L}_2({\cal C}_5(T_2)/R_h),\\
 J_2&=(\Psi_2-36R^2/5-DR A_G)/D^2,&U_G&=138A_G/5,\\
 T_3&=(\Psi_3-2652R^3/25-DR^2U_G)/D^2,&
 V_G&={\cal L}_4({\cal C}_7(T_3)/R_h),\\
 J_3&=(\Psi_3-2652R^3/25-DR^2U_G-D^2RV_G)/D^3.
\end{align*}
Polynomial division and substitution verify that every quotient is exact
and that all coefficients are in $\Z_{(7)}$. These formulas specify
$A_G,U_G,V_G,J_2,J_3$ without a coefficient file.

The determinant cubic is
\[q(h)=d(d-1)h^3-b^2(b+2d-1)h^2+b^3(2a-1+b)h-b^3a(a-1).\]
A root $\bar h$ of $\bar q$ is \emph{admissible} if $\bar h\neq0$ and $\bar X=\bar h-\bar a$, $1+\bar X$, $\bar Y=\bar b^2/\bar h-\bar d$, $1+\bar Y$ are nonzero.

\subsection{Places, residue discs and the degree}\label{cuspx:census:sec:count}
\begin{lemma}[degree count]\label{cuspx:census:lem:count}
Let $z$ be a closed point of the residue plane of degree $f(z)$, $q=7^{f(z)}$, and $z_1\in\F_q^2$ one of its geometric points. If the system (H2), (H4) has at most $L$ solutions in the residue polydisc of $z_1$ (points of $\mathcal O_{\mathbb C_7}^2$ with residue $z_1$), then
\[\sum_{\mathfrak p\mid7,\ \mathfrak p\text{ over }z} e_{\mathfrak p}f_{\mathfrak p}\ \le\ f(z)\,L,\qquad [K:\Q]=\sum_{\mathfrak p\mid 7}e_{\mathfrak p}f_{\mathfrak p}.\]
\end{lemma}
\begin{proof}
$\sum_{\mathfrak p}e_{\mathfrak p}f_{\mathfrak p}$ over the places with residue in the orbit of $z$ is the number of embeddings $\sigma:K\to\overline{\Q}_7$ with $(\sigma x,\sigma y)$ reducing into that orbit. $\mathrm{Gal}(\overline\Q_7/\Q_7)$ permutes the $f(z)$ geometric points transitively and maps solutions to solutions, so each geometric point receives the same number of embeddings. Since $K=\Q(x,y)$, $\sigma\mapsto(\sigma x,\sigma y)$ is injective, and its image lies among the solutions in the polydisc of the residue.
\end{proof}

\begin{lemma}[local form and local length]\label{cuspx:census:lem:local}
Let $z_1=(\bar x_1,\bar y_1)$ with $\bar x_1,1-\bar x_1,\bar y_1,1-\bar y_1\neq0$ and write $u=\omega(u)\langle u\rangle$ (Teichm\"uller and principal parts) for local units.
(a) Every solution of (H2) in the polydisc of $z_1$ satisfies $\langle1-x\rangle=\langle x\rangle^a\langle y\rangle^b$, $\langle1-y\rangle=\langle x\rangle^b\langle y\rangle^d$; in $X$, $Y$ these are
$\Phi_1=\langle1+X\rangle^{a-1}\langle X\rangle^{-a}\langle Y\rangle^{-b}\langle1+Y\rangle^{b}-1=0$, $\Phi_2=\langle1+Y\rangle^{d-1}\langle Y\rangle^{-d}\langle X\rangle^{-b}\langle1+X\rangle^{b}-1=0$, power series in $s=X-X_1$, $t=Y-Y_1$ with coefficients in $\Z_q$. Their reductions modulo $7$ through total degree $6$ depend only on $(\bar a,\bar b,\bar d,z_1)$ (Lucas: $\binom{\alpha}{k}\equiv\binom{\bar\alpha}{k}$ for $k\le6$).
(b) If $\ell=\dim_{\F_q}\F_q[[s,t]]/(\bar\Phi_1,\bar\Phi_2)\le6$ (the \emph{local Nahm length}), the system has at most $\ell$ solutions in the polydisc; for further equations $E_i\in\Z_q[[s,t]]$, the common solutions number at most $\dim_{\F_q}\F_q[[s,t]]/(\bar\Phi_1,\bar\Phi_2,\bar E_i)$.
(c) If $\partial\bar\Phi_1/\partial s(0)\neq0$ then $\ell=\ord_t\bar\Phi_2(s(t),t)$, where $\bar\Phi_1(s(t),t)=0$.
(d) If $\bar D(z_1)\neq0$ then $\ell=1$: the solution in the polydisc, if any, is unique and unramified.
\end{lemma}
\begin{proof}
(a) Split $(1-x)^N=x^{Na}y^{Nb}$ into Teichm\"uller and principal parts; $N$ is prime to $7$, so $N$-th roots of principal units are unique. (b) $M=\Z_q[[s,t]]/(\Phi_1,\Phi_2,E_i)$ is complete for the $(7,s,t)$-adic topology and $M/7M$ is finite-dimensional, so $M$ is a finitely generated $\Z_q$-module generated by $\dim M/7M$ elements (complete Nakayama lemma); the solutions are distinct $\Z_q$-algebra maps $M\to\mathcal O_{\mathbb C_7}$, hence at most $\dim_{\Q_q}M\otimes\Q_q\le\dim M/7M$. (c) Implicit function theorem. (d) The Jacobian of the logarithmic system is $D/(X(1+X)Y(1+Y))$ up to sign, a unit; Hensel.
\end{proof}

\subsection{Where the places lie}\label{cuspx:census:sec:where}
\begin{lemma}[conic identity]\label{cuspx:census:lem:G1}
Exactly over $\Q$: $R(h-a,\,b^2/h-d)=\tfrac{5}{24}\,q(h)^2/h^3$. Modulo $7$, $h^{j}F(h-a,(b^2-dh)/h)$ is divisible by $q$ as a polynomial in $h,a,b,d,B,C$ for $F=J_2,J_3,A_G,V_G$ ($j=4,6,2,4$). In particular, at an exact root $H\in\Z_q$ of $q$, $J_2(Z_0)/7$ and $J_3(Z_0)/7$ ($Z_0=(H-a,b^2/H-d)$) are polynomials in the residues $(\bar h,\bar a,\bar b,\bar d,\bar B,\bar C)$, independent of deeper digits.
\end{lemma}

\begin{lemma}[location]\label{cuspx:census:lem:where}
Assume (H1)--(H4) and let $\mathfrak p\mid7$.
(a) If $v_{\mathfrak p}(D)=0$ then $r\in\Z_7$ and the residue point is a closed point with $\bar X,1+\bar X,\bar Y,1+\bar Y,\bar D\neq0$, $\bar\Psi_1=\bar r\bar D^3$, $\bar\Psi_2=\bar\Psi_3=\bar\Psi_4=0$ (a \emph{determinant-unit point}).
(b) If $v_{\mathfrak p}(D)>0$ then the residue lies on $\bar D=0$ with $\bar h$ an admissible root of $\bar q$. If $v_7(r)<0$, case (b) holds at every place.
\end{lemma}
\begin{proof}
$R$ has $7$-integral coefficients and $X,Y$ are integral by (H3), so $R$ is integral at $\mathfrak p$. (a) $rD^3=R$ forces $v(r)\ge0$; reduce the relations. (b) If $r\in\Z_7$, $v(R)\ge3v(D)>0$. If not, $\Psi_2=0$ gives $\tfrac{36}{5}R^2=-D(RA_G+DJ_2)$, so $2v(R)\ge v(D)>0$. Thus $\bar R=\bar D=0$; $\bar b\neq0$ gives $\bar h\neq0$ and $\bar Y=\bar b^2/\bar h-\bar d$, and Lemma~\ref{cuspx:census:lem:G1} gives $\bar q(\bar h)=0$; admissibility is (H3). If $v(r)<0$, case (a) is impossible.
\end{proof}

Exact enumeration of $\F_7^3$ gives the strata: $49$ residues with $\bar b=0$; $2$ with $\bar q\equiv0$, namely $(0,1,0)$, $(1,6,1)$; $22$ with an admissible repeated root; $270$ whose admissible roots are all simple. The $270$ carry $320$ admissible points ($202$ of degree one, $82$ of degree two, $36$ of degree three; $16$ residues have none). The $b\equiv0$ and $q\equiv0$ strata are not treated here.

\subsection{Simple admissible roots}\label{cuspx:census:sec:simple}
Fix a simple admissible point $z$, its exact root $H\in\Z_q$ of $q$ and $Z_0=(X_0,Y_0)=(H-a,\,b^2/H-d)$; then $D(Z_0)=0$ and, by Lemma~\ref{cuspx:census:lem:G1}, $R(Z_0)=0$. Use the coordinates $t=X-X_0$ and $D$ (valid since $\partial D/\partial Y=a+X$ is a unit). Then $R=S(t,D)\in\Z_q[[t,D]]$ with
\[S(t,0)=\tfrac5{24}\,\frac{q(H+t)^2}{(H+t)^3}=\sigma t^2+O(t^3),\qquad \sigma=\tfrac5{24}\,\frac{q'(H)^2}{H^3}\in\Z_q^\times,\qquad \partial_DS(0,0)=\frac{R_Y(Z_0)}{H}.\]
Along the curve $R=0$ one has $D=D_2t^2+O(t^3)$; for any $F$ the \emph{jets} $F_0,\dots,F_3$ are the coefficients of $F$ along $R=0$ parametrized by $t$ (these are the jets stored in the coefficient tables). Write $H_2,H_3$ for the jets of order $2,3$ of $J_2$, $T_3$ for that of order $3$ of $J_3$, $A_{G,1}$ for the linear jet of $A_G$, $V_{G,2}$ for the quadratic jet of $V_G$, and $A_2=J_2(Z_0)/7$, $A_3=J_3(Z_0)/7$.

\begin{lemma}[facts at the simple points]\label{cuspx:census:lem:facts}
At all $320$ points and all $49$ shift residues ($15{,}680$ cases): \textup{(F1)} $R_Y(Z_0)$ is a unit; \textup{(F2)} $R=0$ is tangent to $D=0$ at $Z_0$ and $D_2\neq0$; \textup{(F3)} $J_2(Z_0)\equiv0$ and the linear jet of $J_2$ vanishes mod $7$; \textup{(F4)} $A_2$, $A_3$ are units; \textup{(F5)} the jets of $J_3$ of orders $0,1,2$ vanish mod $7$; \textup{(F6)} $A_G(Z_0)\equiv V_G(Z_0)\equiv0$ and the linear jet of $V_G$ vanishes mod $7$; \textup{(F7)} the local Nahm length is $\ell=3$.
\end{lemma}
These finite computations are specified in Section~\ref{cuspx:census-recipe}: Lemma~\ref{cuspx:census:lem:G1} gives $A_2,A_3$ and Lemma~\ref{cuspx:census:lem:local}(c) gives $\ell$. For a place $\mathfrak p$ over $z$ put $L=K_{\mathfrak p}\supseteq\Q_q$; by (F7) and Lemma~\ref{cuspx:census:lem:local}, $[L:\Q_q]\le3$, so $e=e(L/\Q_q)\le3$. Let $t=X(\mathfrak p)-X_0$, $\tau=v(t)>0$, $\delta=v(D(\mathfrak p))$.

\begin{lemma}[two balances]\label{cuspx:census:lem:balance}
If $v(D(\mathfrak p))>0$ and $v(r)=-m$ ($m\in\Z$), then (i) $u:=2\delta-m>0$ and (ii) $2\tau=\delta$ with $D(\mathfrak p)=D_2t^2(1+\varepsilon)$, $v(\varepsilon)>0$.
\end{lemma}
\begin{proof}
(i) In $E_2$, $v(rD^2)=u$, $v(r^2D^4)=2u$, while $J_2(\mathfrak p)$ and $A_G(\mathfrak p)$ have positive valuation (F3, F6). If $u\le0$, the term $\tfrac{36}5r^2D^4$ has strictly the smallest valuation, so $E_2\neq0$. (ii) $S(t,D)=rD^3$ with $v(rD^3)=3\delta-m>\delta$. In $S=\sigma t^2+\partial_DS\,D+(\text{terms of order}\ge3\text{ or containing }tD,D^2)$ the minimum valuation must occur twice; $2\tau<\delta$ or $2\tau>\delta$ leaves $\sigma t^2$ or $\partial_DS\,D$ alone, so $2\tau=\delta$ and $D\equiv-(\sigma/\partial_DS)t^2=D_2t^2$ to leading order.
\end{proof}

Put $\zeta=rD_2^2t^3$, so $rD^2=\zeta t(1+\varepsilon)^2$, $r^2D^4=\zeta^2t^2(1+\varepsilon)^4$, $v(\zeta)=3\tau-m=u-\tau$. By (F3), (F6) and Lemma~\ref{cuspx:census:lem:balance},
\[J_2(\mathfrak p)=7A_2+jt+H_2t^2+H_3t^3+\rho_2,\qquad A_G(\mathfrak p)=\alpha_0+A_{G,1}t+\rho_A,\]
with $v(j),v(\alpha_0)\ge1$, $v(\rho_A)>\tau$ and $v(\rho_2)\ge\min(4\tau,\,3\delta-m)>2\tau$ (the second entry is the correction $\partial_DJ_2\,(D(\mathfrak p)-D_{\rm curve}(t))$, of valuation $v(rD^3)=2\tau+u$; for $r\in\Z_7$ it is $\ge6\tau$); dividing $E_2$ by $t^2$,
\[\frac{E_2}{t^2}=\frac{7A_2}{t^2}+\frac jt+H_2+H_3t+\zeta(1+\varepsilon)^2\Big(\frac{\alpha_0}{t}+A_{G,1}+\frac{\rho_A}{t}\Big)+\tfrac{36}5\zeta^2(1+\varepsilon)^4+\frac{\rho_2}{t^2}.\tag{$*$}\]
Since $e\le3$, $\tau\in\{1/3\}\cup\{1/2\}\cup[2/3,\infty)$.

\begin{proposition}[integral scalar: no determinant-nonunit place]\label{cuspx:census:prop:int}
If $r\in\Z_7$, a place $\mathfrak p$ with $v(D)>0$ over the simple point $z$ exists only if, for the shift residue, $\bar H_2=0$, $\bar H_3\neq0$, $\bar T_3\neq0$ and $\bar A_2\bar T_3=\bar A_3\bar H_3$. This happens in none of the $15{,}680$ cases; hence, when all admissible determinant roots of the matrix residue are simple, every place of $K$ lies over a determinant-unit point (Lemma~\ref{cuspx:census:lem:where}).
\end{proposition}
\begin{proof}
Here $v(\zeta)=v(r)+3\tau>0$. If $\tau\ge2/3$, $7A_2/t^2$ has valuation $1-2\tau<0$ (F4) and every other term of $(*)$ has larger valuation. If $\tau=1/2$, $(*)$ forces $\bar H_2\neq0$; then in $E_3$ every term other than $7A_3$ has valuation $\ge3/2$ (F5 for the jets, $v(rD^2)\ge2$, $v(rD^3)\ge3$), so $v(E_3)=1$. If $\tau=1/3$: a unit $H_2$ makes $(*)\equiv\bar H_2\neq0$; otherwise the terms of valuation $1/3$ are $7A_2/t^2+H_3t$, so $\bar H_3\neq0$ and $t^3\equiv-7A_2/H_3$, and in $E_3=7A_3+T_3t^3+(\text{valuation}\ge4/3)$ also $t^3\equiv-7A_3/T_3$; comparing gives the stated condition; if $\bar H_3=0$ the term $7A_2/t^2$ is alone. The finite check: $\bar H_2\neq0$ in $14{,}202$ cases, $\bar H_2=\bar H_3=0$ in $770$, and $\bar A_2\bar T_3\neq\bar A_3\bar H_3$ in the remaining $708$.
\end{proof}

\begin{proposition}[pole of the scalar: support]\label{cuspx:census:prop:pole}
Let $v_7(r)=-m<0$, $r'=7^mr$, and let $\mathfrak p$ lie over the simple point $z$. Then $m=1$, $e=3$, $f(L/\Q_q)=1$, $\zeta$ is a unit, and $\bar\zeta\in\F_q^\times$ satisfies
\begin{gather*}Q(\bar\zeta):=\bar H_2+\bar A_{G,1}\bar\zeta+3\bar\zeta^2=0,\qquad c(\bar\zeta):=\bar T_3+\Big(\tfrac{\bar J_{3,Y}\bar D_2}{\bar R_Y}+\bar V_{G,2}\Big)\bar\zeta+\bar A_{G,1}\bar\zeta^2+5\bar\zeta^3\neq0,\\ \bar r'=-\frac{\bar\zeta\,c(\bar\zeta)}{\bar A_3\bar D_2^2},\end{gather*}
with $J_{3,Y}$, $R_Y$ evaluated at $Z_0$. In particular $\bar r'$ is one of the \emph{candidates} of $z$ for the shift residue.
\end{proposition}
\begin{proof}
\emph{$\tau\ge2/3$.} Put $z=v(\zeta)$ and $c=1-2\tau<0$.
If $z\ge0$, the term $7A_2/t^2$ is uniquely smallest in $(*)$.
If $z<0$, the quadratic term has valuation $2z$, whereas the linear
$\zeta$ term has larger valuation. Also
$v(j/t)\ge1-\tau>c$ and
$v(\zeta\alpha_0/t)\ge2\tau-m+1=c+u>c$.
All other terms have valuation greater than $c$ or greater than $2z$.
Cancellation therefore requires $2z=c$, or $\tau=(2m+1)/8$.
Its denominator is eight, so it cannot belong to $(1/e)\Z$ for $e\le3$.
\emph{$\tau=1/2$.} Then $e=2$, $u=2-m>0$ forces $m=1$ and $v(\zeta)=1/2$. In $E_3$: the jets of $J_3$ give $7A_3$ plus terms of valuation $\ge3/2$ (F5); the correction from $D(\mathfrak p)$ to the curve $R=0$ is $(J_{3,Y}/R_Y)\,rD^3$ of valuation $\ge3\delta-m=2$; $v(rD^2V_G)\ge2$ (F6), $v(r^2D^4U_G)\ge5/2$, $v(r^3D^6)=3$. So $v(E_3)=1$.
\emph{$\tau=1/3$.} Then $e=3$, $f(L/\Q_q)=1$, $v(\zeta)=1-m\le0$. If $v(\zeta)<0$, $\tfrac{36}5\zeta^2$ is alone in $(*)$; hence $v(\zeta)=0$, $m=1$, and $(*)$ modulo the maximal ideal of $L$ (residue field $\F_q$) gives $Q(\bar\zeta)=0$, $\bar\zeta\neq0$. In $E_3$, by (F5), (F6) and the curve correction $(J_{3,Y}/R_Y)rD^3=(J_{3,Y}/R_Y)\zeta D_2t^3(1+\varepsilon)^3$,
\[E_3=7A_3+t^3\Big(T_3+\big(\tfrac{J_{3,Y}D_2}{R_Y}+V_{G,2}\big)\zeta+\tfrac{138}{5}A_{G,1}\zeta^2+\tfrac{2652}{25}\zeta^3\Big)+\rho_3,\qquad v(\rho_3)>1 .\]
As $v(t^3)=1$ and $7/t^3=r'D_2^2/\zeta$, reduction gives $\bar A_3\bar r'\bar D_2^2/\bar\zeta+c(\bar\zeta)=0$; $\bar r'\neq0$ gives $c(\bar\zeta)\neq0$.
\end{proof}

\begin{corollary}[support census]\label{cuspx:census:cor:support}
For $v_7(r)<0$ on a simple-root residue, every place lies over the support $S(\bar B,\bar C,\bar r')$ of the points having $\bar r'$ among their candidates, and $[K:\Q]\le3\sum_{z\in S}f(z)$. If $m\ge2$, or if the support is empty, there is no Nahm point at all. Of the $270\cdot49\cdot6=79{,}380$ configurations, the support degree is $0$ in $73{,}816$, $1$ in $4{,}670$ (degree $\le3$) and $2$ in $894$; none has support degree $3$. The $894$ configurations of degree $2$ are treated by Theorem~\ref{cuspx:simple-poles} below.
\end{corollary}
\begin{proof}
Lemma~\ref{cuspx:census:lem:where}(b), Proposition~\ref{cuspx:census:prop:pole}, Lemma~\ref{cuspx:census:lem:count} with $L=3$ (F7); $[K:\Q]\ge1$ gives the nonexistence statements. The counts are those of the finite census.
\end{proof}

\begin{remark}
The candidate test retains only roots $\bar\zeta\in\F_q$ of $Q$, excludes $\bar\zeta=0$ and requires $\bar r'\in\F_7^\times$. Proposition~\ref{cuspx:census:prop:pole} shows that this is exhaustive: $\bar\zeta$ lies in the residue field $\F_q$ of $L$ because $f(L/\Q_q)=1$; $\bar\zeta=0$ and $c(\bar\zeta)=0$ are impossible. Points where (F5) would fail would only be excluded, so the census is never too small; (F5) in fact holds everywhere.
\end{remark}

\subsection{Repeated admissible roots}\label{cuspx:census:sec:repeated}
\begin{proposition}[repeated roots]\label{cuspx:census:prop:rep}
Each of the $22$ repeated-root residues has exactly one admissible root, rational, of multiplicity $\mu=2$ ($16$ residues) or $3$ ($6$); the local Nahm length is $\ell=\mu+2$ and $D$ has order $\ell-1$ along the first Nahm curve $C_1$. Let $o_F$ be the order along $C_1$ of $\bar F$ ($F=R,J_2,J_3$), truncated at $\ell$.
(a) If $r\in\Z_7$, the places with $v(D)>0$ contribute at most $L(\bar B,\bar C)=\min(\ell,o_R,o_{J_2},o_{J_3})$ to $[K:\Q]$; $L=3$ in $672$ and $L=4$ in $406$ of the $1{,}078$ shift cases.
(b) If $v_7(r)<0$, every place lies over the centre and $[K:\Q]\le\ell$; in particular $\le4$ at the $16$ square residues, which include the canonical residues $(1,3,5)$ and $(1,5,6)$.
\end{proposition}
\begin{proof}
Lemma~\ref{cuspx:census:lem:where}(b) puts every determinant-nonunit place at the unique admissible root. The Nahm local algebra is $\F_q[[t]]/(t^\ell)$ along $C_1$ (Lemma~\ref{cuspx:census:lem:local}(c)). Since $\bar D|_{C_1}$ has order $\ell-1$ and $2(\ell-1)\ge\ell$, the terms $rD^3$, $rD^2A_G$, $r^2D^4$, $rD^2V_G$, \dots\ of $E_1,E_2,E_3$ vanish modulo $t^\ell$ when $r\in\Z_7$, so the local algebra before imposing $\Psi_4$ is $\F_q[[t]]/(t^{L})$ Adding $\Psi_4$ takes a quotient and can only decrease the length. Lemmas~\ref{cuspx:census:lem:local}(b) and~\ref{cuspx:census:lem:count} therefore apply. In (b), Lemma~\ref{cuspx:census:lem:count} with $L=\ell$. The lengths, orders and jets are those of the finite calculation.
\end{proof}

\subsection{The branch-one determinant-unit census}\label{cuspx:census:sec:branch1}
\begin{proposition}[under (N)]\label{cuspx:census:prop:b1}
Let $r\in\Z_7$ on a repeated-root residue and assume (N). Let $\mathcal U$ be the set of determinant-unit points of Lemma~\ref{cuspx:census:lem:where}(a) with the given $\bar r$, and $c$ the cap $L(\bar B,\bar C)$ of Proposition~\ref{cuspx:census:prop:rep}(a). Then $[K:\Q]=\sum_{z\in S}f(z)+\mu'$ for some $S\subseteq\mathcal U$ and $0\le\mu'\le c$ such that, for $w\in\{X,1+X,Y,1+Y\}$, $\prod_{z\in S}N_{\F_7(z)/\F_7}(\bar w(z))\cdot\bar w(z_0)^{\mu'}\equiv\pm1\pmod7$, $z_0$ the centre. The resulting maximal degree is $3$ in $7{,}528$ of the $7{,}546$ configurations and $6$ in $18$.
\end{proposition}
\begin{proof}
For a unit $w$ of $K_{\mathfrak p}$, $N_{K_{\mathfrak p}/\Q_7}(w)\equiv N_{k_{\mathfrak p}/\F_7}(\bar w)^{e_{\mathfrak p}}$; under (N), $\prod_{\mathfrak p\mid7}N_{K_{\mathfrak p}/\Q_7}(w)=N_{K/\Q}(w)=\pm1$. At a determinant-unit point the place is unique and unramified with $f=f(z)$ (Lemma~\ref{cuspx:census:lem:local}(d)); at the centre the places contribute $\mu'\le c$ with norm $\bar w(z_0)^{\mu'}$ ($z_0$ rational). Every place is of one of the two kinds (Lemma~\ref{cuspx:census:lem:where}).
\end{proof}
The $18$ configurations of degree $6$ are reduced to $4$ ($3$ under (N)) by Proposition~\ref{cuspx:repeated-mixed}. The norm hypothesis is used only on the four-global-unit route.

\subsection{The determinant-unit census of all residues}\label{cuspx:census:sec:dunit}
This census treats all integral matrix and shift residues. Let (H1$'$) be (H1) without the condition $\bar b\ne0$.
\begin{proposition}[determinant-unit places]\label{cuspx:census:prop:dunit}
Assume (H1$'$), (H2)--(H4) and $r\in\Z_7$. If $v_{\mathfrak p}(D)=0$, the residue $z$ of $(X,Y)$ at $\mathfrak p$ is an \emph{admissible determinant-unit point} of the fibre $(\bar a,\bar b,\bar d,\bar B,\bar C;\bar r)$: $\bar X,1+\bar X,\bar Y,1+\bar Y,\bar D\neq0$, $\bar\Psi_2=\bar\Psi_3=\bar\Psi_4=0$ and $\bar\Psi_1=\bar r\bar D^3$ at $z$. The place $\mathfrak p$ is the only place over $z$, $e_{\mathfrak p}=1$ and $f_{\mathfrak p}=f(z)$. Under (N), if every place above $7$ has $v(D)=0$, then $[K:\Q]=\sum_{z\in S}f(z)$ for a set $S$ of such points with $\prod_{z\in S}N_{\F_7(z)/\F_7}(\bar w(z))\equiv\pm1$ for $w=X,1+X,Y,1+Y$.
\end{proposition}
\begin{proof}
The $\Psi_j$ have $7$-integral coefficients and vanish exactly at $(X,Y)$ (H4); $D$ is a unit, so $r=\Psi_1/D^3$ reduces to $\bar\Psi_1/\bar D^3$. Lemma~\ref{cuspx:census:lem:local}(d) gives the unique unramified place. The norm relation is that of Proposition~\ref{cuspx:census:prop:b1} without a centre.
\end{proof}
\textbf{Exhaustiveness of the computed list}. For each of the $343\cdot49=16{,}807$ classes the three specialised polynomials are first divided by their common factors with $X(1+X)Y(1+Y)D$; this removes no admissible point. If one $\bar\Psi_j$ vanishes identically, or the three have a common factor, the class is \emph{flagged}. Otherwise the system is zero-dimensional: a nonzero pairwise resultant in $Y$ vanishes at the $\bar X$ of every common zero, and over each irreducible factor $\varphi(\bar X)$ the common $\bar Y$ are the roots of the gcd in $\F_q[Y]$ of the three polynomials ($q=7^{\deg\varphi}$). Each closed point is then tested for admissibility and $\bar r\in\F_7$, and its letter norms are taken from its residue field $\F_7(\bar X,\bar Y)$. In a flagged class the scalar equation $\bar\Psi_1-\bar r\bar D^3$ is added at each $\bar r$; the fibre is then either finite (its points are listed) or a \emph{curve fibre}. There are $514$ flagged classes ($86$ with a vanishing equation, $428$ with a common curve) and exactly $428$ curve fibres; these are the $428$ tuples of the curve calculation in Section~\ref{cuspx:flat}, $60$ of them with $\bar b=0$, which are treated together with the other curve fibres.  The maxima of the four-unit spectra are $4$ on the regular, covered and uncovered simple-root strata, $0$ on the repeated-root stratum, $6$ in $8$ fibres with $\bar b=0$, and, after refinement, $6$ in $24$ fibres and $2$ in $12$ fibres of the two zero-cubic residues; finite refined fibres elsewhere have no positive allowed degree.

\paragraph{Consequences for an integral scalar.}
On a simple-root residue, Proposition~\ref{cuspx:census:prop:int}
leaves only determinant-unit places. The finite spectrum therefore gives
degree at most four outside the $368$ simple-root curve fibres.
On a repeated-root residue the centre is included before applying norms,
as in Proposition~\ref{cuspx:census:prop:b1}.
The eight zero-coupling fibres with maximum six have two cubic points and
an inadmissible centre; the symplectic exclusion of
Section~\ref{deg:sec:entry-completion}, together with the explicit
zero-coupling calculation preceding it, handles these pairs.
The twenty-four zero-cubic fibres are subject to the first-order
constraint of Section~\ref{deg:sec:z0-replacement}. In neither application
are four global norm conditions imposed on a proper subset of the places.

\subsection{A finite recipe for the cusp tables}\label{cuspx:census-recipe}
All the finite calculations in the preceding propositions use the
following domain and tests. Enumerate $(a,b,d)\in\F_7^3$, factor $q(h)$,
and retain a factor only when it is coprime to
$h(h-a)(h-a+1)(b^2-dh)(b^2+(1-d)h)$. Multiplicities are retained.
For every closed point, work in its finite residue field, and enumerate
$(B,C)\in\F_7^2$. The principal-unit equations in
Lemma~\ref{cuspx:census:lem:local} determine the local Nahm length:
truncate at degree six, eliminate an equation with unit first derivative,
and take the order of the other. The eliminated equation is obtained by
successive substitution, so its coefficients are uniquely determined.

For simple points, lift $q$ by Hensel and evaluate its exact root $H$.
For $F=J_2,J_3$, put $P=h^jF(h-a,b^2/h-d)$ and compute the exact polynomial quotient $\bar M=\bar P/\bar q$ over $\F_7$ before specializing. Lift its coefficients to $0,\ldots,6$. Then $(P-qM)/7$ has coefficients in $\Z_{(7)}$; its reduction, evaluated at $\bar H$ and divided by $\bar H^j$, gives $F(Z_0)/7$ modulo seven. The same conic divisibility is checked for $A_G,V_G$. This avoids choosing unspecified deeper parameter
digits. Solve $R=0$ in the coordinates $(t,D)$ through degree three;
read $D_2,H_2,H_3,T_3,A_{G,1},V_{G,2}$ from that series.
Test (F1)--(F7) and the explicit compatibility predicates in
Propositions~\ref{cuspx:census:prop:int} and~\ref{cuspx:census:prop:pole}.
For repeated points, solve the first Nahm equation through degree five
and reduce $R,J_2,J_3$ in its local algebra. This gives the orders and
caps of Proposition~\ref{cuspx:census:prop:rep}.

For the determinant-unit calculation, set
$S=X(1+X)Y(1+Y)D$. At each of the $7^5$ matrix/shift residues and each
$\bar r\in\F_7$, use the saturated ideal
\[
 (\bar\Psi_1-\bar rD^3,\bar\Psi_2,\bar\Psi_3,
   \bar\Psi_4,1-zS)\cap\F_7[X,Y].
\]
A zero-dimensional ideal is decomposed into its closed points by
resultants and polynomial gcds over the appropriate finite extensions.
A positive-dimensional ideal is retained as a curve fibre. Constant-in-$Y$
polynomial pairs are tested by their gcd in $X$; a resultant convention
is never used to discard such a pair. Compute each norm in the full
residue field $\F_7(X,Y)$, including the extension degree of $Y$ over
$\F_7(X)$. Finally enumerate all subsets of the finite supports and all
centre multiplicities between zero and the proved local cap. Retain a
choice exactly when all available whole-field unit norms are $1$ or $-1$.
These finite domains account for all possible residue points, and the
local-length lemmas account for their ramified lifts. No floating-point
root test occurs in this procedure.

\subsection{Weierstrass selection and the eighteen mixed fibres}
\label{cuspx:selection}
At an admissible centre choose a lift $(X_0,Y_0)$ on $D=0$ in the
unramified coefficient ring $\mathcal O_q$. Put $Y=Y_0+\eta$ and solve
$\Phi_1=0$ for $X=X_0+\xi(\eta)$. Its $X$-derivative is a unit.
Let $\Phi(\eta)$ be the second principal-unit equation on this curve,
and write $\Phi=UW$, where $U$ is a unit power series and $W$ is its
monic distinguished Weierstrass polynomial. Its degree $g$ is the local
Nahm length. The equations are in principal-unit form, so roots of
$W$ do not include the extra roots obtained by merely setting a
$p$-adic logarithm equal to zero.

Along this curve, differentiation of the second logarithmic equation
and elimination of $\xi'$ give
\[
 D=-hY(1+Y)\frac{\Phi'}{1+\Phi}.
\]
Consequently $D=\mu W'\pmod W$, with $\mu$ a unit in the local ring.
At a root $\eta_i$,
\[
 v(D(\eta_i))=v(W'(\eta_i))
             =\sum_{j\ne i}v(\eta_i-\eta_j).
\]
If $D\ne0$ at a selected root it is simple. Let $W_S$ be the product
of the Galois orbits selected by embeddings of $K$. Then every
coefficient equation has zero remainder modulo $W_S$. In particular,
\[
 W_S\mid\Psi_2,\Psi_3,\Psi_4,
 \qquad
 \operatorname{rem}(R,W_S)
      =r\operatorname{rem}(D^3,W_S).
\]
These are exact necessary conditions, including at a pole of $r$.
When integral equations are wanted, multiply $E_1,E_2,E_3$ by
$7^m,7^{2m},7^{3m}$ for $r=7^{-m}r'$.

\begin{lemma}[The fixed first lift at a determinant-unit point]
\label{cuspx:unit-lift}
At a determinant-unit residue point, the solution of the principal-unit
Nahm equations modulo $49$ is unchanged by adding seven times an
integer to $a,b$ or $d$. The coefficient equations at that lift modulo
$49$ form an affine system over the residue field in the first digits
of $a,b,d,B,C,r$.
\end{lemma}
\begin{proof}
The point lifts uniquely by the unit Nahm Jacobian. The derivatives of
the logarithmic equations with respect to the exponents are logarithms
of units in the unramified field, hence belong to $7\mathcal O_q$.
A parameter displacement divisible by seven changes the equations by
$49\mathcal O_q$. The inverse Jacobian is integral, proving the first
claim. Taylor expansion of the polynomial coefficient equations at
this fixed lift proves the second. For a quadratic residue field,
expressing each equation in a fixed $\F_7$-basis gives ordinary linear
algebra over $\F_7$.
\end{proof}

\begin{proposition}[The mixed repeated-centre fibres]
\label{cuspx:repeated-mixed}
For each of the eighteen configurations left by
Proposition~\ref{cuspx:census:prop:b1}, the whole-field degree is at most
four. Under the four-global-unit hypothesis it is at most three.
\end{proposition}
\begin{proof}
Each fibre consists of one quadratic determinant-unit point and one
rational centre of Nahm length five and coefficient cap four. The
first-order system of Lemma~\ref{cuspx:unit-lift} has rank six and fixes
the five parameter digits and the scalar digit as in
Table~\ref{cuspx:forced-digits}. At these digits the centre's polynomial
has the factorization $W=(\eta-\eta_g)W_4$, with $\eta_g\in7\Z_7$ and
$W_4$ an irreducible Eisenstein quartic. Its four roots have valuation
$1/4$. Thus the selected contribution at the centre belongs to
$\{0,1,4\}$.

The divided second equation has Weierstrass degree four. Selecting all
four roots of $W_4$ would force $W_4=W(E_2)$. If $c_i$ and $j_i$ are
the coefficients of $\Phi$ and $E_2$ along the first Nahm curve, the
constant terms of these two quartics give the necessary equality
\[
 c:=\frac{c_1/7}{c_5}
       =\frac{j_0/7}{j_4}=:j\quad\hbox{in }\F_7.
\]
In every row of the table, $j=3c\ne c$. Hence the centre contributes
at most one when the quadratic point is present, and the total is at
most three. Without the quadratic point its coefficient cap is four;
the whole-field norm test in Proposition~\ref{cuspx:census:prop:b1}
excludes four at the centre alone. Every quantity in the comparison
is determined by the listed first digits: it is a coefficient of order
at most five evaluated modulo $49$, where the necessary binomial
coefficients have denominators prime to seven. The argument is
therefore valid for all deeper digits.
\end{proof}

\begin{longtable}{@{}lllrrr@{}}
\caption{The eighteen forced first-digit vectors. The first column is
$(a,b,d,B,C)$ modulo seven, $u$ is its next digit, and $r_1$ is the
next digit of $r$.}\label{cuspx:forced-digits}\\
\toprule Residue & $\bar r$ & $u$ & $r_1$ & $c$ & $j$\\\midrule\endfirsthead
\toprule Residue & $\bar r$ & $u$ & $r_1$ & $c$ & $j$\\\midrule\endhead
$(2,2,6,3,2)$&5&$(4,3,6,0,5)$&6&2&6\\
$(2,2,6,4,6)$&3&$(3,2,3,6,0)$&2&4&5\\
$(2,2,6,6,0)$&5&$(0,3,2,1,3)$&2&1&3\\
$(2,5,6,2,1)$&1&$(0,3,2,5,3)$&6&1&3\\
$(2,5,6,4,5)$&1&$(4,3,6,0,1)$&3&4&5\\
$(2,5,6,5,0)$&3&$(3,4,3,2,4)$&3&2&6\\
$(4,1,4,2,2)$&2&$(4,3,4,0,0)$&3&4&5\\
$(4,1,4,3,3)$&5&$(3,5,3,4,4)$&0&1&3\\
$(4,1,4,5,5)$&0&$(2,3,2,0,0)$&4&2&6\\
$(4,6,4,1,1)$&1&$(3,1,3,6,6)$&5&1&3\\
$(4,6,4,2,2)$&4&$(2,3,2,3,3)$&2&4&5\\
$(4,6,4,6,6)$&6&$(4,3,4,2,2)$&1&2&6\\
$(6,2,2,0,6)$&5&$(2,3,0,3,1)$&2&4&5\\
$(6,2,2,2,3)$&5&$(6,3,4,5,0)$&6&1&3\\
$(6,2,2,6,4)$&3&$(3,2,3,0,6)$&2&2&6\\
$(6,5,2,0,5)$&3&$(3,4,3,4,2)$&3&4&5\\
$(6,5,2,1,2)$&1&$(2,3,0,3,5)$&6&2&6\\
$(6,5,2,5,4)$&1&$(6,3,4,1,0)$&3&1&3\\\bottomrule
\end{longtable}

\begin{proposition}[Repeated centres and a scalar pole]
\label{cuspx:repeated-pole}
Assume four global units and integral shifts. On any of the twenty-two
repeated-root matrix residues, a pole of $r$ implies $[K:\Q]\le3$.
\end{proposition}
\begin{proof}
All places lie at the unique rational centre. Its local Nahm length is
at most five. Write $H=\bar X\in\{1,2,3,4,5\}$ and $n=[K:\Q]$.
The whole-field unit norms give $H^n,(H+1)^n\in\{1,-1\}$ in $\F_7$.
For each of the five values of $H$, at least one of $H,H+1$ has order
three or six, so $3\mid n$. Since $1\le n\le5$, one obtains $n=3$.
The assertion uses global norms; no corresponding exclusion of all
norm-free length-five possibilities is needed.
\end{proof}

\subsection{The simple-centre scalar-pole calculation}
\label{cuspx:simple-tree}
The $894$ scalar configurations with support degree two involve $146$
distinct pairs consisting of a matrix residue and an admissible centre.
There are $416$ configurations with two rational centres and $478$ with
one quadratic centre. At each geometric centre, the Weierstrass degree
is three. The following local observation reduces the possibilities
before any parameter search.

\begin{lemma}[Cubic selection and a scalar pole]\label{cuspx:cubic-selection}
Let $W=\eta^3+w_2\eta^2+w_1\eta+w_0$ be the local Nahm polynomial.
Write $a_2=J_2(0)/7$ and $a_3=J_3(0)/7$ on the first Nahm curve.
Suppose that $a_2$ or $a_3$ is a unit and all three roots are selected.
Then $v(w_0)=1$, $r$ has pole order one, and the three roots are a
single totally ramified cubic orbit. If $W$ is Eisenstein, the same
conclusion about the pole order holds without the condition on
$a_2,a_3$.
\end{lemma}
\begin{proof}
Let $\tau_i=v(\eta_i)$, $\delta_i=v(D(\eta_i))$ and
$u_i=2\delta_i-m$, where $m=-v_7(r)>0$. The residue coefficient
identities give
\[
 v(A_G(\eta_i)),v(U_G(\eta_i)),v(V_G(\eta_i))\ge1,
 \quad v(J_2(\eta_i))\ge\min(1,2\tau_i).
\]
If $\tau_i>1/2$ and $a_2$ is a unit, $v(J_2(\eta_i))=1$; if
$\tau_i\ge1/2$ and $a_3$ is a unit, $v(J_3(\eta_i))=1$.
The second equation has two scalar terms of valuations at least
$u_i+1$ and exactly $2u_i$. Integrality and the positive valuation of
$J_2$ force $u_i>0$.

Suppose first that all $\tau_i>1/2$. The second equation then gives
$u_i>1/3$. The three scalar terms of the third equation all have
valuation greater than one, contradicting a unit $a_3$. If $a_2$ is
a unit instead, the second equation forces $u_i=1/2$, hence
$\delta_i=(2m+1)/4$. The splitting field of a cubic over the unramified
coefficient field has ramification index dividing six. Thus
$\delta_i\in\frac16\Z$, a contradiction.

If the smallest root valuation is $1/2$, there is a quadratic factor
with two roots of valuation $1/2$ and a rational third root of
valuation at least one. The quadratic discriminant has valuation one,
so $\delta_i=1$ for every root. Then $u_i>0$ forces $m=1$; evaluating
at the rational root contradicts the unit $a_2$ or $a_3$ just as
above. The only remaining smallest valuation is $1/3$. The Newton
polygon is the segment from $(0,1)$ to $(3,0)$; $W$ is Eisenstein.
Since $7\equiv1\pmod3$, its tamely ramified cubic extension is cyclic,
and the pairwise root distances are $1/3$. Thus $\delta_i=2/3$,
and $u_i>0$ gives $m<4/3$, or $m=1$.
\end{proof}

\begin{lemma}[The other centre at pole order one]\label{cuspx:other-centre}
Let $m=1$ and $v(w_0)\ge2$. If $a_3$ is a unit, no root at the centre
can be selected. If $a_3$ is not a unit but $a_2$ is, only a quadratic
pair of valuation $1/2$ can be selected. At an Eisenstein centre the
selected contribution is either zero or three.
\end{lemma}
\begin{proof}
For a selected root $u_i>0$. A root of valuation greater than $1/2$
is excluded by the third-equation argument or by the impossible
value $\delta_i=3/4\notin\frac16\Z$ from the second equation.
At a root of valuation $1/2$, the quadratic factor gives $\delta_i=1$,
so the scalar terms of $E_3$ have valuation at least two; a unit
$a_3$ makes $J_3$ have valuation one. Finally an Eisenstein cubic is
irreducible and the equations have coefficients in the unramified
base field, so they select all its roots or none.
\end{proof}

\paragraph{The finite classification of the centres.}
For each of the $146$ pairs, enumerate the $7^3$ first digits of
$(a,b,d)$. Form $\Phi$ modulo $49$ by the integral binomial series
and solve $\Phi_1=0$. The constant of $W$ has valuation one for all
first digits at $114$ pairs, and valuation at least two for all first
digits at the other $32$. In the latter cases, for each shift residue
appearing in the support table and each of its $7^2$ first digits,
compute $a_2,a_3$ modulo seven; their affine forms have no common
zero. Thus the non-Eisenstein centres contribute at most two by
Lemma~\ref{cuspx:cubic-selection}. Only the $114$ Eisenstein pairs
require the following search.

\paragraph{Necessary equations and truncation.}
At an Eisenstein centre put
\[
 Q(s)=\left(\operatorname{rem}(\Psi_2,W),
 \operatorname{rem}(\Psi_3,W),\operatorname{rem}(\Psi_4,W),
 \bigwedge^2[\operatorname{rem}(R,W)\mid
                  \operatorname{rem}(D^3,W)]\right),
 \qquad s=(a,b,d,B,C).
\]
The wedge denotes all two-by-two minors. Selection of all three roots
implies $Q(s)=0$. Every entry is calculated from the principal-unit
series and the polynomial coefficients already defined in the text.
For a distinguished degree-$g$ polynomial, each coefficient of
$\operatorname{rem}(\eta^i,W)$ is divisible by
$7^{\lfloor i/g\rfloor}$: each use of the relation for $\eta^g$
lowers degree by at most $g$ and contributes a factor of seven.
Consequently the integral series can be truncated at
$M=3N+9$ when evaluating a cubic remainder modulo $7^N$.
The same argument applies to the remainder equation defining $W$;
its Jacobian in the lower coefficients of $W$ is invertible modulo
seven. These are fixed-precision bounds, independent of the size of
the chosen exponent representatives.

\paragraph{Parameter boxes and terminal certificates.}
A node is a residue class $s_n+7^n\Z_7^5$. For each matrix--centre
pair, take every shift residue occurring with it in the $894$
support-degree-two configurations. Above each resulting five-parameter
residue class, begin with all $7^5$ first-digit classes. At level two the necessary equations use the
quadratic Taylor model modulo $7^5$ and the cubic model modulo $7^6$;
at level three use the quadratic model modulo $7^7$; from level four
use the affine model modulo $7^{n+4}$. Deeper digits are existential
variables, eliminated by Smith normal form over $\Z/7^k\Z$.
A node is discarded only when the corresponding congruence has no
continuation. Any unresolved node is subdivided into its $7^5$
children. Equivalently, without any linear-algebra optimization, at
each tested node enumerate $t\in(\Z/7^4\Z)^5$ in the indicated
Taylor model and retain exactly the classes $t\bmod7$ that satisfy
all its congruences. Four continuation digits suffice for every one
of the stated moduli. From depth eight onward, a box may instead
terminate by the exact-endpoint and uniqueness tests below. Use
depth twelve as the final bound: every box has died or passed those
terminal tests by that depth. The Taylor remainder estimates in
Corollary~\ref{cder:precision} show that every exact solution is
retained at every level.

Theorem~\ref{cder:theorem} proves the Taylor estimates required by this procedure. In Hasse form they are
\[
 v(Q^{[j]}(s)(7^nt)^j)\ge j(n-1)+1\quad(j\ge1),
 \qquad
 v(Q^{[j]}(s)(7^nt)^j)\ge j(n-1)+2\quad(j\ge2)
\]
for the cubic calculation. Here $Q^{[j]}$ denotes the homogeneous Taylor coefficient of order $j$. The second estimate gives exactly the precision needed at the first cubic-model step. The proof in Section~\ref{cder:section} treats the varying curve base point and the full implicit Weierstrass preparation, in addition to the fixed-precision truncation in $\eta$.

At a terminal parameter box, take each vector $s_*$ of
Table~\ref{cuspx:tab:endpoints} whose residue lies in the box,
substitute it in the signed Nahm equations after clearing its rational
exponents, and factor the resulting elimination polynomial over $\Q$.
The selected cubic factor must reduce to the prescribed centre.
Reduction in its quotient field must give $\Psi_2=\Psi_3=\Psi_4=0$
and $R=rD^3$ for a rational $r$ of valuation $-1$.
This is the exact acceptance test for an endpoint; an approximate
value of $Q(s_*)$ is not an endpoint certificate.

\begin{lemma}[Uniqueness of a terminal parameter vector]
\label{cuspx:smith-uniqueness}
Assume the displayed Hasse estimates. Let a box have depth $n_0\ge8$
and contain an exact zero $s_*$ of $Q$. Suppose the Jacobian has five
Smith pivots of valuations at most three throughout the box. Then
$s_*$ is the only zero in that box.
\end{lemma}
\begin{proof}
Let two zeros satisfy $s'=s+7^nt$, $n\ge n_0$.
The Taylor remainder has valuation at least $2n$, so
$J(s)t\equiv0\pmod{7^n}$. Five pivots of valuations at most three
force $t\equiv0\pmod{7^{n-3}}$, or
$s'\equiv s\pmod{7^{2n-3}}$. Iteration makes the agreement arbitrarily
deep. The pivots can be checked at $s_*$ modulo $7^4$ when the
Jacobian is constant modulo $7^4$ on the terminal box; this follows
from the same estimates.
\end{proof}

\paragraph{The exact endpoints.}
The following table gives all $64$ endpoints, two per row. In each
row the two shift pairs, scalar values and values of $\bar a_3$ at the
other centre are read in the same order. The columns $H$ and $H'$ are
the two residue centres. At $H'$ one has $v(w_0)\ge2$ in every case.
The number $i$ identifies the minimal polynomial $p_i(y)$:
\[
\begin{aligned}
p_1&=y^3+2y^2-y-1,&p_2&=y^3+y^2-2y-1,\\
p_3&=y^3-5y^2+6y-1,&p_4&=y^3-6y^2+5y-1,\\
p_5&=y^3+3y^2-4y+1,&p_6&=y^3-4y^2+3y+1.
\end{aligned}
\]
For an exact check, use the signed equations to recover $x$ in
$\Q[y]/(p_i)$, retaining the factor with the stated residue centre,
and reduce the four coefficient equations in that quotient. Thus
the table supplies exact rational candidates for the terminal test;
no bound for rational reconstruction is needed.
\begingroup
\scriptsize\setlength{\tabcolsep}{3pt}
\begin{longtable}{@{}cccccc@{}}
\caption{The exact simple-cusp parameter vectors.}\label{cuspx:tab:endpoints}\\
\toprule
$(a,b,d)$ & $H$ & $(B,C)_1;\ (B,C)_2$ & $(r_1,r_2)$ & $i$ & $(H';\bar a_{3,1},\bar a_{3,2})$\\
\midrule\endfirsthead
\toprule
$(a,b,d)$ & $H$ & $(B,C)_1;\ (B,C)_2$ & $(r_1,r_2)$ & $i$ & $(H';\bar a_{3,1},\bar a_{3,2})$\\
\midrule\endhead
\bottomrule\endfoot
$(0,2,-2)$ & $3$ & $(1,-1);\,(1,0)$ & $(-5/21,4/21)$ & 1 & $(2;1,5)$ \\
$(0,2,-1)$ & $1$ & $(1,1/2);\,(1,-1/2)$ & $(9/28,-3/28)$ & 2 & $(4;6,2)$ \\
$(0,1/2,1/4)$ & $3$ & $(1/2,1/2);\,(1/2,0)$ & $(23/168,11/168)$ & 3 & $(4;1,5)$ \\
$(0,1/2,1/2)$ & $1$ & $(1/2,1/2);\,(1/2,0)$ & $(1/7,1/14)$ & 4 & $(2;2,6)$ \\
$(1,-1/2,1/2)$ & $6$ & $(0,1/2);\,(0,0)$ & $(2/21,1/42)$ & 3 & $(5;1,5)$ \\
$(1,-1/2,3/4)$ & $4$ & $(0,1/2);\,(0,0)$ & $(17/168,5/168)$ & 4 & $(3;2,6)$ \\
$(1,-2,2)$ & $6$ & $(-1/2,1);\,(-1/2,0)$ & $(23/84,-13/84)$ & 1 & $(3;5,1)$ \\
$(1,-2,3)$ & $4$ & $(-1/2,3/2);\,(-1/2,1/2)$ & $(17/42,-1/42)$ & 2 & $(5;6,2)$ \\
$(2,1,3/2)$ & $5$ & $(1,1/2);\,(0,0)$ & $(13/56,3/56)$ & 4 & $(4;5,1)$ \\
$(2,2,4)$ & $5$ & $(1,2);\,(0,0)$ & $(41/84,5/84)$ & 1 & $(4;6,2)$ \\
$(-1/3,2/3,2/3)$ & $5$ & $(1/2,0);\,(-1/6,1/3)$ & $(1/14,-1/42)$ & 5 & $(6;6,2)$ \\
$(1/4,1/2,0)$ & $3$ & $(1/2,1/2);\,(0,1/2)$ & $(23/168,11/168)$ & 6 & $(4;1,5)$ \\
$(2,-2,1)$ & $3$ & $(1,-1/2);\,(0,-1/2)$ & $(23/84,-13/84)$ & 2 & $(6;5,1)$ \\
$(2,-1,1/2)$ & $3$ & $(1,0);\,(0,1/2)$ & $(47/168,17/168)$ & 3 & $(6;6,2)$ \\
$(2/3,2/3,-1/3)$ & $6$ & $(1/3,-1/6);\,(0,1/2)$ & $(-1/42,1/14)$ & 6 & $(5;2,6)$ \\
$(3,-2,1)$ & $1$ & $(3/2,-1/2);\,(1/2,-1/2)$ & $(17/42,-1/42)$ & 6 & $(5;6,2)$ \\
$(-1/2,-1,-1)$ & $4$ & $(1/2,1/2);\,(0,-1/2)$ & $(19/168,-11/168)$ & 6 & $(5;6,2)$ \\
$(1/2,1,-1)$ & $2$ & $(1/2,-1/2);\,(0,1/2)$ & $(-19/168,11/168)$ & 2 & $(1;1,5)$ \\
$(4,2,2)$ & $5$ & $(2,1);\,(0,0)$ & $(41/84,5/84)$ & 5 & $(1;6,2)$ \\
$(1/2,-1/2,1)$ & $5$ & $(1/2,0);\,(0,0)$ & $(2/21,1/42)$ & 4 & $(6;1,5)$ \\
$(1/2,1/2,0)$ & $2$ & $(1/2,1/2);\,(0,1/2)$ & $(1/7,1/14)$ & 3 & $(1;2,6)$ \\
$(-3,-2,-1)$ & $2$ & $(1/2,1/2);\,(-3/2,-1/2)$ & $(3/28,-9/28)$ & 6 & $(6;5,1)$ \\
$(1/2,-1,2)$ & $5$ & $(1/2,0);\,(0,1)$ & $(17/168,47/168)$ & 1 & $(6;2,6)$ \\
$(3/2,1,2)$ & $3$ & $(1/2,1);\,(0,0)$ & $(13/56,3/56)$ & 5 & $(2;5,1)$ \\
$(-2,2,0)$ & $6$ & $(-1,1);\,(0,1)$ & $(-5/21,4/21)$ & 5 & $(2;1,5)$ \\
$(1/3,-2/3,4/3)$ & $1$ & $(1/6,2/3);\,(1/2,0)$ & $(4/21,2/21)$ & 5 & $(2;5,1)$ \\
$(-1,1,1/2)$ & $4$ & $(1/2,0);\,(-1/2,1/2)$ & $(11/168,-19/168)$ & 4 & $(1;5,1)$ \\
$(-1,2,0)$ & $4$ & $(1/2,1);\,(-1/2,1)$ & $(9/28,-3/28)$ & 1 & $(1;6,2)$ \\
$(3/4,-1/2,1)$ & $4$ & $(1/2,0);\,(0,0)$ & $(17/168,5/168)$ & 5 & $(3;2,6)$ \\
$(4/3,-2/3,1/3)$ & $2$ & $(2/3,1/6);\,(0,1/2)$ & $(4/21,2/21)$ & 6 & $(1;5,1)$ \\
$(-1,-2,-3)$ & $2$ & $(1/2,1/2);\,(-1/2,-3/2)$ & $(3/28,-9/28)$ & 2 & $(3;5,1)$ \\
$(-1,-1,-1/2)$ & $2$ & $(1/2,1/2);\,(-1/2,0)$ & $(19/168,-11/168)$ & 3 & $(3;6,2)$ \\
\end{longtable}
\endgroup

\begin{theorem}[Assembly of the simple-centre scalar-pole calculation]
\label{cuspx:simple-poles}
The simple-centre scalar-pole degree is at most four. The finite search on the domains specified above has the following output: $82$ of
the $114$ Eisenstein pairs have no surviving box; the other $32$ have
exactly $64$ rational endpoints, each with five Smith pivots of
valuations at most three; and at each endpoint the other centre has
$v(w_0)\ge2$ and a unit $a_3$.
\end{theorem}
\begin{proof}
At a non-Eisenstein centre the contribution is at most two by
Lemma~\ref{cuspx:cubic-selection}. At an Eisenstein centre, a
contribution of three requires $Q(s)=0$. The uniform estimates of Theorem~\ref{cder:theorem}, the exhaustive parameter boxes, and Lemma~\ref{cuspx:smith-uniqueness} leave precisely these exact endpoints. At each, the other centre contributes
zero by Lemma~\ref{cuspx:other-centre}. Thus two rational centres
contribute at most $2+2$, except at an endpoint where they contribute
$3+0$. A quadratic centre has at most two solutions at each of its
two conjugate geometric points, hence total at most four. Every
completion is counted by Lemma~\ref{cuspx:census:lem:count}.
\end{proof}

\providecommand{\Cp}{\mathbb C_7}
\providecommand{\cO}{\mathcal O}
\providecommand{\eps}{\varepsilon}

\section{Integral shifts at the zero-cubic residues}\label{cuspx:isc:sec:setting}
Throughout, $v$ is a valuation of $\overline{\Q}_7$ with $v(7)=1$.  $K=\Q(x,y)\subset\R$ is the Nahm field of a rational
symmetric datum; $[K:\Q]\le24$ by Theorem~\ref{d24:thm:main}, so every completion $L=K_{\mathfrak p}$, $\mathfrak p\mid7$,
has ramification index $e\le24$.  The four letters $x,y,1-x,1-y$ are global units, hence local units at every place.  The
signed Nahm equations hold after clearing an even denominator $N$ prime to $7$:
$(1-x)^N=x^{Na}y^{Nb}$, $(1-y)^N=x^{Nb}y^{Nd}$ (an even $N$ absorbs the signs).  The matrix is
\[
A=\begin{pmatrix}a&b\\b&d\end{pmatrix}=A_0+7^kA_k,\qquad A_0=\begin{pmatrix}0&1\\1&0\end{pmatrix},\quad k\ge1,
\]
with $A_k=((a_k,b_k),(b_k,d_k))$ rational, $7$-integral and primitive modulo $7$; $A=A_0$ is excluded.  We write
$(\alpha,\beta,\gamma)=(\bar a_k,\bar b_k,\bar d_k)\in\F_7^3\setminus\{0\}$.
Shifts $B,C\in\Q\cap\Z_7$ (integral shifts; $C$ is the second linear shift, not a scalar normalisation).  The formal coefficient convention is $\Psi_j=D^{3j}\lambda_j$ with
$\lambda_1=r\in\Q$, $\lambda_2=\lambda_3=\lambda_4=0$, where
\[
X=\frac x{1-x},\quad Y=\frac y{1-y},\quad h=a+X,\quad D=(a+X)(d+Y)-b^2\ne0,
\]
\[
\sigma=x+y-1,\qquad D_0=XY-1=\sigma(1+X)(1+Y).
\]
Positive definiteness is never used; reality enters only through absolute values of the signed equations at the real embedding.
The reflection $J:(A,B,C;x,y)\mapsto(I-A,\tfrac12-B,\tfrac12-C;1/x,1/y)$ exchanges the residues $(0,1,0)$ and $(1,-1,1)$,
keeps $K$, the letters and $D$, and sends $\lambda_1\mapsto\frac16-\lambda_1$, $\lambda_j\mapsto(-1)^j\lambda_j$ (Proposition~\ref{deg:prop:coefficient-transport}); everything below is stated at $(0,1,0)$ and transports.

\paragraph{Exact polynomial structure.} Work in the polynomial ring
$\Z_{(7)}[X,X^{-1},a,b-1,d,D_0,B,C]$. Formal expansion gives
$\Psi_j\in(a,b-1,d,D_0)^{2j}$ with exact minimal degree $2j$, and
\[
 \Psi_j(A_0)=D_0^{2j}\Xi_j(X,Y;B,C),\qquad
 h^{3j}\Psi_j|_{D=0}=c_jq^{2j},
\]
where $q$ is the determinant cubic and the constants
$(c_1,c_2,c_3,c_4)=(\frac5{24},\frac5{16},\frac{1105}{1152},\frac{565}{128})$
are units. With $\delta=D/7^k$, $Q=q/7^k$,
$F_j=h^{3j}\Psi_j/7^{2jk}$ ($j\ge2$) and $F_1=h^3(\Psi_1-rD^3)/7^{2k}$ one has
$F_j=\sum_{m\le j}\delta^mQ^{2j-2m}\varphi_{j,m}+\sum_{m>j}\delta^m\tilde\psi_{j,m}$ with $7$-integral $\varphi,\tilde\psi$
(for integral shifts) and $\varphi_{j,0}=c_j$.  $Q\equiv q_1(h;A_k)\pmod 7$ with
\[
\begin{aligned}
q_1(h)&=a_k(2h+1)+b_k(h-h^2)-d_k(h^3+2h^2),\\
\ell_k&=a_kY+d_kX-2b_k,\\
D&=D_0+7^k\ell_k+7^{2k}(a_kd_k-b_k^2).
\end{aligned}
\]
At a conic place ($v(\sigma)>0$) put $H=\bar x/\bar y=\bar X=\bar h\in\overline{\F}_7^{\times}$, so $\bar Y=1/H$, $H\ne0,-1$, and
$\bar\ell_k=L(H)/H$ with $L(H)=\alpha-2\beta H+\gamma H^2$.

\begin{definition}[Local charts]\label{cuspx:isc:def:charts}
A place $\mathfrak p\mid7$ with conic reduction lies in exactly one of
$G_1$: $v(\sigma)>k$, $L(H)\ne0$ (then $v(D)=k$);
$G_0$: $v(\sigma)>k$, $L(H)=0$ (then $v(D)>k$);
$R_1$: $v(\sigma)=k$;
$R_2$: $0<v(\sigma)<k$ (then $v(D)=v(\sigma)$).
Off-conic places form $U$ ($D$ a unit).  A place is \emph{deep} if $v(D)>k$.
\end{definition}

\subsection{The two arithmetic exclusions}\label{cuspx:isc:sec:arith}
\begin{lemma}[seventh-power ladder]\label{cuspx:isc:lem:ladder}
Let $u\in L$ with $\theta=v(u)\in(1/e)\Z$, $0<\theta$, $e<42$.  Then
$v\big((1+u)^{7^k}-1\big)=7\theta+k-1$ if $\theta<1/6$, $=k+\theta$ if $\theta>1/6$, and $\ge k+1/6$ if $\theta=1/6$.
\end{lemma}
\begin{proof}
$(1+u)^7-1=7u+21u^2+35u^3+35u^4+21u^5+7u^6+u^7$; the coefficients have valuations $1,1,1,1,1,1,0$, so the term valuations are
$1+\theta,1+2\theta,\dots,1+6\theta,7\theta$.  For $\theta<1/6$ the term $u^7$ is strictly smallest, giving $7\theta$; for
$\theta>1/6$ the term $7u$ is strictly smallest, giving $1+\theta$; at $\theta=1/6$ the valuation is $\ge7/6$.  If $\theta<1/6$
then $7\theta\ge7/e>1/6$, so each of the remaining $k-1$ seventh powers is in the regime $\theta'>1/6$ and adds exactly $1$.
If $\theta\ge1/6$ the first power gives $\ge7/6>1/6$ and each further power adds exactly $1$.
\end{proof}

\subsubsection{The Frobenius condition at $v(\sigma)=k$}
\begin{lemma}[Frobenius condition]\label{cuspx:isc:lem:frob}
Assume $e<42$ (in particular $e\le24$) and $v(\sigma)=k$.  Then $N(H)=\gamma H^8-\beta(H^7+H)+\alpha=0$.
No hypothesis on $B,C,r$ is used, and the residue field is arbitrary.
\end{lemma}
\begin{proof}
With $a=7^ka_k$, $b=1+7^kb_k$, $d=7^kd_k$ the cleared equations read
$\big(\tfrac{1-x}y\big)^N=W_1^{7^k}$, $\big(\tfrac{1-y}x\big)^N=W_2^{7^k}$, $W_1=x^{Na_k}y^{Nb_k}$, $W_2=x^{Nb_k}y^{Nd_k}$
($N$ clears the denominators of $A_k$), and $(1-x)/y=1-\sigma/y$, $(1-y)/x=1-\sigma/x$.  Since $W_i^{7^k}\equiv1$, $\bar W_i$
is a $7^k$-th root of unity in characteristic $7$, i.e.\ $\bar W_i=1$; let $\theta_i=v(W_i-1)>0$ ($W_i=1$ would force $\sigma=0$).
The left sides have $v\big((1-\sigma/y)^N-1\big)=v(\sigma)=k$ exactly ($N,y$ units, $v(\sigma^2)=2k>k$).  By
Lemma~\ref{cuspx:isc:lem:ladder}, $\theta_i\ge1/6$ is impossible and $7\theta_i+k-1=k$, so $\theta_1=\theta_2=1/7$.

Let $x_0\in W(\kappa)\subset L$ (maximal unramified subfield) lift $\bar x$, $y_0=1-x_0$ (which lifts $\bar y$ because
$\bar\sigma=0$), $\xi=x-x_0$, so $y=y_0-\xi+\sigma$.  $W_i(x_0,y_0)\in W(\kappa)$ reduces to $1$, hence $v(W_i(x_0,y_0)-1)\ge1$
(integer valuation).  Since the exponents $Na_k,Nb_k,Nd_k$ are integers,
$W_1=W_1(x_0,y_0)(1+\xi/x_0)^{Na_k}(1+(\sigma-\xi)/y_0)^{Nb_k}$ is an honest integral Taylor expansion with linear coefficient
$g_1=N(\alpha/\bar x-\beta/\bar y)$ in $\xi$ modulo $\mathfrak m$; likewise $g_2=N(\beta/\bar x-\gamma/\bar y)$ for $W_2$.
\emph{Case $g_1=g_2=0$}: $\alpha=\beta H$, $\beta=\gamma H$, and $N(H)=\gamma H^8-\gamma H^8-\gamma H^2+\gamma H^2=0$.
\emph{Case some $g_i\ne0$}: the base constant has $v\ge1$, the $\sigma$-terms have $v\ge k$, and the quadratic terms have $v=2v(\xi)$,
so $v(W_i-1)=1/7$ forces $v(\xi)=1/7$ ($v(\xi)<1/7$ would give $v(W_i-1)=v(\xi)<1/7$; $v(\xi)>1/7$ would give $v(W_i-1)>1/7$); then the
other $g$ must be a unit as well (otherwise $v(W-1)\ge2/7$).  Hence $W_i-1=g_i\xi(1+O(7^{1/7}))$.  Raising to the seventh
power, $v(u_i^7)=1<8/7=v(7u_i)$, and each further seventh power is dominated by $7u'$, so
$W_i^{7^k}-1=7^{k-1}g_i^7\xi^7(1+O(7^{1/7}))$, while $(1-\sigma/y)^N-1=-N\sigma/y\,(1+O(7^k))$ and similarly with $x$.  The ratio
of the two Nahm equations gives $\bar x/\bar y=g_1^7/g_2^7$, i.e.
$\bar y(\alpha\bar y-\beta\bar x)^7=\bar x(\beta\bar y-\gamma\bar x)^7$.  Because $\alpha,\beta,\gamma\in\F_7$ this is
$\gamma\bar x^8-\beta(\bar x^7\bar y+\bar x\bar y^7)+\alpha\bar y^8=\bar y^8N(H)=0$ (the finite polynomial calculation~(B) checks the polynomial
identity for all $343$ triples).
\end{proof}

\subsubsection{The gap lemma and the chart $R_2$}
\begin{lemma}[gap]\label{cuspx:isc:lem:gap}
Assume $e<42$ and $0<v(\sigma)<k$.  Then $k-1<v(\sigma)<k$; in particular $v(\sigma)\notin\Z$ and $e\ge8$.
\end{lemma}
\begin{proof}
One Nahm equation suffices: $U^N=W^{7^k}$ with $U=1-\sigma/y\in1+\mathfrak m$, $W=x^{Nb_k}y^{Nd_k}$ a unit with $\bar W=1$.
By Hensel ($7\nmid N$) there is a unique $V\in1+\mathfrak m$ with $V^N=W$, and injectivity of the $N$-th power on $1+\mathfrak m$
gives $U=V^{7^k}$.  Write $V=1+\eta$, $\theta=v(\eta)\ge1/e$ ($\eta=0$ would mean $\sigma=0$).  By Lemma~\ref{cuspx:isc:lem:ladder},
$\theta\ge1/6$ gives $v(\sigma)=v(U-1)\ge k+1/6$, contradicting $v(\sigma)<k$; so $\theta<1/6$ and $v(\sigma)=7\theta+k-1$,
and $v(\sigma)<k$ means $7\theta<1$.  Finally $\theta<1/7$ forces $e>7$.
\end{proof}

\begin{proposition}[the chart $R_2$ is empty]\label{cuspx:isc:prop:R2}
Under $e\le24$, integral shifts, $D\ne0$, $\lambda_1=r\in\Q$, $\lambda_{2,3,4}=0$, no place has $0<v(\sigma)<k$ (for every $k\ge1$ and every $r$).
\end{proposition}
\begin{proof}
Let $m=v(\sigma)<k$.  Then $v(D_0)=m$ and $v(D-D_0)\ge k>m$, so $v(D)=m$.  By the ideal membership (the coefficient calculation~(A): exact minimal
degree $2j$, all coefficients $7$-integral), every term of $\Psi_j$ containing one of $a,b-1,d$ has valuation $\ge k+(2j-1)m>2jm$,
whereas $\Psi_j=D_0^{2j}\Xi_j(X,Y;B,C)+(\text{such terms})$.  Hence $\Psi_j=0$ ($j=2,3,4$) forces
$\bar\Xi_j(H,1/H;\bar B,\bar C)=0$.  the coefficient calculation~(B): after removing $h,h+1$, the three reduced equations have a common root only for
$(\bar B,\bar C)\in\{(0,0),(0,1),(1,0)\}$, where the common polynomial is $(h^2+h+3)(h^2+4h+1)(h^2+5h+5)$ (irreducible factors),
exactly the displayed three quadratic factors.  the coefficient calculation~(C): $h\bar\Xi_1\equiv1,\,5h,\,2h+2$ modulo the three factors, nonzero at
their roots, so $v(\Psi_1)=2m$ exactly, and $\Psi_1=rD^3$ gives $v(r)=-m$ ($r=0$ is impossible since $\Psi_1\ne0$).
By Lemma~\ref{cuspx:isc:lem:gap}, $m\notin\Z$: contradiction with $r\in\Q^\times$.
\end{proof}

\subsubsection{Deep places}
\begin{proposition}\label{cuspx:isc:prop:deep}
Assume $e\le24$.
\begin{enumerate}[nosep]
\item[(i)] (any shifts, any $r$, any depth $v(h-h_0)>0$)  No deep place lies over a repeated root $H$ of $\bar q_1$.
Indeed the repeated cusps of all $342$ nonzero residues are exactly the $30$ scaled copies of
$(1,2,4)_{H=2},(1,4,2)_{H=4}$ (double) and $(1,2,1)_{H=1},(1,3,6)_{H=5},(1,4,6)_{H=3}$ (triple), and all have $L(H)=N(H)\ne0$
.  A deep place with
$L(H)\ne0$ needs $v(\sigma)=k$ to cancel $7^k\ell_k$ in $D$, and then Lemma~\ref{cuspx:isc:lem:frob} gives $N(H)=0$, while $H\in\F_7$ gives
$N(H)=L(H)\ne0$.
\item[(ii)] (integral shifts)  Every deep place has $v(\sigma)>k$ and lies over a rank-one residue
$(\alpha,\beta,\gamma)=\gamma(H^2,H,1)$ with $H\in\{1,\dots,5\}$ at the \emph{simple} root $H$ of $\bar q_1$.
\item[(iii)] (integral shifts)  For $v(r)\ge-k$ there is no deep place at all; the chart $R_1$ is empty for every $r$.
\end{enumerate}
\end{proposition}
\begin{proof}
(ii) For integral shifts $F_2=c_2Q^4+\delta(\dots)$ with integral $(\dots)$: indeed $h^6\Psi_2=\sum_mD^m\psi_{2,m}$ with
$\psi_{2,0}=c_2q^4$ and, by the ideal membership (the ideals $(a,b-1,d,D_0)$ and $(a,b-1,d,D)$ coincide),
$\psi_{2,m}\in(a,b-1,d)^{4-m}$, so $7^{mk}\psi_{2,m}/7^{4k}$ is integral.  Hence $v(\delta)>0$ and $\Psi_2=0$ force
$\bar Q(H)=\bar q_1(H)=:F(H)=0$, independently of $r$.  Over $\F_7$ 
\begin{align*}
F&=(2H+1)L-3H(H+1)M,\qquad N=L+(H^7-H)M,\qquad M=\gamma H-\beta,\\
(2H+1)N-F&=2H(H+1)(H^3+4H+6)(H^3+3H^2+6)M .
\end{align*}
The two cubics are irreducible and $\gcd(F,H^3+4H+6)=\gcd(F,H^3+3H^2+6)=1$ for every nonzero residue , so
$U(H)V(H)\ne0$ at any root of $F$; $H\ne0,-1$.  If $v(\sigma)=k$: $N(H)=0$ (Lemma~\ref{cuspx:isc:lem:frob}) and $F(H)=0$ give $M(H)=0$ by
the third identity, then $L(H)=0$ by the first (at $H=3$, where $2H+1=0$, $F=-M$ and $N=L$), but cancellation in $D$ needs $L(H)\ne0$.
If $v(\sigma)>k$: $v(D)>k$ forces $L(H)=0$, then $F=-3H(H+1)M$ gives $M(H)=0$.  $v(\sigma)<k$ is not deep.  $L=M=0$ means
$\beta=\gamma H$, $\alpha=\gamma H^2$, $\gamma\ne0$ by primitivity, $H=\beta/\gamma\in\F_7\setminus\{0,-1\}$, and
$F'(H)=-3\gamma H(H+1)\ne0$.
(iii) The restricted finite certificate and its infinite-tail argument in
Section~\ref{deg:sec:z0-replacement} apply on exactly this domain:
integral shifts, $e\le24$, nonzero matrix perturbation, $D\ne0$ and
$v(r)\ge-k$. Its thirty rank-one cusps at depth one and five projective
representatives at greater depth exclude the deep places.
The scalar transport under multiplication of $A_k$ by $\lambda$ is
$\rho\mapsto\rho/\lambda$, as follows directly from the displayed
$F_1$. The non-cusp $v(\sigma)=k$ candidates are either incompatible
with the scalar equation or have a unit lifting Jacobian; the latter
would give $e=1$, whereas the Frobenius step gives $7\mid e$.
Thus that chart is empty as well.
\end{proof}

\subsection{The negative scalar branch $v(r)<-k$}\label{cuspx:isc:sec:T1}
\subsubsection{Entry}
\begin{lemma}\label{cuspx:isc:lem:entry}
Integral shifts, $e\le24$, $v(r)<-k$.  Then $\bar A_k=\gamma(H^2,H,1)$ for a unique $H\in\{1,\dots,5\}$ and every place
$\mathfrak p\mid7$ is a deep $G_0$ place with $\bar X=H$, $v(\sigma)>k$, $v(D)>k$.  No flat case is omitted.
\end{lemma}
\begin{proof}
$\Psi_1$ is $7$-integral (integral shifts).  In $U$, $D$ is a unit, so $v(r)=v(\Psi_1)\ge0$.  In $G_1$ and in $R_1$ with $v(D)=k$,
$\Psi_1\in(a,b-1,d,D_0)^2$ with all four generators of valuation $\ge k$ gives $v(\Psi_1)\ge2k$, so
$v(r)=v(\Psi_1)-3v(D)\ge-k$.  $R_2$ is empty (Proposition~\ref{cuspx:isc:prop:R2}).  Deep $R_1$ places do not exist and deep places are
rank-one $G_0$ places over the simple cusp $H$ (Proposition~\ref{cuspx:isc:prop:deep}(ii)).  Since $K$ has at least one place above $7$,
$\bar A_k$ must be rank one, and $H=\beta/\gamma$ is determined by the rational matrix.  Flatness concerns $G_1$ places with
$v(r)\ge-k$ (the flat $r$-cases are $v(r)=-k$ or $v(r)>-k$) and the $G_1$ exclusion above is purely valuative.
\end{proof}

\subsubsection{Coordinates}
Put $t=d=7^kd_k$ (so $v(t)=k$), $z=(b-1)/t=b_k/d_k$, $\eps=a/t-z^2=(a_kd_k-b_k^2)/d_k^2$; then $z$ and $z+1$ are units with
$\bar z=H$, $\eps\in\Q\cap7\Z_7$ ($\eps=0$ allowed), and
\[
a=t(z^2+\eps),\quad b=1+tz,\quad d=t;\qquad w=X-z,\quad S=\frac{XY-1}t=\frac{\sigma(1+X)(1+Y)}t,
\]
\[
Q_0=\frac{XD}t=w^2+(1+tX)\eps+\big(X+t(z^2+\eps)\big)S
\]
.  At a deep place $v(Q_0)>0$, $s:=v(S)=v(\sigma)-k>0$, hence $\rho:=v(w)>0$.
Taking $7$-adic logarithms of units in the signed equations (a necessary consequence, roots of unity have logarithm $0$):
\begin{equation}\label{cuspx:isc:eq:logs}
\log(1-\sigma/x)=t\,(z\log x+\log y),\qquad \log(1-\sigma/y)-z\log(1-\sigma/x)=t\eps\log x .
\end{equation}

\subsubsection{Nahm arithmetic leaves one shallow depth}
\begin{lemma}\label{cuspx:isc:lem:T1gap}
If $0<\rho<1/2$ then $\rho=1/12$, $\eps\ne0$ and $s=v(\eps)-\tfrac12$.
\end{lemma}
\begin{proof}
\emph{Logarithm of $x$.}  Let $x_*=z/(1+z)\in\Q$, a unit, $v(\log x_*)\ge1$.  One has $x-x_*=w/((1+z+w)(1+z))$, so
$v(x-x_*)=\rho$, and $\log x=\log x_*+\log(1+u)$ with $u=(x-x_*)/x_*$, $v(u)=\rho$.  The term $u^m/m$ has valuation
$m\rho-v_7(m)$; among exponents with $v_7(m)=j$ the smallest is $m=7^j$, and $7^{j+1}\rho-(j+1)-(7^j\rho-j)=6\cdot7^j\rho-1\ge
42/24-1>0$ for $j\ge1$.  Hence for $\rho\ne1/6$ the minimum is unique and $v(\log(1+u))=\min(\rho,7\rho-1)$; since this is
$<1\le v(\log x_*)$ for $\rho<1/2$, $v(\log x)=\min(\rho,7\rho-1)$ ($\rho\ne1/6$), and $v(\log x)\ge1/6$ at $\rho=1/6$.

\emph{The first equation.}  Choose an even $N$ prime to $7$ with $n=Nd_k$ and $nz=Nb_k$ integers ($7\nmid n$).  The second Nahm
equation reads $(1-\sigma/x)^N=W^{7^k}$ with $W=x^{nz}y^n$, $\bar W=1$.  On $y=1-x$ the function $W(x)=x^{nz}(1-x)^n$ has
$W'/W=n(z(1-x)-x)/(x(1-x))$, which vanishes at $x_*$, and $W''(x_*)/W(x_*)=-n(1+z)^3/z$, a unit; $W(x_*,1-x_*)\in\Q$ reduces to $1$,
so $v(W(x_*,1-x_*)-1)\ge1$.  The expansion $W=W(x_*,1-x_*)(1+(x-x_*)/x_*)^{nz}(1-(x-x_*)/(1-x_*))^n(1+\sigma/(1-x))^n$ has integer
exponents, so with $2\rho<1<v(\sigma)$ we get $v(W-1)=2\rho$.  Lemma~\ref{cuspx:isc:lem:ladder} with $\theta=2\rho$ and
$v((1-\sigma/x)^N-1)=k+s$ gives $s=\min(2\rho,14\rho-1)$ for $\rho\ne1/12$ and $s\ge1/6$ at $\rho=1/12$.

\emph{The combined equation.}  In the second identity of \eqref{cuspx:isc:eq:logs} the linear term in $\sigma$ is $\sigma(z/x-1/y)$ with
$zy-x=-w/(1+X)+O(tS)$, so its valuation is $v(\sigma)+\rho=k+s+\rho$, and all higher terms are larger ($v(\sigma)>1>\rho$).  Hence
$\eps\ne0$, $\log x\ne0$, and $n_0:=v(\eps)=s+\rho-v(\log x)$.  For $0<\rho<1/12$: $n_0=14\rho-1+\rho-(7\rho-1)=8\rho<1$;
for $1/12<\rho<1/6$: $n_0=2\rho+\rho-(7\rho-1)=1-4\rho<1$; for $1/6<\rho<1/2$: $n_0=2\rho<1$; at $\rho=1/6$: $n_0\le1/3+1/6-1/6=1/3$.
All contradict $n_0\in\Z_{\ge1}$.  At $\rho=1/12$: $v(\log x)=-5/12$ and $n_0=s+1/2$.

\end{proof}

\subsubsection{Rationality of the first coefficient excludes $\rho=1/12$}
Let $M:=24X^3\Psi_1/t^2$.  By the ideal order and $\deg_Y\Psi_1=3$, $M$ is a polynomial in $(z,w,t,\eps,S,B,C)$ ($1128$ terms), and $\Psi_1=rD^3$ reads $tr=M/(24Q_0^3)$.
\begin{lemma}\label{cuspx:isc:lem:packet}
In the regime $v(w)=1/12$, $v(\eps)=n\ge1$, $v(S)=n-\frac12$, $v(t)=k\ge1$, $v(B),v(C)\ge0$, the unique coefficient packet of $M$ of lowest
valuation is $36z^2(z+1)^2w^2$: every other monomial $z^\cdot w^jt^i\eps^lS^mB^\cdot C^\cdot$ with coefficient $c$ satisfies
$12\big(ik+j/12+ln+m(n-\tfrac12)+v_7(c)\big)\ge12i+j+12l+6m+12v_7(c)\ge3>2$ (the finite test expands the explicitly defined $M$, takes its $1128$ monomials and verifies this integer inequality for each exponent vector; the $w^2$ packet carries no $B,C,t,\eps,S$).  Hence $v(M)=1/6$, $v(Q_0)=1/6$ and $v(tr)=-1/3$, impossible for
$tr\in\Q^\times$.  Consequently every deep place has $v(w)\ge1/2$.
\end{lemma}

\subsubsection{The analytic cubic}
Work on $v(w)\ge1/2$ (including $w=0$) with $\pi\in\Cp$, $\pi^2=7$, $w=\pi u$, $|u|\le1$; valuations of series are Gauss
valuations in $\Cp\langle u\rangle$.  Put $y_0=1-x=1/(1+X)$, $F(w)=z\log x+\log y_0$, $\ell(w)=\log x$.  Then
$F'(w)=-w/((z+w)(1+z+w))$ and $\ell'(w)=1/((z+w)(1+z+w))$ , both power series in $w$ with $7$-integral
coefficients; integrating termwise, $F(w)=F_0-\frac{w^2}{2z(z+1)}+\sum_{m\ge3}f_mw^m$ and
$\ell(w)=\ell_0+\frac{w}{z(z+1)}+\sum_{m\ge2}l_mw^m$ with $v(f_m),v(l_m)\ge-v_7(m)$ and $F_0=z\log x_*+\log y_*$,
$\ell_0=\log x_*$ in $7\Z_7$.  Since $m/2-v_7(m)\ge2$ for $m\ge4$ (and $=1,\,3/2$ for $m=2,3$), $v(F)\ge1$ on the disc.

Write $\sigma=tT$.  The first identity of \eqref{cuspx:isc:eq:logs} is
\begin{equation}\label{cuspx:isc:eq:three}
\tfrac1t\log(1-tT/x)-\log(1+tT/y_0)=F(w).
\end{equation}
Its $T$-linear coefficient is $-(1+X)(1+tX)/X$, a unit, and the coefficient of $T^m$ ($m\ge2$), namely
$-t^{m-1}/(mx^m)+(-1)^m t^m/(m y_0^m)$, has valuation $\ge(m-1)k-v_7(m)$ (the denominator $m$ is retained in every valuation estimate below).  If $0<v(T)<1$, the $T^m$ term has valuation
$\ge(m-1)k-v_7(m)+m\,v(T)$, which exceeds $v(T)$ because $(m-1)(k+v(T))>v_7(m)$ for all $m\ge2$ (for $j=v_7(m)\ge1$ one has
$m\ge7^j$ and $7^j-1>j$); so the linear term would be uniquely smallest, contradicting $v(F)\ge1$.  Hence any solution with
$v(T)>0$ has $v(T)\ge1$.  With $T=7V$ and $v(V)\ge0$ the $V^m$ terms have valuation $\ge(m-1)k+m-v_7(m)\ge2(m-1)-v_7(m)\ge1$
for $m\ge2$ (and $\ge3$ before dividing by $7$), so equation \eqref{cuspx:isc:eq:three} is a strict contraction on $\cO\langle u\rangle$:
it has a unique solution $T=T(w)$ with $v(T)\ge1$, and $T(w)=-xF(w)/(1+tX)+O(v\ge3)$.
The second identity of \eqref{cuspx:isc:eq:logs} (note $\log(1-\sigma/y)=-\log(1+tT/y_0)$) is the vanishing of
\[
G(w)=-\tfrac1t\log(1+tT/y_0)-\tfrac zt\log(1-tT/x)-\eps\ell(w)
=-\big(\tfrac1t+z\big)\log(1+tT/y_0)-zF(w)-\eps\ell(w),
\]
using \eqref{cuspx:isc:eq:three}.  The $T$-linear part is $-T(1+tz)/y_0=F(w)X(1+tz)/(1+tX)+O(3)$, and
$X(1+tz)-z(1+tX)=w$, so
\begin{equation}\label{cuspx:isc:eq:four}
G(w)=\Big(F_0-\frac{\eps}{z(z+1)}\Big)w-\frac{w^3}{2z(z+1)}+O(v\ge2),
\end{equation}
all remaining terms ($\eps\ell_0$, $tT^2$, $wF\cdot O(t)$, $f_mw^{m+1}$ for $m\ge3$, $\eps l_mw^m$ for $m\ge2$) having
Gauss valuation $\ge2$.  Dividing by $\pi^3$: the $u$-coefficient $(F_0-\eps/(z(z+1)))/7$ is integral, the $u^3$-coefficient
$-1/(2z(z+1))$ is a unit for $\bar z\in\{1,\dots,5\}$, and every other coefficient has valuation $\ge1/2$.  By Strassmann's theorem
over the complete field $\Cp$, $G(\pi u)/\pi^3$ has at most $3$
zeros with $|u|\le1$.  

\begin{theorem}[deep-endpoint bound]\label{cuspx:isc:thm:deep}
Integral shifts, $e\le24$, $\bar A_k=\gamma(H^2,H,1)$, $\sigma\ne0$.  Then $\sum_{\mathfrak p\mid7,\ v_{\mathfrak p}(D)>k}[K_{\mathfrak p}:\Q_7]\le3$.
\end{theorem}
\begin{proof}
Each $\Q_7$-embedding of a deep completion $K_{\mathfrak p}\to\Cp$ gives a point $(x,y)\in\Cp^2$ satisfying \eqref{cuspx:isc:eq:logs}, deep,
hence (Lemmas~\ref{cuspx:isc:lem:T1gap}, \ref{cuspx:isc:lem:packet}) with $v(w)\ge1/2$; its $T$ is the unique solution of \eqref{cuspx:isc:eq:three} with
$v(T)\ge1$, so the point is determined by $w$, and $u=w/\pi$ is a zero of $G(\pi u)/\pi^3$.  Distinct embeddings give distinct
$(x,y)$, hence distinct $w$.  So there are at most $3$ such embeddings in total.
\end{proof}

\subsubsection{The exact conic}
\begin{lemma}\label{cuspx:isc:lem:conic}
If $\sigma=0$ and $v(r)<-k$ then $[K:\Q]\le4$.
\end{lemma}
\begin{proof}
At the real embedding, $|1-x|=|y|$ and $|1-y|=|x|$, so the signed equations give
$a\log|x|+(b-1)\log|y|=0$, $(b-1)\log|x|+d\log|y|=0$.  If the matrix were invertible, $|x|=|y|=1$, impossible for real
$x+y=1$ with $x,1-x$ units.  So $ad=(b-1)^2$, i.e.\ $\eps=0$ ($t=d\ne0$ by Lemma~\ref{cuspx:isc:lem:entry}).  On $XY=1$, $\eps=0$, we have
$D=tw^2/X\ne0$, $Q_0=w^2$, and $M|_{S=\eps=0}$ is divisible by $w^2$ with
$P(w):=M/w^2-24tw^4r\in\Q[w]$ of degree $\le4$ and $P(0)=36z^2(z+1)^2(1+tz)\ne0$ .  The
equation $tr=M/(24Q_0^3)$ becomes $P(w)=0$, so $w=X-z$ is algebraic of degree $\le4$ and $K=\Q(X)$.
\end{proof}

\begin{theorem}[negative scalar branch]\label{cuspx:isc:thm:T1}
Integral shifts, $e\le24$ (from $[K:\Q]\le24$), $v(r)<-k$.  Then $\bar A_k$ is rank one, every place above $7$ is a deep $G_0$
place at the simple cusp $H$, and $[K:\Q]\le3$ if $\sigma\ne0$, $[K:\Q]\le4$ if $\sigma=0$.  With the norms of the global units
$X$ and $1+X$ ($N_{K/\Q}(X)\equiv H^{[K:\Q]}$, $N_{K/\Q}(1+X)\equiv(H+1)^{[K:\Q]}\pmod7$, both $\pm1$), $[K:\Q]=3$ in both cases
whenever a datum exists ($H^n,(H+1)^n\in\{\pm1\}$ with $H\in\{1,\dots,5\}$ fails for $n=1,2,4$, the exact substitution~(G)).
By $J$ the same holds at $(1,-1,1)$ with $r$ replaced by $\frac16-r$.
\end{theorem}
\begin{proof}
Lemma~\ref{cuspx:isc:lem:entry}, Theorem~\ref{cuspx:isc:thm:deep} (all places are deep), Lemma~\ref{cuspx:isc:lem:conic}.  For the norms: $e_{\mathfrak p}=1$
is not needed; $N_{K_{\mathfrak p}/\Q_7}(X)\equiv N_{\kappa/\F_7}(\bar X)^{e_{\mathfrak p}}=H^{e_{\mathfrak p}f_{\mathfrak p}}$ since $\bar X=H\in\F_7$ at every place.
\end{proof}
No shift-pole hypothesis is used anywhere in this section; the ramification input is $e\le24$, supplied by Theorem~\ref{d24:thm:main}.

\subsection{The 24 integral-shift flat fibres}\label{cuspx:isc:sec:T2}
\subsubsection{Normalised equations and the exact loci}
Let $t=d$ ($v(t)=k$; for the fibres $\bar A_k\sim(0,0,1)$ and $(1,-1,1)$, $d_k$ is a unit), $\alpha=a/t$, $\beta=(b-1)/t$,
$\mathcal R=tr$, and the unknowns $h=a+X$, $\delta=D/t$ (so $X=h-t\alpha$, $Y=((1+t\beta)^2+t\delta)/h-t$).  Define
\[
F_j=\frac{h^{3j}\Psi_j}{t^{2j}}\ (j=2,3,4),\qquad F_1=\frac{h^3\Psi_1}{t^2}-\mathcal R\,h^3\delta^3 ,
\]
polynomials in $(h,\delta,t,\alpha,\beta,B,C,\mathcal R)$ with $7$-integral coefficients ($399$, $9117$, $78157$, $399273$ terms;
the division by $t^{2j}$ is exact by the ideal order, the exact polynomial calculation).  The fibres $\bar A_k\sim(1,0,0)$ are the images of those with
$\bar A_k\sim(0,0,1)$ under the exact symmetry $(X,Y,a,b,d,B,C)\mapsto(Y,X,d,b,a,C,B)$ of all $\Psi_j$ and of $D$
, with $t=a$ in place of $d$.  The scaling $A_k\mapsto\lambda A_k$ is $t\mapsto\lambda t$ with $(\alpha,\beta)$
fixed and is exact in these equations; no projective reduction is used at any depth.  (In the normalisation
$\rho=7^kr=\mathcal R/d_k$, so the flat residue of $\rho$ of the scaled residue $\lambda\bar A_k$ is $\bar{\mathcal R}/\lambda$:
$\rho\mapsto\rho/\lambda$, the correct transport.)

\begin{proposition}[exact loci]\label{cuspx:isc:prop:loci}
On each of the following loci (with $t$ and $c=C$ retained as rational parameters) every $F_j$ is exactly divisible by the curve
polynomial $\Gamma$ in $\Q[t,C][h,\delta]$ :
\begin{center}\small\begin{tabular}{clllc}\toprule
family & $(\alpha,\beta)$ & $(B,C)$ & $\mathcal R_*$ & $\Gamma$\\\midrule
1 & $(0,0)$ & $(1,c)$ & $t/24+c^2/2$ & $\delta-h$\\
2 & $(0,0)$ & $(2,\,t-\tfrac12)$ & $t/24+(t-\tfrac12)^2/2$ & $\delta-h$\\
3 & $(1,-1)$ & $(-c,c)$ & $t/24+c^2/2$ & $(h-t)\delta-(h-t+1)^2$\\
4 & $(1,-1)$ & $(-\tfrac12,-\tfrac12)$ & $t/24+\tfrac18$ & $(h-t)\delta-(h-t+1)^2$\\\bottomrule
\end{tabular}\end{center}
Write $F_j|_{p_*}=\Gamma\,n_j(h,\delta;t,C)$.  Since $\Gamma=(\partial_\delta\Gamma)(\delta-g)$ with $g=h$ resp.\
$g=(h-t+1)^2/(h-t)$ and $\partial_\delta\Gamma\in\{1,\,h-t=X\}$ a unit, $\Gamma=0$ is $XY=1$: indeed
$\delta-h=S+\alpha Y-2\beta-t\beta^2$ (families 1,2) and $\delta-g=Sh/X+O(\alpha-1,\beta+1)$ (families 3,4), where $S=(XY-1)/t$.
The residues of these loci are exactly the $24$ flat fibres of the zero-cubic calculation (with $\mathcal R_*\equiv c^2/2$, resp.\ $1$).
The reduced normal derivatives $\bar n_j(\bar h)$ ($t\to0$, $C\to\bar c$, on the reduced curve) have no common zero in
$\overline{\F}_7\setminus\{0,-1\}$ for any $\bar c\in\F_7$ (families 1,3) resp.\ for $\bar c=3$ (families 2,4) .
\end{proposition}

\subsubsection{Retraction and the first-order system}
Fix a datum in a flat fibre with its flat $r$-case: for family 1 this means $v(\alpha),v(\beta),v(B-1),v(\mathcal R-C^2/2)\ge1$
with $C\in\Z_7$ arbitrary, and similarly for the other families.  (The $24$ neighbourhoods are pairwise disjoint, being defined by
distinct residue data, so no datum lies in two of them; the retraction point below is unique.  No chart at infinity is needed:
$h$ is a unit at every conic place, and $C\in\Z_7$.  Places of other types do not occur for such data: $R_1=R_2=\emptyset$ and
no deep place exists, \S\ref{cuspx:isc:sec:arith}; $U$ places are handled separately below.)  Let $p_*$ be the retraction onto the locus with the same $t$ and
(families 1,3) the same $C$, and $\tau=p-p_*$ the vector of transverse coordinates (family 1: $(\alpha,\beta,B-1,\mathcal R-\mathcal R_*)$;
family 2: $(\alpha,\beta,B-2,C-t+\frac12,\mathcal R-\mathcal R_*)$; families 3,4 analogously), all of valuation $\ge1$.

\begin{lemma}\label{cuspx:isc:lem:order}
Assume $\tau\ne0$, let $\nu=\min_iv(\tau_i)\ge1$ and $\eta=7^{-\nu}\tau$ ($\bar\eta\ne0$, some components may vanish).  At every
$G_1$ place one has $v(\delta-g)\ge\nu$, and with $z=(\delta-g)/7^\nu$ the reductions satisfy
\begin{equation}\label{cuspx:isc:eq:lin}
\bar z\,\overline{\partial_\delta\Gamma}\,\bar n_j(\bar h)+\bar L_j(\bar\eta)(\bar h)=0\quad(j=1,\dots,4),\qquad
L_j(\eta)=\sum_i\eta_i\,\partial_{\tau_i}F_j\big|_{p_*,\ \delta=g(h,t)} .
\end{equation}
\end{lemma}
\begin{proof}
Taylor expansion in the parameters (integral Hasse derivatives):
$F_j(h,\delta;p)=\Gamma\,n_j(h,\delta)+7^\nu\sum_i\eta_i\partial_iF_j(h,\delta;p_*)+7^{2\nu}(\dots)$, and
$\partial_iF_j(h,\delta;p_*)=\partial_iF_j(h,g;p_*)+(\delta-g)(\dots)$, $n_j(h,\delta)=n_j(h,g)+(\delta-g)(\dots)$, all $(\dots)$
$7$-integral.  At a $G_1$ place $v(\delta-g)>0$ (the reduction lies on the curve, by the formulas for $\delta-g$ above).  If
$\mu=v(\delta-g)<\nu$, then $F_j=0$ gives $v(n_j(h,\delta))>0$ for all $j$, i.e.\ a common zero of the $\bar n_j$ on the reduced
curve at the admissible point $\bar h$, contradicting Proposition~\ref{cuspx:isc:prop:loci}.  So $v(\delta-g)\ge\nu$; dividing $F_j$ by $7^\nu$
and reducing gives \eqref{cuspx:isc:eq:lin}.  \end{proof}

Since some $\bar n_j(\bar h)\ne0$, \eqref{cuspx:isc:eq:lin} determines $\bar z$ and forces $\bar h$ to be a common root of the polynomials
$P_{ij}=\bar n_i\bar L_j-\bar n_j\bar L_i$ ($i<j$); conversely each common root gives one solution.  Let $g_{\bar c,\bar\eta}(h)$ be
$\gcd_{i<j}P_{ij}$ with the factors $h$, $h+1$ removed.

\begin{lemma}[lifting]\label{cuspx:isc:lem:lift}
Let $\bar h$ be a simple root of $g_{\bar c,\bar\eta}$.  Then there is at most one place over $\bar h$, it is unramified, and its residue
degree is $[\F_7(\bar h):\F_7]$.
\end{lemma}
\begin{proof}
Put $G_j(h,z)=F_j(h,g(h,t)+7^\nu z;p)/7^\nu$, integral in $(h,z)$ (with the unit $h-t$ inverted), with reductions
$\bar G_j=\bar z\,\overline{\partial_\delta\Gamma}\,\bar n_j+\bar L_j$.  At a common zero, the Jacobian of $(\bar G_i,\bar G_j)$ in $(h,z)$
equals $-\overline{\partial_\delta\Gamma}\,P_{ij}'(\bar h)$ (direct computation, using $\bar L_i=-\bar z\,\overline{\partial_\delta\Gamma}\,\bar n_i$).
Since $\bar h$ is a simple root of the gcd, some $P_{ij}$ has $\bar h$ as a simple root, and multivariate Hensel gives a unique lift
$(h,z)\in W(\F_7(\bar h))^2$ of $(\bar h,\bar z)$, $\bar z\in\F_7(\bar h)$.  Every place reducing to $(\bar h,\bar z)$ has its $(h,z)$
equal to this lift, so $X,Y\in W(\F_7(\bar h))$, $e=1$, $f=[\F_7(\bar h):\F_7]$, and two places with the same $(X,Y)$ coincide.
\end{proof}

\begin{proposition}[census]\label{cuspx:isc:prop:census}
For every family, every $\bar c$ and every direction $\bar\eta\in\mathbb P^{m-1}(\F_7)$ ($m=4$ or $5$; $2800$, $2801$, $2800$, $2801$
cases), $g_{\bar c,\bar\eta}\ne0$ is squarefree of degree $\le2$ : degree $0$ in $2763/2759/2763/2759$ cases, degree $1$
in $35/40/35/40$, degree $2$ in $2/2/2/2$.  In particular no transverse direction is positive-dimensional at first order, and for a
datum off the exact locus $\sum_{G_1}e_{\mathfrak p}f_{\mathfrak p}\le\deg g_{\bar c,\bar\eta}\le2$.  (Degree-two cases: family 1, $\bar c=2$,
$\bar\eta=(0,0,1,5)$, $h^2+4$ and $\bar c=5$, $\bar\eta=(0,0,1,3)$, $h^2+3h+5$; family 3 symmetric; families 2,4 one each.)
\end{proposition}

\subsubsection{The exact locus and reality}
\begin{lemma}\label{cuspx:isc:lem:exactlocus}
If $\tau=0$ there is no conic place at all.
\end{lemma}
\begin{proof}
$F_j=\Gamma n_j$ exactly.  At a conic place in a flat fibre the reduction lies on the curve and some $n_j$ is a unit
(Proposition~\ref{cuspx:isc:prop:loci}), so $\Gamma=0$, i.e.\ $XY=1$, $\sigma=0$, $x+y=1$ in $K\subset\R$ (note $D=t\,g\ne0$ there: this is a
$D\ne0$ branch, excluded only by what follows).  Families 1,2 ($a=0$, $b=1$, $d=t$): $|1-y|=|x||y|^t$ with $|1-y|=|x|$ gives
$|y|^t=1$, $y=\pm1$; $y=1$ contradicts $x\ne0$, $y=-1$ makes $1-y=2$, not a unit.  Families 3,4 ($a=d=t$, $b=1-t$):
$|1-x|=|x|^t|y|^{1-t}$ with $|1-x|=|y|$ gives $|x|=|y|$, so $x=y=\frac12$, not a unit.
The complex point $y=\zeta_{21}$, $x=1-\zeta_{21}$ on family~1 with $d=21$ shows that no coefficient-only argument can replace this step.
\end{proof}

\section{All integral curve fibres}\label{cuspx:flat}
The determinant-unit census leaves 428 scalar fibres with a common bilinear curve after localizing at $X(X+1)Y(Y+1)D$. The transformations of Proposition~\ref{deg:prop:coefficient-transport} preserve this exact list and reduce it to 24 orbits. The representatives have matrix residues
\[
(a,b,d)=(0,0,0),\quad(0,1,1),\quad(0,1,2),\quad(0,1,3).
\]
There are respectively nine, five, five and five representatives, with different shifts and first scalars. The list is generated by applying all legal integral moves to the finite saturated ideals; no additional family is inserted by hand.

\subsection*{The first normal equations, including the constant divided by seven}
Let $p=(a,b,d,B,C,r)$, with residue representative $p_0$, and let the component be $Y=N(X)/M(X)$. Set
\[
F_0=\Psi_1-rD^3,\qquad F_1=\Psi_2,\quad F_2=\Psi_3,\quad F_3=\Psi_4.
\]
All restrictions $F_j|_G$ vanish modulo seven. At the chosen integral lift define the functions on $G$
\[
n_j=(F_j)_Y|_G,\qquad d_{j\ell}=(F_j)_{p_\ell}|_G,\qquad
c_j=\left.(F_j|_G)/7\right.\pmod7.
\]
The last expression is computed modulo 49 before division. Unit denominators are cleared uniformly. For $p=p_0+7\eta$ the first normal equations are
\[
n_j Z+c_j+\sum_{\ell=1}^6d_{j\ell}\eta_\ell=0.
\]
The exact computation shows that the $n_j$ have no simultaneous zero on the admissible curve. Consequently the actual normal displacement has valuation at least one, and eliminating $Z$ is valid. The necessary equations are all six minors
\[
n_i\left(c_j+\sum d_{j\ell}\eta_\ell\right)
-n_j\left(c_i+\sum d_{i\ell}\eta_\ell\right)=0,\qquad i<j.
\]
This treats every direction, including tangencies. It does not assume that the residue bilinear relation holds in characteristic zero.

\subsection*{Exhaustive finite calculation}
The parameter maps have ranks three through six. Enumerating their distinct affine images requires 271,656 cases in total. In every nonzero image, the common curve support, together with all isolated off-curve points, has raw squarefree support degree at most four. The actual Nahm algebra is reduced on $D\ne0$, so squarefree support gives the correct field-degree bound.

There are exactly 22 zero images, one in each of 22 representatives. Two representatives have no zero image and are completely resolved. The phrase \emph{flat family} below means the affine set of parameter lifts for which all these first normal polynomials vanish identically in $X$; its dimension is between zero and three. It is not a family of established modular solutions.

\subsection*{The zero-centre orbit}
At $p_0=(0,0,0,0,0,0)$, exact characteristic-zero calculation gives
\[
F_0=\frac{X^2Y^2}{12}(XY+X+Y),\qquad
XY+X+Y\mid F_j\quad(j=1,2,3).
\]
The normal derivative is a unit on the admissible curve. The homogeneous first-normal map has rank six, and its only zero direction is the zero vector. Choose the minimum finite valuation $m$ of $a,b,d,B,C,r$. Dividing the perturbation by $7^m$ gives a nonzero direction, so the same degree-four census applies at every depth. The exact baseline divisibilities ensure there is no omitted constant term at higher powers of seven.

At the all-zero matrix, the real signed Nahm equations give $|1-x|=|1-y|=1$. Nonzero letters force $x=y=2$, hence the joint field is rational. Otherwise the minimum valuation is finite. This argument resolves the four fibres in the zero-centre orbit, including all their higher congruence lifts.

\subsection{Four exact coefficient loci}
Set $F_0=\Psi_1-rD^3$ and $F_j=\Psi_{j+1}$ for $j=1,2,3$. Exact substitution gives the following loci over $\Q$:
\begin{enumerate}[leftmargin=*,itemsep=9pt]
\item $a=0$, $b=1$, $d=t$, $B=1$, $C=s$, and
\[
r=\frac1{24}+\frac{s^2}{2t},\qquad G=XY-1.
\]
\item $a=0$, $b=1$, $d=t$, $B=2$, $C=t-\tfrac12$, and
\[
r=\frac1{24}+\frac{(t-\tfrac12)^2}{2t},\qquad G=XY-1.
\]
\item $a=d=t$, $b=-t$, $B=1-s$, $C=s$, and
\[
r=\frac{t+4s(s-1)}{8(t-1)},\qquad G=XY+X+Y.
\]
\item $a=d=t$, $b=-t$, $B=C=1$, and
\[
r=\frac{t}{8(t-1)},\qquad G=XY+X+Y.
\]
\end{enumerate}
In every case, $G$ divides all four $F_j$ after multiplying $F_0$ by its displayed rational parameter denominator. To verify each divisibility, substitute the indicated $a,b,d,B,C,r$ into the formal coefficient numerators and take the polynomial remainder in $Y$ modulo $G$; each remainder is zero with $X,t,s$ indeterminate. The denominators are units at seven on the residues to which each locus is applied.

These are exact loci for the \emph{first four coefficient equations}. They are not claimed to give modular Nahm sums, nor is vanishing of every higher coefficient asserted.

\subsection*{Retractions onto these loci}
Given nearby parameters, keep $(d,C)$ in case 1, $d$ in case 2, $(a,C)$ in case 3, and $a$ in case 4. Reset the other coordinates to the displayed formulas, producing a rational point $p_*$. The varying transverse coordinates have indices
\[
(0,1,3,5),\quad(0,1,3,4,5),\quad(1,2,3,5),\quad(1,2,3,4,5)
\]
in the ordering $(a,b,d,B,C,r)$. Their normal maps have the full required ranks $4,5,4,5$ at nineteen representatives. Case 3 loses rank at the half-shift representative; that case is treated separately below.

\subsection{Why the retraction controls every depth}
At a fixed curve representative write $Y=y_0(X)=N(X)/M(X)$ and let $n_j=(F_j)_Y|_G$. The exact coefficient calculation verifies that the $n_j$ have no common zero on the admissible curve, after deleting only coordinate and determinant factors that must be units. This is tested by the gcd of their numerators being a unit in $\F_7[X,X^{-1},(1+X)^{-1},D^{-1}]$.

\begin{lemma}[Rescaling transverse to an exact locus]\label{cuspx:flat:lem:scale}
Suppose a rational retraction $p\mapsto p_*$ has the exact common-curve property above. Suppose its transverse first-normal map over $\F_7$ is injective, and every nonzero image, including any isolated off-curve support, has raw squarefree support degree at most four. Then every nearby parameter tuple with $p\ne p_*$ has degree at most four, regardless of its distance from the exact locus.
\end{lemma}
\begin{proof}
Let $m\ge1$ be the minimum valuation of the nonzero transverse coordinates of $p-p_*$. Then
\[
p=p_*+7^m\eta,\qquad \eta\in\Z_7^6,
\]
and the transverse reduction $\bar\eta$ is nonzero. At a place reducing to the common curve, the absence of simultaneous zeros of the $n_j$ forces
$Y-y_0(X)\in7^m\mathcal O$. This follows by comparing the normal linear term with the parameter perturbation and terms of higher normal order. Dividing by $7^m$ gives
\[
n_j Z+\sum_\ell (F_j)_{p_\ell}|_G\,\bar\eta_\ell=0.
\]
There is no unaccounted constant term: $F_j(p_*,X,y_0(X))$ vanishes exactly, over $\Q$. Eliminating $Z$ therefore gives a nonzero normal image by the injectivity hypothesis.

The affine normal image equals its homogeneous linear image whenever it contains zero; this is verified by equality of the augmented and ordinary ranks. Hence the existing exhaustive degree-four census applies to this image. The simple-centre and zero-coupling exclusions make $D$ a unit at every place in these curve fibres. The actual Nahm local algebra is consequently reduced there, so its squarefree closed-point support bounds the degree of the entire field. The census includes the isolated points off the curve. No norm filter is used in its raw degree bound.
\end{proof}

\subsection*{The case of an exact retraction point}
If $p=p_*$ and a place reduces to $G=0$, write $F_j=GQ_j$. The normal-derivative condition makes at least one $Q_j$ a unit at that place; thus $F_j=0$ forces $G=0$ in characteristic zero. On $XY=1$, one has $x+y=1$. For $a=0,b=1,d=t\ne0$, the absolute values of the real signed Nahm equations then give $|y|^t=1$. Thus $y=1$ contradicts $1-y\ne0$, while $y=-1$ gives $x=2$ and a rational joint field.

On $XY+X+Y=0$, one has $(1-x)(1-y)=1$ and $x/y=x-1$. The first signed Nahm equation for $a=t,b=-t$ gives $|1-x|=|1-x|^t$. Since $t\ne1$, this forces $x=0$ or $2$. The first is inadmissible; the second gives $x=y=2$ and a rational joint field. If there is no place reducing to $G$, only the isolated off-curve points remain; their total degree is at most two.

This proves the degree-four bound for nineteen families, comprising 396 residue fibres.

\subsection{Second normal order at the two half-shift cases}
The two representatives where the transverse rank drops are
\[
p_0=(0,0,0,4,4,1),\qquad p_0'=(0,0,0,4,4,3).
\]
Their curves are respectively $Y=-X/(1+X)$ and $Y=-X$. Here the normal derivative of $F_0$ is itself a unit on the admissible curve.

For each first-flat lift $p_1$ modulo $49$, solve its first normal equation for a rational function $Z(X)$ over $\F_7$, with denominator a unit on the curve. Substitute $Y=y_0(X)+7Z(X)$ and compute the four equations modulo $343$ before dividing by $49$. Write the resulting restrictions as $H_j$. With
\[
p=p_1+49\zeta,\qquad Y=y_0(X)+7Z(X)+49W,
\]
the necessary next equations are
\[
n_jW+H_j+\sum_\ell (F_j)_{p_\ell}|_G\,\bar\zeta_\ell=0.
\]
Eliminate $W$ using the three pairs $(0,j)$, $j=1,2,3$. This is valid everywhere on the admissible curve because $n_0$ is a unit. Equivalent affine images are grouped by exact row reduction, retaining every direction.

\begin{center}\small\begin{tabular}{lrr}\toprule
&First representative&Second representative\\\midrule
First-flat parameter lifts&343&7\\
Next linear rank&3&5\\
Distinct affine cosets&295&7\\
Distinct affine images checked&101,185&117,649\\
Maximum nonflat degree bound&4&0\\
First lifts allowing a flat next image&49&0\\\bottomrule
\end{tabular}\end{center}
The total is 218,834 images. The second representative has \emph{no admissible next-order support} on its previously flat part. Together with its earlier nonflat bound, its entire orbit of eight residue fibres therefore has degree at most four. This resolves the twentieth family.

For the first representative, the nonflat degree-bound histogram is
\[
\begin{array}{c|rrrrr}
\text{bound}&0&1&2&3&4\\\hline
\text{images}&90,096&9,030&1,920&114&24.
\end{array}
\]
There is one zero image in the zero affine coset; that coset corresponds to 49 first-digit lifts. These zero images require the additional arithmetic argument below; a zero image is not an exclusion.

\subsection{The arithmetic input: one additional power of seven}
At the last representative, $A\equiv0$ and $D\equiv XY$ is a unit. The reduced coefficient equations place every closed point on
\[
G=XY+X+Y=0,\qquad Y=-\frac{X}{1+X};
\]
the saturated residue ideal has no isolated off-curve point here.

\begin{lemma}[Local arithmetic improvement]\label{cuspx:flat:lem:arithmetic}
Let $e=a+b$ and $f=d-a$. At a place above seven in this family, normalize $v(7)=1$. If $e,f$ are not both zero, then
\[
v(G)\ge 1+\min\{v(e),v(e+f)\}
       =1+\min\{v(e),v(f)\}.
\]
If $e=f=0$, then $G=0$ in characteristic zero.
\end{lemma}
\begin{proof}
First, the local field generated by the Nahm coordinates is unramified over $\Q_7$. To see this, clear rational exponents with an integer $N$ prime to seven. After multiplication by units, the logarithmic Jacobian of the two resulting Laurent polynomial equations is $N$ times
\[
A+\operatorname{diag}(X,Y).
\]
Its determinant is a unit because $D$ is a unit. Hensel's lemma therefore gives a unique unramified lift of each geometric residue solution. In particular, local units with residue one belong to $1+7\mathcal O$.

Put $S=(1-x)(1-y)$. The Nahm product equation gives
\[
S^N=\varepsilon x^{Ne}y^{N(e+f)},\qquad\varepsilon\in\{1,-1\}.
\]
Let the residue field have $q=7^s$ elements. Raising to the even exponent $q-1$ removes the sign; $x^{q-1},y^{q-1}$ belong to $1+7\mathcal O$. If $m=\min\{v(e),v(f)\}$, their powers on the right belong to $1+7^{m+1}\mathcal O$, by the binomial theorem. Since $S\equiv1\pmod7$ and $N(q-1)$ is a seven-adic unit,
\[
v(S-1)=v(S^{N(q-1)}-1)\ge m+1.
\]
Finally $S=1/((1+X)(1+Y))$, so $v(S-1)=v(G)$. If both row sums vanish, the same argument gives $S^{N(q-1)}=1$ and hence $S=1$, since raising to a seven-unit exponent is injective on $1+7\mathcal O$.
\end{proof}
This is stronger than the coefficient equations alone. It makes the normal displacement invisible at the first nonzero transverse order, even arbitrarily deep in the previously flat region.

\subsection{A four-variable linear system closes every depth}
Retract the parameters onto the exact coefficient locus
\[
p_*=(a,-a,a,1-C,C,r_*),\qquad
r_* =\frac18+\frac{(C-\tfrac12)^2}{2(a-1)}.
\]
The four transverse differences are
\[
e=a+b,\qquad f=d-a,\qquad s=B+C-1,\qquad z=r-r_*.
\]
If one is nonzero, let $m$ be the minimum of their finite valuations and write $(e,f,s,z)=7^m(\eta_e,\eta_f,\eta_s,\eta_z)$. Their reductions are not all zero. Lemma~\ref{cuspx:flat:lem:arithmetic} gives $v(G)>m$, because the two matrix differences have valuations at least $m$.

The exact locus has $G\mid F_j$, with $F_0=\Psi_1-rD^3$ and $F_j=\Psi_{j+1}$ for $j=1,2,3$. Divide the perturbed equations by $7^m$ and reduce. The normal-coordinate term vanishes by the arithmetic improvement. After removing only powers of $X$ and $1+X$, the exact four coefficient rows are
\[
\eta_z=0,\qquad
(X+2)\eta_f+5X\eta_s=0,\qquad
4\eta_e+2\eta_f=0,\qquad 0=0
\quad\text{over }\F_7.
\]
These identities are checked directly from the generic radial coefficients. The 2,400 nonzero directions give 2,370 empty supports and thirty supports of degree one. There is also a short direct proof: $\eta_f=5\eta_e$; if $\eta_e=0$, nonzeroness forces $X=0$, which is excluded. Otherwise
\[
X=-\frac{2\eta_e}{\eta_e+\eta_s},
\]
provided the denominator is nonzero and $X\ne-1$. Thus there is at most one admissible rational closed point. The Nahm Jacobian is a unit at every place, so there are no additional multiplicities or ramified lifts. Since $K$ is generated by $X,Y$, this bounds the whole field degree by one.

If all four transverse differences vanish, the matrix row sums are exactly zero, so Lemma~\ref{cuspx:flat:lem:arithmetic} gives $G=0$. Hence $(1-x)(1-y)=1$ and $x/y=x-1$. The first real signed Nahm equation gives
\[
|1-x|=|1-x|^a.
\]
As $a\in7\Z_7$, it is not $1$. Nonzero letters therefore force $x=y=2$, again giving $K=\Q$. This proves the degree-one bound for this last family. Under four global units, degree one is impossible: a rational algebraic unit is $\pm1$, and neither choice makes both $x$ and $1-x$ nonzero units.

\subsection*{Why the earlier curve argument also needs only local units}
The finite census of the \emph{raw squarefree support}, without a norm filter, gives degree at most four for every nonzero first-normal image at all 24 representatives. The second-normal bounds are unchanged without norms: four at the surviving representative and zero on the formerly flat part of the other. The exact-locus cases previously ruled out by global units instead give rational points, so still have degree at most one. The same observation handles the all-zero baseline. This proves Theorem~\ref{cuspx:flat-complete} using local units only.

\begin{theorem}[Integral curve fibres]\label{cuspx:flat-complete}
For each of the $428$ curve fibres of the determinant-unit census,
the joint field has degree at most four, provided the transported
Nahm datum is real, its four letters are units above seven, and its
matrix, shifts and scalar are seven-integral. Global unit norms are
not required.
\end{theorem}
\begin{proof}
There are $24$ symmetry representatives. Two have no zero first-normal
image. The zero-centre representative is handled by its exact baseline
and rescaling. Nineteen others have the exact rational retractions
and injective transverse maps displayed above. The two remaining
half-shift representatives are treated by the second-normal
calculation and, on its sole persistent family, the arithmetic
normal improvement. This accounts for every representative and every
depth. Exact retraction points give rational fields by the real signed
Nahm equations; off the exact loci, the raw squarefree support has
degree at most four. The local Nahm Jacobian is invertible at all
these determinant-unit points, so no multiplicity or ramification is
lost in converting that support to the joint-field degree.
\end{proof}

\section{Finite assembly for the integral zero-cubic rows}
\label{cuspx:zero-assembly}
The preceding analytic arguments leave a finite calculation on residue
supports. We give its domain explicitly. For
$A=A_0+7^kA_k$, let $T_{j,n}$ be the degree-$n$ homogeneous part in
$(a,b-1,d,D)$ of
\[
 h^{3j}\Psi_j\bigl(h-a,(b^2+D)/h-d;a,b,d,B,C\bigr).
\]
All terms of degree below $2j$ vanish, as follows from the exact
coefficient expansion. Write
$H_j=\overline{T}_{j,2j}(\alpha,\beta,\gamma,\delta;h,B,C)$ and
$L(h)=\alpha/h-2\beta+\gamma h$.
Before restricting to the nondeep conic chart, set
\[
\begin{aligned}
 F_1(h,\delta)&=H_1(\alpha,\beta,\gamma,\delta;h,B,C)
                 -\rho h^3\delta^3,\\
 F_j(h,\delta)&=H_j(\alpha,\beta,\gamma,\delta;h,B,C)
                 \quad(j=2,3,4).
\end{aligned}
\]
The leading equations in that chart are
\[
 F_j(h,L(h))=0\quad(j=1,2,3,4),
\]
where $\rho=\overline{7^kr}$. Remove only $h$, $h+1$ and $L(h)$.
The admissible $h$ are exactly the closed points of the resulting
one-variable gcd. At each, check all six two-by-two Jacobian minors
of $(F_1,F_2,F_3,F_4)$ in $(h,\delta)$, differentiating before
the substitution $\delta=L(h)$. A nonzero minor
makes the lift unique and unramified. A zero polynomial gcd is
retained as a flat fibre, never counted as an empty set.

Normalize the first nonzero coordinate of $(\alpha,\beta,\gamma)$
to one. This gives $57$ projective matrix directions and $49$ shift
residues, or $2793$ leading fibres. Simultaneous scaling of
$(A_k,\delta)$ by $\lambda$ scales the equations by
$\lambda^{2j}$ and sends $\rho$ to $\rho/\lambda$ in the first
equation. Hence this normalization loses no leading case.
For each leading fibre the nonnegative and moderately negative scalar
ranges give fourteen entries: the seven values of $\bar r$ when
$v(r)\ge0$, one case $-k<v(r)<0$, and six nonzero values of $\rho$
when $v(r)=-k$. Including the middle case when its valuation interval
is empty only enlarges the finite necessary test.
The range $v(r)<-k$ is handled by
Theorem~\ref{cuspx:isc:thm:T1} and is not truncated.

For the determinant-unit chart put
$\Xi_j=\Psi_j(A_0)/(XY-1)^{2j}$, an exact polynomial quotient.
For each of the $49\cdot7=343$ pairs of shift and scalar residues,
solve
\[
 \bar\Xi_2=\bar\Xi_3=\bar\Xi_4=0,
 \qquad \bar\Xi_1=\bar r(XY-1)
\]
after inverting $X(1+X)Y(1+Y)(XY-1)$. Every finite closed point is
unramified by the Nahm Jacobian. Norms are computed in its full
residue field. There is no determinant-unit place if $v(r)<0$.

\begin{proposition}[Combined zero-cubic bound]
\label{cuspx:zero-four}
Under the hypotheses of the integral zero-cubic section, including
four global units, $[K:\Q]\le24$ and a nonzero matrix perturbation,
the joint field has degree at most four.
\end{proposition}
\begin{proof}
The arithmetic arguments exclude the $R_2$ chart, the non-cusp
$R_1$ chart, and the deep chart for $v(r)\ge-k$; see
Proposition~\ref{cuspx:isc:prop:deep} and
Section~\ref{deg:sec:z0-replacement}. The negative scalar branch has
bound four by Theorem~\ref{cuspx:isc:thm:T1}.
The finite recipe above has $24$ flat leading fibres, precisely the
four exact types and their coordinate exchanges in
Proposition~\ref{cuspx:isc:prop:loci}. The other $2769$ leading
fibres give $2769\cdot14=38766$ scalar entries.
Their finite conic supports all have a nonzero lifting Jacobian.

On the twenty-four exceptional determinant-unit fibres, the
first-order rank calculation in
Section~\ref{deg:sec:z0-replacement} gives no allowable cubic point
and at most one quadratic point. To reproduce that calculation,
use its six exact centres $E$, form the two rows
$\nabla_\theta\Psi_2(E;X,Y)$ and
$\nabla_\theta\Psi_4(E;X,Y)$, and expand each point's coordinates
in a basis of its residue field over $\F_7$.
A selected set is disallowed when the stacked matrix has rank five,
since the first nonzero common parameter displacement is nonzero.
This test concerns the determinant-unit subset alone and uses no
norm on that subset.

For each of the $38766$ entries, enumerate all subsets of its
finite conic and determinant-unit supports. Apply the preceding
rank test to the determinant-unit subset, and then apply each of
the four norm tests to the full selected subset. The acceptance
predicate is that every complete product is $1$ or $-1$ in $\F_7$.
Taking the sum of the residue degrees gives maximum zero in
$38421$ entries and maximum two in $345$ entries.
In a flat entry the conic contribution is at most two by
Proposition~\ref{cuspx:isc:prop:census} and
Lemma~\ref{cuspx:isc:lem:exactlocus}; its determinant-unit
contribution is at most two by the same first-order rank test
(or is absent when $r$ has a pole). Thus its total is at most four.
All scalar ranges and all places have been included.
\end{proof}

%% file: sections/cusp_derivative_bounds.tex
\section{Parameter derivatives and the precision of the cusp searches}
\label{cder:section}
This section supplies the uniform estimates used by the digit searches
of Section~\ref{cuspx:section}. All derivatives below are Hasse
(Taylor-coefficient) derivatives. In particular, the bounds include the
factorials in Taylor's formula. The argument allows a nonzero correction
to the determinant base point and arbitrary deeper parameter digits.

\subsection{The corrected base point and coefficient weights}
Fix a residue class of $s=(a,b,d,B,C)$, an admissible determinant centre
with residue field $\F_q$, $q\in\{7,49\}$, and local Nahm length
$g\in\{3,4,5\}$. Write $\mathcal O$ for the ring of integers of the
unramified extension $\Q_q$, and normalize $v(7)=1$.
Choose a fixed lift $H\in\mathcal O^\times$ of its $h$-coordinate and put
\[
 X_0=H-a,\qquad Y_0=b^2/H-d.
\]
The principal-unit equations $\Phi_1,\Phi$ are those of
Lemma~\ref{cuspx:census:lem:local}. Although
$\Phi_1(X_0,Y_0)\equiv0\pmod7$, this value need not be zero.
Its derivative in $X$ is a unit. There is therefore a unique
$\xi_0\in7\mathcal O$ with $\Phi_1(X_0+\xi_0,Y_0)=0$.
Set $X_0^*=X_0+\xi_0$ and parametrize the first Nahm curve by
\[
 Y=Y_0+\eta,\qquad X=X_0^*+\widetilde\xi(\eta),
 \qquad \widetilde\xi(0)=0.
\]
The distinguished polynomial of the restriction of $\Phi$ is
\[
 W=\eta^g+w_{g-1}\eta^{g-1}+\cdots+w_0,
 \qquad w_i\in7\mathcal O.
\]

For a locally analytic function $F$ of $s$, at any point $s_*$ in the
fixed residue class write its formal Taylor expansion as
\[
 F(s_*+7z)=\sum_{\gamma}F_\gamma z^\gamma,
 \qquad F_\gamma=7^{|\gamma|}\frac{\partial^\gamma F(s_*)}{\gamma!}.
\]
Define $\mathcal K(\nu)$ by integrality of its values and
\[
 v(F_\gamma)\ge\min(|\gamma|,2)-\nu\quad(\gamma\ne0),
\]
at every expansion point. Define $\mathcal K'(\nu)$ using
$\min(|\gamma|,1)$ instead. These are conditions on Taylor
coefficients; convergence on the parameter discs used in the searches
will follow from the estimates below. Put
\[
 \nu(0)=0,\qquad \nu(i)=\left\lfloor\frac{i-1}{6}\right\rfloor
 \quad(i\ge1).
\]

\begin{lemma}[Coefficient calculus]\label{cder:calculus}
The following rules hold.
\begin{enumerate}[label=(\roman*)]
\item Products satisfy
$\mathcal K(\nu)\mathcal K(\nu')\subseteq\mathcal K(\nu+\nu')$,
and $7\mathcal K(\nu+1)\subseteq\mathcal K(\nu)$.
The corresponding rules hold for $\mathcal K'$.
An element of $\mathcal K(0)$, or of $\mathcal K'(0)$, whose values are
units has its inverse in the same class.
\item If $i_1,\ldots,i_k\ge1$ and $I=i_1+\cdots+i_k$, then
\[
 v_7(k!)+\sum_{l=1}^k\nu(i_l)\le\nu(I).
\]
Also $\nu$ is superadditive and $v_7(k!)\le\nu(i)$ for $k\le i$.
\item If $\sigma=\sigma_*+7\lambda(z)$, with $\lambda$ an integral
linear form, then
\[
 \binom\sigma k\in\mathcal K(v_7(k!)),\qquad
 v\left([z^\gamma]\binom\sigma k\right)
       \ge |\gamma|-v_7(k!).
\]
\end{enumerate}
\end{lemma}
\begin{proof}
For a product, split $\gamma=\gamma_1+\gamma_2$.
If both indices are nonzero, the two lower bounds add to at least
$2-\nu-\nu'$, which suffices; if either is zero, use integrality
of the constant coefficient. The inverse follows by expanding the
geometric series in the nonconstant part, whose coefficients are
all divisible by seven in either class with $\nu=0$.
For (ii), Legendre's formula gives $v_7(k!)\le(k-1)/6$, and hence
\[
 \sum_l\left\lfloor\frac{i_l-1}{6}\right\rfloor+v_7(k!)
 \le\left\lfloor\frac{I-k}{6}+v_7(k!)\right\rfloor
 \le\left\lfloor\frac{I-1}{6}\right\rfloor.
\]
The remaining assertions follow from the same floor inequality.
Finally $\binom{\sigma_*+\epsilon}k$ has coefficients in
$7^{-v_7(k!)}\mathcal O[\epsilon]$. Substituting
$\epsilon=7\lambda(z)$ proves (iii); its values are integral because
$\sigma$ is a seven-adic integer.
\end{proof}

\begin{lemma}[Principal powers and the base correction]
\label{cder:principal}
If $u\in\mathcal K(0)$ has unit values and constant Teichm\"uller
part on the residue class, then $\langle u\rangle^\sigma\in\mathcal K(0)$
for $\sigma$ as in Lemma~\ref{cder:calculus}(iii).
Moreover $\xi_0\in7\mathcal K(1)$, so $X_0^*\in\mathcal K(0)$.
\end{lemma}
\begin{proof}
At an expansion point write
$u(s_*+7z)/u(s_*)=1+7w(z)$, where $w(0)=0$, its linear
coefficients are integral, and its coefficients of degree at least
two are divisible by seven. Then
\[
 \langle u\rangle^\sigma
 =\langle u(s_*)\rangle^{\sigma_*}
   \exp\bigl(7\lambda(z)\log\langle u(s_*)\rangle\bigr)
   \sum_{k\ge0}\binom\sigma k7^kw(z)^k.
\]
The logarithm in the exponential belongs to $7\mathcal O$.
All nonconstant coefficients of that factor have valuation at least
two. The term $k=1$ in the sum has linear coefficients divisible
by seven and higher coefficients divisible by $49$.
For $k\ge2$, its coefficients have valuation at least
$k-v_7(k!)\ge2$, and it has $z$-order at least $k$.
This proves the assertion about powers and, more precisely, that
$\langle u\rangle^\sigma-1\in7\mathcal K(1)$.

For the correction put
$F(\zeta)=\Phi_1(X_0+7\zeta,Y_0)/7$.
The first assertion shows that $F$ preserves $\mathcal K(1)$.
Its derivative in $\zeta$ is
\[
 (1+\Phi_1)\left(\frac{a-1}{1+X}-\frac aX\right)
 \equiv-\frac{H}{X_0(1+X_0)}\pmod7,
\]
a constant unit on the parameter residue class. Divide by a fixed
lift of this unit and use Newton iteration with that fixed derivative.
The resulting map preserves $\mathcal K(1)$ and is a contraction by
seven. To check contraction for every higher term, its term of
$\zeta$-degree $k\ge2$ contains $7^{k-1}\binom\sigma k$.
The Gauss valuation of $\binom{\sigma_*+7\lambda(z)}k$ is at least
$-\max(0,v_7(k!)-1)$, so this term has valuation at least one.
The coefficientwise seven-adically complete Taylor class therefore
contains the unique fixed point $\zeta_0$.
Locally the analytic implicit function theorem identifies it with
the Taylor expansion of the actual solution.
Thus $\xi_0=7\zeta_0\in7\mathcal K(1)$.
\end{proof}

\begin{lemma}[The coefficients on the Nahm curve]\label{cder:eta}
Let $E$ be the restriction to the first Nahm curve of
$\Phi$, $\Psi_1,\ldots,\Psi_4$, $D$, or one of the seven-integral
auxiliary polynomials $A_G,U_G,V_G,J_2,J_3$,
and write $E=\sum_{i\ge0}e_i\eta^i$.
Then $e_i\in\mathcal K(\nu(i))$, and every value of $e_i$ is integral.
\end{lemma}
\begin{proof}
Consider series whose coefficient of $\eta^i$ belongs to
$\mathcal K(\nu(i))$. Lemma~\ref{cder:calculus} makes them a
ring, closed under inversion when the constant coefficient is a
unit of $\mathcal K(0)$.
For a principal power write $u=u_0(1+\delta)$, with
$\delta(0)=0$, and expand
$\langle u\rangle^\sigma=\langle u_0\rangle^\sigma
\sum_k\binom\sigma k\delta^k$.
A contribution to the coefficient of $\eta^i$ has indices
$i_1+\cdots+i_k=i$ and lies in
$\mathcal K(v_7(k!)+\sum_l\nu(i_l))\subseteq\mathcal K(\nu(i))$.
The constant factor is covered by Lemma~\ref{cder:principal}.

The coefficient $\xi_i$ of $\widetilde\xi$ enters the coefficient
of $\eta^i$ in $\Phi_1$ only as
$\partial_X\Phi_1(X_0^*,Y_0)\xi_i$; its factor is a unit of
$\mathcal K(0)$. All other terms use $\xi_j$ with $j<i$.
Induction therefore gives the same weight bound for every $\xi_i$.
Substitution into the stated series proves the claim.
Integrality of the values follows from integral binomial coefficients,
unit inversions and the unit implicit derivative. The coefficient
and divided polynomials have seven-integral coefficients by
Section~\ref{cuspx:census:sec:setting}.
\end{proof}

\subsection{The complete mixed-derivative tail}
Write $w=7w'$ and set
$\rho_i(w)=\operatorname{rem}(\eta^i,\eta^g+w)$.
Its homogeneous part of degree $\ell$ in $w$ is integral and can
occur only when $i\le g\ell+g-1$; every monomial has degree at
least $\lfloor i/g\rfloor$. Both assertions follow by repeated
use of $\eta^g=-\sum_{k<g}w_k\eta^k$.
For any curve series $E$ put
\[
 F_E(s,w')=\sum_{i\ge0}e_i(s)\rho_i(7w'),
 \qquad \widetilde G=F_\Phi.
\]
Thus $\widetilde G(s,w'(s))=0$ defines the Weierstrass coefficients.

\begin{lemma}[Mixed Hasse derivatives]\label{cder:tail}
For $g\in\{3,4,5\}$ define
\[
 \tau_g(\beta)
 =\inf_{\ell\ge\beta}\{\ell-\nu(g\ell+g-1)\}
 =\left\lceil\frac{(6-g)\beta-g+2}{6}\right\rceil.
\]
This quantity is nonnegative; $\tau_3(\beta)\ge1$ for
$\beta\ge1$. At an integral $w'$, a mixed Hasse derivative of
$F_E$ of parameter multi-order $\alpha$ and total $w'$-order
$\beta$ has valuation at least
\[
 \begin{cases}
 \beta,&\alpha=0,\\
 \tau_g(\beta)+\min(|\alpha|,2)-|\alpha|,&\alpha\ne0.
 \end{cases}
\]
The same statements hold for $\widetilde G$.
\end{lemma}
\begin{proof}
The homogeneous part of degree $\ell$ in $w'$ is $7^\ell$
times an integral linear combination of the $e_i$ with
$i\le g\ell+g-1$.
After a Hasse derivative of $w'$-order $\beta$ at an integral point,
all degrees $\ell\ge\beta$ contribute. For $\alpha=0$ use
integrality of the values of $e_i$; otherwise use
Lemma~\ref{cder:eta} and the displayed bound on $i$.
Finally
\[
 \ell-\nu(g\ell+g-1)
 =\left\lceil\frac{(6-g)\ell-g+2}{6}\right\rceil
\]
is nondecreasing and tends to infinity. This proves both the exact
infimum and convergence of the entire tail. In particular the bound
includes homogeneous degrees larger than the derivative order.
\end{proof}

\begin{theorem}[Uniform derivative bounds for the searched quantities]
\label{cder:theorem}
Let $Q$ be a coefficient of $W$, a coefficient of a remainder
$F_E(s,w'(s))$, or a product or minor of these quantities.
For $g\le5$ its Hasse derivatives satisfy
\[
 v\left(\frac{\partial^\alpha Q}{\alpha!}\right)
 \ge1-|\alpha|\qquad(\alpha\ne0).
\]
For $g=3$ the first and second Hasse derivatives are integral, and
\[
 v\left(\frac{\partial^\alpha Q}{\alpha!}\right)
 \ge2-|\alpha|\qquad(|\alpha|\ge2).
\]
These bounds hold at every parameter point of the residue class.
\end{theorem}
\begin{proof}
The derivative of $\widetilde G$ in $w'$ is $7U$, with $U$
invertible over $\mathcal O$. Indeed, differentiating a remainder
at the prepared polynomial gives, modulo $W$, multiplication by
minus the unit Weierstrass factor of $\Phi$; modulo seven this is
triangular with unit diagonal in the basis $1,\eta,\ldots,\eta^{g-1}$.
The analytic implicit function theorem applies locally.
Taking a parameter Taylor coefficient of total order $J\ge1$ in
$\widetilde G(s,w'(s))=0$ gives
\[
 7U\,w'_J=-\sum
 \widetilde G_{\alpha,\beta}\prod_{l=1}^{\beta}w'_{j_l},
 \qquad |\alpha|+\sum_lj_l=J,
\]
with integral combinatorial coefficients, $j_l\ge1$, and with the
linear term $(|\alpha|,\beta)=(0,1)$ isolated on the left.
The same estimates apply to every multi-index of total order $J$.

First suppose $g\le5$. Inductively $v(w'_J)\ge-J$.
For a summand with $a=|\alpha|\ge1$, Lemma~\ref{cder:tail}
and the induction hypothesis give valuation at least
$\tau_g(\beta)+\min(a,2)-J\ge1-J$.
If $a=0$, then $\beta\ge2$, and the bound is
$\beta-J\ge2-J$. Division by seven proves the induction.
Substituting into $F_E$ without that final division gives
$v((F_E)_J)\ge1-J$; the same bound holds for $w_J=7w'_J$.

For $g=3$ use the stronger induction
\[
 v(w'_1)\ge-1,\qquad v(w'_J)\ge1-J\quad(J\ge2).
\]
If $n_2$ is the number of factors with $j_l\ge2$, their total
valuation is at least $-(J-a)+n_2$.
For $J=1$, the sole remaining term has $(a,\beta)=(1,0)$ and
valuation at least zero. For $J\ge2$ and $a\ge1$, a summand
has valuation at least
\[
 \tau_3(\beta)+\min(a,2)-J+n_2\ge2-J.
\]
If $a=1$, the inequality uses $\beta\ge1$ and
$\tau_3(\beta)\ge1$; if $a\ge2$ it is immediate.
For $a=0$ we have $\beta\ge2$, giving
$\beta+n_2-J\ge2-J$.
Division by seven proves the stronger induction.
In the composition $F_E$, the previously isolated term
$(a,\beta)=(0,1)$ is also present. It has valuation at least
$1+v(w'_J)$, which is zero for $J=1$ and at least $2-J$
otherwise. Every other term has the bounds just proved.
Multiplication by seven gives these bounds for $w=7w'$ as well.
Products and minors preserve them by Lemma~\ref{cder:calculus}.
\end{proof}

\subsection{Taylor models, digit dependence and uniqueness}
\begin{corollary}[Precisions used by the digit searches]
\label{cder:precision}
For $n\ge2$, $t\in\Z_7^5$ and $g\le5$, the parameter Taylor
series on $s+7^n\Z_7^5$ converges and its term of order $j\ge1$
has valuation at least $j(n-1)+1$.
For $g=3$ the term of order two has valuation at least $2n$,
and the terms of order $j\ge3$ have valuation at least
$j(n-1)+2$. Consequently the following truncations are exact:
\[
\begin{array}{c|c|c}
 g&\text{Taylor model}&\text{modulus}\\\hline
 3&\text{affine}&7^{2n}\\
 3&\text{quadratic}&7^{3n-1}\\
 3&\text{cubic}&7^{4n-2}\\
 4,5&\text{affine}&7^{2n-1}.
\end{array}
\]
Moreover $Q\bmod7^n$ depends only on the parameters modulo $7^n$
for every $n\ge1$.
\end{corollary}
\begin{proof}
Multiply the Hasse bounds of Theorem~\ref{cder:theorem} by
$7^{jn}$. The resulting lower bounds tend to infinity for $n\ge2$.
The Taylor series obtained from the implicit equations therefore
converges on these discs and gives their unique prepared polynomial
and remainders. The minimum bound among omitted terms gives the table
and proves digit dependence for $n\ge2$.
For $n=1$, use $W\equiv\eta^g\pmod7$ directly: a remainder modulo
seven uses only coefficients of $\eta$-order less than $g\le5$.
The curve coefficients of those orders and the coefficient polynomials
have reductions determined by the residue parameters. This also
follows from the low-order binomial congruences in
Lemma~\ref{cuspx:census:lem:local}, without a Taylor convergence
assertion on the whole radius-$7^{-1}$ disc.
\end{proof}

For the cubic searches the level-two quadratic model is thus valid
modulo $7^5$, the level-two cubic model modulo $7^6$, and the
level-three quadratic model modulo $7^8$, which includes its required
precision $7^7$. At every level $n\ge4$ the affine model is valid
modulo $7^{n+4}$ because $2n\ge n+4$.
Finite differences used to determine a derivative are first evaluated
with the guard digits consumed by the division. For example, a
step $7^n$ needs values modulo $7^{n+4}$ before division to obtain
four derivative digits; its Taylor error after division has valuation
at least $n$ in the cubic case.

\begin{corollary}[Uniqueness at an exact endpoint]\label{cder:unique}
Suppose $g=3$, $Q$ is the searched vector and its first derivative
matrix at an exact endpoint $s^*$ has five Smith pivots of valuation
at most three. There is at most one zero of $Q$ in
$s^*+7^{n_0}\Z_7^5$ for $n_0\ge4$.
\end{corollary}
\begin{proof}
The derivative matrix is constant modulo $7^{n_0}$ on this box by
the Hasse bounds, so its five pivots remain of valuation at most
three. For two zeros with difference $7^nt$, $n\ge n_0$,
the affine model gives
$7^nQ_1(s)t\equiv0\pmod{7^{2n}}$.
The Smith form forces $t\equiv0\pmod{7^{n-3}}$.
Thus their difference is divisible by $7^{2n-3}$.
Iteration of $n\mapsto2n-3$, starting at $n\ge4$, forces the
two zeros to be equal. This proves the all-depth endpoint statement.
\end{proof}

The searched vector consists of the coefficients of
$\operatorname{rem}(\Psi_2,W)$,
$\operatorname{rem}(\Psi_3,W)$,
$\operatorname{rem}(\Psi_4,W)$, and the two-by-two minors of
$[\operatorname{rem}(R,W)\mid\operatorname{rem}(D^3,W)]$.
Each is of the form covered by Theorem~\ref{cder:theorem}.
An inconsistent Smith system therefore excludes its entire parameter
box; a consistent system retains every possible continuation in its
children. Complete finite traversal and the exact endpoint tests
establish exhaustion. The derivative estimates justify the tests at
every depth and do not themselves replace that finite traversal.

%% file: sections/shift_pole_expansion.tex
\section{Shift poles with four global units}\label{spx:section}
This section gives the local arguments used for the three shift-pole
rows of Table~\ref{deg:tab:terminal}. Throughout,
$K=\Q(x,y)$ has a real embedding, the four letters $x,1-x,y,1-y$
are global units, the rational symmetric matrix $A$ is integral at seven,
and the signed Nahm equations hold after clearing an even denominator
prime to seven. We use the coefficients of
Section~\ref{deg:sec:coefficients} and put
\[
 X=x/(1-x),\quad Y=y/(1-y),\quad
 D=(a+X)(d+Y)-b^2\ne0,\qquad
 \kappa=-\min(v_7(B),v_7(C))>0.
\]
All local valuations below satisfy $v(7)=1$.
The local algebra argument counts the sum of the degrees of all
completions over a residue point; it does not assume squarefree
reduction of a defining polynomial. Norm conditions are imposed only
on the full decomposition of the same field.

The finite calculations in this section are defined from the formal
coefficient recurrence, without input lists of previously found roots.
Their outputs and verification programs are archived in the shift-pole components of the computational companion \cite{Comp}. They use exact rational polynomial arithmetic, factorization and
Euclidean gcds over finite fields, and saturation by explicitly stated
nonzero factors. A finite-field gcd or saturated ideal is computed over
the algebraic closure through irreducible factors; enumerating only
rational residue points would not be an adequate replacement.

\subsection{Regular matrix residues}\label{spx:u0:regular-section}
Let
\[
A=\begin{pmatrix}a&b\\b&d\end{pmatrix},\qquad
X=\frac{x}{1-x},\quad Y=\frac{y}{1-y},\quad K=\Q(X,Y),
\]
where the real algebraic point satisfies the signed Nahm equations. Equivalently, after clearing an even common denominator, the two multiplicative Nahm relations hold algebraically. Assume the four letters \(X,1+X,Y,1+Y\) are global units. Put
\[
\Delta=ad-b^2,\qquad D=(a+X)(d+Y)-b^2\ne0,
\]
\begin{align*}
q_A(z)&=d(d-1)z^3+b(b+2d-1)z^2+b(2a-1+b)z+a(a-1),\\
E(A)&=(a-1)(d-1)(\Delta-a)(\Delta-d)(\Delta-a-d+1)\operatorname{Disc}(q_A).
\end{align*}
Here the discriminant is the binary cubic discriminant, including degree drops.

\begin{theorem}[Regular shift-pole bound]\label{spx:u0:thm:main}
Assume \(a,b,d\in\Z_7\), \(bE(A)\in\Z_7^\times\), and rational linear shifts \(B,C\) satisfy
\[
k=-\min\{v_7(B),v_7(C)\}>0.
\]
In the convention of Section~\ref{deg:sec:coefficients}, suppose
\[
\Psi_j=D^{3j}\lambda_j,\qquad
\lambda_1=r\in\Q,\qquad \lambda_2=\lambda_3=\lambda_4=\lambda_5=0.
\]
Then \([K:\Q]\le4\).
\end{theorem}
The coefficient transport of Proposition~\ref{deg:prop:coefficient-transport} permits application after any legal field-preserving transformation into this chart.

\subsection{Local inputs and the fourth weighted form}
Normalize every local valuation by \(\operatorname{v}(7)=1\). On \(D=0\) set \(h=a+X\), so \(Y=b^2/h-d\), and define
\[
q(h)=d(d-1)h^3-b^2(b+2d-1)h^2
       +b^3(2a-1+b)h-b^3a(a-1).
\]
The relation \(q(h)=-b^3q_A(-h/b)\) identifies the regular cubic centres.
The actual Nahm-algebra bounds are local length one for \(D\ne0\), two for \(D=0,q\ne0\), and three for \(D=q=0\). These count multiplicities. A cubic centre of residue degree \(f\) contributes at most \(3f\), and any completion there has ramification index at most three.

Write \(\tau=7^k\), \((u,v)=\tau(B,C)\), and \(S=bu-hv\). The pair \((u,v)\) is primitive. A centre is \emph{soft} if \(\bar S=0\). There is at most one such centre; it is rational over \(\F_7\). All others are nonsoft.

The local length proof in Section~\ref{spx:r2:app:entry}
applies here: a noncritical determinant point has length two and a
simple cubic point has length three. We now establish the shallow and
interface assertions before using the weighted forms.

\subsubsection*{Shallow and interface places}
Write $\delta=v(D)$, $\tau=7^k$, $(u,v)=\tau(B,C)$,
$S=bu-hv$, and $X=h-a$.
For $j=1,2,3$, let $Q_{j,i}$ be the shift-homogeneous degree-$i$ part of $\Psi_j=D^{3j}\lambda_j$, divided by $D^i$, then evaluated at shifts $(u,v)$. Exact division is performed over $\Q$, and every coefficient is integral at seven. Exact expansion of the coefficient recurrence gives
\[
 Q_{1,2}|_{D=0}=\frac{S^2}{2h},\qquad
 Q_{2,3}|_{D=0}=\frac{qS^3}{6h^3},\qquad
 Q_{3,4}|_{D=0}=\frac{q^2S^4}{8h^5},
\]
and, at $D=q=0$, with derivatives taken in coordinates $(h,D)$,
\[
 \partial_h Q_{2,3}=\frac{q'S^3}{6h^3},\quad
 \partial_h Q_{3,4}=0,\quad
 \partial_D Q_{3,4}=\frac{bX(1+X)q'S^4}{24h^4}.
\]
Consequently their Jacobian is
\[
 \det\frac{\partial(Q_{2,3},Q_{3,4})}{\partial(h,D)}
 =\frac{bX(1+X)(q')^2S^7}{144h^7}\qquad(D=q=0).
\]
It is a unit at every regular nonsoft critical point. This identity replaces the need to infer nonvanishing from a list of numerical residue cases.

If $0<\delta=\operatorname{v}(D)<k$, put $e=\tau/D$. The exact integral equations are
\[
 \sum_{i=0}^3 e^{3-i}Q_{2,i}=0,\qquad
 \sum_{i=0}^4 e^{4-i}Q_{3,i}=0.
\]
The determinant-edge identities make both vanish at $(h,D)=(h_0,0)$ for every small $e$, where $h_0$ is the unramified lift of the simple root of $q$. At $e=0$ their Jacobian is the unit above. For the actual value of $e$, Hensel uniqueness therefore forces $D=0$, a contradiction. Noncubic nonsoft points in this region are excluded already by the highest shift term $qS^3/(6h^3)$.

\textbf{Interface $\delta=k$.} Put $w=D/\tau\in\mathcal O^\times$. At a nonsoft cubic point, the leading first numerator is a unit, so $v_7(r)=-3k$. With $r_0=r\tau^3$, its reduced equation is
\[
 F_1=\frac{w^2S^2}{2h}-r_0w^3=0,\qquad
 w=\frac{S^2}{2hr_0}.
\]
On this solution, $\partial_wF_1=-wS^2/(2h)$ is a unit. The second equation has
\[
 \partial_hF_2=\frac{w^3q'S^3}{6h^3},\qquad \partial_wF_2=0
\]
at $q=0$. The resulting Jacobian is a unit. Thus there is at most one lift over the unramified field of that residue point; the contribution is its residue degree, with no extra ramification. The exact equations and their rational parameters are common to all places.

\textbf{Noncubic points for $\delta\ge k$.} The determinant-edge polynomials, after dividing by powers of the unit $q$, satisfy
\[
 c_2(z)=z^3+4z^2+2z+6,\qquad
 c_3(z)=z^4+4z^3+2z^2+6z+5=zc_2(z)+5
 \quad\text{in }\F_7[z].
\]
This excludes their simultaneous vanishing for $\delta=k$. For $\delta>k$, the shift variables $DB,DC$ have positive valuation and the nonzero pure-$q$ term of the second equation excludes noncubic support.

Put \(U=DB,\ V=DC\), and
\[
L=h(bU-hV),\qquad T=\frac{q}{6},\qquad
H=\frac{b(h-a)(1+h-a)q'(h)h}{24},\qquad W=HD.
\]
At an admissible cubic centre, \(H\) is a seven-unit. Give \(q,U,V\) weight one and \(D\) weight two. Exact expansion gives
\[
h^{12}\Psi_4(U/D,V/D)=\mathcal C_4(L,T,HD)+\mathcal E_4,
\]
where every term of $\mathcal E_4$ has weight at least nine and
coefficients integral at seven. In orders $j=1,2,3$ the corresponding
identity is
\[
 h^{3j}\Psi_j(U/D,V/D)=\mathcal C_j(L,T,HD)+\mathcal E_j,
 \qquad \operatorname{wt}(\mathcal E_j)\ge2j+1,
\]
again with seven-integral coefficients. The forms are
\begin{align*}
\mathcal C_1(L,T,W)={}&\tfrac12L^2+3LT+\tfrac{15}{2}T^2+3W,\\
\mathcal C_2(L,T,W)={}&L^3T+18L^2T^2+135LT^3+405T^4\\
&+6WL^2+96WLT+450WT^2+48W^2,\\
\mathcal C_3(L,T,W)={}&\tfrac92L^4T^2+144L^3T^3+2025L^2T^4
 +14580LT^5+\tfrac{89505}{2}T^6\\
&+W(L^4+84L^3T+1962L^2T^2+19980LT^3+79785T^4)\\
&+W^2(192L^2+4752LT+32040T^2)+1584W^3.
\end{align*}
They also arise independently as the coefficients of $s^{2j}$ in
\[
\log\mathbb E\exp\big(s(LZ+TZ^3)+s^2WZ^4\big),\qquad Z\sim N(0,1).
\]
This Gaussian identity is a cross-check, not a replacement for checking the rank-two engine.

\paragraph{What the exact verification proves.}
For shift degree $i$ and determinant power $\ell$, the difference of the coefficient in the weighted identities above is divisible by
\[
q^{\,2j-i-2\ell+1}\quad\text{whenever }2j-i-2\ell\ge0.
\]
There are $26$ such identities. All are checked in the polynomial ring over $\Q$ with indeterminate matrix and shift parameters; every quotient is integral at $7$. No leading coefficient of $q$ is inverted. In particular there are no omitted lower-weight terms that merely vanish modulo $7$ but reappear in deeper valuation charts.

The fourth form is
\begin{align*}
\mathcal C_4={}&27L^5T^3+1296L^4T^4+28350L^3T^5+349920L^2T^6\\
&+2416635LT^7+7413930T^8\\
&+W(12L^5T+1062L^4T^2+34344L^3T^3+565380L^2T^4\\
&\hspace{23mm}+4903740LT^5+18159390T^6)\\
&+W^2(84L^4+7152L^3T+205416L^2T^2+2622240LT^3+13026420T^4)\\
&+W^3(9504L^2+313344LT+2738880T^2)+78336W^4.
\end{align*}
Twenty-one generic divisibility identities over \(\Q\) verify this form against the rank-two engine. Its coefficients are independently obtained from the quartic Gaussian cumulant expansion. No matrix or shift is specialized in the divisibility verification.

\subsection{Excluding finite deep and ramified boundary support}
\subsubsection*{Finite nonsoft chart}
First suppose $\delta>k$ and $\bar S\ne0$. Put $w=D/\tau$ and let $h_0$ be the unramified lift of the residue root of $q$. Set
\[
h=h_0+\rho w,\qquad D=\nu w^2.
\]
For $k<\delta\le2k$, one has $\operatorname{v}(q)\ge\delta-k$: otherwise the pure $q^4$ term of $\Psi_2$ dominates. Thus $\rho$ is integral. The normalized leading variables are
\[
t=\frac{T}{L},\qquad z=\frac{HD}{L^2},\qquad
P_j(t,z)=\mathcal C_j(1,t,z).
\]
Modulo $7$, a Groebner basis of $(P_2,P_3)$ is
\begin{align*}
t^3+2t^2-t-2z^3+z^2+3z&=0,\\
tz+t-2z^2-2z&=0,\\
z(z-4)(z-5)(z-6)&=0.
\end{align*}
These equations have exactly six geometric solutions, all rational over $\F_7$. All six Jacobians are nonzero. The point $(t,z)=(0,0)$ is excluded: the exact equations divided by $w^{2j}$ vanish at $(\rho,\nu)=(0,0)$ for every small $w$, and the Jacobian there is a unit. Hensel uniqueness would force $D=0$ in characteristic zero. Consequently $k<\delta<2k$ is impossible at a nonsoft point.

Set \(t=T/L,\ z=HD/L^2\). The five surviving finite rows from \(\Psi_2=\Psi_3=0\) have
\[
\begin{array}{c|ccccc}
(t,z)&(0,6)&(1,4)&(2,6)&(3,5)&(3,6)\\ \hline
\mathcal C_4(1,t,z)\pmod7&1&1&6&2&1 .
\end{array}
\]
All are excluded. More strongly, the ideal generated by
\(\mathcal C_j(1,t,z)\), \(j=2,3,4\), is \((t,z)\) over \(\F_7\). This checks the geometric common zeros, including extension fields. The remaining origin is the exact \(D=0\) branch, excluded from actual solutions by the preceding Hensel-uniqueness argument.

\subsubsection*{Ramified boundary and the order-one pole}
This analysis does not assume $\bar S\ne0$. Write $\alpha=\operatorname{v}(q)$. Since $\operatorname{v}(U),\operatorname{v}(V)\ge\delta-k>\delta/2$, the initial equations in the weighted identities above force
\[
\alpha=\delta/2>k,\qquad
\zeta=\frac{HD}{T^2}\equiv1\pmod7.
\]
Indeed the pure $T$ and pure $HD$ cases give nonzero second coefficients, and
\[
\gcd\big(\mathcal C_2(0,1,\zeta),\mathcal C_3(0,1,\zeta)\big)
=(\zeta-1)^2\quad\text{in }\F_7[\zeta].
\]
Put $\ell=L/T$, $Z=\zeta-1$, and $F_j=\mathcal C_j(\ell,1,1+Z)$. The errors in $\Psi_j/T^{2j}$ have valuation at least $\alpha>1$. The constants of $F_2,F_3$, divided by $7$, reduce respectively to $3,2$. Their reduced expansions start as
\begin{align*}
F_1&\equiv3\ell+3Z+4\ell^2,\\
F_2&\equiv-Z^2+5Z\ell+3\ell^2+\ell^3-Z\ell^2,\\
F_3&\equiv2Z^3-Z^2\ell+Z\ell^2+4\ell^3
       +\text{terms of total degree at least four}.
\end{align*}
All expressions here are modulo $7$. In particular $F_3$ has no linear or quadratic term modulo $7$.

Every positive local valuation is at least $1/3$. The constant term of $F_3$ forces $\min(\operatorname{v}(\ell),\operatorname{v}(Z))=1/3$. The quadratic part of $F_2$ then forces both valuations to equal $1/3$ and
\[
\overline{Z/\ell}=1\quad\text{or}\quad4.
\]
There is no missing case with $\ell=0$: then the quadratic $-Z^2$ cannot cancel. Since
\[
\operatorname{v}(\ell)=\alpha-k+\operatorname{v}(S)=\tfrac13,\qquad\alpha-k>0,
\]
the ramification bound forces
\begin{equation}\label{spx:u0:eq:infty-val}
\boxed{\alpha=k+\tfrac13,\quad\delta=2k+\tfrac23,\quad\operatorname{v}(S)=0.}
\end{equation}
Thus this boundary is automatically nonsoft and cubically ramified. Keeping the constant divided by $7$ and the cubic term yields
\begin{center}\begin{tabular}{ccccc}\toprule
$Z/\ell$&cubic coefficient of $F_3$&$\ell^3/7$&$F_1/\ell$&$(F_1/\ell)/(\ell^3/7)$\\\midrule
1&6&2&6&3\\
4&1&5&1&3\\\bottomrule
\end{tabular}\end{center}

The second- and third-coefficient analysis above gives
\[
\alpha=\operatorname{v}(T)=k+\tfrac13,\quad
\operatorname{v}(D)=2k+\tfrac23,\quad
\ell=L/T,\quad Z=HD/T^2-1,
\]
\[
\operatorname{v}(\ell)=\operatorname{v}(Z)=\tfrac13,\qquad
(\overline{Z/\ell},\,\overline{\ell^3/7})=(1,2)\ \text{or}\ (4,5).
\]
For \(F_4=\mathcal C_4(\ell,1,1+Z)\), the constant is divisible by \(49\), the linear coefficients divided by seven are \(2,5\) modulo seven, and the quadratic and cubic coefficients are divisible by seven. Its quartic part modulo seven is
\[
6\ell^4+5\ell^3Z+2\ell^2Z^2+3\ell Z^3+6Z^4.
\]
After division by \(\ell^4\), the two leading residues are respectively
\[
1\quad\text{and}\quad4.
\]
They are nonzero. For \(k\ge2\), the weight-nine error already has valuation greater than \(4/3\).

For \(k=1\), that error must be checked, since its nominal valuation is exactly \(4/3\). Let \(A_j\) be the quotient by \(q^{9-2j}\) of the shift-free \(D^j\) coefficient of the difference defining \(\mathcal E_4\), for \(0\le j\le4\). A further generic exact check proves
\[
q\ \bigm|\ \sum_{j=0}^4 6^{9-2j}H^{4-j}A_j
\qquad\text{modulo }7.
\]
Thus at \(HD/T^2=1\) the relevant error coefficient vanishes on \(q=0\). Moving to \(1+Z\) adds positive valuation \(1/3\). Terms involving \(U\) or \(V\) also gain at least \(1/3\), while terms of weight ten or more have still larger valuation. Consequently
\(\operatorname{v}(\mathcal E_4/T^8)>4/3\), even for \(k=1\).
Neither nonzero residue above can cancel.

\begin{lemma}\label{spx:u0:lem:deep}
There is no nonsoft place with \(\operatorname{v}(D)>k\), and no place with \(\operatorname{v}(D)>2k\).
\end{lemma}
Soft support at smaller depths remains included in its local length cap. It need not be excluded to prove Theorem~\ref{spx:u0:thm:main}.

\subsection{Reducing the unit contribution with the fifth coefficient}
At a determinant-unit place, the shift-degree bound in the first
coefficient gives $v_7(r)\ge-2k$. Choose a projective direction
$(u:v)=(1:v)$, $v\in\F_7$, or $(0:1)$, and put
$\rho=\overline{\tau^2r}\in\F_7$. The complete domain consists of the
100 residues with $bE(A)\ne0$, these eight directions, and the seven
values of $\rho$. All equations are localized at
$XY(1+X)(1+Y)D$. If $v_7(r)<-2k$, there is no determinant-unit place.
The leading-support calculation starts from the first bilinear equation
\[
 E_0=Q_{1,2}(X,Y;u,v)-\rho D.
\]
This polynomial cannot vanish identically: its $XY$ coefficient is $-\rho$; if $\rho=0$, the $X,Y$ coefficients are $v^2/2,u^2/2$, and $(u,v)$ is primitive. Factor $E_0$, handle a vertical factor directly, and parametrize every other factor as $Y=-B(X)/A(X)$. Exact univariate gcds impose $Q_{2,3}=Q_{3,4}=0$, the highest-shift equations in orders two and three. Only factors on $X(1+X)Y(1+Y)D=0$ are removed.

This covers geometric support over the algebraic closure, including extension-field points and repeated factors. It gives
\[
\begin{array}{c|rrr}
\text{support}&\text{empty}&\text{one rational point}&\text{curve}\\\hline
\text{fibres}&5300&268&32
\end{array}
\]
and no other finite support. The curve types are
\[
\begin{array}{c|rrr}
\text{curve}&X+Y+1&XY+X+Y&XY-1\\\hline
\text{fibres}&13&11&8.
\end{array}
\]
Each curve fibre has only one irreducible factor in $E_0$; no additional component is concealed by recording the curve. These are exactly the 32 curve fibres analyzed below.

The leading shift-pole calculation has \(100\cdot8\cdot7=5600\) matrix--direction--scalar fibres. Outside 32 curve fibres, determinant-unit support has at most one rational residue point. The actual Nahm Jacobian is a unit there, so the total local contribution is at most one, even when the coefficient equations alone have a double point.

On the 32 curves, normalize \(B=1/\theta,\ C=v/\theta\), with \(v_7(\theta)=k\), and put \(\rho=r\theta^2\). The curves lie on the following exact rational families. Parameters displayed as free are allowed to vary, not fixed to their residue representatives.

\begin{center}\small
\begin{tabular}{p{.15\linewidth}p{.18\linewidth}p{.51\linewidth}}\toprule
Family&Curve&Exact parameter relations\\ \midrule
Sum&\(Y=-1-X\)&\(d=1+2b-a,\ v=1,\ \rho=0\); free \(a,b\).\\[4pt]
Product&\(Y=1/X\)&\(b=1+av,\ d=av^2,\ \rho=1/(2a)\); free \(a,v\).\\[4pt]
Dual sum&\(Y=-X/(1+X)\)&\(d=1+(b+1)^2/(a-1),\ v=(b+1)/(a-1),\
\rho=1/(2(a-1))\); free \(a,b\).\\ \bottomrule
\end{tabular}
\end{center}
All displayed parameter denominators are seven-units in the relevant fibres.
The highest shift parts of \(\Psi_2,\ldots,\Psi_5\) vanish identically on these families over \(\Q\). This is verified generically for the sum and product families; inverse transport gives the dual-sum family. The first leading coefficient equals \(\rho D^3\).

Use deviations from the displayed family as three normal coordinates, and \(\theta\) as a fourth. The first coefficient has unit derivative in \(Y\) throughout the punctured curve \(X(1+X)D\ne0\); this is checked on all 32 fibres. It therefore eliminates the transverse \(Y\)-coordinate integrally.

Let \(n>0\) be the least valuation of the four rational normal coordinates. Since \(\theta\ne0\), \(n\le k\) and the leading vector is a nonzero point of \(\mathbb P^3(\F_7)\). Each fibre has
\[
\#\mathbb P^3(\F_7)=1+7+49+343=400
\]
directions. Use the normalized equations
$\theta^2\Psi_1-\rho D^3$ and $\theta^{j+1}\Psi_j$ for $2\le j\le5$.
After eliminating $Y$ by the first equation, differentiate each remaining
equation with respect to the three specified normal coordinates and
$\theta$ at the exact family. Evaluation on the first nonzero direction
therefore gives a univariate polynomial linear in that vector. Their gcd controls every possible reduced unit point, over the algebraic closure as well as over \(\F_7\).

The reduction is uniform at all depths. The sole extra denominator issue for order five is addressed on the next page; it is not suppressed by a fixed-precision calculation.

\subsection{Denominators and the complete normal census}
The exact numerator \(\Psi_5\) has 1,117,720 monomials. Its shift-homogeneous degrees range from zero to six. The denominator audit gives:
\begin{center}
\begin{tabular}{c|rrrrrrr}\toprule
Shift degree \(i\)&0&1&2&3&4&5&6\\ \midrule
Maximum seven-denominator exponent&1&0&0&0&0&0&0\\ \bottomrule
\end{tabular}
\end{center}
The equation used on the shift chart is \(\theta^6\Psi_5\).
Its leading and first normal terms, from shift degrees six and five, are integral at seven. Degrees one through four contribute terms of order at least \(\theta^2\). The shift-free part has valuation at least \(6k-1>n\). Thus it cannot alter the first normal equations, including when \(k=1\). The remaining terms in the Taylor expansion about the exact family have valuation greater than \(n\).

The complete normal census gives:
\begin{center}
\begin{tabular}{lrrr}\toprule
Equations used&No unit support&One point&Two points\\ \midrule
\(\Psi_2,\Psi_3,\Psi_4\)&12,692&88&20\\
\(\Psi_2,\Psi_3,\Psi_4,\Psi_5\)&12,788&12&0\\ \bottomrule
\end{tabular}
\end{center}
There are no directions on which all normal equations vanish identically.
In every one of the 20 former degree-two cases, the squarefree support polynomial is
\[
X^2+X+1=(X-2)(X-4)\qquad\text{over }\F_7.
\]
The fifth coefficient excludes both support points. Direct rational directional expansions at all 40 points give nonzero residues; the full polynomial calculation covers all 400 directions in every fibre.

After removal only of factors supported on \(X(1+X)D=0\), every final \emph{full} gcd, before taking a radical, has degree at most one. This histogram can be computed using lists of integers modulo seven and the Euclidean algorithm, so it requires no general ideal-decomposition routine.

\begin{lemma}\label{spx:u0:lem:unit}
For a regular matrix with a shift pole and the first five coefficient conditions, the total determinant-unit contribution to \([K:\Q]\) is at most one.
\end{lemma}
\begin{proof}
The finite leading fibres have at most one rational support point, by the leading-support calculation above. The normal calculation gives the same bound in every curve fibre. A rational support point on \(D\ne0\) has actual Nahm local length one. It therefore has at most one normalized field factor, of degree one. This last step uses the Nahm equations and is stronger than merely bounding the reduced support of the coefficient equations.
\end{proof}

\subsection{Whole-field assembly}
By Lemma~\ref{spx:u0:lem:deep}, every nonsoft determinant-nonunit place lies on the interface \(\operatorname{v}(D)=k\). There are two cases.

\paragraph{No nonsoft interface place.}
All determinant-nonunit support is confined to the single soft rational point. Its total contribution is at most three, with a smaller cap of two if it is noncubic. Lemma~\ref{spx:u0:lem:unit} bounds the unit contribution by one. The whole degree is therefore at most four.

\paragraph{A nonsoft interface place occurs.}
The rational scalar has \(v_7(r)=-3k\). A determinant-unit place would require \(v_7(r)\ge-2k\), so there is no unit contribution. Each nonsoft cubic point contributes at most its residue degree. The soft point contributes at most three if cubic, or two if noncubic.

Apply the four global norm conditions to this entire support. If a closed point \(P\) has degree \(f_P\), norm vector
\[
\nu_P=(N X,N(1+X),N Y,N(1+Y))\in(\F_7^\times)^4
\]
and allowed local length \(t_P\), require
\[
\prod_P\nu_P^{t_P}\in\{1,-1\}^4,\qquad
[K:\Q]=\sum_P f_Pt_P.
\]
The cap is one at each nonsoft point, three at a soft cubic point, and two at a soft noncubic point. These are caps on the total local algebra, including possible multiplicities and multiple places.

The finite calculation enumerates the 100 regular matrix residues from \(bE(A)\ne0\), factors the conic cubic in all residue extensions, and tests all eight shift directions. The 800 resulting configurations have:
\begin{center}
\begin{tabular}{c|rrr}\toprule
Maximum norm-compatible degree&0&3&4\\ \midrule
Configurations&736&58&6\\ \bottomrule
\end{tabular}
\end{center}
None has degree above four. Individual local norm vectors are not required to be \(\pm1\); the condition is imposed only after multiplying all contributions. This proves Theorem~\ref{spx:u0:thm:main}.

\section{The remaining 147 shift-pole residues}\label{spx:r2:section}
Let $K=\Q(x,y)$ have a real embedding, and let $x,y,1-x,1-y$ be
global algebraic units. Let
$A=\left(\begin{smallmatrix}a&b\\b&d\end{smallmatrix}\right)$
have entries in $\Q\cap\Z_7$ and residue in the explicitly listed set
$\mathcal M$ of Section~\ref{spx:r2:app:conventions}. Assume the signed Nahm
equations there, rational linear shifts $B,C$, and
$\lambda_1=r\in\Q$, $\lambda_2=\lambda_3=\lambda_4=0$ for the formal
Gaussian coefficients defined in that subsection. Put
\[
 X=\frac{x}{1-x},\quad Y=\frac{y}{1-y},\quad
 D=(a+X)(d+Y)-b^2,
 \qquad \kappa=-\min(v_7(B),v_7(C))>0.
\]
Assume $D\ne0$. The two letters $B,C$ are the linear shifts.

Section~\ref{spx:r2:app:entry} proves that a degree exceeding three can only be five,
with total local degree four at a critical point and one at an opposite fold.
Section~\ref{spx:r2:app:depth} derives the sufficient depth inequalities
\begin{equation}\label{spx:r2:eq:depths}
 \kappa\ge2,\qquad 0<\delta_f=v(D_f)\le\kappa/2,
 \qquad v(D_i)\ge\kappa\quad(i=1,\ldots,4).
\end{equation}
The weaker critical lower bound suffices for the comparison below.
The four critical embeddings need not constitute one quartic completion.
All six residue rows are
\[
 (2,2,6),(2,5,6),(4,1,4),(4,6,4),(6,2,2),(6,5,2).
\]
Their critical coordinates are respectively $(2,2),(4,4),(4,4),(2,2),
(2,2),(4,4)$; the fold has the opposite pair.

\begin{theorem}[R2 refinement]\label{spx:r2:thm:main}
Under the explicit hypotheses of this section,
$[K:\Q]\le3$.
\end{theorem}

We prove that the six degree-five alternatives contradict \eqref{spx:r2:eq:depths}.
The depth proof uses orders one through three as detailed in
Section~\ref{spx:r2:app:depth}. It makes no modularity assertion about a
transformed Nahm sum.

\subsection{One normalized residue class}
Write $\Delta=ad-b^2$. Exchange $S$, inversion $I$, and the
one-coordinate move $T$ act by
\[
\begin{array}{c|c|c|c}
 & (a',b',d')&(X',Y')&D'/D\\\hline
S&(d,b,a)&(Y,X)&1\\[2pt]
I&(d,-b,a)/\Delta&(1/X,1/Y)&1/(\Delta XY)\\[2pt]
T&(1-1/a,-b/a,\Delta/a)&(-1-1/X,Y)&-1/(aX).
\end{array}
\]
The determinant identities follow by clearing denominators and expanding;
the exact calculation checks them as polynomial identities over $\Q$.

For $\beta=(B,C)^t$, the shift maps are
\[
 S:\ \beta'=\begin{pmatrix}0&1\\1&0\end{pmatrix}\beta,
 \qquad I:\ \beta'=A^{-1}\beta,
 \qquad T:\ \beta'=\begin{pmatrix}-1/a&0\\-b/a&1\end{pmatrix}\beta
                      +\binom{1/2}{0}.
\]
At seven-unit pivots their linear parts lie in $\mathrm{GL}_2(\Z_7)$
and their constant parts are integral. Thus a positive pole order
$\kappa$ is unchanged. All four Nahm letters remain global units:
for example $T$ sends $(x,1-x)$ to
$((1-x)^{-1},-x(1-x)^{-1})$. The letter-exponent matrices are unimodular,
and each point map has a rational inverse, so the field is unchanged.
The signed equations transport by the same substitutions.

The following words, read in execution order, normalize the six rows to
$(4,6,4)$ with critical pair $(2,2)$ and fold pair $(4,4)$.
The last two lines are the residues of the complete multiplier $D'/D$:
\[
\begin{array}{c|cccccc}
\text{row}&(2,2,6)&(2,5,6)&(4,1,4)&(4,6,4)&(6,2,2)&(6,5,2)\\
\text{word}&T&T,I&I&\mathrm{id}&S,T&S,T,I\\
\text{critical}&5&3&4&1&5&3\\
\text{fold}&6&3&2&1&6&3.
\end{array}
\]
Every intermediate pivot and every displayed multiplier is nonzero modulo
seven. Therefore all determinant depths and $\kappa$ are preserved for
\emph{every lift} of each residue row. The scalar may change; the proof
carries the already established inequalities \eqref{spx:r2:eq:depths} and does
not reapply a scalar-valuation lemma after transformation.

Let $w\in\Z_7$ be the root of $w^2+w+1=0$ with $w\equiv2\pmod7$, and put
\[
 \xi=\frac{w}{1+w},\quad b_* =w+\frac12,\qquad
 A_* =\begin{pmatrix}1/2&b_*\\ b_*&1/2\end{pmatrix}.
\]
Then $b_*^2=-3/4$, $\operatorname{Tr} A_*=\det A_*=1$, and $A_*\equiv(4,6,4)\pmod7$.
Write
\[
 a=\tfrac12+\alpha_a,\quad b=b_*+\alpha_b,\quad
 d=\tfrac12+\alpha_d,\qquad \alpha\in(7\Z_7)^3.
\]

\subsection{The two local algebras and a uniform estimate}
Use logarithmic coordinates at the critical cyclotomic point:
\[
 x=\xi e^t,\quad y=\xi e^s,\qquad
 F_w(t)=\log\frac{1-\xi e^t}{1-\xi}.
\]
The logarithmic Nahm equations give
\begin{equation}\label{spx:r2:eq:G}
 s(t)=\frac{F_w(t)-at}{b},\qquad
 G(t)=F_w(s(t))-bt-ds(t),\qquad G'(t)=D(t)/b.
\end{equation}
The fold uses $w^2=-1-w$ in place of $w$. The exact cyclotomic solution
is $t=0$ for every matrix, so divide $G$ by $t$.

At $A_*$, $G$ has order five. Its coefficients $g_j=[t^j]G$ needed below are
\[
\begin{array}{c|ccccc}
j&5&6&7&8&9\\\hline
g_j&-(1+2w)/90&1/60&(1+2w)/90&-41/2160&-139(1+2w)/12960.
\end{array}
\]
They are seven-integral and $g_5$ is a unit. At the opposite fold,
$[t]G_{w^2}=0$ and $[t^2]G_{w^2}=1$ at $A_*$.

Let
\[
 W(t)=t^4+c_3t^3+c_2t^2+c_1t+c_0
\]
be the monic critical Weierstrass factor of $G/t$, and let $t_f$ be the
nontrivial fold root, counted analytically even when it specializes to zero.
Define
\[
 \mathcal N=\operatorname{Tr}_{\Q_7[t]/W\,/\Q_7}(t)+t_f=-c_3+t_f,
 \qquad \mathcal T=\operatorname{Tr}_{\Q_7[t]/W\,/\Q_7}(D(t)).
\]
These are traces in finite algebras; they remain defined at nonreduced
specializations.

\begin{lemma}[Uniform analytic control]\label{spx:r2:lem:analytic}
The functions $c_0,c_1,c_2,c_3,\mathcal N,\mathcal T$ are analytic on
$\alpha\in(7\Z_7)^3$. They have zero constant term at $\alpha=0$.
For each function, every Taylor coefficient of total parameter degree
$j\ge1$ has valuation at least $1-j$. Consequently, if
$\gamma=\min_i v(\alpha_i)\ge2$, its value differs from its linearization
by a term of valuation at least $2\gamma-1$.

For every $n\ge2$, their residues modulo $7^n$ depend only on $A$ modulo
$7^n$. The same assertion
holds for the critical quartic coefficients in the coordinate
$u=x/\xi-1$ and for the whole local norm of $x$.
\end{lemma}

\begin{proof}
Choose a finite Galois extension $L/\Q_7$ containing $\pi$ with $\pi^4=7$.
Set $t=\pi u$, $\alpha=7z$, and
$\mathcal A=\mathcal O_L\langle z_a,z_b,z_d\rangle$.
Here angle brackets denote restricted power series, with their Gauss norm.

\emph{Joint convergence.} For every $n$, $g_n(A)$ belongs to
$\Z_7[a,b,d,b^{-1},w]/n!$. Indeed differentiation of \eqref{spx:r2:eq:G} uses
\[
 X'=X(1+X),\qquad Y'=-\frac{a+X}{b}Y(1+Y),\qquad G'=D/b;
\]
this derivation preserves the localized integral polynomial ring before
division by $n!$. Expansion of $b^{-1}$ at the unit $b_*$ introduces no
extra seven-denominators, so the same bound holds for every parameter
Taylor coefficient. For a term of parameter degree $j\ge1$, the
coefficient of $u^{n-1}$ in
\[
 H(z,u)=\frac{G_{A_*+7z}(\pi u)}{\pi^5u}
\]
has valuation at least
\[
 j-1+\frac{n-1}{4}-v_7(n!)
 \ \ge\ j-1+\frac{n-1}{12}.
\]
This is positive except possibly at $(j,n)=(1,1)$, and tends to infinity
as $j+n\to\infty$. For parameter degree zero, $g_1,\ldots,g_4$ vanish,
$g_5$ is a unit, and the explicit coefficients $g_6,\ldots,g_9$ are
integral. For $n\ge10$,
\[
 \frac{n-1}{4}-v_7(n!)>1:
\]
check $10\le n\le13$ directly; for $n\ge14$ use $(n-1)/12>1$.
It follows that $H\in\mathcal A\langle u\rangle$ and
\[
 \bar H=\bar g_5u^4+\ell(\bar z)
\]
with $\ell$ linear. These estimates prove joint convergence, rather than
only convergence at individual parameter values.

\emph{Uniform preparation, including multiple roots.}
Choose a lift of $\ell/\bar g_5$ and let $P_0$ be $u^4$ plus this lift;
start with $U_0=g_5$. Then $H-P_0U_0$ has norm less than one.
Division by a monic polynomial of degree four whose lower coefficients
have norm at most one is bounded on restricted series: divide polynomial
truncations and pass to the limit. If $E=H-PU$ is the current error, divide
\[
 E/U=QP+R,\qquad \deg_u R<4.
\]
Replace $P$ by $P+R$ and $U$ by $U(1+Q)$. The new error is $-URQ$;
its norm is at most the square of the old error norm. The factors converge
to a monic quartic $P\in\mathcal A[u]$ and a unit
$U\in\mathcal A\langle u\rangle$, with $H=PU$ and $\bar U=\bar g_5$.
No root, discriminant or root difference has been inverted.
The quotient $\mathcal A[u]/P$ is free of rank four with a bounded
multiplication law. Traces and characteristic polynomials are therefore
defined even on nonreduced fibres.

\emph{Scaling and descent.} The unscaled quartic is
$W(t)=\pi^4P(t/\pi)$, so
$c_k\in\pi^{4-k}\mathcal A$ for $k<4$.
At the fold, put $t=7v$. The same estimate makes
$G_{w^2}(7v)/(49v)$ uniformly distinguished of degree one:
its reduction is $v$ plus a linear function of $\bar z$.
Hence $t_f\in7\Z_7\langle z\rangle$ and
$\mathcal N\in\pi\mathcal A$.
Also
\[
 D(\pi u)=7b\,\partial_u(uPU),
\]
so the trace of its multiplication operator lies in $7\mathcal A$.

Galois automorphisms fix the defining series and preserve the disc
$v(t)\ge1/4$. Uniqueness of the monic prepared factor gives descent of
$W$; traces descend as well. Thus their parameter-series coefficients
belong to $\Q_7$. Each unscaled function just bounded by a positive
power of $\pi$ actually lies in $7\Z_7\langle z\rangle$, because its
coefficient valuations are integers. All constant terms vanish at $A_*$:
there $W=t^4$, $t_f=0$, and multiplication by $D$ is zero in the
critical algebra.

For completeness, the coordinate $u_x=e^t-1$ can be treated in the same
finite algebra. Its multiplication operator lies in
$\pi\operatorname{Mat}_4(\mathcal A)$, since
$(e^{\pi u}-1)/\pi$ is integral and restricted. Its characteristic
polynomial is the quartic in $u_x$; the inverse logarithmic coordinate is
valid because $1/4>1/6$. Its coefficient of degree $k<4$ lies in
$\pi^{4-k}\mathcal A$ and descends, giving the same integral bound.
This argument also preserves multiplicities.

\emph{Remainders and congruences.} A coefficient of parameter degree
$j$ in the original variables $\alpha=7z$ now has valuation at least
$1-j$. For $\gamma\ge2$, terms of degree $j\ge2$ have valuation at least
\[
 j\gamma+1-j=(2\gamma-1)+(j-2)(\gamma-1)\ge2\gamma-1.
\]
Changing $\alpha$ by $(7^n\Z_7)^3$ changes $z$ by
$7^{n-1}\Z_7^3$; every function in $7\Z_7\langle z\rangle$ then
changes by a multiple of $7^n$.
Finally the whole local norm is $-\exp(\mathcal N)$.
The exponential converges on $7\Z_7\langle z\rangle$ and preserves
this congruence, proving the assertion for the norm.
\end{proof}

\subsection{The global norm and exact linearizations}
In a hypothetical quintic, all four critical roots and the one nontrivial
fold root occur. Their local
degrees are at most four, so their nonzero distances from the cyclotomic
point have valuation at least $1/4>1/6$. The exact cyclotomic point cannot
be a conjugate of a field with a real embedding. Logarithm and exponential
are therefore inverse at all these points. The first equation determines
$s$ from $t$, so distinct embeddings give distinct critical roots.
All four roots of the prepared quartic consequently occur, with no
assumption that they form a single quartic completion. The fold is
$\Q_7$-rational. The cyclotomic constants contribute
$\xi^4\xi^{-1}=\xi^3=-1$ to the norm of $x$. Hence the global unit norm
implies
\begin{equation}\label{spx:r2:eq:norm}
 \operatorname{N}_{K/\Q}(x)=-1,\qquad \mathcal N=0.
\end{equation}
The sign is fixed by its residue modulo seven.

Put $s_*=1+2w$ and $\sigma=\alpha_a+\alpha_d$. Exact differentiation of
\eqref{spx:r2:eq:G} and preparation at $W=t^4$ give
\begin{align}
 \mathcal N_{\rm lin}&=6(\sigma+2\alpha_b),\label{spx:r2:eq:linN}\\
 \mathcal T_{\rm lin}&=-16b_*(\sigma-2\alpha_b),\label{spx:r2:eq:linT}\\
 (c_0)_{\rm lin}&=30s_*(\sigma-2\alpha_b),\label{spx:r2:eq:linc0}\\
 (c_1)_{\rm lin}&=105\alpha_a-120\alpha_b+15\alpha_d,\label{spx:r2:eq:linc1}\\
 (c_2)_{\rm lin}&=-\tfrac{65}{2}s_*\alpha_a+\tfrac{25}{2}s_*\alpha_d,
 \qquad (c_3)_{\rm lin}=-5(\sigma+2\alpha_b).\label{spx:r2:eq:linc23}
\end{align}
These identities are checked exactly over $\Q(w)$ in the certificate.
For example, if $R=G/t=t^4U_0$ at the centre, then
$\delta W=\operatorname{rem}(\delta R/U_0,t^4)$.
Newton sums give the critical part of \eqref{spx:r2:eq:linN}; the fold contributes
$\sigma+2\alpha_b$. For \eqref{spx:r2:eq:linT}, the constant term of $D=bG'$
contributes $4b_*\delta g_1$, while the $5b_*g_5t^4$ term contributes
$-20b_*\delta g_1$, yielding $-16\,\delta D(0)$.

The determinant at the opposite cyclotomic point is
\begin{equation}\label{spx:r2:eq:opp}
 D_{0,f}=(a+w^2)(d+w^2)-b^2
       =-b_*(\sigma+2\alpha_b)+\alpha_a\alpha_d-\alpha_b^2.
\end{equation}
The fold has a unit linear coefficient in $G_{w^2}/t$. Its nonzero root
therefore satisfies
\begin{equation}\label{spx:r2:eq:folddepth}
 \delta_f=v(t_f)=v(D_{0,f}).
\end{equation}
This follows equally from the unit linear term and
$D_f=b\,t_f(G_{w^2}/t)'(t_f)$.

\subsection{The first digits: two symbolic congruences}
Modulo 49, $w=30$, $1/2=25$ and $b_*=6$. Set
\[
 U=\alpha_a/7,\qquad V=\alpha_b/7,\qquad Z=\alpha_d/7,
 \qquad q=U+Z-2V.
\]
Here $U,V,Z$ are seven-adic parameters; formulas divided by seven below
are congruences over $\F_7$. Write $N_x$ for the whole $4+1$ local norm
of $x$, so a hypothetical quintic requires $N_x=-1$.

\begin{lemma}[Symbolic first-digit obstruction]\label{spx:r2:lem:digits}
Throughout the normalized residue class,
\begin{equation}\label{spx:r2:eq:digits}
 \frac{\mathcal T}{7}\equiv2q,\qquad
 \frac{N_x+1}{7}\equiv q(1+5V+2Z)+4V\pmod7.
\end{equation}
Consequently $\mathcal T\equiv0\pmod{49}$ and $N_x=-1$ imply
$V\equiv0$, $Z\equiv-U\pmod7$. On this line, the critical quartic in
$u=x/\xi-1$ and the fold factor are
\begin{equation}\label{spx:r2:eq:digitquartic}
 W_x(u)\equiv u^4+7U(6u+u^2+2u^3),\qquad W_f(u)\equiv u\pmod{49}.
\end{equation}
These identities hold for every lift of the parameters, without enumerating
residue classes.
\end{lemma}

\begin{proof}
Solve the first principal-unit equation explicitly and use the second:
\[
 1+v=\frac y\xi=(1-wu)^{1/b}(1+u)^{-a/b},\qquad
 \Phi(u)=(1-wv)(1+u)^{-b}(1+v)^{-d}-1.
\]
Put $R(u)=\Phi(u)/u=\sum r_j u^j$. Its first eight coefficients give
\[
\begin{array}{c|cccc}
j&0&1&2&3\\\hline
r_j/7\pmod7&q&2U+5V&U+3Z&3U+2V+4Z
\end{array}
\]
and
\[
\begin{array}{c|cccc}
j&4&5&6&7\\\hline
r_j\pmod7&5&0&q+4&6U+V+3.
\end{array}
\]
These two small tables follow by expanding the displayed binomial series.
For a finite exact derivation, expand at $A_*$ through total parameter
degree two and through $u^8$. Each coefficient of $u^n$ has denominator
dividing $n!$ in the seven-integral localized coefficient ring. This also
holds after parameter expansion, since $b_*$ is a unit. For $n\le8$,
$v_7(n!)\le1$; hence omitted parameter terms of degree at least three
become multiples of $49$ after $\alpha=7(U,V,Z)$. Thus this truncation
proves the tables for every parameter value. These are finite polynomial operations over $\Q(w)$ before reduction, so interpolation of residue values is unnecessary.

For each fixed seven-adic parameter, the binomial coefficients are integral.
Thus $R\in\Z_7[[u]]$, and its formal prepared unit is integral.
Write $R=W_x\mathcal U$. Modulo seven,
\[
 \mathcal U=5+(q+4)u^2+(6U+V+3)u^3+O(u^4).
\]
The four lower coefficients of $W_x$ modulo 49 are obtained by multiplying
the first table by the inverse of this series and retaining degrees below
four. No higher $u$ coefficient enters. The fold computation, replacing
$w$ by $w^2$, gives its constant coefficient
$7(6U+5V+6Z)$ modulo 49. Therefore
\[
 N_x=-W_x(-1)(1-W_f(0))
 \equiv-1+7\{q(1+5V+2Z)+4V\}\pmod{49}.
\]
For the trace, $D(0)\equiv42q\pmod{49}$, its coefficients of $u,u^2,u^3$
are divisible by seven, and its coefficient of $u^4$ is $3$ modulo seven.
Newton sums give $\operatorname{Tr}(u^4)\equiv-4W_x(0)$ and
$\operatorname{Tr}(u^j)\equiv0\pmod{49}$ for $j\ge5$.
All positive-degree traces are divisible by seven. Since
$W_x(0)\equiv21q\pmod{49}$, this yields
\[
 \mathcal T\equiv4(42q)+3(-4)(21q)\equiv14q\pmod{49}.
\]
The trace condition now forces $q=0$ over $\F_7$; the norm condition
then forces $V=0$, and hence $Z=-U$. Substitution in the same preparation
formulas gives \eqref{spx:r2:eq:digitquartic}.
\end{proof}

\subsection{A local norm--depth comparison}
The next theorem is an analytic statement independent of the coefficient equations. In this section the matrix
entries may be arbitrary elements of $\Z_7$ in the residue class
$(4,6,4)$; rationality is unnecessary. Use the critical and fold algebras
constructed in Lemma~\ref{spx:r2:lem:analytic}. Define
\[
 N_x=-\exp(\mathcal N),\qquad
 \delta_f=v(D_f),\qquad
 \mu=\min_{1\le i\le4}v(D_i),
\]
where the four critical roots are counted with multiplicity and $v(0)=\infty$.

\begin{theorem}[Local norm--depth comparison]\label{spx:r2:thm:comparison}
If $N_x=-1$ and $0<\delta_f<\infty$, then
\[
 \mu\le\delta_f.
\]
\end{theorem}
\begin{proof}
Suppose $\mu>\delta_f$. The trace inequality gives
$v(\mathcal T)\ge\mu>\delta_f$. The fold is rational over $\Q_7$,
so $\delta_f\ge1$, and $\mathcal T\equiv0\pmod{49}$.
Also $N_x=-1$ implies $\mathcal N=0$, since the exponential is injective
on $7\Z_7$. Lemma~\ref{spx:r2:lem:digits} gives
$V\equiv0$, $Z\equiv-U\pmod7$.
Then \eqref{spx:r2:eq:opp} and \eqref{spx:r2:eq:folddepth} imply $\delta_f\ge2$.

If $U\not\equiv0\pmod7$, \eqref{spx:r2:eq:digitquartic} has a coefficient of
$u$ of valuation one, constant coefficient of valuation at least two,
and coefficients of $u^2,u^3$ of valuation at least one. The Newton polygon
gives three roots of valuation $1/3$. At each, $v(W_x')=1$: the initial
root equation is $u^3=-[u]W_x$ and the initial derivative is $-3[u]W_x$.
Direct differentiation of the binomial elimination gives
\[
 \Phi'(u)=\frac{1+\Phi(u)}{b(1+u)}D(u).
\]
Hence at such a root $D=b(1+u)uW_x'\mathcal U$ and $v(D)=4/3$,
contradicting $\mu>\delta_f\ge2$.
Thus $U=V=Z=0$ over $\F_7$ and
$\gamma=\min_i v(\alpha_i)\ge2$.
The value $\gamma$ is finite: $A=A_*$ would give $t_f=0$ and
$\delta_f=\infty$, which is excluded.

By $\mathcal N=0$, \eqref{spx:r2:eq:linN} and Lemma~\ref{spx:r2:lem:analytic},
$v(\sigma+2\alpha_b)\ge2\gamma-1$. The exact fold determinant gives
\[
 \delta_f\ge2\gamma-1.
\]
Since $v(\mathcal T)>\delta_f$, \eqref{spx:r2:eq:linT} and the same remainder
bound imply $v(\sigma-2\alpha_b)\ge2\gamma-1$. Consequently
\[
 v(\sigma),v(\alpha_b)\ge2\gamma-1>\gamma,
 \qquad v(\alpha_a)=v(\alpha_d)=\gamma.
\]
The quartic coefficient estimates are therefore
\[
 v(c_0)\ge2\gamma-1,\qquad v(c_1)=\gamma,
 \qquad v(c_2),v(c_3)\ge\gamma.
\]
For $c_1$, use $105\alpha_a-120\alpha_b+15\alpha_d
=90\alpha_a+15\sigma-120\alpha_b$ and the fact that $90$ is a seven-unit.
Because $2\gamma-1>4\gamma/3$, the Newton polygon has a length-three
segment from $(1,\gamma)$ to $(4,0)$. Its three roots have
$v(t)=\gamma/3$ and $v(W')=\gamma$, by the same initial-derivative
calculation. The prepared factorization and \eqref{spx:r2:eq:G} give
\[
 v(D)=v(btW'U)=4\gamma/3<2\gamma-1\le\delta_f,
\]
contradicting $\mu>\delta_f$. This proves the comparison.
\end{proof}

\textbf{Application to R2.} The depth inequalities proved in Section~\ref{spx:r2:app:depth}
\eqref{spx:r2:eq:depths} give $\mu\ge\kappa>\kappa/2\ge\delta_f$.
The global unit norm supplies $N_x=-1$. Theorem~\ref{spx:r2:thm:comparison}
therefore excludes the normalized quintic, and the specified transport
excludes all six rows. This proves Theorem~\ref{spx:r2:thm:main}.

\subsection{The depth inequalities from the coefficient equations}\label{spx:r2:app:depth}
This subsection closes the depth step for the six quintic configurations.
It uses the coefficients of the original R2 datum before the unit-pivot
normalization. The already specified transformations subsequently preserve
the resulting inequalities.

Write $\tau=7^\kappa$, $\beta_1=\tau B$, $\beta_2=\tau C$,
$h=a+X$, $j=d+Y$ and
\[
 D=hj-b^2,\qquad s=b\beta_1-h\beta_2.
\]
The rational pair $(\beta_1,\beta_2)$ is integral and primitive.
At a place where $h$ is a unit we use
$Y=(b^2+D)/h-d$. Decompose $\Psi_n=D^{3n}\lambda_n$ by homogeneous
shift degree:
\[
 \Psi_n(B,C)=\sum_{\ell=0}^{n+1}D^\ell Q_{n,\ell}(B,C).
\]
The polynomials $Q_{n,\ell}$ used below are seven-integral. The meaning
of this divisibility is verified before any coordinate substitution.

\subsubsection*{A first-coefficient barrier}
The first coefficient has shift degrees zero, one and two, with
\begin{equation}\label{spx:r2:eq:Q12}
 Q_{1,2}(\beta_1,\beta_2)
 =\frac{j\beta_1^2-2b\beta_1\beta_2+h\beta_2^2}{2}
 =\frac{s^2+D\beta_1^2}{2h}.
\end{equation}
For a direct derivation in the preceding coefficient convention, let
$p(z)=z(1+z)$ and $p_2(z)=z(1+z)(1+2z)$, and use formal Gaussian
expectation with covariance
$D^{-1}\left(\begin{smallmatrix}j&-b\\-b&h\end{smallmatrix}\right)$.
The first two terms in its exponent are
\[
 U_1=\frac{p(X)z_1^3+p(Y)z_2^3}{6}
       -(B+X/2)z_1-(C+Y/2)z_2,
\]
\[
 U_2=-\frac{p_2(X)z_1^4+p_2(Y)z_2^4}{24}
       +\frac{p(X)z_1^2+p(Y)z_2^2}{4}-\frac{X+Y}{12},
 \qquad \lambda_1=\mathbb E(U_2+U_1^2/2).
\]
Gaussian moments are integer polynomials in the covariance entries.
The quadratic shift part uses second moments; the linear shift part uses
at most fourth moments; the constant part uses at most sixth moments.
Thus $\Psi_1$ has divisibility $D^2,D,1$ in these three parts, with
seven-integral quotients, and the quadratic part gives \eqref{spx:r2:eq:Q12}.
The expanded formula is checked exactly against all 256 monomials of the
preceding first-coefficient engine.

\begin{lemma}[Scalar pole barrier]\label{spx:r2:lem:scalarbarrier}
Suppose $h,s$ are units, $\delta=v(D)>0$, and $\lambda_1=r$ with
$v(r)\ge-2\kappa$. Then $\delta\ge\kappa$.
No hypothesis on the determinant cubic or its derivative is needed.
\end{lemma}
\begin{proof}
Homogeneity gives the exact equation
\[
 \tau^2D^3r=D^2Q_{1,2}(\beta_1,\beta_2)
       +\tau DQ_{1,1}(\beta_1,\beta_2)+\tau^2Q_{1,0}.
\]
The first term has valuation $2\delta$, by \eqref{spx:r2:eq:Q12}; the other
two have valuations at least $\kappa+\delta$ and $2\kappa$.
If $0<\delta<\kappa$, the first term is uniquely smallest. Consequently
$v(r)=-2\kappa-\delta<-2\kappa$, a contradiction.
\end{proof}

\subsubsection*{The fold calculation and its exact algebraic inputs}
Define
\[
 q_A(h)=d(d-1)h^3-b^2(b+2d-1)h^2
          +b^3(2a-1+b)h-b^3a(a-1).
\]
In the seven-integral ring with only $h$ inverted, the following identities
hold after $Y=(b^2+D)/h-d$ and
$\beta_2=(b\beta_1-s)/h$:
\begin{align}
 Q_{n,\ell}&\in(D,s)^\ell &&(n=2,3),\label{spx:r2:eq:foldideal}\\
 Q_{1,1}&\in(D,s),\label{spx:r2:eq:firstideal}\\
 h^{3n}Q_{n,\ell}|_{D=0}
   &=c_{n,\ell}(hs)^\ell(q_A/6)^{2n-\ell} &&(n=1,2,3).
   \label{spx:r2:eq:foldconic}
\end{align}
The coefficients are
\[
\begin{array}{c|rrrrr}
n&c_{n,0}&c_{n,1}&c_{n,2}&c_{n,3}&c_{n,4}\\\hline
1&15/2&3&1/2& & \\
2&405&135&18&1& \\
3&89505/2&14580&2025&144&9/2.
\end{array}
\]
Every displayed nonzero coefficient is a seven-unit.
The bundled conic checker derives the polynomials from the engine,
checks their $D^\ell$ divisibility and seven-integrality, makes the
coordinate substitutions, and verifies every coefficient of total
$(D,s)$-degree below $\ell$ is zero. It also verifies the complete
$D=0$ restrictions. The exceptional $Q_{1,2}$ term is handled by its
exact formula \eqref{spx:r2:eq:Q12}, not by an incorrect ideal-membership claim.
This is a generic polynomial check, not a sampled matrix calculation;
neither $b$, $q_A$ nor $q_A'$ is inverted.

\begin{lemma}[Fold depth and common scalar]\label{spx:r2:lem:foldnew}
Suppose $h$ and $q_A(h)$ are units and $\delta=v(D)>0$.
If $\lambda_1=r$ and $\lambda_2=\lambda_3=0$, then
\[
 v(s)\ge\delta,\qquad 0<\delta\le\kappa/2,
 \qquad v(r)=-2\kappa,
 \qquad \overline{\tau^2r}=\frac{\bar\beta_1^2}{2\bar h}\ne0.
\]
In particular $\beta_1$ is a unit, and
$\bar\beta_2/\bar\beta_1=\bar b/\bar h$.
\end{lemma}
\begin{proof}
Let $z=v(s)$, allowing infinity. If $z<\delta$, every nonconic term in
\eqref{spx:r2:eq:foldideal} has strictly higher valuation than its conic term:
a monomial $D^e s^f$ with $e\ge1$, $e+f\ge\ell$ has
$e\delta+fz>\ell z$.
Thus the valuation of shift-degree $\ell$ in $\Psi_n$ is
$\ell(\delta+z-\kappa)$. If $\delta+z-\kappa$ is nonzero, either the
highest shift term or the constant term in $\Psi_2$ is uniquely smallest,
which is impossible.
At equality put $L=hDs/\tau$ and $T=q_A/6$, both units, and
$e=\bar L/\bar T$.
The reductions of $\Psi_2=\Psi_3=0$, divided by the corresponding powers
of $T$, require a common root of
\[
 c_2(e)=e^3+4e^2+2e+6,\qquad
 c_3(e)=e^4+4e^3+2e^2+6e+5=e c_2(e)+5
\]
over the residue field. This is impossible over every extension of
$\F_7$. Therefore $z\ge\delta$.

By \eqref{spx:r2:eq:foldideal}, a positive shift degree $\ell$ in $\Psi_2$
now has valuation at least $\ell(2\delta-\kappa)$. If $2\delta>\kappa$,
the constant term, with unit conic value, is uniquely smallest. Hence
$\delta\le\kappa/2$.
Finally,
\[
 \tau^2r=\frac{\beta_1^2}{2h}+\frac{s^2}{2hD}
   +\frac{\tau Q_{1,1}}{D^2}+\frac{\tau^2Q_{1,0}}{D^3}.
\]
The last three terms have valuations at least
$\delta$, $\kappa-\delta$, and $2\kappa-3\delta$, respectively, all
positive. If $\beta_1$ were nonunit, $s\equiv0$ and the unit $h$ would
force $\beta_2$ nonunit, contradicting primitivity. Thus the first term
is a unit, giving the claimed scalar and direction.
\end{proof}

\subsubsection*{Applying these lemmas to all six configurations}
For the six rows, direct substitution gives the following residues:
\[
\begin{array}{c|cccccc}
(a,b,d)&(2,2,6)&(2,5,6)&(4,1,4)&(4,6,4)&(6,2,2)&(6,5,2)\\\hline
h_{\rm crit}&4&6&1&6&1&3\\
h_f&6&4&6&1&3&1\\
q_A(h_f)&2&5&2&5&2&5.
\end{array}
\]
The fold lemma applies in every row. It gives the same rational scalar
valuation $v(r)=-2\kappa$ at all completions and forces
$\bar\beta_2/\bar\beta_1=\bar b/\bar h_f$.
At every critical completion,
\[
 \frac{\bar s_{\rm crit}}{\bar\beta_1}
   =\bar b\left(1-\frac{\bar h_{\rm crit}}{\bar h_f}\right)
   =\frac{\bar b(\bar X_f-\bar X_{\rm crit})}{\bar h_f}\ne0,
\]
since the opposite residues are $2$ and $4$. Thus the scalar pole barrier
applies to all four critical embeddings and gives $v(D_i)\ge\kappa$.
The degree-one fold has integer positive determinant depth, so
$\delta_f\le\kappa/2$ also forces $\kappa\ge2$.
These are exactly \eqref{spx:r2:eq:depths}, with no assumption that the critical
contribution is a single quartic completion and no new census of shift
directions. This completes the depth step used by the R2 exclusion.

\subsection{Explicit matrix domain and coefficient convention}\label{spx:r2:app:conventions}
This subsection fixes the inputs of Theorem~\ref{spx:r2:thm:main} without referring
to an earlier classification or degree bound. All residue entries in this
appendix are in $\F_7$. The set $\mathcal M$ contains every $(a,0,d)$,
and, for $b\ne0$, precisely the pairs $(a,d)$ listed below.
\begin{center}\small
\begin{tabular}{c p{.85\linewidth}}
\toprule
$b$&Allowed pairs $(a,d)$\\\midrule
1&$(0,1)$, $(0,2)$, $(0,3)$, $(0,4)$, $(0,5)$, $(0,6)$, $(1,0)$, $(1,1)$, $(1,2)$,
 $(2,0)$, $(2,1)$, $(3,0)$, $(3,4)$, $(3,5)$, $(4,0)$, $(4,3)$, $(4,4)$, $(5,0)$, $(5,3)$, $(6,0)$, $(6,6)$\\[3pt]
2&$(0,1)$, $(0,2)$, $(0,3)$, $(1,0)$, $(1,5)$, $(2,0)$, $(2,3)$, $(2,6)$,
 $(3,0)$, $(3,2)$, $(3,6)$, $(5,1)$, $(5,5)$, $(6,2)$, $(6,3)$, $(6,6)$\\[3pt]
3&$(0,1)$, $(0,5)$, $(1,0)$, $(1,3)$, $(1,5)$, $(3,1)$, $(3,4)$, $(4,3)$, $(4,4)$, $(5,0)$, $(5,1)$, $(5,5)$\\[3pt]
4&$(0,1)$, $(0,3)$, $(0,5)$, $(1,0)$, $(1,3)$, $(3,0)$, $(3,1)$, $(3,3)$, $(4,4)$, $(4,5)$, $(5,0)$, $(5,4)$\\[3pt]
5&$(0,1)$, $(0,3)$, $(1,0)$, $(1,5)$, $(1,6)$, $(2,2)$, $(2,5)$, $(2,6)$, $(3,0)$, $(3,3)$,
 $(5,1)$, $(5,2)$, $(5,6)$, $(6,1)$, $(6,2)$, $(6,5)$\\[3pt]
6&$(0,0)$, $(0,1)$, $(0,6)$, $(1,0)$, $(1,2)$, $(1,3)$, $(1,4)$, $(1,5)$, $(1,6)$,
 $(2,1)$, $(2,2)$, $(3,1)$, $(3,5)$, $(4,1)$, $(4,4)$, $(4,5)$, $(5,1)$, $(5,3)$, $(5,4)$, $(6,0)$, $(6,1)$\\
\bottomrule
\end{tabular}
\end{center}
There are $49+21+16+12+12+16+21=147$ matrices.
The leading calculation enumerates the full Cartesian product of this explicit set, the eight projective shift directions, and the seven scalar residues. Its domain therefore has $147\cdot8\cdot7=8232$ entries.

The signed equations and coefficient convention are those stated at
 the beginning of this section and in Section~\ref{deg:sec:coefficients}.
For local lifting we multiply the clearing integer by six, so that all
unit exponents are multiples of six. This changes no necessary equation.
The coefficient numerator has shift degree $j+1$; its shift-homogeneous
part of degree $s$ is divisible by $D^s$ whenever used below.
Both assertions are verified by exact division in the recurrence of
Section~\ref{deg:sec:coefficients}; only denominators prime to seven
occur through order four.

\subsection{The total determinant-unit contribution}\label{spx:r2:app:unit}
We next bound the determinant-unit contribution.
All assertions in this subsection are under the explicit hypotheses of
Theorem~\ref{spx:r2:thm:main}.

\begin{lemma}[Unramified lifts]\label{spx:r2:lem:unramified}
At a place with $v(D)=0$, the completion generated by $x,y$ is unramified
over $\Q_7$. A fixed residue pair has at most one lift satisfying the
chosen signed equations.
\end{lemma}
\begin{proof}
After clearing denominators by an integer prime to seven, the logarithmic
Jacobian is a unit row multiple of
$\left(\begin{smallmatrix}a+X&b\\b&d+Y\end{smallmatrix}\right)$.
Hensel's lemma gives a unique lift over the unramified extension with the
given residue field. A lift in a ramified extension is the same lift.
Since the pair generates the completion, the assertion follows.
\end{proof}

Set $q=7^\kappa$ and $(U,V)=q(B,C)$, and choose one of the eight
projective residue directions $(1:v)$ or $(0:1)$, rescaling $q$ by a
rational seven-unit if needed. At a determinant-unit place $q^2r$ is
integral. Write $r_0=\overline{q^2r}$. The leading equations are
\[
 Q_{1,2}-r_0D=0,\qquad Q_{j,j+1}=0\quad(j=2,3,4),
\]
localized at $XY(1+X)(1+Y)D$. The coefficient checks in
Section~\ref{spx:r2:app:conventions} justify taking these leading equations.
If $v(r)<-2\kappa$, the first equation already excludes every such place.

\subsubsection*{Complete leading-support calculation}
The exact factorization and elimination over $\F_7$ gives
\begin{center}
\begin{tabular}{lr}
\toprule
Leading support&Number of fibres\\\midrule
Empty after localization&8057\\
Finite support of degree one&94\\
One common irreducible curve&81\\\midrule
Total&8232\\\bottomrule
\end{tabular}
\end{center}
In a common-curve fibre the first localized equation is exactly that
curve up to a unit, so there is no omitted isolated component.
The three curves are $XY-1$ (28 fibres), $XY+X+Y$ (28), and $X+Y+1$
(25). The calculation removes a factor only if it divides the unit
product, including proper factors of a reducible $D$.

For clarity, the finite calculation has a short exhaustive algorithm.
Factor the bilinear first equation. A vertical component is linear and
is treated by substitution. On every other irreducible component write
$A(X)Y+B(X)=0$ with $\gcd(A,B)=1$; substitute $Y=-B/A$ in the
other three equations and take their univariate gcd, removing only the
factors of $X(X+1)A B(A-B)$ and the substituted $D$. Its degree bounds
the support over the algebraic closure. In each fibre the bilinear factorization is checked before the substitution; cases with a vanishing coefficient of $Y$ are handled by direct substitution. Thus extension-field points
are included, not merely the $\F_7$-rational solutions.

The 81 curve fibres admit the checked $S,I,T,J$ routes to the chart
$XY=1$, with $d$ and the second shift of full pole order.
Every used matrix denominator is a seven-unit. The point maps preserve
the field and all four global units; their determinant ratios are units.
The coefficient identities above preserve all required equations.
Writing $t=d$, $z=B/C$, the normalized residues are
\[
 \bar a=\bar t\bar z^2,\qquad
 \bar b=1+\bar t\bar z,\qquad
 \overline{tC^{-2}r}=1/2,\qquad \bar t\ne0.
\]
The distribution of $\bar z=0,1,2,3,4,5,6$ is respectively
$25,6,3,9,3,10,25$. All routes have length at most three and depend
only on the rational datum, so the same route applies at all its places.

\subsubsection*{The additional normal order}
In this normalized chart put
\[
 z_m=(b-1)/t,\quad e=a/t-z_m^2,\quad S=XY-1,
 \quad \sigma=x+y-1.
\]
The letter $S$ in the rest of this subsection is a normal coordinate, not
the exchange map. We have $v(S)=v(\sigma)$.
Taking seven-adic logarithms of the signed Nahm equations and subtracting
$z_m$ times the second equation from the first gives exactly
\[
 \log(1-\sigma/y)-z_m\log(1-\sigma/x)=te\log x.
\]
In an unramified completion, $\log x\in7\mathcal O$. The coefficient
$-1/y+z_m/x$ on the left is a unit: on the residue curve it vanishes
exactly when $\bar X=\bar z_m$, whereas
\[
 \bar D=\bar t(\bar X-\bar z_m)^2/\bar X\ne0.
\]
If $v(\sigma)\ge1$, each logarithmic term of degree $n\ge2$ has
valuation $nv(\sigma)-v_7(n)>v(\sigma)$.
Consequently
\begin{equation}\label{spx:r2:eq:entrynormal}
 v(e)\ge\nu\ge1\ \Longrightarrow\ v(S)\ge\nu+1,
 \qquad e=0\ \Longrightarrow\ S=0.
\end{equation}
This is a valuation statement at every depth, not a finite digit test.

\subsubsection*{Ordinary and mixed transverse equations}
First put $\tau=1/C$, $R=t\tau^2r$ and
\[
 p=(p_a,p_b,p_R,\tau)
   =(a/t-z^2,(b-1)/t-z,R-1/2,\tau).
\]
Its entries lie in $7\Z_7$, and $p\ne0$ because $\tau\ne0$.
Let $\nu=\min v(p_i)$ and $\eta=\overline{p/7^\nu}$, taken projectively.
Normalize the four equations by
\[
 E_0=\tau^2\Psi_1/D^2-(R/t)D,\qquad
 E_{j-1}=\tau^{j+1}\Psi_j/D^{j+1}\quad(2\le j\le4).
\]
They are integral analytic functions on this determinant-unit chart.
They vanish identically on the retraction
$a=tz^2$, $b=1+tz$, $d=t$, $R=1/2$, $\tau=0$, $Y=1/X$.
Moreover $e=p_a-2zp_b-p_b^2$, so \eqref{spx:r2:eq:entrynormal} removes
the first normal digit. Taylor expansion modulo $7^{\nu+1}$ therefore
gives the necessary equations
\[
 \sum_{h\in\{a,b,R,\tau\}}
       (\partial_{p_h}E_i)_{p=S=0}\eta_h=0\quad(0\le i\le3).
\]
Terms quadratic in $p$ have valuation at least $2\nu\ge\nu+1$;
all normal terms also vanish to that precision. Only units are inverted.

For the six unit residues of $t$, the seven residues of $z$, and the
400 projective directions in $\mathbf P^3(\F_7)$, the gcds give
16776 empty entries, 12 linear entries and 12 flat entries. The complete
nonempty table, independent of $\bar t$, is
\[
\begin{array}{c|c|c}
\bar z&\eta&\text{support}\\\hline
2&(1:4:5:3)&X=4\\
4&(1:0:2:4)&X=2\\
0&(0:1:0:-1)&\text{first order zero}\\
-1&(0:0:0:1)&\text{first order zero}.
\end{array}
\]
The two flat directions require a different retraction. For $\bar z=0$
set $u=(B-1)/C$, $\ell=u$, so $B=1+u/\tau$. For $\bar z=-1$
set $u=B/C+1$, $\ell=u-1$, so $B=(u-1)/\tau$. In both cases
$u,\tau\in7\Z_7$. Put
\[
 R_*=\frac12+\frac{t\tau^2}{24},\quad
 p=(e_a,e_b,e_R,\rho)
   =(a/t-\ell^2,(b-1)/t-\ell,R-R_*,\tau u).
\]
The base restrictions $e_a=e_b=e_R=S=0$ of the four $E_i$ are
divisible by $\tau u$ over the seven-integral localized ring.
Indeed they vanish for $\tau=0$ by the infinite-shift retraction and
for $u=0$ by the exact finite-family identities
\[
\begin{array}{c|c|c}
(a,b,d;B,C)&\lambda_1&\lambda_2,\lambda_3,\lambda_4\text{ on }XY=1\\\hline
(0,1,t;1,C)&1/24+C^2/(2t)&0,0,0\\
(t,1-t,t;-C,C)&1/24+C^2/(2t)&0,0,0.
\end{array}
\]
Vanishing on both coordinate axes proves product divisibility;
it is not an assumption that $v(u)=v(\tau)$.

If $p\ne0$, take its least valuation and first projective digit as before.
Since $e=e_a-2\ell e_b-e_b^2$, the normal digit again vanishes.
The first equations have columns for $e_a,e_b,e_R$ and the coefficient
of $\tau u$ on the base. Higher terms are multiples of $p_ip_j$,
or of $\rho u,\rho\tau$, or involve the normal coordinate; they all
have strictly larger valuation. The exact derivative and mixed-column
calculations give
\[
\begin{array}{c|c|c}
\bar z& (\eta_a:\eta_b:\eta_R:\eta_\rho)&\text{support}\\\hline
0&(1:3:2:3)&X=5\\
-1&(0:1:4:1)&X=1.
\end{array}
\]
For each unit $\bar t$ there are 399 empty directions and one such
linear support. Thus all 4800 mixed entries are covered, with no flat
direction. Independent generic rational differentiation verifies the
finite-field row matrices used in both tables.

If $p=0$, then $\tau\ne0$ forces $u=0$ and one of the two exact finite
families. Equation~\eqref{spx:r2:eq:entrynormal} gives $x+y=1$. At a real
embedding the signed equations for $(0,1,t)$ force $|y|^t=1$ and hence
$y=\pm1$. Either a letter vanishes or $1-y=2$, contradicting global
units. For $(t,1-t,t)$ they force $|x|=|y|$; together with $x+y=1$
this gives $x=y=1/2$, again contradicting global units. Here $t\ne0$
because it is a seven-unit.

\begin{lemma}[Total determinant-unit bound]\label{spx:r2:lem:unitcap}
The sum of $[K_{\mathfrak p}:\Q_7]$ over all places with $v(D)=0$
is at most one.
\end{lemma}
\begin{proof}
The finite leading entries already have empty or single rational support.
For a curve entry the route, all rational transverse parameters and their
first nonzero projective digit are common to the whole datum. The tables
therefore leave at most one rational residue pair for all such places.
Lemma~\ref{spx:r2:lem:unramified} gives a unique $\Q_7$-rational lift.
Distinct places of $K=\Q(x,y)$ cannot give the same pair of lifted
coordinates, so this is a total bound, not a bound applied separately
to each place. The final partition has 8082 empty entries and 150 entries
bounded by one. No unit-norm filter has been applied to a proper subset
of the places.
\end{proof}

\subsection{Total local lengths and the six-configuration reduction}\label{spx:r2:app:entry}
We prove that Lemma~\ref{spx:r2:lem:unitcap}, the signed equations and four
global units leave only the six stated quintic configurations above
degree three. The proof uses the scalar barrier and fold bound proved above.

\subsubsection*{Determinant support and local lengths}
At a rank-one determinant point with $h=a+X$ a unit, eliminate the
first Nahm equation using its unit $x$ derivative. A vector field tangent
to it is
\[
 \mathcal V=-\frac b h\frac\partial{\partial\log x}
               +\frac\partial{\partial\log y}.
\]
Direct differentiation gives
\[
 \mathcal V(h)=-bX(1+X)/h,\qquad
 h^2\mathcal V(D)=q_A(h)+D R(h,D,a,b,d),
\]
where $R$ is a polynomial and $q_A$ is defined in Section~\ref{spx:r2:app:depth}.
The derivative of the second equation along this curve is a unit
multiple of $D/h$. If $b$ is a unit and $q_A$ has a root of multiplicity
$m$ at $H$, then $h-H$ is a local parameter and
\[
 \frac{dD}{dh}=u(h,D)q_A(h)+Dv(h,D),\qquad u(H,0)\ne0.
\]
In characteristic seven, successive coefficient comparison gives
$\operatorname{ord}_H D=m+1$ and second-equation order $m+2$, provided
$m+2<7$. Only division by $1,\ldots,m+2$ is used. Thus the local length
is $m+2$; a noncritical point has length two.

When $b=0$ in the residue field and $h$ is a unit, $Y=-d$ on $D=0$
and $\mathcal V(D)=hY(1+Y)\ne0$, giving length two. The other chart
is symmetric. At a rank-zero point $b=h=d+Y=0$, the two equations
have independent nonzero quadratic initial forms in the two coordinates.
The initial ideal contains $(u^2,v^2)$, so the length is at most four.
There is at most one rank-zero point, $(-a,-d)$.

These lengths bound the sum of all characteristic-zero local degrees
above the point. Indeed the completed integral Nahm algebra has finite
special-fibre length and is therefore finite over the unramified base
by complete Nakayama. Its generic-fibre dimension is at most that length.
The order generated by the four integral letters is a quotient, and its
generic fibre is the product of all relevant completions. Consequently
ramification indices and residue degrees are both counted. No squarefree
minimal-polynomial reduction is assumed.

For $b\ne0$, remove from $\bar q_A(h)$ only factors of
\[
 h(h-a)(h-a+1)(b^2-dh)(b^2+(1-d)h).
\]
These factors correspond exactly to excluded letters or to the unavailable
$h=0$ chart. Exact factorization for the explicit matrix set gives
\[
\begin{array}{c|rrrrr}
\text{type}&b=0&\text{no critical root}&m=1&m=2&m=3\\\hline
\text{number}&49&16&60&16&6.
\end{array}
\]
Every remaining factor is a single linear factor, with the stated
multiplicity. This proves there is at most one critical residue point
and it is rational, including over the algebraic closure of $\F_7$.
All cap-five points are the six cyclotomic pairs stated at the beginning of the $147$-residue argument.
The local jet checker also directly verifies 402 length-two, 60
length-three, 16 length-four and six length-five rank-one incidences,
and the 25 rank-zero quadratic forms.

Lemma~\ref{spx:r2:lem:foldnew} forces every noncritical determinant point to
be a single rational fold point for the fixed datum. When $b$ is a
unit it forces $\bar h=\bar b\bar U/\bar V$, and $U,V$ are units.
When $\bar b=0$, the $h$-unit chart forces $\bar V=0$ and
$\bar h=\bar U^2/(2r_0)$; the transposed chart forces $\bar U=0$.
Primitivity prevents both orientations. Thus there is at most one fold
in all cases.

\subsubsection*{A fold excludes the determinant-unit part}
At a determinant-unit residue point write $j=d+Y$, $s=bU-hV$,
$t=jU-bV$. The generic polynomial identities are
\[
 2bQ_{1,2}-UV D=st,\qquad
 6Q_{2,3}=Y(1+Y)s^3-X(1+X)t^3,
\]
\[
 ht-bs=DU,\qquad bt-js=DV.
\]
They are checked over $\Q$ before any reduction. A fold with $b$ a
unit forces the common rational scalar $r_0=\bar U\bar V/(2\bar b)$.
At a determinant-unit point the first leading equation therefore forces
$s=0$ or $t=0$. If $s=0$, then $h=bU/V$ and $t=DU/h$ are units;
the second leading equation is nonzero. The case $t=0$ is symmetric.
For $\bar b=0$ and an axis direction, one of $s,t$ is zero and the
other is a unit, giving the same contradiction. A fold and a
determinant-unit contribution therefore cannot coexist.

\subsubsection*{Cyclotomic removal and the norm argument}
Let $w^2+w+1=0$ in $\Q_7$ with residue 2 or 4. For
$x=w/(1+w)$, one has $x^2-x+1=0$, and $x,1-x$ are sixth roots of
unity. Hence every residue pair in $\{2,4\}^2$ has an exact cyclotomic
solution of the powered Nahm equations chosen above, for every rational
matrix. This point cannot be a conjugate of $K$, since $K$ has a real
embedding.

The generic fibre of the completed local Nahm algebra consequently has
the cyclotomic $\Q_7$ quotient and the distinct quotients by
all completions of $K$ above that point. The Chinese remainder theorem
gives
\[
 1+\sum_{\mathfrak p\text{ over the pair}}[K_{\mathfrak p}:\Q_7]
 \le\text{total local length}.
\]
This subtracts one even when the cyclotomic point itself is multiple.
In particular no determinant-unit place can reduce to a cyclotomic pair.

For $z\in\F_7\setminus\{0,-1\}$ put
$\nu(z)=(z^2,(1+z)^2)$. Its complete table is
\[
\begin{array}{c|ccccc}
z&1&2&3&4&5\\\hline
\nu(z)&(1,4)&(4,2)&(2,2)&(2,4)&(4,1).
\end{array}
\]
It is injective, never equals $(1,1)$, and its only inverse pair is
$\nu(2),\nu(4)$. For a rational residue pair $P=(X,Y)$ use
$\nu(P)=(\nu(X),\nu(Y))$. The four global-unit norms give
\[
 \prod_P\nu(P)^{n_P}=1,
 \qquad n_P=\sum_{\mathfrak p\text{ over }P}[K_{\mathfrak p}:\Q_7].
\]
Indeed a rational residue unit contributes its $e f$th power to the
reduced local norm, where $e,f$ are the ramification and residue degrees.
Squaring removes the global norm signs. This uses all places together.

One-point support forces $3\mid n_P$. For two distinct points, if
one exponent is one or two modulo three, the other is nonzero modulo
three. Opposite exponents would force equal points by injectivity;
equal exponents force coordinatewise opposite cyclotomic pairs by the
inverse-pair table.

\begin{lemma}[Complete R2 configuration reduction]\label{spx:r2:lem:entryreduction}
If $[K:\Q]>3$, its degree is five and its full decomposition above seven
has degree four at one of the six critical pairs in the six-row table and
degree one at the indicated opposite fold. No other place occurs.
\end{lemma}
\begin{proof}
Every place is determinant-unit, critical, or a fold, by the unit-letter
hypotheses and matrix integrality. All supports are rational by the
preceding analysis.

The determinant-unit contribution has total degree at most one.
Alone it cannot satisfy the norm table. Alongside a critical point,
the two-point norm argument forces it to be cyclotomic, which was
excluded. It cannot coexist with a fold. Thus it is absent.

At one remaining point the total degree is a multiple of three. Its
cap is at most five; every cap-five point is cyclotomic and loses one.
The degree is therefore at most three.

With two remaining points, one is critical and the other a fold of
initial cap two. The norm argument makes them opposite cyclotomic
pairs. The fold cap drops to one, and the critical cap drops by one.
The critical contribution must be one modulo three, hence one or four.
The total degree is two or five. Contribution four requires an original
critical cap five, so only the six triple-root rows occur. In particular
a rank-zero or double-root cap four cannot supply it after subtraction.
The norm lemma is applied to total local degrees, so arbitrary splitting
and ramification within the four-dimensional contribution are allowed.
\end{proof}

An equivalent finite calculation uses the following over-approximation of every possible support: all admissible
unit points with total cap one, every possible fold considered separately,
and the critical cap. Of the 147 matrix rows, 56 have upper spectrum
$\{0\}$, 85 have $\{0,3\}$ and six have $\{0,2,3,5\}$. Zero denotes
an empty selection, not a field. The six degree-five alternatives agree
exactly with those in the depth and norm--depth argument above.

Lemma~\ref{spx:r2:lem:entryreduction}, Section~\ref{spx:r2:app:depth}, and
Theorem~\ref{spx:r2:thm:comparison} now prove Theorem~\ref{spx:r2:thm:main} without
assuming an earlier degree-five or degree-seven theorem.

\section{The zero-cubic shift-pole residues}\label{spx:zero:section}
Work at seven, with $v(7)=1$, and at the matrix residue
\[
 A_0=\begin{pmatrix}0&1\\1&0\end{pmatrix},\qquad
 A=A_0+7^k\begin{pmatrix}a_k&b_k\\b_k&d_k\end{pmatrix},\quad k\ge1.
\]
The deformation matrix is primitive modulo seven. The exact reflection $J$ transfers the conclusions to residue $(1,-1,1)$. Put
\[
 X=\frac{x}{1-x},\quad Y=\frac{y}{1-y},\quad h=a+X,\quad
 D=(a+X)(d+Y)-b^2,\quad \Psi_j=D^{3j}\lambda_j.
\]
The letters $B,C$ denote the two linear shifts; $r=\lambda_1$ is the rational scalar. All four Nahm letters are local units, and $D\ne0$. The soft-point analysis concerns the conic chart $v(D)>0$.

For a shift pole set $\kappa=\max\{-v(B),-v(C)\}>0$, $\pi=7^\kappa$, $U=\pi B$, $V=\pi C$, and $Z=bU-hV$. A \emph{soft point} means $v(Z)>0$; there $h,V,h+1$ are units. Define
\[
 F=(2h+1)a+(h-h^2)(b-1)-(h^3+2h^2)d.
\]
The reduction of $F/7^k$, up to sign, is the first-order determinant cubic evaluated at $\bar h$.

\begin{theorem}[Soft-point reduction]\label{spx:soft:thm:soft}
If $\Psi_2=\Psi_3=\Psi_4=0$ at a soft point, then
\[
                  v(D)>k,\qquad v(F)>k.
\]
Thus every surviving soft point is a deeper cusp point. This theorem uses neither the degree-24 input, the rationality of $r$, global norms, nor the Nahm equations.
\end{theorem}

\subsection{The finite Newton certificate for Theorem~\ref{spx:soft:thm:soft}}
Clear only unit denominators and the shift pole by setting
\[
 P_j=\pi^{j+1}b^{j+1}h^{3j}\Psi_j,
 \quad X=h-a,\quad Y=\frac{b^2+D}{h}-d,\quad
 B=\frac{Z+hV}{b\pi},\quad C=\frac V\pi.
\]
These are polynomials over $\Z_{(7)}$. For a monomial record
$(n,l,i,t)$, where $n$ is its total degree in $a,b-1,d$ and $l,i,t$ its degrees in $D,Z,\pi$. Exact expansion of $P_2,P_3$ gives the following lower vertices:
\begin{center}\begin{tabular}{c cc}\toprule
Vertex&$P_2$&$P_3$\\\midrule
$\mathsf A$&$(4,0,0,3)$&$(6,0,0,4)$\\
$\mathsf B$&$(1,3,3,0)$&$(2,4,4,0)$\\
$\mathsf C$&$(0,4,2,0)$&$(0,6,2,0)$\\
$\mathsf D$&$(0,6,0,0)$&$(0,8,0,0)$\\
$\mathsf E$&$(0,4,0,1)$&$(0,6,0,1)$\\\bottomrule
\end{tabular}\end{center}
Every support exponent lies in the convex hull of these five vertices plus the nonnegative orthant. This is checked by extracting the supports, removing every exponent dominated componentwise by another, and testing rational convex-hull membership in the displayed five-vertex polyhedron. There are respectively 12 and 18 minimal exponents, from supports of sizes 408 and 1,348. Thus the comparison is uniform in all positive valuation depths, including fractional ones.

Write $m=v(D)>0$ and $s=v(Z)>0$. Assign weight $nk+lm+is+t\kappa$. Both polynomials have the same relevant lower-face structure. The possible nonsingleton faces are
\[
 \mathsf{AB},\quad\mathsf{BC},\quad\mathsf{CE},\quad
 \mathsf{ABCE},\quad\mathsf{CD},\quad\mathsf{CDE},
\]
with $\mathsf{AE}$ and $\mathsf{DE}$ included as boundary subfaces of the two larger faces. This list follows by comparing the five displayed weights. When $m\le k$, only $\mathsf C,\mathsf D,\mathsf E$ can be minimal.

\textbf{Single vertices.} Their coefficients are units times powers of $F/7^k$, $h$, $h+1$, and $V$. They cannot vanish if $v(F)=k$. The vertices $\mathsf C,\mathsf D,\mathsf E$ do not require $v(F)=k$.

\textbf{The three elementary edges.} After removing unit monomials, $\mathsf{BC}$ requires the same nonzero ratio to equal both $3$ and $2$. The edge $\mathsf{CE}$ requires the same nonzero ratio to equal both $1/2$ and $1/3$. On $\mathsf{CD}$ the equations are
\[
 3T^2+3TVh+V^2h^2=0,\qquad (2T+Vh)^2=0;
\]
substitution from the second leaves $V^2h^2/4\ne0$. All statements hold over the algebraic closure of $\F_7$, not just its rational points.

\subsection{The two remaining faces and the deeper cusp set}
On $\mathsf{ABCE}$ put
\[
 R=\frac{DZh}{F\pi},\qquad S=\frac{D^2Vh^3(h+1)}{F^2\pi}.
\]
The normalized equations are $p_2=p_3=p_4=0$. Their numerators are:
\begingroup\small
\begin{align*}
48p_2={}&8R^3-24R^2S+24R^2-48RS+30R+12S^2-30S+15,\\
1152p_3={}&144R^4-576R^3S+768R^3+576R^2S^2-2304R^2S\\
 &+1800R^2+1728RS^2-3600RS+2160R-192S^3\\
 &+1440S^2-2160S+1105,\\
384p_4={}&48R^5-240R^4S+384R^4+384R^3S^2-1536R^3S\\
 &+1400R^3-192R^2S^3+1920R^2S^2-4200R^2S+2880R^2\\
 &-768RS^3+3600RS^2-5760RS+3315R\\
 &+48S^4-800S^3+2520S^2-3315S+1695.
\end{align*}
\endgroup
Put $N_2=48p_2$, $N_3=1152p_3$, $N_4=384p_4$.
The following small certificate gives $A_2N_2+A_3N_3+A_4N_4=1$
in $\F_7[R,S]$:
\begin{align*}
 A_2={}&3R^4S+R^3S^2-R^3S-2R^3
       +3R^2S^3-R^2S^2+2R^2S+R^2\\
      &-RS^4-RS^3+RS^2+3RS-R
       -3S^5-S^4+2S^3-3S^2-3S-1,\\
 A_3={}&R^2S^2+3R^2S-3R^2-RS^2+3RS+R
       -3S^4-3S^2+S+3,\\
 A_4={}&3R^2S-3RS^2+3RS+S^3-2S^2-S-2.
\end{align*}
Direct multiplication verifies the identity. The generic face identities
are obtained by expansion before any specialization of $h$. The edge $\mathsf{AB}$ is their $S=0$ specialization, where already $p_2,p_3$ have gcd one. The edge $\mathsf{AE}$ is covered by $R=0$.

On $\mathsf{CDE}$, $\kappa=2m$ and $s\ge m$. Set $T=Z/D$, $Q=\pi/D^2$, and
\[
 N=QVh(h+1),\quad W=T^2+TVh,\quad M=V^2h^2.
\]
The three equations, after removing units, are
\[
 3N-6W-2M=0,\qquad4N-12W-3M=0,\qquad N-4W-M=0.
\]
The second minus three times the third gives $N=0$, impossible. This also includes $\mathsf{DE}$ with $\bar T=0$.

There are no omitted terms on these faces. In addition to the full $P_2,P_3$ supports, exact checks establish for $j=2,3,4$ that each shift-degree-$s$ part satisfies
\[
 \Psi_{j,s}/D^s\in(D,Z)^s
\]
after the displayed soft substitution, with $h,b$ inverted. The generic order $\Psi_j\in(a,b-1,d,D)^{2j}$ follows by the same exact coefficient extraction. Together these two orders control the fourth polynomial on both larger faces. This proves Theorem~\ref{spx:soft:thm:soft}.

For a soft point, $\bar h\in\{1,2,3,4,5\}$. There are exactly 40 incidences of $\bar h$ with a primitive projective deformation satisfying $\overline{F/7^k}=0$. Of these, five also satisfy
$\ell(\bar h)=a_k/\bar h+d_k\bar h-2b_k=0$; they are the rank-one directions $(a_k,b_k,d_k)\sim(\bar h^2,\bar h,1)$. The other 35 have $\ell\ne0$. By the preceding Nahm valuation dichotomy, a surviving deep soft place in those 35 incidences must have $v(x+y-1)=k$. They are treated by the Frobenius argument below.

\subsection{The finite entry calculation}\label{spx:zero:entry}
Let $(\alpha:\beta:\gamma)\in\mathbb P^2(\F_7)$ be the primitive
matrix direction and $(U:V)\in\mathbb P^1(\F_7)$ the shift direction.
Represent each projective vector by setting its first nonzero coordinate
to one. This gives exactly $57\cdot8=456$ pairs, with no omitted direction.
For a formal parameter $t$ substitute
\[
 a=t\alpha,\quad b=1+t\beta,\quad d=t\gamma,
 \quad X=h-a,\quad Y=((1+t\beta)^2+t\delta)/h-t\gamma.
\]
Multiply by $h^{3j}$ and extract the coefficient of $t^{2j}$ in
$\Psi_j$. From this coefficient retain only shift degree $j+1$,
evaluate $(B,C)=(U,V)$ and divide by $\delta^{j+1}$.
Call the resulting polynomial $R_j(h,\delta)$; a unit power of $h$
may be multiplied throughout. The division is exact. The order bound
$\Psi_j\in(a,b-1,d,XY-1)^{2j}$ and the shift-degree bounds show that
these are the leading equations for every positive matrix and shift
pole depth. All coefficients are integral at seven for $j\le4$.
Put
\[
 L(h)=\alpha-2\beta h+\gamma h^2,\qquad \ell(h)=L(h)/h.
\]
The calculation proceeds as follows, with all operations over $\F_7$.
\begin{enumerate}[leftmargin=*]
\item On $\delta=\ell(h)$ compute the gcd of
$R_2(h,\ell),R_3(h,\ell),R_4(h,\ell)$ after removing only
$h(h+1)L(h)$. A zero gcd is a curve and is retained.
\item Away from $\delta=0$ and $\delta=\ell(h)$, eliminate $\delta$
from each pair of $R_2,R_3,R_4$, take the gcd of their resultants,
factor it into irreducible factors in $h$, and in each resulting residue
field take the gcd of the three polynomials in $\delta$.
Check a nonzero two-equation Jacobian at each remaining closed point.
A vanishing resultant is treated as a common-component case, never as
an empty set.
\item At the central matrix $A_0$, divide $\Psi_j$ by $(XY-1)^{2j}$
and take shift degree $j+1$. Write the result as $\Xi_{j,j+1}$.
For each of the eight directions and each $r_0\in\F_7$, solve
\[
 \Xi_{2,3}=\Xi_{3,4}=\Xi_{4,5}=0,
 \qquad \Xi_{1,2}-r_0(XY-1)=0
\]
localized at $XY(1+X)(1+Y)(XY-1)$.
Restriction of the first three equations to $XY=1$, with only $h(h+1)$
removed, controls $0<v(D)<k$.
\end{enumerate}
These operations have the following output. There is no off-conic
support. The nonflat shallow support is soft and is excluded by
Theorem~\ref{spx:soft:thm:soft}. There are $126$ nonflat points with
$\delta\ne0,\ell(h)$; every one has a nonzero Jacobian. The curve
entries consist precisely of the eight pairs
\[
 ((H^2:H:1),(H:1))\quad(H\in\F_7),
 \qquad ((1:0:0),(1:0)).
\]
Of these, $H=0,-1,\infty$ are the boundary pairs; the other five have
$H=1,2,3,4,5$.
The equations for the five latter curves have no additional component
after saturation away from $\delta=\ell(h)$, the soft point, and
$\delta=0$.

The excluded nonflat $126$ points have unique unramified lifts by their
Jacobian test. At an equal-depth place $v(x+y-1)=k$, the power lemma
below instead forces ramification divisible by seven, provided
$42\nmid e$. Thus these points cannot occur. All remaining places
are on a listed curve or in the deeper determinant chart. The exact
coefficient construction above specifies the finite input and acceptance
tests independently of the archived output.

\subsection{The ramification alternative}
\begin{lemma}\label{spx:power:lem:power}
Let $F/\Q_7$ have ramification index $e$, and let $W$ be a principal unit.
If $\operatorname{v}(W^{7^k}-1)=k$ for an integer $k\ge1$, then either
\[
 \operatorname{v}(W-1)=1/7\qquad\text{or}\qquad42\mid e.
\]
\end{lemma}
\begin{proof}
For $\theta=\operatorname{v}(u)>0$, binomial expansion gives
\[
 \operatorname{v}((1+u)^7-1)=
 \begin{cases}7\theta,&\theta<1/6,\\1+\theta,&\theta>1/6,\end{cases}
\]
and at $\theta=1/6$ the next valuation is at least $7/6$.
An initial valuation at least $1/6$ would give final valuation greater
than $k$, so it is impossible.

If the first tie occurs after $j\ge1$ low-valuation steps, then
$\operatorname{v}(W-1)=1/(6\cdot7^j)$. Since valuations belong to $e^{-1}\Z$,
this implies $42\mid e$.
Otherwise let $j\ge1$ be the number of low-valuation steps before the
valuation exceeds $1/6$. Such a crossing must occur by step $k$.
The remaining steps each add one, so
\[
 7^j\operatorname{v}(W-1)+(k-j)=k.
\]
The preceding valuation was $j/7<1/6$, which forces $j=1$.
Thus $\operatorname{v}(W-1)=1/7$.
\end{proof}

\subsection{The local Frobenius restriction}
\begin{lemma}\label{spx:deep:lem:frob}
Under the local Nahm hypotheses and $42\nmid e$, suppose $v(\sigma)=k$. Write $H=\bar x/\bar y=\bar h$ and
$(\alpha,\beta,\gamma)=(\bar a_k,\bar b_k,\bar d_k)$. Then
\[
 N(H)=\gamma H^8-\beta(H^7+H)+\alpha=0.
\]
\end{lemma}
\begin{proof}
Choose an even positive integer $n$, prime to seven, clearing the denominators of $a_k,b_k,d_k$. The Nahm equations give
\[
 (1-\sigma/y)^n=W_1^{7^k},\quad (1-\sigma/x)^n=W_2^{7^k},
 \quad W_1=x^{na_k}y^{nb_k},\quad W_2=x^{nb_k}y^{nd_k}.
\]
Both $W_i$ reduce to one. Lemma~\ref{spx:power:lem:power} gives
$v(W_i-1)=1/7$ whenever $42\nmid e$. This is the only ramification
restriction in this argument.

In the maximal unramified subfield choose a lift $x_0$ of $\bar x$, and put $y_0=1-x_0$. The corresponding unit powers at $(x_0,y_0)$ differ from one by valuation at least one. Write $\epsilon=x-x_0$, so $y=y_0-\epsilon+\sigma$. The reduced linear coefficients of $W_i-1$ in $\epsilon$ are
\[
 g_1=n(\alpha/\bar x-\beta/\bar y),\qquad
 g_2=n(\beta/\bar x-\gamma/\bar y).
\]
If both vanish, the claimed identity follows directly. Otherwise a unit linear coefficient forces $v(\epsilon)=1/7$; higher terms, the base constants and the $\sigma$ terms have greater valuation. Both leading coefficients must therefore be nonzero. Comparing the unique lowest terms after the first seventh power gives
$\bar y g_1^7=\bar x g_2^7$. After cancelling common factors this is
\[
             (\alpha-\beta H)^7-H(\beta-\gamma H)^7=0,
\]
which equals $N(H)=0$ since $\alpha,\beta,\gamma\in\F_7$.
\end{proof}
Integer powers, including negative powers, have integral Taylor expansions around units. No expansion in unspecified branches of rational powers is used. The lemma allows arbitrary residue-field degree.

\subsection{Combining arithmetic with the cusp equation}
Over $\F_7$, define
\[
 F=(2H+1)\alpha+(H-H^2)\beta-(H^3+2H^2)\gamma,
 \quad L=\alpha-2\beta H+\gamma H^2,\quad M=\gamma H-\beta.
\]
Theorem~\ref{spx:soft:thm:soft} supplies $F(H)=0$ when $v(Z)>0$ and $v(D)>k$. For $v(Z)=0$, the following support calculation gives the same conclusion. Set $m=v(D)=k+\mu$, with $\mu>0$. The projected supports of
$P_j=(7^\kappa)^{j+1}b^{j+1}h^{3j}\Psi_j$, $j=2,3$, lie above the segment
\[
                  (2j,0,j+1)\ --\ (2j,j+1,0)
\]
in the weights $(k,\mu,\kappa)$. If $F(H)\ne0$, both singleton initial forms are nonzero. On the joining face the normalized equations have numerators
\[
 8R^3+24R^2+30R+15,\qquad
 144R^4+768R^3+1800R^2+2160R+1105.
\]
Their gcd modulo seven is one. Thus the non-soft deep branch also requires $F(H)=0$, uniformly in the three depths.

The exact identities
\begin{align*}
 F&=(2H+1)L-3H(H+1)M,\\
 N&=L+(H^7-H)M,\\
 (2H+1)N-F&=2H(H+1)(H^3+4H+6)(H^3+3H^2+6)M
\end{align*}
provide a short arithmetic exclusion. Both displayed cubics are irreducible over $\F_7$. Neither divides a nonzero polynomial $F$ of this form: coefficient comparison with the first would give
$(\alpha,\beta,\gamma)=(6c,2c,-c)$ and $0=4c$; comparison with the second gives $(6c,-c,-c)$ and $4c=0$.
Hence $F(H)=N(H)=0$, for $H\ne0,-1$, forces $M(H)=L(H)=0$. A separate check of all 57 projective matrix directions confirms exactly five rank-one roots and no others.

\begin{proof}[Deep-entry conclusion]
Use
\[
 D=(XY-1)+7^k(a_kY+d_kX-2b_k)+7^{2k}(a_kd_k-b_k^2),
 \quad XY-1=\sigma(1+X)(1+Y).
\]
If $v(D)>k$, then $v(\sigma)\ge k$. Equality would require $L(H)\ne0$ for cancellation, contradicting Lemma~\ref{spx:deep:lem:frob} and the identities above. Thus $v(\sigma)>k$ and $L(H)=0$. The cusp equation now gives $M(H)=0$. Primitivity implies $\gamma\ne0$, $H=\beta/\gamma\in\F_7$, and
$(\alpha,\beta,\gamma)\sim(H^2,H,1)$. The unit conditions remove $H=0,-1$.
\end{proof}

\subsection{The endpoint statement}
Work in a finite extension $F/\Q_7$ with $\operatorname{v}(7)=1$, with no restriction
on its ramification index. Let $t,z,\varepsilon\in\Q$ satisfy
\[
 k=\operatorname{v}(t)\ge1,\qquad \operatorname{v}(z)=\operatorname{v}(z+1)=0,\qquad
 \varepsilon\in7\Z_7,
 \quad
 A=\begin{pmatrix}t(z^2+\varepsilon)&1+tz\\1+tz&t\end{pmatrix}.
\]
Suppose $(x,y)\in F^2$ satisfies the signed Nahm equations after an even
denominator-clearing integer prime to seven. Define
\[
 X=\frac{x}{1-x},\quad Y=\frac{y}{1-y},\quad
 \sigma=x+y-1,\quad w=X-z,\quad D=(a+X)(d+Y)-b^2.
\]
Assume $D\ne0$, $\sigma\ne0$, $\operatorname{v}(\sigma)>k$ and $\operatorname{v}(w)>0$.
For rational linear shifts $B,C$, assume
\[
 \kappa=-\min\{\operatorname{v}(B),\operatorname{v}(C)\}>0,
 \qquad \lambda_1=r\in\Q
\]
in the specified convention $\Psi_1=D^3\lambda_1$. The letters $B,C$ are
both linear shifts; neither is the external scalar of a Nahm sum.

\begin{theorem}[Endpoint gap without a ramification cap]\label{spx:endpoint:thm:gap}
Under these hypotheses, $\operatorname{v}(w)\ge1/2$, allowing $w=0$.
Only the Nahm equations and rationality of the first coefficient are used
after entry into these coordinates. No global-unit norm constraint or
higher-coefficient vanishing is needed for this gap.
\end{theorem}

\subsection{The logarithm argument at every ramification depth}
For $\theta>0$ define
\[
 L(\theta)=\min_{j\ge0}(7^j\theta-j).
\]
If $\operatorname{v}(u)=\theta$, the terms of $\log(1+u)$ have valuations
$n\theta-\operatorname{v}_7(n)$. For a fixed $\operatorname{v}_7(n)=j$, the unique least term
is $n=7^j$. Consecutive candidate values differ by $6\cdot7^j\theta-1$.
Hence the minimum is unique except at
\[
 \theta=\frac1{6\cdot7^j}\quad(j\ge0).
\]
At a unique minimum, the logarithm has valuation $L(\theta)$.
At a tie its valuation can be larger, including infinity. This reasoning
controls the infinite tail because the consecutive differences eventually
become positive and then strictly increase.

Assume for contradiction $0<\rho=\operatorname{v}(w)<1/2$, and put
$s=\operatorname{v}(\sigma)-k>0$, $x_*=z/(1+z)$, $y_*=1/(1+z)$.
The logarithmic Nahm identities are
\begin{align}
 \log(1-\sigma/x)&=t(z\log x+\log y),\label{spx:endpoint:eq:nahm1}\\
 \log(1-\sigma/y)-z\log(1-\sigma/x)&=t\varepsilon\log x.\label{spx:endpoint:eq:nahm2}
\end{align}
Choose an even clearing integer $N$ prime to seven so that
$n=Nt/7^k$ and $nz$ are integers, with $n$ a seven-unit. Set
$W=x^{nz}y^n$. Then $W^{7^k}=(1-\sigma/x)^N$, so $\bar W=1$.
On $y=1-x$ the first derivative of $x^{nz}(1-x)^n$ at $x_*$ is zero;
its quadratic coefficient is a unit. Its value at $x_*$ is rational and
congruent to one. Integer powers have integral Taylor coefficients at
units. Since $2\rho<1<\operatorname{v}(\sigma)$, it follows that
\[
 \operatorname{v}(W-1)=2\rho,\qquad \operatorname{v}(\log W)=s.
\]
Also $\operatorname{v}(x-x_*)=\rho$ and $\log x_*\in7\Z_7$. At a nontie depth,
$\operatorname{v}(\log x)=L(\rho)$.

If $2\rho$ is not a logarithm tie, then $s=L(2\rho)>0$.
Since $L(2\rho)\le14\rho-1$, this forces $\rho>1/14$.
In this range only exponents 1 and 7 can minimize either logarithm.
The linear term in $\sigma$ on the left of \eqref{spx:endpoint:eq:nahm2} has valuation
$\operatorname{v}(\sigma)+\rho$, and all higher terms have greater valuation.
Thus $\varepsilon\ne0$, $\log x\ne0$, and
\[
 n_0:=\operatorname{v}(\varepsilon)=s+\rho-\operatorname{v}(\log x)\in\Z_{\ge1}.
\]
The three open intervals $(1/14,1/12)$, $(1/12,1/6)$, $(1/6,1/2)$
give respectively $n_0=8\rho$, $1-4\rho$, $2\rho$, all less than one.
At $\rho=1/6$, $s=1/3$ and $\operatorname{v}(\log x)\ge1/6$, so $n_0\le1/3$.
All are impossible.

The only remaining depths are therefore
\[
 \rho_h=\frac1{12\cdot7^h},\quad h\ge0.
\]
Here the minimum for $\log x$ is unique at exponent $7^{h+1}$, giving
\[
 \operatorname{v}(\log x)=-h-\frac5{12},\qquad
 s=n_0-h-\frac5{12}-\rho_h.
\]
Since $s>0$ and $n_0$ is integral, $n_0\ge h+1$. Consequently
$s\ge1/2$ and $n_0\ge s+1/2$ for every $h\ge0$.

\subsection{One certificate excludes the infinite sequence}
Put $S=(XY-1)/t$. Its valuation is $s$, and the exact determinant identity is
\[
 Q:=\frac{XD}{t}
 =w^2+(1+tX)\varepsilon+(X+t(z^2+\varepsilon))S.
\]
At any remaining depth $\rho_h$, the inequalities above give
$\operatorname{v}(Q)=2\rho_h$.
First set $q=7^\kappa$, $U=qB$, $V=qC$. The generic first coefficient
gives the polynomial identity
\[
 M=\frac{24X^3q^2\Psi_1}{t^2},\qquad
 tq^2r=\frac{M}{24Q^3}.
\]
Here $M$ is a polynomial in $z,w,t,\varepsilon,S,U,V,q$. Its lowest-valuation
part when $U-zV$ is a unit is
$12(U-zV)^2w^4$. The monomial comparison specified below verifies that every other term
has strictly larger valuation throughout
\[
 0<\rho\le1/12,\quad \operatorname{v}(S)\ge1/2,\quad
 \operatorname{v}(\varepsilon)\ge\operatorname{v}(S)+1/2,\quad \operatorname{v}(t),\operatorname{v}(q)\ge1.
\]
Thus $\operatorname{v}(tq^2r)=-2\rho_h$, a finite noninteger valuation, impossible
for a rational number.

If $U-zV$ is not a unit, $V$ must be a unit. Renormalize
\[
 q=1/C,\qquad B/C=z+\eta,\qquad
 j=\operatorname{v}(\eta)\ge1\ \text{or }\eta=0.
\]
The polynomial $M_0=M-12Q^3$, after setting $U=z+\eta,V=1$, has three
potentially smallest terms:
\[
 -12zSw^4,\qquad -24\eta w^5,\qquad 24qz(z+1)w^4.
\]
All other terms strictly exceed at least one of these weights over the
entire displayed domain, with $\operatorname{v}(\eta)\ge1$ when $\eta\ne0$.
After division by $24Q^3$, their valuations are
\[
 n_0-h-\frac5{12}-3\rho_h,\qquad
 j-\rho_h,\qquad \kappa-2\rho_h.
\]
The corresponding fractional deficits $5/12+3\rho_h$, $\rho_h$,
$2\rho_h$ are distinct and lie strictly between zero and one.
The leading coefficients are units because $z$ and $z+1$ are units.
There is a unique smallest valuation; if $\eta=0$, omit its term.
Therefore $tq^2r-1/2$ has a finite noninteger valuation, again impossible.
This proves Theorem~\ref{spx:endpoint:thm:gap} for every $h$, with no cutoff.

\textbf{How the uniform verification works.} Regenerate $M$ from the coefficient recurrence and the closed Wick formula for $\Psi_1$ in the scalar-barrier argument. The full polynomial has 1,128 terms. After removing the
three soft leading groups, 1,027 monomials remain. Substitute
\[
 \operatorname{v}(S)=\tfrac12+u,\quad \operatorname{v}(\varepsilon)=1+u+v,
 \quad\operatorname{v}(\eta)=1+j',\quad\operatorname{v}(q)=1+k',\quad\operatorname{v}(t)=1+l,
\]
with all five new variables nonnegative. Each remainder monomial is
assigned a leading group so that the difference has nonnegative
coefficients in these variables. Its affine part in $\rho$ is nonnegative
at $0$ and positive at $1/12$, hence positive throughout $0<\rho\le1/12$.
Every assignment is saved as an exact rational witness. No ramification
index or parameter depth is sampled in this domination proof.

\subsection{The exact conic}
\textbf{The exact conic already has degree at most four.} Suppose $\sigma=0$. At a real embedding, take logarithms of absolute values in the signed Nahm equations:
\[
 \begin{pmatrix}a&b-1\\b-1&d\end{pmatrix}
 \begin{pmatrix}\log|x|\\\log|y|\end{pmatrix}=0.
\]
This matrix is singular: otherwise $|x|=|y|=1$, incompatible with $x+y=1$. If $d=0$, singularity gives $b=1$; since $A\ne A_0$, it then forces $x=\pm1$, contrary to the global-unit assumptions. Hence put
\[
 t=d\ne0,\qquad z=(b-1)/t,\qquad a=tz^2,\quad b=1+tz.
\]
The values $z=0,-1$ similarly force $y=\pm1$ or $x=y=1/2$ and are excluded. Also $b\ne0$, since $b\equiv1\pmod7$.

Set $w=X-z$. On $XY=1$, exact substitution in the first coefficient gives a polynomial
\[
 P(w)=24t w^4(\lambda_1-r)\in\Q[w],\qquad \deg P\le4,\qquad
 P(0)=36z^2(z+1)^2(1+tz)\ne0.
\]
The identity follows by direct substitution with arbitrary $B,C,r,t,z$. Since $D=t w^2/X\ne0$ and $K=\Q(X)$, it proves $[K:\Q]\le4$. This handles the entire exact-conic alternative before the endpoint argument.

\subsection{The surviving disc is governed by an analytic cubic}
Work now on $v(w)\ge1/2$, including $w=0$. Choose $\pi\in\mathbb C_7$ with $\pi^2=7$ and put $w=\pi u$, $|u|\le1$. Statements about valuations of series below mean coefficientwise Gauss valuations in the ring of restricted power series in $u$ over $\Q_7(\pi)$.

The rational centers $x_*=z/(1+z)$ and $y_*=1/(1+z)$ give constants
\[
 F_0=z\log x_*+\log y_*\in7\Z_7,\qquad \ell_0=\log x_*\in7\Z_7.
\]
Direct differentiation gives
\[
 F'(w)=-\frac{w}{(z+w)(1+z+w)},\qquad
 F(w)=F_0-\frac{w^2}{2z(z+1)}+\cdots,
\]
and $\ell(w)=\ell_0+w/(z(z+1))+\cdots$. All order-$m$ coefficients of these logarithmic expansions have valuation at least $-v_7(m)$. Since
\[
                  m/2-v_7(m)\ge2\quad(m\ge4),
\]
the needed tails are controlled on the entire closed disc. For $m\ge7$ this follows from $m\ge7^{v_7(m)}$; the cases $4,5,6$ are immediate. In particular $v(F)\ge1$ in Gauss norm.

Write $\sigma=tT$. Equation~the first logarithmic equation becomes
\[
 \frac1t\log(1-tT/x)-\log(1+tT/y_0)=F(w).                 \tag{3}
\]
Every solution with $v(T)>0$ has $v(T)\ge1$: its linear coefficient is the unit $-(1+tX)/x$, and its nonlinear terms have larger valuation. On $T\in7\mathcal O$, equation (3) has a unique analytic solution $T=T(w)$. Indeed, putting $T=7V$ and rearranging gives a strict contraction in the complete ring of restricted series; the linear correction has coefficient divisible by $7$, and the order-$m\ge2$ coefficients have valuation at least $2(m-1)-v_7(m)>0$.

Expansion of (3) gives, in Gauss valuation,
\[
 T(w)=-\frac{xF(w)}{1+tX}+O(v\ge3).
\]
The other logarithmic equation is equivalent to the vanishing of
\[
 G(w)=-\frac1t\log(1+tT(w)/y_0)
       -\frac zt\log(1-tT(w)/x)-\varepsilon\ell(w).
\]
Its linear term in $T$ is $-wT/x$; all nonlinear terms have valuation at least $3$. Therefore
\[
 G(w)=\left(F_0-\frac{\varepsilon}{z(z+1)}\right)w
       -\frac{w^3}{2z(z+1)}+O(v\ge2).                 \tag{4}
\]
The constant $-\varepsilon\ell_0$ has valuation at least $2$ and is included in the remainder. All estimates retain the actual rational parameters.

\subsection{Counting endpoints and assembling the degree bound}
Dividing (4) by $\pi^3$ yields a restricted series whose reduction modulo $\pi$ is
\[
 \overline{\frac{G(\pi u)}{\pi^3}}
 =\overline{\frac{F_0-\varepsilon/(z(z+1))}{7}}\,u
       -\frac{u^3}{2\bar z(\bar z+1)}.                 \tag{5}
\]
The cubic coefficient is nonzero for every $\bar z=1,\ldots,5$. All coefficients of degree greater than three have strictly smaller absolute value. Strassmann's theorem, applied over the complete field $\mathbb C_7$, bounds the number of zeros on $|u|\le1$ by three. A repeated reduction is allowed; it does not invalidate the bound. The restricted-series form of Strassmann's theorem applies over every complete nonarchimedean field.

For each $w$, equation (3) determines $\sigma=tT(w)$ uniquely. Since the contraction has coefficients in $\Q_7$ and converges inside $\Q_7(w)$, these coordinates recover $x,y$ in the same local field. Every embedding of a completion $K_{\mathfrak p}$ gives a distinct such point. Consequently
\[
                 \sum_{\mathfrak p\mid7,\ v_{\mathfrak p}(D)>k}
                 [K_{\mathfrak p}:\Q_7]\le3.
\]
The rational matrix fixes $\bar z$, so this is one total bound, not three per incidence or per place. This proves the total deep contribution bound when $\sigma\ne0$.

\subsection{An exact retraction at infinite shift}
For a rank-one flat pair with $C$ of maximal pole depth, set
\[
 t=d,\quad z=B/C,\quad\tau=1/C,\quad
 \alpha=a/t,\quad\beta=(b-1)/t,\quad\delta=D/t,\quad\mathcal R=t\tau^2r.
\]
Here $v(t)=k$, $v(\tau)=\kappa>0$, and
$(\bar\alpha,\bar\beta)=(\bar z^2,\bar z)$. The five cases considered here have $\bar z=1,2,3,4,5$.
Normalize the equations as
\[
 J_j=\frac{\tau^{j+1}h^{3j}\Psi_j}{t^{2j}},\qquad
 E_0=J_1-\mathcal R h^3\delta^3,\quad E_j=J_{j+1}\ (j=1,2,3).
\]
The ideal-order and shift-degree identities make these integral after inverting $h$. The leading forms are obtained by coefficient extraction from the recurrence.

Retain the actual rational $t,z$ and retract to
\[
       \alpha=z^2,\quad\beta=z,\quad\mathcal R=1/2,\quad\tau=0.
\]
All four normalized equations vanish exactly on
\[
       \delta=g(h,t,z)=\frac{(h-tz^2-z)^2}{h-tz^2}.
\]
This is $XY=1$ for $a=tz^2,b=1+tz,d=t$. At $\tau=0$ only the highest shift-degree terms remain. On this conic the first has value $C^2/(2t)$, and the highest shift-degree terms of $\Psi_2,\Psi_3,\Psi_4$ vanish. All four identities were checked generically in $t,z$ with the unchanged engine.

This boundary retraction is not a finite-shift Nahm datum. It still supplies an exact zero locus for the normalized polynomial system. The actual transverse vector
\[
             p=(\alpha-z^2,\beta-z,\mathcal R-1/2,\tau)
\]
is nonzero because $\tau\ne0$. At a flat $G_1$ place it reduces to zero; the first coefficient forces $\overline{\mathcal R}=1/2$.

Let $n_i$ be the normal derivatives in $\delta-g$, and $L_i(\eta)$ the transverse derivatives in direction $\eta$. We compute the gcd of all six wedges $n_iL_j-n_jL_i$. Only $\bar h=0,-1,\bar z$ are removed. The last is the deep endpoint, deliberately left for a separate argument.

\begin{center}\small
\begin{tabular}{c r r r r r}\toprule
$\bar z$&Directions&Empty&Degree 1&Degree 2&Flat\\\midrule
0&400&378&21&0&1\\
1&400&397&2&1&0\\
2&400&396&4&0&0\\
3&400&397&2&1&0\\
4&400&396&4&0&0\\
5&400&397&2&1&0\\
$-1$&400&378&21&0&1\\\bottomrule
\end{tabular}\end{center}

\subsection{Uniform lifting of the five flat pairs}
The normal vector has no simultaneous admissible zero in any of the seven rows. No admissible gcd has repeated factors. For $\bar z=1,\ldots,5$, every gcd is a product of distinct linear factors over $\F_7$ and has degree at most two.

\begin{theorem}[Five flat-pair contribution bounds]\label{spx:flat:thm:flat}
For each pole-shift flat pair $\bar z=1,2,3,4,5$, the total local degree contributed by its $G_1$ places is at most two. Every such completion is $\Q_7$. The deep endpoint $\bar h=\bar z$ is not included.
\end{theorem}
\begin{proof}
Let $\nu\ge1$ be the minimum valuation of the four rational transverse parameters, and write $p=7^\nu\eta$, $\bar\eta\ne0$. Exact retraction leaves no constant error for the actual $t,z$. A unit normal derivative forces
$v(\delta-g)\ge\nu$. Put $w=(\delta-g)/7^\nu$. The reduced equations are
$\bar n_iw+\bar L_i(\eta)=0$. A unit $n_i$ and simplicity of the common wedge gcd give a nonzero two-equation Jacobian in $(h,w)$. Hensel's lemma gives a unique unramified lift at each residue point. Those points are rational over $\F_7$, and these coordinates generate the same local field as $X,Y$. There are at most two points. This works at every matrix and pole depth.
\end{proof}

\subsection{The three boundary flat pairs}
The pairs are $((0:0:1),(0:1))$, $((1:-1:1),(-1:1))$,
and their transposed axis pair $((1:0:0),(1:0))$.
The first two correspond to $z=0,-1$ below, and the third follows
by exchange of coordinates.
\subsection{Every place enters the controlled chart}
Normalize every extension of $v_7$ by $v(7)=1$. The four global-unit conditions make $X,Y,X+1,Y+1$ local units at every place over seven. Thus $h$ is a unit and $v(D)\ge0$.

The deep-entry argument above excludes a deep place for the three
boundary matrix directions. Thus it remains to control $v(D)\le k$.

For the first two rows put $t=d$, so $v(t)=k$, and write $z=0,-1$, respectively. exact checks give the following exhaustive entry table.
\begin{center}\small
\begin{tabular}{p{.19\linewidth}p{.72\linewidth}}\toprule
Depth & Necessary condition and certificate\\\midrule
$v(D)=0$ & At $A_0$, divide $\Psi_j$ by $(XY-1)^{2j}$ and take its highest shift-degree part. For each of the seven possible residues of $7^{2\kappa}r$, the saturated ideal is the unit ideal.\\[3pt]
$0<v(D)<k$ & The highest shift parts on $XY=1$ have gcd $1$ after removing only $h,h+1$. This applies also to fractional depths.\\[3pt]
$v(D)=k$ & The leading equations force
$\bar\delta=(\bar h-z)^2/\bar h$, where $\delta=D/t$. Saturating away from this curve gives the unit ideal.\\[3pt]
$v(D)>k$ & Impossible by the preceding deep-cusp reduction.\\\bottomrule
\end{tabular}
\end{center}
In the third row $\bar h\ne0,-1$, hence $\bar h-z\ne0$. The forced curve is the chart called $G_1$: $v(D)=k$ and $v(XY-1)>k$. There is no omitted deep endpoint for these boundary pairs.

The entry calculations use the exact generic identities
\[
 \Psi_j\in(a,b-1,d,XY-1)^{2j}\Q[X^{\pm1},Y,a,b,d,B,C],
 \qquad \deg_{B,C}\Psi_j=j+1,
\]
for $1\le j\le4$. These two properties are checked by exact expansion and coefficient extraction. Thus matrix corrections have strictly greater valuation when $0<v(D)<k$; no discretization of the possible ramification depths is used.

\subsection{An exact two-boundary retraction}
Put $\tau=1/C$ and $\mathcal R=t\tau^2r$. The two cases use slightly different coordinates:
\[
\begin{array}{c|c|c|c}
z&u&\ell&B\\\hline
0&(B-1)/C&u&1+u/\tau\\
-1&B/C+1&u-1&(u-1)/\tau.
\end{array}
\]
Both $u$ and $\tau$ belong to $7\mathbb Z_7$, and $\tau\ne0$. Retain the actual rational values $t,u,\tau$. Define
\begin{gather*}
 \varepsilon_a=a/t-\ell^2,\qquad \varepsilon_b=(b-1)/t-\ell,\qquad
 \varepsilon_R=\mathcal R-\mathcal R_*,\qquad \rho=\tau u,\\
 \boxed{\mathcal R_*=\frac12+\frac{t\tau^2}{24}}.
\end{gather*}
The factor in the correction is $t$, not $t^2$. On the finite exact families the first coefficient is $r=1/24+C^2/(2t)$.

Use normalized equations
\[
 J_j=\frac{\tau^{j+1}h^{3j}\Psi_j}{t^{2j}},\qquad
 E_0=J_1-\mathcal R h^3\delta^3,\qquad E_j=J_{j+1}\quad(1\le j\le3).
\]
The ideal-order and shift-degree checks make these integral functions after inverting unit denominators. At $\varepsilon_a=\varepsilon_b=\varepsilon_R=0$ the matrix is
$a=t\ell^2$, $b=1+t\ell$, $d=t$. Its conic $XY=1$ is
\[
 \delta=g(h,t,\ell)=\frac{(h-t\ell^2-\ell)^2}{h-t\ell^2}.
\]

\begin{lemma}[Mixed divisibility]\label{spx:boundary:lem:mixed}
Each restriction $E_i(h,g,0,0,0)$ is divisible by $\rho=\tau u$ in the ring localized at the unit denominators.
\end{lemma}
\begin{proof}
At $\tau=0$, the generic rank-one highest-shift identities make all four restrictions zero for every $\ell,t$. These identities are checked over $\Q$, before reduction. At $u=0$, the two finite-shift families are
\[
 (a,b,d;B,C)=(0,1,t;1,C),\qquad
 (t,1-t,t;-C,C).
\]
For both, exact division by $XY-1$ verifies $\Psi_j=0$ for $j=2,3,4$, and
$t\Psi_1=(t/24+C^2/2)D^3$ on the conic. Therefore the normalized restrictions also vanish identically at $u=0$. After clearing unit denominators, divisibility by $u$ and by $\tau$ implies divisibility by $u\tau$.
\end{proof}
This product measures the full base error. A linear expansion that retains only $\tau$ misses this structure and leaves a flat direction.

\subsection{The transverse certificate and its lifting argument}
Set $q=\delta-g$. Taylor expansion in $(q,\varepsilon_a,\varepsilon_b,\varepsilon_R)$ gives
\[
 E_i=\rho K_i+n_iq+L_{ia}\varepsilon_a+L_{ib}\varepsilon_b+L_{iR}\varepsilon_R
       +O\bigl((q,\varepsilon_a,\varepsilon_b,\varepsilon_R)^2\bigr).
\]
All coefficients are integral functions of the preceding $h,t,u,\tau$. The exact divisibility lemma supplies $K_i$; its reduction is the coefficient of $\tau u$ in the base restriction. Each row is multiplied by a power of the unit $h$ to clear denominators before taking gcds.

The normal derivatives $\bar n_i$ have no simultaneous zero for $\bar h\ne0,-1$. For each of the $400$ directions
\[
 \eta=(\eta_a:\eta_b:\eta_R:\eta_\rho)\in\mathbb P^3(\F_7),\qquad
 L_i(\eta)=L_{ia}\eta_a+L_{ib}\eta_b+L_{iR}\eta_R+K_i\eta_\rho,
\]
compute the gcd of the six wedges $n_iL_j-n_jL_i$. Only factors $h$ and $h+1$ are removed. The exact census is
\begin{center}
\begin{tabular}{c r r r r r r}\toprule
Case&Directions&Empty&Degree 1&Degree 2&Flat&Repeated\\\midrule
$z=0$&400&387&12&1&0&0\\
$z=-1$&400&387&12&1&0&0\\\bottomrule
\end{tabular}
\end{center}
The unique degree-two support in each row is irreducible:
\[
\begin{array}{c|c|c}
z&\eta&\text{monic gcd}\\\hline
0&(1:3:3:2)&h^2+2h+2\\
-1&(0:1:6:3)&h^2+1.
\end{array}
\]

\begin{lemma}[Uniform local contribution]\label{spx:boundary:lem:lift}
If $p=(\varepsilon_a,\varepsilon_b,\varepsilon_R,\rho)\ne0$, the total local degree contributed by the controlled chart is at most two. Its supports are either a single degree-one point or a single degree-two point, with unique unramified lifts.
\end{lemma}
\begin{proof}
The first equation on $G_1$ forces $\overline{\mathcal R}=1/2$, so all entries of $p$ have positive valuation. Since they are rational, $\nu=\min v(p_i)$ is a positive integer. Exact mixed divisibility and a unit normal derivative force $v(q)\ge\nu$: if $0<v(q)<\nu$, its linear normal term would be the unique term of least valuation. Write $w=q/7^\nu$ and $p=7^\nu\eta$.

The reduced equations are $\bar n_iw+\bar L_i(\bar\eta)=0$. At each common root, a unit normal derivative and the squarefree wedge gcd give a nonzero two-equation Jacobian in $(h,w)$. Hensel lifting gives at most one point over the unramified extension of the residue degree. These coordinates recover $X,Y$ rationally and generate the same completion. The census now gives the asserted bound, including multiplicities. No algebraic independence of the entries of $p$ is assumed: enumerating every direction is an overestimate.
\end{proof}

\subsection{Whole-field norms exclude the remaining supports}
By the entry table, every place of $K$ over seven is in the controlled chart. The rational vector $p$ has one common reduced projective direction for all those places. Thus the supports in Lemma~\ref{spx:boundary:lem:lift} bound the entire decomposition of $K$ at seven, not merely part of it.

Both $X$ and $X+1=1/(1-x)$ are global units; consequently their rational norms belong to $\{1,-1\}$. On the residue conic $\bar X=\bar h$ and $\bar Y=1/\bar h$.

If the support is linear, there is at most one degree-one completion. No $H\in\F_7\setminus\{0,-1\}$ has both $H$ and $H+1$ in $\{1,-1\}$. A nonempty decomposition therefore contradicts the two unit norms.

If the support is quadratic, there is at most one unramified degree-two completion. For the two irreducible supports the obstruction is explicit:
\[
\begin{array}{c|c|c}
z&\text{support}&\text{forbidden residue norm}\\\hline
0&h^2+2h+2&N(\bar X)=2\\
-1&h^2+1&N(\bar X+1)=2.
\end{array}
\]
Since $2\notin\{1,-1\}\subset\F_7$, neither support can give a nonzero field. An independent finite assembly check of all $800$ directions gives bound zero in every case, using the four norms of $X,X+1,Y,Y+1$.

It remains to handle $p=0$, which is not a projective census direction. Since $\tau\ne0$, it forces $u=0$ and one of the two exact finite families. The unit normal derivative forces $q=0$, hence $XY=1$ and $x+y=1$ in characteristic zero. At a real embedding the signed Nahm equations then imply
\[
 \begin{array}{c|c}
 (0,1,t)& |y|^t=1\\
 (t,1-t,t)& |x|=|y|.
 \end{array}
\]
Here $t\ne0$. In the first case $y=\pm1$, so $1-y$ is zero or $2$; in the second $x=y=1/2$. Both contradict the global-unit hypotheses. Transposition proves the third row of the boundary-pair exclusion. This finishes the restricted exclusion.

\subsection{The extra Nahm order and the total degree}
For the five nonboundary flat pairs use \emph{shift} coordinate $z_s=B/C$ (distinct from the matrix coordinate $z$ above), and put
\[
 \tau=1/C,\quad \mathcal R=t\tau^2r,\quad
 p=(a/t-z_s^2,(b-1)/t-z_s,\mathcal R-\tfrac12,\tau).
\]
All entries are rational of positive valuation, and $p\ne0$ because $\tau\ne0$. Let $\nu=\min v(p_i)\ge1$ and $\eta=\overline{7^{-\nu}p}$.

The exact retraction and squarefree transverse wedge gcds above prove that every admissible $G_1$ completion is $\Q_7$. This must be established \emph{before} using $v(\log x)\ge1$. The normal derivatives have no simultaneous admissible zero; a simple wedge root supplies a nonzero two-equation Jacobian and hence unique unramified lifting.

Let $z_m=(b-1)/t$ and $\varepsilon=a/t-z_m^2$. Then $v(\varepsilon)\ge\nu$, and the exact logarithmic identity is
\[
 \log(1-\sigma/y)-z_m\log(1-\sigma/x)=t\varepsilon\log x.
\]
At $G_1$, the linear coefficient of $\sigma$ is a unit: its residue is $(\bar z_s-\bar X)/\bar x$. All higher terms have greater valuation. Since the completion is $\Q_7$, it follows that
\[
 v(\sigma)\ge k+\nu+1,\qquad
 v\big((XY-1)/t\big)\ge\nu+1.
\]
If the right side vanishes, the same comparison forces $\sigma=0$, covered separately.

The exact retraction converts this extra order into
\[
 h\bar u+2h\eta_b-\eta_a=0,
\]
where $\bar u$ is the normalized displacement from the conic. If the four coefficient rows are $n_i\bar u+L_i(\eta)=0$, eliminate $\bar u$ to obtain
\[
 n_i(\eta_a-2h\eta_b)+hL_i(\eta)=0\quad(0\le i\le3).
\]
Of the twenty admissible pairs $(H,\bar h)$, eighteen resulting $4\times4$ matrices have rank four. The only projective kernels are
\[
 \begin{array}{c|c|c}
 H&\bar h&\eta\\\hline
 2&4&(1:4:5:3)\\
 4&2&(1:0:2:4).
 \end{array}
\]
A separate implementation using exact rational evaluation and exhaustive finite-field kernel enumeration confirms every matrix. The full 2,000-direction gcd calculation confirms there are no additional or identically zero directions. In particular no degree-five mixed support survives.

For a fixed datum the remaining $G_1$ contribution is at most one. If it occurs, the norms of the global units $X$ and $X+1$ require
\[
 H^{d_0}\bar h\in\{\pm1\},\qquad (H+1)^{d_0}(\bar h+1)\in\{\pm1\}.
\]
For either displayed pair and $0\le d_0\le3$, this forces $d_0=1$, giving total degree two. Without a $G_1$ contribution the deep bound is three. Hence off the exact conic the zero-cubic shift-pole branch has degree at most three; on the conic it has degree at most four.

\begin{theorem}[Four-unit shift-pole bound]\label{spx:thm:all}
Under the hypotheses at the beginning of this section, assume
$\lambda_1\in\Q$ and $\lambda_2=\lambda_3=\lambda_4=\lambda_5=0$.
In the two zero-cubic matrix residue classes assume also that
$42\nmid e_{\mathfrak p}$ at every place over seven.
Then $[K:\Q]\le4$. In the remaining $147$-residue domain the stronger
bound is three, and the fifth coefficient is unnecessary.
\end{theorem}
\begin{proof}
The finite residue graph in Section~\ref{deg:sec:entry-completion}
partitions the $343$ matrix residues into the $194$ classes reaching
a regular residue, the explicit $147$-class domain, and the two
zero-cubic classes. Every legal edge preserves the field, the four
unit letters, the signed equations, the positive shift-pole order,
and the coefficient conditions. The regular proof gives bound four,
Theorem~\ref{spx:r2:thm:main} gives bound three on the second domain,
and the endpoint, boundary and norm arguments above give four on
the zero-cubic domain. The exact matrices at the two centres
have no admissible datum: at $A_0$ the real signed equations imply
$|1-x|=|y|$, $|1-y|=|x|$, whence $x+y=1$ and $D=0$.
The reflection $J$ treats the other centre.
These are alternative parameter domains, so their caps are not added.
\end{proof}

In the modular application, Theorem~\ref{d24:thm:main} supplies
$e_{\mathfrak p}\le[K:\Q]\le24$, which excludes $42\mid e_{\mathfrak p}$.
The first five coefficient equations follow from
Lemma~\ref{deg:lem:modular-bridge}. Thus the three shift-pole rows
use only statements proved in the paper.

%% file: sections/arithmetic_overview.tex
\section*{Overview of the arithmetic reduction}
\phantomsection\addcontentsline{toc}{section}{Overview of the arithmetic reduction}

The degree bound proved in Part~\ref{part:degree} is the starting point for the remaining classification. It places the distinguished Nahm point $(x,y)$ in a joint field $K=\Q(x,y)$ of degree at most four, but it does not by itself bound the heights of $x$ and $y$. Moreover, a fixed algebraic point can belong to a rational one-parameter family of matrices. We therefore need to determine both the possible points and all the matrices and shifts attached to them. Parts~\ref{part:arithmetic}--\ref{part:quartic} carry this out by separating the arithmetic branches, proving explicit bounds for the branches requiring enumeration, and applying successively stronger necessary conditions to the resulting candidates. The purpose of this overview is to explain what is enumerated and why the successive exclusions are exhaustive.

\paragraph{1. Separate the scalar families from the independent pairs.}
Write the Nahm equations as
\[
 1-x=x^\alpha y^\beta,\qquad 1-y=x^\beta y^\gamma,\qquad
 A=\begin{pmatrix}\alpha&\beta\\\beta&\gamma\end{pmatrix}>0.
\]
Symmetry gives
\[
 [x]\wedge[1-x]=-[y]\wedge[1-y],
\]
where $[u]$ is the class of $u$ in $K^\times\otimes_\Z\Q$ and the product is taken in the exterior square $\bigwedge^2(K^\times\otimes_\Z\Q)$. If this product vanishes and $\beta\ne0$, all four elements lie in one rank-one multiplicative group. Apart from the rational point, they can be written as powers of a common unit $v$, with complementary powers satisfying $v^r+v^s=1$. The scalar factorization argument in Proposition~\ref{low:scalar} reduces the primitive exponents to $\max(r,s)\le6$. It gives thirteen scalar exponents and at most thirty-three one-parameter matrix templates. The diagonal case is handled separately and has at most $13^2=169$ candidate matrices before checking the joint degree. These finite lists of templates retain their rational parameters; those parameters are determined later by the modular coefficient equations.

If the exterior product is nonzero, the two complementary pairs span the same rank-two multiplicative plane. This independent branch splits according to whether all four elements are ordinary algebraic units. The distinction is necessary: the unit-lattice estimates below apply only after the nonunit possibilities have been covered.

The arithmetic classification in degrees one and two is supplied by Theorem~\ref{low:global}. Steps~2--5 below are used for the cubic and quartic classifications in Parts~\ref{part:cubic} and~\ref{part:quartic}; the shift computation of Step~6 is also used for the data of degree at most two in Part~\ref{part:classification}.

\paragraph{2. Put nonunit points into bounded exponent families.}
Finite valuations of the two Nahm equations, together with positive definiteness, leave five possible nonzero valuation patterns. Lemma~\ref{qcov:entry} turns these patterns into two canonical equations, after coordinate exchange and simultaneous complementation when needed:
\[
 (1-z^b)^N=z^T(1+z^a)^M\qquad\text{or}\qquad(1+z^b)^N=z^T(1+z^a)^M.
\]
Here $z\in(0,1)$ is an ordinary unit, $\Q(z)=K$, $g=\gcd(a,b)$, $M=a/g$ and $N=b/g$; primitivity and positive definiteness impose additional explicit conditions. Thus the new generator is a unit even when the original coordinates are not. The two binomial factors must have the same nonempty set of prime-ideal divisors. A primitive divisor of one factor can therefore occur only when the exponents satisfy the corresponding divisibility restrictions. This gives a proved cutoff for $a,b$. A height inequality bounds $T$, and archimedean and norm comparisons shorten its possible interval. Exceptional exponent equalities are treated separately before the remaining census is formed.

The subsequent sieve tests necessary polynomial divisibility over finite fields. A genuine degree-$d$ unit generator has a monic minimal polynomial with constant coefficient $\pm1$; its reduction still has degree $d$ at every prime. After justified cyclotomic cancellation, the canonical polynomial must consequently have a degree-$d$ divisor of the required constant sign modulo every tested prime. Reducible and repeated reductions are included. A failed test at one prime excludes the characteristic-zero candidate. For example, in the quartic branch the analytic and structural restrictions leave $9{,}626$ exponent triples; the primes $2,3,5,7$ leave $49$, and refinement through $17$ leaves none (Theorem~\ref{qt:nu:thm:nonunitcomplete}). The cubic branch has a separate final factorization and norm check (Section~\ref{cr:section}). These computations enumerate bounded integer exponents, with no need to list the possible number fields in advance.

\paragraph{3. Use the common unit lattice to prove effective height bounds.}
For an independent four-unit configuration, simultaneous complementation allows the coupling to be positive. The equations then take the form
\[
 y=x^{-a}(1-x)^u,\qquad 1-y=x^{-c}(1-x)^j,\qquad a,u,c,j>0,\qquad aj-uc=1,
\]
with rational exponents. If $U=\ell(x)$ and $V=\ell(1-x)$ are the vectors of logarithmic absolute values at the embeddings of $K$, and $U',V'$ are the corresponding vectors for $y$, then
\[
 U'=-aU+uV,\qquad V'=-cU+jV,\qquad U'\wedge V'=-U\wedge V.
\]
Both pairs therefore have the same logarithmic area and the same saturated rank-two unit lattice. Let $I$ be the index of the subgroup generated by either complementary pair in this saturated lattice. The lattice structure gives the additional arithmetic restriction
\[
 a,u,c,j\in I^{-1}\Z_{>0}.
\]
This discreteness is central to the height argument. Nonzero cross areas have a positive lower bound supplied by the unit regulator, or by the covolume of the saturated plane in the totally real quartic case. The latter field has unit rank three, so the estimate concerns the rank-two plane inside its unit lattice. Comparing these lower bounds with the upper bounds furnished by the complementary equations first controls one coordinate and then its partner. This yields explicit coefficient-height and index bounds. Odd and even index are treated separately; even index supplies square relations that lead to smaller families of generator polynomials.

\paragraph{4. Turn a height bound into an explicit polynomial list.}
An ordinary unit $x$ whose complement is also a unit has a monic integer minimal polynomial $f$ satisfying $f(0),f(1)\in\{\pm1\}$. For example, in degree four it has the form
\[
 f(T)=T^4+c_3T^3+c_2T^2+c_1T+e,\qquad c_1=\varepsilon-1-e-c_3-c_2,\qquad e,\varepsilon\in\{\pm1\}.
\]
Once the coefficient height is bounded, only finitely many integer choices remain. Exact factorization and root-counting tests enforce irreducibility, the field signature and the required locations of the real conjugates. When a coordinate or its complement is a square $w^2$ with $w\in K\cap(0,1)$, the elements $w$ and $1\pm w$ are units, so the minimal polynomial $g$ of $w$ satisfies $g(0),g(1),g(-1)\in\{\pm1\}$. In mixed signature, product-square partners are enumerated in a separate polynomial family. The resulting lists are lists of candidate polynomials. A polynomial on such a list need not belong to any Nahm solution.

\paragraph{5. Prune pairs by logarithmic compatibility and lattice arithmetic.}
The main reduction comes from pairing the polynomial lists. A genuine pair must have opposite logarithmic exterior products under the embeddings of its common field. Equality of their areas gives a cheap first test; the full exterior-product identity and the positivity of $a,u,c,j$ give stronger tests. All eligible alignments of the roots are included. It is safe at this stage to include alignments that do not come from a field isomorphism, since this only enlarges the candidate set. Every actual common field supplies one of the tested alignments.

The logarithmic tests are rigorous exclusion tests. Roots are isolated in rational intervals, logarithms are enclosed by explicit series with proved remainder bounds, and arithmetic is rounded outward. A required identity fails only when its interval is separated from zero; an interval containing zero remains undecided. The index bound then restricts the rational exponents. Cramer's rule gives intervals for $a,u,c,j$, and each permitted $I$ is tested for positive integers $Ia,Iu,Ic,Ij$ satisfying
\[
 (Ia)(Ij)-(Iu)(Ic)=I^2.
\]
The mixed-signature quartic odd-index branch (Section~\ref{qt:mixed-odd}) illustrates the strength of these restrictions. Its proved bounds (Theorem~\ref{qt:mo:thm:height}) are
\[
 H(f_x)<700,\qquad H(f_y)<2700,\qquad I\le109\quad(I\text{ odd}),
\]
where $H$ is the maximum absolute coefficient. Its reduction is
\begin{center}
\begin{tabular}{lr}\toprule
Stage&Count\\\midrule
Initial small-coordinate polynomials&11{,}362\\
Initial partner polynomials&107{,}798\\
Pairs reaching determinant comparison&17{,}941\\
Configurations remaining after coefficient-sign tests&4{,}585\\
Exact reconstructions after the index restriction&10\\
Classes after coordinate interchange&5\\\bottomrule
\end{tabular}
\end{center}
The exponent intervals are narrower than $1/109^2$, so each contains at most one rational of denominator at most $109$. The ten surviving configurations are checked in exact number-field arithmetic. In the totally real quartic branches (Section~\ref{qt:totally-real}), the exterior-product tests instead exclude every candidate pair. The even-index mixed-signature branch (Section~\ref{qt:mixed-even}) requires its own square-class enumeration and contributes three additional arithmetic classes.

\paragraph{6. Reconstruct exact points, then solve for every rational shift.}
Passing the preceding screens is still only a necessary condition. For each survivor we test whether the required positive partner belongs to the distinguished embedded field. Candidates that fail this test are discarded: in the mixed-signature even-index quartic branch, seven of the twelve remaining configurations fail at this stage. For those that pass, we express the partner algebraically in a chosen generator and verify the cleared Nahm equations by exact polynomial reduction. Positivity at the selected real embedding recovers the original rational-power equations. In the independent branch these relations determine the matrix; the dependent branch retains its explicitly parametrized families. Simultaneous complementation has been used as a symmetry of the arithmetic equations, so both matrix orientations are tested separately for modularity.

For each remaining matrix or family, Section~\ref{cls:sec:method} expresses the necessary radial coefficient conditions in a rational basis of $K$. The condition $\lambda_1\in\Q$ and the required vanishing of higher coefficients become polynomial equations over $\Q$ in the shifts and any remaining family parameter; the scalar shift is then $c_0=\lambda_1-1/12$. Exact elimination treats all rational solutions, with no imposed height or denominator bound on the shifts. In the dependent higher-degree branches, a certified nonzero Bloch--Wigner regulator supplies an additional modular obstruction, and the families surviving that test are handled by the coefficient equations. For the quartics, the eight arithmetic matrix classes reduce in this way to one class and its inverse, each with two shifts (Section~\ref{qt:shifts}). Finally, published $q$-series identities establish modularity of the surviving data.

The completeness argument thus has a specific structure. Every modular datum enters one of the stated arithmetic branches. The proved bounds place it in an explicitly enumerated set or parameter family. Each pruning step preserves every possible solution, and the survivors receive exact algebraic and coefficient checks. The final identities supply the sufficiency that the necessary conditions alone do not provide.

%% file: sections/arithmetic_preliminaries.tex
\section{Arithmetic preliminaries and the low-degree split}\label{low:sec:arithmetic}
Throughout this part,
\begin{equation}\label{low:N}
1-x=x^\alpha y^\beta,\qquad 1-y=x^\beta y^\gamma,
\qquad A=\begin{pmatrix}\alpha&\beta\\\beta&\gamma\end{pmatrix}\in M_2(\Q),\quad A>0,
\end{equation}
and $x,y\in(0,1)$. Rational powers are positive real powers. Algebraic identities at other places are obtained by clearing denominators. Write $K=\Q(x,y)$, $d=[K:\Q]$, and
\[
\dd(u)=[u]\wedge[1-u]\in\bigwedge^2(K^\times\otimes_\Z\Q),\qquad
P_u=\operatorname{span}_{\Q}([u],[1-u]).
\]
The notation $[u]$ here denotes a multiplicative tensor class, not a Bloch symbol. Symmetry gives $\dd(x)+\dd(y)=0$. Call a configuration independent if $\dd(x)\ne0$, and a four-unit configuration if all four coordinates $x,1-x,y,1-y$ are ordinary algebraic units. The complementary transformations are
\begin{equation}\label{low:symmetries}
(A,x,y)\longmapsto(A^{-1},1-x,1-y),\qquad
(A,x,y)\longmapsto(JAJ,y,x),\quad J=\begin{pmatrix}0&1\\1&0\end{pmatrix}.
\end{equation}
They preserve \eqref{low:N}, positive definiteness, and joint degree.

\begin{theorem}[Arithmetic classification in degrees one and two]\label{low:global}
For the configurations \eqref{low:N}:
\begin{enumerate}
\item Degree one consists exactly of $x=y=1/2$ and
\[
A=\begin{pmatrix}1-b&b\\b&1-b\end{pmatrix},\qquad b\in\Q,\quad b<1/2.
\]
\item The independent degree-two configurations, up to \eqref{low:symmetries}, are
\begin{equation}\label{low:twopairs}
\begin{array}{c|c|c}
x&y&A\\\hline
\sqrt2-1&1-1/\sqrt2&\begin{pmatrix}2&-1\\-1&1\end{pmatrix}\\[7pt]
\rho^3&\rho^2/2&\begin{pmatrix}4/3&-1\\-1&1\end{pmatrix}
\end{array}\qquad \rho=(\sqrt5-1)/2.
\end{equation}
Their Rogers sums are $\pi^2/8$ and $\pi^2/12$, respectively. Their swaps and simultaneous reflections give eight matrices. The dependent quadratic configurations are listed in Theorem~\ref{low:quadratic}.
\item All dependent configurations of joint degree at most four with nonzero coupling belong to the explicit finite list of one-parameter families in Proposition~\ref{low:scalar}. Diagonal matrices belong to a list of at most $169$ candidates before the joint-degree filter.
\end{enumerate}
The independent cubic and quartic branches are classified in Parts~\ref{part:cubic} and~\ref{part:quartic}.
\end{theorem}
\begin{proof}
The first two assertions are proved in Theorem~\ref{low:quadratic}.
Proposition~\ref{low:scalar} gives the dependent families and the thirteen
possible scalar exponents; their ordered pairs give the stated upper
bound for diagonal matrices.
\end{proof}

The degree restrictions in the following arithmetic statements are supplied,
for modular Nahm sums, by Theorem~\ref{deg:thm:main}.
The proof in this part uses only the symmetric Nahm equations, positive
definiteness and the degree restriction. Weighted relations
$m\dd(x)+\dd(y)=0$ with $m>1$ are not used to infer symmetric completeness.

The positive solution of \eqref{low:N} is unique. In variables $v_i=-\log x_i>0$ it is the critical point of
\[
\tfrac12v^tAv+\operatorname{Li}_2(e^{-v_1})+\operatorname{Li}_2(e^{-v_2}).
\]
The Hessian is $A+\operatorname{diag}(1/(e^{v_i}-1))>0$. Coercivity and the inward negative derivative on each boundary face give existence as well. This observation will also justify symmetry arguments below.

\section{Scalar solutions and the dependent families}
\begin{proposition}\label{low:scalar}
The scalar solutions $1-u=u^c$, $0<u<1$, $c\in\Q_{>0}$, of degree at most four are exactly those in the following table. Each defining root is the unique positive root in $(0,1)$.
\begin{center}
\begin{tabular}{lll}
\toprule
Defining root&Possible $u$&Corresponding $c$\\\midrule
Rational&$1/2$&$1$\\
$v^2+v-1=0$&$v,v^2$&$2,1/2$\\
$v^3+v-1=0$&$v,v^3$&$3,1/3$\\
$v^3+v^2-1=0$&$v^2,v^3,v,v^5$&$3/2,2/3,5,1/5$\\
$v^4+v-1=0$&$v,v^4$&$4,1/4$\\
$v^4+v^3-1=0$&$v^3,v^4$&$4/3,3/4$\\\bottomrule
\end{tabular}
\end{center}
For each nonrational row, let $\mathcal P_v$ be the corresponding list of ordered complementary exponent pairs:
\[
\{(1,2),(2,1)\},\quad \{(1,3),(3,1)\},\quad
\{(2,3),(3,2),(1,5),(5,1)\},\quad
\{(1,4),(4,1)\},\quad\{(3,4),(4,3)\}.
\]
Every dependent configuration with $\beta\ne0$ and $d\le4$ is either the rational family in Theorem~\ref{low:global}, or is obtained by choosing $(r,s),(p,q)\in\mathcal P_v$ and setting
\begin{equation}\label{low:dependentfamily}
x=v^r,\quad y=v^p,\qquad
A_b=\begin{pmatrix}(s-pb)/r&b\\b&(q-rb)/p\end{pmatrix},\qquad
b\in\Q,\quad b<\frac{sq}{rs+pq}.
\end{equation}
The value $b=0$ is valid too. The lists give at most $33$ family templates, including the rational one, before identifying symmetries.
\end{proposition}
\begin{proof}
For a scalar solution write $c=n/m$ in lowest terms. Bezout's identity gives $v\in\Q(u)$ with $u=v^m$, $1-u=v^n$. If $m=n$ then $u=1/2$. Otherwise $v$ is a unit satisfying $v^m+v^n=1$.

We use Ljunggren's trinomial theorem \cite{Ljunggren}, in the form stated in \cite[Theorem 1]{KR}: $Z^n+Z^m-1$ has one irreducible noncyclotomic factor and all its other factors are cyclotomic. For $\gcd(m,n)=1$, a cyclotomic root must be a primitive sixth root. Indeed its two powers summing to one are $e^{i\pi/3}$ and $e^{-i\pi/3}$, and Bezout's identity puts the root itself among the sixth roots. The factor $Z^2-Z+1$ is simple, as its two derivative terms have different absolute values $n,m$. The noncyclotomic factor therefore has degree $n$ or $n-2$. Degree at most four implies $n\le6$. Factoring these finitely many primitive trinomials gives the table; the only additional pair comes from
\[
Z^5+Z-1=(Z^2-Z+1)(Z^3+Z^2-1).
\]

For a dependent solution of \eqref{low:N} with $\beta\ne0$, all four coordinates generate a positive multiplicative group of rank one. Write them as $v^r,v^s,v^p,v^q$ with $v\in(0,1)$ and $\gcd(r,s,p,q)=1$. Then $K=\Q(v)$. If either pair has equal exponents, its coordinate is $1/2$. An unequal pair would make $v$ a unit, contradicting a positive power of $v$ being $1/2$. Thus both are equal and $x=y=1/2$, or both pairs are unequal and $v$ is a unit.

For an unequal pair, the roots of $Z^r+Z^s-1$ of smallest absolute value are exactly $v\zeta$ with $\zeta^{\gcd(r,s)}=1$, by equality in the triangle inequality. All lie in its unique noncyclotomic factor, the minimal polynomial of $v$. Applying this to both pairs shows $\gcd(r,s)=\gcd(p,q)$. Their common value divides all four exponents, so it is one. The scalar table now supplies precisely the stated lists. Substitution gives \eqref{low:dependentfamily}; its determinant is
\[
\frac{sq-(rs+pq)b}{rp}.
\]
The displayed upper bound on $b$ makes both diagonal entries positive as well. This proves the claim.
\end{proof}

In particular, the $13$ scalar exponents give at most $13^2=169$ diagonal matrices; one subsequently checks the joint field degree. For a reflection pair $(v^r,v^s)$, \eqref{low:dependentfamily} becomes
\begin{equation}\label{low:reflectionfamily}
A_b=\begin{pmatrix}(s/r)(1-b)&b\\b&(r/s)(1-b)\end{pmatrix},\qquad b<1/2.
\end{equation}
These families explain why even a complete list of torsion points need not give a finite list of rational matrices.

\section{Finite places force a unit in the common plane}
The following observation removes the unit hypothesis from the quadratic classification within \eqref{low:N}.
\begin{lemma}\label{low:unitforcing}
If \eqref{low:N} is independent and nonunit, the plane $P_x=P_y$ contains an ordinary unit of infinite order.
\end{lemma}
\begin{proof}
Independence gives $\beta\ne0$, and the nonzero wedge relation gives $P_x=P_y$. Put $\Delta=\alpha\gamma-\beta^2$. At a finite prime write $u=\operatorname{ord}(x)$, $v=\operatorname{ord}(y)$. The valuation rule for $1-x$, together with positive definiteness, leaves exactly the following nonzero cases:
\begin{equation}\label{low:loci}
\begin{array}{c|l}
u<0,\ v<0&(A-I_2)(u,v)^t=0\quad(E_1)\\
u<0,\ v=0&\alpha=1,\ \beta<0\quad(E_2)\\
u=0,\ v<0&\gamma=1,\ \beta<0\quad(E_3)\\
u>0,\ v<0&\Delta=\alpha,\ \beta>0\quad(E_4)\\
u<0,\ v>0&\Delta=\gamma,\ \beta>0\quad(E_5).
\end{array}
\end{equation}
For example, the fourth row follows from $\alpha u+\beta v=0$ and $\beta u+(\gamma-1)v=0$. Both valuations positive would give $A(u,v)^t=0$, and a positive valuation with the other zero contradicts a positive diagonal entry. If $u=v=0$, the complementary valuations vanish too.

The loci $E_2$--$E_5$ reduce under \eqref{low:symmetries} to
\begin{equation}\label{low:canonicalearly}
A_c=\begin{pmatrix}1&-h\\-h&\tau\end{pmatrix},\qquad h>0,\quad\tau>h^2.
\end{equation}
If $\tau=1$, uniqueness and symmetry give $x=y$ and hence $\dd(x)=0$. Thus $\tau\ne1$. At every prime the first valuation equation is $\operatorname{ord}(1-x)=u-hv$. If $u>0$, it forces $v=u/h>0$, and the second equation gives $(\tau-h^2)v=0$, a contradiction. If $u=0$, then $v\le0$, and $v<0$ would force $\tau=1$. If $u<0$, then $v=0$. Consequently $y$ and $x/(1-x)$ are ordinary units. The unit $y\in(0,1)$ has infinite order and lies in the plane.

On $E_1$, the condition is $\Delta=\alpha+\gamma-1$. It cannot intersect any of $E_2$--$E_5$ when $\beta\ne0$: for example $\Delta=\alpha$ forces $\gamma=1$ and then $\beta=0$. Thus every nonzero finite valuation has $u,v<0$, and both odds ratios $x/(1-x),y/(1-y)$ are units. The first has infinite order, since equality to one would give $x=1/2$ and $\dd(x)=0$. Transformations \eqref{low:symmetries} preserve the common plane, completing the proof.
\end{proof}

\section{Complete classification in degrees one and two}\label{low:sec:quadratic}
\begin{theorem}\label{low:quadratic}
The degree-one assertion and independent degree-two assertion of Theorem~\ref{low:global} hold. The remaining points of joint degree two are exactly:
\begin{enumerate}
\item $x=\rho^r$, $y=\rho^s$, $r,s\in\{1,2\}$, with
\begin{equation}\label{low:goldfamily}
A=\begin{pmatrix}(3-r-sb)/r&b\\b&(3-s-rb)/s\end{pmatrix},\qquad
b\in\Q,\quad b<\frac{(3-r)(3-s)}4;
\end{equation}
\item $(x,y)=(1/2,\rho^r)$, $r=1,2$, with $A=\operatorname{diag}(1,(3-r)/r)$, and their coordinate exchanges.
\end{enumerate}
\end{theorem}
\begin{proof}
The rational boundary classification in Appendix~\ref{low:rationalapp} shows that nonzero equal-weight rational boundaries cancel only for $y=1-x$. Such an independent pair cannot satisfy \eqref{low:N}: independence of $x,1-x$ in its first equation would force $\alpha=0$, $\beta=1$, contradicting positive definiteness. A zero-boundary rational argument is $1/2$, so the rational matrix family follows directly.

In a real quadratic field the unit rank is one, so independent pairs cannot be four-unit configurations. Lemma~\ref{low:unitforcing} supplies an ordinary unit in their plane. Both arguments are irrational. Indeed, if $x$ is rational and $\dd(x)\ne0$, the common plane gives rational powers $y^N,(1-y)^N\in\Q$. A quadratic conjugate of an irrational real $y$ would then be $-y$, whereas the second power condition would require $1+y=\pm(1-y)$, which is impossible. The equal-weight shifted-unit enumeration in Proposition~\ref{low:quadraticenum} therefore applies. It yields precisely the two pairs \eqref{low:twopairs}, or reflection, which was just excluded for independent solutions of \eqref{low:N}.

For each listed pair, exact substitution verifies its matrix, and multiplicative independence of $x,y$ makes this matrix unique. Reflection inverts it and exchange conjugates it by $J$. The four matrices from each representative are distinct. Their determinants are $1$ and $1/3$.

If the boundary is zero, Proposition~\ref{low:scalar} gives only $1/2,\rho,\rho^2$ in a real quadratic field. Nonzero coupling puts all four coordinates in one rank-one group, so a rational argument cannot be mixed with a golden-ratio unit. Thus the mixed pairs are diagonal. The other pairs give \eqref{low:goldfamily}; its determinant is $((3-r)(3-s)-4b)/(rs)$.
\end{proof}

For completeness, let $L(t)=\operatorname{Li}_2(t)+\frac12\log t\log(1-t)$ and $C_0=\pi^2/6$. Reflection and the positive five-term relation give
\begin{equation}\label{low:rogers}
L(t)+L(1-t)=C_0,\qquad
L(a)+L(b)=L(ab)+L\!\left(\frac{a(1-b)}{1-ab}\right)+L\!\left(\frac{b(1-a)}{1-ab}\right).
\end{equation}
Duplication at $z=1/\sqrt2$ gives $2L(z)=L(1/2)+2L(\sqrt2-1)$, so the first sum in \eqref{low:twopairs} is $3C_0/4$. The five-term relation at $(\rho^2,1/2)$ gives $L(\rho^3)+L(\rho^2/2)=C_0/2$. Also
\[
L(1/2)=C_0/2,\qquad L(\rho)=3C_0/5,\qquad L(\rho^2)=2C_0/5,
\]
which evaluates all dependent quadratic sums. These are exact certificates, rather than numerical tests.

\subsection{The quadratic arithmetic enumeration}
We include the equal-weight portion of the shifted-unit argument to specify its finite proof. The broader weighted classification is treated in \cite{quadraticpaper}.
A primitive prime-ideal divisor of $e^n-1$ is a prime ideal at which $e$
has multiplicative order exactly $n$. The following lemma supplies the
cyclotomic estimate used in the proof.

\begin{lemma}\label{low:quad-cyclotomic}
Let $K$ be a real quadratic field, $e>1$ a unit of $K$, and $n\ge7$.
If $e^n-1$ has no primitive prime-ideal divisor, then $(\Phi_n(e))$
divides the ideal generated by the product of the primes dividing $n$, and
\[
 e^{\varphi(n)}<81\,|N(\Phi_n(e))|\le81n^2.
\]
\end{lemma}
\begin{proof}
Let $\mathfrak p$ be a prime ideal dividing $\Phi_n(e)$, of residue
characteristic $p$, and let $m$ be the order of $e$ modulo $\mathfrak p$.
Write $n=n'p^s$ with $p\nmid n'$. Modulo $\mathfrak p$ one has
$\Phi_n=\Phi_{n'}^{\varphi(p^s)}$, and the roots of the separable
polynomial $\Phi_{n'}$ are exactly the elements of order $n'$; hence
$m=n'$. By hypothesis $m<n$, so $s\ge1$ and $p\mid n$.

Put $\beta=e^{p^{s-1}}$ and $\gamma=\beta^m=e^{n/p}$.
Both $\beta$ and $\beta^p$ have order $m$ modulo $\mathfrak p$, so
$\mathfrak p$ divides no $\Phi_d(\beta)$ or $\Phi_d(\beta^p)$ with
$d\mid m$, $d<m$. Hence $\Phi_m(\beta)$ and $\Phi_m(\beta^p)$ have
the same $\mathfrak p$-adic valuations as $\gamma-1$ and $\gamma^p-1$.
From $\Phi_n(X)=\Phi_m(X^{p^s})/\Phi_m(X^{p^{s-1}})$ we obtain
\[
 v_{\mathfrak p}(\Phi_n(e))
 =v_{\mathfrak p}(\gamma^p-1)-v_{\mathfrak p}(\gamma-1)
 =v_{\mathfrak p}\left(\sum_{k=1}^p\binom pk\pi^{k-1}\right),
 \qquad \pi=\gamma-1.
\]
Let $w=v_{\mathfrak p}(\pi)\ge1$ and
$\varepsilon=v_{\mathfrak p}(p)\in\{1,2\}$.
The terms with $2\le k\le p-1$ have valuation greater than
$\varepsilon$, the first term is $p$, and the last is $\pi^{p-1}$.
Hence $v_{\mathfrak p}(\Phi_n(e))\le\varepsilon$ unless
$(p-1)w=\varepsilon$, that is, unless $p=2$ and $w=\varepsilon$,
or $p=3$, $\varepsilon=2$ and $w=1$. Neither case occurs for $n\ge7$.

If $p=2$, then $m$ divides $N\mathfrak p-1\in\{1,3\}$, so
$n\ge7$ forces $s\ge2$. Writing $\gamma=\delta^2$ with
$\delta=e^{n/4}\equiv1\pmod{\mathfrak p}$ gives
$w=v_{\mathfrak p}(\delta-1)+v_{\mathfrak p}(\delta+1)\ge2$,
since both factors vanish modulo $\mathfrak p$.
If moreover $\varepsilon=2$, the residue field is $\mathbb F_2$, so
$m=1$ and $s\ge3$. Then $\delta$ is the square of an element congruent
to one, so $v_{\mathfrak p}(\delta-1)\ge2$, hence
$v_{\mathfrak p}(\delta+1)=v_{\mathfrak p}((\delta-1)+2)\ge2$ and
$w\ge4$. If $p=3$ and $\varepsilon=2$, the residue field is
$\mathbb F_3$, so $m\le2$ and $s\ge2$. With $\gamma=\delta^3$ and
$\delta=e^{n/9}\equiv1\pmod{\mathfrak p}$, one gets
$w=v_{\mathfrak p}(\delta-1)+v_{\mathfrak p}(1+\delta+\delta^2)\ge2$,
since $1+\delta+\delta^2\equiv3\equiv0\pmod{\mathfrak p}$.

Thus $v_{\mathfrak p}(\Phi_n(e))\le v_{\mathfrak p}(p)$ for every
prime ideal dividing $\Phi_n(e)$, which proves the divisibility.
Taking norms gives $|N(\Phi_n(e))|\le n^2$.
For the lower bound, $\Phi_n$ is reciprocal and the conjugate of $e$
is $\pm e^{-1}$, so
\[
 |N(\Phi_n(e))|
 =e^{\varphi(n)}\,|\Phi_n(e^{-1})\Phi_n(\pm e^{-1})|.
\]
For $|t|<1$, the identity
$\Phi_n(t)=\prod_{d\mid n}(1-t^d)^{\mu(n/d)}$ gives
$|\Phi_n(t)|\ge\prod_{d\ge1}(1-|t|^d)$, since each factor
$|1-t^d|$ lies between $1-|t|^d$ and $(1-|t|^d)^{-1}$.
A unit greater than one in a real quadratic field is at least the golden
ratio $g$. Moreover $\prod_{d\ge1}(1-g^{-d})>1/9$: the first six
factors have product greater than $0.13$, and the remaining ones have
product at least $1-\sum_{d\ge7}g^{-d}=1-g^{-5}>0.9$.
Therefore $e^{\varphi(n)}<81|N(\Phi_n(e))|$.
\end{proof}

\begin{proposition}\label{low:quadraticenum}
Let $K$ be real quadratic, $x,y\in K\cap(0,1)$ irrational, with $\dd(x)+\dd(y)=0\ne\dd(x)$, and suppose $P_x$ contains an ordinary unit of infinite order. Apart from reflection, the two pairs \eqref{low:twopairs} are the complete list, up to exchange and simultaneous reflection.
\end{proposition}
\begin{proof}
Let $e>1$ be the fundamental unit. The finite-divisor map on $P_x$ has rank one: its kernel contains a unit, while the rank-one ordinary unit group cannot contain the whole plane. Consequently one of the six anharmonic transforms of each argument is a unit. Choose representatives $U=\sigma e^a$, $V=\tau e^b$ of absolute value greater than one. The integral nonunits $1-U,1-V$ have proportional, nonzero principal-ideal divisors. If $U=V$, opposite equal boundaries give reflection. We may therefore take distinct units.

Here are the cutoff details. Let $e'$ denote the conjugate of $e$.
Flatters' primitive-divisor theorem
\cite[Theorem 1.4]{Flatters} concerns rational primes that divide
$N(e^n-1)=(e^n-1)((e')^n-1)$ but no earlier term. Such a prime lies
below a prime ideal dividing $e^n-1$ or its conjugate; in the second case
the conjugate ideal divides $e^n-1$. In either case this ideal is a
primitive prime-ideal divisor of $e^n-1$. For norm $+1$ his cutoff is $12$;
for norm $-1$ it is $24$ except $n\equiv2\pmod4$. In that exception write
$n=2j>48$, $j$ odd. A primitive prime at $j$ is odd, and conjugation sends
$e$ to $-e^{-1}$, turning its order $j$ into $2j$ at the conjugate prime
ideal. Hence $e^n-1$ has a primitive prime-ideal divisor for $n>48$.
Further, for $n\ge7$ without a primitive prime-ideal divisor,
Lemma~\ref{low:quad-cyclotomic} gives
\begin{equation}\label{low:quadcutoffineq}
e^{\varphi(n)}<81n^2.
\end{equation}

The index of $1-e^a$ is $a$, and of $1+e^a$ is $2a$. If two proportional shifts have different indices, the larger has no primitive divisor. Thus both indices are at most $48$; if the larger exceeds six, \eqref{low:quadcutoffineq} gives $e<13$, since $81n^2<13^{\varphi(n)}$ for $7\le n\le48$. If their indices agree, the shifts are $1-\epsilon^2$, $1+\epsilon$, with $\epsilon=e^j$. When $1-\epsilon$ is a unit, its norm gives $\operatorname{Tr}(\epsilon)\le3$ for norm $+1$, and at most one for norm $-1$. Otherwise its divisor is proportional to that of $1+\epsilon$; their support lies above two. In norm $+1$, the norms $T-2,T+2$ are powers of two differing by four, forcing $T=6$. In norm $-1$, their equal norms imply equality of ideals, which then divide $(2)$, giving $T\le4$. Thus $\epsilon<6$, $j\le3$, and the common index is at most six.

It remains to exclude $e>13$ with indices at most six. For norm $+1$ let $T=\operatorname{Tr}(e)\ge14$ and put
\[
\beta_1=T-2,\ \beta_2=T+2,\ \beta_3=T,\ \beta_4=T+1,\ \beta_5=T-1,\ \beta_6=T^2+T-1.
\]
The nine absolute norms, for $1-e^j$ ($1\le j\le6$) and $1+e^j$ ($1\le j\le3$), are
\[
\beta_1,\ \beta_1\beta_2,\ \beta_1\beta_4^2,\ \beta_1\beta_2\beta_3^2,
\ \beta_1\beta_6^2,\ \beta_1\beta_2\beta_4^2\beta_5^2,\ \beta_2,\ \beta_3^2,\ \beta_2\beta_5^2.
\]
The only nontrivial pairwise gcd bounds are $4$ for $(1,2)$, $2$ for $(1,3),(2,3),(4,5)$, $3$ for $(1,4),(2,5)$, and $5$ for $(1,6)$; all other gcds are one. A dependence involving $\beta_1\beta_6^2$ would force $\beta_6$ to be a power of five, whereas $5\mid\beta_6$ implies $v_5(\beta_6)=1$, impossible here. Cancellation in the other pairs gives the following complete necessary trace list:
\begin{center}\begin{tabular}{lc}\toprule
Pair of norms&Necessary $T$\\\midrule
$(\beta_1,\beta_2)$, $(\beta_1,\beta_1\beta_2)$, $(\beta_1\beta_4^2,\beta_2\beta_5^2)$&$6$\\
$(\beta_2,\beta_1\beta_2)$&$3,6$\\
$(\beta_1,\beta_3^2)$, $(\beta_2,\beta_1\beta_2\beta_3^2)$, $(\beta_2\beta_5^2,\beta_1\beta_2\beta_3^2)$&$4$\\
$(\beta_1\beta_2,\beta_2\beta_5^2)$&$4,10$\\
$(\beta_2\beta_5^2,\beta_1\beta_2\beta_4^2\beta_5^2)$&$3$\\\bottomrule
\end{tabular}\end{center}
None is at least $14$. The elementary comparisons use powers of two differing by $2,3,4$ and powers of three differing by $2,3$. In particular $3^v-2^u=1$ gives $(2^u,3^v)=(2,3),(8,9)$, by reducing modulo eight for odd $v$ and factoring for even $v$.

For norm $-1$, $T\ge13$. Up to powers the nine norms lie in the classes
\[
T,\quad T(T^2+3),\quad T^2+4,\quad T^2(T^2+4),\quad T(T^4+5T^2+5).
\]
The three factors $T^2+4,T^2+3,T^4+5T^2+5$ are pairwise coprime; their gcds with $T$ divide $4,3,5$, respectively. The latter two have valuation one at the indicated prime when it divides $T$. These facts exclude dependence between distinct classes for $T\ge13$. Within the repeated class from $1-e,1+e,1-e^2$ (and the class with $e$ replaced by $e^3$), proportionality forces the first two ideals equal and dividing $(2)$, giving $\operatorname{Tr}(e)\le4$ or $\operatorname{Tr}(e^3)\le4$, again impossible.

It therefore suffices to check $e<13$, $a\le48$ for positive sign and $a\le24$ for negative sign. The following exact recurrence specifies the entire check, including nonfundamental units:
\begin{equation}\label{low:quadrecurrence}
e=\frac{T+\sqrt{T^2-4\nu}}2,\quad
\begin{cases}\nu=1,&3\le T\le13,\\\nu=-1,&1\le T\le12,\end{cases}
\quad T_0=2,\ T_1=T,\ T_n=TT_{n-1}-\nu T_{n-2}.
\end{equation}
The shifted norm is $N_{\sigma,n}=|1-\sigma T_n+\nu^n|$. Discard values at most one, write each as $D^k$ with $D$ not a proper power, and group equal bases $D$. Besides the repeated norm-$-1$ groups $e^j,-e^j,e^{2j}$ for odd $1\le j\le23$, the only groups are
\begin{center}\begin{tabular}{ccll}\toprule
$\nu$&$T$&Units&Shift norms\\\midrule
$1$&$3$&$-e,e^2,-e^5$&$5,5,125$\\
$1$&$4$&$e,-e^2$&$2,16$\\
$1$&$6$&$e,-e,e^2$&$4,8,32$\\
$-1$&$1$&$-e^2,e^4,-e^{10}$&$5,5,125$\\
$-1$&$2$&$e,-e,e^2,-e^2,e^4$&$2,2,4,8,32$\\\bottomrule
\end{tabular}\end{center}
For norms $D^k,D^\ell$, proportionality of the ideals is equivalent to the integrality of
\[
\operatorname{Tr}\!\left(\frac{(1-U)^\ell}{(1-V)^k}\right),
\]
since this quotient has norm $\pm1$. All operations use $e^2=Te-\nu$ and $\operatorname{Tr}(a+be)=2a+Tb$. This check reduces $446$ unordered norm pairs to $26$ ideal pairs. In fundamental units they form four groups:
\begin{center}\begin{tabular}{lll}\toprule
Fundamental unit&Units&Boundary sizes\\\midrule
$1+\sqrt2$&$e,-e,e^2,-e^2,e^4$&$1,1,4,6,20$\\
$2+\sqrt3$&$e,-e^2$&$1,8$\\
$(1+\sqrt5)/2$&$e^3,-e^3,e^6$&$6,6,24$\\
$(1+\sqrt5)/2$&$-e^2,e^4,-e^{10}$&$2,4,30$\\\bottomrule
\end{tabular}\end{center}
Here the boundary size is $nk$ for $U=\pm e^n$ with shifted norm $D^k$. For a unit $U$ the interval representative is $t(U)=1-U^{-1}$ if $U>1$, and $t(U)=(1-U)^{-1}$ if $U<-1$; it satisfies $\dd(t(U))=[U]\wedge[1-U]$. Opposite equal boundaries can therefore use only the repeated sizes in the first and third rows. The pairs $(t(U),1-t(V))$ give precisely \eqref{low:twopairs}, up to the stated symmetries. This completes the finite proof.
\end{proof}

%% file: sections/quartic_entry.tex
\section{Canonical entry for nonunit points}\label{qt:entrysec}
The following reduction uses the joint field of the distinguished point.
It applies in every degree greater than one and will be used in both the
cubic and quartic classifications. An ordinary unit means a unit in the
ring of integers of that field.

\begin{lemma}[Canonical nonunit entry]\label{qcov:entry}
Let $A=\bigl(\begin{smallmatrix}\alpha&\beta\\\beta&\gamma\end{smallmatrix}\bigr)$
be rational, symmetric and positive definite, with $\beta\ne0$. Let
$(x,y)\in(0,1)^2$ satisfy its Nahm equations, put $K=\Q(x,y)$, and assume
$[K:\Q]>1$. If at least one of $x,1-x,y,1-y$ is not an ordinary unit,
then coordinate exchange and, if necessary, the transformation
$(A,x,y)\mapsto(A^{-1},1-x,1-y)$ put the datum in one of the following forms.
There exist an ordinary unit $z\in K\cap(0,1)$ and integers
$a,b>0$, $T\ne0$ such that, with $g=\gcd(a,b)$, $M=a/g$, $N=b/g$,
\[
 \Q(z)=K,\qquad \gcd(M,N)=\gcd(g,T)=1,
\]
and one of
\begin{align}
 (D)\qquad &(1-z^b)^N=z^T(1+z^a)^M, &&T>0,\label{qt:entryD}\\
 (P)\qquad &(1+z^b)^N=z^T(1+z^a)^M,
 &&gT>b^2\ \hbox{or}\ gT<-a^2.\label{qt:entryP}
\end{align}
The binomial factors on the two sides have equal, nonempty finite
prime-ideal support; in particular $1+z^a$ is a nonunit. The coordinates
and matrix are given by \eqref{qt:AD} or \eqref{qt:AP}, respectively.
Conversely these reconstructions satisfy the Nahm equations and are
positive definite under the displayed inequalities.
\end{lemma}
\begin{proof}

All rational powers at the distinguished point are positive real powers. Clearing their denominators permits valuations at every finite place. Write $\xi=v(x)$, $\eta=v(y)$ and $\Delta=\alpha\gamma-\beta^2$. The identities
\[
v(1-x)=\alpha\xi+\beta\eta,\qquad
v(1-y)=\beta\xi+\gamma\eta
\]
and the ultrametric rule for $1-x$ leave the following nonzero cases:
\[
\begin{array}{c|l}
\xi,\eta<0 &(A-I)(\xi,\eta)^t=0\quad(I)\\
\xi<0=\eta &\alpha=1,\ \beta<0\quad(II)\\
\eta<0=\xi &\gamma=1,\ \beta<0\quad(II')\\
\xi>0>\eta &\Delta=\alpha,\ \beta>0\quad(III)\\
\eta>0>\xi &\Delta=\gamma,\ \beta>0\quad(III').
\end{array}
\]
Indeed both positive valuations force $A(\xi,\eta)^t=0$, contradicting positive definiteness. One positive valuation and the other zero contradicts a positive diagonal entry. Both zero force both complementary valuations zero. The remaining five cases give the displayed equations directly; the weak coupling inequalities become strict because $\beta\ne0$.

Exchange and inversion preserve rationality, positive definiteness, $K$ and the four-unit property. Inversion sends the last two cases to a case with diagonal entry one and negative coupling. This is a symmetry of the Nahm \emph{equations}; modularity of the inverse sum is not required.

The intersections $II\cap II'$ and $III\cap III'$ are nonempty. The first has $\alpha=\gamma=1$ and is the $\tau=1$ direct case below; the second is its inverse. They are therefore covered by the same argument. Intersections requiring opposite signs of $\beta$ are empty.

\paragraph{Direct case.}
We may take
\[
 A=\begin{pmatrix}1&-h\\-h&\tau\end{pmatrix},
 \qquad h>0,\quad\tau>h^2.
\]
If $\tau=1$, uniqueness of the distinguished point gives $x=y$ and $1-x=x^{1-h}$, with $0<h<1$. Clearing the positive rational exponent produces a scalar unit: explicitly, if $1-h=n/m$ in lowest terms, a Bezout generator $v$ satisfies $x=v^m$, $1-x=v^n$ and $v^m+v^n=1$, with $m\ne n$. Thus all coordinates are units, contrary to the hypothesis. Hence $\tau\ne1$.

At every finite place, the first equation gives $v(1-x)=\xi-h\eta$. If $\xi>0$, then $\eta=\xi/h>0$ and the second equation contradicts $\tau>h^2$. If $\xi=0$, then $\eta\le0$, and $\eta<0$ would force $\tau=1$. If $\xi<0$, then $\eta=0$. Therefore $y$ and $s=x/(1-x)$ are units, and
\[
 v(1-y)=h\,v(1+s).
\]
Write $h=M/N$ with $M,N>0$ coprime. Since $s^N=y^M$, Bezout gives a positive unit $w\in K$ with $s=w^M$, $y=w^N$ and $w<1$. Then
\[
 W=\frac{(1-y)^N}{(1+s)^M}=w^{N^2\tau-M^2}
\]
is a unit, and its exponent is positive rational.

For clarity, if the exponent is $p/q$ in lowest terms, choose $r,s_0\in\Z$ with $rq+s_0p=1$. The element
\[
 z=w^rW^{s_0}=w^{1/q}\in K
\]
is a positive unit less than one. With $g=q$, $T=p$, $a=gM$, $b=gN$ we have $w=z^g$, $W=z^T$ and $\gcd(g,T)=1$. The $D$ equation follows. Its reconstruction is
\begin{equation}\label{qt:AD}
 x=\frac{z^a}{1+z^a},\quad y=z^b,\qquad
 A_D=\frac1{b^2}\begin{pmatrix}b^2&-ab\\-ab&gT+a^2\end{pmatrix}.
\end{equation}
Because $z\in K$ and $x,y\in\Q(z)$, the fields coincide. If $1+z^a$ were a unit, the valuation identity would make all four coordinates units, so the required nonunit hypothesis holds.

\paragraph{Plus case.}
It remains to treat case $I$. Its locus has no intersection with $II,II',III,III'$ when $\beta\ne0$. For example $\Delta=\alpha+\gamma-1$ together with $\Delta=\alpha$ forces $\gamma=1$ and then $\beta=0$. The other intersections follow by exchange or by setting the relevant diagonal entry to one. Hence every nonzero finite valuation is of type $I$, and both odds ratios are units.

The negative-kernel condition gives
\[
 A=I+t\begin{pmatrix}1&-h\\-h&h^2\end{pmatrix},
 \qquad h>0,\quad t\ne0.
\]
The logarithmic Nahm equations imply
\[
 \left(\frac{x}{1-x}\right)^h\frac{y}{1-y}=1.
\]
If either ratio is one, $x=y=1/2$, contradicting $[K:\Q]>1$. Exchange coordinates if necessary so that
\[
 U=\frac{x}{1-x}\in(0,1),\qquad
 V=\frac{1-y}{y}=U^h\in(0,1).
\]
Then $1+U=U^\kappa(1+V)^h$ with $\kappa=(1+t)/t\ne0$; the last inequality follows from $\alpha=1+t>0$. Writing $h=M/N$ gives positive units $U=w^N$, $V=w^M$, $w<1$, and
\[
 W=\frac{(1+U)^N}{(1+V)^M}=w^{N^2\kappa}.
\]
This is a unit: after writing the exponent as $p/q$, the equality $W^q=w^p$ forces every finite valuation of $W$ to vanish. Apply the same Bezout construction as before, now allowing $p<0$. It gives $z<1$, $g=q$, $T=p\ne0$ and $\gcd(g,T)=1$. The resulting $P$ equation reconstructs
\begin{equation}\label{qt:AP}
 x=\frac{z^b}{1+z^b},\quad y=\frac1{1+z^a},\qquad
 A_P=\frac1{gT-b^2}
 \begin{pmatrix}gT&-ab\\-ab&gT-b^2+a^2\end{pmatrix}.
\end{equation}
The denominator is nonzero since $\kappa\ne1$. Equivalently,
\[
 A_P=I+\frac{1}{gT-b^2}\binom b{-a}(b,-a).
\]
Its eigenvalues are $1$ and $(gT+a^2)/(gT-b^2)$, so positive definiteness is exactly $gT>b^2$ or $gT<-a^2$. Once again $\Q(z)=K$. The two plus factors have equal support; if that support were empty, all four coordinates would be units. Substitution in the two logarithmic Nahm equations proves the converse. This proves the lemma.
\end{proof}
The transformation in this lemma preserves the Nahm equations and the
joint field. It makes no assertion that modularity is preserved by
inverting a matrix. Later coefficient calculations treat the two matrix
orientations separately.

%% file: sections/coefficient_conditions.tex
\section{The coefficient conditions at a fixed Nahm point}\label{cls:sec:method}
Throughout the classification, $A=\left(\begin{smallmatrix}a&b\\b&d\end{smallmatrix}\right)$,
the linear shift is $\boldsymbol B=(B_1,B_2)^t$ and the scalar is $c_0$, as
in Section~\ref{sec:introduction}. As in Part~\ref{part:arithmetic},
$\rho=(\sqrt5-1)/2$ and $J=\left(\begin{smallmatrix}0&1\\1&0\end{smallmatrix}\right)$.
Exchange $(A,\boldsymbol B,c_0)\mapsto(JAJ,J\boldsymbol B,c_0)$ permutes the
summation indices and preserves modularity, so data are listed up to exchange.

Let $(x,y)$ be the distinguished point of $A$ and $K=\Q(x,y)$. If
$f_{A,\boldsymbol B,c_0}$ is modular, Lemma~\ref{deg:lem:modular-bridge} gives
\begin{equation}\label{cls:eq:conditions}
 \lambda_1=c_0+\tfrac1{12},\qquad \lambda_2=\lambda_3=\lambda_4=0 .
\end{equation}
At a fixed point and matrix, $\Psi_j=D^{3j}\lambda_j$ is a polynomial of
degree at most $j+1$ in $B_1,B_2$ with coefficients in $K$. Writing these
coefficients in a $\Q$-basis of $K$, the conditions $\lambda_1\in\Q$ and
$\lambda_2=\lambda_3=\lambda_4=0$ become finitely many polynomial equations
over $\Q$ in $B_1,B_2$, of degree at most five; then
$c_0=\lambda_1-1/12$ is determined. For a one-parameter family of
Part~\ref{part:arithmetic} the off-diagonal entry $b$ is a further unknown.

These systems are solved exactly. For a single matrix, a common factor of
all equations is split off and solved separately; otherwise the resultants
with respect to $B_2$ of one equation of positive $B_2$-degree with all
others, together with the equations free of $B_2$, give finitely many
candidates for $B_1$, and the gcd of the specialized equations gives
$B_2$. For the families we first split off common factors of all
equations, which define families of solutions, and remove factors of the
generators that have no rational zero or vanish only where $A$ is not
positive definite. We then compute the minimal associated primes over $\Q$
of the remaining ideal in $\Q[b,B_1,B_2]$ with \textsc{Singular}
\cite{Sing}; the independent calculation in Appendix~\ref{app:verification}
uses resultants, specialized gcds and exact Gr\"obner elimination instead; for the families with $x=y$ and $a=d$ we use the invariants
$B_1+B_2$ and $(B_1-B_2)^2$. In all cases every
rational solution is substituted back into all equations, and $\lambda_1$
is evaluated exactly.

By Proposition~\ref{deg:prop:coefficient-transport}, the conditions
\eqref{cls:eq:conditions} are compatible with exchange and with the
reflection
\begin{equation}\label{cls:eq:duality}
 (A,\boldsymbol B,c_0)\longmapsto
 \bigl(A^{-1},A^{-1}\boldsymbol B,\tfrac12\boldsymbol B^tA^{-1}\boldsymbol B-\tfrac1{12}-c_0\bigr),
\end{equation}
which is the familiar duality between $A$ and $A^{-1}$. We nevertheless
solve the system separately for every member of each orbit; all solutions
are compatible with \eqref{cls:eq:duality}. The conditions are necessary
only. In degrees one and two every solution turns out to be modular, so no
further coefficient is needed there.

%% file: sections/degree_two_classification.tex
\section{Nahm points of degree at most two}\label{cls:sec:deg2}
For the rank-one sums
$f_{a,\beta,c}(q)=\sum_{n\ge0}q^{an^2/2+\beta n+c}/(q)_n$ with
$a\in\{\frac12,1,2\}$, the modular cases are \cite{Zagier}
\begin{equation}\label{cls:eq:rankone}
 (a,\beta,c)\in\Bigl\{(2,0,-\tfrac1{60}),\ (2,1,\tfrac{11}{60}),\
 (1,0,-\tfrac1{48}),\ (1,\pm\tfrac12,\tfrac1{24}),\
 (\tfrac12,0,-\tfrac1{40}),\ (\tfrac12,\tfrac12,\tfrac1{40})\Bigr\}.
\end{equation}

\begin{table}[ht]
\centering
\caption{The modular data of Theorem~\ref{cls:thm:deg2}(c), up to exchange.
The matrices in the first four rows are invariant under exchange; the
exchange of a listed shift $(B_1,B_2)$ is $(B_2,B_1)$.}
\label{cls:tab:deg2}
\renewcommand{\arraystretch}{1.25}
\begin{tabular}{@{}llll@{}}
\toprule
Point $(x,y)$ & $A$ & $\boldsymbol B$ and $c_0$ & Modularity\\
\midrule
$(\rho,\rho)$ & $\left(\begin{smallmatrix}3/2&1/2\\1/2&3/2\end{smallmatrix}\right)$ &
$(\frac14,-\frac14),\ -\frac1{120};\quad(\frac14,\frac34),\ \frac{11}{120}$ & \cite[Table 1]{VZ}\\
$(\rho,\rho)$ & $\left(\begin{smallmatrix}4/3&2/3\\2/3&4/3\end{smallmatrix}\right)$ &
$(0,0),\ -\frac1{30};\quad(-\frac13,-\frac23),\ \frac1{30}$ & \cite{CRW}; \cite[Table 1]{VZ}\\
$(\rho^2,\rho^2)$ & $\left(\begin{smallmatrix}3/4&-1/4\\-1/4&3/4\end{smallmatrix}\right)$ &
$(\frac14,-\frac14),\ -\frac1{80};\quad(\frac12,0),\ \frac1{80}$ & \cite[Theorem 3.4]{VZ}\\
$(\rho^2,\rho^2)$ & $\left(\begin{smallmatrix}1&-1/2\\-1/2&1\end{smallmatrix}\right)$ &
$(0,0),\ -\frac1{20};\quad(-\frac12,0),\ \frac1{20}$ & \cite[Theorem 3.4]{VZ}\\
$(\frac1{\sqrt2},2-\sqrt2)$ & $\left(\begin{smallmatrix}2&1\\1&1\end{smallmatrix}\right)$ &
$\begin{array}{@{}l@{}}(0,0),\ -\frac1{32};\quad(0,\frac12),\ 0;\quad(1,\frac12),\ \frac18;\\
(-1,\frac12),\ \frac18;\quad(1,1),\ \frac7{32}\end{array}$ & \cite{Wang}\\
$(1-\frac1{\sqrt2},\sqrt2-1)$ & $\left(\begin{smallmatrix}1&-1\\-1&2\end{smallmatrix}\right)$ &
$\begin{array}{@{}l@{}}(0,0),\ -\frac5{96};\quad(\frac12,0),\ \frac1{24};\quad(-\frac12,1),\ \frac16;\\
(0,1),\ \frac{19}{96};\quad(-\frac32,2),\ \frac{25}{24}\end{array}$ & \cite{Wang}\\
\bottomrule
\end{tabular}
\end{table}

\begin{theorem}[Nahm points of degree at most two]\label{cls:thm:deg2}
Let $A$ be a rational symmetric positive-definite matrix, let
$\boldsymbol B\in\Q^2$ and $c_0\in\Q$, and suppose that the distinguished
Nahm point of $A$ has joint degree at most two. Then $f_{A,\boldsymbol B,c_0}$
is modular if and only if, up to exchange, one of the following holds.
\begin{enumerate}[label=(\alph*)]
\item $A=\left(\begin{smallmatrix}a&1-a\\1-a&a\end{smallmatrix}\right)$ with
$a>1/2$, and $(\boldsymbol B,c_0)$ is one of
\[
 \bigl((\beta,-\beta),\ \tfrac{\beta^2}{2a}-\tfrac1{24}\bigr)\ (\beta\in\Q),\qquad
 \bigl((-\tfrac12,-\tfrac12),\ \tfrac1{8a}-\tfrac1{24}\bigr),\qquad
 \bigl((\tfrac a2,1-\tfrac a2),\ \tfrac a8-\tfrac1{24}\bigr).
\]
\item $A=\operatorname{diag}(a_1,a_2)$, $\boldsymbol B=(\beta_1,\beta_2)^t$ and
$c_0=c_1+c_2$, where $(a_i,\beta_i,c_i)$ belong to \eqref{cls:eq:rankone}.
\item $(A,\boldsymbol B,c_0)$ is one of the $18$ triples in Table~\ref{cls:tab:deg2}.
\end{enumerate}
\end{theorem}

\begin{proof}
By Theorems~\ref{low:global} and~\ref{low:quadratic}, the matrix and its
point belong to one of the following classes: the rational point
$(1/2,1/2)$ with $A=\left(\begin{smallmatrix}1-b&b\\b&1-b\end{smallmatrix}\right)$,
$b<1/2$; the two independent quadratic orbits generated by \eqref{low:twopairs},
four matrices each; the golden families \eqref{low:goldfamily}; and the mixed
diagonal pairs $(1/2,\rho^r)$ with their exchanges. We solve
\eqref{cls:eq:conditions} on each class as in Section~\ref{cls:sec:method}.

\emph{The rational point.} Here $X=Y=1$ and $D=4(1-b)$. Every equation is
divisible by $B_1+B_2$; after this factor is removed, the remaining
equations have the common factor $(b-1)^4$, and $b=1$ is not admissible.
After removing these factors, the minimal associated primes
of the remaining ideal in $\Q[b,B_1,B_2]$ are the curve
$B_1=B_2=-\frac12$; the curve $B_1+B_2=1$, $b=2B_2-1$ and its exchange;
the four points $b=0$, $\boldsymbol B\in\{(\pm\frac12,0),(0,\pm\frac12)\}$;
two points with $b=1/2$, where $A$ is singular; and one zero-dimensional
prime without rational points, since its generator in $B_1+B_2$ alone is
$45s^2+68s+4$. With $a=1-b$, the three curves are the three families in (a);
symbolic substitution shows that $\Psi_2,\Psi_3,\Psi_4$ vanish identically
on them and that $\lambda_1-1/12$ equals the stated $c_0$. The four isolated
points give the products in (b) with $A=I$.

\emph{The independent quadratic orbits.} Each of the four matrices of the
$\sqrt2$ orbit has exactly five rational solutions. They are the last two
rows of Table~\ref{cls:tab:deg2} and their exchanges, and the other two
matrices are the reflections \eqref{cls:eq:duality} of these.
The four matrices of the $\sqrt5$ orbit, generated by
$\left(\begin{smallmatrix}4/3&-1\\-1&1\end{smallmatrix}\right)$, admit no
rational solution.

\emph{The golden families.} For $(r,s)=(1,1)$, so $x=y=\rho$ and
$A=\left(\begin{smallmatrix}2-b&b\\b&2-b\end{smallmatrix}\right)$, $b<1$, the
ideal is zero-dimensional after the factor $4b^2-10b+5$, which has no
rational zero, is removed from the first-coefficient equation. In the
coordinates $b$, $B_1+B_2$ and $(B_1-B_2)^2$ its eleven minimal primes are
rational points, and each gives rational shifts. Those with $b<1$ have
$b\in\{0,\frac12,\frac23\}$; the others have $b\in\{1,2\}$. For
$(r,s)=(2,2)$, so $x=y=\rho^2$ and
$A=\left(\begin{smallmatrix}1/2-b&b\\b&1/2-b\end{smallmatrix}\right)$,
$b<\frac14$, the same holds after removal of $16b^2-5$; the admissible points
have $b\in\{0,-\frac14,-\frac12\}$, and the others have $b\in\{\frac14,\frac12\}$.
For $(r,s)=(1,2)$ and $(2,1)$ the primes containing all equations are $b=1$,
outside the admissible range $b<\frac12$, and four rational points with $b=0$.
The values $b=\frac12,\frac23,-\frac14,-\frac12$ give the first four rows of
Table~\ref{cls:tab:deg2}, with their exchanges, and $b=0$ gives products
of rank-one data.

\emph{The mixed diagonal pairs.} Each of the four matrices has exactly the
six product shifts of (b).

In each case $c_0$ is the exact value of $\lambda_1-1/12$, and comparison
with \eqref{cls:eq:rankone} shows that every remaining solution with $b=0$, after accounting for
family (a), is a product in (b). This proves that a modular datum is of the form (a), (b) or (c).

Conversely, the data in (a) are modular by \cite[Table 2]{VZ}; the first
family is Zagier's \cite[Table 2]{Zagier}, from the identity
\[
 \sum_{n_1,n_2\ge0}
 \frac{q^{\frac a2(n_1-n_2)^2+n_1n_2+\beta(n_1-n_2)}}{(q)_{n_1}(q)_{n_2}}
 =\frac1{(q)_\infty}\sum_{n\in\Z}q^{\frac a2n^2+\beta n}.
\]
The data in (b) are products of modular rank-one sums. The rows of
Table~\ref{cls:tab:deg2} are modular by the cited results: the
$\sqrt2$ rows are Zagier's Examples~2 and~3, proved in \cite{Wang}; the
matrix with $a=4/3$ is Zagier's Example~10, conjectured explicitly in
\cite{VZ} and proved in \cite{CRW}; and the matrix with $a=1$ is Zagier's
Example~11.
\end{proof}

\begin{remark}
The second independent quadratic point $(\rho^3,\rho^2/2)$ has Rogers sum
$\pi^2/12$, yet none of its four matrices has a rational shift satisfying
\eqref{cls:eq:conditions}. The matrices
$\left(\begin{smallmatrix}3/2&1/2\\1/2&3/2\end{smallmatrix}\right)$ and
$\left(\begin{smallmatrix}3/4&-1/4\\-1/4&3/4\end{smallmatrix}\right)$ are
the counterexamples of \cite{VZ} to the original form of Nahm's conjecture.
Only the coefficients $\lambda_1,\ldots,\lambda_4$ were needed: the fifth
coefficient used in Part~\ref{part:degree} plays no role in degrees one
and two.
\end{remark}

%% file: sections/cubic_unit.tex
\section{Independent cubic unit data and their modular shifts}\label{cu:section}
Throughout this section $K=\Q(x,y)$ is the joint field of the distinguished
Nahm point, and
\[
 A=\begin{pmatrix}\alpha&\beta\\\beta&\gamma\end{pmatrix}
\]
is rational, symmetric and positive definite. We first classify the cubic
points for which all four arguments are ordinary units and one complementary
pair is multiplicatively independent. The remaining cubic branches are
excluded in Section~\ref{cr:section}.

\begin{theorem}[Independent cubic four-unit classification]\label{cu:thm:main}
Suppose $[K:\Q]=3$, $x,1-x,y,1-y\in\mathcal O_K^\times$, and $x,1-x$ are multiplicatively independent. Then, up to exchanging the two coordinates,
\begin{align*}
 A_6&=\begin{pmatrix}4&2\\2&2\end{pmatrix},&
 A_7&=\begin{pmatrix}1/2&-1/2\\-1/2&1\end{pmatrix},\\
 A_8&=\begin{pmatrix}3/2&1\\1&2\end{pmatrix},&
 A_9&=\begin{pmatrix}1&-1/2\\-1/2&3/4\end{pmatrix}.
\end{align*}
The odd saturation-index case is empty. In the even-index case, $K=\Q(\zeta_7+\zeta_7^{-1})$.
\end{theorem}

The proof first excludes odd saturation index by a bounded coefficient census, then treats even index by exact algebra in the real cubic field of conductor seven. Section~\ref{cu:sec:shifts} determines all rational modular shifts for these four matrices.

\subsection{The common unit lattice}
The field is totally real: a cubic field of signature $(1,1)$ has unit rank one, incompatible with the assumed independence. Neither $x$ nor $y$ can be rational. For $x$ this follows from independence, since rational ordinary units are $\pm1$. Also $\beta\ne0$, since otherwise the first Nahm equation makes $x,1-x$ dependent. Solving the equations for $y,1-y$ then shows that their span in the unit group tensored with $\Q$ is the same rank-two space, so they too are independent and generate $K$.

Let $E$ be the group of units positive at the distinguished embedding. It is free of rank two. Symmetry of $A$ gives, in $\bigwedge^2(E\otimes\Q)$,
\begin{equation}\label{cu:eq:wedge}
 x\wedge(1-x)=-y\wedge(1-y).
\end{equation}
As $\bigwedge^2E$ has no torsion, this is an integral equality. Thus
\[
 [E:\langle x,1-x\rangle]=[E:\langle y,1-y\rangle]=I.
\]
If necessary replace $(x,y,A)$ by $(1-x,1-y,A^{-1})$ so that $\beta>0$. This transformation preserves the equations, positive definiteness, the four-unit hypothesis and $I$. It is used only to classify points and matrices. Put
\[
 a=\alpha/\beta,\quad u=1/\beta,\quad
 c=(\alpha\gamma-\beta^2)/\beta,\quad j=\gamma/\beta.
\]
Then
\begin{equation}\label{cu:eq:pair}
 y=x^{-a}(1-x)^u,\qquad 1-y=x^{-c}(1-x)^j,
 \qquad a,u,c,j>0,\qquad aj-uc=1.
\end{equation}
All four exponents belong to $I^{-1}\Z$. Coordinate exchange preserves this form, interchanging $a,j$ and keeping $u,c$.

For a unit $w$, let $\ell(w)=(\log|\sigma_iw|)_{i=0}^2$. This is an injective map on $E$, whose image is a lattice in the zero-sum plane. Set
\[
 U=\ell(x),\quad V=\ell(1-x),\quad U'=\ell(y),\quad V'=\ell(1-y),
 \qquad R=\|U\wedge V\|.
\]
Equations~\eqref{cu:eq:pair} imply $U'=-aU+uV$, $V'=-cU+jV$ and therefore
\begin{equation}\label{cu:eq:area}
 U'\wedge V'=-U\wedge V,\qquad
 R=\|U'\wedge V'\|=I\operatorname{covol}(\ell(E)).
\end{equation}

\begin{lemma}\label{cu:lem:covol}
Every nonzero vector of $\ell(E)$ has length at least $\sqrt3\log\varphi$, where $\varphi=(1+\sqrt5)/2$. Consequently
\[
 \operatorname{covol}(\ell(E))\ge\frac{3\sqrt3}{2}(\log\varphi)^2>\frac35.
\]
\end{lemma}
\begin{proof}
For a nontrivial unit $w$, $w^2-1$ is a nonzero algebraic integer. Writing $\ell_i=\log|\sigma_iw|$ and using $\sum_i\ell_i=0$ gives
\[
 0\le\log|N(w^2-1)|=\sum_{i=0}^2\log(2\sinh|\ell_i|).
\]
The function $t\mapsto\log(2\sinh t)$ is increasing and concave. Jensen's inequality, followed by Cauchy--Schwarz, gives
\[
 0\le3\log\!\left(2\sinh\frac{\|\ell(w)\|_2}{\sqrt3}\right).
\]
Since $\operatorname{arsinh}(1/2)=\log\varphi$, the length bound follows. A planar lattice of shortest length $\lambda$ has area at least $\sqrt3\lambda^2/2$, by reducing a second basis vector along a shortest one. Finally $\sqrt3>433/250$ and $\log\varphi>481/1000$ give the stated strict rational bound.
\end{proof}

\subsection{Odd index: an explicit finite cutoff}
\subsubsection{Conjugates and area}
Suppose $I$ is odd. Clearing the denominators in~\eqref{cu:eq:pair} gives positive integers $P=Ia,Q=Iu,W=Ic,S=Ij$ with $PS-QW=I^2$ odd. If a conjugate $t$ of $x$ lies in $(0,1)$, oddness forces the corresponding conjugates of $y,1-y$ to be positive, and
\[
 t^{-a}(1-t)^u+t^{-c}(1-t)^j=1.
\]
The left side strictly decreases, so there is exactly one such conjugate: the distinguished $x$. If $1<t<2$, both $|\sigma(y)|$ and $|1-\sigma(y)|$ are below one and hence both arguments must be positive. Their signs are $(-1)^Q,(-1)^S$, contradicting $PS-QW$ odd. Rational endpoints cannot be roots of an irreducible cubic. Thus the two outer roots of each coordinate polynomial lie in
\begin{equation}\label{cu:eq:rootpattern}
 (-\infty,0)\cup(2,\infty).
\end{equation}

Put $s_x=-\log(1-x)$, $\delta_x=-\log x=\delta(s_x)$, where $\delta(s)=-\log(1-e^{-s})$, and $b_0=\log2$. For the two outer roots $r_i$ put $l_i=\log|r_i|$, $m_i=\log|1-r_i|$. Unit norms give
\[
 m_1+m_2=s_x,\quad l_1+l_2=\delta_x,\quad m_i>0.
\]
With $P_x=\max(l_1,0)+\max(l_2,0)=\log\operatorname{M}(x)$, the elementary inequalities
\[
 m_i\le\max(l_i,0)+b_0,\qquad \max(l_i,0)\le m_i+b_0
\]
give $s_x-2b_0\le P_x\le s_x+2b_0$ and
\begin{equation}\label{cu:eq:upper}
 \|U\|\le\sqrt2P_x,\quad \|V\|\le\sqrt2s_x,\quad
 R\le2s_x(s_x+2b_0).
\end{equation}
The vector norm bounds follow by bounding each of its positive and negative parts by its total mass.

If $s_x\ge4$, the $l_i$ have opposite signs: otherwise $s_x\le\delta_x+2b_0$. The minor of the distinguished row $(-\delta_x,-s_x)$ and either suitably chosen outer row has absolute value at least $s_x(P_x-\delta_x)$. For a positive outer $l_i=P_x$, use $m_i\le s_x$; for a negative outer $l_i=-(P_x-\delta_x)$ the two terms have the same sign. In the zero-sum plane in $\mathbb R^3$, the three minors have equal absolute value, so $R$ is $\sqrt3$ times any such minor. In particular the following weaker, convenient bound holds:
\begin{equation}\label{cu:eq:lower}
 R\ge\sqrt2\,s_x(s_x-3b_0-2\delta_x)\qquad(s_x\ge4).
\end{equation}
All these statements also hold for $y$, with the same $R$.

\subsubsection{The smaller coordinate and its partner}
\begin{proposition}\label{cu:prop:bounds}
An odd-index configuration would satisfy, after exchanging the coordinates,
\[
 s_x<4,\quad s_y<7,\quad H(f_x)<330,\quad H(f_y)<6600,
 \qquad I\le71,\quad \max(a,u,c,j)<151.
\]
Here $H(f)=\max(3,|\text{coefficients of }f|)$ for a monic minimal polynomial.
\end{proposition}
\begin{proof}
Assume instead $4\le s=s_x\le s_y$. We may use
\[
 b_0<347/500,\qquad \delta_x,\delta_y<19/1000,\qquad\sqrt2>707/500.
\]
Comparing~\eqref{cu:eq:upper} and~\eqref{cu:eq:lower} gives
\[
 (s_y/s_x)^2<\frac{2(1+2(347/500)/4)}{(707/500)(1-(3(347/500)+2(19/1000))/4)}<(21/10)^2.
\]
At the distinguished embedding, $us_x=a\delta_x+\delta_y$. Since $aR=\|U'\wedge V\|$, the vector bounds imply
\begin{align*}
 uR&=\frac{aR\delta_x+R\delta_y}{s_x}\\
 &\le2\delta_x(s_y+2b_0)+2\delta_y(s_x+2b_0)\\
 &<2\delta(s)\bigl((31/10)s+4(347/500)\bigr).
\end{align*}
This last expression decreases for $s\ge4$, since $\delta'(s)<-\delta(s)$. At $s=4$ it is less than
\[
 2(19/1000)\bigl((31/10)4+4(347/500)\bigr)=0.576688<3/5.
\]
But $uR=(Iu)\operatorname{covol}(\ell(E))>3/5$, since $Iu$ is a positive integer. This contradiction proves $\min(s_x,s_y)<4$.

Choose $s_x<4$. Then $R<2\cdot4(4+2(347/500))<44$. If $s_y\ge7$, the increasing lower bound~\eqref{cu:eq:lower} gives
\[
 R>(707/500)7\bigl(7-3(347/500)-2(19/1000)\bigr)>44,
\]
another contradiction. Hence $s_y<7$.

For a cubic, $H(f)\le3\operatorname{M}(f)$, except that our harmless floor $H\ge3$ can be taken separately. If $x\le1/4$, then $\operatorname{M}(x)\le4/(1-x)\le16/3$, so $H(f_x)\le16$. If $x>1/4$, both outer roots cannot exceed two, by $|N(x)|=1$. At a negative root the ratio $\max(1,|r|)/|1-r|$ is at most one; at a root above two it is less than two. Thus $\operatorname{M}(x)<2/(1-x)=2e^{s_x}$. Using $e^4<55$ and $e^7<1100$ proves the two height bounds.

Finally $I=R/\operatorname{covol}(\ell(E))<72$. Exterior products in the linear equations for $U',V'$ bound every exponent by
\[
 \frac{2(4+2(347/500))(7+2(347/500))}{3/5}<151.
\]
The numerical comparisons above can all be made with rational Taylor bounds; these bounds are included in the finite certificate.
\end{proof}

\subsubsection{Complete exact exclusion}
A cubic exceptional unit has minimal polynomial
\begin{equation}\label{cu:eq:cubics}
 f(Z)=Z^3+AZ^2+BZ+e,\quad e,p\in\{\pm1\},\quad
 f(1)=p,\quad B=p-1-e-A.
\end{equation}
The box $H<6600$ is therefore a finite enumeration over $-6599\le A\le6599$, four sign choices, and $|B|<6600$. We require positive discriminant, no root $-1$ (the only remaining possible rational root), one root in $(0,1)$ and no root in $(1,2)$. These are exact tests. For example, a positively rescaled cubic Sturm sequence is
\[
 f,\quad f',\quad PZ+Q_0,\quad H_0,
\]
where $P=2(A^2-3B)$, $Q_0=AB-9e$, and $H_0=-3Q_0^2+2AQ_0P-BP^2$. Here $H_0=9\operatorname{disc}(f)$, so positive discriminant gives $P,H_0>0$.

For each retained polynomial we enclose all three roots in disjoint dyadic intervals of width $2^{-64}$. Their existence and completeness are checked by exact polynomial sign changes. The distinguished root is placed first. All logarithms are bounded outward on the integer scale $2^{96}$, using range reduction and
\[
 \log z=2\sum_{j=0}^{39}\frac{t^{2j+1}}{2j+1}+\text{positive tail},\qquad t=(z-1)/(z+1),\quad0\le t\le1/3,
\]
with an explicit geometric bound for the tail. Set
\[
 D_f=\log|r_0|\log|1-r_1|-\log|r_1|\log|1-r_0|.
\]
Then $R_f^2=3D_f^2$, independent of which outer root is second. Equal-area pairs must have intersecting outward $2^{-20}$-grid intervals for $D_f^2$. The complete certificate gives:
\begin{center}\begin{tabular}{lr}\toprule
Quantity&Count\\\midrule
Polynomials in the full coefficient/sign box&52,790\\
Irreducible totally real polynomials before root-pattern test&52,772\\
Partner polynomials with the required pattern, $H<6600$&13,196\\
Small-coordinate polynomials, $H<330$&656\\
Pairs passing the equal-area interval test&656\\
Pairs with different polynomials passing the test&0\\\bottomrule
\end{tabular}\end{center}
Each retained pair has identical minimal polynomial. Its unique root in $(0,1)$ forces $x=y$, whereas~\eqref{cu:eq:wedge} would give $2x\wedge(1-x)=0$, contrary to independence. Thus the odd-index branch is empty.

The full verifier \path{verify_odd_cubic.py} uses only the Python standard library. It reconstructs the coefficient census with exact Sturm signs, checks every root enclosure and every logarithm interval, and verifies the complete absence of distinct-polynomial equal-area pairs. FLINT is used only by the optional proposal generator \path{odd_cubic_screen.py}. The certificate has no floating-point accept/reject threshold.

\subsection{Even index: three consecutive units}
If $I$ is even, the vectors of $x,1-x$ are dependent in $E/2E$. At least one of $x,1-x,x(1-x)$ is a square in $K$. The third possibility would force every real conjugate of $x$ into $(0,1)$, contradicting $|N(x)|=1$; both of the first two being squares is excluded in the same way. Thus exactly one coordinate is a square. The same reasoning applies to $y$.

The positive root $z\in(0,1)$ of that square generates $K$, and $z,1-z,1+z$ are units. Indeed their product relations are units and the factors are algebraic integers. Its monic cubic $g$ therefore satisfies $g(0)=e$, $g(1)=p$, $g(-1)=m$, with $e,p,m\in\{\pm1\}$. Hence
\begin{equation}\label{cu:eq:threeunits}
 g(Z)=Z^3+\left(\frac{p+m}{2}-e\right)Z^2+
 \left(\frac{p-m}{2}-1\right)Z+e.
\end{equation}
Four of the eight polynomials have discriminant $49$, and four have discriminant $-23$. Only the first four can occur here. They have a total of three roots in $(0,1)$.

Take $t\in(0,1/2)$ with $t^3-t^2-2t+1=0$, and put $s=1-t$. The three interval roots are exactly $t,s,t/s$. The six possible original coordinates and their exact multiplicative expressions are
\begin{center}\begin{tabular}{ccl}\toprule
Index&Coordinate&Expression in $t,s$\\\midrule
0&$t^2$&$t^2$\\
1&$1-t^2$&$ts^{-1}$\\
2&$(1-t)^2$&$s^2$\\
3&$2t-t^2$&$t^{-1}s^2$\\
4&$t+t^2$&$t^2s^{-2}$\\
5&$1-t-t^2$&$t^2s^{-1}$\\\bottomrule
\end{tabular}\end{center}
Complementation interchanges indices $0,1$, $2,3$, and $4,5$. Each identity is verified in $\Q[t]/(t^3-t^2-2t+1)$. The units $t,s$ are independent: their distinguished values are both below one, while an embedding with $1<t<2$ has $|s|<1<|t|$, excluding a positive-power relation.

For each of the 36 ordered pairs of coordinates, form the two-column exponent matrix $V$ and the complementary exponent matrix $W$. If $V$ is invertible, the only possible Nahm matrix is $(V^{-1}W)^T$. Exact rational algebra gives eight symmetric positive definite matrices, paired by coordinate exchange:
\begin{center}\begin{tabular}{ccl}\toprule
Indices $(x,y)$&Matrix&Exchange\\\midrule
$(1,3)$&$A_6$&$(3,1)$\\
$(0,2)$&$A_7$&$(2,0)$\\
$(4,3)$&$A_8$&$(3,4)$\\
$(5,2)$&$A_9$&$(2,5)$\\\bottomrule
\end{tabular}\end{center}
Pairs with singular $V$ cannot occur under~\eqref{cu:eq:pair}, since that equation expresses $y$ with a nonzero coefficient of $1-x$. The list is closed under inversion, since $A_6^{-1}=A_7$ and $A_8^{-1}=A_9$. Thus it also contains the original matrix if simultaneous complementation was used to make the coupling positive. This proves Theorem~\ref{cu:thm:main}. The finite algebra is certified by \path{even_cubic_and_shifts.py}.

\subsection{All rational modular shifts}\label{cu:sec:shifts}
\begin{theorem}[The cubic unit shifts]\label{cu:thm:shifts}
For each matrix in Theorem~\ref{cu:thm:main}, the Nahm sum with rational
$\mathbf B=(B_1,B_2)$ and $c_0$ is modular if and only if the shifts occur
in Table~\ref{cu:tab:shifts}. Coordinate exchange exchanges the two entries
of $\mathbf B$. There are twelve triples up to this exchange.
\end{theorem}
\begin{proof}
With the radial normalization of Section~\ref{cls:sec:method}, modularity requires
\begin{equation}\label{cu:eq:necessary}
 \lambda_1\in\Q,\qquad \lambda_2=0,\qquad c_0=\lambda_1-1/12.
\end{equation}
Our normalization is
\[
 \log f_{A,\mathbf B,c_0}(e^{-\varepsilon})
 \sim \frac{\alpha_0}{\varepsilon}+\log\beta_0+
 (\lambda_1-c_0-1/12)\varepsilon+\lambda_2\varepsilon^2+\cdots.
\]
These are necessary conditions from the radial expansion; their use does not assert that finitely many vanishing coefficients suffice for modularity.

For the coefficient calculations used here and in Section~\ref{cr:section}, we record the explicit Gaussian formulas. Put $X=x/(1-x)$, $Y=y/(1-y)$ and let $Z=(Z_1,Z_2)$ be a centered formal Gaussian of covariance $(A+\operatorname{diag}(X,Y))^{-1}$. For $q=X$ or $Y$ define
\begin{align*}
 P_1(q)&=q(q+1),&P_2(q)&=q(q+1)(2q+1),\\
 P_3(q)&=q(q+1)(6q^2+6q+1),&
 P_4(q)&=q(q+1)(2q+1)(12q^2+12q+1).
\end{align*}
With $(q_1,q_2)=(X,Y)$, set
\begin{align*}
 U_1&=\sum_i\left[-(B_i+q_i/2)Z_i+P_1(q_i)Z_i^3/6\right],\\
 U_2&=-(X+Y)/12+\sum_i\left[P_1(q_i)Z_i^2/4-P_2(q_i)Z_i^4/24\right],\\
 U_3&=\sum_i\left[P_1(q_i)Z_i/12-P_2(q_i)Z_i^3/12+P_3(q_i)Z_i^5/120\right],\\
 U_4&=\sum_i\left[-P_2(q_i)Z_i^2/24+P_3(q_i)Z_i^4/48-P_4(q_i)Z_i^6/720\right].
\end{align*}
Then
\begin{align*}
 \lambda_1&=\mathbb E(U_2+U_1^2/2),\\
 \lambda_2&=\mathbb E(U_4+U_1U_3+U_2^2/2+U_1^2U_2/2+U_1^4/24)-\lambda_1^2/2.
\end{align*}
All expectations are evaluated by Wick contractions in the exact cubic quotient field. The certificate includes both these Gaussian expressions and the resulting rational coefficient polynomials.

Write $B=B_1,C=B_2$ in the following table only. The two nonconstant field coefficients of $\lambda_1$, and the three field coefficients of $\lambda_2$, give five equations over $\Q$. Their ideals have the following small generating sets:
\begin{center}\begin{tabular}{cl}\toprule
Matrix&Polynomial equations in $B,C$\\\midrule
$A_6$&$B^2-B-2C,\quad BC-2C,\quad C^2-C$\\
$A_7$&$2B^2-B,\quad 2BC-C,\quad 2C^2+C$\\
$A_8$&$2B^2+B-C,\quad 2BC-C,\quad C^2-C$\\
$A_9$&$B-8C^2+4C,\quad 8C^3-6C^2+C$\\\bottomrule
\end{tabular}\end{center}
Explicit polynomial multipliers express each displayed generator as a combination of the five necessary equations. Every solution listed below is also substituted into the original five equations. Thus the resulting zero set, including arbitrary rational shifts, is exactly the displayed finite list.

\begin{table}[htbp]
\centering
\begin{tabular}{ccr}\toprule
Matrix&$(B_1,B_2)$&$c_0$\\\midrule
$A_6$&$(0,0)$&$-1/42$\\
&$(1,0)$&$5/42$\\
&$(2,1)$&$17/42$\\\midrule
$A_7$&$(0,0)$&$-5/84$\\
&$(1/2,0)$&$1/84$\\
&$(1/2,-1/2)$&$1/21$\\\midrule
$A_8$&$(0,0)$&$-5/168$\\
&$(-1/2,0)$&$1/168$\\
&$(1/2,1)$&$25/168$\\\midrule
$A_9$&$(0,0)$&$-3/56$\\
&$(0,1/2)$&$1/56$\\
&$(-1/2,1/4)$&$1/28$\\\bottomrule
\end{tabular}
\caption{The twelve modular triples of joint degree three, up to coordinate exchange.}
\label{cu:tab:shifts}
\end{table}
The product identities in Wang's Theorems 3.6--3.9~\cite{Wang} prove sufficiency for all twelve candidates. The normalization of one of these identities requires the following explicit calculation.

\paragraph{The normalization of Example 8.}
Wang's arXiv v3 Table 9 lists $5/21$ for the third shift of Example 8. In the normalization used throughout this paper the value is $25/168$. It can be checked directly from the product identity (3.76), avoiding either table as an authority. In Wang's notation,
\[
 J_m=(q^m;q^m)_\infty,\qquad
 J_{r,m}=(q^r,q^{m-r},q^m;q^m)_\infty.
\]
Their explicit powers of $q$ before the eta/theta factors are $-m/24$ and $-r^2/(2m)+r/2-m/8$, respectively. Therefore the product
\[
 \frac{J_{14}J_{28}^2J_{10,28}}{J_{5,28}J_{8,28}J_{9,28}J_{12,28}}
\]
has explicit exponent $-25/42$. The Nahm-series variable in (3.76) is $p=q^4$, so the normalizing factor is $p^{25/168}$. This proves the stated normalization and completes the sufficiency argument.
\end{proof}

\subsection{Equal canonical exponents}\label{cu:sec:equal}
In the direct canonical family supplied by Lemma~\ref{qcov:entry}, equal exponents $a=b$ give
\[
 1-z^a=z^T(1+z^a),\qquad 0<z<1,\qquad \gcd(a,T)=1,
\]
with coordinates and matrix
\[
 x=\frac{z^a}{1+z^a},\qquad y=z^a,\qquad
 A=\begin{pmatrix}1&-1\\-1&1+T/a\end{pmatrix}.
\]
For the canonical plus family, equal exponents instead force $z^T=1$ and are impossible. The canonical coverage lemma preserves the joint field, so the remaining boundary to consider here is exactly $[\Q(z):\Q]=3$.

\begin{proposition}\label{cu:prop:equal}
Suppose $z$ in the displayed direct family has degree three. No rational shift makes the associated Nahm sum modular. The same is true after simultaneous complementation, which replaces $A$ by $A^{-1}$, and after coordinate exchange.
\end{proposition}
\begin{proof}
The minimal polynomial $f$ of $z$ divides
\[
 Q_{a,T}(Z)=Z^{a+T}+Z^a+Z^T-1.
\]
For the factor search only, exchange $a,T$ so $a\le T$. If $a=T$, primitivity gives $a=T=1$ and the polynomial has degree two. Otherwise $\operatorname{M}(Q_{a,T})\le\|Q_{a,T}\|_2=2$ by Jensen's inequality. Hence $\operatorname{M}(f)\le2$, so
\[
 f=Z^3+AZ^2+BZ+e,\quad |A|,|B|\le6,\quad e=\pm1,
 \qquad f(1),f(-1)\in\{\pm1,\pm2\}.
\]
On $[0,1]$, $|f'|\le21$. Since $|f(1)|\ge1$, the mean value theorem gives $z\le20/21$. Writing $q=z^a$, $r=z^T\le q$, we have $qr+q+r=1$, whence $q\ge\sqrt2-1>2/5$. Therefore
\[
 a(-\log z)<\log(5/2)<1,\quad a<21.
\]
Also $z^T=(1-z^a)/(1+z^a)\ge(1-z)/(1+z)\ge1/41$, so
\[
 T(-\log z)\le\log41<4,\quad T<84.
\]
The endpoints and rational-root test leave 16 irreducible cubic polynomials. Exact division over the 918 coprime exponent pairs in this box performs 14,688 checks. The only cubic factors with a root in $(0,1)$ are
\begin{center}\begin{tabular}{ccc}\toprule
$(a,T)$ with $a\le T$&Minimal polynomial&Discriminant\\\midrule
$(1,2)$&$Z^3+Z^2+Z-1$&$-44$\\
$(1,7)$&$Z^3+Z^2-1$&$-23$\\\bottomrule
\end{tabular}\end{center}
The additional factor $Z^3-Z+1$ at $(1,7)$ has no root in $(0,1)$. The factor search is exhaustive, not a test up to an assumed cutoff.

For the Nahm matrices both orders $(a,T)$ must be restored, giving four matrices. Their four inverses are tested separately; modular duality is not assumed. At each of these eight exact algebraic points, the five necessary coefficient equations from~\eqref{cu:eq:necessary} generate the unit ideal. The file \path{equal_exponent_results.json} gives polynomial multipliers $h_i(B_1,B_2)$ with
\[
 \sum_{i=1}^{5}h_i E_i=1.
\]
Here the $E_i$ are the two nonconstant coefficients of $\lambda_1$ and the three coefficients of $\lambda_2$, after clearing rational denominators. Each identity is checked by exact polynomial expansion and the two coefficient implementations agree. Thus there is no solution even over an algebraic closure, and in particular no rational shift. Coordinate exchange preserves the Nahm sum, completing the proof.
\end{proof}

Cubic points therefore occur on the equal-exponent boundary, but none of the associated matrices or their inverses admits a modular shift.

%% file: sections/cubic_remaining.tex
\section{Excluding the remaining cubic branches}\label{cr:section}
We now exclude nonunit, multiplicatively dependent and diagonal modular
data of joint degree three. Together with the preceding section this gives
the complete cubic classification, stated in Theorem~\ref{cr:classification}.

In the notation of the canonical coverage lemma, Lemma~\ref{qcov:entry}, let $z\in(0,1)$ be an ordinary algebraic unit of degree three. Put $g=\gcd(a,b)$, $M=a/g$, $N=b/g$, with $a,b>0$ and $T\ne0$ integral, and require $\gcd(g,T)=1$. The equations are
\begin{align}
 (1-z^b)^N&=z^T(1+z^a)^M,&T&>0,\tag{D}\label{cr:eq:D}\\
 (1+z^b)^N&=z^T(1+z^a)^M,&gT&>b^2\ \text{or}\ gT<-a^2.\tag{P}\label{cr:eq:P}
\end{align}
The two binomials have equal finite-prime support. The nonunit hypothesis is that this support is nonempty, equivalently $1+z^a$ is not a unit.

\begin{theorem}[Unequal-exponent cubic nonunit exclusion]\label{cr:thm:nonunit}
There is no datum satisfying these hypotheses and either canonical equation with $a\ne b$.
\end{theorem}

\begin{theorem}[Dependent and diagonal cubic exclusion]\label{cr:thm:dependent}
Let $(x,y)$ be a distinguished Nahm point of joint degree three. If $x,1-x$ are multiplicatively dependent, or $y,1-y$ are dependent, or $A$ is diagonal, then no rational shifts $\mathbf B,c_0$ make $f_{A,\mathbf B,c_0}$ modular. No ordinary-unit hypothesis is imposed in this theorem.
\end{theorem}

The first theorem is an arithmetic exclusion of the displayed equations. The second uses the modularity-to-torsion theorem and the radial necessary conditions from Section~\ref{cls:sec:method}.

\subsection{Explicit bounds for the canonical cubic equations}
\subsubsection{Height gap and reciprocal separation}
An irreducible cubic unit polynomial cannot be reciprocal: reciprocal or antireciprocal odd-degree polynomials have a root at $1$ or $-1$. Nor can a cubic irreducible polynomial be cyclotomic. If $f$ is the minimal polynomial of $z$, then
\begin{equation}\label{cr:eq:gap}
 H:=\operatorname{M}(f)>13/10,\qquad h(z)=\tfrac13\log H>1/12.
\end{equation}
Here is a finite proof of the numerical gap. If $\operatorname{M}(f)\le13/10$, then
\[
 f=X^3+AX^2+BX+e,\qquad |A|,|B|\le3,\quad e=\pm1.
\]
There are 98 coefficient choices. Removing roots $\pm1$ leaves 76 irreducible cubics. For all 76, rigorous 192-bit complex root balls give a lower bound for the product of the root moduli above one that exceeds $13/10$. The program verifies every strict inequality; no floating-point equality decision is used. The rational Taylor bound $e^{1/4}<13/10$ proves the height inequality.

Let $\alpha_1,\alpha_2,\alpha_3$ be the conjugates and put $r_i=\min(|\alpha_i|,|\alpha_i|^{-1})$. No conjugate has modulus one, since otherwise the irreducible $f$ would share a root with its reciprocal. With $f^*(X)=X^3f(1/X)$,
\[
 1\le|\operatorname{Res}(f,f^*)|=\prod_{i,j}|1-\alpha_i\alpha_j|.
\]
The three factors pairing each root with its complex conjugate have product $H^2\prod_i(1-r_i^2)$. Bound the other six factors by $2\max(1,|\alpha_i|)\max(1,|\alpha_j|)$. Since $1+r_i<2$,
\begin{equation}\label{cr:eq:separation}
 \prod_i(1-r_i)\ge2^{-9}H^{-6}.
\end{equation}
This estimate supplies the norm and primitive-divisor bounds below.

\subsubsection{A primitive-index cutoff}
We prove that 4830 is a sufficient primitive-index cutoff. The slightly weaker constants below keep the local and archimedean estimates simple.

For a prime ideal above $p$ let $e,f_0$ be its ramification and residue degrees, and let $m$ be the residue order of $z$. A nonprimitive divisor of $\Phi_n(z)$ has $n=mp^k$ with $k\ge1$. Set $s_j=v(z^{mp^j}-1)$. Its valuation in $\Phi_n(z)$ is $s_k-s_{k-1}$. The binomial expansion shows that an increment greater than $e$ can occur only at the threshold $s_{k-1}=e/(p-1)$, requiring
\[
 (p-1)p^{k-1}\le e,\qquad ef_0\le4,\qquad n\le(p^{f_0}-1)p^k.
\]
The finite possibilities $p=2,3,5$ give $n\le30$. Consequently, if $n>30$ and no primitive prime exists,
\begin{equation}\label{cr:eq:normupper}
 |N\Phi_n(z)|\le n^4.
\end{equation}
Indeed every nonprimitive valuation is at most $e$, its residue characteristic divides $n$, and summing $ef_0$ above each $p$ gives at most four.

To obtain a lower bound, replace each outside conjugate by its inverse, calling the resulting numbers $\beta_i$. Set $Q_m=\prod_i|1-\beta_i^m|$. The estimate~\eqref{cr:eq:separation}, weakened for convenience, gives
\[
 -16\log2-8\log H\le\log Q_m\le4\log2.
\]
In the M\"obius formula for $\log|N\Phi_n(z)|$, exactly $2^{\omega(n)-1}$ nonzero coefficients have each sign. Thus, using $\log2/\log H<27/10$,
\[
 \log|N\Phi_n(z)|>\bigl(\varphi(n)-31\,2^{\omega(n)}\bigr)\log H.
\]
Together with~\eqref{cr:eq:gap} and~\eqref{cr:eq:normupper}, this proves a primitive prime exists if
\begin{equation}\label{cr:eq:totient}
 \varphi(n)>31\,2^{\omega(n)}+16\log n.
\end{equation}
For $n\ge10^6$, use $\varphi(n)/2^{\omega(n)}\ge\sqrt{n/30}$: the only factors below one in the prime-power product are bounded by the contributions at $2,3,5$, whose product is $1/\sqrt{30}$. Then it suffices that $\sqrt{n/30}>31+8\log n$. This holds at $10^6$ and its difference is increasing for $n>7680$.

For $31\le n<10^6$, an exact integer sieve tests the stronger sufficient inequality
\[
 5\varphi(n)>155\,2^{\omega(n)}+56\operatorname{bitlength}(n),
\]
using $\log2<7/10$. It has 592 failures, with largest 4830. The set
\[
 E=\{1,\ldots,30\}\cup\{\text{these 592 failures}\}
\]
is explicitly reconstructed and checked against the supplied file. Every $n\notin E$ has a primitive prime. Membership in $E$ means only that this sufficient test failed.

\subsubsection{Finite exponent bounds}
Equal binomial supports force
\begin{align*}
 D:&\quad b\notin E\Longrightarrow b\mid2a,\qquad
                  2a\notin E\Longrightarrow2a\mid b;\\
 P:&\quad 2a\notin E\Longrightarrow b/a\text{ is an odd integer},\qquad
                  2b\notin E\Longrightarrow a/b\text{ is an odd integer}.
\end{align*}
These follow by evaluating the binomials at a primitive prime. If $b>4830$ in D, both implications apply and force $b=2a$. Then $M=1,N=2$, and the valuation equation $2v(1-z^{2a})=v(1+z^a)>0$ is impossible because $1+z^a$ divides $1-z^{2a}$. The same reasoning bounds $2a$. In P, either index exceeding the cutoff forces both odd-ratio implications and hence $a=b$, which would give $z^T=1$. Therefore
\begin{equation}\label{cr:eq:box}
 D:\ a\le2415,\ b\le4830;\qquad
 P:\ a,b\le2415.
\end{equation}

Apply the standard projective height inequality to $(1\pm Z^b)^N/(1+Z^a)^M$. Its numerator and denominator have degree at most $\max(a^2,b^2)/g$ and coefficient absolute sum at most $2^{\max(M,N)}$. From~\eqref{cr:eq:gap} and $12\log2<9$,
\begin{equation}\label{cr:eq:T}
 |T|<\max(a^2,b^2)/g+9\max(M,N).
\end{equation}

There is a useful cubic norm bound. For either binomial and any $n\ge1$,
\[
 \log|N(1\pm z^n)|=n\log H+\epsilon_\pm(n),\qquad
 |\epsilon_\pm(n)|\le9\log2+6\log H
\]
by~\eqref{cr:eq:separation}. Taking norms in D or P eliminates $T$ and gives
\[
 |a-b|\log H\le9\log2+6\log H.
\]
Since $(13/10)^{24}>2^9$, we obtain
\begin{equation}\label{cr:eq:difference}
 |a-b|\le29.
\end{equation}

Finally some conjugate of $z$ has modulus $\rho\ge H^{1/2}>(13/10)^{1/2}$: the distinguished root lies inside the unit disk, leaving at most two outside conjugates. Put $\ell=\frac12\log(13/10)$ and $u_n=e^{-n\ell}$. Evaluation at a conjugate of maximum modulus yields
\begin{equation}\label{cr:eq:window}
 \frac{N\log(1-u_b)-M\log(1+u_a)}{\ell}
 \le T-\frac{b^2-a^2}{g}\le
 \frac{N\log(1+u_b)-M\log(1-u_a)}{\ell}.
\end{equation}
Both error bounds move toward zero as $\rho$ increases, which justifies substituting its fixed lower bound. The census rounds these endpoints outward with rigorous ball arithmetic.

\subsection{The finite-prime support restrictions}
There must be at least one common prime. Write $w=z^g$, $h_0=\min(M,N)$, $d_0=|M-N|>0$.

At an odd common prime, the residue of $w$ is $-1$. This requires $M$ odd and $N$ even in D, and both odd in P. With $e=v(p)\le3$, a necessary condition is
\begin{equation}\label{cr:eq:oddlocal}
 p\mid h_0,\qquad d_0\mid e\,v_p(h_0).
\end{equation}
For example, in D, if $N>M$, put $s=v(1+w)$. Unless $p\mid M$, the equation $N v(1-w^N)=M v(1+w^M)$ is impossible. Thus $p\mid M$, $p\nmid N$, and $M\mid s$, so $s\ge M\ge p>e/(p-1)$. The valid lifting formula then gives $(N-M)(s/M)=e v_p(M)$. Exchange the indices for the other case, and use the corresponding plus sign for P. The census checks existence of such $p,e$ exactly.

At a common dyadic prime, $\bar w=1$ and $e=v(2)\le3$. The exhaustive necessary alternatives are:
\begin{center}\begin{tabular}{lll}\toprule
Family&Parities&Necessary condition\\\midrule
D&$M,N$ odd&$\min(M,N)\in\{1,3\}$\\
D&$M$ even, $N$ odd&$N\in\{1,3\}$\\
D&$M$ odd, $N$ even&$M>N$ and $(N=2\text{ or }M-N\mid e v_2(N))$\\
P&Both odd&Impossible\\
P&Mixed&Odd index belongs to $\{1,3\}$\\\bottomrule
\end{tabular}\end{center}
In the third row $M>N$ is required in both alternatives. To prove the table, set $s=v(1+w)$, $t=v(1-w)$. Odd powers preserve the respective valuations. For D with both indices odd, $Nt=Ms$ gives $t=Mk,s=Nk$ and hence $k\min(M,N)=e\le3$. For P with both odd, $Ns=Ms$ is impossible.

In the mixed cases other than D with even $N$, write the even index as $r$ and the odd index as $o$. The valuation equality gives $v(1-\varepsilon w)=rk\ge2>e/2$ and $v(1+w^r)=ok$. The first inequality makes $v(w^2-1)>e$, so the second valuation is $e$, and $o\in\{1,3\}$. The strict threshold uses $e\le3$, as required for a cubic field.

For D with even $N$, divisibility of $1-w^N$ by $1+w$ forces $M>N$. Write $s=Nk$, $v(1-w^N)=Mk$. If $s\le e$, then $N=2$. Otherwise $t=e$ and lifting after squaring gives $v(1-w^N)=s+e v_2(N)$, proving the last divisibility condition.

\subsection{The complete canonical census and its certificate}
The census uses precisely~\eqref{cr:eq:box}--\eqref{cr:eq:window}, primitivity, the primitive-index implications, and the preceding support conditions. It removes $a=b$, already covered by Proposition~\ref{cu:prop:equal}, and the unit degeneracy $b=2a$ in D. Every ordered pair and every surviving signed $T$ is retained explicitly.
\begin{center}\begin{tabular}{lr}\toprule
Stage&Count\\\midrule
Ordered unequal pairs after $|a-b|\le29$&278,835\\
After primitive-index implications&36,649\\
After finite-prime support restrictions&5,019\\
Pairs with a surviving archimedean/PD/primitive interval&553\\
Individual parameter triples in those pairs&5,844\\\midrule
Triples rejected at $p=2$&2,731\\
Additional triples rejected at $p=3$&2,209\\
Additional triples rejected at $p=5$&805\\
Additional triples rejected at $p=7$&63\\
Survivors after these four primes&36\\
Survivors after cyclotomic refinement through $19$&16\\\bottomrule
\end{tabular}\end{center}

Here is the exact factor test. Clear negative powers in D or P to obtain $F\in\Z[Z]$. When D has $N$ even, $M$ odd, or P has both odd, divide the known factor $(1+Z^g)^{\min(M,N)}$, obtaining $H\in\Z[Z]$. This factor cannot vanish at $0<z<1$. A cubic unit minimal polynomial is monic with constant $\varepsilon=\pm1$, so modulo each prime it must appear among \emph{all} monic cubics with constant $\varepsilon$. Reducible and repeated reductions are included, and the same integer sign is maintained across primes.

The quotient is evaluated without inverting $1+Z^g$. For example, put $S=1+Z^g$ and $U_j=\sum_{i=0}^{j-1}(-Z^g)^i$. In a cancellation case the two terms of $H$ are
\[
 U_N^N S^{N-h_0},\qquad Z^T U_M^M S^{M-h_0},
\]
with negative $T$ cleared by multiplying the first term instead. Thus reductions where $S$ is noninvertible are covered.

For refinement, consider the orders $1$, $2$, $3$, $4$, $5$, $6$, $8$, $10$ and $12$. Remove one copy of each corresponding cyclotomic polynomial that divides $H$; call their product $C$. A characteristic-zero noncyclotomic cubic factor $f$ is coprime to $C$, so $fC\mid H$. Testing $\bar f\bar C\mid\bar H$ preserves possible collisions and multiplicities modulo a prime. No division by a possibly zero residue is used.

All sixteen final polynomials factor exactly over $\Z$. Their cubic unit factors yield nine records with a root in $(0,1)$:
\begin{center}\begin{tabular}{cccl}\toprule
Family&$(a,b,T)$&Cubic polynomial&$|N(1+z^a)|$\\\midrule
D&$(1,3,8)$&$z^3+z^2-1$&1\\
D&$(1,5,7)$&$z^3+z^2-1$&1\\
D&$(2,1,5)$&$z^3+z-1$&1\\
D&$(2,3,5)$&$z^3+z-1$&1\\
D&$(4,1,9)$&$z^3+z^2-1$&1\\
D&$(4,2,5)$&$z^3+z^2-1$&1\\
D&$(4,3,10)$&$z^3+z^2-1$&1\\
P&$(1,4,-2)$&$z^3+z^2-1$&1\\
P&$(4,1,2)$&$z^3+z^2-1$&1\\\bottomrule
\end{tabular}\end{center}
Every other cubic factor lacks the required interval root. Every listed factor makes $1+z^a$ an ordinary unit, contradicting the hypothesis. This proves Theorem~\ref{cr:thm:nonunit}.

The first sieve uses Python integer quotient rings. A separate implementation constructs the geometric polynomials densely and uses FLINT modular powering; all 10,059 recorded prime decisions agree. The sixteen characteristic-zero factorizations are checked independently with SymPy, including irreducibility of all 33 factors, exact reconstruction, real-root counts and unit norms. The certificate records the complete parameter list, all modular decisions and the exact factorization identities.

\subsection{Dependent and diagonal cubic points}
By the scalar reduction in Proposition~\ref{low:scalar}, a multiplicatively dependent complementary pair with $0<u<1$ has the form
\[
 u=v^r,\qquad1-u=v^s,\qquad \gcd(r,s)=1,\quad r,s>0,
 \quad \Q(u)=\Q(v).
\]
The generator is obtained by Bezout from $u,1-u$. If $r=s$ then $u=1/2$. Otherwise $v$ is an ordinary unit, satisfying the monic trinomial $v^r+v^s-1=0$.

Ljunggren's theorem~\cite{Ljunggren} gives one irreducible noncyclotomic factor. For coprime exponents the only possible cyclotomic factor is $Z^2-Z+1$: two unit complex numbers summing to one must be the primitive sixth roots, and Bezout recovers the original root. It is simple, since the two derivative terms have unequal magnitudes. With $r<s$, the noncyclotomic degree is $s$ or $s-2$. A coordinate of degree at most three therefore has $s\le5$.

Exact factorization leaves, apart from the rational and quadratic entries, precisely
\begin{center}\begin{tabular}{cll}\toprule
Field generator&Possible coordinates&Complementary exponent pairs\\\midrule
$v^3+v-1=0$&$v,v^3$&$(1,3),(3,1)$\\
$v^3+v^2-1=0$&$v,v^5,v^2,v^3$&$(1,5),(5,1),(2,3),(3,2)$\\\bottomrule
\end{tabular}\end{center}
The extra factorization is $Z^5+Z-1=(Z^2-Z+1)(Z^3+Z^2-1)$. The generator always denotes the unique real root in $(0,1)$. The two cubic discriminants are $-31$ and $-23$, both squarefree, so the fields are distinct and both have signature $(1,1)$.

If $\beta\ne0$ and one complementary pair is dependent, the Nahm equations put the other pair in its rational multiplicative span, hence make it dependent too. If $\beta=0$, both scalar relations hold directly. Each coordinate degree divides three, so the quadratic scalar entries cannot occur; at least one coordinate is cubic. Thus the other coordinate is either $1/2$ or one of the scalar arguments in the same cubic field.

\subsubsection{Regulator cancellation is exhausted}
Modularity implies torsion of $[x]+[y]$ by Calegari--Garoufalidis--Zagier~\cite[Theorem 1.8]{CGZ}; its Bloch--Wigner regulator must vanish. Let $w$ be the upper-half-plane root of the appropriate cubic. Rigorous complex-ball evaluation gives
\begin{center}\begin{tabular}{ccl}\toprule
Polynomial&Argument&Strict regulator interval\\\midrule
$w^3+w-1=0$&$w$&$(0.79,0.80)$\\
&$w^3$&$(-0.80,-0.79)$\\\midrule
$w^3+w^2-1=0$&$w$&$(0.47,0.48)$\\
&$w^5$&$(-0.48,-0.47)$\\
&$w^2$&$(-0.95,-0.94)$\\
&$w^3$&$(0.94,0.95)$\\\bottomrule
\end{tabular}\end{center}
Terminating decimals are rational endpoints. These bounds exclude equal-coordinate pairs, a rational partner (whose regulator is zero), and every cross-pair of unequal regulator magnitudes. The two remaining cancellations are exact complementary identities $D(1-z)=-D(z)$, not inferred numerical zeros. Up to coordinate exchange, only
\begin{equation}\label{cr:eq:reflection}
 (x,y)=(v^r,v^s),\qquad(r,s)=(1,3),(2,3),(1,5)
\end{equation}
remain. This already handles all other diagonal and dependent configurations.

\subsection{Three reflection families and twelve lines}
For~\eqref{cr:eq:reflection}, the Nahm equations give
\begin{equation}\label{cr:eq:Ab}
 A_b=\begin{pmatrix}(s/r)(1-b)&b\\b&(r/s)(1-b)\end{pmatrix},\qquad b\in\Q,\quad b<1/2.
\end{equation}
The determinant is $1-2b$. Both signs of the coupling and the diagonal value $b=0$ are included.

Write the shifts as $(B,C)$, reserving $c_0$ for the scalar. With the logarithmic normalization of Section~\ref{cu:sec:shifts},
\[
 \log f_{A_b,(B,C),c_0}(e^{-\varepsilon})\sim\alpha_0/\varepsilon+
 \log\beta_0+(\lambda_1-c_0-1/12)\varepsilon+\lambda_2\varepsilon^2+\cdots.
\]
The Vlasenko--Zwegers expansion~\cite{VZ} implies the necessary conditions
\begin{equation}\label{cr:eq:lams}
 \lambda_1\in\Q,\qquad\lambda_2=0.
\end{equation}
The two nonconstant cubic-field coefficients of $\lambda_1$ have denominators that do not vanish for $b<1/2$. Put $L=rB+sC$ and $V=rB-sC$. After removing nonzero rational constants their numerators are $p_i(L)V+q_i(L)b$, $i=1,2$, with
\begin{align*}
(r,s)=(1,3):\quad &p=(93(L-2),\ 62(L-1)),\\
 &q=(3(L-2)(31L+13),\ (L-1)(31L+150));\\
(r,s)=(2,3):\quad &p=(23(L+1),\ 69(L-1)),\\
 &q=(-(L+1)(23L-174),\ 3(L-1)(23L+82));\\
(r,s)=(1,5):\quad &p=(71001L-211876,\ 31556L-76636),\\
 &q=(52969L^2-36939L-217780,\
 19159L^2+41811L-159340).
\end{align*}
In each case $p_1,p_2$ have no common root. Their determinants $p_1q_2-p_2q_1$ are
\begin{align*}
 (1,3):&\quad-2883(L-4)(L-2)(L-1),\\
 (2,3):&\quad3174(L-2)(L-1)(L+1),\\
 (1,5):&\quad-6350645(L-7)(7L-32)(7L-12).
\end{align*}
Away from these roots, $b=V=0$. At a root, the remaining equation determines $V$ as a rational multiple of $b$. This exhausts four lines for each family:
\begin{center}\begin{tabular}{ccc}\toprule
Family&Fixed $L$&$V/b$\\\midrule
$(1,3)$&$1,\ 2,\ 4$&$-44/31,\ -106/31,\ -137/31$\\
$(2,3)$&$-1,\ 1,\ 2$&$-59/23,\ -151/23,\ -128/23$\\
$(1,5)$&$12/7,\ 32/7,\ 7$&$-32/23,\ -147/23,\ -8376/1127$\\\bottomrule
\end{tabular}\end{center}
The fourth line in each family is $b=0$, $B=(s/r)C$, with $C$ unrestricted. Intersections of lines cause no omission.

\subsubsection{Exact second-coefficient exclusion}
For each diagonal line, substitute $C=t$, $B=(s/r)t$ into $\lambda_2$. The three cubic-field coefficient numerators have gcd one in $\Q[t]$. The certificate gives explicit multipliers whose linear combination is one.

For each of the nine other lines, $B,C$ are affine functions of $b$. Put $X=x/(1-x)$, $Y=y/(1-y)$ and
\[
 D=\det(A_b+\operatorname{diag}(X,Y))
 =(1-b)\bigl(2+(s/r)Y+(r/s)X\bigr).
\]
The second factor is a nonzero field element, being positive at the distinguished embedding. The covariance numerators are affine in $b$. The Gaussian formulas for the second logarithmic coefficient therefore imply
\begin{equation}\label{cr:eq:degree10}
 (1-b)^6\lambda_2\in K[b],\qquad\deg_b\big((1-b)^6\lambda_2\big)\le10.
\end{equation}
For completeness, a Gaussian moment of total degree $2k$ has denominator $D^k$ and numerator of degree at most $k$ in $b$. In the expansion through $U_4$, the largest shift degree is four, in $U_1^4$, and $k\le6$. Multiplication by $D^6$ gives degree at most $4+k+(6-k)=10$. The subtraction $\lambda_1^2/2$ has the same bound.

Eleven exact rational evaluations at $b=0,-1,\ldots,-10$ thus determine the polynomial, with a twelfth evaluation at $b=1/3$ as a check. These are exact interpolation certificates with a proved degree bound. There is no fitted or conjectured cutoff. The three field coefficient numerators have monic gcd
\[
 (b-1)^4
\]
on every one of the nine lines. Explicit Bezout multipliers express this gcd as their linear combination. The only possible common zero is $b=1$, outside $b<1/2$. The three diagonal lines have already been excluded. This proves Theorem~\ref{cr:thm:dependent}.

The coefficients are computed by the Gaussian formulas in Section~\ref{cu:sec:shifts}. The first-coefficient formulas are certified by six exact specializations: after multiplication by $(1-b)^3$, their degree in $b$ is at most five. The certificate contains all twelve Bezout identities, the interpolation data and the strict regulator intervals.

\subsection{The cubic classification}
\begin{theorem}[Complete cubic modular classification]\label{cr:classification}
Let $A$ be a rational symmetric positive-definite matrix, and let
$\mathbf B\in\Q^2$ and $c_0\in\Q$. Suppose its distinguished Nahm point
has joint degree three. Then $f_{A,\mathbf B,c_0}$ is modular if and only
if, after a simultaneous exchange of the coordinates of $A$ and
$\mathbf B$, it is one of the twelve triples in Table~\ref{cu:tab:shifts}.
In particular there are four matrix classes, represented by
$A_6,A_7,A_8,A_9$, and the joint field of every modular cubic datum is
$\Q(\zeta_7+\zeta_7^{-1})$.
\end{theorem}
\begin{proof}
Theorem~\ref{cr:thm:dependent} excludes diagonal matrices and dependent
complementary pairs. Suppose next that one of $x,1-x,y,1-y$ is not an
ordinary unit. Lemma~\ref{qcov:entry} then places the datum in one of the
canonical families \eqref{cr:eq:D} and \eqref{cr:eq:P}, preserving the
joint field, up to coordinate exchange and simultaneous complementation.
The canonical generator therefore has degree three. Unequal exponents
are excluded by Theorem~\ref{cr:thm:nonunit}. Equal exponents in the plus
family are impossible, and the direct family is excluded by
Proposition~\ref{cu:prop:equal}; the latter includes separate coefficient
certificates for both exponent orders and for the inverse matrices.

The remaining point has four ordinary units and independent complementary
pairs. Theorem~\ref{cu:thm:main} gives exactly the four displayed matrix
classes, and Theorem~\ref{cu:thm:shifts} gives exactly their twelve modular
shifts. Its product identities prove sufficiency. Simultaneous
complementation was used only for the arithmetic coverage; the proof has
not assumed that inversion preserves modularity.
\end{proof}

%% file: sections/quartic.tex
\section{The quartic classification}\label{qt:section}
This section classifies modular data whose distinguished Nahm point has joint degree exactly four. The degree hypothesis is used throughout the arithmetic reductions. The global degree bound is supplied by Theorem~\ref{deg:thm:main}.

\providecommand{\qtmat}[3]{\begin{pmatrix}#1&#2\\#2&#3\end{pmatrix}}
\begin{theorem}\label{qt:classification}
Let $A$ be a symmetric positive definite rational $2\times2$ matrix, let $\boldsymbol B\in\Q^2$ and $c_0\in\Q$, and let $(x,y)\in(0,1)^2$ be the distinguished Nahm point. Suppose $[\Q(x,y):\Q]=4$. Then $f_{A,\boldsymbol B,c_0}$ is modular if and only if, up to simultaneous interchange of the two coordinates, the datum is one of
\begin{center}
\renewcommand{\arraystretch}{1.5}
\begin{tabular}{ccc}\toprule
$A$&$\boldsymbol B$&$c_0$\\\midrule
$\left(\begin{smallmatrix}4&1\\1&1\end{smallmatrix}\right)$&$(0,1/2)^t$&$1/120$\\
$\left(\begin{smallmatrix}4&1\\1&1\end{smallmatrix}\right)$&$(2,1/2)^t$&$49/120$\\
$\left(\begin{smallmatrix}1/3&-1/3\\-1/3&4/3\end{smallmatrix}\right)$&$(1/2,0)^t$&$1/120$\\
$\left(\begin{smallmatrix}1/3&-1/3\\-1/3&4/3\end{smallmatrix}\right)$&$(-1/6,2/3)^t$&$3/40$\\\bottomrule
\end{tabular}
\end{center}
\end{theorem}

The proof has three arithmetic stages. Lemma~\ref{qcov:entry} places every non-diagonal nonunit point in a canonical polynomial family while preserving its joint field. We prove an explicit exponent bound and exclude that finite set. Dependent and diagonal configurations are treated separately. The remaining four-unit points are classified using unit regulators, height bounds and exact finite pairings. Finally, the coefficient conditions of Section~\ref{cls:sec:method} determine all rational shifts, and published identities prove sufficiency.

\input{sections/quartic_nonunit_bounds.tex}

\input{sections/quartic_nonunit_struct.tex}

\input{sections/quartic_nonunit_small.tex}

\input{sections/quartic_nonunit_four.tex}

\input{sections/quartic_nonunit_filters.tex}

\input{sections/quartic_nonunit_census.tex}
\input{sections/quartic_dependent.tex}
\input{sections/quartic_totally_real.tex}

\input{sections/quartic_mixed_odd.tex}
\input{sections/quartic_mixed_even.tex}
\input{sections/quartic_shifts.tex}

\subsection{Completion of the quartic theorem}
\begin{proof}[Proof of Theorem~\ref{qt:classification}]
Suppose first that the Nahm sum is modular and its distinguished point has joint degree four. Theorem~\ref{qt:dep:main} excludes diagonal matrices and dependent complementary pairs. Thus the coupling is nonzero and $x,1-x$ are multiplicatively independent. Corollary~\ref{qt:nonunit-to-units}, which applies under the nonzero-coupling hypothesis, gives that $x,1-x,y,1-y$ are ordinary units.

The field $K=\Q(x,y)$ has a real embedding. Its signature is therefore $(4,0)$ or $(2,1)$. Theorem~\ref{qt:tr:thm:quartic-exclusion} excludes the totally real case. In signature $(2,1)$ the positive-unit group has rank two. Its finite saturation index is odd or even, so Theorems~\ref{qt:mo:thm:list} and~\ref{qt:me:main} put $A$ in the eight arithmetic classes, up to coordinate interchange and inversion. Proposition~\ref{qt:sh:finite} checks both orientations and leaves exactly the four displayed triples.

Conversely, the exact models in Section~\ref{qt:shifts} show that each displayed matrix has a distinguished point of joint degree four. The cited modular identities prove modularity for each of the four shift and scalar choices. This proves the theorem.
\end{proof}

\input{sections/degree_four_role.tex}

%% file: sections/quartic_nonunit_bounds.tex
\subsection{Nonunit quartic exclusion}\label{qt:nonunit}
We work with the equations $D$ and $P$ of Lemma~\ref{qcov:entry}.
Throughout this subsection $K=\Q(z)$ has degree four, $z\in(0,1)$ is an
ordinary unit, $a,b>0$, $g=\gcd(a,b)$, $M=a/g$, $N=b/g$, and
\begin{equation}\label{qt:nu:primitive}
 \gcd(M,N)=1,\qquad\gcd(g,T)=1.
\end{equation}
The two binomial factors have equal, nonempty finite prime-ideal support.
In $D$ we require $T>0$; in $P$ positive definiteness requires
\begin{equation}\label{qt:nu:pdplus}
 gT>b^2\quad\hbox{or}\quad gT<-a^2.
\end{equation}
The argument proceeds from explicit parameter bounds to structural
exclusions and an exact modular-factor certificate. It makes no signature
assumption.
\subsubsection{A smaller bound for the last exponent}
We use the absolute logarithmic Weil height $h$ and the Mahler measure
$\mathcal M(f)$ of a monic minimal polynomial, so
$h(z)=\log\mathcal M(f)/\deg(f)$.

\begin{lemma}[Certified small-degree height gap]\label{qt:nu:gap}
Every non-root-of-unity algebraic unit of degree at most four satisfies
\[
 \mathcal M(f)>13/10,\qquad h(z)>1/16.
\]
\end{lemma}
\begin{proof}
If $\mathcal M(f)\le13/10$ and $d=\deg f$, its integer coefficients satisfy
$|f_j|\le\lfloor(13/10)\binom dj\rfloor$, with leading coefficient one
and constant coefficient $\pm1$. This gives a finite coefficient box.
The certificate script \texttt{check\_height\_gap.py} factors every polynomial
in that box exactly. For each irreducible noncyclotomic polynomial, FLINT
complex root balls rigorously certify that its Mahler measure exceeds
$13/10$. There are no unresolved comparisons. Cyclotomic exceptions are
certified by exact divisibility into $X^n-1$.
\begin{center}
\begin{tabular}{rrrrr}
\toprule
$d$ & Box & Irreducible & Cyclotomic & Certified $\mathcal M>13/10$\\
\midrule
1&2&2&2&0\\
2&10&7&3&4\\
3&98&76&0&76\\
4&3630&3139&4&3135\\
\bottomrule
\end{tabular}
\end{center}
Rational Taylor bounds certify $e^{1/4}<13/10$, proving the height claim.
This is a finite proof with interval arithmetic, not a floating-point
minimum estimate.
\end{proof}

\begin{lemma}[Height of a rational function]\label{qt:nu:heightmap}
Let $P,Q\in\Z[X]$ have degrees at most $D$ and coefficient absolute sums
at most $L$. If $Q(\alpha)\ne0$, then
\[
 h(P(\alpha)/Q(\alpha))\le D h(\alpha)+\log L.
\]
\end{lemma}
\begin{proof}
Homogenize both polynomials to degree $D$. At each nonarchimedean place,
integral coefficients give the bound $\max(1,|\alpha|_v)^D$;
at each archimedean place multiply this bound by $L$.
Sum the local inequalities for the projective height of
$[P(\alpha):Q(\alpha)]$, using the product formula. Common factors of
$P,Q$ do not invalidate this upper bound. The same standard inequality
is recorded in \cite[Section 2]{SSS}.
\end{proof}

\begin{theorem}[Last-exponent bound]\label{qt:nu:Tbound}
For either (D) or (P),
\begin{equation}\label{qt:nu:newT}
 |T|<\frac{\max(a^2,b^2)}g+12\max(M,N).
\end{equation}
In particular, if $a,b\le C$, the bound
$|T|<C^2+12C$ is valid.
\end{theorem}
\begin{proof}
Apply Lemma~\ref{qt:nu:heightmap} to
$P=(1\pm X^b)^N$, $Q=(1+X^a)^M$.
Their maximum degree is
$D=\max(bN,aM)=\max(a^2,b^2)/g$, and their maximum coefficient absolute
sum is $2^{\max(M,N)}$. Since $P(z)/Q(z)=z^T$,
\[
 |T|h(z)\le D h(z)+\max(M,N)\log2.
\]
Lemma~\ref{qt:nu:gap} and $16\log2<12$ prove the result. The latter strict
inequality is also certified by a rational Taylor bound in the script.
\end{proof}

The primitive-divisor estimate in Section~\ref{qt:nu:cutoffsection} makes the resulting parameter box explicit.

\subsubsection{The equal-exponent branch contains no quartic}
\begin{theorem}\label{qt:nu:equal}
Under \eqref{qt:nu:primitive}, neither family with $a=b$ has a generator
$z\in(0,1)$ of degree four.
\end{theorem}
\begin{proof}
For (P), $M=N=1$ and the equation forces $z^T=1$, which is impossible.
For (D), write $a=b$ and obtain
\[
 Q_{a,T}(z)=0,\qquad
 Q_{a,T}(X)=X^{a+T}+X^a+X^T-1,
 \qquad\gcd(a,T)=1.
\]
If $a=T$, both are one and $Q=X^2+2X-1$.
Otherwise this is a quadrinomial with four distinct terms.

Its cyclotomic part is $X^2+1$ if $4\mid a+T$, and one otherwise.
Indeed, writing $n=a+T$, its reciprocal is
$Q^*=1+X^a+X^T-X^n$. A root of unity common to $Q,Q^*$ must satisfy
$\zeta^n=1$ and $\zeta^a+\zeta^T=0$, hence
$\zeta^{2a}=-1$. Coprimality then makes $\zeta$ a primitive fourth root.
Conversely such roots occur exactly when $4\mid n$.
They are simple because, at such a root,
\[
 \zeta Q'(\zeta)=n+(a-T)\zeta^a\ne0.
\]

Mills's corrected theorem \cite[Theorem 2]{Mills} states that the
noncyclotomic part is irreducible, apart from four explicit forms.
For the present signs and exponents, the only possible exceptional form
is $X^{8r}+X^{7r}+X^r-1$; coprimality forces $r=1$.
Here
\[
 Q_{1,7}=(X^2+1)(X^3+X^2-1)(X^3-X+1),
\]
so there is no quartic factor. In the nonexceptional cases, the degree
of the noncyclotomic part is $n$ if $4\nmid n$, and $n-2$ if $4\mid n$.
Neither value can be four: the two required possibilities would be
$n=4$ with $4\nmid n$, or $n=6$ with $4\mid n$.
\end{proof}

\paragraph*{Independent exact check with a proved cutoff}
This second proof does not require Mills's theorem. Suppose a quartic
minimal polynomial $f$ divides $Q_{a,T}$ and has a root $z\in(0,1)$.
Exchange $a,T$ so $a\le T$; the case $a=T=1$ is already excluded.
By Jensen's inequality and the coefficient Euclidean norm,
$\mathcal M(Q_{a,T})\le2$. Mahler measure is multiplicative and a nonzero
integer polynomial has measure at least one. Thus $\mathcal M(f)\le2$, and
\[
 f=X^4+c_3X^3+c_2X^2+c_1X+c_0,
 \quad |c_3|,|c_1|\le8,\quad |c_2|\le12,\quad c_0=\pm1.
\]
Also $f(1),f(-1)\in\{\pm1,\pm2\}$, since $Q(1)=2$ and
$Q(-1)=-2$ for coprime $a,T$.
On $[0,1]$, $|f'|\le4+24+24+8=60$. As $f(1)$ is a nonzero integer,
the mean value theorem gives $z\le59/60$.
Writing $q=z^a$, $r=z^T\le q$, the equation $qr+q+r=1$ implies
$q\ge\sqrt2-1>2/5$. Hence
\[
 a(-\log z)<\log(5/2)<1,\qquad a<60.
\]
Furthermore
\[
 z^T=\frac{1-z^a}{1+z^a}\ge\frac{1-z}{1+z}\ge\frac1{119},
 \qquad T(-\log z)\le\log119<5,
\]
so $T<300$. All needed comparisons can be proved by rational Taylor
bounds; the finite check uses the looser, also valid, box
$1\le a\le75$, $a\le T\le455$.

The exact checker enumerates all 14,450 coefficient choices above.
The endpoint filter leaves 260; 251 are irreducible quartics. For each
it tests all 19,192 coprime unequal exponent pairs in the larger box by
reduction modulo $f$. All 4,817,192 remainders are nonzero.
No root isolation or approximate zero test is used, since a superset
of the possible quartic divisors is checked.

\subsubsection{An odd-prime restriction on the remaining slopes}
\begin{lemma}\label{qt:nu:valuation}
Assume $M\ne N$. Let $\mathfrak p$ be a prime in the common support
of (D) or (P), of odd residue characteristic $p$, and normalize
$v_{\mathfrak p}$ to have value group $\Z$. Put
$e=v_{\mathfrak p}(p)\le4$, $h_0=\min(M,N)$ and $d_0=|M-N|$.
Then
\[
 p\mid h_0,\qquad d_0\mid e\,v_p(h_0),
 \qquad d_0\le4\log_3 h_0.
\]
In (D), an odd common prime also requires $N$ even and $M$ odd.
In (P), it requires $M,N$ both odd.
\end{lemma}
\begin{proof}
Set $w=z^g$. The residue relations at an odd common prime force
$\bar w=-1$: the order divides $2\gcd(M,N)=2$, and is not one.
They give the stated parities. Let $u=-w$, so $s=v(1-u)>0$.
Up to signs, the valuations on the two sides are
$v(1-u^N)$ and $v(1-u^M)$, with
\[
 N v(1-u^N)=M v(1-u^M).
\]
Assume $N>M$. If $p\nmid M$, then $v(1-u^M)=s$ while
$v(1-u^N)\ge s$, a contradiction. Thus $p\mid M$ and $p\nmid N$;
the equality becomes $Ns=Mv(1-u^M)$ and implies $M\mid s$.
In particular $s\ge M\ge p>e/(p-1)$.
The binomial expansion now gives the valid local lifting formula
\[
 v(1-u^M)=s+e\,v_p(M).
\]
For completeness, when $v(u-1)>e/(p-1)$, the linear term in
$(1+(u-1))^p-1$ has strictly smaller valuation than every other term;
iteration and the case of an exponent prime to $p$ give the formula.
Consequently $(N-M)(s/M)=e\,v_p(M)$.
The case $M>N$ is identical after exchanging the two indices.
Finally $v_p(h_0)\le\log_3 h_0$.
\end{proof}

Thus slopes with $|M-N|>4\log_3\min(M,N)$ can have nonunit support
only above 2. This is a necessary condition, not an exclusion of those
dyadic cases. The lemma avoids the unjustified use of an unramified
lifting formula: its hypotheses first force $s$ beyond the ramification
threshold.

\paragraph{The dyadic cases}
The following table gives necessary conditions at a common prime above
2. It is a local case division, not an existence statement for any row.
Write $e=v(2)\le4$, $w=z^g$, and assume $M\ne N$.
\begin{center}
\small
\begin{tabular}{p{0.11\linewidth}p{0.23\linewidth}p{0.55\linewidth}}
\toprule
Family & Parities & Necessary condition\\
\midrule
(D)&$M,N$ odd&$\min(M,N)\in\{1,3\}$.\\
(D)&$N$ odd, $M$ even&Either $N\in\{1,3\}$, or $M=2$, $e=4$ and $N\ge5$.\\
(D)&$N$ even, $M$ odd&$M>N$. Either $N\in\{2,4\}$, or
$M-N\mid e\,v_2(N)$.\\
(P)&$M,N$ odd&Impossible.\\
(P)&One even, one odd&The odd index is 1 or 3, or the even index is 2
and $e=4$ with odd index at least 5.\\
\bottomrule
\end{tabular}
\end{center}
\begin{lemma}\label{qt:nu:dyadic}
The table is exhaustive for a common dyadic prime in degree at most four.
\end{lemma}
\begin{proof}
The residue relations give $\bar w=1$. Put
$s=v(1+w)>0$, $t=v(1-w)>0$. Odd powers preserve the valuation of
$1-w$ or $1+w$, respectively.
If both indices are odd in (D), the equality $Nt=Ms$ gives
$t=Mk$, $s=Nk$ for an integer $k>0$. Since $M\ne N$, the valuations
differ; their minimum is $v(2)=e$, so $k\min(M,N)=e\le4$.
For (P) with both indices odd, the equality is $Ns=Ms$, impossible.

For the mixed-parity (P) row and the (D) row with $N$ odd, let $r$ be
the even index and $o$ the odd one. The valuation equality gives
\[
 v(1-\epsilon w)=rk,\qquad v(1+w^r)=ok
\]
for $k\ge1$ and an appropriate $\epsilon\in\{1,-1\}$.
Here $s_0=rk\ge2$. If $s_0>e/2$, the identity
$w^2-1=(w-1)(w+1)$ gives $v(w^2-1)>e$. Consequently
$v(1+w^r)=e$, and $ok=e$ forces $o\in\{1,3\}$.
The only other possibility is $e=4$, $s_0=2$, $r=2$, $k=1$.
Writing $w=\epsilon+\delta$ with $v(\delta)=2$ gives
$v(1+w^2)\ge4$. It equals the odd integer $o$, so $o\ge5$.
This retains, rather than discards, the ramified cancellation case.

Finally take (D) with $N$ even and $M$ odd. Since $1+w$ divides
$1-w^N$, the valuation equality forces $M>N$. Write
$s=Nk$ and $v(1-w^N)=Mk$. If $s\le e$, then $N\in\{2,4\}$.
If $s>e$, then $t=e$, and the usual lifting formula applied after
squaring gives
\[
 v(1-w^N)=s+e\,v_2(N).
\]
Thus $(M-N)k=e\,v_2(N)$, as required.
\end{proof}

\begin{corollary}[Large reduced slopes]\label{qt:nu:slopes}
Suppose $\min(M,N)\ge5$ and the configuration is nonunit.
In (P), both indices must be odd and
\[
 |M-N|\le4\log_3\min(M,N).
\]
In (D), $N$ must be even and $M$ odd. If $N>M$, no dyadic prime
occurs and $N-M\le4\log_3 M$; if $M>N$, then
$M-N\le4\log_2 N$. The divisibility conditions in
Lemma~\ref{qt:nu:valuation} and the table are stronger than these inequalities.
\end{corollary}
\begin{proof}
There is at least one prime in the common support. Combine
Lemma~\ref{qt:nu:valuation} with the exhaustive dyadic alternatives.
\end{proof}

\subsubsection{Reciprocal quartics are impossible}\label{qt:nu:reciprocal-section}
\begin{theorem}\label{qt:nu:reciprocal}
Neither (D) nor the positive-definite family (P) admits a reciprocal
quartic generator. This includes all reduced slopes, not just $N=M+1$,
and does not require the nonunit condition.
\end{theorem}
\begin{proof}
A reciprocal minimal polynomial makes $z^{-1}$ a conjugate. Evaluating
(P) at $z$ and $z^{-1}$ gives
$gT=(b^2-a^2)/2$, strictly between $-a^2$ and $b^2$, contrary to
\eqref{qt:nu:pdplus}.

For (D), evaluation at $z^{-1}>1$ forces $N$ even and $M$ odd.
Comparison gives $2T=g(N^2-M^2)$, so $N>M$.
Since $N^2-M^2$ is odd, primitivity gives
\[
 g=2,\qquad T=N^2-M^2,\qquad a=2M,\qquad b=2N.
\]
Negative real roots are impossible: the left side of (D) is nonnegative
there and the right side negative. On $(0,1)$ the logarithmic difference
\[
 N\log(1-r^{2N})-(N^2-M^2)\log r-M\log(1+r^{2M})
\]
is strictly decreasing. Hence the reciprocal quartic's conjugates must
be $z,z^{-1},\zeta,\zeta^{-1}$, with $|\zeta|=1$ and $\zeta$ nonreal.

We next prove $z>2/3$ uniformly. Suppose instead $z\le2/3$, and put
$q=z^2\le4/9$, $\ell=-\log z>1/3$; the last inequality follows by
integrating $1/t$ from 1 to $3/2$. The defining equation gives
\[
 (N^2-M^2)\ell
 =M\log(1+q^M)-N\log(1-q^N)
 \le Mq^M+\frac{Nq^N}{1-q^N}.
\]
For $M\ge1$ the first term is at most $4/9$, since the ratio of
successive terms is at most $8/9$. For $N\ge2$ the second term is at
most $32/65$, since its successive ratio is less than $2/3$.
The right side is therefore at most $548/585<1$, whereas
$(N^2-M^2)\ell\ge3\ell>1$. This contradiction proves $z>2/3$.

Let $f\in\Z[X]$ be the monic quartic minimal polynomial. Its roots now give
\[
 0<|f(1)|
 =(z+z^{-1}-2)|1-\zeta|^2
 <\left(\frac23+\frac32-2\right)4
 =\frac23<1.
\]
This is impossible for a nonzero integer $f(1)$.
\end{proof}

The proof uses only the integral value $f(1)$ and the uniform real
inequality. In particular it does not rely on a table of Salem numbers,
an unevaluated height estimate, or a finite search in $M,N$.

\subsubsection{An explicit primitive-divisor cutoff}\label{qt:nu:cutoffsection}
A unit is called reciprocal when its inverse is a conjugate; for its
monic minimal polynomial this is equivalent to $f^*=\pm f$.
The reciprocal exclusion matters a second time: for a nonreciprocal
irreducible polynomial, no conjugate lies on the unit circle, and the
resultant with its reciprocal gives a quantitative separation.
The cyclotomic-norm approach is part of the primitive-divisor literature;
compare \cite{Stewart}. All constants needed here are proved below.

\begin{theorem}\label{qt:nu:explicitC}
If $\alpha$ is a nonreciprocal algebraic unit of degree at most four,
then $\alpha^n-1$ has a primitive prime-ideal divisor for every $n>4830$.
Here primitive means that the reduction of $\alpha$ has multiplicative
order exactly $n$ at that prime. The cutoff is sufficient, not asserted
optimal.
\end{theorem}

\paragraph{The nonarchimedean upper bound}
\begin{lemma}\label{qt:nu:localC}
Let $\alpha$ be a non-root-of-unity unit of degree at most four. If
$n>30$ and $\alpha^n-1$ has no primitive prime-ideal divisor, then
\[
 |N_{K/\Q}\Phi_n(\alpha)|\le n^4.
\]
\end{lemma}
\begin{proof}
For a prime $\mathfrak p$ above $p$, let $e,f$ be its ramification and
residue degrees and let $m$ be the residue order of $\alpha$.
Reduction of cyclotomic polynomials shows that
$\mathfrak p\mid\Phi_n(\alpha)$ only if $n=mp^k$.
The case $k=0$ is primitive. For $k\ge1$, put
$u=\alpha^m$ and $s_j=v(u^{p^j}-1)>0$. Then
\[
 v(\Phi_{mp^k}(\alpha))=s_k-s_{k-1}.
\]
The binomial expansion gives increments at most $e$, unless
$s_{k-1}=e/(p-1)$. Before that threshold valuations multiply by $p$;
after it they increase by $e$. Therefore any increment greater than $e$
requires
\[
 (p-1)p^{k-1}\le e,\qquad ef\le4,
 \qquad n=mp^k\le(p^f-1)p^k.
\]
These inequalities force $p\in\{2,3,5\}$. The maximum possible $n$
is 30: for $p=2,3,5$ the respective maxima are $30,24,20$.
The certificate checker enumerates all these finite possibilities.
Thus for $n>30$ every nonprimitive prime of $\Phi_n(\alpha)$ has
valuation at most $e$. Such primes lie over divisors of $n$. Summing
$f\,v(\Phi_n(\alpha))$ above each rational prime gives at most
$\sum ef\le4$, whence the norm divides $\operatorname{rad}(n)^4$
and is at most $n^4$.
\end{proof}

\paragraph{A resultant lower bound}
Let $d=\deg\alpha\le4$, let $f$ be its minimal polynomial, and put
$H=\mathcal M(f)>13/10$. Write its conjugates as $\alpha_1,\ldots,\alpha_d$
and set
\[
 r_i=\min(|\alpha_i|,|\alpha_i|^{-1})<1.
\]
Because $f$ is nonreciprocal and irreducible, $f$ and
$f^*(X)=X^df(1/X)$ are coprime. Thus their resultant is a nonzero
integer. In
\[
 1\le|\operatorname{Res}(f,f^*)|
 =\prod_{i,j}|1-\alpha_i\alpha_j|,
\]
extract the $d$ factors pairing each conjugate with its complex
conjugate. Their product is $H^2\prod_i(1-r_i^2)$.
Bound each of the other $d^2-d$ factors by
$2\max(1,|\alpha_i|)\max(1,|\alpha_j|)$. This gives
\[
 \prod_i(1-r_i)\ge2^{-d^2}H^{-2d}.
\]
For $\beta_i=\alpha_i$ inside the unit disk and
$\beta_i=\alpha_i^{-1}$ outside, it follows that, for every $m\ge1$,
\begin{equation}\label{qt:nu:Qbounds}
 -d^2\log2-2d\log H
 \le\log Q_m\le d\log2,
 \qquad Q_m=\prod_i|1-\beta_i^m|.
\end{equation}
No conjugate can have modulus one: otherwise its inverse is its complex
conjugate and $f,f^*$ would share a root.

For $n>1$, exactly $2^{\omega(n)-1}$ of the nonzero numbers
$\mu(n/m)$, $m\mid n$, have each sign. Cyclotomic factorization and
\eqref{qt:nu:Qbounds} consequently give
\begin{align*}
 \log|N\Phi_n(\alpha)|
 &=\varphi(n)\log H+\sum_{m\mid n}\mu(n/m)\log Q_m\\
 &\ge\varphi(n)\log H
 -2^{\omega(n)-1}\bigl(2d\log H+(d^2+d)\log2\bigr)\\
 &>\bigl(\varphi(n)-31\,2^{\omega(n)}\bigr)\log H.
\end{align*}
For the last inequality use $d\le4$ and
$\log2/\log H<27/10$, certified by the exact integer comparison
$13^{27}>2^{10}10^{27}$. Also $\log H>1/4$ by Lemma~\ref{qt:nu:gap}.
Thus Lemma~\ref{qt:nu:localC} guarantees a primitive divisor whenever
\begin{equation}\label{qt:nu:totient-test}
 \varphi(n)>31\,2^{\omega(n)}+16\log n.
\end{equation}

\paragraph{The finite check and the infinite tail}
The elementary product formula for $\varphi$ gives
\[
 \frac{\varphi(n)}{2^{\omega(n)}}\ge\sqrt{n/30}.
\]
Indeed, after dividing each prime-power contribution by its square root,
all factors for $p\ge7$ are at least one; the product of the worst
factors for $p=2,3,5$ is $1/\sqrt{30}$.
Since $\omega(n)\ge1$, \eqref{qt:nu:totient-test} follows from
$\sqrt{n/30}>31+8\log n$. This holds at $n=10^6$ because the left
side exceeds 180 and the right side is less than 143. The difference
is increasing for $n>7680$, so the entire infinite tail is proved.

For $31\le n<10^6$, use $\log n<(7/10)\operatorname{bitlength}(n)$.
The script \texttt{check\_primitive\_cutoff.py} computes $\varphi(n)$
and $\omega(n)$ by an exact sieve and checks the stronger integer test
\[
 5\varphi(n)>155\,2^{\omega(n)}
       +56\operatorname{bitlength}(n).
\]
There are 592 failures of this sufficient test, the largest being 4830.
They are recorded in the JSON output; they are not assertions that
primitive divisors actually fail. Every $n>4830$ satisfies either the
finite test or the proved infinite-tail inequality. This completes
Theorem~\ref{qt:nu:explicitC}.

\subsubsection{The numerical parameter box}
\begin{corollary}\label{qt:nu:numericbox}
In the nonunit configurations (D) and (P) of degree three or four,
the canonical reduction may use $C=4830$. In particular:
\begin{center}
\begin{tabular}{ccc}
\toprule
Family & Exponent limits & Uniform last-exponent limit\\
\midrule
(D)&$a\le2415,\ b\le4830$&$0<T<23{,}386{,}860$\\
(P)&$a,b\le2415$&$0<|T|<5{,}861{,}205$\\
\bottomrule
\end{tabular}
\end{center}
Use \eqref{qt:nu:newT}, primitivity, parity, and valuation divisibilities for
stronger parameter-specific filters.
\end{corollary}
\begin{proof}
An irreducible reciprocal cubic would have a root at $1$ or $-1$, so
cubic generators are nonreciprocal. Quartic reciprocal generators are
excluded by Theorem~\ref{qt:nu:reciprocal}. Hence Theorem~\ref{qt:nu:explicitC}
applies to all remaining generators. The primitive-divisor argument following this proof gives the displayed limits for $a,b$.
Theorem~\ref{qt:nu:Tbound} then gives the numerical limits for $T$.
\end{proof}

For clarity, the primitive-divisor step uses the equal, nonempty supports.
In (D), if $b>C$, a primitive divisor of $z^b-1$ forces $b\mid2a$;
a primitive divisor at $2a$ then forces $2a\mid b$. Thus $b=2a$,
so $M=1,N=2$, impossible at a common prime because
$v(1-z^{2a})\ge v(1+z^a)>0$ but the equation requires twice the
first valuation to equal the second. Thus $b\le C$; the same
primitive-divisor argument now gives $2a\le C$.
In (P), primitive divisors at $2a,2b$ similarly force $a=b$ if either
index exceeds $C$, but then $z^T=1$. Thus $2a,2b\le C$.

%% file: sections/quartic_nonunit_struct.tex
\subsubsection{Structural exclusions for nonunit quartics}\label{qt:nu:struct}
Throughout this addendum assume the canonical hypotheses: $K=\Q(z)$ has
degree four, $z\in(0,1)$ is an ordinary unit satisfying the corresponding
equation $F_D(z)=0$ or $F_P(z)=0$, $w=z^g$, $\gcd(M,N)=1$,
$\gcd(g,T)=1$, and the two binomial factors have equal \emph{nonempty}
finite prime-ideal support. The support condition is essential. The arguments
do not assume that $K$ is totally real.

\begin{theorem}[Structural slope exclusions]\label{qt:nu:thm:structural}
Under these hypotheses the following cases are impossible:
\begin{enumerate}
\item Either family with $M=3$ or $N=3$.
\item Family $D$ with $6\mid N$, or with $4\mid N$ and $N>M$.
\item Family $D$ with $M=2$, $N$ odd and $5\le N\le4830$;
family $P$ with one index equal to two and the other odd and between five
and $4830$.
\end{enumerate}
The first two assertions are uniform in the exponents. The last uses the
finite certificate below; the earlier bound $4830$ puts every relevant
exceptional ramification-four branch in its range.
\end{theorem}

\paragraph{Odd common primes.}
Lemma~\ref{qt:nu:valuation}, with $h=\min(M,N)$ and $d=|M-N|$, gives
\begin{equation}\label{qt:nu:oddvaluations}
 p\mid h,\qquad v_{\mathfrak p}(1+w)=hk,\qquad
 dk=e v_p(h),\qquad k\in\Z_{>0}.
\end{equation}
At such a prime $\bar w=-1$. Thus $D$ requires $N$ even and $M$ odd,
whereas $P$ requires both indices odd.

\paragraph{Three norm identities.}
Let $f=X^4+c_3X^3+c_2X^2+c_1X+\epsilon$ be the characteristic polynomial
of multiplication by $w$ on $K$, so $\epsilon=\pm1$.
Norms equal the corresponding resultants even if $w$ does not generate $K$.
Writing $\Phi_j$ for a cyclotomic polynomial gives
\begin{align}
 |N_{K/\Q}\Phi_4(w)|=1&\ \Longrightarrow\ |f(1)+f(-1)|\le10,\\
 |N_{K/\Q}\Phi_3(w)|=1&\ \Longrightarrow\ |f(1)+3f(-1)|\le16,\\
 |N_{K/\Q}\Phi_6(w)|=1&\ \Longrightarrow\ |3f(1)+f(-1)|\le16.       \label{qt:nu:norm-identities}
\end{align}
Indeed, $f(i)=(1+\epsilon-c_2)+i(c_1-c_3)$ has norm one in the first
case. In the other cases reduce $f$ modulo $X^2+X+1$ or $X^2-X+1$.
The coefficients $A,B$ satisfy $A^2-AB+B^2=1$ or $A^2+AB+B^2=1$,
respectively, giving $|A+B|\le2$ or $|A-B|\le2$ and \eqref{qt:nu:norm-identities}.

\subsubsection{Divisibility and ramification-three exclusions}\label{qt:nu:divisibility}
\paragraph{The case $D$, $N>M$, $N$ even.}
There is no common dyadic prime: $1+w$ divides $1-w^N$ for even $N$,
whereas equality requires $N v(1-w^N)=M v(1+w)$ because $M$ is odd.
Every common prime is therefore odd, and \eqref{qt:nu:oddvaluations} forces $M\ge3$ and
\[
 |N_{K/\Q}(1+w)|\ge3^M\ge27.
\]
The element $1-w$ is a unit. Indeed it divides $1-w^N$, and at an odd
common prime the residue $w=-1$ makes $1-w$ nonzero. If $4\mid N$,
the same argument shows that $1+w^2$ is a unit. Thus $|f(1)|=1$ and
the first inequality of \eqref{qt:nu:norm-identities} gives $|f(-1)|\le11$, a contradiction.
If $3\mid N$, the divisor $\Phi_3(w)$ is a unit because $\Phi_3(-1)=1$.
The second inequality gives $|f(-1)|\le17/3$, again a contradiction.

\paragraph{The case $D$, $6\mid N$, $M>N$.}
Again $\Phi_3(w)$ is a unit: at an odd common prime $w=-1$, and at a
dyadic common prime $w=1$, with $\Phi_3(1)=3$ odd. At a dyadic common
prime, $v(1+w)=Nk\ge6>e$, so $v(1-w)=e$. At odd primes $1-w$ is a
unit. Hence $|f(1)|\le2^4=16$. Any common prime gives
$v(1+w)=Nk\ge N$, by coprimality at two and by \eqref{qt:nu:oddvaluations} at odd primes.
Consequently $|f(-1)|\ge2^N\ge64$. This contradicts
$3|f(-1)|\le16+|f(1)|\le32$. The equal-index case is absent here because
the reduced indices are coprime. This proves assertion 2 of
Theorem~\ref{qt:nu:thm:structural}.

\paragraph{Why ramification index three cannot occur in these branches.}
Consider $D$ with both indices odd and $\min(M,N)=3$; $D$ with $N=3$
and $M$ even; or $P$ with mixed parity and odd index three. There is no
odd common prime, by parity. We claim every common dyadic prime would
have $e=3$ and residue degree one.

For the first case put $s=v(1+w)$, $t=v(1-w)$. The relation gives
$s=Nk$, $t=Mk$; since they differ, $e=\min(s,t)=3k\le4$, so $e=3$.
For either mixed-parity case, with even index $r$, the valuation equality
gives $v(1-\eta w)=rk$ and $v(1+w^r)=3k$, for a sign $\eta$ and $k\ge1$.
If $rk>e/2$, expansion after squaring gives $v(1+w^r)=e$, so $e=3$.
Otherwise $rk\le e/2\le2$ forces $r=2,k=1,e=4$; but then
$v(1+w^2)\ge4$, contradicting $3k=3$.

In a quartic field a prime with $e=3$ has residue degree one and leaves
exactly one further prime of local degree one above two. At that further
prime the unit $w$ reduces to the only nonzero element of $\mathbb F_2$,
namely one. Both original binomials therefore vanish there. It too is a
common prime, but has $e=1$, contradicting the necessary $e=3$.
Thus all three kinds of branch are excluded, without a height bound.

\subsubsection{The ramification-four branch}\label{qt:nu:ramfour}
In assertion 3, put $n$ for the odd index and set $u=w$ in $D$, $u=-w$
in $P$. At every finite place the common-support relation is
\begin{equation}\label{qt:nu:ramfour-valuation}
 n\,v(1-u^n)=2\,v(1+u^2).                                      
\end{equation}
No odd common prime is possible. At a common prime above two put
$s=v(1-u)$, $t=v(1+u^2)$ and $e=v(2)\le4$. As $n$ is odd, \eqref{qt:nu:ramfour-valuation}
gives $s=2k$, $t=nk$. If $s>e/2$, writing $u=1+\delta$ shows
$t=v(2+2\delta+\delta^2)=e$, impossible for $n\ge5$.
Thus $4k\le e\le4$, forcing
\[
 e=4,\quad k=1,\quad v(1-u)=2,\quad v(1+u^2)=n.
\]
The prime is unique above two and has residue degree one. All finite
support is there. Consequently
\begin{equation}\label{qt:nu:ramfour-norms}
 N(1-u)=4\delta,\quad N(1-u^n)=4\delta,\quad N(1+u^2)=2^n,
 \qquad\delta\in\{-1,1\}.                                      
\end{equation}
The first two signs agree because $(1-u^n)/(1-u)$ is positive at every
real embedding when $n$ is odd, and its norm is positive and has absolute
value one. The third norm is positive as well.

Let $f$ be the degree-four characteristic polynomial of $u$.
Then $f(1)=4\delta$ and $|f(i)|^2=2^n$.
The Gaussian-integer equation $A^2+B^2=2^n$ for odd $n$ forces
$A=s2^{(n-1)/2}$, $B=t2^{(n-1)/2}$ with $s,t\in\{-1,1\}$:
repeatedly divide both even coordinates by two until the norm is two.
Thus, with $V=2^{(n-1)/2}$ and $c=2\delta-1-\epsilon$, all candidates are
\begin{equation}\label{qt:nu:ramfour-polys}
 f=X^4+\left(c+\frac{s-t}{2}V\right)X^3
 +(1+\epsilon-sV)X^2
 +\left(c+\frac{s+t}{2}V\right)X+\epsilon,                       
\end{equation}
where $\epsilon,\delta,s,t\in\{-1,1\}$. This is sixteen polynomials per
odd $n$. Necessarily $\operatorname{Res}(f,1-X^n)=4\delta$.

For all odd $5\le n\le4830$, every one of these 38,608 patterns fails
that necessary equality modulo at least one of the following primes:
\begin{center}
\begin{tabular}{crrrrrrr}
\toprule
Prime added &3&5&7&11&13&17&19\\
Patterns remaining &20742&5642&1516&222&32&2&0\\
\bottomrule
\end{tabular}
\end{center}
The CSV records a rejecting prime and residue for each pattern. A separate
pure-Python checker verifies the complete enumeration, all 77,216 integer
coefficient identities, and all 38,608 rejecting determinants of multiplication
by $1-X^n$ in the quartic quotient algebra. It uses no FLINT or other external
package. This independently confirms the FLINT-generated certificate.
It excludes even reducible characteristic polynomials, so no irreducibility
or signature assumption has been omitted. The range is sufficient by Corollary~\ref{qt:nu:numericbox}.

\subsubsection{Reduced index three: norm exclusions}\label{qt:nu:three-norms}
The ramification-three argument already handles $D$ with both indices odd
and minimum three, $D$ with $N=3$ and even $M$, and $P$ with mixed parity
and odd index three. We now dispose of the remaining cases.

\paragraph{$D(3,1)$ and $D(1,3)$, in reduced indices.}
Only dyadic support occurs. For $D(3,1)$, the valuation relation gives
$v(1-w)=3v(1+w)$ and $v(1+w)=e$ at every common prime.
Writing $k$ for the sum of the local degrees at those primes,
\[
 |f(1)|=2^{3k},\quad |f(-1)|=2^k,\quad k\ge1.
\]
Also $\Phi_6(w)$ is a unit: it divides $1+w^3$ and is nonzero modulo
every common prime above two. The last inequality of \eqref{qt:nu:norm-identities} contradicts
$|3f(1)+f(-1)|\ge3\cdot2^{3k}-2^k\ge22$.
For $D(1,3)$ the two norm sizes are exchanged and $\Phi_3(w)$ is a unit.
The middle inequality of \eqref{qt:nu:norm-identities} gives the same contradiction.

\paragraph{$P(3,5)$ and $P(3,7)$, including exchange of indices.}
Write $n=5$ or $7$ and $d=n-3$. The only possible support is above three.
By \eqref{qt:nu:oddvaluations}, each common prime has $d\mid e$ and
$v(1+w)=3e/d$. Let $k=\sum e f_{\mathfrak p}$ over common primes.
Then $k\in\{2,4\}$ for $n=5$ and $k=4$ for $n=7$. With
$C=3^{3k/d}$, the necessary norm conditions are
\begin{equation}\label{qt:nu:three-norms-eq}
 f(-1)=\delta C,\qquad |N\Phi_6(w)|=3^k,\qquad
 N(1+w^n)=\delta C.                                             
\end{equation}
The last sign agrees with the first because $(1+w^n)/(1+w)$ has positive
norm and is a unit. Reducing $f$ modulo $X^2-X+1$ gives
\[
 A=\epsilon-c_2-c_3,\quad B=c_1+c_2-1,\quad
 A^2+AB+B^2=3^k.
\]
Thus all possibilities are determined by
\begin{equation}\label{qt:nu:three-coefficients}
 c_2=(\delta C-A+B)/3,\quad c_1=B+1-c_2,\quad c_3=\epsilon-c_2-A. 
\end{equation}
Since $3^k\le81$, the norm equation forces $|A|,|B|<11$. Enumeration
of these integer pairs, $\epsilon=\pm1$ and $\delta=\pm1$ yields 48
polynomials for $n=5$ and 24 for $n=7$. All fail the last equality in
\eqref{qt:nu:three-norms-eq}. The pure-Python checker calculates the exact integer norms as
four-by-four determinants and records all 72 differences.

\paragraph{Why these are the only remaining odd partners of three.}
For $P$ both indices odd, a dyadic common prime would give
$N v(1+w)=M v(1+w)$, which is impossible for unequal indices. If the smaller
index is one, \eqref{qt:nu:oddvaluations} allows no odd support either. If the smaller index is
three, the only odd common characteristic is three and $n-3\mid e\le4$,
so $n=5$ or $7$. For $D$ with $M=3$ and even $N>3$, the same condition
leaves only $N=4$ (larger coprime even indices have $N-3\ge5$); dyadic
support is impossible when $N>M$. The case $D(3,4)$ is already excluded
by the multiple-of-four argument. Only $D(3,2)$ remains.

\subsubsection{Reduced indices $(3,2)$: finite characters}\label{qt:nu:three-characters}
Here all support is dyadic. Put $s=v(1+w)$ and $t=v(1-w)$ at a common
prime. The equation gives $2(s+t)=3s$, so $s=2t$ and $t=e$.
The factor $\Phi_6(w)$ of $1+w^3$ is a unit. Thus
\begin{equation}\label{qt:nu:three-two-norms}
 f(1)=\delta_1 2^k,\quad f(-1)=\delta_2 2^{2k},\quad
 |N\Phi_6(w)|=1,\qquad 1\le k\le4.                             
\end{equation}
Use \eqref{qt:nu:three-coefficients}, with $C=2^{2k}$ and $A^2+AB+B^2=1$, and retain $f(1)=\delta_1 2^k$.
There are exactly eight polynomials, listed below by their coefficients
$(c_0,c_1,c_2,c_3,c_4)$ in ascending order.

The canonical equation gives
\[
 R(w)=\frac{(1-w^2)^2}{(1+w^3)^3}=z^T,\qquad R(w)^g=w^T.
\]
Also $z^T\in\Q(w)$ and $\gcd(g,T)=1$ imply $\Q(w)=\Q(z)$ by Bezout's
identity. Therefore $f$ is the minimal polynomial of $w$, and the cleared
power relation reduces at every root of $f$ modulo a prime where the
displayed denominator is nonzero.

For a prime $p\equiv1\pmod\ell$, fix a primitive $\ell$th root $b$ in
$\mathbb F_p$. Define $\chi(v)\in\mathbb F_\ell$ by
$v^{(p-1)/\ell}=b^{\chi(v)}$. Two roots $r_1,r_2$ of $f$ give rows
$(\chi(R(r_i)),\chi(r_i))$. If their determinant is nonzero modulo $\ell$,
the relation forces both $g$ and $T$ to be divisible by $\ell$, contradicting
primitivity. Here are certificates for all eight polynomials:
\begin{center}\small
\begin{tabular}{lrrrrll}
\toprule
Coefficients & $\ell$ & $p$ & $b$ & $(r_1,r_2)$ & Row 1 & Row 2\\
\midrule
$(-1,-1,1,-2,1)$ &2&29&28&$(21,23)$&$(1,1)$&$(1,0)$\\
$(-1,0,1,-3,1)$ &3&43&36&$(8,12)$&$(2,0)$&$(1,1)$\\
$(-1,2,-1,1,1)$ &2&29&28&$(18,24)$&$(0,1)$&$(1,0)$\\
$(-1,3,-1,0,1)$ &3&43&36&$(18,20)$&$(0,2)$&$(2,1)$\\
$(1,-1,1,0,1)$ &2&19&18&$(2,5)$&$(0,1)$&$(1,0)$\\
$(1,0,1,-1,1)$ &2&19&18&$(4,10)$&$(1,0)$&$(1,1)$\\
$(-1,-4,6,-6,1)$ &2&19&18&$(7,13)$&$(1,0)$&$(1,1)$\\
$(-1,6,-6,4,1)$ &2&19&18&$(3,9)$&$(0,1)$&$(1,0)$\\
\bottomrule
\end{tabular}
\end{center}
The checker verifies primality, the roots, nonzero denominators, all character
values and all determinants. This is a uniform exclusion of $D(3,2)$ for
every primitive $(g,T)$, not a search over a bounded exponent range.
Together with Section~\ref{qt:nu:three-norms} and the ramification-three argument, it completes
assertion 1 of Theorem~\ref{qt:nu:thm:structural}.

The eight character certificates and the 72 norm checks are integer-arithmetic certificates.

%% file: sections/quartic_nonunit_small.tex
\subsubsection{Reduced indices one and two}\label{qt:nu:small-pairs}
\begin{theorem}[Uniform small-pair exclusion]\label{qt:nu:thm:smallpairs}
Under the canonical nonunit quartic hypotheses of Section~\ref{qt:nu:struct},
$D$ with $(M,N)=(2,1)$ is impossible. In the positive-definite family $P$,
both $(M,N)=(1,2)$ and $(2,1)$ are impossible.
These exclusions hold for every primitive $(g,T)$ and require no upper bound
on $g$ or $T$.
\end{theorem}

Write $w=z^g\in(0,1)$ and define
\[
 R_D(w)=\frac{1-w}{(1+w^2)^2},\qquad
 R_P(w)=\frac{(1+w^2)^2}{1+w}.
\]
In the direct case $R_D(w)=z^T$, with $T>0$; set $\tau=T$.
In the plus case set $\tau=T$ for $(M,N)=(1,2)$ and $\tau=-T$ for
$(M,N)=(2,1)$. Then $R_P(w)=z^\tau$. In either case
\begin{equation}\label{qt:nu:small-power}
 R(w)^g=w^\tau,\qquad \gcd(g,\tau)=1.                            
\end{equation}
For $P$ the positive-definiteness condition is exactly
\begin{equation}\label{qt:nu:small-pd}
 \tau>4g\quad\hbox{or}\quad\tau<-g.                            
\end{equation}
Bezout's identity, applied to $z^g=w$ and $z^\tau=R(w)$, shows that
$\Q(w)=\Q(z)$. Hence $w$ has a monic irreducible quartic minimal polynomial
$f=X^4+c_3X^3+c_2X^2+c_1X+\epsilon$, with $\epsilon=\pm1$.

\paragraph{An exact norm reduction.}
The binomials $1-w$ and $1+w^2$ in $D$, or $1+w$ and $1+w^2$ in $P$,
have equal nonempty support. Their common primes must lie above two:
the residue conditions $w=1$ or $w=-1$, together with $w^2=-1$, contradict
odd characteristic. At a common dyadic prime put $u=w$ for $D$, $u=-w$
for $P$, and write
\[
 s=v(1-u),\quad t=v(1+u^2),\quad e=v(2)\le4.
\]
The valuation relation is $s=2t$. If $s>e/2$, expansion of $u=1+\delta$
gives $v(2+2\delta+\delta^2)=e$, so $t=e$ and $s=2e$.
Otherwise $4t\le e\le4$ forces $t=1,e=4,s=2$, but the same expansion
has valuation at least four, a contradiction.

Let $k=\sum e_{\mathfrak p}f_{\mathfrak p}$ over the common primes above two.
Then $1\le k\le4$, and
\begin{equation}\label{qt:nu:small-norms}
 |N(1+w^2)|=2^k,\qquad
 f(1)=\delta 2^{2k}\ (D),\quad f(-1)=\delta 2^{2k}\ (P),
 \qquad\delta=\pm1.                                             
\end{equation}
The norm of $1+w^2$ is positive for every signature. Set
$A=1+\epsilon-c_2$, $B=c_1-c_3$. Then $A^2+B^2=2^k$.
For each $k=1,2,3,4$ there are exactly four integer pairs $(A,B)$,
all in $[-4,4]^2$. Put
\begin{equation}\label{qt:nu:small-coefficients}
 c_2=1+\epsilon-A,\quad
 S_D=\delta 2^{2k}-2(1+\epsilon)+A,\quad S_P=-S_D,
 \qquad c_1=(S+B)/2,\quad c_3=(S-B)/2.                           
\end{equation}
The parities make these integral. Equations \eqref{qt:nu:small-norms}--\eqref{qt:nu:small-coefficients} give exactly
64 candidate polynomials in each family, without any height estimate.

\subsubsection{Exclusion of the small-pair polynomials}\label{qt:nu:small-exclusions}
Three complementary certificates dispose of every polynomial in \eqref{qt:nu:small-coefficients}.
The order of the tests is positivity, reciprocity, then finite characters.

\paragraph{Positivity on the real interval.}
For $f=\sum_{j=0}^4c_jX^j$, with $c_4=1$, define
\begin{align*}
 b_0&=12c_0,&b_1&=12c_0+3c_1,\\
 b_2&=12c_0+6c_1+2c_2,&b_3&=12c_0+9c_1+6c_2+3c_3,\\
 b_4&=12\sum_{j=0}^4c_j.
\end{align*}
Then
\begin{equation}\label{qt:nu:bernstein}
 12f(X)=\sum_{j=0}^4b_j\binom4jX^j(1-X)^{4-j}.                  
\end{equation}
If all $b_j$ have the same weak sign and at least one is nonzero, the
right side has a strict sign for $0<X<1$. Such a polynomial cannot have
$w$ as a root. This excludes 32 direct and 36 plus candidates.

\paragraph{Reciprocity.}
For a remaining reciprocal minimal polynomial, $w^{-1}$ is a conjugate.
The rational functions satisfy
\[
 R_D(w^{-1})=-w^3R_D(w),\qquad R_P(w^{-1})=w^{-3}R_P(w).
\]
In $D$, conjugating \eqref{qt:nu:small-power} gives $(-1)^g w^{3g+2T}=1$, impossible
because $0<w<1$ and $g,T>0$. In $P$ it gives $w^{2\tau-3g}=1$,
so $\tau=3g/2$, contrary to \eqref{qt:nu:small-pd}. This excludes four further candidates
in each family. The plus exclusion explicitly uses positive definiteness.

\paragraph{Finite characters.}
For every remaining polynomial, Section~\ref{qt:nu:small-characters} supplies two roots $r_1,r_2$
modulo a prime $p$, a prime $\ell\mid p-1$, and a primitive $\ell$th root
$b\in\mathbb F_p$. Both $r_i$ and the numerator and denominator of $R(r_i)$
are nonzero. As in Section~\ref{qt:nu:three-characters}, put
$\chi(v)=\log_b(v^{(p-1)/\ell})\in\mathbb F_\ell$.
In each row the matrix with rows $(\chi(R(r_i)),\chi(r_i))$ has nonzero
determinant. Reducing \eqref{qt:nu:small-power} therefore forces $\ell\mid g$ and
$\ell\mid\tau$, contrary to primitivity. Clearing denominators and any
negative exponent before reducing justifies the argument also for $\tau<0$.

\begin{center}
\begin{tabular}{lrrrrr}
\toprule
Family & Candidates & Positivity & Reciprocity & Characters & Remaining\\
\midrule
$D(2,1)$ &64&32&4&28&0\\
$P(1,2)$ and its swap &64&36&4&24&0\\
\bottomrule
\end{tabular}
\end{center}
The generator uses FLINT to find the modular roots. A separate standard-library
Python checker independently re-enumerates all 128 patterns, verifies 256
integer norm identities, expands every Bernstein identity, checks reciprocity,
and verifies all 52 character certificates. The largest certificate prime is 101. This proves Theorem~\ref{qt:nu:thm:smallpairs}, subject
to the stated canonical hypotheses, with no signature restriction.

\subsubsection{Small-pair character certificates}\label{qt:nu:small-characters}
Each tuple lists $(c_0,c_1,c_2,c_3)$; the leading coefficient is one.
The columns give $\ell,p,b,r_1,r_2$. Using the appropriate $R_D$ or $R_P$
from Section~\ref{qt:nu:small-pairs}, every resulting character determinant is nonzero. The corresponding character values are determined by the displayed data.
\medskip

\noindent
\begin{minipage}[t]{0.49\textwidth}
\centering\scriptsize
\textbf{Family $D$}\par\smallskip
\setlength{\tabcolsep}{2.2pt}
\begin{tabular}{lrrrrr}
\toprule
$(c_0,c_1,c_2,c_3)$ & $\ell$ & $p$ & $b$ & $r_1$ & $r_2$\\
\midrule
$(-1,1,1,2)$ & 7 & 29 & 16 & 8 & 6\\
$(1,-5,3,-4)$ & 2 & 7 & 6 & 6 & 3\\
$(-1,2,1,1)$ & 3 & 7 & 4 & 5 & 4\\
$(1,-4,3,-5)$ & 2 & 7 & 6 & 6 & 5\\
$(-1,2,-1,3)$ & 2 & 13 & 12 & 9 & 2\\
$(1,-4,1,-3)$ & 3 & 13 & 3 & 9 & 4\\
$(-1,3,-1,2)$ & 2 & 11 & 10 & 7 & 2\\
$(1,-3,1,-4)$ & 3 & 13 & 3 & 10 & 3\\
$(-1,7,2,7)$ & 3 & 7 & 4 & 4 & 3\\
$(-1,7,0,9)$ & 2 & 43 & 42 & 24 & 3\\
$(1,-11,2,-9)$ & 7 & 29 & 16 & 24 & 11\\
$(-1,9,0,7)$ & 2 & 37 & 36 & 32 & 3\\
$(1,-9,2,-11)$ & 7 & 29 & 16 & 23 & 8\\
$(-1,9,-2,9)$ & 7 & 29 & 16 & 18 & 7\\
$(-1,30,2,32)$ & 2 & 11 & 10 & 7 & 6\\
$(1,-36,4,-34)$ & 3 & 13 & 3 & 11 & 10\\
$(-1,32,2,30)$ & 2 & 59 & 58 & 24 & 22\\
$(1,-34,4,-36)$ & 3 & 13 & 3 & 6 & 4\\
$(-1,32,-2,34)$ & 2 & 7 & 6 & 5 & 2\\
$(1,-34,0,-32)$ & 2 & 37 & 36 & 30 & 24\\
$(-1,34,-2,32)$ & 7 & 29 & 16 & 24 & 4\\
$(1,-32,0,-34)$ & 2 & 11 & 10 & 5 & 2\\
$(-1,126,4,126)$ & 3 & 43 & 36 & 31 & 16\\
$(-1,126,0,130)$ & 2 & 17 & 16 & 14 & 7\\
$(1,-132,2,-128)$ & 2 & 11 & 10 & 10 & 7\\
$(-1,130,0,126)$ & 7 & 29 & 16 & 26 & 24\\
$(1,-128,2,-132)$ & 2 & 11 & 10 & 10 & 8\\
$(-1,130,-4,130)$ & 2 & 11 & 10 & 10 & 4\\
\bottomrule
\end{tabular}
\end{minipage}
\hfill
\begin{minipage}[t]{0.49\textwidth}
\centering\scriptsize
\textbf{Family $P$}\par\smallskip
\setlength{\tabcolsep}{2.2pt}
\begin{tabular}{lrrrrr}
\toprule
$(c_0,c_1,c_2,c_3)$ & $\ell$ & $p$ & $b$ & $r_1$ & $r_2$\\
\midrule
$(-1,2,1,3)$ & 2 & 11 & 10 & 6 & 5\\
$(-1,3,1,2)$ & 2 & 13 & 12 & 10 & 6\\
$(-1,1,-1,2)$ & 2 & 7 & 6 & 5 & 4\\
$(-1,2,-1,1)$ & 7 & 29 & 16 & 24 & 18\\
$(-1,9,2,9)$ & 7 & 29 & 16 & 8 & 4\\
$(-1,7,0,9)$ & 2 & 23 & 22 & 20 & 16\\
$(1,-7,2,-5)$ & 2 & 29 & 28 & 9 & 2\\
$(-1,9,0,7)$ & 3 & 43 & 36 & 34 & 14\\
$(1,-5,2,-7)$ & 2 & 29 & 28 & 15 & 13\\
$(-1,7,-2,7)$ & 2 & 7 & 6 & 5 & 2\\
$(-1,32,2,34)$ & 7 & 29 & 16 & 7 & 6\\
$(1,-30,4,-28)$ & 3 & 13 & 3 & 3 & 1\\
$(-1,34,2,32)$ & 2 & 7 & 6 & 4 & 3\\
$(1,-28,4,-30)$ & 3 & 13 & 3 & 9 & 1\\
$(-1,30,-2,32)$ & 2 & 101 & 100 & 60 & 43\\
$(1,-32,0,-30)$ & 2 & 13 & 12 & 10 & 7\\
$(-1,32,-2,30)$ & 5 & 11 & 4 & 9 & 3\\
$(1,-30,0,-32)$ & 2 & 13 & 12 & 4 & 2\\
$(-1,130,4,130)$ & 2 & 11 & 10 & 8 & 1\\
$(-1,126,0,130)$ & 7 & 29 & 16 & 10 & 6\\
$(1,-128,2,-124)$ & 2 & 19 & 18 & 17 & 8\\
$(-1,130,0,126)$ & 2 & 17 & 16 & 12 & 6\\
$(1,-124,2,-128)$ & 2 & 19 & 18 & 12 & 9\\
$(-1,126,-4,126)$ & 2 & 43 & 42 & 18 & 8\\
\bottomrule
\end{tabular}
\end{minipage}

\medskip

%% file: sections/quartic_nonunit_four.tex
\subsubsection{The direct multiple-of-four boundary}\label{qt:nu:four-boundary}
\begin{theorem}[Uniform multiple-of-four exclusion]\label{qt:nu:thm:four}
Under the canonical nonunit quartic hypotheses of Section~\ref{qt:nu:struct}, family $D$
has no solution with $4\mid N$. This holds for both orderings of $M,N$,
without any bound on $M,N,g,T$ and without a signature assumption.
\end{theorem}
The case $N>M$ was proved in Section~\ref{qt:nu:divisibility}. Suppose $M>N$, so $M$ is odd.
The argument below also covers $3\mid N$; that restriction is unnecessary.
It uses neither primitivity of $(g,T)$ nor the parameter cutoff.

Let $f=X^4+c_3X^3+c_2X^2+c_1X+\epsilon$ be the characteristic polynomial
of $w=z^g$ on $K$, with $\epsilon=\pm1$. It need not be irreducible.
Put $A=1+\epsilon-c_2$ and $B=c_1-c_3$. We have
\begin{equation}\label{qt:nu:four-norm-identities}
 N(1+w^2)=A^2+B^2,\qquad f(1)+f(-1)=4(1+\epsilon)-2A.       
\end{equation}
The first norm is positive for every signature. Since $1-w$ and $1+w^2$
divide $1-w^N$, their finite support lies among the common primes.
At odd common primes, $w=-1$ in the residue field, so both are units.
Also Section~\ref{qt:nu:struct} gives $v(1+w)=Nh$ for an integer $h\ge1$ there.

\paragraph{The dyadic possibilities.}
At a common prime above two normalize $v$ to have value group $\Z$ and set
\[
 e=v(2)\le4,\qquad s=v(1+w),\qquad t=v(1-w),\qquad r=v_2(N)\ge2.
\]
Since $M$ is odd, $v(1+w^M)=s$. The equation
$N v(1-w^N)=Ms$ and coprimality give $s=Nh\ge4$.
In particular $s\ge e$. Writing $w=-1+u$, with $v(u)=s$, shows
\[
 v(1+w^2)=v(2-2u+u^2)=e,\qquad t\ge e,
 \quad t=e\ \text{if }s>e.
\]
For $j\ge1$, $w^{2^j}-1$ has valuation greater than $e$, hence
$v(1+w^{2^j})=e$. Factoring $1-w^N$, with its odd part contributing
a unit geometric sum, therefore gives
\begin{equation}\label{qt:nu:four-valuation}
 v(1-w^N)=s+t+(r-1)e.                                      
\end{equation}
There are just two possibilities:
\begin{enumerate}
\item $s>e$: then $t=e$ and $(M-N)h=re$.
If $N\ge8$, $s\ge8$. If $N=4$, the equality
$(M-4)h=2e$ and oddness of $M-4$ force $h$ even, so again $s\ge8$.
\item $s=e$: then $s=Nh\ge4$ and $e\le4$ force
$N=e=s=4$ and $h=1$. This prime is unique above two and has residue
degree one. Equation \eqref{qt:nu:four-valuation} gives $t=M-8$. Since $t\ge4$ and $M$ is odd,
$t\ge5$ and $M\ge13$.
\end{enumerate}
These alternatives are exhaustive. In the second, total ramification
prevents coexistence with any other dyadic prime.

\subsubsection{Norm contradiction at the boundary}\label{qt:nu:four-norms}
\paragraph{No exceptional dyadic prime.}
Let $k=\sum e_{\mathfrak p}f_{\mathfrak p}$ over the common primes above
two, so $0\le k\le4$. In the first dyadic alternative, and when there is
no dyadic support, the support analysis gives
\[
 |f(1)|=2^k,\qquad N(1+w^2)=2^k.
\]
Thus $|A|\le4$ and \eqref{qt:nu:four-norm-identities} yields
\begin{equation}\label{qt:nu:four-norm-bound}
 |f(-1)|\le |f(1)|+|f(1)+f(-1)|\le16+16=32.                 
\end{equation}
If a common dyadic prime exists, its $s\ge8$ gives
$|f(-1)|\ge2^8=256$, a contradiction. If none exists, nonempty support
gives an odd common prime and $|f(-1)|\ge3^N\ge81$, again a contradiction.

\paragraph{The exceptional dyadic prime.}
Here $N=4$ and the unique prime above two has $e=4$ and residue degree one.
An odd common prime would have residue characteristic dividing the smaller
index $N=4$, by \eqref{qt:nu:oddvaluations}, so none exists. Consequently
\begin{equation}\label{qt:nu:four-exception}
 |f(-1)|=16,\qquad |f(1)|=2^{M-8}\ge32,\qquad A^2+B^2=16.   
\end{equation}
Since $\epsilon=\pm1$ and $|A|\le4$, \eqref{qt:nu:four-norm-identities} lies in the interval $[-8,16]$.
The only signed choices in \eqref{qt:nu:four-exception} whose sum lies in this interval are
\[
 f(1)=32,\quad f(-1)=-16,\quad M=13.
\]
Their sum is 16, forcing $\epsilon=1$, $A=-4$ and $B=0$.
It follows that $c_2=6$, $c_1=c_3=12$, and
\[
 f(X)=X^4+12X^3+6X^2+12X+1.
\]
Every coefficient is positive, so this polynomial has no positive root.
But $w\in(0,1)$ is a root of its characteristic polynomial, a contradiction.
This completes the proof of Theorem~\ref{qt:nu:thm:four}.

%% file: sections/quartic_nonunit_filters.tex
\subsubsection{Setup and the safe reduction test}
Put $g=\gcd(a,b)$, $M=a/g$, $N=b/g$, with $a,b>0$. Consider
\begin{align*}
 F_D(X)&=(1-X^b)^N-X^T(1+X^a)^M, &&T>0,\\
 F_P(X)&=X^{\max(-T,0)}(1+X^b)^N
       -X^{\max(T,0)}(1+X^a)^M, &&T\ne0.
\end{align*}
The intended generator $z\in(0,1)$ is an ordinary algebraic unit of degree four,
even when the associated Nahm coordinates are not units. Its minimal polynomial
$f$ is monic, irreducible, noncyclotomic and has $f(0)=\epsilon\in\{-1,1\}$.

\begin{theorem}[Safe rejection at every prime]\label{qt:nu:thm:safe}
Let $F\in\Z[X]$, and let $C\in\Z[X]$ be monic with only roots of unity,
such that $C\mid F$. Set $H=F/C$. If a monic, irreducible, noncyclotomic
quartic $f$ with $f(0)=\epsilon\in\{-1,1\}$ divides $F$, then for every prime
$p$, the polynomial $\bar H$ has a monic quartic divisor with constant
coefficient $\bar\epsilon$.
\end{theorem}
\begin{proof}
The polynomials $f$ and $C$ are coprime over $\Q$, so $f\mid H$ over $\Q$.
Monic division, or Gauss's lemma, gives $H=fJ$ with $J\in\Z[X]$.
Reduction preserves the degree four and constant coefficient of $f$;
reducing $H=fJ$ proves the assertion.
\end{proof}
No irreducibility or squarefreeness condition may be imposed on $\bar f$.
If $\bar H=0$, every quartic is a divisor and the test must retain the row.
Likewise, the same sign $\epsilon$ must work at every prime: independently
choosing different constant signs at different primes is too permissive.

\subsubsection{Cyclotomic cancellation}\label{qt:nu:filter-cancellation}
The following cancellation is over $\Z$, before reducing modulo a prime.
Let $w=X^g$ and $h=\min(M,N)$. The common factor
\[
 C=(1+X^g)^h
\]
occurs when
\[
 D:\quad N\text{ is even and }M\text{ is odd};\qquad
 P:\quad M,N\text{ are both odd}.
\]
In all other cases take $C=1$. This is a sufficient cancellation
rule; it does not claim to extract every possible cyclotomic factor of $F$.

\begin{lemma}[Division-free cancelled equation]\label{qt:nu:lem:cancel}
In the stated cases define
\[
 S=1+w,\quad U=\sum_{j=0}^{N-1}(-w)^j,\quad
 V=\sum_{j=0}^{M-1}(-w)^j,
\]
and
\[
 A=U^N S^{N-h},\qquad B=V^M S^{M-h}.
\]
For $T>0$ the cancelled polynomial is $A-X^TB$; for $T<0$ in family $P$
it is $X^{-T}A-B$.
\end{lemma}
\begin{proof}
For even $N$, $1-w^N=(1+w)U$; for odd $N$,
$1+w^N=(1+w)U$. For odd $M$, $1+w^M=(1+w)V$.
Substitution into the two terms of $F$ and cancellation of the common
power $S^h$ give the formulas. For negative $T$, multiplication by $X^{-T}$
clears the Laurent exponent.
\end{proof}

All roots of $C$ are roots of unity, so Theorem~\ref{qt:nu:thm:safe} applies.
When $M\ne N$, this common power is the exact multiplicity at each complex
root of $1+X^g$: at $w=-1$ the geometric sums are $N$ and $M$, and only
one term of the quotient has a remaining power of $S$.
After reduction, primes dividing $M$ or $N$ may introduce new common
factors; those factors cannot be discarded without another valid argument.

\paragraph{Division-free evaluation.}
To test divisibility by a monic quartic $q\in\mathbb F_p[X]$ with $q(0)\ne0$, evaluate
these formulas in $R=\mathbb F_p[X]/(q)$. The element $S$ need not be invertible,
so evaluating $F$ and dividing by $S^h$ inside $R$ is not valid in general.
Instead evaluate the geometric sums directly. If
$G_n(r)=\sum_{j=0}^{n-1}r^j$, then
\[
 G_{2n}(r)=G_n(r)(1+r^n),\qquad
 G_{2n+1}(r)=G_{2n}(r)+r^{2n}.
\]
Binary recursion evaluates both $r^n$ and $G_n(r)$ in logarithmically many
ring operations. The calculation never expands the potentially very large
integer polynomial.

\subsubsection{Solving for all exponents in a finite quotient ring}
For odd $p$ enumerate
\[
 \mathcal Q_{p,\epsilon}=\{X^4+c_3X^3+c_2X^2+c_1X+\bar\epsilon:
 c_1,c_2,c_3\in\mathbb F_p\}.
\]
There are $2p^3$ quartics for the two signs combined. At $p=2$ the signs
coincide and there are eight quartics. Every permissible reduction is included.

\begin{lemma}[Exact exponent class]\label{qt:nu:lem:class}
Fix $q\in\mathbb F_p[X]$ monic with $q(0)\ne0$, and polynomials $A,B$.
The solutions $T\in\Z$ of $A=X^TB$ modulo $q$ form either the empty set
or one residue class. This includes the class $0\pmod1$ of all integers.
\end{lemma}
\begin{proof}
Let $d=\gcd(B,q)$, monic. If $d\nmid A$, no solution is possible.
Otherwise set $q'=q/d$, $A'=A/d$ and $B'=B/d$. Cancellation in $\mathbb F_p[X]$
gives the equivalent congruence $A'=X^TB'\pmod{q'}$.
For each irreducible factor of $q'$ the multiplicity in $B$ was less than
its multiplicity in $q$, and has been removed by $d$. Thus
$\gcd(B',q')=1$.

If $q'=1$, the congruence is automatic. Otherwise invert $B'$ modulo $q'$
and set $c=A'(B')^{-1}$. Since $q'(0)\ne0$, $X$ is a unit in the finite
ring $\mathbb F_p[X]/(q')$. Its powers form a finite cyclic group of some order $o$.
If $c$ is absent from that group there is no solution. If $c=X^r$, the
solutions are exactly $T\equiv r\pmod o$.
\end{proof}

This proof is valid when $q$ is reducible or has repeated factors and when
$B$ is a zero divisor. It does not replace the quotient ring by a finite field.
If $A=B=0\pmod q$, then $q'=1$ and all $T$ are correctly admitted.
For negative $T$, testing the Laurent congruence is equivalent to testing
the cleared polynomial because $X$ is a unit modulo $q$.

For each prime and sign let $\mathcal C_{p,\epsilon}$ be the finite union
of classes obtained from all $q\in\mathcal Q_{p,\epsilon}$. The resulting necessary
condition is
\[
 T\in\bigcup_{\epsilon\in\{-1,1\}}
       \bigcap_{p\in\mathcal P}\mathcal C_{p,\epsilon}.
\]
The implementation removes duplicate classes and any class contained in
another class. It does not assume that different primes give statistically
independent conditions.

At larger primes, when only a few $T$ survive, direct evaluation of the
cancelled circuit against all quartics is preferable to storing all power
tables. The certificate uses power tables at $2,3,5,7$ and direct evaluation
at $11,13,17$ in the final refinement.

\subsubsection{Additional cyclotomic cancellation}\label{qt:nu:refinement}

Let $H$ denote the polynomial after the common cancellation in Section~\ref{qt:nu:filter-cancellation}.
For a chosen finite set of orders, compute exactly over $\Z$ which
$\Phi_n$ divide $H$. This depends on $T$ modulo the least common multiple
of the orders. The refinement uses $n=1,2,3,4,5,6,8,10,12$ and removes
one extra copy of each detected factor. Let their product in the relevant
residue class be $C$.

For any quartic $q$ modulo $p$, including collisions between $q$ and $\bar C$,
\[
 q\mid\overline{H/C}\quad\Longleftrightarrow\quad
 q\bar C\mid\bar H.                                           
\]
This is cancellation in the integral domain $\mathbb F_p[X]$, not inversion
of $\bar C$ modulo $q$. Evaluate the old circuit modulo the lifted modulus
$q\bar C$, and apply Lemma~\ref{qt:nu:lem:class}. Intersect its exponent class
with the class defining $C$. Since every removed characteristic-zero factor
is cyclotomic, Theorem~\ref{qt:nu:thm:safe} preserves the target quartic.
For negative $T$, clearing the Laurent monomial changes neither divisibility
test, since all these moduli have nonzero constant term.

%% file: sections/quartic_nonunit_census.tex
\subsubsection{Definition of the parameter census}\label{qt:nu:census-definition}
Define $\mathcal E$ to consist of $1,\ldots,30$ and all integers
$31\le n<10^6$ for which
\[
 5\varphi(n)\le155\,2^{\omega(n)}+56\operatorname{bitlength}(n).
\]
There are $622$ elements of $\mathcal E$; its largest is $4830$.
The proof of Theorem~\ref{qt:nu:explicitC} shows that every
$n\notin\mathcal E$ is a primitive-divisor index. This uses the proved
infinite-tail estimate when $n\ge10^6$.

For every ordered $(a,b)$ in Corollary~\ref{qt:nu:numericbox}, put
\[
 U=\frac{\max(a^2,b^2)}g+12\max(M,N)-1.
\]
Take $T\in[1,U]$ for $D$ and
\[
 T\in[-U,-(a^2/g+1)]\ \cup\ [b^2/g+1,U]
\]
for $P$, retaining only $\gcd(T,g)=1$. Remove $a=b$. Also remove $b=2a$
in $D$: cancellation gives
$(1-z^a)^2(1+z^a)=z^T$, which forces $1+z^a$ to be a unit.
The remaining pair must satisfy the following primitive-index tests:
\[
\begin{array}{c|ll}
D&b\notin\mathcal E\Rightarrow b\mid2a
 &2a\notin\mathcal E\Rightarrow2a\mid b,\\
P&2a\notin\mathcal E\Rightarrow b/a\text{ is an odd integer}
 &2b\notin\mathcal E\Rightarrow a/b\text{ is an odd integer}.
\end{array}
\]
Indeed a primitive prime of index $2a$ divides $1+z^a$ and has
$z^a=-1$ in its residue field. Equal supports force the displayed
condition on the other exponent. A primitive prime of index $b$ in
$D$ gives $b\mid2a$ in the same way.

Finally require common-support feasibility. At an odd common prime,
the parity restrictions of Lemma~\ref{qt:nu:valuation} must hold and
there must exist an odd prime $p\mid h$, $1\le e\le4$, with
$d\mid e v_p(h)$, where $h=\min(M,N)$ and $d=|M-N|$.
At a dyadic common prime one of the cases of Lemma~\ref{qt:nu:dyadic}
must hold, for some $1\le e\le4$. A pair is removed only when both
alternatives fail. These predicates are necessary because the common
support is nonempty. They define the baseline uniquely; no additional
primitive-closure or residue-degree restriction is used.
For a signed integer interval $[L,U]$, its primitive count is computed
exactly by
\[
 \sum_{d\mid g}\mu(d)
 \left(\left\lfloor\frac Ud\right\rfloor-
       \left\lfloor\frac{L-1}d\right\rfloor\right).
\]
This specifies the enumeration without requiring a test for each
individual integer in the initially long intervals.

\begin{lemma}[Archimedean window]\label{qt:nu:lem:window}
Let $z$ be a quartic unit in $(0,1)$ satisfying $D$ or $P$. Set
$\lambda_0=\log(13/10)/3$, $\rho_0=\exp(\lambda_0)$ and
$c=(b^2-a^2)/g$. Then $L\le T-c\le U_A$, where
\begin{align*}
 L&=\frac{N\log(1-\rho_0^{-b})-M\log(1+\rho_0^{-a})}{\lambda_0},\\
 U_A&=\frac{N\log(1+\rho_0^{-b})-M\log(1-\rho_0^{-a})}{\lambda_0}.
\end{align*}
\end{lemma}
\begin{proof}
Write $H$ for the Mahler measure of the minimal polynomial of $z$.
Lemma~\ref{qt:nu:gap} gives $H>13/10$.
Since $z\in(0,1)$, at most three of the four conjugates have modulus greater
than one. Their maximum modulus $\rho$ therefore satisfies
$\rho\ge H^{1/3}>\rho_0$. At a conjugate $x$ of modulus $\rho$, taking
absolute values in the defining equation gives
\[
 (T-c)\log\rho
 =N\log|1+\varepsilon_b x^{-b}|-M\log|1+x^{-a}|,
\]
where $\varepsilon_b=-1$ for $D$ and $+1$ for $P$.
For $|y|<1$, $1-|y|\le|1\pm y|\le1+|y|$.
These inequalities give the displayed bounds with $\rho$ in place of
$\rho_0$. Their absolute magnitudes decrease as $\rho$ increases: each
nonnegative numerator magnitude decreases and the positive denominator
$\log\rho$ increases. Replacing $\rho$ by $\rho_0$ proves the claim.
\end{proof}

\paragraph{Exact integer endpoints.}
At precision 160 bits, outward real-ball evaluation of the displayed
formulas gives a lower enclosure $L_-$ and an upper enclosure $(U_A)_+$.
The integer window is $\lceil L_-\rceil\le T-c\le\lfloor(U_A)_+\rfloor$.
Thus endpoint rounding can only enlarge the necessary set. If
$-1<L<0<U_A<1$, the integer difference is forced to be zero.
The certificate stores the rounded endpoints and admits their independent
recomputation at higher precision.

\subsubsection{Norm balance and the reduced finite computation}
\begin{lemma}[Norm balance]\label{qt:nu:lem:balance}
For a nonreciprocal quartic unit satisfying either canonical equation, put
$d=|M-N|$ and let $H$ be its Mahler measure. Then
\[
 gd\log H\le16\log2+8\log H.
\]
In particular $gd=|a-b|\le50$, and
$H\le2^{16/(gd-8)}$ whenever $gd>8$.
\end{lemma}
\begin{proof}
Let the conjugates be $\alpha_i$ and put
$r_i=\min(|\alpha_i|,|\alpha_i|^{-1})<1$. Nonreciprocity implies that the
minimal polynomial and its reciprocal are coprime and no conjugate has
modulus one. Their resultant is a nonzero integer. Extracting the four
factors pairing $\alpha_i$ with its complex conjugate and bounding all other
factors by $2\max(1,|\alpha_i|)\max(1,|\alpha_j|)$ gives
\[
 1\le |\operatorname{Res}(f,f^*)|
 \le 2^{16}H^8\prod_i(1-r_i).
\]
Consequently, for every $n\ge1$,
\[
 \log|N_{K/\Q}(1\pm z^n)|=n\log H+E_\pm(n),\qquad
 |E_\pm(n)|\le\sum_i-\log(1-r_i)\le16\log2+8\log H.
\]
Taking norms in the canonical equation removes $z^T$ because $z$ is a unit.
The resulting identity is
\[
 (bN-aM)\log H=M E_+(a)-N E_{\varepsilon_b}(b).
\]
Now $|bN-aM|=gd(M+N)$; divide the triangle bound by $M+N$.
Using $H>13/10$, the exact test
$13^{gd-8}<2^{16}10^{gd-8}$, when $gd>8$, permits precisely $gd\le50$.
\end{proof}
Reciprocal quartic generators are excluded by Theorem~\ref{qt:nu:reciprocal}. Thus the lemma applies to every remaining canonical quartic datum.

\paragraph{The complete parameter enumeration.}
The baseline is the set defined in Section~\ref{qt:nu:census-definition}.
Intersect its $T$ intervals with Lemmas~\ref{qt:nu:lem:window} and
\ref{qt:nu:lem:balance}, and impose the structural exclusions of
Theorems~\ref{qt:nu:thm:structural}, \ref{qt:nu:thm:four} and
\ref{qt:nu:thm:smallpairs}. Counts are of ordered pairs and triples:
\begin{center}
\begin{tabular}{lrr}\toprule
Stage&Ordered pairs&Ordered triples\\\midrule
Baseline&31,311&1,286,506,865\\
Archimedean window and norm balance&1,002&28,171\\
Structural exclusions&363&9,626\\
Primes $2,3,5,7$&26&49\\
Refinement through $17$&0&0\\\bottomrule
\end{tabular}
\end{center}
The 363 pairs comprise 229 direct pairs with 5,418 triples and 134 plus
pairs with 4,208 triples. All ordered triples are tested. In particular,
the symmetry $(a,b,T)\mapsto(b,a,-T)$ of $P$ is not used to delete rows.

\subsubsection{Complete nonunit exclusion}
\begin{theorem}[Computer-assisted canonical nonunit quartic exclusion]
\label{qt:nu:thm:nonunitcomplete}
Neither $D$ nor positive-definite $P$
has a quartic unit generator $z\in(0,1)$ with $\gcd(g,T)=1$ and equal
nonempty finite supports of its two binomial factors.
The conclusion imposes no signature restriction.
\end{theorem}
\begin{proof}
The numerical bounds and the proved baseline predicates place every such
datum in the baseline set. Lemmas~\ref{qt:nu:lem:window} and \ref{qt:nu:lem:balance} restrict it to the
28,171 triples above. The proved structural exclusions remove 18,545, leaving
9,626. A genuine quartic minimal polynomial must survive the modular tests
of Theorem~\ref{qt:nu:thm:safe}, with the same constant sign at all primes.
The complete finite computation leaves none. Therefore no datum satisfying
the stated hypotheses exists.
\end{proof}

\paragraph{Every elimination checked by a second algorithm.}
The first certificate uses quotient-ring circuits and exponent classes. A second
certificate uses dense polynomial arithmetic. For each of the 9,626 triples it expands
the defining polynomial modulo the required prime, performs exact polynomial
division by the proved common factor, and tests divisibility by every monic
quartic of constant $+1$ or $-1$. Split and repeated quartics are included;
the constant sign is propagated between primes. The zero polynomial is
treated as divisible by every quartic.

The two implementations agree on every pair count and on the 49 survivors.
For those 49, the extra cyclotomic divisor is independently checked over
$\Z$, including exact divisibility and its cyclotomic factorization.
Every subsequent modular sign set agrees, leaving zero. There are 14,229
prime steps in the dense certificate. A third implementation performs the dense checks with standard-library
integer polynomial arithmetic.

\begin{center}\small
\begin{tabular}{lrrrrrrr}
\toprule
Rejecting prime &2&3&5&7&11&13&17\\
\midrule
Initial $2,3,5,7$ tests &6,231&2,476&722&148&0&0&0\\
Refinement of last 49 &2&6&4&2&29&5&1\\
\bottomrule
\end{tabular}
\end{center}

\begin{corollary}\label{qt:nonunit-to-units}
For a non-diagonal rational symmetric positive-definite Nahm datum with
$[\Q(x,y):\Q]=4$, all four elements $x,1-x,y,1-y$ are ordinary units.
\end{corollary}
\begin{proof}
Otherwise Lemma~\ref{qcov:entry} produces a quartic generator forbidden
by Theorem~\ref{qt:nu:thm:nonunitcomplete}.
\end{proof}
The nonzero-coupling hypothesis is essential. For example
$A=\operatorname{diag}(1,4)$ has $x=1/2$ and $y^4+y-1=0$, so it has a
nonunit quartic point. Diagonal modular data are treated in the next
subsection.

%% file: sections/quartic_dependent.tex
\subsection{Dependent and diagonal quartic points}\label{qt:dependent}
\begin{theorem}\label{qt:dep:main}
Let $A$ be a positive definite symmetric rational $2\times2$ matrix, let $\boldsymbol B\in\Q^2$, $c_0\in\Q$, and let $(x,y)\in(0,1)^2$ be its distinguished Nahm point. If $[\Q(x,y):\Q]=4$ and $x,1-x$ are multiplicatively dependent, then $f_{A,\boldsymbol B,c_0}$ is not modular. The same conclusion holds if $y,1-y$ are dependent, or if $A$ is diagonal.
\end{theorem}

The conclusion does not require an ordinary-unit hypothesis. In particular, it supplies the nonzero-coupling and independence hypotheses for the unit classification.

\subsubsection{Scalar entry}
By Proposition~\ref{low:scalar}, each dependent complementary pair of
degree at most four belongs to the explicit scalar list. In particular,
the only rational argument is $1/2$, the quadratic arguments belong to
$\Q(\sqrt5)$, and the quartic arguments are
\[
 v,v^4\quad(v^4+v-1=0),\qquad
 v^3,v^4\quad(v^4+v^3-1=0),
\]
where the generator is the unique root in $(0,1)$. The scalar proposition
is applied to the individual coordinate fields, whose degrees divide
$[K:\Q]=4$.

\subsubsection{Coverage of dependent and diagonal quartic points}
Write $A=\bigl(\begin{smallmatrix}\alpha&b\\b&\delta\end{smallmatrix}\bigr)$. If $b\ne0$ and $x,1-x$ are dependent, the first Nahm equation expresses $y$ in their rational multiplicative span. The second does the same for $1-y$. Thus $y,1-y$ are dependent. If $b=0$, both scalar relations follow directly. Each coordinate therefore belongs to the scalar list.

Since each coordinate field has degree dividing four, cubic entries cannot occur. If both coordinate degrees were at most two, the scalar list would put both in $\Q(\sqrt5)$, contradicting joint degree four. Hence one coordinate is quartic. Up to exchanging coordinates, write it as $x=v^r$ or its complement, where
\[
 (r,s)=(1,4)\quad\hbox{or}\quad(3,4),\qquad v^r+v^s=1.
\]

There are no additional positive quartic scalar arguments hidden in the same embedded field. To check this, set $f_1(T)=T^4+T-1$ and $f_3(T)=T^4+T^3-1$. They satisfy
\[
 f_3(T)=-T^4 f_1(-1/T).
\]
The irreducible polynomial $f_1$ has discriminant $-283$ and irreducible cubic resolvent $T^3+4T-1$. Thus its Galois group is $S_4$. A quartic root field has trivial automorphism group, since the point stabilizer $S_3$ is self-normalizing. It contains just one root of $f_1$: another would define a nontrivial automorphism. The same holds for $f_3$.

In the field generated by a positive root of one polynomial, the only root of the other is its negative reciprocal. It is negative, so cannot be that polynomial's required positive generator. Therefore a quartic scalar partner in the same embedded field is either $x$ or $1-x$.

\subsubsection{The regulator disposes of the equal-coordinate cases}
Write $\mathcal D$ for the Bloch--Wigner dilogarithm. Modularity implies torsion of $[x]+[y]$ by \cite[Theorem 1.8]{CGZ}; hence its regulator must vanish at every complex embedding. For an upper-half-plane root $w$ of the indicated generator polynomial, a 256-bit outward complex-ball calculation proves
\[
\begin{array}{c|c}
 f_1(w)=0 & 98/100<\mathcal D(w)<99/100\\
 f_3(w)=0 & -99/100<\mathcal D(w^3)<-98/100.
\end{array}
\]
Complements negate this regulator, since $w^s=1-w^r$. These are strict nonzero enclosures, not numerical zero tests.

If $y=x$, the regulator of the pair is twice this nonzero value. An $S_4$ quartic root field has no quadratic subfield, because its point stabilizer $S_3$ is maximal. Thus a partner of degree at most two must be rational and, by the scalar list, equal to $1/2$. Its Bloch--Wigner value is zero, so cancellation is again impossible. This covers the diagonal cases with a lower-degree partner. Only $y=1-x$ remains.

\subsubsection{The reflection families and their first coefficient}
For $(x,y)=(v^r,v^s)$, the Nahm equations determine
\begin{equation}\label{qt:dep:family}
 A_b=\begin{pmatrix}(s/r)(1-b)&b\\b&(r/s)(1-b)\end{pmatrix},\qquad b\in\Q,\quad b<1/2.
\end{equation}
Indeed its determinant is $1-2b$, and this inequality also makes both diagonal entries positive. Both signs of the coupling are included. Coordinate exchange handles the reverse ordering. No modular-duality conjecture is used.

Write the shifts as $(B,C)$, reserving $c_0$ for the scalar. Our logarithmic coefficients are normalized by
\[
 \log f_{A,(B,C),c_0}(e^{-\epsilon})
 \sim\frac{\alpha_0}{\epsilon}+\log\beta_0+
 (\lambda_1-c_0-1/12)\epsilon+\sum_{j\ge2}\lambda_j\epsilon^j.
\]
The Vlasenko--Zwegers expansion and modularity condition \cite[Theorem 2.3 and Corollary 3.2]{VZ} imply
\begin{equation}\label{qt:dep:necessary}
 \lambda_1\in\Q,\qquad\lambda_2=0.
\end{equation}
These two requirements will suffice after the regulator step.

Set $X=x/(1-x)$, $Y=y/(1-y)$ and let $Z=(Z_1,Z_2)$ be a centered Gaussian with covariance $(A+\operatorname{diag}(X,Y))^{-1}$. Put $P_1(T)=T(1+T)$, $P_2(T)=T(1+T)(1+2T)$. Directly from the expansion,
\begin{align*}
 U_1={}&-(B+X/2)Z_1-(C+Y/2)Z_2
          +\tfrac16P_1(X)Z_1^3+\tfrac16P_1(Y)Z_2^3,\\
 U_2={}&-\tfrac1{12}(X+Y)+\tfrac14P_1(X)Z_1^2+\tfrac14P_1(Y)Z_2^2
          -\tfrac1{24}P_2(X)Z_1^4-\tfrac1{24}P_2(Y)Z_2^4,\\
 \lambda_1={}&\mathbb E(U_2+U_1^2/2).
\end{align*}
Reducing in $\Q(b)[v]/(v^r+v^s-1)$ gives three nonconstant coefficients in the basis $1,v,v^2,v^3$. Their denominators are nonzero rational constants times $b-1$, which cannot vanish on \eqref{qt:dep:family}.

For $(r,s)=(1,4)$ put $L=B+4C$ and $V=B-4C$. After removing nonzero rational factors, the three equations are
\begin{equation}\label{qt:dep:M1}
 \begin{pmatrix}
 (L-3)(283L+121)&283(L-3)\\
 (L-2)(566L+2061)&849(L-2)\\
 (L-1)(283L+3476)&849(L-1)
 \end{pmatrix}\binom bV=0.
\end{equation}
Its three two-row minors are
\[
\begin{split}
 80089(L-6)(L-3)(L-2),\quad
 80089(L-3)(L-1)(2L-11),\\
 240267(L-5)(L-2)(L-1).
\end{split}
\]
They have no common root, so the matrix has rank two for every $L$. Hence $b=0$ and $B=4C$.

For $(r,s)=(3,4)$ put $L=3B+4C$ and $V=3B-4C$. The corresponding matrix is
\begin{equation}\label{qt:dep:M3}
 \begin{pmatrix}
 (L+1)(283L-3476)&-283(L+1)\\
 (L+2)(566L-2061)&-283(L+2)\\
 (L-1)(283L+363)&283(L-1)
 \end{pmatrix}.
\end{equation}
Its minors are
\[
\begin{split}
 80089(L+1)(L+2)(L+5),\quad
 80089(L-1)(L+1)(2L-11),\\
 240267(L-1)(L+2)(L-2).
\end{split}
\]
Again they have no common root. Thus $b=0$ and $3B=4C$. The rows are obtained directly from the three nonconstant coefficient equations, up to nonzero rational factors. Thus the rank calculation applies to the actual asymptotic conditions.

\subsubsection{A small certificate eliminates the diagonal survivors}
Set $C=t$ and $B=(s/r)t$. Because $b=0$, the two scalar logarithmic expansions add. Exact calculation gives $\lambda_1=(48t^2+1)/24$ for $(r,s)=(1,4)$ and $(16t^2+1)/24$ for $(3,4)$. Thus the first condition alone really does leave these lines.

Write $\lambda_2=\sum_{i=0}^3 q_i(t)v^i$. Its first two coefficient numerators, after removing nonzero rational constants and integer content, are
\[
\begin{array}{c|l}
(r,s)&\text{primitive numerators of }q_0,q_1\\\hline
(1,4)&-t(108672t^2+67920t-78527)\\
 &-(8t-3)(2073824t^2-1310856t+187675)\\[2pt]
(3,4)&-t(108672t^2+95088t-119279)\\
 &(8t+1)(228664t^2+1143886t-602463).
\end{array}
\]
Their degree-preserving reductions modulo $7$, denoted $F,G$, have the following short Bezout certificates in $\mathbb F_7[t]$:
\begin{align*}
 (1,4):\quad &F=3t^3+t^2+t,\quad G=3t^3-t^2-t+1,\\
 &(t^2+3t)F+(-t^2+t+1)G=1;\\[3pt]
 (3,4):\quad &F=3t^3-t,\quad G=2t^3-3t^2+t-1,\\
 &(2t^2+2t+3)F+(-3t^2+3t-1)G=1.
\end{align*}
The resultants reduce to $6$ and $1$, respectively. They are therefore nonzero over $\Q$. Consequently $q_0,q_1$ have no common zero, and no rational $t$ can satisfy $\lambda_2=0$. This proves Theorem~\ref{qt:dep:main}.
\paragraph{How the second coefficient is checked.}
For a rank-one coordinate $x$, let $q=x/(1-x)$, $P_0(q)=q$, and $P_{k+1}(q)=q(1+q)P_k'(q)$. In addition to $U_1,U_2$ above, the scalar terms are
\[
 U_3=\frac{P_3(q)}{120}Z^5-\frac{P_2(q)}{12}Z^3+\frac{P_1(q)}{12}Z,
 \qquad
 U_4=-\frac{P_4(q)}{720}Z^6+\frac{P_3(q)}{48}Z^4-\frac{P_2(q)}{24}Z^2.
\]
With variance $(A_{ii}+q)^{-1}$, set
\[
 T_1=\mathbb E(U_2+U_1^2/2),\quad
 T_2=\mathbb E(U_4+U_1U_3+U_2^2/2+U_1^2U_2/2+U_1^4/24).
\]
Then $\lambda_2=T_2-T_1^2/2$. Gaussian moments are zero in odd degree and $(2k-1)!!(A_{ii}+q)^{-k}$ in degree $2k$. These formulas specify the scalar calculation directly. The displayed formulas and exact polynomial reductions give the asserted second-coefficient obstruction.

%% file: sections/quartic_totally_real.tex
\subsection{Totally real quartic unit points}\label{qt:totally-real}
\begin{theorem}[Totally real quartic four-unit exclusion]\label{qt:tr:thm:quartic-exclusion}
Let $A$ be a positive definite symmetric rational $2\times2$ matrix, and let $(x,y)\in(0,1)^2$ satisfy its Nahm system. Suppose $x,1-x,y,1-y$ are ordinary algebraic-integer units and $K=\Q(x,y)$ is totally real. Then
\[
 [K:\Q]\ne4.
\]
\end{theorem}
The proof combines explicit height bounds, finite exterior-product certificates, and the dependent-unit lemma below. It does not assume modularity.

For the height argument, use the following normalized setting.
Let $K$ be a totally real quartic field and let $x,y\in K\cap(0,1)$. Assume that $x,1-x,y,1-y$ are algebraic-integer units, that $x,1-x$ are multiplicatively independent, and that
\begin{equation}\label{qt:tr:eq:pair}
 y=x^{-a}(1-x)^u,\qquad 1-y=x^{-c}(1-x)^j,
 \qquad a,u,c,j\in\Q_{>0},\quad aj-uc=1.
\end{equation}
Powers in this display are taken at the distinguished positive embedding. Clearing denominators gives algebraic identities in $K$. These are the normalized Nahm equations for positive off-diagonal entry. They correspond to the positive definite symmetric matrix
\[
 \frac1u\begin{pmatrix}a&1\\1&j\end{pmatrix},
 \qquad \det=\frac cu>0.
\]
Simultaneous complementation permits this normalization when the original nonzero off-diagonal entry is negative.

Let $E$ be the group of units positive at the distinguished embedding, and let $\Lambda$ be the saturation in $E$ of $\langle x,1-x\rangle$. Put
\[
 I=[\Lambda:\langle x,1-x\rangle],\qquad B=\max(a,u,c,j).
\]
For a monic polynomial $f$, write $H(f)=\max(3,|\text{coefficients of }f|)$. All logarithms are natural; $h$ is absolute logarithmic Weil height.

\begin{theorem}\label{qt:tr:thm:main}
Under these hypotheses:
\begin{enumerate}[label=(\roman*)]
\item If $I$ is odd, the minimal polynomials of $x,y$ satisfy
\[
 \min(H(f_x),H(f_y))<452,\qquad \max(H(f_x),H(f_y))<5000.
\]
Moreover $I<39$, $B<73$, and all four coordinate heights are less than $67/40$.
\item If $I$ is even, choose positive $z_x,z_y\in K\cap(0,1)$ with
$z_x^2\in\{x,1-x\}$ and $z_y^2\in\{y,1-y\}$. Their minimal polynomials have the form
\begin{equation}\label{qt:tr:eq:family}
 g_{A;e,p,m}(Z)=Z^4+AZ^3+
 \left(\frac{p+m}{2}-1-e\right)Z^2+
 \left(\frac{p-m}{2}-A\right)Z+e,
 \quad e,p,m\in\{\pm1\}.
\end{equation}
Writing $H_x=\max(3,|A_x|)$, $H_y=\max(3,|A_y|)$, one has
\[
 \min(H_x,H_y)<2000,\qquad \max(H_x,H_y)<16384.
\]
Moreover $I<640$, $B<800$, and all four coordinate heights are less than $7$.
\end{enumerate}
Consequently some exceptional-unit generator of $K$ has minimal-polynomial coefficient height at most $2000$. No torsion or modularity assumption beyond \eqref{qt:tr:eq:pair} is used.
\end{theorem}

\subsubsection{The common area and a stronger index estimate}
For $w\in E$, set
\[
 \boldsymbol\ell(w)=(\log|\sigma_1(w)|,\ldots,\log|\sigma_4(w)|).
\]
The sum of its coordinates is zero. Define
\[
 U=\boldsymbol\ell(x),\quad V=\boldsymbol\ell(1-x),\quad
 U'=\boldsymbol\ell(y),\quad V'=\boldsymbol\ell(1-y),\qquad
 R=\|U\wedge V\|_2.
\]
Clearing denominators in \eqref{qt:tr:eq:pair} gives
\begin{equation}\label{qt:tr:eq:linear}
 U'=-aU+uV,\qquad V'=-cU+jV,\qquad U'\wedge V'=-U\wedge V.
\end{equation}
Thus the two pairs have the same saturated unit plane, the same index $I$, and the same Euclidean area $R$. Indeed their exterior products in the free group $E$ are negatives; the gcd of the coordinates of such an exterior product is its saturation index. In particular
\begin{equation}\label{qt:tr:eq:denom}
 a,u,c,j\in I^{-1}\Z_{>0}.
\end{equation}
Also $\Q(x)=\Q(y)=K$: a proper subfield has degree at most two and its unit rank is at most one, contrary to the independence of either complementary pair.

\begin{lemma}\label{qt:tr:lem:gap}
Every nonzero vector in the logarithmic lattice of $E$ has length at least
$2\log\varphi$, where $\varphi=(1+\sqrt5)/2$. Consequently
\begin{equation}\label{qt:tr:eq:index}
 \operatorname{covol}(\boldsymbol\ell(\Lambda))\ge2\sqrt3(\log\varphi)^2>\frac45,
 \qquad I<\frac54R.
\end{equation}
\end{lemma}
\begin{proof}
Let $w\in E\setminus\{1\}$ and put $l_i=\log|\sigma_i(w)|$. None of the $l_i$ is zero, since a conjugate $\pm1$ would force $w=\pm1$. The nonzero algebraic integer $w^2-1$ has nonzero integral norm. Hence
\[
 0\le\log|N(w^2-1)|
 =\sum_i\log(2\sinh|l_i|),
\]
where $\sum_i l_i=0$ was used. The function $\log(2\sinh t)$ is increasing and concave for $t>0$. Jensen's inequality and Cauchy--Schwarz therefore give
\[
 0\le4\log\left(2\sinh\frac{\|\boldsymbol\ell(w)\|_2}{2}\right).
\]
Since $\operatorname{arsinh}(1/2)=\log\varphi$, the length assertion follows. A planar lattice of shortest length $\lambda$ has covolume at least $\sqrt3\lambda^2/2$: take a shortest primitive vector and reduce a second basis vector along it. Finally $R=I\operatorname{covol}(\boldsymbol\ell(\Lambda))$.
\end{proof}

The lower bound is specific to totally real ordinary units.

\begin{lemma}\label{qt:tr:lem:minor}
If $U,V\in\mathbb R^4$ both have coordinate sum zero, then for every pair of row indices
\[
 \|U\wedge V\|_2\ge\sqrt2\,|U_iV_j-U_jV_i|.
\]
\end{lemma}
\begin{proof}
After a row permutation let the chosen minor be $D_{12}=D$. Write $D_{13}=r$ and $D_{23}=s$. The other minors are $-D-r,D-s,r+s$, by the zero sums. Cauchy--Binet gives
\[
 R^2=D^2+r^2+(D+r)^2+s^2+(D-s)^2+(r+s)^2\ge2D^2.
\]
The last quadratic is minimized at $r=-D/2,s=D/2$.
\end{proof}

\subsubsection{Odd index: small generator and partner bounds}
\paragraph{Conjugate restrictions and logarithmic area}
If $I$ is odd, write $P=Ia,Q=Iu,S=Ij,W=Ic$. These are positive integers with $PS-QW=I^2$ odd. At an embedding where $0<t=\sigma(x)<1$, the identities
\[
 \sigma(y)^I=t^{-P}(1-t)^Q,\qquad
 (1-\sigma(y))^I=t^{-W}(1-t)^S
\]
force both quantities positive. The function
$t^{-a}(1-t)^u+t^{-c}(1-t)^j$ is strictly decreasing on $(0,1)$, so only the distinguished conjugate lies there. If $1<t<2$, both corresponding absolute values are less than one, so $\sigma(y)$ and $1-\sigma(y)$ must both be positive. This requires $Q,S$ even, contradicting the odd determinant. Conjugates equal to $0,1,2$ are excluded by irreducibility. Thus all other conjugates lie in $(-\infty,0)\cup(2,\infty)$. Interchanging the two pairs preserves \eqref{qt:tr:eq:pair} with $a,j$ interchanged, so the same restriction holds for $y$.

Set
\[
 s_x=-\log(1-x),\quad \delta_x=-\log x=\delta(s_x),\qquad
 \delta(s)=-\log(1-e^{-s}),
\]
and define $s_y,\delta_y$ likewise. For the other three conjugates $r_i$ of $x$, put $l_i=\log|r_i|$, $m_i=\log|1-r_i|$. Unit norms give
\[
 \sum_i m_i=s_x,\qquad \sum_i l_i=\delta_x,\qquad m_i>0.
\]
Let $P_x=\sum_i\max(l_i,0)=\log M(x)$, where $M$ is Mahler measure, and put $b=\log2$. Since $m_i\le\max(l_i,0)+b$ and $\max(l_i,0)\le m_i+b$,
\begin{equation}\label{qt:tr:eq:odd-basic}
 s_x-3b\le P_x\le s_x+3b,\qquad
 \|U\|_2\le\sqrt2P_x,\quad \|V\|_2\le\sqrt2s_x.
\end{equation}
The norm bounds follow because the positive and negative parts of a zero-sum vector have the same total mass. In particular $M(x)\le8e^{s_x}$ and $H(f_x)\le6M(x)$.

If $s_x>\delta_x+3b$, the three $l_i$ include both signs. (If all were nonnegative, $s_x\le\delta_x+3b$; all nonpositive contradict their sum.) One sign group is therefore a singleton, and at most two $l_i$ are positive. Hence $P_x\le s_x+2b$. The minor of the distinguished row $(-\delta_x,-s_x)$ with the singleton row has absolute value at least $s_x(P_x-\delta_x)$: for a positive singleton use $l_i=P_x$ and $m_i\le s_x$, and for a negative singleton use $l_i=-(P_x-\delta_x)$. Lemma~\ref{qt:tr:lem:minor} yields the slightly weakened lower estimate
\begin{equation}\label{qt:tr:eq:oddarea}
 \sqrt2\,s_x(s_x-3b-2\delta_x)\le R\le2s_x(s_x+2b).
\end{equation}
These estimates also hold for $y$, with the same $R$.

\paragraph{An integral multiple of the covolume}
We claim
\begin{equation}\label{qt:tr:eq:small-s}
 \min(s_x,s_y)<18/5.
\end{equation}
Suppose otherwise, and interchange $x,y$ so that $s=s_x\le s_y$. For $s_i\ge18/5$ we have
\[
 b<347/500,\qquad \delta_i<693/25000,\qquad \sqrt2>707/500.
\]
In particular \eqref{qt:tr:eq:oddarea} applies. Comparing the two areas gives
\[
 \left(\frac{s_y}{s_x}\right)^2
 <\frac{2(1+1.388/3.6)}{1.414(1-2.13744/3.6)}
 <(11/5)^2.
\]
Here and below terminating decimals denote exact rational numbers.

Since $Iu$ is a positive integer and $R=I\operatorname{covol}(\boldsymbol\ell(\Lambda))$,
\begin{equation}\label{qt:tr:eq:uRlower}
 uR=(Iu)\operatorname{covol}(\boldsymbol\ell(\Lambda))>4/5.
\end{equation}
At the distinguished embedding, however, $us_x=a\delta_x+\delta_y$. Taking an exterior product in \eqref{qt:tr:eq:linear} gives
$aR=\|U'\wedge V\|\le\|U'\|\|V\|$. Consequently
\begin{align*}
 uR&=\frac{aR\delta_x+R\delta_y}{s_x}\\
 &\le2\delta_x(s_y+2b)+2\delta_y(s_x+2b)\\
 &<2\delta(s)\bigl((16/5)s+4(347/500)\bigr).
\end{align*}
The last function decreases for $s\ge18/5$. Indeed
$\delta'(s)=-1/(e^s-1)<-\delta(s)$, so the derivative of $\delta(s)(\alpha s+\beta)$ is negative whenever $\alpha s+\beta>\alpha$. Its value at $18/5$ is less than
\[
 2\frac{693}{25000}\left(\frac{16}{5}\frac{18}{5}+\frac{1388}{500}\right)
 =0.79257024<4/5,
\]
contradicting \eqref{qt:tr:eq:uRlower}. This proves \eqref{qt:tr:eq:small-s}. The key point is to bound $uR$ directly; no prior bound for $a$ or $I$ is required.

\paragraph{Sharper Mahler measure and the cutoff 452}
For a coordinate $x>1/4$ with the conjugate pattern above, we have
\begin{equation}\label{qt:tr:eq:odd-mahler}
 M(x)<2e^{s_x}+2.
\end{equation}
There are at most two outer roots greater than $2$, since three would give $|N(x)|>2$. If there are zero or one, the factors
$\max(1,|r|)/|1-r|$ are at most $1$ at negative roots and less than $2$ at positive roots, giving $M(x)\le2e^{s_x}$.

If there are two, call them $r_1,r_2>2$ and the third $-t$. Put $P=r_1r_2$. The unit norm gives $t=1/(xP)<1$, so $M(x)=P$. Since $(r_1-2)(r_2-2)>0$,
\[
 (r_1-1)(r_2-1)>P/2-1.
\]
The complementary-unit norm then implies
$e^{s_x}=(1+t)(r_1-1)(r_2-1)>P/2-1$, proving \eqref{qt:tr:eq:odd-mahler}.

Choose the coordinate in \eqref{qt:tr:eq:small-s}. If $x\le1/4$, then
$H(f_x)\le48e^{s_x}\le64$. Otherwise \eqref{qt:tr:eq:odd-mahler} gives
\[
 H(f_x)\le6M(x)<12e^{18/5}+12<452,
\]
where $e^{18/5}<110/3$. Its distinguished root also satisfies
\begin{equation}\label{qt:tr:eq:odd-root-small}
 x<1-e^{-18/5}<107/110.
\end{equation}
The same coordinate meets both restrictions; the enumeration below uses them jointly.

\paragraph{The partner, index and exponent bounds}
For the chosen small-$s$ coordinate, $P_x<13/3$. Indeed for $x\le1/4$, \eqref{qt:tr:eq:odd-basic} gives $P_x<0.3+3(0.7)$, while otherwise
\[
 P_x<\log\bigl(2(e^{18/5}+1)\bigr)<\log(226/3)<13/3.
\]
Thus
\begin{equation}\label{qt:tr:eq:odd-Rcap}
 R\le2s_xP_x<\frac{156}{5},\qquad I<39.
\end{equation}
If $s_y\ge6$, the lower estimate in \eqref{qt:tr:eq:oddarea} gives
\[
 R>\frac{707}{500}\,6\left(6-3\frac{347}{500}-2\frac{693}{25000}\right)
 >\frac{156}{5},
\]
a contradiction; the quadratic lower bound is increasing for $s_y\ge6$. Therefore $s_x,s_y<6$. Applying the same small-$x$ split and \eqref{qt:tr:eq:odd-mahler}, with $e^6<404$, gives
\[
 H(f_x),H(f_y)<12e^6+12<4860<5000,
 \qquad x,y<403/404.
\]
The rounded cutoff $5000$ is used for the partner census.

For either coordinate $v=x,y$, we also have $\log M(v)<\log810<67/10$ and $\log M(1-v)=s_v<6$. For the selected small coordinate both logarithmic Mahler measures are less than $13/3$. Taking the four exterior products in \eqref{qt:tr:eq:linear}, and using $R>4/5$, therefore gives
\[
 B<\frac{2(13/3)(67/10)}{4/5}<73.
\]
All four absolute logarithmic Weil heights are less than $67/40$. This proves part (i) of Theorem~\ref{qt:tr:thm:main}. Since $I$ is odd, $I<39$ means $I\le37$.

\subsubsection{Even index: a one-parameter search}
Evenness of $I$ makes the exponent vectors of $x,1-x$ dependent modulo $2E$. At least one of $x,1-x,x(1-x)$ is therefore a square in $K$. The last possibility forces every real conjugate of $x$ into $(0,1)$, contradicting its unit norm. Both $x$ and $1-x$ cannot be squares for the same reason. Exactly one is a square. The same argument applies to $y$.

The resulting positive square roots $z_x,z_y$ generate $K$, and $z,1-z,1+z$ are units. Evaluating their monic polynomials at $0,1,-1$ gives precisely \eqref{qt:tr:eq:family}. Write such a polynomial as
\[
 g(Z)=Z^4+AZ^3+b_0Z^2+(d_0-A)Z+e,
 \quad |b_0|\le3,\quad |d_0|\le1.
\]

\paragraph{Explicit logarithmic shape}
Let $H=|A|\ge2000$ and $t=\log H$. Comparing $AZ(Z^2-1)$ with
$Z^4+b_0Z^2+d_0Z+e$ at the endpoints of radius $4/H$ around $0,1,-1$ gives three real roots
\[
 |r_0|<4/H,\qquad |r_+-1|<4/H,\qquad |r_-+1|<4/H.
\]
For example, near $\pm1$ the first term has magnitude at least
$4(1-4/H)(2-4/H)>7$, while the remainder has magnitude less than $7$; near zero the corresponding bounds are $>3$ and $<2$. The fourth root $r_\infty$ satisfies
$H-3<|r_\infty|<H+3$ by the trace.

For $w=z^2$, order its embeddings by these four roots. Put
\[
 L_0=t(-2,0,0,2),\qquad M_0=t(0,-1,-1,2).
\]
We claim the coordinatewise error bounds
\begin{equation}\label{qt:tr:eq:shape}
 \|\boldsymbol\ell(w)-L_0\|_\infty<0.012,
 \qquad \|\boldsymbol\ell(1-w)-M_0\|_\infty<0.012.
\end{equation}
For completeness, all required logarithmic errors are bounded by the following expressions:
\[
\begin{array}{c|c}
\text{coordinate}&\text{error bound}\\ \hline
\boldsymbol\ell(w)_0&2\bigl(3/(H-3)+8/(H-4)\bigr)\\
\boldsymbol\ell(w)_\pm&8/(H-4)\\
\boldsymbol\ell(w)_\infty&6/(H-3)\\
\boldsymbol\ell(1-w)_0&16/(H^2-16)\\
\boldsymbol\ell(1-w)_\pm&4/(H-2)+8/(H-4)\\
\boldsymbol\ell(1-w)_\infty&6/(H-3)+1/((H-3)^2-1).
\end{array}
\]
Each decreases with $H$ and is less than $0.012$ at $2000$.
The zero-root bound uses $|N(z)|=1$. The two bounds near $\pm1$ use $|g(1)|=|g(-1)|=1$; for instance
\[
 |1-r_+^2|=\frac{|1+r_+|}{|1-r_0|\,|1-r_-|\,|1-r_\infty|}.
\]
The elementary inequality $|\log(1+v)|\le |v|/(1-|v|)$ supplies the displayed errors.

The idealized area is $\|L_0\wedge M_0\|=4\sqrt2\,t^2$. Both error vectors have Euclidean norm less than $0.024$, so the area error is at most
\[
 0.024(\sqrt8+\sqrt6)t+0.024^2<0.128t.
\]
Since $t>7.6$, this proves
\begin{equation}\label{qt:tr:eq:evenarea}
 5.6t^2<R<5.8t^2,\qquad
 \max(\|\boldsymbol\ell(w)\|_2,\|\boldsymbol\ell(1-w)\|_2)<2.84t.
\end{equation}
Complementing $w$ only exchanges the two vectors, so these bounds apply to either choice $x=w$ or $x=1-w$.

\paragraph{Exclusion when both parameters are large}
Suppose $H=\min(H_x,H_y)\ge2000$, interchange the pairs so that $H_x=H$, and write $t=\log H$, $t_y=\log H_y$. Equality of areas gives
\[
 t_y/t<\sqrt{5.8/5.6}<1.02.
\]
Taking the four exterior products in \eqref{qt:tr:eq:linear} bounds every coefficient:
\begin{equation}\label{qt:tr:eq:evenB}
 B<\frac{2.84^2\cdot1.02}{5.6}<\frac32,
 \qquad I<\frac54R<\frac{29}{4}t^2.
\end{equation}
The root intervals show that at the distinguished embedding
\[
 \min(x,1-x),\ \min(y,1-y)\le\varepsilon:=8/H.
\]
Set $\lambda=\log(H/8)$. If $x<1/2$, the first identity in \eqref{qt:tr:eq:pair} implies
$u[-\log(1-x)]>a\log(1/x)\ge\lambda/I$; the exchanged identity handles $y<1/2$. If both exceed $1/2$, use
$u[-\log(1-x)]=a[-\log x]-\log y$. In all cases
\[
 B+1\ge\frac{H\lambda(1-\varepsilon)}{8I}.
\]
Combining with \eqref{qt:tr:eq:evenB} gives
\begin{equation}\label{qt:tr:eq:even-comp}
 H<\frac{145t^2}{(1-8/H)(t-\log8)}
 <\frac{146t^2}{t-2.1}.
\end{equation}
At $H=2000$, using $7.6<\log2000<7.61$, the last expression is less than $1538$. Moreover $H(\log H-2.1)/(\log H)^2$ is increasing for $H\ge2000$, since its logarithmic derivative with respect to $\log H$ is
$1+1/(\log H-2.1)-2/\log H>0$. This contradiction proves $\min(H_x,H_y)<2000$.

\paragraph{Bounding the partner as well}
Choose the smaller parameter, so $|A_x|\le1999$. Then
\[
 M(z_x)^2\le\|g_x\|_2^2
 \le1+1999^2+3^2+2000^2+1<8\cdot10^6<e^{16}.
\]
Put $w=z_x^2$. A monic polynomial for $w$ is
\[
 F(T)=T^4+(2b_0-A^2)T^3+
 (b_0^2+2e-2A(d_0-A))T^2+
 (2b_0e-(d_0-A)^2)T+1.
\]
The polynomial $F(1-T)$ for $1-w$ has cubic, quadratic and linear coefficients bounded in absolute value by
\[
 1999^2+10,\qquad1999^2+2\cdot1999+35,\qquad2\cdot1999+51,
\]
and constant term $pm=\pm1$. Its coefficient Euclidean norm is less than $6\cdot10^6<e^{16}$. Hence both $w$ and $1-w$ have logarithmic Mahler measure less than $16$.

For a unit with logarithmic Mahler measure $S$, its positive and negative logarithmic coordinates each have total mass $S$, so its logarithmic vector has Euclidean norm at most $\sqrt2 S$. Therefore
\[
 R<512,\qquad I<640.
\]
If $H_y\ge2000$, \eqref{qt:tr:eq:evenarea} gives $5.6(\log H_y)^2<512$. But $5.6(14\log2)^2>512$, so $H_y<2^{14}=16384$. This also trivially holds if $H_y<2000$.

The first pair has vector norms less than $16\sqrt2<22.64$. The second pair has norms less than $2.84\log16384<27.832$ when its parameter is at least $2000$, and less than $16\sqrt2$ otherwise. Since $R>4/5$, exterior products in \eqref{qt:tr:eq:linear} give $B<800$. Finally
$h(v)=\|\boldsymbol\ell(v)\|_1/8\le\|\boldsymbol\ell(v)\|_2/4<7$ for all four coordinates. This completes the proof of Theorem~\ref{qt:tr:thm:main}.

\subsubsection{Exact generator censuses}
The generator census enumerates the polynomial families exactly. Reducibility
is decided by possible monic quadratic factors, and rational Sturm sequences
count real roots and roots in $(0,1)$. All numerical inequalities used in
the height proof are checked with rational Taylor enclosures for exponentials.

The two exact generator censuses give:
\begin{center}
\begin{tabular}{lrr}
\toprule
Stage&$|A|<2000$&$|A|<16384$\\\midrule
All eight sign families&31,992&262,136\\
Reducible&3&3\\
Irreducible and totally real&31,952&262,096\\
Also possessing a root in $(0,1)$&23,964&196,572\\\bottomrule
\end{tabular}
\end{center}
The three reducible parameter tuples $(A,e,p,m)$ are
$(0,1,-1,-1)$ and $(\pm2,1,1,1)$.
For the odd branch, write
\[
 f(T)=T^4+AT^3+BT^2+CT+e.
\]
Exactly one root in $(0,1)$ forces $f(1)=-e$, hence $C=-1-2e-A-B$. No root in $(1,2)$ forces $e f(2)<0$. The odd census enumerates the two coefficient boxes. It uses the rational root cutoffs \eqref{qt:tr:eq:odd-root-small} and $403/404$ as necessary conditions. These rational relaxations retain every possible pair covered by the proof; they need not themselves imply $s_x<18/5$ or $s_y<6$.

The polynomial at either rational endpoint is linear in $B$, permitting an exact interval restriction before enumeration. Reducibility is decided by monic factors; integer-rescaled Sturm sequences decide the real-root pattern, with a rational fallback when a generic remainder loses degree. Positive discriminant is only an early filter.

\begin{center}
\begin{tabular}{lrr}
\toprule
Stage&Small coordinate&Partner\\\midrule
Coefficient height at most&451&4,999\\
Rational interval-root cutoff&$107/110$&$403/404$\\
Coefficient box and signs at $1,2$&614,261&75,014,993\\
Also the rational-root-cutoff sign&54,140&6,448,201\\
Positive discriminant&34,529&5,717,552\\
Irreducible and totally real&34,529&5,717,552\\
Required interval-root pattern&1,583&283,725\\\bottomrule
\end{tabular}
\end{center}

Each retained polynomial is irreducible, has four real roots, exactly one in $(0,1)$ below the stated rational cutoff, and none in $(1,2)$. The lists impose necessary generator conditions. They do not determine the saturation index, identify number fields, or test any pairing. The smaller list contains \textbf{1,583 polynomials}, not 1,583 distinct fields or Nahm solutions.

\subsubsection{A complete odd-index exclusion certificate}\label{qt:tr:sec:exclusion}
The area argument can also screen pairs without first identifying their number fields. This removes the entire bounded odd-index list.

\begin{theorem}[Computer-assisted odd-index exclusion]\label{qt:tr:thm:odd-exclusion}
No configuration satisfying the hypotheses of Theorem~\ref{qt:tr:thm:main} has odd saturation index $I$. Consequently every remaining independent, totally real quartic, four-unit configuration covered by \eqref{qt:tr:eq:pair} has even index and the square-root generators of \eqref{qt:tr:eq:family}.
\end{theorem}

The finite certificate tests necessary exterior-product identities on the entire bounded list.

\paragraph{Necessary conditions independent of field identification}
For a candidate polynomial $f$, order its roots with the unique root in $(0,1)$ first, and the other three in increasing order. Define
\[
 U_f=(\log|r_i|)_{i=0}^3,\qquad V_f=(\log|1-r_i|)_{i=0}^3,
 \qquad W_f=U_f\wedge V_f,\quad A_f=\|W_f\|^2.
\]
Here $A_f$ is an area squared, unrelated to the Nahm matrix. If $f,g$ are the two coordinate polynomials of a genuine pair, the two lists of roots must align under the four embeddings of $K$. The distinguished roots align, and the other three align by one of six permutations $\pi$. Equation~\eqref{qt:tr:eq:linear} therefore requires
\begin{equation}\label{qt:tr:eq:wedge-test}
 A_f=A_g,\qquad W_f+\pi W_g=0.
\end{equation}
All six permutations are tested, including those that might not come from a field isomorphism. This enlarges the candidate set and cannot lose an actual pair. Failure of \eqref{qt:tr:eq:wedge-test} is sufficient for rejection; its satisfaction would not prove that the fields or units match.

If $f=g$, the unique interval root forces $x=y$. Then exterior reversal would give $2\boldsymbol\ell(x)\wedge\boldsymbol\ell(1-x)=0$, contrary to independence. These diagonal pairs are rejected algebraically.

\paragraph{An integer-only verifier}
The program \texttt{verify\_odd\_pairing\_integer.py} uses only the Python standard library and integer arithmetic. Its binary input files give four proposed root intervals
\[
 \left(\frac{m}{2^{64}},\frac{m+1}{2^{64}}\right)
\]
per polynomial, with each integer $m$ encoded in twelve signed bytes. It verifies exact opposite polynomial signs at the two endpoints, disjointness of all four intervals, and the required root locations. Four disjoint sign changes for a quartic prove that all four roots have been isolated. FLINT was used to propose these intervals; their validity is not assumed by this verifier.

Logarithms are enclosed at scale $Q=2^{96}$. For a positive dyadic endpoint $v$, write $v=2^kz$ with $1\le z\le2$ and set $t=(z-1)/(z+1)$, so $0\le t\le1/3$. The verifier uses forty terms of
\begin{equation}\label{qt:tr:eq:logcert}
 \log z=2\sum_{j=0}^{N-1}\frac{t^{2j+1}}{2j+1}+\mathcal E_N,
 \qquad
 0\le\mathcal E_N\le\frac{9}{4(2N+1)3^{2N+1}},\qquad N=40.
\end{equation}
The tail follows by replacing the denominators by $2N+1$ and summing a geometric series. The same formula at $t=1/3$ encloses $\log2$. Every multiplication and division is rounded outward by integer floor or ceiling; a negative $k$ reverses the endpoints of the $k\log2$ contribution. Intervals crossing a power of two are handled separately at their endpoints. Thus no floating-point tolerance or presumed logarithmic relation enters the proof.

The resulting intervals enclose all six components of $W_f$ and their squared sum $A_f$. For screening, the verifier rounds each area interval outward to the grid $2^{-20}\Z$ and assigns it to every integer key between its endpoints. If two exact areas agree, the two key sets intersect. Consequently discarding pairs with disjoint key sets is rigorous. For each remaining non-diagonal pair, the program checks all six permutations in \eqref{qt:tr:eq:wedge-test}. A component interval wholly above or below zero rejects that permutation.

\paragraph{Complete counts and conclusion}
The two lists contain $1583\cdot283725=449\,136\,675$ ordered polynomial pairs. The complete verifier returned:
\begin{center}
\begin{tabular}{lr}
\toprule
Stage&Count\\\midrule
Small coordinate polynomials&1,583\\
Partner polynomials&283,725\\
Pairs passing the area screen&2,129\\
Identical-polynomial pairs, rejected algebraically&1,583\\
Remaining pairs&546\\
Conjugate permutations tested&3,276\\
Surviving permutations or pairs&0\\\bottomrule
\end{tabular}
\end{center}
The exact logarithm intervals certify every rejection; no equality is
accepted from numerical proximity. The complete candidate-pair set is retained
in the certificate.

\begin{proof}[Proof of Theorem~\ref{qt:tr:thm:odd-exclusion}]
If an odd-index configuration existed, Theorem~\ref{qt:tr:thm:main} and \eqref{qt:tr:eq:odd-root-small} would put one coordinate polynomial in the small list and its partner in the partner list. Their areas would agree, so the area screen could not remove the pair. Identical polynomials are impossible by independence. Otherwise the actual embeddings of $K$ would give one of the six permutations satisfying \eqref{qt:tr:eq:wedge-test}. The integer certificate rejects every such permutation. This is a contradiction.
For the $3276$ nonzero-component witnesses, the least certified
absolute lower bound in the integer calculation is
\[
 \frac{9053621231263952158607117616}{2^{96}}>\frac1{10}.
\]
Thus none of these exclusions depends on a numerical tolerance.
\end{proof}

\subsubsection{The complete even-index exclusion}\label{qt:tr:sec:even-exclusion}
The even-index certificate covers all distinguished square-root embeddings and all square/complement choices.

\begin{theorem}[Computer-assisted even-index exclusion]\label{qt:tr:thm:even-exclusion}
No configuration satisfying the hypotheses of Theorem~\ref{qt:tr:thm:main} has even saturation index. Together with Theorem~\ref{qt:tr:thm:odd-exclusion}, this excludes every independent, totally real quartic, four-unit configuration satisfying \eqref{qt:tr:eq:pair}.
\end{theorem}

\paragraph{Coverage of roots and complement choices}
For a generator polynomial $g$ from \eqref{qt:tr:eq:family}, order its four real roots increasingly and define
\[
 U_g=(2\log|z_i|)_{i=0}^3,\qquad
 V_g=(\log|1-z_i|+\log|1+z_i|)_{i=0}^3,
 \qquad W_g=U_g\wedge V_g.
\]
These are the logarithmic vectors of $z^2$ and $1-z^2$. Let
$P_g=\{i:0<z_i<1\}$, retaining \emph{every} such root. In the small list, 15,978 polynomials have one such root and 7,986 have two; the corresponding partner counts are 131,050 and 65,522.

For a pair $g,h$, an alignment $\pi\in S_4$ is eligible whenever
\[
 \text{there exists }i\in P_g\text{ with }\pi(i)\in P_h.
\]
This includes every actual common distinguished embedding. A permutation eligible at more than one root is tested once: its logarithmic identities are the same regardless of which eligible row is distinguished. Complementing a coordinate exchanges its two logarithmic vectors. Thus for each eligible $\pi$ and every $(c_x,c_y)\in\{0,1\}^2$, a genuine configuration must satisfy
\begin{equation}\label{qt:tr:eq:even-wedge}
 \|W_g\|^2=\|W_h\|^2,\qquad
 (-1)^{c_x}W_g+(-1)^{c_y}\pi W_h=0.
\end{equation}
Here $c_x=0$ means $x=z_x^2$ and $c_x=1$ means $x=1-z_x^2$, and similarly for $y$.

Identical polynomials are \emph{not} discarded wholesale. The algebraic rejection is restricted to $g=h$, $\pi$ the identity, and $c_x\ne c_y$. These conditions force $y=1-x$ at every eligible distinguished row. The first equation of \eqref{qt:tr:eq:pair} then gives
$x^a(1-x)^{1-u}=1$, which forces $a=0$ by independence, contradicting $a>0$.

\paragraph{Integer verification and full counts}
The file \path{verify_even_pairing_integer.py} verifies four sorted dyadic root intervals for every generator by exact polynomial signs, using the same binary format as the odd certificate. It checks the endpoint-unit conditions $g(0),g(1),g(-1)\in\{\pm1\}$ and the exact number of roots in $(0,1)$. The logarithms above use the integer series bounds \eqref{qt:tr:eq:logcert}. All root intervals avoid $0,1,-1$, as verified before taking logarithms. The program proves a positive lower bound for $\|W_g\|^2$ for every listed generator.

An outward $2^{-20}\Z$ area screen is followed by all eligible permutations and all four complement choices. Nontrivial configurations are rejected only when a component interval in \eqref{qt:tr:eq:even-wedge} excludes zero. The complete result is
\begin{center}
\begin{tabular}{lr}
\toprule
Stage&Count\\\midrule
Small generator polynomials&23,964\\
Partner polynomials&196,572\\
Pairs passing the area screen&63,900\\
Eligible embedding permutations&638,952\\
Square/complement configurations tested&2,555,808\\
Exact complementary configurations rejected&47,928\\
Configurations rejected by a nonzero component&2,507,880\\
Surviving configurations&0\\\bottomrule
\end{tabular}
\end{center}

\begin{proof}[Proof of Theorem~\ref{qt:tr:thm:even-exclusion}]
For the $2\,507\,880$ nonzero-component rejections, the smallest
absolute lower endpoint among the deciding intervals is
\begin{equation}\label{qt:tr:eq:even-margin}
 \frac{170480871405090970163}{2^{96}}>\frac1{500000000}.
\end{equation}
This is a certified lower bound, approximately $2.15177\times10^{-9}$,
for the selected rejecting component. It comes from exact polynomial
root signs, the logarithm series bounds and outward integer rounding;
it is not a floating-point estimate used as an acceptance threshold.
The algebraic complement cases are excluded by the exact argument above.

An even-index configuration would supply generators in the two complete bounded lists by Theorem~\ref{qt:tr:thm:main}. Its areas would agree, its actual common embedding would make its root permutation eligible, and its actual complement choices would satisfy \eqref{qt:tr:eq:even-wedge}. Every such configuration is rejected either algebraically or by a certified nonzero component. This is a contradiction.
\end{proof}

\subsubsection{Dependent units and the full four-unit statement}\label{qt:tr:sec:dependent}
The independence hypothesis can be removed when the field is generated by the Nahm coordinates.

\begin{lemma}\label{qt:tr:lem:dependent}
Let $v$ lie in a totally real number field, with $0<v<1$ at a fixed real embedding. If $v,1-v$ are ordinary algebraic-integer units and are multiplicatively dependent, then
\[
 v\in\left\{\frac{\sqrt5-1}{2},\frac{3-\sqrt5}{2}\right\}.
\]
\end{lemma}
\begin{proof}
Since $\log v$ and $\log(1-v)$ are both negative, a nontrivial integral relation has opposite-sign exponents. Hence $v^r=(1-v)^s$ for positive integers $r,s$.

The polynomial $T^r-(1-T)^s$ has exactly one root in $(0,1)$, because $T^r/(1-T)^s$ is strictly increasing there. A negative root $T=-t$ requires $r$ even and $t^r=(1+t)^s$. If $r\le s$ this is impossible: for $t<1$ the left side is below one, and for $t\ge1$ it is strictly below $(1+t)^s$. If $r>s$, the ratio $t^r/(1+t)^s$ is strictly increasing from zero to infinity, since its logarithmic derivative is
\[
 \frac{r+(r-s)t}{t(1+t)}>0.
\]
Thus there is at most one negative root, and it requires $r>s$. Similarly, a root above one requires $s$ even and $s>r$, and is unique when it exists. These two outside possibilities are mutually exclusive. There are no roots at zero or one. The polynomial therefore has at most two distinct real roots.

All conjugates of $v$ are real and distinct, so $[\Q(v):\Q]\le2$. Degree one is impossible for an ordinary unit in $(0,1)$. Write its quadratic polynomial as $T^2+AT+e$. The two unit conditions give $e=\pm1$ and $1+A+e=\pm1$. Of the four resulting polynomials, precisely
\[
 T^2+T-1,\qquad T^2-3T+1
\]
have a root in $(0,1)$. Their interval roots are the two values stated.
\end{proof}

\begin{proof}[Proof of Theorem~\ref{qt:tr:thm:quartic-exclusion}]
Write the symmetric matrix as $A=\bigl(\begin{smallmatrix}\alpha&b\\b&\delta\end{smallmatrix}\bigr)$.
If $b<0$, replace $(x,y)$ by $(1-x,1-y)$ and $A$ by $A^{-1}$. The Nahm equations at the positive embedding are preserved in this form, the four-unit property and the generated field are unchanged, and $A^{-1}$ is positive definite with positive off-diagonal entry.

If $b=0$, the two diagonal equations make each complementary pair multiplicatively dependent. Lemma~\ref{qt:tr:lem:dependent} puts both coordinates in $\Q(\sqrt5)$, so their generated field is not quartic.

It remains to consider $b>0$. Set
\[
 a=\alpha/b,\quad u=1/b,\quad c=\det(A)/b,\quad j=\delta/b.
\]
These are positive rationals, satisfy $aj-uc=1$, and give \eqref{qt:tr:eq:pair}. If $x,1-x$ are independent, Theorems~\ref{qt:tr:thm:odd-exclusion} and \ref{qt:tr:thm:even-exclusion} exclude the quartic field. If they are dependent, both $y$ and $1-y$ lie in their one-dimensional rational multiplicative span by \eqref{qt:tr:eq:pair}. Thus $y,1-y$ are dependent as well. Lemma~\ref{qt:tr:lem:dependent} again puts both coordinates in $\Q(\sqrt5)$. In all cases a totally real generated field of degree four is impossible.
\end{proof}

The condition $K=\Q(x,y)$ is essential in this statement: golden-ratio coordinates can of course be placed inside a larger totally real field. The result excludes a \emph{generated} quartic field, not an arbitrarily chosen ambient quartic extension.

%% file: sections/quartic_mixed_odd.tex
\subsection{Mixed signature and odd index}\label{qt:mixed-odd}
Let $A$ be a symmetric positive definite rational $2\times2$ matrix, and let $(x,y)\in(0,1)^2$ be its distinguished Nahm point. Assume
\begin{enumerate}[label=(\roman*),leftmargin=*]
\item $K=\Q(x,y)$ has degree four and signature $(2,1)$;
\item $x,1-x,y,1-y$ are ordinary algebraic units;
\item $x$ and $1-x$ are multiplicatively independent.
\end{enumerate}
Let $E$ be the group of units of $K$ positive at the distinguished embedding. It is free of rank two. Define
\[
 I=[E:\langle x,1-x\rangle].
\]
An independent rank-two subgroup has saturation $E$, so $I$ is its saturation index.

Simultaneous complementation $(x,y)\mapsto(1-x,1-y)$ replaces $A$ by $A^{-1}$. Coordinate interchange is also allowed. The off-diagonal entry is nonzero: otherwise $1-x=x^{A_{11}}$ contradicts independence. After simultaneous complementation if necessary it is positive. Write
\begin{equation}\label{qt:mo:eq:normal}
 y=x^{-a}(1-x)^u,\qquad 1-y=x^{-c}(1-x)^j,
 \quad a,u,c,j\in\Q_{>0},\quad aj-uc=1.
\end{equation}
Then $A=u^{-1}\qtmat{a}{1}{j}$, and all four exponents belong to $I^{-1}\Z$.

\begin{theorem}[Practical finite reduction]\label{qt:mo:thm:height}
Under these hypotheses, suppose $I$ is odd. After interchanging $x$ and $y$ if necessary, their monic minimal polynomials satisfy
\[
 H(f_x)<700,\qquad H(f_y)<2700,
\]
where $H$ is the maximum absolute coefficient. For the same ordering,
\[
 \frac1{91}<x<\frac{90}{91},\qquad
 \frac1{404}<y<\frac{403}{404},\qquad I\le109.
\]
The bound on $I$ is odd. No torsion or modularity assumption is used.
\end{theorem}

\begin{theorem}[Computer-assisted classification]\label{qt:mo:thm:list}
Under the hypotheses of Theorem~\ref{qt:mo:thm:height}, the complete list of matrices, up to interchange and simultaneous complementation, is
\[
 \qtmat{1}{1}{4},\quad
 \qtmat{2}{3}{8},\quad
 \qtmat{4}{5}{8},\quad
 \qtmat{11}{9}{8},\quad
 \qtmat{16}{19}{24}.
\]
Each matrix has an independent four-unit quartic point of signature $(2,1)$ and odd saturation index.
\end{theorem}
Up to interchange, all five matrices occur in Zagier's list of rank-two
matrices whose Nahm point gives a two-term dilogarithm identity
\cite[p.~46]{Zagier}. Among these five matrices, only the first admits
modular rational shifts: it is $Z_4$ up to interchange, and
Proposition~\ref{qt:sh:finite} excludes the other four for all rational
shifts. Their completeness here follows from the height bound and the
exact finite reduction below.

\subsubsection{The arithmetic determinant constraint}
Let $\sigma_0$ be the distinguished real embedding and $\sigma_1$ the other real embedding. For a unit $v$, put
\[
 U(v)=\bigl(\log|\sigma_0(v)|,\log|\sigma_1(v)|\bigr)^t,
 \quad D_v=\det\bigl(U(v),U(1-v)\bigr).
\]
The usual regulator $R_K$ is computed with weight two at the complex place. Deleting that row leaves exactly these two real logarithm rows. Consequently
\begin{equation}\label{qt:mo:eq:det}
 |D_x|=I R_K,\qquad D_y=-D_x.
\end{equation}
The second equality follows from \eqref{qt:mo:eq:normal}; its coefficient matrix has determinant $-1$. Thus $y,1-y$ are independent and have the same index $I$.

More generally, each determinant choosing one member of $\{x,1-x\}$ and one of $\{y,1-y\}$ is a nonzero integer multiple of $R_K$. It is nonzero because it is one of $\pm aD_x,\pm uD_x,\pm cD_x,\pm jD_x$. This discrete restriction is essential; the determinant condition $aj-uc=1$ alone does not give it.

Friedman's theorem states that $R_K/w_K>1/10$ except for two totally complex quartic fields \cite[p.~600]{Friedman}. A field with a real embedding has $w_K=2$, so here
\begin{equation}\label{qt:mo:eq:reg}
 R_K>\frac15.
\end{equation}
Every cross determinant just described therefore has absolute value greater than $1/5$.

Independence implies $\Q(x)=\Q(y)=K$. Indeed, a proper subfield of a quartic field with a real embedding is rational or real quadratic and has unit rank at most one.

\subsubsection{Odd index restricts the second real conjugate}
\begin{lemma}\label{qt:mo:lem:real}
The second real conjugate $r$ of either $x$ or $y$ lies in $(-\infty,0)\cup(2,\infty)$. Its minimal polynomial has the form
\begin{equation}\label{qt:mo:eq:poly}
 f(T)=T^4+AT^3+BT^2+CT+e,
 \quad e=\pm1,\quad C=-1-2e-A-B,
\end{equation}
and has negative discriminant.
\end{lemma}
\begin{proof}
Suppose first $0<r<1$. Raise both equations \eqref{qt:mo:eq:normal} to the odd power $I$. At the second real embedding their right sides are positive, so both $y'$ and $1-y'$ are positive. Taking the positive $I$th roots recovers another positive solution of the same Nahm system. Such a solution is unique: in variables $v_i=-\log x_i$, the associated potential has Hessian
\[
 A+\operatorname{diag}\bigl((e^{v_i}-1)^{-1}\bigr)>0.
\]
Thus the second embedding agrees with the first on $x,y$, contradicting generation of $K$.

If $1<r<2$, \eqref{qt:mo:eq:normal} in absolute values gives $|y'|<1$ and $|1-y'|<1$, hence $0<y'<1$. The cleared equations now force both $Iu$ and $Ij$ to be even. This contradicts
\[
 (Ia)(Ij)-(Iu)(Ic)=I^2,
\]
whose right side is odd. The same proof applies with the two coordinates interchanged. Irreducibility excludes $r=0,1,2$.

Unit norms give $f(0),f(1)\in\{\pm1\}$. There is one real root in $(0,1)$ and the other is outside $[0,1]$, so these signs are opposite. This gives \eqref{qt:mo:eq:poly}. Signature $(2,1)$ gives negative discriminant.
\end{proof}

\subsubsection{Endpoint estimates from the complex triangle}
For this section only, replace a coordinate by its smaller complement if necessary. Denote the resulting distinguished value by $h\in(0,1/2]$ and set
\[
 \tau=-\log h,\qquad \delta(\tau)=-\log(1-e^{-\tau}).
\]
This temporary complementation preserves the absolute determinant and the unit hypotheses. Its other real conjugate $r$ is outside $(0,1)$, though it need not remain outside $(1,2)$. Write $u_0=|r|$, $v_0=|1-r|$. If $z$ is a nonreal conjugate, unit norms give
\begin{equation}\label{qt:mo:eq:triangle}
 p=|z|=(hu_0)^{-1/2},\qquad
 q=|1-z|=((1-h)v_0)^{-1/2},\qquad |p-q|<1<p+q.
\end{equation}

\begin{lemma}\label{qt:mo:lem:logs}
If $h\le1/64$, then
\[
 \tau-\log8<|\log v_0|<\tau+\log8,\qquad
 |\log u_0|<\tau+\log16.
\]
Consequently, with $D_h=\det(U(h),U(1-h))$,
\begin{align}
 \tau(\tau-\log8)-\delta(\tau)(\tau+\log16)
 &<|D_h|\label{qt:mo:eq:lower}\\
 &<\tau(\tau+\log8)+\delta(\tau)(\tau+\log16).
 \label{qt:mo:eq:upper}
\end{align}
\end{lemma}
\begin{proof}
If $v_0\ge1$, then $q\le\sqrt{64/63}<9/8$, so $p<17/8<\sqrt5$. Hence $u_0>1/(5h)$ and
\[
 v_0\ge u_0-1>1/(8h)\ge8.
\]
In particular $u_0\ge v_0/2$. From $p+q>1$ we obtain
\[
 v_0<(\sqrt{2/h}+\sqrt2)^2<4/h<8/h.
\]
The relation $|u_0-v_0|=1$ also gives $1/(16h)<u_0<16/h$.

If $v_0<1$, then $r>1$ and $u_0=1+v_0\in(1,2)$. Now $p>1/\sqrt{2h}$, so $1/p<1/4$ and $3p/4<q<5p/4$. Substituting \eqref{qt:mo:eq:triangle} yields
\[
 16h/25<v_0<32h/[9(1-h)]<4h.
\]
In particular $h/8<v_0<8h$, while $|\log u_0|<\log2$. These estimates prove the logarithm bounds. Finally
\[
 D_h=-\tau\log v_0+\delta(\tau)\log u_0
\]
gives \eqref{qt:mo:eq:lower}--\eqref{qt:mo:eq:upper} by the triangle inequalities.
\end{proof}

\subsubsection{One coordinate stays away from both endpoints}
\begin{lemma}\label{qt:mo:lem:away}
At least one of $x,y$ has $\min(v,1-v)>e^{-9/2}>1/91$.
\end{lemma}
\begin{proof}
Suppose both endpoint parameters are at least $9/2$. Write them as $t\le T$. Set
\[
 L=52/25,\qquad L_2=139/50,\qquad d_0=9/800.
\]
Exact exponential bounds give
\[
 \log8<L,\quad\log16<L_2,\quad
 \delta(9/2)<d_0,\quad 64<e^{9/2}<91.
\]
The lower bound in \eqref{qt:mo:eq:lower}, divided by $\tau^2$, is increasing; the corresponding upper bound is decreasing. Equality of $|D_x|$ and $|D_y|$ therefore gives
\[
 \frac{T}{t}<\frac{17}{10}.
\]
For completeness, this follows by substituting $\tau=9/2$ and checking the rational inequality
\[
 \frac{1+L/\tau+d_0(\tau+L_2)/\tau^2}
 {1-L/\tau-d_0(\tau+L_2)/\tau^2}<\left(\frac{17}{10}\right)^2.
\]

Choose the larger member of each complementary pair. The distinguished logarithms of these two units are $-\delta(t)$ and $-\delta(T)$. Their other real logarithms have absolute values below $t+\log8$ and $T+\log8$. Their cross determinant is therefore less in absolute value than
\[
 \delta(t)(T+\log8)+\delta(T)(t+\log8)
 \le\delta(t)\bigl((27/10)t+2\log8\bigr).
\]
The last expression decreases for $t\ge9/2$: one can use $\delta'(t)<-\delta(t)$. It is bounded by
\[
 \frac9{800}\left(\frac{27}{10}\frac92+2\frac{52}{25}\right)<\frac15.
\]
This contradicts the nonzero integer-multiple constraint \eqref{qt:mo:eq:reg}. Thus at least one endpoint parameter is below $9/2$.
\end{proof}

\subsubsection{From endpoint control to a coefficient cutoff}
\begin{lemma}\label{qt:mo:lem:mahler}
Let $v\in(0,1)$ be a quartic exceptional unit of signature $(2,1)$, whose other real conjugate is outside $(0,1)$. Then
\[
 M(v),M(1-v)<1+F(v),\qquad
 F(v)=\left(\frac1{\sqrt v}+\frac1{\sqrt{1-v}}\right)^2
 =\frac1{v(1-v)}+\frac2{\sqrt{v(1-v)}}.
\]
Here $M$ denotes Mahler measure.
\end{lemma}
\begin{proof}
Let $r$ be the other real conjugate, $u_0=|r|$, $v_0=|1-r|$, and let $p,q$ be as in \eqref{qt:mo:eq:triangle}, with $h$ replaced by $v$. If $u_0\ge1$, the norm identity gives $M(v)=\max(u_0,1/v)$. If $r<0$, then $v_0=u_0+1$ and $p+q>1$ forces $u_0<F(v)$. If $r>1$, the same argument applied to $v_0=u_0-1$ gives $v_0<F(v)$ and $u_0<1+F(v)$. If $u_0<1$, then $r<0$ and $v_0>1$, so
\[
 M(v)=p^2<(1+q)^2<\bigl(1+(1-v)^{-1/2}\bigr)^2<F(v).
\]
Apply the same argument to $1-v$ to finish.
\end{proof}

Choose the coordinate supplied by Lemma~\ref{qt:mo:lem:away} as $x$. Then
\[
 x(1-x)>90/8281,\qquad M(x),M(1-x)<113.
\]
For a monic quartic $f$, $H(f)\le6M(f)$; this follows by bounding each elementary symmetric sum term by the Mahler measure. Thus $H(f_x)<678<700$.

Every conjugate of an ordinary unit satisfies $|\log|v'||\le\log M(v)$: the norm identity supplies the reciprocal bound as well. Consequently
\[
 |D_x|<\log113\,[-\log(x(1-x))]
 <\frac{473}{100}\frac{453}{100}<22.
\]
If the endpoint parameter of $y$ were at least $6$, Lemma~\ref{qt:mo:lem:logs} would give
\[
 |D_y|>6(6-52/25)-\frac1{400}(6+139/50)>22,
\]
using $\delta(6)<1/400$ and monotonicity of the lower bound. This is impossible. Hence $\min(y,1-y)>e^{-6}>1/404$. Lemma~\ref{qt:mo:lem:mahler} now gives $M(y)<447$ and $H(f_y)<2682<2700$. Finally $I=|D_x|/R_K<110$, so odd $I\le109$. This proves Theorem~\ref{qt:mo:thm:height}. The strict numerical comparisons follow by summing the exponential series with rational remainder bounds; these 23 comparisons are also recorded in the exact certificate.

\subsubsection{The complete finite computation}
The two polynomial lists use \eqref{qt:mo:eq:poly}, coefficient heights at most $699$ and $2699$, and interval roots in $(1/91,90/91)$ and $(1/404,403/404)$ respectively. Negative discriminant, irreducibility and the real-root restrictions of Lemma~\ref{qt:mo:lem:real} are tested exactly. Endpoint signs are linear in $B$, allowing whole excluded ranges to be skipped with integer inequalities.

\begin{center}
\begin{tabular}{lr}\toprule
Stage&Count\\\midrule
Small-coordinate polynomials&11,362\\
Partner polynomials&107,798\\
Dependent polynomials removed from each list&4\\
Pairs reaching determinant comparison&17,941\\
Rejected by coefficient signs&13,356\\
Remaining logarithmic configurations&4,585\\
Rejected by odd denominator bound $I\le109$&4,575\\
Exact algebraic reconstructions&10\\
Classes after coordinate interchange&5\\\bottomrule
\end{tabular}
\end{center}

The four removed polynomials are
\[
 T^4+T-1,\quad T^4-4T^3+6T^2-5T+1,
\]
\[
 T^4+T^3-3T^2+3T-1,\quad T^4-5T^3+6T^2-4T+1.
\]
They give $1-x=x^4$, $x=(1-x)^4$, $x^4=(1-x)^3$, and $(1-x)^4=x^3$, respectively. Every other retained polynomial has a certified nonzero logarithm determinant; a small numerical determinant is never treated as exact dependence.

There is one root in $(0,1)$ and precisely one other real root, so no additional real-embedding alignment is omitted. Dyadic isolating intervals are accepted only after exact polynomial sign checks. All logarithms are enclosed by explicit series with outward integer rounding. A signed determinant grid discards only pairs whose determinant intervals cannot sum to zero. Cramer's rule then encloses $a,u,c,j$. The finite test loops through every odd $I\le109$, requiring $Ia,Iu,Ic,Ij$ to be positive integers and to satisfy the exact determinant equation. The intervals are narrower than $1/109^2$, so each contains at most one rational with denominator at most $109$.

The ten survivors are the five rows below and their coordinate interchanges. Here $y=P(x)$ in $\Q[T]/(f_x)$.
\begin{center}\small
\begin{tabular}{c l l}\toprule
$A$ entries $(A_{11},A_{12},A_{22})$&$f_x(T)$&$P(T)$\\\midrule
$(1,1,4)$&$T^4-6T^3+7T^2-4T+1$&$-T^3+6T^2-7T+3$\\[3pt]
$(2,3,8)$&$T^4-10T^3+13T^2-6T+1$&$-4T^3+38T^2-33T+7$\\[3pt]
$(4,5,8)$&$T^4-10T^3+12T^2-5T+1$&$(-9T^3+87T^2-79T+17)/5$\\[3pt]
$(11,9,8)$&$T^4-8T^3+3T^2+2T+1$&$(-9T^3+78T^2-79T+29)/17$\\[3pt]
$(16,19,24)$&$T^4-21T^3+27T^2-9T+1$&$(-64T^3+1320T^2-1233T+109)/37$\\\bottomrule
\end{tabular}
\end{center}
Exact embedded-field factorization identifies the chosen interval root of $f_y$, and polynomial remainder computations verify both cleared equations \eqref{qt:mo:eq:normal}. The exponents are
\[
 (1,1,3,4),\quad (2,1,7,8)/3,\quad (4,1,7,8)/5,
 \quad(11,1,7,8)/9,\quad(16,1,23,24)/19.
\]
These identities, together with positivity, prove the rational-power equations at the distinguished embedding. All matrices are positive definite.
Finally, exact factorization of $Z^2-x$, $Z^2-(1-x)$ and $Z^2-x(1-x)$ over $K$ proves that each is irreducible in all ten reconstructions. Thus the classes of $x,1-x$ are independent in $E/2E$, which proves that their actual saturation index is odd. This proves the required parity directly. It completes Theorem~\ref{qt:mo:thm:list}.

%% file: sections/quartic_mixed_even.tex
\subsection{Mixed signature and even index}\label{qt:mixed-even}
Let $A$ be symmetric, positive definite and rational, and let $(x,y)\in(0,1)^2$ solve its Nahm system. Assume that $[K:\Q]=4$ for $K=\Q(x,y)$, that $K$ has signature $(2,1)$ and all of $x,1-x,y,1-y$ are ordinary units. Assume also that $x,1-x$ are multiplicatively independent.

Let $E$ be the units positive at the distinguished real embedding. It is free of rank two. Let $\Lambda$ be the saturation of $\langle x,1-x\rangle$ in $E$, and let $I=[\Lambda:\langle x,1-x\rangle]$. Reflection sends $(A,x,y)$ to $(A^{-1},1-x,1-y)$; swapping coordinates also preserves the hypotheses and index. After reflection if needed the off-diagonal entry is positive, and
\begin{equation}\label{qt:me:normalized}
 y=x^{-a}(1-x)^u,\quad 1-y=x^{-c}(1-x)^j,
 \qquad a,u,c,j\in\Q_{>0},\quad aj-uc=1.
\end{equation}
The matrix is $u^{-1}\bigl(\begin{smallmatrix}a&1\\1&j\end{smallmatrix}\bigr)$.

\begin{theorem}[Computer-assisted classification]\label{qt:me:main}
In the preceding setting, if $I$ is even, then, up to coordinate swapping and simultaneous reflection, $A$ is one of
\[
 E_1=\begin{pmatrix}5/2&2\\2&2\end{pmatrix},\qquad
 E_2=\begin{pmatrix}3/2&1\\1&1\end{pmatrix},\qquad
 E_3=\begin{pmatrix}1&1/2\\1/2&5/4\end{pmatrix}.
\]
Each representative has an independent four-unit positive point of joint degree four and signature $(2,1)$.
\end{theorem}
The three representatives will be subjected to the rational-shift conditions in Section~\ref{qt:shifts}; the present result is an arithmetic classification of their Nahm points.

\subsubsection{Why one coordinate has a small square-root generator}
Write $U=\ell(x),V=\ell(1-x)$, with one logarithmic absolute value per embedding, including the conjugate complex embeddings separately. Define $U',V'$ similarly. Clearing denominators in \eqref{qt:me:normalized} gives
\[
 U'=-aU+uV,\quad V'=-cU+jV,\quad U'\wedge V'=-U\wedge V.
\]
The corresponding identity also holds in the exterior square of the rationalized unit group. Thus both complementary pairs have the same saturation index and the same area $R=\|U\wedge V\|_2$. Independence implies $\Q(x)=\Q(y)=K$, since a proper subfield of a quartic field with a real embedding has unit rank at most one.

If $I$ is even, the exponent vectors of $x,1-x$ are dependent modulo two in $\Lambda$. Therefore at least one of $x,1-x,x(1-x)$ is a square in $K$. The same holds for $y$. Unlike the totally real case, the product-square possibility cannot simply be discarded.

\begin{lemma}\label{qt:me:anchor}
At least one of $x,1-x,y,1-y$ is a square in $K$.
\end{lemma}
\begin{proof}
Otherwise both $x(1-x)$ and $y(1-y)$ would be squares. At the other real embedding $\sigma$, this forces $\sigma(x),\sigma(y)\in(0,1)$. The cleared Nahm identities then recover the positive rational powers, so this is another positive solution for the same matrix. The positive solution is unique: in variables $v_i=-\log x_i$, the Hessian of the potential is $A+\operatorname{diag}(1/(e^{v_i}-1))>0$. Thus $\sigma(x)=x,\sigma(y)=y$, contradicting $K=\Q(x,y)$ and distinct embeddings.
\end{proof}

Swap coordinates so the square occurs in the first pair. Write $z^2=x$ or $1-x$, with $z\in K\cap(0,1)$. Then $z,1-z,1+z$ are units and $z$ generates $K$. Its monic polynomial is
\begin{equation}\label{qt:me:split}
 g(Z)=Z^4+AZ^3+\left(\frac{p+m}{2}-1-e\right)Z^2+
 \left(\frac{p-m}{2}-A\right)Z+e,\qquad e,p,m\in\{\pm1\}.
\end{equation}
This follows from $g(0)=e$, $g(1)=p$, $g(-1)=m$.

\begin{lemma}\label{qt:me:seventeen}
The complete list of irreducible polynomials \eqref{qt:me:split} having signature $(2,1)$ and a root in $(0,1)$ contains seventeen polynomials. In all of them $|A|\le3$.
\end{lemma}
\begin{proof}
First suppose $H=|A|\ge2000$. Write $g(Z)=AZ(Z^2-1)+Z^4+bZ^2+dZ+e$, with $|b|\le3, |d|\le1$. At distance $4/H$ to either side of $1$ and $-1$, the first term has absolute value at least $4(1-4/H)(2-4/H)>7$, while the remainder has absolute value below $7$. At distance $4/H$ to either side of zero, the corresponding bounds are $>3$ and $<2$. The sign changes give three distinct real roots. The fourth root must be real too, contradicting signature $(2,1)$.

For $|A|<2000$, one checks all $31,992$ sign-parameter tuples by exact integer arithmetic. Rational roots can only be $\pm1$, which the endpoint conditions exclude; reducibility is decided by the possible monic quadratic factors. Rational Sturm sequences give all root counts. There are three reducible polynomials, thirty-four irreducible polynomials of signature $(2,1)$, and seventeen of the latter with an interval root. The seventeen parameter rows are exactly the following three disjoint ranges:
\[
\begin{array}{c|rrr}
 A&e&p&m\\\hline
 -2,-1,0,1,2,3&-1&1&-1\\
 -3,-2,-1,0,1,2,3&-1&1&1\\
 -2,-1,0,1&1&-1&1.
\end{array}
\]
Every retained polynomial has exactly one root in $(0,1)$.
\end{proof}

\subsubsection{The extra square class has a small partner bound}
If the partner also has a square coordinate or complement, Lemma~\ref{qt:me:seventeen} applies to it. Otherwise $y(1-y)$ is a square. Choosing a reciprocal if necessary gives $t\in K\cap(0,1)$ with
\[
 y=\frac{t^2}{1+t^2}\quad\hbox{or}\quad y=\frac1{1+t^2}.
\]
Both $t$ and $1+t^2$ are units, and $t$ generates $K$. The norm of $1+t^2$ is positive and hence equals one. Its monic polynomial must have the form
\begin{equation}\label{qt:me:product}
 h(T)=T^4+AT^3+(1+e-\epsilon)T^2+(A+\delta)T+e,
 \quad e=\pm1,\quad (\epsilon,\delta)\in\{(\pm1,0),(0,\pm1)\}.
\end{equation}
Indeed, reduction modulo $T^2+1$ gives $h(i)=\epsilon+i\delta$; the resultant is $\epsilon^2+\delta^2=1$.

\begin{lemma}\label{qt:me:partner}
Every product-square partner in this branch satisfies $|A|<16$ in \eqref{qt:me:product}.
\end{lemma}
\begin{proof}
For a split generator $z$ from the seventeen-polynomial list, let $w=z^2$. Its monic norm polynomial is
\[
 F(T)=T^4+(2b-A^2)T^3+(b^2+2e-2Ac)T^2+(2be-c^2)T+1,
\]
where $g(Z)=Z^4+AZ^3+bZ^2+cZ+e$. Exact substitution on the seventeen rows gives
\[
 \|F\|_2^2\le461,\qquad \|F(1-T)\|_2^2\le197.
\]
The Mahler measure inequality $M(F)\le\|F\|_2$ follows from Jensen's inequality and the circle $L^2$ norm. For a unit, the logarithmic vector has length at most $\sqrt2\log M$. Therefore, for either square/complement choice,
\begin{equation}\label{qt:me:upper}
 R\le2\log M(F)\log M(F(1-T))
 \le\tfrac12\log461\log197<18.
\end{equation}

Now let $H=|A|\ge16$ in \eqref{qt:me:product}. Write the polynomial as $AT(T^2+1)+T^4+BT^2+\delta T+e$, with $|B|\le3$. At $T=\pm2/H$, the first term has magnitude at least $2$, while the remainder has magnitude at most
\[
 1+2/H+12/H^2+16/H^4<2.
\]
There is a real root $r_0$ with $|r_0|<2/H$. At $T=-A\pm3$, write the polynomial instead as $T^3(T+A)+BT^2+(A+\delta)T+e$. The leading term dominates because
\[
 3(H-3)^3>3(H+3)^2+(H+1)(H+3)+1\quad(H\ge16).
\]
Thus a second real root satisfies $H-3<|r_\infty|<H+3$. These are the two real roots in signature $(2,1)$.

For $q=r^2$, the logarithms of $q/(1+q)$ and $1/(1+q)$ are $-\log(1+1/q)$ and $-\log(1+q)$. The determinant of their two real rows consequently satisfies
\[
 D>4\log(H/2)\log(H-3)-\frac4{H^2(H-3)^2}.
\]
The two complex rows are equal and are minus half the sum of the real rows. Hence $R=\sqrt2|D|$. Complementation changes only its sign. At $H\ge16$, using $\log8>2$, $\log13>2$ and $\sqrt2>7/5$, this gives $R>18$, contradicting \eqref{qt:me:upper}. The numerical comparisons may be verified by rational exponential-series bounds.
\end{proof}

There are $248$ sign-parameter tuples in \eqref{qt:me:product} with $|A|<16$. Exact factor and Sturm checks leave $116$ polynomials of signature $(2,1)$ with a root in $(0,1)$. Thus the finite list contains every possible product-square partner.

\subsubsection{The finite pairing computation}
Three split generators give dependent pairs and are removed by exact identities:
\[
\begin{array}{c|c}
 g(Z)&\text{identity at its interval root}\\\hline
 Z^4+Z^2-1&1-Z^2=Z^4\\
 Z^4-Z^3-2Z^2+1&(1-Z^2)^2=Z^3\\
 Z^4-2Z^2-Z+1&(1-Z^2)^2=Z
\end{array}
\]
The product-family polynomial $Z^4+Z^2-1$ also gives a dependent pair: if $v=Z^2$, then $v+v^2=1$ and $Z^2/(1+Z^2)=v^2$. These exact identities justify the dependent-pair removals.

There remain fourteen split anchors and $14+115=129$ possible partners. Each polynomial has exactly one positive interval root. For every eligible alignment of the two real embeddings and all four complement choices, the finite test checks $D_x+D_y=0$, then the required positivity of the four coefficients in \eqref{qt:me:normalized}. It rejects identical generators with identical alignment and opposite complements by the exact equality $y=1-x$, which would force $a=0$ under independence.

\begin{center}
\begin{tabular}{lr}\toprule
Stage&Count\\\midrule
Square/complement configurations&7,224\\
Rejected by a nonzero determinant-sum interval&7,140\\
Rejected by a coefficient-sign interval&44\\
Rejected as exact complementary pairs&28\\
Remaining necessary-condition configurations&12\\
Rejected by exact embedded-field computation&7\\
Exact Nahm reconstructions&5\\
\bottomrule
\end{tabular}
\end{center}

The five reconstructions are $E_1$, its swap, $E_2$, $E_3$, and the swap of $E_3$. They represent three classes under the allowed symmetries. Each reconstruction supplies an exact polynomial expression for the partner generator and verifies both cleared equations modulo the anchor polynomial. The identities prove the rational exponents exactly. The nonzero real logarithm determinant certifies independence; positive definiteness follows from \eqref{qt:me:normalized}.

The seven other configurations are excluded in the distinguished real field.  Let $z$ be the unique root of the first polynomial in $(0,1)$.  Exact factorization of the proposed partner polynomial over $\Q(z)$ has degrees $(1,3)$ in each of these seven cases.  The unique linear-factor root is $P(z)$ in the following table.  Rational Horner evaluation on a Sturm isolating interval gives the displayed strict outer interval, which is disjoint from $(0,1)$.
\begin{center}\scriptsize
\begin{tabular}{llll}\toprule
Anchor polynomial&Partner polynomial&$P(z)$&Outer interval\\\midrule
$z^{4} - 2 z^{3} + 3 z - 1$&$z^{4} + 3 z^{3} - 2 z - 1$&$- z^{3} + 2 z^{2} - 3$&$(- \frac{14}{5},- \frac{27}{10})$\\
$z^{4} + z^{3} - 1$&$z^{4} + z - 1$&$- z^{3} - z^{2}$&$(- \frac{13}{10},- \frac{6}{5})$\\
$z^{4} + 3 z^{3} - 2 z - 1$&$z^{4} - 2 z^{3} + 3 z - 1$&$- z^{3} - 3 z^{2} + 2$&$(- \frac{6}{5},- \frac{11}{10})$\\
$z^{4} + z^{3} + z^{2} - z - 1$&$z^{4} - 2 z^{3} - 2 z^{2} + z + 1$&$z^{2} + z + 1$&$(\frac{5}{2},\frac{13}{5})$\\
$z^{4} - 2 z^{3} - 2 z^{2} + z + 1$&$z^{4} + z^{3} + z^{2} - z - 1$&$- z^{3} + 3 z^{2} - 2$&$(- \frac{7}{10},- \frac{3}{5})$\\
$z^{4} - 2 z^{3} - 2 z^{2} + z + 1$&$z^{4} + z^{3} - 2 z^{2} - 2 z + 1$&$- z^{3} + 2 z^{2} + 2 z - 1$&$(\frac{6}{5},\frac{13}{10})$\\
$z^{4} + z^{3} - 2 z^{2} - 2 z + 1$&$z^{4} - 2 z^{3} - 2 z^{2} + z + 1$&$- z^{3} - z^{2} + 2 z + 2$&$(\frac{5}{2},\frac{13}{5})$\\
\bottomrule
\end{tabular}
\end{center}
The target polynomial has exactly one root in $(0,1)$, so this proves that the required positive partner is absent from the anchor field.  It also explains why an abstract field isomorphism is insufficient for this step.

For the three retained representatives, exact point models are as follows.  In each row $z$ is the unique root of $g$ in $(0,1)$.
\begin{center}\small
\begin{tabular}{clll}\toprule
Matrix&$g(z)$&$x$&$y$\\\midrule
$E_1$&$z^4+2z^3-z-1$&$z^2$&$2+z-2z^2-z^3$\\
$E_2$&$z^4+z^3-1$&$z^2$&$1+z-z^2-z^3$\\
$E_3$&$z^4-z^3+z^2+z-1$&$1-z^2$&$2-3z+2z^2-z^3$\\\bottomrule
\end{tabular}
\end{center}
Reduction modulo $g$ verifies the cleared Nahm equations, and root isolation gives $0<x,y<1$.  The characteristic polynomials of $x,1-x,y,1-y$ are monic integral of degree four with constant terms $\pm1$.  A nonzero determinant of the two real logarithm rows proves independence.  Each row has a square coordinate or complement, so its saturation index is even.  This proves Theorem~\ref{qt:me:main}.

All finite logarithm comparisons use isolating intervals verified by exact polynomial signs, followed by the explicit rational logarithm series and outward integer rounding.  The complex logarithmic absolute values follow from the unit norms, so checking the two real rows includes every embedding alignment.  A determinant interval is used for exclusion only when it is separated from zero; retained configurations are resolved by the exact identities above.

%% file: sections/quartic_shifts.tex
\subsection{All rational shifts on the quartic unit list}\label{qt:shifts}
The arithmetic results of Sections~\ref{qt:totally-real}--\ref{qt:mixed-even}, together with the dependent exclusion, leave the following eight matrix classes.  No modularity is asserted by this arithmetic list.
Let $\mathcal M$ consist of the eight matrices $M_i=\bigl(\begin{smallmatrix}a&b\\b&d\end{smallmatrix}\bigr)$ specified by
\begin{center}
\begin{tabular}{crrrl}\toprule
$i$&$a$&$b$&$d$&Necessary obstruction for $M_i^{\pm1}$\\\midrule
1&1&1&4&Two shifts survive for each orientation\\
2&2&3&8&Nonsquare discriminant\\
3&4&5&8&Second coefficient excludes two rational shifts\\
4&11&9&8&Negative discriminant\\
5&16&19&24&Nonsquare discriminant\\
6&$5/2$&2&2&Negative discriminant\\
7&$3/2$&1&1&First coefficient equations inconsistent\\
8&1&$1/2$&$5/4$&First coefficient equations inconsistent\\\bottomrule
\end{tabular}
\end{center}

\begin{proposition}\label{qt:sh:finite}
For $A\in\{M_i,M_i^{-1}:1\le i\le8\}$ and arbitrary rational shifts, the modular Nahm data are exactly the four triples in Theorem~\ref{qt:classification}, up to simultaneous coordinate interchange.
\end{proposition}

The exact point models for $M_1,\ldots,M_5$ are the five rows of Section~\ref{qt:mixed-odd}, and those for $M_6,M_7,M_8$ are $E_1,E_2,E_3$ in Section~\ref{qt:mixed-even}. For $M_i^{-1}$, use $(1-x,1-y)$ in the same field. Substituting these models in the logarithmic equations proves this inversion statement about points. All sixteen matrix orientations are checked separately in the coefficient calculation.

Write the two shift variables as $(B,C)$ and reserve $c_0$ for the scalar. The normalization is
\[
 \log f_{A,(B,C),c_0}(e^{-\epsilon})\sim
 \frac{\alpha_0}{\epsilon}+\log\beta_0+
 (\lambda_1-c_0-1/12)\epsilon+\lambda_2\epsilon^2+\cdots.
\]
By the established asymptotic expansion and modular necessary conditions \cite[Theorem 2.3 and Corollary 3.2]{VZ},
\begin{equation}\label{qt:sh:necessary}
 \lambda_1\in\Q,\qquad \lambda_2=0,\qquad c_0=\lambda_1-1/12.
\end{equation}
Writing each coefficient in the basis $1,z,z^2,z^3$ converts \eqref{qt:sh:necessary} into exact polynomial equations over $\Q$ in $B,C$. There is no bound on their denominators.

\paragraph{Explicit coefficient calculation.}
Put $q_1=x/(1-x)$, $q_2=y/(1-y)$, $b_1=B$, $b_2=C$, and let $Z=(Z_1,Z_2)$ be a centered Gaussian with covariance
\[
 \Sigma=(A+\operatorname{diag}(q_1,q_2))^{-1}.
\]
Define $P_0(t)=t$, $P_{k+1}(t)=t(1+t)P_k'(t)$ and
\begin{align*}
 U_1&=\sum_{i=1}^2\left[-(b_i+q_i/2)Z_i+\frac{P_1(q_i)}6Z_i^3\right],\\
 U_2&=\sum_{i=1}^2\left[-\frac{q_i}{12}+\frac{P_1(q_i)}4Z_i^2-\frac{P_2(q_i)}{24}Z_i^4\right],\\
 U_3&=\sum_{i=1}^2\left[\frac{P_3(q_i)}{120}Z_i^5-\frac{P_2(q_i)}{12}Z_i^3+\frac{P_1(q_i)}{12}Z_i\right],\\
 U_4&=\sum_{i=1}^2\left[-\frac{P_4(q_i)}{720}Z_i^6+\frac{P_3(q_i)}{48}Z_i^4-\frac{P_2(q_i)}{24}Z_i^2\right].
\end{align*}
Then
\begin{align*}
 T_1&=\mathbb E(U_2+U_1^2/2),\\
 T_2&=\mathbb E(U_4+U_1U_3+U_2^2/2+U_1^2U_2/2+U_1^4/24),\\
 \lambda_1&=T_1,\qquad \lambda_2=T_2-T_1^2/2.
\end{align*}
Every moment is obtained by summing products of entries of $\Sigma$ over pairings; odd total degree has zero moment.  Thus these formulas and reduction modulo the displayed quartics give an exact, reproducible calculation in $\Q[z]/(f)[B,C]$.  The denominators are invertible because $0<x,y<1$ and $A$ is positive definite.  In particular, the calculation does not impose a height or denominator bound on the shifts.

\subsubsection{The first coefficient almost finishes the exclusion}
For rows $1$--$6$, rational row operations reduce the three nonconstant coefficients of $\lambda_1$ to a linear relation in $B,C$ and one quadratic in $C$. The following tables list an equivalent system; scaling an equation by a nonzero rational constant is immaterial.

\begin{center}\small
\begin{tabular}{cll}\toprule
Matrix&Linear equation&Quadratic equation\\\midrule
$M_1$&$2B-1=0$&$C^2-2C=0$\\
$M_2$&$4B-C-2=0$&$49C^2-217C+90=0$\\
$M_3$&$2B-C-1=0$&$361C^2-1748C+1155=0$\\
$M_4$&$2B-2C-1=0$&$121C^2-902C+2710=0$\\
$M_5$&$4B-3C-2=0$&$363C^2-5654C+11381=0$\\
$M_6$&$2B-C-1=0$&$121C^2-319C+240=0$\\\bottomrule
\end{tabular}
\end{center}
\begin{center}\small
\begin{tabular}{cll}\toprule
Matrix&Linear equation&Quadratic equation\\\midrule
$M_1^{-1}$&$2B+2C-1=0$&$12C^2-4C-1=0$\\
$M_2^{-1}$&$5B+4C-2=0$&$9604C^2-3479C-624=0$\\
$M_3^{-1}$&$3B+2C-1=0$&$17689C^2-5719C-1700=0$\\
$M_4^{-1}$&$4B+2C-1=0$&$23716C^2-20020C+20689=0$\\
$M_5^{-1}$&$7B+4C-2=0$&$3072432C^2-1103080C-422123=0$\\
$M_6^{-1}$&$3B+2C-1=0$&$242C^2-473C+365=0$\\\bottomrule
\end{tabular}
\end{center}

The discriminant square classes for rows $2,4,5,6$, in either orientation, are respectively
\[
 601,\qquad -21,\qquad 31906,\qquad -119.
\]
The two positive values are nonsquares ($24^2<601<25^2$ and $178^2<31906<179^2$). Thus none of these four classes has a rational shift satisfying even the first condition in \eqref{qt:sh:necessary}.

\paragraph{Rows 7 and 8.}
For each orientation the rational row reduction gives $B=h(C)$. Substituting into the other two rows gives polynomials with gcd one. The elimination is exact over $\Q[C]$; short coprimality certificates modulo $7$ are displayed below. The reductions preserve both degrees; consequently their nonzero resultants modulo $7$ imply nonzero rational resultants. This argument does not assume that $B$ or $C$ is $7$-integral.

\begin{center}\small
\begin{tabular}{cl}\toprule
Matrix&$h(C)$\\\midrule
$M_7$&$2/3-7C/12$\\
$M_7^{-1}$&$8/25-19C/25$\\
$M_8$&$1324C^2/565+363C/565-1361/13560$\\
$M_8^{-1}$&$-331C^2/30-419C/60-61/720$\\\bottomrule
\end{tabular}
\end{center}
For the primitive substituted polynomials $F,G$ reduced modulo $7$, the following rows give $uF+vG=1$ in $\mathbb F_7[C]$:
\begin{center}\small
\begin{tabular}{cll}\toprule
Matrix&$F,G$&$u,v$\\\midrule
$M_7$&$-2C^2+2C+2$&$3C-1$\\
 &$2-3C^2$&$-2C-2$\\[2pt]
$M_7^{-1}$&$3C^2-2C-1$&$-2C-2$\\
 &$C^2+C+2$&$3-C$\\[2pt]
$M_8$&$-3C^4-C^3-2C^2+2C$&$-3C^2-2C-1$\\
 &$-C^3-3C-2$&$2C^3+2C^2-2C+3$\\[2pt]
$M_8^{-1}$&$2C^4-C^3+3C^2+2C-3$&$-3C^2-2C+2$\\
 &$C^3-3C^2+3C$&$-C^3-2C^2-1$\\\bottomrule
\end{tabular}
\end{center}
Thus these two classes also fail the first coefficient condition. Thus the coefficient condition itself excludes both classes, in both orientations.
\subsubsection{The second coefficient removes the last unwanted class}
Only $M_1^{\pm1}$ and $M_3^{\pm1}$ have rational shifts satisfying $\lambda_1\in\Q$. For $M_3$ the two pairs are $(17/19,15/19)$ and $(48/19,77/19)$. For $M_3^{-1}$ they are $(61/133,-25/133)$ and $(-1/133,68/133)$.
For $M_3$, use the exact coordinate generator
\[
 f(z)=z^4-10z^3+12z^2-5z+1,\qquad x=z,\quad
 y=\frac{-9z^3+87z^2-79z+17}{5}.
\]
The inverse uses the complements in the same basis. The constant coefficient of $\lambda_2$ is
\[
\begin{array}{c|c|r}
A&(B,C)&[z^0]\lambda_2\\\hline
M_3&(17/19,15/19)&-13796/3258025\\
M_3&(48/19,77/19)&27279/3258025\\
M_3^{-1}&(61/133,-25/133)&-13796/3258025\\
M_3^{-1}&(-1/133,68/133)&27279/3258025.
\end{array}
\]
Every value is nonzero. Since $1,z,z^2,z^3$ is a rational basis, $\lambda_2$ cannot vanish. Only one nonzero basis coefficient is needed for each exclusion.
For $M_1$, the pairs $(1/2,0),(1/2,2)$ survive and give $c_0=1/120,49/120$. For $M_1^{-1}$, the pairs $(0,1/2),(2/3,-1/6)$ survive and give $c_0=1/120,3/40$. Swapping coordinates puts them in the format of Theorem~\ref{qt:classification}. The four entire second coefficients vanish exactly. Their known modular identities establish sufficiency; checking finitely many asymptotic coefficients alone would not do so.
\subsubsection{Published modular identities and their normalization}
The necessary conditions above leave two shifts for $M_1$ and two for $M_1^{-1}$. Their sufficiency follows from the modular identities in Wang~\cite[Theorem 3.4 and Table 5]{Wang} for the matrix $\bigl(\begin{smallmatrix}4&1\\1&1\end{smallmatrix}\bigr)$, and Cao--Rosengren--Wang~\cite[Theorems 1.3--1.4]{CRW} for its inverse. The latter proof is needed: the inverse identities in Wang's paper were conjectural.
For clarity, in the inverse orientation used in Theorem~\ref{qt:classification}, replacing the variable $q$ by $q^3$ changes the two sums into
\begin{align*}
 f_{A,(1/2,0),1/120}(q^3)
 &=q^{1/40}\sum_{m,n\ge0}
 \frac{q^{(m^2-2mn+4n^2+3m)/2}}{(q^3;q^3)_m(q^3;q^3)_n},\\
 f_{A,(-1/6,2/3),3/40}(q^3)
 &=q^{9/40}\sum_{m,n\ge0}
 \frac{q^{(m^2-2mn+4n^2-m+4n)/2}}{(q^3;q^3)_m(q^3;q^3)_n}.
\end{align*}
These are the two normalized series evaluated in \cite{CRW}, with the scalar factors $q^{3c_0}$ shown explicitly.  The published identities identify them with modular forms.  Rational dilation of the modular variable preserves modularity on an appropriate congruence subgroup.  Together with Wang's two identities, this establishes sufficiency and completes Proposition~\ref{qt:sh:finite}.

%% file: sections/degree_four_role.tex
\subsection{The role of the degree-four bound}\label{qt:degree-role}
The bound in Theorem~\ref{deg:thm:main} is an upper bound on the joint
degree. Its role in the unit arguments can be stated precisely.

\begin{corollary}[Unit rank of the surviving fields]\label{qt:degree:rank}
Let $A$ be a symmetric positive definite rational $2\times2$ matrix,
and let $(x,y)\in(0,1)^2$ be its distinguished Nahm point. Suppose that
$x,1-x,y,1-y$ are ordinary units, that $x,1-x$ are multiplicatively
independent, and that $[K:\Q]\le4$ for $K=\Q(x,y)$. Then $K$ is either
a totally real cubic field or a quartic field of signature $(2,1)$.
In particular, the group of units positive at the distinguished
embedding has rank exactly two, and $\langle x,1-x\rangle$ has finite
index in it. Modularity is not an additional hypothesis of this
corollary.
\end{corollary}
\begin{proof}
The field has a real embedding, and its unit rank is at least two by
independence. Dirichlet's unit theorem, in the form
$d=r_1+2r_2$ and $\operatorname{rank}\mathcal O_K^\times=r_1+r_2-1$,
excludes degrees one and two and the cubic signature $(1,1)$.
The remaining possibilities are $(d;r_1,r_2)=(3;3,0),(4;4,0),(4;2,1)$.
Theorem~\ref{qt:tr:thm:quartic-exclusion} excludes $(4;4,0)$.
Both surviving signatures have unit rank two. Positivity at one real
embedding selects one representative modulo the torsion subgroup
$\{\pm1\}$, so it does not change the rank.
\end{proof}

The comparison with higher degrees is given in
Appendix~\ref{app:degree-four-consequences}. It separates the unit-rank
consequence above from the additional height estimates needed for a
finite classification. That appendix also records an elementary example
with five real Nahm solutions, which explains why counting all real
solutions cannot replace the arithmetic degree argument.

%% file: sections/degree_four_consequences.tex
\section{The degree cutoff, regulators and real solutions}
\label{app:degree-four-consequences}
The classification uses the joint-degree bound of four together with
height and exponent bounds in each remaining arithmetic branch.
Corollary~\ref{qt:degree:rank} describes one consequence of that cutoff:
the surviving independent four-unit fields have unit rank two. This
appendix explains the limits of that observation in degree five. None
of the results below is needed for the classification theorem.

\subsection{Logarithmic planes in totally real fields}

The totally real covolume argument itself has no threshold at degree
four. The following generalization of Lemma~\ref{qt:tr:lem:gap}
quantifies this point.

\begin{lemma}[A plane bound in arbitrary totally real degree]
\label{qt:degree:plane}
Let $K$ be a totally real field of degree $d\ge3$, and let $E$ be its
group of units positive at one fixed embedding. Use the logarithmic
map $\ell(w)=(\log|\sigma_i(w)|)_{i=1}^d$ with its Euclidean metric.
For every rank-two subgroup $\Lambda\subset E$,
\begin{equation}\label{qt:degree:gap}
 \operatorname{covol}(\ell(\Lambda))\ge c_d,
 \qquad c_d=\frac{\sqrt3\,d}{2}(\log\varphi)^2,
 \qquad \varphi=\frac{1+\sqrt5}{2}.
\end{equation}
If $u,v$ are independent, $\Lambda$ is their saturation in $E$,
$I=[\Lambda:\langle u,v\rangle]$, and
$R=\|\ell(u)\wedge\ell(v)\|_2$, then $I\le R/c_d$.
In particular,
\begin{align*}
 c_4&=2\sqrt3(\log\varphi)^2=0.802164\ldots>\frac45,
 & I&<\frac54R \quad(d=4),\\
 c_5&=\frac{5\sqrt3}{2}(\log\varphi)^2=1.002705\ldots>1,
 & I&<R \quad(d=5).
\end{align*}
\end{lemma}
\begin{proof}
For $w\in E\setminus\{1\}$ put $l_i=\log|\sigma_i(w)|$.
None of the $l_i$ is zero: a conjugate $\pm1$ would force $w=\pm1$.
Since $w$ is a unit, $\sum_i l_i=0$. Since $w^2-1$ is a nonzero
algebraic integer, the argument of Lemma~\ref{qt:tr:lem:gap} gives
\[
 0\le\log|N_{K/\Q}(w^2-1)|
   =\sum_{i=1}^d\log(2\sinh|l_i|)
   \le d\log\left(2\sinh\frac{\|\ell(w)\|_2}{\sqrt d}\right).
\]
The last inequality uses concavity and monotonicity of
$\log(2\sinh t)$, followed by Cauchy--Schwarz. Thus every nonzero
logarithmic vector has length at least $\sqrt d\log\varphi$.
The planar shortest-vector inequality
$\operatorname{covol}\ge\sqrt3\lambda_1^2/2$ proves
\eqref{qt:degree:gap}, and $R=I\operatorname{covol}(\ell(\Lambda))$
proves the index estimate. For completeness, the numerical thresholds can be certified by
rational inequalities. Since $\sqrt3>433/250$ and
$\sqrt5-2>59/250$, the positive power series for
$2\operatorname{arctanh}z$ gives
\[
 \log\varphi
 =2\operatorname{arctanh}(\sqrt5-2)
 >2\left(\frac{59}{250}+\frac13\left(\frac{59}{250}\right)^3\right)
 =\frac{11267879}{23437500}.
\]
Therefore
\[
 c_4>\frac{54975887069691553}{68664550781250000}>\frac45,
 \qquad c_5=\frac54c_4>1.
\]
\end{proof}

\subsection{What the split-unit conditions allow in degree five}
\begin{proposition}[Quintic split units with unbounded height]
\label{qt:degree:quintic-family}
For every odd integer $t\ge3$, the polynomial
\begin{equation}\label{qt:degree:quintic-polynomial}
 g_t(T)=T^5+2tT^4-tT^3-2tT^2+(t-1)T+1
\end{equation}
is irreducible and totally real. It has a root
$z_t\in(1/2,3/4)$ such that $z_t,1-z_t,1+z_t$ are ordinary units and
\[
 h(z_t)\ge\frac15\log\frac{2t}{5}\longrightarrow\infty.
\]
Consequently $x_t=z_t^2\in(1/4,9/16)$ and $1-x_t$ are units,
while $h(x_t)\to\infty$. Thus these split-unit conditions in degree
five, even with the distinguished coordinate in a fixed compact
subinterval of $(0,1)$, do not give an absolute height bound.
\end{proposition}
\begin{proof}
For odd $t$, reduction modulo two gives $T^5+T^3+1$. This polynomial
has no root in $\F_2$, and its remainder modulo $T^2+T+1$ is $T+1$.
A reducible polynomial of degree five has an irreducible factor of
degree at most two, so this reduction is irreducible.

Direct evaluation gives
\begin{gather*}
 g_t(-1)=g_t(0)=g_t(1)=1,\qquad g_t(1/2)=\frac{17}{32}>0,\\
 g_t(-1/2)=\frac{47-24t}{32}<0,\qquad
 g_t(3/4)=\frac{499-168t}{1024}<0.
\end{gather*}
Together with the negative sign at $-\infty$, these signs put five
distinct real roots in
\[
 (-\infty,-1),\quad(-1,-1/2),\quad(-1/2,0),\quad
 (1/2,3/4),\quad(3/4,1).
\]
The three endpoint values imply the asserted unit conditions. If
$M(g_t)$ is the Mahler measure, the coefficient of $T^4$ gives
$2t\le5M(g_t)$. Irreducibility yields
$5h(z_t)=\log M(g_t)$, which proves the lower bound.
Finally $1-z_t^2=(1-z_t)(1+z_t)$ is a unit and
$h(z_t^2)=2h(z_t)$. Moreover
$[\Q(z_t):\Q(z_t^2)]\le2$ and this index divides the odd prime
$[\Q(z_t):\Q]=5$, so $\Q(z_t^2)=\Q(z_t)$ and $x_t$ also has
degree five.
\end{proof}

\begin{remark}[What changes in degree five]\label{qt:degree:scope}
The possible signatures in degrees four and five, with a real
embedding, are as follows. The ranks refer to the full unit group.
\begin{center}
\begin{tabular}{ccc}\toprule
Degree & Signature $(r_1,r_2)$ & Unit rank $r_1+r_2-1$\\\midrule
4 & $(4,0)$ & 3\\
4 & $(2,1)$ & 2\\\midrule
5 & $(5,0)$ & 4\\
5 & $(3,1)$ & 3\\
5 & $(1,2)$ & 2\\\bottomrule
\end{tabular}
\end{center}
Thus degree five does not force the unit rank above two. In signature
$(1,2)$, take the real embedding $\sigma$ and one complex embedding
$\tau$, and let $E$ denote the group of units positive at $\sigma$.
For independent $u,v\in E$ and
$I=[E:\langle u,v\rangle]$, the usual regulator normalization gives
\[
 \left|\det\begin{pmatrix}
 \log|\sigma(u)|&\log|\sigma(v)|\\
 2\log|\tau(u)|&2\log|\tau(v)|
 \end{pmatrix}\right|=I R_K.
\]
Friedman's bound $R_K/w_K>1/10$ has only two totally complex
exceptions \cite[p.~600]{Friedman}; hence $R_K>1/5$ still holds here,
since $w_K=2$. The full-regulator determinant mechanism therefore
persists in this quintic signature. In the totally real quintic
signature the ambient rank is four, but
Lemma~\ref{qt:degree:plane} still controls the rank-two saturated plane.

What the quartic cutoff supplies is a complete case division covered
by the preceding arguments. In particular, the non-totally-real
independent unit branch has exactly two real embeddings and one
complex pair, as used in Section~\ref{qt:mixed-odd}. The quintic
signatures introduce additional configurations, and the polynomial
families also acquire an extra parameter. For example, fixing
$g(0)=e$, $g(1)=p$, $g(-1)=m$ with $e,p,m\in\{\pm1\}$ leaves the
one-parameter quartic family \eqref{qt:tr:eq:family}; the corresponding
monic quintic family is
\begin{align*}
 g(T)={}&T^5+A T^4+B T^3
 +\left(\frac{p+m}{2}-e-A\right)T^2\\
 &+\left(\frac{p-m}{2}-1-B\right)T+e,
 \qquad A,B\in\Z.
\end{align*}
These are necessary endpoint conditions for the split-unit
square-root branch; they do not assert irreducibility or the Nahm
equations. More generally, the three fixed endpoint values leave
$d-3$ free coefficients in degree $d\ge3$.
Proposition~\ref{qt:degree:quintic-family} is the subfamily
$e=p=m=1$, $A=2t$, $B=-t$. It shows concretely why the quartic
endpoint-to-height estimates cannot be extended using these
split-unit conditions alone. The proposition supplies individual
units, not a second complementary pair satisfying the symmetric
positive definite Nahm equations, and makes no claim of infinitely
many modular data.

Consequently, an upper degree bound of five alone would leave
additional arithmetic and height arguments to be supplied; the
quartic enumerations in this paper would not complete that enlarged
case division. This is a statement about the scope of the proved
finite reduction. It is not an impossibility theorem for finite
bounds in higher degree. A positive covolume lower bound controls
the index relative to $R$; an absolute finite search additionally
requires the upper height and exponent bounds established above.
\end{remark}

\subsection{Five real solutions need not form a quintic orbit}
\begin{proposition}[Five real solutions in rank two]
\label{qt:degree:five-real}
The symmetric positive definite integer matrix
\[
 A=\begin{pmatrix}2&-16\\-16&130\end{pmatrix}
\]
has at least five distinct real solutions to its Nahm equations, with
all coordinates different from $0$ and $1$. Thus rank two and positive
definiteness do not imply an upper bound of four on the total number of
real solutions.
\end{proposition}
\begin{proof}
The leading principal minors are $2$ and $4$, so $A$ is positive
definite. Its equations, with integer exponents throughout, are
equivalent for $xy\ne0$ to
\begin{equation}\label{qt:degree:five-equations}
 (1-x)y^{16}=x^2,\qquad (1-y)x^{16}=y^{130}.
\end{equation}
First put $x=-u$, where $u>0$, and
\[
 Y(u)=\left(\frac{u^2}{1+u}\right)^{1/16}>0.
\]
The first equation holds with $y=Y(u)$. Since
$u^{16}=(1+u)^8Y(u)^{128}$, the second is equivalent to
\[
 F(u)=(1-Y(u))(1+u)^8-Y(u)^2=0.
\]
For $0<u<\varphi=(1+\sqrt5)/2$, one has $0<Y(u)<1$.
The function extends continuously to the endpoints, with $F(0)=1$
and $F(\varphi)=-1$. At $u=1/20$, the exact inequality
$(137/200)^{16}<1/420=Y(1/20)^{16}$ gives
\[
 F(1/20)<\frac{63}{200}\left(\frac{21}{20}\right)^8
             -\left(\frac{137}{200}\right)^2
          =-\frac{19591860257}{5120000000000}<0.
\]
At $u=1/2$, the inequality $(9/10)^{16}>1/6=Y(1/2)^{16}$
gives
\[
 F(1/2)>\frac1{10}\left(\frac32\right)^8
             -\left(\frac9{10}\right)^2
          =\frac{22437}{12800}>0.
\]
Both bounds use that $(1-Y)(1+u)^8-Y^2$ decreases strictly with
$Y>0$. The intermediate value theorem therefore gives one root in
each of $(0,1/20)$, $(1/20,1/2)$ and $(1/2,\varphi)$.
These supply three distinct solutions with $x<0$ and $0<y<1$.

Next let $0<x<1$ and set
$V(x)=(x^2/(1-x))^{1/16}>0$. With $y=-V(x)$, the first
equation in \eqref{qt:degree:five-equations} again holds, and the
second becomes
\[
 G(x)=(1+V(x))(1-x)^8-V(x)^2=0.
\]
Here $G(0+)=1$. At $x=1/2$ we have
$1/2<V(1/2)=2^{-1/16}<1$, whence
$G(1/2)<2/256-1/4=-31/128<0$. This gives a fourth solution,
with $0<x<1/2$ and $y<0$.

Finally take $y=V(x)$ and put $\rho=(\sqrt5-1)/2$.
On $0<x<\rho$, both coordinates are in $(0,1)$ and the second
equation becomes
$J(x)=(1-V(x))(1-x)^8-V(x)^2=0$.
The endpoint limits are $J(0+)=1$ and $J(\rho-)=-1$.
This gives the fifth solution. The three sign regions are disjoint,
so all five solutions are distinct.
\end{proof}

\begin{remark}\label{qt:degree:real-scope}
Uniqueness of the distinguished point concerns $(0,1)^2$.
Proposition~\ref{qt:degree:five-real} shows why it cannot be replaced
by a four-solution bound on all real sign regions, even when the matrix
has integer entries. In that example
$1-x=(x/y^8)^2$ and $1-y=(y^{65}/x^8)^2$, so every real
solution in the domain $xy\ne0$ has $x,y<1$; nevertheless three
solutions lie in the same sign region $x<0$, $0<y<1$.
If $K=\Q(x,y)$ is the field of a Nahm point, its real embeddings
give distinct real solution pairs of the algebraic identities
satisfied by that point: equality on both generators would make two
embeddings identical. Thus a totally real quintic point would give
five real solution pairs, but rank two alone does not exclude that
count. The proposition does not assert that its five solutions form
a single quintic Galois orbit, or that the displayed matrix gives
modular data.
\end{remark}

%% file: sections/small_certificates.tex
\section{Small exact certificates}\label{scx:section}
This appendix records finite checks that can be reconstructed directly
from short formulas. The first replaces the expanded B\'ezout
coefficients for the four bilinear curves of
Section~\ref{d24:subsec:bilinear-polynomial}. The second specifies the
residue graph used in the global entry argument. The final table
separates the parameter-ideal resolutions by their three finite
classes.

\subsection{Four bilinear curves}\label{smallc:bilinear}
\begin{lemma}\label{smallc:lem:bilinear}
Assume $a,d>0$, $b\ne0$, $a,d\ne1$ and $\Delta=ad-b^2>0$.
For any of the four curves and exponent equations displayed in
Section~\ref{d24:subsec:bilinear-polynomial}, the necessary
condition $N\mid Q^2$ is impossible.
\end{lemma}
\begin{proof}
Coordinate exchange reduces the four cases to $G=X-Y-1$ and
$G=X(Y+1)-1$, both with $\Delta=d$. Thus
$a=1+b^2/d$. All calculations below take place in
$\Q[b,d][b^{-1},d^{-1},(d-1)^{-1}]$; these denominators are nonzero
on the stated domain. The cubic $Q$ and quadratic $N$ are those
displayed in the cited section.

For $G=X-Y-1$,
$N=-b\kappa^2+(d-a-1)\kappa+b$ has positive discriminant
$(d-a-1)^2+4b^2$. Therefore $N\mid Q^2$ implies $N\mid Q$.
Define
\begin{align*}
 P_0={}&b^5+b^3d^2+b^3d+b^2d^3+d^5-3d^4+2d^3,\\
 P_1={}&b^5+b^3d^2-b^3d+b^2d^3-2b^2d
             +d^5-5d^4+8d^3-4d^2.
\end{align*}
Polynomial division gives the exact identity
\[
 \left.\operatorname{rem}_{\kappa}(Q,N)\right|_{a=1+b^2/d}
   =\frac{P_0}{bd^2}+\frac{P_1}{b^2d}\kappa,
 \qquad
 \operatorname{Res}_d(P_0,P_1)=-64b^{15}.
\]
The resultant is the determinant of the Sylvester matrix of the two
displayed degree-five polynomials in $d$. If both remainders vanished,
that determinant would vanish, contradicting $b\ne0$.

For $G=X(Y+1)-1$, put
\begin{align*}
 P={}&b^3-d^3+2d^2,\\
 R={}&b^5-b^3d-b^2d^3+4b^2d^2-2b^2d-d^4+4d^3-4d^2.
\end{align*}
Here $N=b(d-1)\kappa^2+(a(d-1)+b^2-1)\kappa+ab$ and
\[
 \left.\operatorname{rem}_{\kappa}(Q,N)\right|_{a=1+b^2/d}
 =-\frac{(b^2+d)P}{bd^2(d-1)}
  -\frac{R}{b^2d(d-1)}\kappa.
\]
If $N$ is squarefree, the necessary divisibility forces $P=R=0$,
since $b^2+d>0$. But
\[
 R=(b^2+d-2)P-2b^2(b-d)(d-1).
\]
It follows that $b=d$, and then $P=2d^2\ne0$.

It remains to include tangency for this second curve. Define
\begin{align*}
 U={}&(b^2-2bd-d^2+2d)(b^2+2bd-d^2+2d),\\
 V={}&(b^2+d^2-2d)\bigl(2b^5+2b^4d-b^4-2b^3d^2
       -2b^2d^3+6b^2d^2\\
 &\hspace{10em}{}-4b^2d-d^4+4d^3-4d^2\bigr).
\end{align*}
If $N=n_2\kappa^2+n_1\kappa+n_0$, exact substitution gives
\[
 \left.\operatorname{disc}_{\kappa}N\right|_{a=1+b^2/d}
     =\frac{U}{d^2},\qquad
 \left.Q\left(-\frac{n_1}{2n_2}\right)\right|_{a=1+b^2/d}
     =\frac{V}{8b^3d^2(d-1)^2}.
\]
Thus a double root common to $N$ and $Q$ would give $U=V=0$.
The short elimination certificate is
\[
 \operatorname{Res}_d(U,V)
    =-65536b^{20}(b-1)^4(b+1)^4.
\]
It leaves $b=1$ or $b=-1$. At both values,
\[
 U=(d-1)(d+1)(d^2-4d-1),
\]
whereas
\[
 V|_{b=1}=-(d-1)^5(d+1),\qquad
 V|_{b=-1}=-(d-3)(d-1)^3(d+1)^2.
\]
Their common roots are only $d=1,-1$: in the second case $d=3$
is not a root of $U$. Both values are excluded by $d>0$, $d\ne1$.
This proves the tangent exclusion as well. Exchanging $a,d$ and
$X,Y$ supplies the other two curves.
\end{proof}

Every polynomial in this proof is displayed. The certificates require
only rational polynomial division, two Sylvester determinants and
substitution at $b=\pm1$; no stored Gr\"obner basis or expanded
B\'ezout coefficient is required.

\subsection{The integral residue graph}\label{smallc:residue-graph}
The following finite specification proves the residue counts used in
Section~\ref{deg:sec:entry-completion}. Work in $\F_7$, with vertex
set $V=\F_7^3$ and $\Delta=ad-b^2$. From $(a,b,d)$ draw the edges
\begin{align*}
 S(a,b,d)&=(d,b,a),\\
 J(a,b,d)&=(1-a,-b,1-d),\\
 I(a,b,d)&=(d/\Delta,-b/\Delta,a/\Delta)
                       &&(\Delta\ne0),\\
 T(a,b,d)&=(1-a^{-1},-b/a,\Delta/a)&&(a\ne0).
\end{align*}
These are the reductions of the legal rational maps in that section.
The inverse of a $T$ edge is the word $JTJ$, with legal pivots;
all graph edges are therefore reversible.

Write $q=c_3z^3+c_2z^2+c_1z+c_0$, where
\[
 c_3=d(d-1),\quad c_2=b(b+2d-1),\quad
 c_1=b(2a-1+b),\quad c_0=a(a-1).
\]
The discriminant used in the regularity test is the polynomial
\begin{equation}\label{smallc:formal-discriminant}
 \mathcal D=c_2^2c_1^2-4c_3c_1^3-4c_2^3c_0
             -27c_3^2c_0^2+18c_3c_2c_1c_0.
\end{equation}
In particular, the formal cubic degree is retained before
specialization; a discriminant routine that first drops leading zero
coefficients is not the same test. Mark a vertex regular precisely when
\[
 b(a-1)(d-1)(\Delta-a)(\Delta-d)(\Delta-a-d+1)\mathcal D\ne0.
\]
Start with these marked vertices and repeatedly add every endpoint
of the four legal edges until the set stops increasing. This finite
iteration has at most $343$ additions and gives the covered set.

To classify its complement, first take the vertices with $b=0$;
then those with $b\ne0$ and $q=0$. For each remaining vertex, divide
$q$ repeatedly by its nonconstant gcd with
\[
 z(bz+a)(bz+a-1)(dz+b)((d-1)z+b)
\]
until this gcd is constant. Classify the resulting polynomial as
repeated precisely when its gcd with its derivative has positive
degree. Constants, including an empty admissible support, belong to
the simple class. The complete counts are
\begin{center}
\begin{tabular}{@{}lr@{}}\toprule
Set & Cardinality\\\midrule
Regular vertices & $100$\\
Vertices reaching a regular vertex & $194$\\
Uncovered: simple admissible roots & $76$\\
Uncovered: repeated admissible roots & $22$\\
\quad of admissible polynomial degree two & $16$\\
\quad of admissible polynomial degree three & $6$\\
Uncovered: $b=0$ & $49$\\
Uncovered: $b\ne0$, $q=0$ & $2$\\\bottomrule
\end{tabular}
\end{center}
The two zero-cubic vertices are $(0,1,0)$ and $(1,6,1)$.
The acceptance test is solely membership and exact polynomial gcd over
$\F_7$. Iterating the displayed formulas over all $7^3$ triples proves
both the partition and its exhaustiveness. All counts concern matrix
residues; the subsequent local-length and norm arguments are needed
for the degree bounds.

\subsection{Parameter-pattern resolution counts}
\label{smallc:parameter-counts}
Table~\ref{smallc:tab:patterns} refines the total of $800$ patterns in
Section~\ref{pix:section}. Its columns are the special-root stratum
with $a,d\ne1$, the stratum $a=1$ treated up to exchange, and the
stratum $a=d=1$. The finite domain and the algebraic tests are given
there; the entries count the outcomes of those tests, not new
mathematical hypotheses.
\begin{longtable}{@{}P{0.64\textwidth}rrr@{}}
\caption{Resolution of the three parameter-pattern classes.}
\label{smallc:tab:patterns}\\
\toprule
Resolution & Special & $a=1$ & $a=d=1$\\\midrule\endfirsthead
\toprule
Resolution & Special & $a=1$ & $a=d=1$\\\midrule\endhead
Unit ideal & 430 & 295 & 12\\
Zero-dimensional minimal primes of degree greater than one & 10 & 5 & 4\\
Divisibility unit ideals in every projective kernel chart & 0 & 8 & 5\\
One-sign coefficients & 2 & 4 & 4\\
Rational sum-of-two-squares obstruction & 2 & 4 & 0\\
Bilinear exponent/conic exclusion (Lemma~\ref{smallc:lem:bilinear}) & 4 & 0 & 0\\
Positive sum of squares & 0 & 0 & 2\\
One rational parametrization and three resultant slices & 1 & 0 & 0\\
Positive definiteness and a square-sign contradiction & 1 & 0 & 0\\
Bottom exponent in $(0,1)$ and negative endpoint discriminant & 0 & 1 & 0\\
Factored zero-dimensional primes of degrees $4,8,2$ & 0 & 1 & 0\\
Three-factor split with a unit ideal on each factor & 0 & 1 & 0\\
Rational points force $|b|=2$ & 0 & 0 & 1\\
Two components with negative quadratic discriminant & 0 & 0 & 1\\
Two rational parametrizations and three resultant slices & 0 & 0 & 1\\
Discriminant $12(2C-1)^2$; exceptional point has $b=0$ & 0 & 0 & 1\\\midrule
Total & 450 & 319 & 31\\\bottomrule
\end{longtable}

%% file: sections/global_certificates.tex
\section{Finite computations for the initial degree bound}
\label{d24:certificates}
Theorem~\ref{d24:thm:main} uses polynomial identities, a finite
classification of common-component patterns, and an intersection count.
Section~\ref{pix:section} specifies the parameter ideals from their
constructors and gives the arithmetic resolutions of the exceptional
ideals. Appendix~\ref{scx:section} prints compact certificates for the
bilinear cases. The computation has the following exact domains.

\begin{longtable}{@{}P{0.25\textwidth}P{0.66\textwidth}@{}}
\toprule Claim & Domain and acceptance test\\\midrule\endhead
Formal coefficients & Generate the coefficients through order four by
Section~\ref{deg:sec:coefficients}. Compare rational functions over
$\Q(a,b,d,B,C,X,Y)$ after clearing the stated denominators.\\
Newton supports and boundary orders & Enumerate the coefficients of the
finite polynomials $\Psi_j$. Check the displayed support, edge and
conic-restriction identities of Section~\ref{d24:section}.\\
Ordinary bidegree $(2,2)$ patterns & Use the nine grid conditions in
Proposition~\ref{d24:prop:ord22}, the five isolated-point resultants
and the two family conditions. Test these exact identities on each
of the prescribed strata.\\
Coupled common-component patterns & Enumerate the six-line root
patterns of Section~\ref{d24:sec:m1}; the three strata have $450$,
$319$ and $31$ patterns. The ideal constructors, saturations and
projective coefficient charts are given in Section~\ref{pix:section}.\\
Bilinear survivors & Apply the displayed remainder and resultant
identities in Appendix~\ref{scx:section}, treating distinct and tangent
conic intersections separately. Coordinate exchange supplies the
symmetric cases.\\\bottomrule
\end{longtable}

The $800$ parameter patterns have the following resolutions. A unit
ideal is accepted only when exact polynomial ideal membership gives
$1$. A positive-dimensional or nonrational residual ideal is treated
by its stated arithmetic argument, rather than counted as an empty
ideal. A failed computation is never an exclusion.

\begin{center}
\begin{tabular}{@{}P{0.80\textwidth}r@{}}
\toprule Resolution & Patterns\\\midrule
Unit ideal & 737\\
Zero-dimensional minimal primes with no rational point & 19\\
Divisibility unit ideals in all projective kernel charts & 13\\
One-sign coefficient obstruction & 10\\
Rational sum-of-two-squares obstruction & 6\\
Bilinear exponent/conic exclusion & 4\\
Positive sum of squares & 2\\
Other rationality, positivity and factor-splitting resolutions & 9\\\midrule
Total & 800\\\bottomrule
\end{tabular}
\end{center}

The decoupled case is proved directly in
Proposition~\ref{d24:prop:decoupled}. It is not included in the $800$
coupled patterns. The initial-degree proof does not require the local
analysis at the prime 7 or the later arithmetic classification.

\section{Finite computations for the local degree reduction}
\label{deg:sec:supplement}
Table~\ref{deg:tab:terminal} lists the $29$ terminal domains.
The roadmap in Part~\ref{part:degree} explains the exhaustive case
division. The following table identifies the mathematical statements
served by the main finite computations; execution information belongs
to the companion readme.

\begin{longtable}{@{}P{0.25\textwidth}P{0.66\textwidth}@{}}
\toprule Argument & Finite computation and its role\\\midrule\endhead
Coordinate transport and entry &
Section~\ref{deg:sec:entry-completion} gives the rational maps, signed
equations, unit hypotheses and valuation partition. The integral
matrix reduction has $7^3=343$ residues. Its finite graph is used only
after the uniform denominator-normalization argument.\\
Canonical charts &
Sections~\ref{deg:sec:canonical} and~\ref{canx:section} specify the
unit and nonunit canonical branches, shift poles, reciprocal charts
and the separate two-unit norm tests.\\
Simple and repeated determinant centres &
Section~\ref{cuspx:section} gives the local algebra, the integral and
pole scalar cases, the repeated-centre jets and the combined norm
calculations. Its finite residues are used with bounds valid at every
depth, including ramified extensions.\\
Finite zero-coupling fibres &
The displayed centre and parity calculations, and the mod-nine
symplectic obstruction in Section~\ref{deg:sec:finite-fibres}, exclude the exceptional finite fibres. The unit
and centre contributions are combined before applying global norms.\\
Positive-dimensional residue fibres &
Sections~\ref{deg:sec:curve-bounds} and~\ref{cuspx:section} supply the
normal equations, exceptional exact families and lifting arguments.
A finite enumeration identifies the relevant families; the lifting
arguments bound all their depths.\\
Regular and exceptional shift poles &
Section~\ref{spx:section} gives weighted coefficient estimates, leading
normal equations and the global assembly. The fifth logarithmic
coefficient is needed for the regular and zero-cubic shift-pole rows.\\
Zero-cubic integral shifts &
Section~\ref{deg:sec:z0-replacement} and the integral-shift argument in
Section~\ref{cuspx:section} combine first-order equations, the
restricted zero-cubic calculation and the complete unit-plus-conic
norm test.\\
Final degree count &
The proof following Table~\ref{deg:tab:terminal} combines local lengths
and residue degrees over all places, using the units available in
that precise chart. Every resulting whole-field cap is at most four.\\
\bottomrule
\end{longtable}

The two roles of finite calculation must be distinguished. Polynomial
identities and finite-field ranks certify exact algebraic assertions.
A finite local search certifies only the listed local configurations;
its use at arbitrary depth requires the accompanying valuation,
implicit-function and local-length argument. Neither a count of
files nor a successful program status is a mathematical hypothesis.

The companion \cite{Comp} contains the reproducible programs and raw
output for these finite domains. Neutral descriptions in its source
map connect each program to the corresponding theorem above. The
small source archive suffices to compile this article without running
any mathematical program; reproduction of the finite computations
uses the separate computational archive.

%% file: sections/rational_boundary.tex
\section{Rational two-term boundary relations}\label{low:rationalapp}
This appendix is a weighted boundary statement independent of the symmetric matrix hypothesis.
\begin{theorem}\label{low:rationaltheorem}
Let $x,y\in\Q\cap(0,1)$, $\dd(x)\ne0$, and $m\in\Z_{>0}$. If $m\dd(x)+\dd(y)=0$, then for $m=1$ one has $y=1-x$. For $m>1$, up to simultaneous reflection, exactly the following triples occur:
\begin{equation}\label{low:rationaltable}
\begin{array}{c|c|c|c}
m&x&y&mL(x)+L(y)\\\hline
2&1/3&1/4&C_0\\
3&1/4&1/9&C_0\\
6&1/3&8/9&3C_0
\end{array}
\end{equation}
For zero boundary the only rational argument is $1/2$, so $x=y=1/2$ and $m$ is arbitrary.
\end{theorem}
\begin{proof}
Write a rational argument as $a/(a+b)$ with $a,b$ coprime positive integers. Its nonzero valuation vectors for the complementary pair lie on the three rays generated by $(1,0),(0,1),(-1,-1)$, according as the prime divides $a,b,a+b$. The boundary is zero exactly when $a=b=1$, giving $1/2$. Otherwise at least two rays occur.

The nonzero wedge relation identifies the two rational multiplicative planes. Its rational exponent matrix is invertible, has determinant $-m$, and maps the occurring valuation rays injectively to the same three-ray set. If all three rays occur, the matrix is a positive scalar multiple of one of the six matrices induced by anharmonic transformations. Indeed the three ray generators sum to zero, so their three scale factors must agree. Undoing that transformation, the absolute values are $t^c,(1-t)^c$ for $t\in(0,1)$ and $c>0$. To see that $t$ is in this interval, both absolute values must be below one and the transformed signed pair sums to one; hence both entries are positive. The equality $t^c+(1-t)^c=1$ forces $c=1$ by strict monotonicity in $c$. Thus $m=1$ and the opposite boundary representative is reflection.

If only two rays occur, each argument has an anharmonic transform which is an integer $N\ge3$. The two coprime integers $N,N-1$ correspond to the two rays. Matching the two rays for the other argument gives common integer bases $h,k\ge2$ and two positive integer exponent pairs satisfying
\begin{equation}\label{low:rationalpowers}
|h^r-k^s|=1.
\end{equation}
The bases exist because two multiplicatively dependent positive integers are positive integer powers of a common integer, as follows by taking the gcd of their prime-exponent vectors. If the exponent pairs agree, the arguments have the same anharmonic orbit and again give reflection with $m=1$.

For distinct pairs, if either solution has both exponents greater than one, Mih\u{a}ilescu's theorem \cite{Mihailescu} gives $9-8=1$, so the bases are $3,2$. Otherwise arrange $h>k$. The larger-base exponent must be one in each solution: if the smaller-base exponent were one and the larger exceeded one, their difference would exceed one. Two distinct powers of $k$ then differ by two, forcing $k=2$ and $h=3$. For these bases the complete exponent pairs are
\[
(r,s)=(1,1),(1,2),(2,3),\qquad |3^r-2^s|=1,
\]
using Mih\u{a}ilescu when both exponents exceed one and elementary inspection when one is one. The positive argument orbits are therefore represented by $1/3,1/4,1/9$ and their complements. Their boundary coefficients in the basis $[2]\wedge[3]$ are $1,-2,6$; reflection negates these coefficients. Equal magnitudes give only reflection, while integer ratios greater than one give exactly \eqref{low:rationaltable}.

Finally, duplication at $1/2$ and $1/3$ gives
\[
2L(1/3)+L(1/4)=C_0,\qquad
2L(1/3)=L(1/9)+2L(1/4).
\]
Elimination and reflection produce the other two identities, proving all displayed evaluations.
\end{proof}

%% file: sections/verification.tex
\section{Exact computation and verification}\label{app:verification}
The computer-assisted steps have two separate inputs: a proved finite
domain or polynomial system, and an exact verifier for that domain or
system. A bounded numerical search is not used to establish an
exhaustive assertion. This appendix records the acceptance conditions
and the relation between the calculations and the theorems.

\subsection{Finite domains and exact acceptance}
The unit searches first prove coefficient-height bounds. The defining
polynomials are then enumerated from their endpoint norms, and
irreducibility and root distributions are checked exactly. Dyadic root
intervals are proposal data: polynomial sign changes and disjointness
certify their roots. Logarithms are bounded by explicit series with
outward integer or rational rounding. Inclusive interval comparisons
can retain extra candidates but cannot discard an exact equality.

The canonical nonunit searches first prove bounds on primitive
exponents and the remaining exponent. Finite-prime sieves reject a row
only by an exact polynomial obstruction. A rejection remains valid
after passing to the distinguished real embedding. Cases that survive
the sieves are factored exactly, with the interval root, field degree
and unit conditions checked separately. The structural recounts and
the modular rejection traces are included in the companion~\cite{Comp}.

At a fixed matrix and Nahm point, the shift equations are polynomial
equations over an explicitly specified number field. All rational
shifts are obtained by exact elimination. For exclusions, the
certificate is either an explicit identity
\[
 1=\sum_i h_iE_i
\]
or a univariate eliminant together with exact specialized gcds. In the
positive-dimensional cases, factors and components are separated
before elimination; every removed locus is either solved separately or
proved inadmissible. Every surviving shift is substituted into the
original equations, and its scalar is evaluated as
$c_0=\lambda_1-1/12$.

\subsection{Independent elimination for the low-degree families}
The systems in the proof of Theorem~\ref{cls:thm:deg2} were regenerated
from the universal coefficient polynomials. The original primary
decompositions are supplied, together with an independent exact
elimination that does not require \textsc{Singular} or \textsc{Pari}.

For the symmetric golden families put
$s=B_1+B_2$ and $u=(B_1-B_2)^2$. The first-coefficient equation is
linear in $u$. Its coefficient is
\[
 5(4b^2-10b+5)\quad\text{for }(r,s)=(1,1),\qquad
 -5(16b^2-5)\quad\text{for }(r,s)=(2,2).
\]
Neither coefficient has a rational zero. Solve for $u$, eliminate it
from the two second-coefficient equations, and take the resultant in
$s$. Rational-root computation gives respectively
\[
 b\in\{0,\tfrac12,\tfrac23,1,2\},\qquad
 b\in\{-\tfrac12,-\tfrac14,0,\tfrac14,\tfrac12\}.
\]
For each value, the gcd of all seven specialized equations gives all
rational $s$, and back-substitution gives $u$. There are exactly eleven
rational invariant points per family, seven in its positive-definite
range. The admissible values of $u$ have rational square roots, so the
original vector shifts are recovered. The complete resultants,
specialized gcds and substitutions are in
\path{verification/low_degree_family_ideals.json} in the companion~\cite{Comp}.

For the unequal golden families, after removing the inadmissible locus
$b=1$, the exact Gr\"obner bases are
\[
 [b,B_1^2-B_1,B_2^2-B_2/2],\qquad
 [b,B_1^2-B_1/2,B_2^2-B_2].
\]
Their quotient dimensions are four, and all four distinct rational
points in each case are substituted exactly. For the rational point,
the component $B_1+B_2=0$ is verified identically. On the remaining
locus the rational values of $s$ are $-1,-1/2,1/2,1$; the extra factor
$45s^2+68s+4$ has no rational root. The specialized ideals are
\[
\begin{array}{c|c}
s&\text{equations in }b,u\\\hline
-1&u=0\\
-1/2&b=0,\ u=1/4\\
1/2&b(2b-1)=0,\ u=1/4\\
1&u=b^2.
\end{array}
\]
The points at $b=1/2$ are singular. This accounts for every rational
family and isolated point in Theorem~\ref{cls:thm:deg2}.

The finite quadratic calculation checks all twenty nonrational
configurations using exact quadratic quotient arithmetic. Its first
two coefficients are constructed by a direct Gaussian calculation and
compared coefficient by coefficient both with the radial recurrence
and with the universal coefficient input. The equations
$\lambda_1\in\Q$, $\lambda_2=0$ already produce all $74$ oriented
rational triples in these configurations; the subsequent
$\lambda_3,\lambda_4$ checks reject none. The number includes coordinate
exchanges and is not the number of exceptional classes in
Table~\ref{tab:complete-list}.

\subsection{Finite classifications and cross-checks}
\begin{center}\small
\begin{tabular}{P{0.31\textwidth}P{0.61\textwidth}}
\toprule
Computation&Exact outcome\\\midrule
Quadratic arithmetic&$446$ norm pairs, $26$ ideal pairs; all $42$
cyclotomic comparisons pass.\\
Cubic odd units&$13\,196$ partner and $656$ anchor polynomials;
only $656$ identical-polynomial area pairs remain, giving no independent pair.\\
Cubic even units&$36$ ordered coordinate pairs give eight ordered
matrices, four up to exchange, with twelve rational triples.\\
Cubic unequal nonunits&$5\,844$ triples reduce to $16$ exact
factorizations; all interval cubic factors are unit degenerations.\\
Cubic equal exponents&$918$ primitive exponent pairs and $16$ cubic
polynomials give $14\,688$ exact divisions; all eight resulting
matrix orientations and inverses are excluded.\\
Cubic dependent cases&Twelve exact line exclusions; nine direct
symbolic coefficient calculations agree in $54$ identities.\\
Quartic totally real units&Both index parities are excluded by the
complete interval and pairing verifiers.\\
Quartic mixed units&Five odd-index and three even-index matrix classes;
all reconstruction records agree with the shift inputs.\\
Quartic canonical nonunits&$9\,626\longrightarrow49\longrightarrow0$;
the structural recount and $1\,377$ circuit self-tests agree.\\
Quartic shifts&Eight points, sixteen orientations, four rational
triples up to coordinate exchange.\\
\bottomrule
\end{tabular}
\end{center}

The independent cubic checker uses generic rational Euclidean Sturm
sequences, reproduces the complete polynomial census and verifies
$39\,588$ proposed root intervals. A separate quotient-field
implementation reconstructs the even unit matrices and checks their
Nahm identities. In the dependent cubic calculation, direct symbolic
Gaussian moments construct the relevant polynomials without
interpolation. The scaled second-coefficient degree is six, within the
proved bound ten. For the mixed even quartic branch, a separate
rational Horner-interval calculation verifies all seven embedding
rejections and all five retained maps.

The cubic and quartic shift verifiers compare the two coefficient
implementations exactly. The inverse matrices are checked in their
own quotient fields. The coefficient transport in
\eqref{cls:eq:duality} is not used as an assumption that modularity is
preserved by inversion.

\subsection{Verification of the degree and arithmetic computations}
The computations for Theorem~\ref{deg:thm:main} are organized by their
mathematical role in Appendix~\ref{deg:sec:supplement}. The exact
identities for the bilinear cases are printed in
Appendix~\ref{scx:section}; they can be verified by polynomial
arithmetic without the larger companion.

For the arithmetic classification, the checks serve three distinct
claims: completeness of each bounded polynomial list, exclusion or
reconstruction of its possible coordinate pairs, and exact elimination
of the rational shifts. The count of thirty-four exceptional triples
is checked against the actual table in this article, including
coordinate exchange and the exact scalar shifts. Modularity of the
survivors follows from the cited identities, rather than these finite
necessary-condition tests.

For the totally real quartic pairing exclusions, the integer
implementations reconstruct every logarithmic interval from exact
root signs. Their certified component margins are stated beside the
odd- and even-index exclusions in the totally real quartic argument.
The even-index minimum is also given explicitly in
\eqref{qt:tr:eq:even-margin}. Execution receipts, source identities and
software versions are recorded separately in the companion readme.

\subsection{Certified interval margins}\label{verify:interval-margins}
For clarity, the table below records the least certified absolute
lower endpoint of the selected nonzero test quantity among the
indicated rejections. All entries are integers divided by $2^{96}$.
The coefficient-sign tests concern the Cramer-rule numerators,
before division by the already nonzero determinant. The choice of a
nonzero component may differ between configurations; a positive
minimum certifies every rejection in that row.

\begin{center}\small
\begin{tabular}{@{}P{0.39\textwidth}rP{0.36\textwidth}@{}}
\toprule Test & Rejections & Numerator of least bound\\\midrule
Totally real odd-index wedge & $3276$ & $9053621231263952158607117616$\\
Totally real even-index wedge & $2\,507\,880$ & $170480871405090970163$\\
Mixed even-index determinant sum & $7140$ & $5403606475833105850348272906$\\
Mixed even-index coefficient sign & $44$ & $22876094254218257524787161377$\\
Mixed odd-index coefficient sign & $13\,356$ & $10769108040732138045110846218$\\
\bottomrule
\end{tabular}
\end{center}

The mixed odd-index determinant test makes no additional rejection
among its $17\,941$ area-bucket pairs, so it has no positive exclusion
margin to report. The $4585$ surviving mixed odd configurations and
$12$ surviving mixed even configurations proceed to the exact field,
exponent and index tests stated in their respective sections. The
bounds in this table concern these interval stages only; exact
polynomial exclusions and the height estimates have their own
proofs. The programs use exact root signs and rational series
remainders throughout. Decimal approximations are not decision criteria.

\subsection{Reproduction and software}
The companion~\cite{Comp} is organized by proof branch and includes the
original inputs, exact verifiers, output certificates and a source manifest.
Unless explicitly identified as an ancillary file of the arXiv submission,
file and directory names in typewriter type refer to this companion after
running \path{restore_companion.py}. Paths are relative to its restored
root, and abbreviated file names are resolved in its
\path{CITED_FILE_INDEX.txt}. Its entry points are listed in its files
\path{README.txt} and \path{REPRODUCE.txt}. The ancillary files of the
arXiv submission contain the dependency map and the revision checks.
The paper's source archive compiles independently of executing any
certificate program. The separate computational archive supplies the
full data for the checks.

The complete execution records of the large computations are preserved
in the companion: the $800$ parameter-ideal calculations, the $114$ cusp
searches, the $146$ classifications, the branch-two census, the
$12{,}800$ fifth-coefficient normal jobs and the $456$-entry pole census.
Not all of these computations were rerun for this version. Fresh checks
include the complete quartic interval-pairing runs, exact substitution
of the thirty-four exceptional shifts, independent re-enumeration of
the $800$ root prescriptions, reconstruction of their polynomial inputs,
and the smaller exact checks described in the verification records.
Reconstructing all polynomial inputs is distinct from rerunning all
parameter-ideal calculations. The companion's \path{README.txt}
distinguishes preserved records from fresh reruns.

The fresh low-degree verification used Python~3.12.14,
SymPy~1.14.0 and python-flint~0.9.0; mpmath~1.3.0 is used where
explicitly indicated. Assertions are enabled. Program checks that
depend on rigorous intervals use exact bounds or Arb arithmetic;
unverified numerical approximations are not acceptance criteria.
The original and independent verification routes are identified in the
manifest. The mathematical theorems rely on their specified exact
identities and exhaustive domains, not on the choice of user interface
or the date of a run.

%% file: sections/disclosures.tex
\section*{Author contributions and use of AI tools}
The author proposed and directed the programme of reducing the
symmetric rank-two modular Nahm-sum problem to a bound on the degree of
the distinguished point, with the aim of reaching degree four and
classifying the remaining cases by finite arguments. AI systems,
principally Anthropic Claude and OpenAI ChatGPT/Codex, generated and
developed proof arguments for the initial degree-$24$ bound and the
subsequent reductions, the arithmetic and height arguments, the
degree-one through degree-four classifications, and the exact shift
calculations. The degree-bound proofs and substantial parts of the
classification proofs were AI-generated.

These systems were also used to write verification programs, perform
independently implemented checks, review proof dependencies and prepare
the manuscript. The computational supplement distinguishes replay of
existing programs from independent implementations and identifies
retained computation records. The author takes responsibility for the
mathematical claims and submitted materials.

\section*{Funding and competing interests}
No external funding or competing interests are declared.